\documentclass[10pt, a4paper, notitlepage]{article}
\pdfoutput=1

\usepackage{biblatex}
\usepackage{bbm}

\usepackage{etoolbox}
\usepackage{tikz}
\usetikzlibrary{calc}
\usetikzlibrary{cd}
\usetikzlibrary{decorations.markings}
\usetikzlibrary{decorations.pathreplacing}
\usetikzlibrary{decorations.pathmorphing}
\usetikzlibrary{decorations.text}
\usetikzlibrary{arrows.meta}
\usetikzlibrary{arrows}
\usetikzlibrary{positioning}

\usepackage{geometry} 
\usepackage{tocloft} 
\usepackage{fancyhdr} 
\usepackage{longtable}
\usepackage{subcaption}
\usepackage{enumitem}
\usepackage{amssymb}
\usepackage{amsmath}
\usepackage{stackrel} 
\usepackage{uniinput}
\usepackage{amsthm}
\usepackage{stmaryrd}
\usepackage{mathtools}
\usepackage{fouridx}
\usepackage{xparse}
\usepackage{bm}
\usepackage{aliascnt}
\usepackage{import}

\usepackage{hyperref}
\theoremstyle{definition}
\newtheorem{theorem}{Theorem}
\let\emph\textbf

\newaliascnt{remark}{theorem}
\newaliascnt{definition}{theorem}
\newaliascnt{proposition}{theorem}
\newaliascnt{lemma}{theorem}
\newaliascnt{corollary}{theorem}
\newaliascnt{example}{theorem}
\newaliascnt{convention}{theorem}
\newaliascnt{todo}{theorem}

\newtheorem{remark}[remark]{Remark}
\newtheorem{definition}[definition]{Definition}
\newtheorem{proposition}[proposition]{Proposition}
\newtheorem{lemma}[lemma]{Lemma}
\newtheorem{corollary}[corollary]{Corollary}
\newtheorem{example}[example]{Example}
\newtheorem{convention}[convention]{Convention}

\aliascntresetthe{remark}
\aliascntresetthe{definition}
\aliascntresetthe{proposition}
\aliascntresetthe{lemma}
\aliascntresetthe{corollary}
\aliascntresetthe{example}
\aliascntresetthe{convention}
\aliascntresetthe{todo}

\numberwithin{equation}{section}
\numberwithin{figure}{section}
\numberwithin{table}{section}
\makeatletter\let\c@table\c@figure\makeatother

\numberwithin{theorem}{section}
\numberwithin{remark}{section}
\numberwithin{definition}{section}
\numberwithin{proposition}{section}
\numberwithin{lemma}{section}
\numberwithin{corollary}{section}
\numberwithin{example}{section}
\numberwithin{convention}{section}
\numberwithin{todo}{section}

\DeclareFieldFormat*{url}{}
\DeclareFieldFormat*{doi}{}
\DeclareFieldFormat*{issn}{}
\DeclareFieldFormat[misc]{url}{\mkbibacro{URL}\addcolon\space\url{#1}}
\DeclareFieldFormat[online]{url}{\mkbibacro{URL}\addcolon\space\url{#1}}

\setlist{topsep=0.5em, itemsep=0em}
\renewcommand{\arraystretch}{1.3}

\SetLabelAlign{parleftalign}{\parbox[t]\textwidth{#1\par\mbox{}}}

\renewcommand{\k}{\mathbbm{k}}
\newcommand{\bbmu}{\mbox{$\raisebox{-0.59ex}
  {$l$}\hspace{-0.21em}\mu\hspace{-0.88em}\raisebox{-0.98ex}{\scalebox{2}
  {$\color{white}.$}}\hspace{-0.416em}\raisebox{+0.88ex}
  {$\color{white}.$}\hspace{0.46em}$}{}}

\newcommand{\Id}{\operatorname{Id}}
\newcommand{\id}{\mathrm{id}}
\newcommand{\coid}{\mathrm{id}^*}
\newcommand{\Hom}{\operatorname{Hom}}
\newcommand{\Ker}{\operatorname{Ker}}
\newcommand{\End}{\operatorname{End}}
\newcommand{\Ob}{\operatorname{Ob}}
\newcommand{\cat}[1]{\mathcal{#1}}

\newcommand{\subsum}[1]{\sum_{\substack{#1}}}
\newcommand{\htensor}{\widehat{\otimes}}
\newcommand{\pre}{\operatorname{pre}}

\newcommand{\Tw}{\operatorname{Tw}}
\newcommand{\Mod}{\operatorname{Mod}}
\newcommand{\Add}{\operatorname{Add}}
\newcommand{\MC}{\operatorname{MC}}
\newcommand{\MCb}{\overline{\operatorname{MC}}}
\newcommand{\Gtl}{\operatorname{Gtl}}
\newcommand{\Jac}{\operatorname{Jac}}
\newcommand{\Fuk}{\operatorname{Fuk}}
\newcommand{\wFuk}{\operatorname{wFuk}}
\newcommand{\preFuk}{\operatorname{Fuk}^\mathrm{pre}}
\newcommand{\relFuk}{\operatorname{relFuk}}
\newcommand{\relpreFuk}{\operatorname{relFuk}^{\mathrm{pre}}}

\newcommand{\MF}{\operatorname{MF}}
\def\H{\operatorname{H}}
\newcommand{\HTw}{\operatorname{H}\operatorname{Tw}}

\newcommand{\mf}{\operatorname{mf}}
\newcommand{\HH}{\operatorname{HH}}
\newcommand{\HC}{\operatorname{HC}}

\newcommand{\pmat}[1]{\begin{pmatrix} #1 \end{pmatrix}}
\newcommand{\vspan}{\operatorname{span}}

\newcommand{\isoto}{\xrightarrow{\sim}}
\newcommand{\verylongisoto}{\xrightarrow{~~~\sim~~~}}
\newcommand{\verylongmapsto}{\xmapsto{~~~\phantom{\sim}~~~}}
\newcommand{\verylongto}{\xrightarrow{~~~\phantom{\sim}~~~}}
\newcommand{\embeds}{\hookrightarrow}
\def\Im{\operatorname{Im}}
\newcommand{\smooth}[1]{\tilde{#1}}
\newcommand{\restr}[1]{|_{#1}}
\newcommand{\Subdisk}{\mathsf{D}}
\newcommand{\Subtree}{\mathsf{T}}
\newcommand{\Subresult}{\mathsf{R}}
\newcommand{\SLd}{\mathsf{Disk}_{\operatorname{SL}}}
\newcommand{\IDd}{\mathsf{Disk}_{\operatorname{ID}}}
\newcommand{\IDr}{\mathsf{Result}_{\operatorname{ID}}}
\newcommand{\PiTr}{\mathsf{Result}_{\operatorname{π}}}
\newcommand{\CRr}{\mathsf{Result}_{\operatorname{CR}}}
\newcommand{\CRd}{\mathsf{Disk}_{\operatorname{CR}}}
\newcommand{\DSr}{\mathsf{Result}_{\operatorname{DS}}}
\newcommand{\DSd}{\mathsf{Disk}_{\operatorname{DS}}}
\newcommand{\DWr}{\mathsf{Result}_{\operatorname{DW}}}
\newcommand{\DWd}{\mathsf{Disk}_{\operatorname{DW}}}
\newcommand{\PiTrM}{\mathsf{Result}_{\operatorname{πM}}}
\newcommand{\MDd}{\mathsf{Disk}_{\operatorname{MD}}}
\newcommand{\MTd}{\mathsf{Disk}_{\operatorname{MT}}}
\newcommand{\MDr}{\mathsf{Result}_{\operatorname{MD}}}
\newcommand{\MTr}{\mathsf{Result}_{\operatorname{MD}}}
\newcommand{\disktarget}{\operatorname{t}}
\newcommand{\Abouzaid}{\operatorname{Abou}}
\newcommand{\punctures}{\operatorname{Punc}}
\newcommand{\sgn}{\operatorname{sgn}}
\newcommand{\qparam}{\operatorname{qparam}}
\newcommand{\disjoint}{\mathbin{\dot{\cup}}}
\newcommand{\ZigzagCat}{\mathbb{L}}
\newcommand{\DefZigzagCat}{\mathbb{L}_q}

\newcommand{\running}{~|~}
\newcommand{\mirQ}{\check{Q}}

\newcommand{\rectified}{\mathsf{R}}
\newcommand{\transversal}{\mathsf{tr}}

\newcommand{\Res}{\operatorname{Res}}

\newcommand{\Perf}{\operatorname{Perf}}
\newcommand{\Stab}{\operatorname{Stab}}
\newcommand{\Def}{\operatorname{Def}}
\newcommand{\Coh}{\operatorname{Coh}}

\newcommand{\sphereodd}{{\mathsf{odd}}}
\newcommand{\sphereeven}{{\mathsf{even}}}
\newcommand{\landau}{\mathcal{O}}

\newcommand{\coeff}[1]{|#1|_B}

\newcommand{\chlQ}{Q^{\ZigzagCat}}
\newcommand{\RefObjects}{\mathbb{L}}
\newcommand{\DefRefObjects}{\mathbb{L}_q}

\NewDocumentCommand{\tensor}{t_}
 {%
  \IfBooleanTF{#1}
   {\tensorop}
   {\otimes}%
 }
\NewDocumentCommand{\tensorop}{m}
 {%
  \mathbin{\mathop{\otimes}\displaylimits_{#1}}%
 }

\newcommand{\paperone}{\cite{Paper-I}}
\newcommand{\papertwoA}{\cite{Paper-IIA}}
\newcommand{\paperthree}{\cite{Paper-III}}
\newcommand{\cA}{\mathcal{A}}

\title{Relative Fukaya Categories via Gentle Algebras}
\author{Jasper van de Kreeke}
\date{15 March 2023}

\begin{document}

\maketitle

\begin{abstract}
Relative Fukaya categories are hard to construct. In this paper, we provide a very explicit construction in the case of punctured surfaces.

The starting point is the gentle algebra $ \Gtl Q $ associated with a punctured surface $ Q $. Our model for the relative Fukaya category is a deformation $ \Gtl_q Q $ which has one formal deformation parameter per puncture. We verify that $ \Gtl_q Q $ is a model for the relative Fukaya category by computing a large part of the $ A_∞ $-structure on the derived category $ \H\Tw\Gtl_q Q $. The core method is the deformed Kadeishvili theorem from \papertwoA.

The paper's technical contributions include (1) a construction for removing curvature from band objects, (2) a method for analyzing Kadeishvili trees, (3) a matching between deformed Kadeishvili trees and relative Fukaya disks. This paper is the second in a series of three, whose aim it is to deform mirror symmetry for punctured surfaces.
\end{abstract}

\newcommand{\papertwoacknowledgements}{%
\subsection*{Acknowledgements}
This paper is part of the author's PhD thesis supervised by Raf Bocklandt. Crucial insights were gained at the 2020 Summer Camp on Derived Categories, Stability Conditions and Deformations, which the author organized in his garden in Berlin. The author thanks all participants for sharing their knowledge, in particular Severin Barmeier for advice and expertise on deformation theory through the lens of $ L_\infty $-algebras. The author's PhD project was supported by the NWO grant “Algebraic methods and structures in the theory of Frobenius manifolds and their applications” (TOP1.17.012).}

\tableofcontents

\section{Introduction}
Fukuya categories capture the global geometry of manifolds. They are complicated even to define and hard to study. In the case of punctured surfaces, gentle algebras were introduced as a remedy by Bocklandt \cite{Bocklandt}. They have become a successful discrete version of the Fukaya category, having already served as A-side in mirror symmetry \cite{Bocklandt} and as standard model to study homological properties and stability conditions \cite{HKK}.
\begin{center}
\begin{tikzpicture}
\path (0, 0) node[align=center] {\textbf{Smooth} \\ Fukaya category $ \Fuk Q $} (4, 0) node {\Large $ \longleftrightarrow $} node[above] {\cite{Bocklandt}} (8, 0) node[align=center] {\textbf{Discrete} \\ Gentle algebra $ \Gtl Q $};
\end{tikzpicture}
\end{center}
On the smooth side, Seidel introduced in 2002 a procedure to deform Fukaya categories \cite{Seidel-relative}. The result is now known as the “relative Fukaya category” and has already been used e.g.~as A-side in deformed mirror symmetry for the $ n $-punctured torus by Lekili-Perutz-Polishchuk \cite{Lekili-Perutz, Lekili-Polishchuk}. Surprisingly, a rigorous construction of the relative Fukaya category was only finished in 2022 by Perutz and Sheridan \cite{Sheridan-Perutz}.

On the discrete side, Bocklandt and the author recently proposed an analog of Seidel's procedure in \paperone. Our analog consists of an explicit deformation of the gentle algebra. Ultimately, these “deformed gentle algebras” will serve as A-side in our proof of deformed mirror symmetry for arbitrary punctured surfaces.

In this paper, we prove that the smooth and discrete deformations are equivalent. More precisely, we show that on the subcategories of zigzag curves, the smooth and discrete deformations have the same $ A_∞ $-structure:
\begin{center}
\begin{tikzpicture}
\path (0, 0) node[align=center] {\textbf{Deformed smooth} \\ Relative Fukaya category $ \relFuk Q $} (4, 0) node {\Large $ \longleftrightarrow $} (8, 0) node[align=center] {\textbf{Deformed Discrete} \\ Deformed gentle algebra $ \Gtl_q Q $};
\end{tikzpicture}
\end{center}

\subsection*{Assembling deformed mirror symmetry}
This paper is the second part in a series of three papers aimed at proving deformed noncommutative mirror symmetry for punctured surfaces. The results of this paper seem to be interesting enough to stand on their own, but their full value becomes visible when viewed in the context of the series. Here we recall the overall aim of the series and the special relevance that this paper has to the workings of the series.

\paragraph*{Mirror symmetry for punctured surfaces} Original mirror symmetry of punctured surfaces due to Bocklandt \cite{Bocklandt} considers as A-side the gentle algebra $ \Gtl Q $ of a dimer $ Q $ and as B-side a category of matrix factorizations $ \mf(\Jac \mirQ, ℓ) $ of the dual dimer $ \mirQ $. Under the assumption that $ \mirQ $ is zigzag consistent, Bocklandt proves the existence of an $ A_∞ $-quasi-isomorphism
\begin{equation*}
\Gtl Q ≅ \mf(\Jac\mirQ, ℓ).
\end{equation*}
This is known as noncommutative mirror symmetry for punctured surfaces. A natural question is which deformation $ \Gtl Q $ corresponds to which deformation of $ \mf(\Jac\mirQ ℓ) $. We shall focus on one broad deformation $ \Gtl_q Q $ from the class of deformations defined in \paperone\ and ask which deformation $ \mf_q (\Jac\mirQ, ℓ) $ of $ \mf(\Jac\mirQ, ℓ) $ corresponds to $ \Gtl_q Q $ such that there is still a quasi-isomorphism of deformed $ A_∞ $-categories $ \Gtl_q Q ≅ \mf_q (\Jac\mirQ, ℓ) $.

\paragraph*{The Cho-Hong-Lau construction}
Proving mirror symmetry is inherently difficult because two categories only vaguely resembling each other need to be matched. More precisely, there exists typically no strict $ A_∞ $-isomorphism between the categories involved. To construct a non-strict functor as in \cite{Bocklandt} requires one to recognize that two given $ A_∞ $-structures are equal up to a kind of homotopy. The analogous question in case of deformations of $ A_∞ $-categories is how to decide whether two given $ A_∞ $-deformations are gauge-equivalent. There are apparently very few tools available to decide this question.

The game changes as soon as we take the work of Cho, Hong and Lau \cite{CHL} into account. They explain how to construct a mirror equivalence for punctured surfaces by a version of Koszul duality. Their paper shows how to systematically obtain both the dual dimer $ \mirQ $ and Bocklandt's mirror equivalence from a systematic construction:

\begin{center}
\begin{tikzpicture}
\path (0, 0) node[align=center] (A) {$ \ZigzagCat ⊂ \cat C $ \\ $ A_∞ $-category with subcategory} (9, 0) node[align=center] (B) {$ F: \cat C → \MF(\Jac(\chlQ, W), ℓ) $ \\ Mirror functor};
\path[draw, ->, decorate, decoration={snake, amplitude=0.2em, post length=0.5em}] ($ (A.east) + (right:1) $) to ($ (B.west) + (left:1) $);
\end{tikzpicture}
\end{center}

The mirror category $ \MF(\Jac(\chlQ, W), ℓ) $ is a category of matrix factorizations. The Jacobi algebra $ \Jac(\chlQ, W) = ℂ\chlQ / (∂_a W) $ and the potential $ ℓ ∈ \Jac(\chlQ, W) $ are determined by a kind of Koszul transform of the $ A_∞ $-structure on the subcategory $ \ZigzagCat ⊂ \cat C $. Mirror symmetry of punctured surfaces is a special case of the Cho-Hong-Lau construction: We set $ \cat C = \HTw\Gtl Q $ and let $ \H\ZigzagCat ⊂ \HTw\Gtl Q $ be the subcategory of zigzag paths (\cite[Chapter 10]{CHL}). This way, the Jacobi algebra $ \Jac(\chlQ, W) $ becomes the Jacobi algebra $ \Jac \mirQ $ of the dual dimer and $ ℓ $ becomes the standard central element $ ℓ ∈ \Jac \mirQ $. The mirror functor $ F: \Gtl Q → \MF(\Jac \mirQ, ℓ) $ one obtains this way is an explicit incarnation of Bocklandt's mirror symmetry for punctured surfaces.

\paragraph*{Deformed Cho-Hong-Lau construction}
We provide a deformed version of the Cho-Hong-Lau construction. The starting point for the deformed construction is a given deformation $ \cat C_q $ of $ \cat C $. We denote by $ \DefRefObjects ⊂ \cat C_q $ the subcategory with the same objects as $ \RefObjects $. The result of our deformed Cho-Hong-Lau construction is a functor of $ A_∞ $-deformations $ F_q: \cat C_q → \MF(\Jac(\chlQ, W_q), ℓ_q) $. The deformed superpotential $ W_q $ and potential $ ℓ_q $ are computed from the deformed $ A_∞ $-structure on $ \DefRefObjects $.

\begin{center}
\begin{tikzpicture}
\path (0, 0) node[align=center] (A) {$ \DefRefObjects ⊂ \cat C_q $ \\ Deformed category with subcategory} (9, 0) node[align=center] (B) {$ F_q: \cat C_q → \MF(\Jac(\chlQ, W_q), ℓ_q) $ \\ Deformed mirror functor};
\path[draw, ->, decorate, decoration={snake, amplitude=0.2em, post length=0.5em}] ($ (A.east) + (right:1) $) to ($ (B.west) + (left:1) $);
\end{tikzpicture}
\end{center}

\paragraph*{Application to deformed gentle algebras}
The strategy to obtain our deformed mirror functor for $ \Gtl_q Q $ is to apply the deformed Cho-Hong-Lau construction to $ \Gtl_q Q $. Starting form the category $ \Gtl_q Q $ and following the definitions outlined in our auxiliary paper \papertwoA, we construct the twisted completion $ \Tw\Gtl_q Q $ and its minimal model $ \H\Tw\Gtl_q Q $. We apply the deformed Cho-Hong-Lau construction to the subcategory of zigzag paths $ \H\DefZigzagCat ⊂ \H\Tw\Gtl_q Q $. The result is a functor of $ A_∞ $-deformations $ \H\Tw\Gtl_q Q → \MF(\Jac_q \mirQ, ℓ_q) $. By virtue of the deformed Cho-Hong-Lau construction, the deformed Jacobi algebra $ \Jac_q \mirQ $ and the deformed potential $ ℓ_q $ are determined from the deformed $ A_∞ $-structure on $ \H\DefZigzagCat $.

\begin{center}
\begin{tikzpicture}
\path (0, 0) node[align=center] (A) {$ \H\DefZigzagCat ⊂ \HTw\Gtl_q Q $ \\ Deformed category of zigzag paths} (9, 0) node[align=center] (B) {$ F_q: \Gtl_q Q → \mf(\Jac_q \mirQ, ℓ_q) $ \\ Deformed mirror functor};
\path[draw, ->, decorate, decoration={snake, amplitude=0.2em, post length=0.5em}] ($ (A.east) + (right:1) $) to ($ (B.west) + (left:1) $);
\end{tikzpicture}
\end{center}

\paragraph*{Assembly of the materials}
The assembly of the total package is divided into the three papers as follows: In the first paper \paperone, we classify all deformations of $ \Gtl Q $ up to gauge equivalence. In the present second paper, we select one certain broad deformation $ \Gtl_q Q $ of $ \Gtl Q $. This deformation induces a deformation $ \HTw\Gtl_q Q $ of the derived category $ \H\Tw\Gtl Q $. With the help of our auxiliary paper \papertwoA, we calculate the deformed $ A_∞ $-structure on the subcategory $ \H\DefZigzagCat ⊂ \HTw\Gtl_q Q $ given by zigzag paths. In the third paper \paperthree, we prove the deformed Cho-Hong-Lau construction. Simply plugging in the description of $ \H\DefZigzagCat $ from the present paper gives the desired mirror functor $ F_q: \Gtl_q Q → \mf(\Jac_q \mirQ, ℓ_q) $. This result amounts to a wide range of deformed mirror equivalences for punctured surfaces.

\subsection*{Results}
We present here the results of this paper in a non-technical manner. The precise statements can be found in \autoref{th:uncurving-th}, \ref{th:subdisk-minmodel-th} and \ref{th:subdisk-main-th}.

\paragraph*{Uncurving of band objects}
The objects of the derived category $ \HTw\Gtl Q $ have been classified up to isomorphism in \cite{HKK}. They fall into two classes: the string objects and the band objects. Geometrically, a string corresponds to a curve $ γ: [0, 1] → |Q| $ which starts and ends at punctures. A band object corresponds to a closed curve $ γ: S^1 → |Q| $ which does not hit any punctures. Both string objects and band objects can also be interpreted as objects in the deformed twisted completion $ \Tw\Gtl_q Q $. In \autoref{th:uncurving-th}, we show that for the typical band object this curvature can be gauged away.

\paragraph*{Analysis of deformed Kadeishvili trees}
We explicitly describe the minimal model $ \H\DefZigzagCat ⊂ \HTw\Gtl_q Q $ in terms of immersed disks. Our strategy is to compute the minimal deformed $ A_∞ $-structure on $ \H\DefZigzagCat ⊂ \HTw\Gtl_q Q $ by means of the deformed Kadeishvili theorem from \papertwoA. The deformed Kadeishvili theorem instructs us to choose a homological splitting for $ \ZigzagCat $, gauge away the curvature of $ \DefZigzagCat $ and perform an infinitesimal “change of basis” on the homological splitting. After this procedure, the $ A_∞ $-structure of $ \H\DefZigzagCat $ is described in terms of Kadeishvili trees. In the present paper, we analyze all these trees in detail.

\paragraph*{Minimal model of the deformed category of zigzag paths}
We succeed in describing the minimal model $ \H\DefZigzagCat ⊂ \HTw\Gtl_q Q $ in terms of immersed disks. Our strategy is to match all result components of Kadeishvili trees with immersed disks. We find four types of immersed disks, the CR, ID, DS and DW disks. The precise description is stated in \autoref{th:subdisk-minmodel-th}. Once we restrict to the transversal part of $ \H\DefZigzagCat $, the description reduces to the smooth immersed disks used for the definition of the relative Fukaya category. Explicitly, the transversal part of $ \H\DefZigzagCat $ agrees with the subcategory of $ \relpreFuk Q $ given by zigzag curves.

\subsection*{The minimal model calculation}
The main storyline of this paper is the calculation of the minimal model $ \H\DefZigzagCat $ by means of our deformed Kadeishvili theorem from \papertwoA. We review here a special case of the deformed Kadeishvili construction, in which the construction reduces to three steps. We describe how each of the three steps plays out during the specific calculation of $ \H\DefZigzagCat $. This description covers the materials contained in \autoref{sec:splitting} till \ref{sec:subdisk}.

\paragraph*{The deformed Kadeishvili theorem \papertwoA}
The classical Kadeishvili theorem states that every $ A_∞ $-category has a minimal model. By definition, a minimal model of an $ A_∞ $-category $ \cat C $ is any $ A_∞ $-category $ \H\cat C $ with vanishing differential $ μ^1_{\H\cat C} $ such that $ \cat C $ and $ \H\cat C $ are quasi-isomorphic. When $ \cat C_q $ is an (infinitesimally curved) deformation of $ \cat C $, it is not clear a priori what a minimal model should be and whether it exists. In our auxiliary paper \papertwoA, we fix a definition of minimal models for deformed $ A_∞ $-categories and show that every deformed $ A_∞ $-category $ \cat C_q $ has a minimal model. We also provide a deformed Kadeishvili construction to compute a minimal model $ \H\cat C_q $ explicitly.

In the present paper, we apply the deformed Kadeishvili construction in order to compute the minimal model $ \H\Tw\Gtl_q Q $. For \autoref{sec:splitting} till \ref{sec:subdisk}, a simplified version of our deformed Kadeishvili construction suffices. This simplified version starts from the datum of a homological splitting $ H ⊕ I ⊕ R $ for $ \cat C $ and the assumption that $ \cat C_q $ has vanishing curvature and satisfies a technical condition with respect to the homological splitting. In this case, the minimal model $ \H\cat C_q $ can be computed as follows:
\begin{enumerate}
\item Perform an infinitesimal change on the homological splitting in order to adapt it to $ μ^1_q $.
\item Calculate the deformed codifferential $ h_q $ and projection $ π_q $.
\item Define the structure of $ \H\cat C_q $ by sums over deformed Kadeishvili trees.
\end{enumerate}

\paragraph*{Gauging away curvature of $ \DefZigzagCat $} Before we can apply the deformed Kadeishvili construction to $ \DefZigzagCat $, we have to gauge away the curvature of $ \DefZigzagCat $. The gauge consists of applying our “complementary angle trick” which we define and treat in the larger generality of band objects. The idea is that the twisted complexes contained in $ \DefZigzagCat $ consist of sums of arcs, with twisted differential given by angles between those arcs. Every angle comes with a certain complementary angle. Our “complementary angle trick” consists of adding the complements of these angles, weighted by deformation parameters, to the twisted differential. This trick succeeds at removing the curvature of $ \DefZigzagCat $.

\paragraph*{Choice of homological splitting}
As a preparatory step for the deformed Kadeishvili construction, we are supposed to choose a homological splitting $ H ⊕ I ⊕ R $ for $ \ZigzagCat $. There are many possible homological splittings, but not all make sense from a geometric point of view. In \autoref{sec:splitting}, we choose one specific homological splitting. To define our splitting, we have to choose an explicit basis for the cohomology of every hom spaces in $ \ZigzagCat $. Since we expect to obtain the relative Fukaya category as minimal model, we choose cohomology basis elements which are geometrically located as close as possible to the intersection points of the zigzag curves.

\paragraph*{Step 1} For the first step of the Kadeishvili construction, we are supposed to calculate the infinitesimal base change. It requires from us that we evaluate the deformed differential $ μ^1_{\DefZigzagCat} $. In \autoref{sec:deformed}, we execute this by investigating all possible contributions to products $ μ^1_{\DefZigzagCat} (ε) = μ^1_q (δ, …, ε, …, δ) $. We introduce the notions of “E, F, G, H disks” and “tails” as bookkeeping tool to systematically construct the required infinitesimal base change.

Already at the present stage after the first step, we have a strong indication that we will obtain the relative Fukaya category as result of the calculation. Take for granted that the basis elements of $ H $ can be identified with intersection points of zigzag curves. The infinitesimal base change of the first step adds infinitesimal amounts of $ R $ to the cohomology basis elements in $ H $. Visually, the interpretation is that the intersection points “grow tails” in all directions where they could possibly bound disks. This is a strong indication that we will obtain the relative Fukaya category as a result.

\paragraph*{Step 2} For the second step, we are supposed to calculate the deformed codifferential $ h_q $ and deformed projection $ π_q $. In \autoref{sec:deformed-codif}, we calculate the deformed codifferential $ h_q (ε) $ and $ π_q (ε) $ for the most important morphisms $ ε $ between zigzag paths. This requires a detailed analysis of the surroundings of $ ε $, which we capture in terms of what we call situations of type A, B, C and D. It turns out also $ h_q (ε) $ comes with an infinitesimal “tail” pointing in all possible directions that can bound disks with $ ε $. We end up with expressions for the deformed codifferential of any morphism.

Very specifically, the reader will see recurrent use of the codifferential expression $ h_q (βα) $ throughout the paper. In this context, the morphism $ β $ always denotes a $ β $-angle associated with a “type A situation”. The analysis shows that $ β $-angles act as extending link between multiple portions of a relative Fukaya disk, which renders them the most powerful angles in this paper.

\paragraph*{Step 3} For the third step, we have to evaluate deformed Kadeishvili trees. In \autoref{sec:resultcomp}, we start with a careful characterization of all results that can possibly come out of the Kadeishvili trees. The simplest Kadeishvili trees can of course be translated into relative Fukaya disks directly. Results of all other Kadeishvili trees are instead results of iterated applications of the deformed product $ μ_q $ of $ \Gtl_q $. We introduce “result components” as a bookkeeping tool for evaluating Kadeishvili trees and provide a full characterization of how result components are derived from each other within a Kadeishvili tree.

\paragraph*{Interpretation of the result} In \autoref{sec:subdisk}, we show how to turn result components of Kadeishvili trees into disks by an explicit method. Due to our characterization, every result component comes with a history, a way in which it was derived from simpler result components. We devise an inductive procedure to draw a disk from a result component. The type of the outcome is not exactly the same as a relative Fukaya disk, but is what we call an SL disk (shapeless disk). The drawing procedure works as follows: Departing from the leaves of the Kadeishvili tree, we start drawing a small portion of the SL disk. As multiple result components are merged into one at any node in the tree, we glue together their small portions. When we reach the root of the tree, we conclude the drawing by closing the SL disk with an output mark. All in all, we have assigned this way an SL disk to a result component.

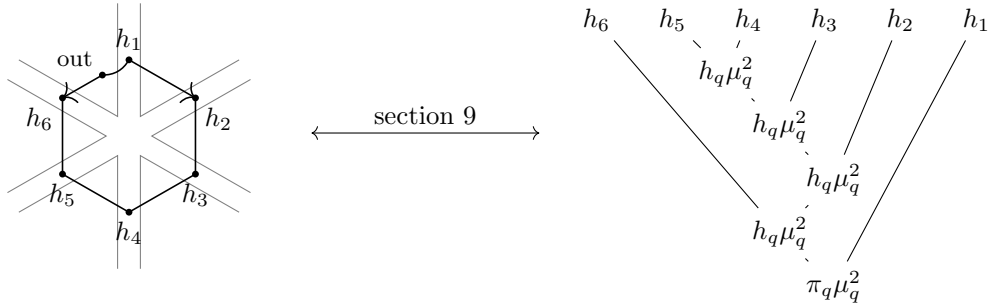
\begin{center}
\begingroup
\newcommand{\hexagondraw}{\path[draw, gray] (0, 0) -- ++(up:1) coordinate[midway] (A) coordinate[pos=0.7] (Ai) -- ++(330:1) coordinate[midway] (B) coordinate[pos=0.3] (Bi) ++(60:0.2) -- ++(150:1) coordinate[midway] (C) coordinate[pos=0.7] (Ci) -- ++(30:1) coordinate[midway] (D) coordinate[pos=0.3] (Di) ++(120:0.2) -- ++(210:1) coordinate[midway] (E) coordinate[pos=0.7] (Ei) -- ++(up:1) coordinate[midway] (F) coordinate[pos=0.3] (Fi) ++(left:0.2) -- ++(down:1) coordinate[midway] (G) coordinate[pos=0.7] (Gi) -- ++(150:1) coordinate[midway] (H) coordinate[pos=0.3] (Hi) ++(240:0.2) -- ++(330:1) coordinate[midway] (I) coordinate[pos=0.7] (Ii) -- ++(210:1) coordinate[midway] (J) coordinate[pos=0.3] (Ji) ++(300:0.2) -- ++(30:1) coordinate[midway] (K) coordinate[pos=0.7] (Ki) -- ++(down:1) coordinate[midway] (L) coordinate[pos=0.3] (Li);
\path ($ (F)!0.5!(G) $) coordinate (1);
\path ($ (D)!0.5!(E) $) coordinate (2);
\path ($ (B)!0.5!(C) $) coordinate (3);
\path ($ (L)!0.5!(A) $) coordinate (4);
\path ($ (J)!0.5!(K) $) coordinate (5);
\path ($ (H)!0.5!(I) $) coordinate (6);}

\begin{tikzpicture}[scale=1.5]
\hexagondraw
\path (1) -- (6) coordinate[pos=0.4] (out);
\path (out) node[above left] {\small out};
\foreach \i in {1, 2, 3, 4, 5, 6, out} \path[fill] (\i) circle[radius=0.03];
\path (1) node[above] {$ h_1 $};
\path (2) node[below right] {$ h_2 $};
\path (3) node[below] {$ h_3 $};
\path (4) node[below] {$ h_4 $};
\path (5) node[below] {$ h_5 $};
\path (6) node[below left] {$ h_6 $};
\path[draw, semithick, -{To[scale=1.5]}] (1) -- (2);
\path[draw, semithick] (2) -- (3);
\path[draw, semithick] (3) -- (4);
\path[draw, semithick] (4) -- (5);
\path[draw, semithick] (5) -- (6);
\path[draw, semithick, {To[scale=1.5]}-] (6) -- (out) to[bend right] (1);
\path (1) -- (2); 
\path (2) -- (3); 
\path (3) -- (4); 
\path (4) -- (5); 
\path (5) -- (6); 
\path (6) -- (out); 
\path[draw, <->] (1.5, 1.2) -- (3.5, 1.2) node[midway, above] {\autoref{sec:subdisk}};
\begin{scope}[shift={(4, 2.2)}]
\path node (A) {$ h_6 $} node[right of=A] (B) {$ h_5 $} node[right of=B] (C) {$ h_4 $} 
node[right of=C] (D) {$ h_3 $} node[right of=D] (E) {$ h_2 $} node[right of=E] (F) {$ h_1 $}
node[below right of=B] (G) {$ h_q μ^2_q $} edge (B) edge (C)
node[below right of=G] (H) {$ h_q μ^2_q $} edge (G) edge (D)
node[below right of=H] (I) {$ h_q μ^2_q $} edge (H) edge (E)
node[below left of=I] (J) {$ h_q μ^2_q $} edge (I) edge (A)
node[below right of=J] (K) {$ π_q μ^2_q $} edge (J) edge (F);
\end{scope}
\end{tikzpicture}
\endgroup

\end{center}

As a final step, we classify all SL disks we have obtained this way. It turns out that the SL disks obtained are of four types, which we call CR, ID, DS and DW disks. In other words, the higher products on $ \H\DefZigzagCat $ are precisely computed by SL disks in the surface that belong to one of those four types. It is useful to know in advance that many, but not all of these disks are transversal. In fact, the transversal ones among them are all of CR type and match exactly the (transversal) relative Fukaya disks. This finishes the computation of the minimal model $ \H\DefZigzagCat $.

\subsection*{Context and philosophical highlights}
The results of our paper are very specific. To get a sense of their general meaning, we put the results into context. We comment on the following philosophical highlights:
\begin{itemize}
\item All deformed $ A_∞ $-categories have derived categories.
\item Computational techniques carry over to the deformed case.
\item Minimal model calculations are possible if one is sensitive to the result.
\item The deformed gentle algebra $ \Gtl_q Q $ is a relative wrapped Fukaya category.
\item Hamiltonian deformations arise naturally from representation theory.
\end{itemize}

\paragraph*{Derived categories of $ A_∞ $-deformations}
In formal $ A_∞ $-deformation theory, one regards infinitesimal deformations of a given $ A_∞ $-structure such that the $ A_∞ $-relations are preserved. As we review in \papertwoA, one needs to permit the deformations to have curvature in order to obtain a homologically invariant notion of deformations. More precisely, one allows infinitesimal curvature in the sense that it lies in a multiple of the maximal ideal of the local ring.

Dealing with curvature is often regarded as tedious, because the curvature prevents the differential from squaring to zero. The presence of curvature is often referred to as the “curvature problem”. A main question is how to gauge away the curvature or otherwise how to deal with the remaining curvature. An instance of the uncurving problem has been studied by Lowen and Van den Bergh \cite{Lowen-vdB}, on which we comment in \autoref{sec:literature-LvdB}.

Part of the “curvature problem” is that curved $ A_∞ $-categories do not have derived categories since their differential does not square to zero. In \papertwoA, we show that the game changes when the curvature is only infinitesimal. We define a twisted completion $ \Tw\cat C_q $ and minimal model $ \H\cat C_q $ for any $ A_∞ $-deformation $ \cat C_q $. The derived category of $ \cat C_q $ can then be written as $ \H\Tw\cat C_q $. It is a deformation of $ \H\Tw\cat C $.

\paragraph*{Minimal model techniques}
Minimal models of $ A_∞ $-deformations are defined abstractly by inducing the deformation afterwards. More precisely, let $ \cat C $ be an $ A_∞ $-category and $ \cat C_q $ a deformation. Let $ \H\cat C $ be a minimal model for $ \cat C $ with quasi-isomorphism $ π: \cat C → \H\cat C $. In \papertwoA\ we define the minimal model $ \H\cat C_q $ to be the deformation of $ \H\cat C $ induced from $ \cat C_q $ via $ π $. This abstract approach means that the “minimal model” $ \H\cat C_q $ may have infinitesimal curvature as well as a residual infinitesimal differential. This is not the same as taking cohomology of the hom complexes $ (\Hom_{\cat C_q} (X, Y), μ^1_q) $. In fact, the deformed differential of an $ A_∞ $-deformation need not even square to zero because of the curvature.

Minimal models of $ A_∞ $-categories can classically be computed by means of homological splittings and Kadeishvili trees. In our deformed Kadeishvili theorem, we show that this method carries over to the deformed case. The starting point is an $ A_∞ $-category together with an deformation $ \cat C_q $. The difficulties encountered in constructing the minimal model are the presence of curvature $ μ^0_q $, the fact that the deformed differential $ μ^1_q $ does not square to zero and the fact that $ μ^1_q $ is not compatible with the homological splitting chosen for $ \cat C $. In \autoref{sec:2Bkadeishvili}, we review how to adapt the Kadeishvili construction to these special circumstances. We view our deformed Kadeishvili theorem as evidence that computational techniques which apply to $ A_∞ $-deformations can be tweaked in order to apply to $ A_∞ $-deformations as well.

\paragraph*{Discrete relative Fukaya category}
Many different constructions of Fukaya categories are available in the literature. The most general approach is the reference work of Seidel \cite{Seidel}. For relative Fukaya categories, a new reference is the construction of Sheridan and Perutz \cite{Sheridan-Perutz}. For the case of punctured surfaces, there are many further specific models available. One can distinguish whether they depart from the discrete side of gentle algebras or from the smooth side of actual Fukaya categories, and whether they consider the punctured surface alone or whether they consider Seidel's deformation. The following is a non-exhaustive overview:

\begin{center}
\begin{tabular}{c|cc}
starting point & non-deformed & deformed \\\hline
geometric & \cite[Appendix B]{Bocklandt} & \cite{Seidel-relative}, \cite{Lekili-Perutz}, \cite{Lekili-Polishchuk} \\
discrete & \cite{Bocklandt}, \cite{Bocklandt-book}, \cite{HKK} & \emph{this paper}
\end{tabular}
\end{center}

The most important reference for us is the construction of the gentle algebras of \cite{Bocklandt}. In \paperone, we proposed a candidate deformation $ \Gtl_q Q $ with the intention to provide a “relative wrapped Fukaya category” for punctured surfaces. Verifying that the deformed gentle algebra $ \Gtl_q Q $ deserves this name would at least entail proving that the transversal part of its derived category is equivalent to the relative Fukaya category, see \autoref{sec:whyshould-candidate}. It is however quite difficult to actually compute the derived category, as we witness in the present paper. If a “relative wrapped Fukaya category” existed already, this would be greatly eased, see \autoref{sec:whyshould-alternative}.

In the present paper, we succeed in showing that at least on the subcategory of zigzag paths, the derived category $ \HTw\Gtl_q Q $ and the relative Fukaya category $ \relFuk Q $ agree. Although our calculation is limited to zigzag paths, we consider our calculation strong evidence that $ \HTw\Gtl_q Q $ indeed contains $ \relFuk Q $. It is a crude verification that we have correctly transported Seidel's vision to gentle algebras and that $ \Gtl_q Q $ can be considered a relative wrapped Fukaya category.

\paragraph*{Minimal model calculations}
A highlight in this paper is our explicit computation of an entire minimal model. Such computations are scarce in the $ A_∞ $-literature and often considered tedious. Indeed, minimal model calculations are hard because of the large amount of Kadeishvili trees involved. For some calculations in the liturature, it is not necessary to perfom the calculation until the end. In the case of \cite{Bocklandt} it suffices to calculate only part of the minimal model because the rest is determined up to homotopy. In the case of \paperone, we also cut short the calculation of an $ L_∞ $-minimal model by means of grading arguments. In the present paper, we perform the minimal model calculation of $ \DefZigzagCat $ until the end.

Minimal model calculations are the core connecting bridge between the discrete and the smooth world. They are regarded as tedious, but we contend that minimal model calculation need not hurt if one has a clue regarding the outcome. The minimal model calculation in the present paper succeeds precisely because we recognize in every step the inherent geometric meaning of the terms that appear. This concerns both the choice of the homological splitting for $ \ZigzagCat $ and the evaluation of the Kadeishvili trees. In \autoref{sec:whydoes} and \ref{sec:reuse-resultcomp}, we offer further explanation on why our method of “result components” works and how to apply it in other situations.

\paragraph*{Hamiltonian deformations}
Implementing Hamiltonian deformations is one of the difficulties one encounters when defining smooth Fukaya categories. In the discrete world, one circumvents this problem by choosing such a small set of generators that the Hamiltonian deformations can be chosen canonically and disappear completely from the picture. When passing to the derived category $ \HTw\Gtl Q $, we however expect the full generality of the smooth Fukaya category to reappear. In particular, we expect to find $ A_∞ $-products on some non-transversal and expect that we can explain these products as an incarnation of Hamiltonian deformation.

In the present paper, we compute the precise $ A_∞ $-products on the deformed category $ \H\DefZigzagCat $. The starting point is our deformed Kadeishvili theorem, whose essential ingredient is a choice of homological splitting for $ \ZigzagCat $. As expected, the products of $ \H\DefZigzagCat $ agree with the products of the relative Fukaya category on transversal sequences. We however also obtain an explicit description of the products on non-transversal sequences. We show how to interpret even the products on non-transversal sequences geometrically as disks being bounded by zigzag curves and their Hamiltonian deformations. While Hamiltonian deformations have to be incorporated as an ingredient into the definition of smooth Fukaya categories from the beginning, they appear naturally through the Kadeishvili construction of the minimal model $ \H\DefZigzagCat $.

The precise shape of the products of $ \H\DefZigzagCat $ depends on the choice of homological splitting for $ \ZigzagCat $. Nevertheless, different homological splittings give quasi-equivalent minimal models $ \H\DefZigzagCat $. We have selected one specific splitting which makes it particularly easy to identify the minimal model as the relative Fukaya category. When choosing a slightly different splitting, we still expect to obtain the same products on transversal sequences, but the products on non-transversal sequences will typically change. These changed products can be interpreted geometrically as products in the relative Fukaya category under application of a different Hamiltonian deformation. While homological splittings for $ \ZigzagCat $ are a discrete and representation-theoretic notion, Hamiltonian deformations are a smooth and geometric notion. Highly simplified, we may say that choices of homological splittings correspond to choices of Hamiltonian deformations. See also \autoref{sec:whydoes-choice}.

\paragraph*{Strings and bands}
Gentle algebras $ \Gtl Q $ were originally introduced in \cite{Bocklandt} to provide a combinatorial description of the wrapped Fukaya category of the punctured surface $ |Q| \setminus Q_0 $. In contrast to the Fukaya category of $ |Q| \setminus Q_0 $, the wrapped Fukaya category also includes curves which start and end at punctures. Haiden, Katzarkov and Kontsevich \cite{HKK} classified the objects of $ \HTw\Gtl Q $ under the additional datum of a $ ℤ $-grading. Their classification indeed finds those types of curves expected from the wrapped Fukaya category. Explictly, their classification divides the objects into two classes, known as string objects and band objects. Roughly speaking, a string object is a non-closed curve running between two punctures of $ Q $ and a band object is a closed curve that avoids the punctures of $ Q $.

A string object or band object given by a curve in $ |Q| $ can be explicitly realized as a twisted complex in $ \Tw\Gtl Q $. The procedure entails approximating the curve by arcs $ a_1, …, a_k $ of $ Q $ together with angles $ α_i $ between the arcs. One then forms a twisted complex $ (\bigoplus_i a_i [s_i], δ=\sum_i α_i) $ by summing up the arcs and using the angles as twisted differential.

Seidel \cite{Seidel-relative} describes which objects in the relative Fukaya category should have curvature according to his vision. For the case of punctured surfaces, his criterion states that curves which bound a so-called teardrop should have curvature. Also those curves which are contractible in the surface $ |Q| $ should have curvature. All other objects in the relative Fukaya category should be curvature-free according to Seidel.

In the present paper, we approach Seidel's vision from the starting point of the deformed gentle algebras $ \Gtl_q Q $ instead of the relative Fukaya category. Translated to our setting, a band object $ (\bigoplus_i a_i [s_i], δ=\sum_i α_i) $ should be uncurvable if its underlying curve in the closed surface $ |Q| $ is not contractible and does not bound a teardrop. In order to make his vision true, we devise a trick to gauge away the curvature of these band objects. Our “complementary angle trick” consists essentially of adding infinitesimal multiples of the complementary angle of $ α_i $ to $ δ $ for all $ i $. In \autoref{sec:uncurving}, we verify that our trick successfully uncurves all band objects whose underlying curve in the closed surface $ |Q| $ is not contractible and does not bound a teardrop, making true Seidel's vision.

\subsection*{Data structures}
\begin{figure}
\centering
\begin{tikzpicture}
\path (-3, -3) node[align=center] (Subdisk) {$ \Subdisk(r) $ \\ subdisks};
\path (-1, -5) node (CR) {CR, ID, DS, DW disks} edge (Subdisk);
\path (2, -5) node (equals) {\Large $ = $};
\path (5, -5) node (Imm) {immersed disks};
\path (-2, 0) node[align=center] (ResComp) {$ r ∈ \PiTr $ \\ result components} edge (Subdisk);
\path (-4, 3) node[align=center] (Hq) {$ φ^{-1} (h) $ \\ deformed cohomology \\ basis elements} edge (ResComp);
\path (-2, 3) node (Piq) {$ π_q $} edge (ResComp);
\path (0, 3) node[align=center] (hq) {$ h_q (ε) $ \\ deformed codifferential} edge (ResComp);
\path (-2, 5) node[align=center] (tails) {$ T(ε) $ \\ tails of morphisms} edge (Hq) edge (Piq) edge (hq);
\path (-2, 7) node[align=center] (muq) {$ μ^1_q (ε) $ \\ deformed differential} edge (tails); 
\path (-4, 8) node[align=center] (H) {$ h: L_1 → L_2 $ \\ cohomology basis \\ elements} edge (Hq);
\path (0, 9) node[align=center] (angles) {$ ε: L_1 → L_2 $ \\ angles} edge (muq);
\path (0, 10.5) node[align=center] (zigzagpaths) {$ L $ \\ zigzag paths} edge (angles) edge (H);
\path (6, 3) node[align=center] (smInt) {$ p ∈ \smooth L_1 ∩ \smooth L_2 $ \\ intersection points} edge (Imm);
\path (6, 10) node[align=center] (zigzaglags) {$ \smooth L $ \\ zigzag curves} edge (smInt);
\path (3, 12) node[align=center] (dimer) {$ Q $ \\ dimer} edge (zigzagpaths) edge (zigzaglags);
\begin{scope}[shift={($ (equals.north) + (-0.4, 1.5) $)}, scale=0.7] 
\path (0, 0) coordinate (A);
\path (1, 0) coordinate (B);
\path (B) ++(288:1) coordinate (C);
\path (C) ++(216:1) coordinate (D);
\path (D) ++(144:1) coordinate (E);
\path[fill=lightgray!0.5] (A) -- (B) -- (C) -- (D) -- (E) -- cycle;
\foreach \i in {A, B, C, D, E} \path[fill] (\i) circle[radius=0.07];
\path[draw, semithick] ($ (A)!-0.3!(B) $) -- ($ (B)!-0.3!(A) $);
\path[draw, semithick] ($ (B)!-0.3!(C) $) -- ($ (C)!-0.3!(B) $);
\path[draw, semithick] ($ (C)!-0.3!(D) $) -- ($ (D)!-0.3!(C) $);
\path[draw, semithick] ($ (D)!-0.3!(E) $) -- ($ (E)!-0.3!(D) $);
\path[draw, semithick] ($ (E)!-0.3!(A) $) -- ($ (A)!-0.3!(E) $);
\end{scope}
\begin{scope}[shift={($ (ResComp.east) + (0, 0.5) $)}] 
\begin{pgflowlevelscope}{\pgftransformscale{0.7}}
\path node (A) {$ r_4 $} node[right of=A] (B) {$ r_3 $} node[right of=B] (C) {$ r_2 $} node[right of=C] (D) {$ r_1 $}
node[below right of=A] (E) {$ r_6 $} edge (A) edge (B) node[below left of=D] (F) {$ r_5 $} edge (C) edge (D)
node[below right of=E] {$ r_7 $} edge (E) edge (F);
\end{pgflowlevelscope}
\end{scope}
%
\begin{scope}[shift={($ (H.west) + (-1, -0.5) $)}, scale=0.5] 
\path[draw] (0, 0) -- ++(45:1) -- ++(up:1) coordinate[midway] (m1) -- ++(45:1);
\path[draw] (1.1, 0) -- ++(135:1) -- ++(up:1) coordinate[midway] (m2) -- ++(135:1);
\path[fill] ($ (m1)!0.5!(m2) $) circle[radius=0.05];
\end{scope}
%
\begin{scope}[shift={($ (Hq.west) + (-0.7, 0.5) $)}, scale=0.5] 
\path[draw] (0, 0) -- ++(30:1) coordinate[pos=0.3] (1-start) -- ++(up:1) coordinate[midway] (m1) -- ++(45:1);
\path[draw] (1.3, 0) -- ++(150:1) coordinate[pos=0.3] (1-end) -- ++(up:1) coordinate[midway] (m2) -- ++(135:1);
\path[draw] (0, 0) -- ++(down:1) -- ++(330:0.6) -- ++(210:0.6) coordinate[pos=0.6] (2-start) -- ++(down:0.6) -- ++(330:0.6) -- ++(210:0.6) coordinate[pos=0.6] (3-start);
\path[draw] (1.3, 0) -- ++(down:1) -- ++(210:0.6) -- ++(330:0.6) coordinate[pos=0.6] (2-end) -- ++(down:0.6) -- ++(210:0.6) -- ++(330:0.6) coordinate[pos=0.6] (3-end);
\foreach \i in {1, 2, 3} \path[draw, ->, bend right=60] (\i-start) to (\i-end);
\path[fill] ($ (m1)!0.5!(m2) $) circle[radius=0.05];
\end{scope}
\begin{scope}[shift={($ (muq.east) + (0.2, 0.6) $)}, scale=0.5] 
\path[draw] (0, 0) -- ++(330:1) -- ++(210:1) coordinate[pos=0.6] (1-start) -- ++(down:1) -- ++(330:1) -- ++(210:1);
\path[draw] (2, 0) -- ++(210:1) -- ++(330:1) coordinate[pos=0.6] (1-end) -- ++(down:1) -- ++(210:1) -- ++(330:1);
\path[draw, bend right=60, ->] (1-start) to node[midway, above] {$ ε $} (1-end);
\path (1, -1.7) node {E};
\end{scope}
\begin{scope}[shift={($ (tails.east) + (1, 4.5) $)}] 
\begin{pgflowlevelscope}{\pgftransformscale{0.6}}
\path node (A) {$ ε $}
node[below of=A] (B) {E} edge (A)
node[below of=B] (C) {E} edge (B)
node[below of=C] (D) {F} edge (C)
node[left of=D] {E} edge (C) node[right of=D] {H} edge (C);
\end{pgflowlevelscope}
\end{scope}
\begin{scope}[shift={($ (Subdisk.east) + (-0.5, -0.5) $)}] 
\begin{pgflowlevelscope}{\pgftransformscale{0.5}}
\path (0, 0) to[bend right=30] ++(-0.5, -1) coordinate (start);
\path[draw, thick, rounded corners] (start) to[bend left=50, looseness=1.5] coordinate[pos=0.5] (1-l) ++(-0.3, -1) to[bend right=100, looseness=4] ++(0.7, 0) to[bend left=50, looseness=1.5] coordinate[pos=0.6] (1-r) ++(0.1, 0.5) to[bend left=50, looseness=1.5] coordinate[pos=0.4] (2-l) ++(0.1, -0.5) to[bend right=100, looseness=4] ++(0.5, 0) to[bend left=50, looseness=1.5] coordinate[pos=0.6] (2-r) ++(0.1, 0.5) to[bend left=50, looseness=1.5] coordinate[pos=0.4] (3-l) ++(0.1, -0.5) to[bend right=100, looseness=4] ++(0.5, 0) to[bend left=50, looseness=1] coordinate[pos=0.42] (3-r) coordinate[pos=0.9] (7-l) ++(0, 0.8) to[bend left=10, looseness=1] coordinate[pos=0.55] (6-l) ++(1, -0.2) to[bend left=40, looseness=1.5] coordinate[pos=0.5] (4-l) ++(0, -0.5) to[bend right=100, looseness=4] ++(0.5, 0) to[bend left=50, looseness=1.5] coordinate[pos=0.42] (4-r) ++(0.1, 0.5) to[bend left=50, looseness=1.5] coordinate[pos=0.5] (5-l) ++(0.1, -0.5) to[bend right=100, looseness=4] ++(0.5, 0) to[bend left=30, looseness=1] coordinate[pos=0.6] (5-r) ++(-0.3, 0.6) to[bend right=80, looseness=1.5] ++(-1, 0.5) to[bend left=90, looseness=3] coordinate[pos=0.5] (6-r) ++(-0.4, 0.1) to[bend right=80, looseness=2] ++(-0.5, 0.1)
to[bend left=100, looseness=4] coordinate[pos=0.3] (7-r) ++(-0.5, 0.1) coordinate (end);
\path[draw, thick, dashed] (start) -- (end);
\end{pgflowlevelscope}
\end{scope}
\begin{scope}[shift={($ (hq.east) + (0,  1) $)}, scale=0.5] 
\path[draw] (0, 0) -- ++(330:1) coordinate[pos=0.4] (1-start) -- ++(210:1) -- ++(down:1) -- ++(330:1) -- ++(210:1) coordinate[pos=0.6] (2-start) -- ++(down:1) -- ++(330:1) -- ++(210:1) coordinate[pos=0.6] (3-start);
\path[draw] (2, 0) -- ++(210:1) -- ++(330:1) coordinate[pos=0.6] (1-end) -- ++(down:1) -- ++(210:1) -- ++(330:1) coordinate[pos=0.6] (2-end) -- ++(down:1) -- ++(210:1) -- ++(330:1) coordinate[pos=0.6] (3-end);
\path[draw, ->, bend right=90, looseness=2] (1-start) to node[pos=0.6, above] {$ ε $} (1-end);
\foreach \i in {2, 3} \path[draw, ->, bend right=60] (\i-start) to (\i-end);
\end{scope}
\begin{scope}[shift={($ (angles.east) + (0, 0.5) $)}, scale=0.7] 
\path[draw] (0, 0) -- ++(330:1) -- ++(210:1) coordinate[pos=0.6] (1-start);
\path[draw] (2, 0) -- ++(210:1) -- ++(330:1) coordinate[pos=0.6] (1-end);
\path[draw, ->, bend right=60] (1-start) to node[midway, above] {$ ε $} (1-end);
\end{scope}
\begin{scope}[shift={($ (zigzagpaths.west) + (-1, 1) $)}, scale=0.5] 
\path[draw] (0, 0) -- ++(300:1) -- ++(40:1) -- ++(250:1) -- ++(15:1) -- ++(230:1) -- ++(320:1) -- ++(180:1) -- ++(300:1);
\end{scope}
\begin{scope}[shift={($ (zigzaglags.east) + (-1, -0.5) $)}, scale=0.5] 
\path[draw] (0, 0) -- ++(300:1) coordinate[midway] (A) -- ++(40:1) coordinate[midway] (B) -- ++(250:1) coordinate[midway] (C) -- ++(15:1) coordinate[midway] (D) -- ++(230:1) coordinate[midway] (E) -- ++(320:1) coordinate[midway] (F) -- ++(180:1) coordinate[midway] (G) -- ++(300:1) coordinate[midway] (H);
\path[draw, rounded corners, thick] ($ (A)!-0.3!(B) $) -- (A) -- (B) -- (C) -- (D) -- (E) -- (F) -- (G) -- (H) -- ($ (H)!-0.3!(G) $);
\end{scope}
\begin{scope}[shift={($ (smInt.east) + (-1.3, -2.3) $)}, scale=0.4] 
\path[draw, color=lightgray] (0, 0) -- ++(30:2) coordinate[pos=0.5] (target-tail) -- ++(up:1.5) coordinate[pos=0.5] (target-middle) coordinate (5-head);
\path[draw, color=lightgray] (5-head) -- ++(30:2) coordinate[pos=0.5] (target-head) node[near end, below] {};
\path[draw, color=lightgray] (3.1, 0) -- ++(150:2) coordinate[pos=0.5] (source-tail) -- ++(up:1.5) coordinate[pos=0.5] (source-middle) coordinate (2-head);
\path[draw, color=lightgray] (2-head) -- ++(150:2) node[near end, below] {} coordinate[pos=0.5] (source-head);
\path (source-middle) -- (target-middle) coordinate[pos=0.5] (middle);
\path[draw, rounded corners, semithick] ($ (middle)!2!(target-tail) $) -- (middle) -- ($ (middle)!1!(target-head) $) coordinate[pos=0.7] (target);
\path[draw, rounded corners, semithick] ($ (middle)!1!(source-tail) $) -- (middle) -- ($ (middle)!2!(source-head) $) coordinate[pos=0.7] (source);
\path[fill] (middle) circle[radius=0.15];
\end{scope}
\end{tikzpicture}
\caption{This graph depicts the essential data structures used in the paper. The left part of the graph depicts data structures used for the computation of the minimal model $ \H\DefZigzagCat $. The right part depicts the construction of the Fukaya category. At the end of the paper, the minimal model $ \H\DefZigzagCat $ is described by means of immersed disks. The relative Fukaya category is defined in terms of immersed disks as well. Ultimately, we conclude that $ \H\DefZigzagCat $ is equivalent to the subcategory of the relative Fukaya category given by zigzag curves.}
\label{fig:intro-data}
\end{figure}
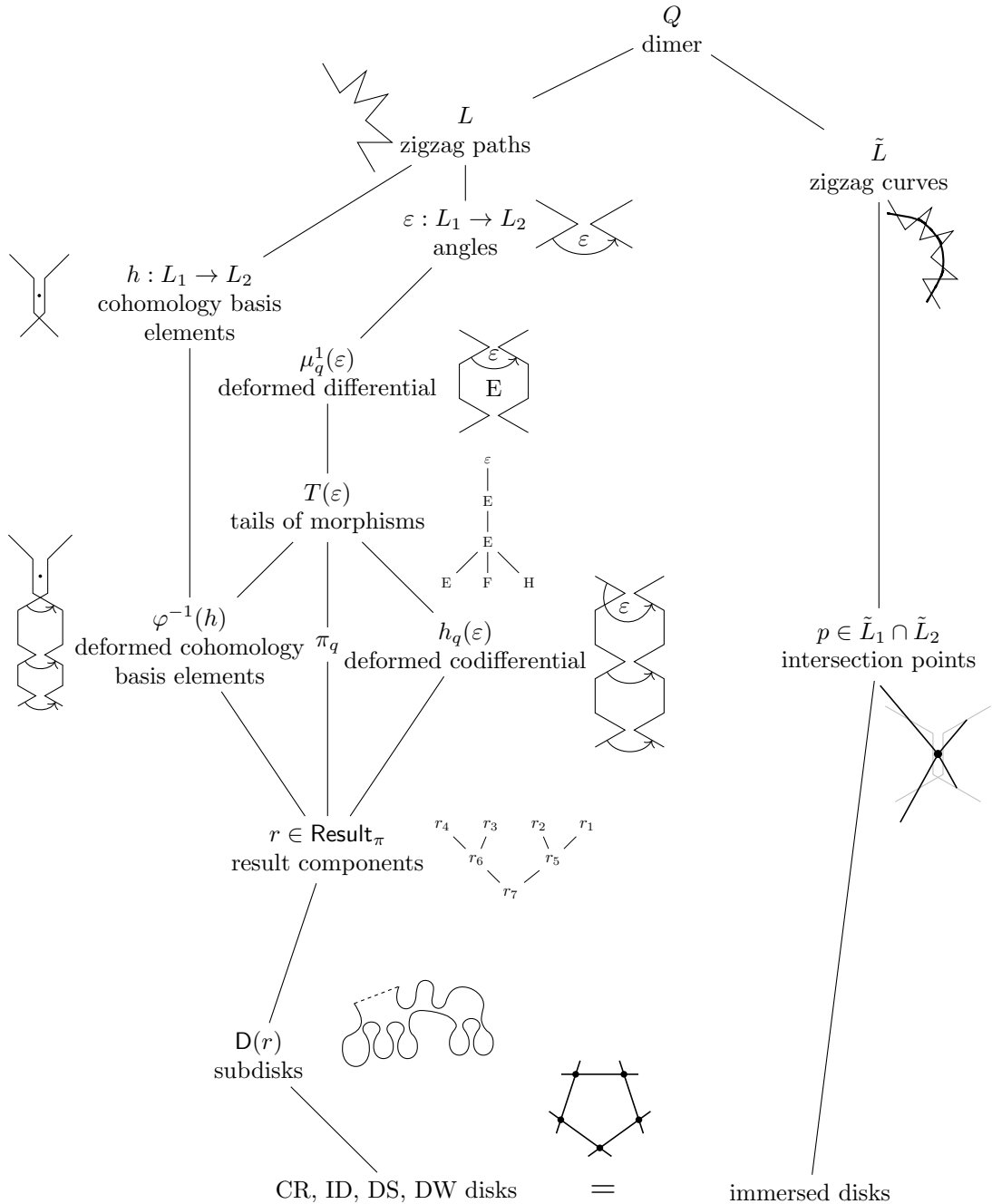

Most of our calculation does not go beyond simple inspection of arc systems and linear algebra. However, organizing result components and matching them with disks requires us to devise a large amount of data structures and fill them with data. For an overview, we depict in \autoref{fig:intro-data} the essential data structures. We shall here explain the purpose and development of these datastructures and which data flows from which structure into which one.

The starting point is a dimer $ Q $, which is a specific type of quiver embedded in a surface. It gives rise to the discrete notion of zigzag paths and the smooth notion of zigzag curves. On the smooth side, the zigzag curves give rise to the notion of intersection points and smooth immersed disks, the foundations of Fukaya categories.

On the discrete side, we regard the category $ \ZigzagCat $ of zigzag paths. A morphism $ ε: L_1 → L_2 $ between two zigzag paths consists of an angle between arcs of $ L_1 $ and $ L_2 $. We determine a basis of cohomology elements for $ \ZigzagCat $. We also define the category $ \DefZigzagCat $ of deformed zigzag paths. Examining the deformed differential $ μ^1_{\DefZigzagCat} $ of this category gives rise to four types of disks which we call E, F, G, H disks. We introduce the auxiliary notion of “tails”. The tail of an angle $ ε: L_1 → L_2 $ is a tree whose nodes are decorated with E, F, G, H disks.

The deformed Kadeishvili theorem gives rise to notions of deformed cohomology basis elements, a deformed projection $ π_q $ and a deformed codifferential $ h_q $. We can describe them explicitly by means of tails. According to the deformed Kadeishvili theorem, the product structure of the minimal model $ \H\DefZigzagCat $ is described in terms of sums over trees. We define a notion of “result components” which serves to systematically track the results of evaluations of trees. From a result component we build a “subdisk” by drawing zigzag curve segments and intersection points. Subdisks of result components are immersed disks and fall into four classes which we call the CR, DS, ID and DW disks.

From the perspective of data structures, this finishes the construction of the minimal model $ \H\DefZigzagCat $. Both $ \H\DefZigzagCat $ and the relative Fukaya category are described by immersed disks. Therefore the category $ \DefZigzagCat $ of deformed zigzag paths is quasi-isomorphic to the subcategory of zigzag curves of the relative Fukaya category. On the level of data structures, this finishes the main theorem.

\subsection*{Structure of the paper}
In \autoref{sec:2Bainfty} we recall $ A_∞ $-categories, their deformations and our deformed Kadeishvili theorem. In \autoref{sec:prelim-gtl}, we recall gentle algebras and deformed gentle algebras. In \autoref{sec:fukaya}, we recall basics of Fukaya categories and explain their subcatgories of zigzag curves. In \autoref{sec:uncurving} we present the uncurving procedure for band objects. In \autoref{sec:splitting}, we exhibit the category of zigzag paths $ \ZigzagCat $ together with a homological splitting. In \autoref{sec:deformed}, we present the deformed version of this homological splitting, together with reference material for the rest of the paper. In \autoref{sec:resultcomp}, we introduce the tool of result components to enumerate products in the minimal model $ \H\DefZigzagCat $. In \autoref{sec:subdisk}, we devise a simple drawing method to transform these result components into immersed disks. In \autoref{th:subdisk-minmodel-th}, we describe explicitly the structure of the minimal model $ \H\DefZigzagCat $ in terms of immersed disks. Our main result \autoref{th:subdisk-main-th} states that $ \H\DefZigzagCat $ has the same products on transversal sequences as the relative Fukaya category.

This paper contains several appendices which are devoted to technical proofs and additional explanation. In \autoref{sec:2examples}, we provide examples of immersed disks together with their corresponding result components. The aim is to facilitate understanding of how disks arise from the minimal model $ \H\DefZigzagCat $. In \autoref{sec:trick}, we complete the proof of uncurvability of band objects. In \autoref{sec:classification}, we finish the proof of the main result by providing an explicit inverse construction which maps CR, DS, ID and DS disks to their corresponding result components. In \autoref{sec:sphere}, we study the case of specific sphere dimers, among which the pair of pants. These dimers are not geometrically consistent and fall outside of the scope of the rest of the paper, but we have included their calculation due to their relevance in mirror symmetry. In \autoref{sec:mirobjects}, we compute a small class of products in the category $ \HTw\Gtl_q Q $ which go beyond zigzag paths. Specifically, this concerns products of morphisms between arcs and zigzag paths from which we determine the mirror objects $ F_q (a) ∈ \MF(\Jac_q \mirQ, ℓ_q) $ in the third paper. In \autoref{sec:literature}, we discuss the relation with the literature in more detail. In \autoref{sec:whyshould}, we explain why one is led to believe from an a priori perspective that $ \H\DefZigzagCat $ agrees with the relative Fukaya category. In \autoref{sec:whydoes}, we summarize from an a posteriori perspective why the very technical calculation of $ \H\DefZigzagCat $ contained in this paper succeeds. In \autoref{sec:reuse}, we share insight on how to reuse the constructions in this paper for other purposes. In \autoref{sec:notation}, we collect notation specific to this paper.

\subsection*{Conventions}
During the course of the paper, we play in two different contexts. In \autoref{sec:uncurving}, the starting point is an arc system $ \cA $ which has no monogons or digons in the closed surface. We summarize this in the [NMDC] condition. In \autoref{sec:splitting} till \ref{sec:subdisk}, the starting point is a geometrically consistent dimer $ Q $. Every zigzag path is supposed to come with a chosen spin structure and locations of identity and co-identity endomorphism. We summarize this setup in \autoref{conv:alpha0-direction}.

\papertwoacknowledgements

\section{Preliminaries on $ A_∞ $-categories}
\label{sec:2Bainfty}
In this section, we recollect background material on $ A_∞ $-categories and fix notation. In \autoref{sec:2Binfty-ainfty}, we recall $ A_∞ $-categories, their functors, twisted completion and minimal models. In \autoref{sec:2Bainfty-defo}, we recall completed tensor products, deformations of $ A_∞ $-categories and their functors. We very briefly comment on the construction of twisted completion and minimal model for $ A_∞ $-deformations from \papertwoA. In \autoref{sec:2Bkadeishvili}, we recall our deformed Kadeishvili theorem.

\subsection{$ A_∞ $-categories}
\label{sec:2Binfty-ainfty}
In this section we recall $ A_∞ $-categories, their twisted completion and functors. The material is standard and can for instance be found in \cite{Bocklandt-book}. Throughout we work over an algebraically closed field of characteristic zero and write $ ℂ $.

\begin{definition}
\label{def:2Binfty-ainfty-def}
A ($ ℤ $- or $ ℤ/2ℤ $-graded, strictly unital) \emph{$ A_∞ $-category} $ \cat C $ consists of a collection of objects together with $ ℤ $- or $ ℤ/2ℤ $-graded hom spaces $ \Hom(X, Y) $, distinguished identity morphisms $ \id_X ∈ \Hom^0(X, X) $ for all $ X ∈ \cat C $, together with multilinear higher products
\begin{equation*}
μ^k: \Hom(X_k, X_{k+1}) ¤ … ¤ \Hom(X_1, X_2) → \Hom(X_1, X_{k+1}), \quad k ≥ 1
\end{equation*}
of degree $ 2-k $ such that the $ A_∞ $-relations and strict unitality axioms hold: For every compatible morphisms $ a_1, …, a_k $ we have
\begin{align*}
& \sum_{0 ≤ j < i ≤ k} (-1)^{‖a_n‖ + … + ‖a_1‖} μ(a_k, …, μ(a_i, …, a_{j+1}), a_j, …, a_1) = 0, \\
& μ^2 (a, \id_X) = a, ~ μ^2 (\id_Y, a) = (-1)^{|a|} a, ~ μ^{≥3} (…, \id_X, …) = 0.
\end{align*}
\end{definition}


Next we recall the additive completion $ \Add\cat C $ of an $ A_∞ $-category $ \cat C $. This category consists of formal sums of shifted objects. The hom space between two objects consists of matrices of morphisms between the summands.

\begin{definition}
\label{def:2Binfty-ainfty-Add}
Let $ \cat C $ be an $ A_∞ $ category with product $ μ_{\cat C} $. The additive completion $ \Add \cat C $ of $ \cat C $ is the category of formal sums of shifted objects of $ \cat C $:
\begin{equation*}
A_1 [k_1] ⊕ … ⊕ A_n [k_n].
\end{equation*}
The hom space between two such objects $ X = \bigoplus A_i [k_i] $ and $ Y = \bigoplus B_i [m_i] $ is
\begin{equation*}
\Hom_{\Add\cat C} (X, Y) = \bigoplus_{i, j} \Hom_{\cat C} (A_i, B_j) [m_j - k_i].
\end{equation*}
Here $ [-] $ denotes the right-shift. The products on $ \Add\cat C $ are given by multilinear extensions of
\begin{equation*}
μ_{\Add \cat C}^k (a_k, …, a_1) = (-1)^{\sum_{j < i} \Vert a_i \Vert l_j} μ_{\cat C}^k (a_k, …, a_1).
\end{equation*}
Here each $ a_i $ lies in some $ \Hom(X_i[k_i], X_{i+1} [k_{i+1}]) $. The integer $ l_i $ denotes the difference $ k_{i+1} - k_i $ between the shifts and the degree $ ‖a_i‖ $ is the degree of $ a_i $ as element of $ \Hom_{\cat C} (X_i, X_{i+1}) $.
\end{definition}

Next we recall the twisted completion $ \Tw\cat C $ of an $ A_∞ $-category $ \cat C $. The objects of this category are virtual chain complexes of objects of $ \cat C $:

\begin{definition}
\label{def:2Binfty-ainfty-Tw}
A \emph{twisted complex} in $ \cat C $ is an object $ X ∈ \Add\cat C $ together with a morphism $ δ ∈ \Hom^1_{\Add\cat C} (X, X) $ of degree $ 1 $ such that $ δ $ is strictly upper triangular and satisfies the Maurer-Cartan equation:
\begin{equation*}
\MC(δ) ≔ μ^1 (δ) + μ^2 (δ, δ) + … = 0.
\end{equation*}
We may refer to the morphism $ δ $ as the \emph{twisted differential}. Note that the upper triangularity ensures that this sum is well-defined. The \emph{twisted completion} of $ \cat C $ is the $ A_∞ $-category $ \Tw\cat C $ whose objects are twisted complexes. Its hom spaces are the same as for the additive completion:
\begin{equation*}
\Hom_{\Tw\cat C} (X, Y) = \Hom_{\Add\cat C} (X, Y).
\end{equation*}
The products on $ \Tw \cat C $ of $ \cat C $ are given by embracing with $ δ $'s:
\begin{equation*}
μ_{\Tw \cat C}^k (a_k, …, a_1) = \sum_{n_0, …, n_k ≥ 0} μ_{\Add \cat C} (\underbrace{δ, …, δ}_{n_k}, a_k, …, a_1, \underbrace{δ, …, δ}_{n_0}).
\end{equation*}
\end{definition}


A functor between two $ A_∞ $-categories is a mapping which matches the products of the two categories:

\begin{definition}
\label{def:2Bprelim-ainfty-functor}
Let $ \cat C $ and $ \cat D $ be $ A_∞ $-categories. Then a \emph{functor} $ F: \cat C → \cat D $ of $ A_∞ $-categories consists of a map $ F: \Ob(\cat C) → \Ob(\cat D) $ together with for every $ k ≥ 1 $ a degree $ 1-k $ multilinear map
\begin{equation*}
F^k: \Hom_{\cat C} (X_k, X_{k+1}) ¤ … ¤ \Hom_{\cat C} (X_1, X_2) → \Hom_{\cat C} (FX_1, FX_{k+1})
\end{equation*}
such that the $ A_∞ $-functor relations hold:
\begin{multline*}
\sum_{0 ≤ j < i ≤ k} (-1)^{‖a_j‖ + … + ‖a_1‖} F(a_k, …, a_{i+1}, μ(a_i, …, a_{j+1}), a_j, …, a_1) \\
= \sum_{\substack{l ≥ 0 \\ 1 = j_1 < … < j_l ≤ k}} μ(F(a_k, …, a_{j_l}), …, F(…, a_{j_2}), F(…, a_{j_1})).
\end{multline*}
The functor $ F $ is an \emph{isomorphism} if $ F: \Ob(\cat C) → \Ob(\cat D) $ is a bijection and $ F^1: \Hom_{\cat C} (X, Y) → \Hom_{\cat D} (FX, FY) $ is an isomorphism for all $ X, Y ∈ \cat C $. The functor $ F $ is a \emph{quasi-isomorphism} if $ F: \Ob(\cat C) → \Ob(\cat D) $ is a bijection and $ F^1: \Hom_{\cat C} (X, Y) → \Hom_{\cat D} (FX, FY) $ is a quasi-isomorphism of complexes for every $ X, Y ∈ \cat C $.
\end{definition}

\begin{definition}
When $ F: \cat C → \cat D $ and $ G: \cat D → \cat E $ are $ A_∞ $-functors, then their composition is given by $ GF: \Ob(\cat C) → \Ob(\cat E) $ on objects and
\begin{equation*}
(GF)(a_k, …, a_1) = \sum G(F(a_k, …), …, F(…, a_1)).
\end{equation*}
\end{definition}


Let us recall minimal models and their notation as follows:

\begin{definition}
An $ A_∞ $-category $ \cat C $ is \emph{minimal} if $ μ^1_{\cat C} = 0 $. A \emph{minimal model} of $ \cat C $ is any minimal $ A_∞ $-category $ \cat D $ together with a quasi-isomorphism $ F: \cat D → \cat C $. A minimal model of $ \cat C $ is generically denoted $ \H\cat C $.
\end{definition}

By the famous Kadeishvili theorem, every $ A_∞ $-category has a minimal model. In fact, a minimal model can be constructed semi-explicitly by sums over trees.

\subsection{Deformations of $ A_∞ $-categories}
\label{sec:2Bainfty-defo}
In this section, we recall deformations of $ A_∞ $-categories. We follow \papertwoA\ where also more detail can be found. We start by recalling completed tensor products. Then we recall $ A_∞ $-deformations and their functors. We comment very briefly on the construction of the twisted completion and minimal models for $ A_∞ $-deformations from \papertwoA.

We recall now completed tensor products $ B \htensor X $ with $ B $ a local ring and $ X $ a vector space. The letter $ B $ will always denote a local ring with extra properties. We have decided to give this a name:

\begin{definition}
\label{def:2Bainfty-defo-base}
A \emph{deformation base} is a complete local Noetherian unital $ ℂ $-algebra $ B $ with residue field $ B/\mathfrak{m} = ℂ $. The maximal ideal is always denoted $ \mathfrak{m} $.
\end{definition}

\begin{remark}
By the Cohen structure theorem, every deformation base is of the form $ ℂ⟦x_1, …, x_n⟧ / I $ with $ I $ denoting some ideal.
\end{remark}

If $ X $ is a vector spaces, then $ B \htensor X = \lim (B/\mathfrak{m}^k \tensor X) $ denotes the completed tensor product over $ ℂ $. For simplicity, we write $ \mathfrak{m}^k X $ to denote the infinitesimal part $ \mathfrak{m}^k X = \mathfrak{m}^k \htensor X ⊂ B \htensor X $. Recall that $ B \htensor X $ is a $ B $-module and comes with the $ \mathfrak{m} $-adic topology, which turns $ B \htensor X $ into a sequential Hausdorff space. For convenience, we may from time to time use expressions like $ x = \landau(\mathfrak{m}^k) $ to indicate $ x ∈ \mathfrak{m}^k X $.

\begin{definition}
A map $ φ: B \htensor X → B \htensor Y $ is \emph{continuous} if it is continuous with respect to the $ \mathfrak{m} $-adic topologies. A map $ φ: (B \htensor X_k) ¤ … ¤ (B \htensor X_1) → B \htensor Y $ is \emph{continuous} if for every $ 1 ≤ i ≤ k $ and every sequence of elements $ x_1, …, \hat x_i, …, x_k $ the map
\begin{equation*}
μ(x_k, …, -, …, x_1): B \htensor X_i → B \htensor Y
\end{equation*}
is continuous.
\end{definition}

\begin{remark}
Every element in $ B \htensor X $ can be written as a series $ \sum_{i = 0}^∞ m_i x_i $. Here $ m_i $ is a sequence of elements $ m_i ∈ \mathfrak{m}^{→∞} $ and $ x_i $ is a sequence of elements $ x_i ∈ X $. We have used the notation $ m_i ∈ \mathfrak{m}^{→∞} $ to indicate that $ m_i ∈ \mathfrak{m}^{k_i} $ for some sequence $ (k_i) ⊂ ℕ $ with $ k_i → ∞ $.
\end{remark}

\begin{remark}
\label{rem:2Bprelim-ainfty-continuityautomatic}
Every $ B $-linear map $ B \htensor X → B \htensor Y $ is automatically continuous (see \papertwoA~for the argument), so is every every $ B $-multilinear map $ (B \htensor X_k) ¤ … ¤ (B \htensor X_1) → B \htensor Y $. Linear maps $ X → B \htensor Y $ can be uniquely extended to $ B $-linear maps $ B \htensor X → B \htensor Y $ and multilinear maps $ X_k ¤ … ¤ X_1 → B \htensor Y $ can be uniquely extended to $ B $-multilinear maps $ (B \htensor X_k) ¤ … ¤ (B \htensor X_1) → B \htensor Y $ (see \papertwoA).
\end{remark}

\begin{remark}
\label{rem:2Bprelim-ainfty-leadingterm}
The \emph{leading term} of a $ B $-linear map $ φ: B \htensor X → B \htensor Y $ is the map $ φ_0: X → Y $ given by the composition $ φ_0 = π φ \restr{X} $, where $ π: B \htensor Y → Y $ denotes the standard projection. If the leading term $ φ_0 $ is injective or surjective, then $ φ $ is injective or surjective itself (see \papertwoA~for the argument).
\end{remark}

We recall now $ A_∞ $-deformations. When $ \cat C $ is an $ A_∞ $-category, the idea is to model its $ A_∞ $-deformations on the collection of enlarged hom spaces $ \{B \htensor \Hom_{\cat C} (X, Y)\}_{X, Y ∈ \cat C} $. Any $ B $-multilinear product on these hom spaces is automatically continuous. Similarly, functors of $ A_∞ $-deformations will be defined as maps between tensor products of the enlarged hom spaces and will be automatically continuous as well.

$ A_\infty $-deformations of $ \cat C $ will always be allowed to have infinitesimal curvature. The reason is that only this way we get a homologically sensible notion: Whenever $ μ_q $ is an (infinitesimally) curved deformation, then $ ν = μ - μ_q $ is a Maurer-Cartan element of the Hochschild DGLA $ \HC(\cat C) $. We comment on this in more detail in \papertwoA.

\begin{definition}
\label{def:2Bprelim-ainfty-defo}
Let $ \cat C $ be an $ A_∞ $ category with products $ μ $ and $ B $ a deformation base. An \emph{$ A_\infty $-deformation} of $ \cat C_q $ of $ \cat C $ consists of
\begin{itemize}
\item The same objects as $ \cat C $,
\item Hom spaces $ \Hom_{\cat C_q} (X, Y) = B \htensor \Hom_{\cat C} (X, Y) $ for $ X, Y ∈ \cat C $,
\item $ B $-multilinear products of degree $ 2 - k $
\begin{equation*}
μ_q^k: \Hom_{\cat C_q} (X_k, X_{k+1}) \tensor … \tensor \Hom_{\cat C_q} (X_1, X_2) → \Hom_{\cat C_q} (X_1, X_{k+1}), ~ k ≥ 1
\end{equation*}
\item Curvature of degree $ 2 $ for every object $ X ∈ \cat C $
\begin{equation*}
μ_{q, X}^0 ∈ \mathfrak{m} \Hom_{\cat C_q}^2 (X, X),
\end{equation*}
\end{itemize}
such that $ μ_q $ reduces to $ μ $ once the maximal ideal $ \mathfrak{m} $ is divided out, and $ μ_q $ satisfies the curved $ A_∞ $ ($ cA_∞ $) relations
\begin{equation*}
\sum_{k ≥ l ≥ m ≥ 0} (-1)^{‖a_m‖ + … + ‖a_1‖} μ_q (a_k, …, μ_q (a_l, …), a_m, …, a_1) = 0.
\end{equation*}
The deformation is \emph{unital} if the deformed higher products still satisfy the unitality axioms
\begin{equation*}
μ^2_q (a, \id_X) = a, ~ μ^2_q (\id_Y, a) = (-1)^{|a|} a, ~ μ^{≥3}_q (…, \id_X, …) = 0.
\end{equation*}
\end{definition}

We are now ready to explain the natural extension of $ A_∞ $-functors to the deformed case.

\begin{definition}
\label{def:2Bprelim-defo-functor}
Let $ \cat C, \cat D $ be two $ A_∞ $-categories and $ \cat C_q, \cat D_q $ deformations. A \emph{functor of deformed $ A_∞ $-categories} consists of a map $ F_q: \Ob(\cat C) → \Ob(\cat D) $ together with for every $ k ≥ 1 $ a $ B $-multilinear degree $ 1-k $ map
\begin{equation*}
F^k_q: \Hom_{\cat C_q} (X_k, X_{k+1}) ¤ … ¤ \Hom_{\cat C_q} (X_1, X_2) → \Hom_{\cat D_q} (F_q X_1, F_q X_{k+1})
\end{equation*}
and infinitesimal curvature $ F^0_{q, X} ∈ \mathfrak{m} \Hom^1_{\cat D} (F_q X, F_q X) $ for every $ X ∈ \cat C $, such that the curved $ A_∞ $-functor relations hold:
\begin{multline*}
\sum_{0 ≤ j ≤ i ≤ k} (-1)^{‖a_j‖ + … + ‖a_1‖} F_q (a_k, …, a_{i+1}, μ_q (a_i, …, a_{j+1}), a_j, …, a_1) \\
= \sum_{\substack{l ≥ 0 \\ 1 = j_1 < … < j_l ≤ k}} μ_q (F_q (a_k, …, a_{j_l}), …, F_q (…, a_{j_2}), F_q (…, a_{j_1})).
\end{multline*}
If $ \cat C_q $ and $ \cat D_q $ are strictly unital, then we say $ F_q $ is strictly unital if $ F_q^1 (\id_X) = \id_{F_q X} $ for every $ X \in \cat C $ and $ F_q^{≥2} (…, \id_X, …) = 0 $.
\end{definition}

\begin{remark}
Note that the functor $ F_q $ itself is allowed to have a curvature component. The first two curved $ A_∞ $-functor relations read
\begin{align*}
F^0_q + F^1_q (μ^0_{\cat C_q, X}) &= μ^1_{\cat D_q} (F^0_{q, X}), \\
F^1_q (μ^1_{\cat C_q} (a)) + (-1)^{‖a‖} F^2_q (μ^0_{\cat C_q, Y}, a) + F^2_q (a, μ^0_{\cat C_q, X}) &= μ^1_{\cat D_q} (F^1_q (a)) + μ^2_{\cat D_q} (F^0_{q, Y}, F^1_q (a)) \\
& \quad + μ^2_{\cat D_q} (F^1_q (a), F^0_{q, X}), \quad ∀ ~ a: X → Y.
\end{align*}
If $ F_q: \cat C_q \to \cat D_q $ is a functor of $ A_\infty $-deformations, then its leading term $ F: \cat C \to \cat D $ is automatically a functor of $ A_\infty $-categories.
\end{remark}

\begin{definition}
\label{def:2Bainfty-quasiiso}
Let $ F_q: \cat C_q → \cat D_q $ be a functor of $ A_∞ $-deformations. Then $ F_q $ is a \emph{quasi-isomorphism} if its leading term $ F: \cat C → \cat D $ is a quasi-isomorphism of $ A_∞ $-categories.
\end{definition}


\emph{Uncurving} refers to the process of removing curvature from a deformation $ \cat C_q $ by means of gauging. More precisely, uncurving refers to finding a functor $ F_q: \cat C_q' → \cat C_q $ where $ \cat C_q' $ is another deformation of $ \cat C $ with possibly less curvature and the functor $ F_q $ is the identity when dividing out the maximal ideal $ \mathfrak{m} $ of the deformation base. For instance, one may regard the \emph{uncurving} $ \cat C_q' $ of $ \cat C_q $ by an element $ r = \{r_X\}_{X ∈ \cat C} $ consisting of $ r_X ∈ \mathfrak{m} \End^1_{\cat C} (X) $ for every $ X ∈ \cat C $. The categories $ \cat C_q' $ and $ \cat C_q $ are related by the functor $ F_q $ of $ A_∞ $-deformations given by
\begin{equation*}
F_q: \cat C_q' \isoto \cat C_q, \quad \text{given by} \quad F_q^0 ≔ r, \quad F_q^1 ≔ \Id, \quad F_q^{≥2} ≔ 0.
\end{equation*}
The curvature of $ \cat C_q' $ is then
\begin{equation*}
μ^0_{\cat C_q', X} = μ^0_{\cat C_q, X} + μ^1_{\cat C_q} (r_X) + μ^2_{\cat C_q} (r_X, r_X) + ….
\end{equation*}
We record the following definition from \papertwoA:
\begin{definition}
Let $ \cat C $ be an $ A_∞ $ category and $ \cat C_q $ a deformation. Let $ X ∈ \cat C $. Then $ X $ is \emph{uncurvable} if there exists an $ r_X ∈ \mathfrak{m} \End^1(X) $ such that
\begin{equation*}
μ^0_X + μ^1_{\cat C_q} (r_X) + μ^2_{\cat C_q} (r_X, r_X) + … = 0.
\end{equation*}
\end{definition}

\begin{remark}
\label{rem:2Bprelim-ainfty-defo-uncurvingqi}
In \papertwoA\ we explain that uncurvability of objects is invariant under quasi-isomorphism. More precisely, if $ X, Y ∈ \cat C $ are quasi-isomorphic and $ X $ is uncurvable in $ \cat C_q $, then also $ Y $ is uncurvable in $ \cat C_q $. We also show that functors of $ A_∞ $-deformations send uncurvable objects to uncurvable objects. More precisely, if $ F_q: \cat C_q → \cat D_q $ is a functor of $ A_∞ $-deformations and $ X ∈ \cat C $ is uncurvable, then $ F_q (X) $ is uncurvable.
\end{remark}


We recall in \autoref{def:2B-ainfty-defo-twisted} that a deformation $ \cat C_q $ has a twisted completion $ \Tw\cat C_q $. This category $ \Tw\cat C_q $ is a deformation of $ \Tw\cat C $. Its objects are defined in terms of twisted differentials as well, but the twisted differentials do not satisfy the Maurer-Cartan equation with respect to the deformed product $ μ_{\cat C_q} $. Instead, the failure to satisfy the Maurer-Cartan equation is captured in the object's curvature. We recall the definition from \papertwoA\ as follows:

\begin{definition}
\label{def:2B-ainfty-defo-twisted}
Let $ \cat C $ be an $ A_∞ $ category with products $ μ_{\cat C} $ and $ \cat C_q $ a deformation with products $ μ_{\cat C_q} $. Then the \emph{twisted completion} $ \Tw \cat C_q $ is the (possibly curved) deformation of $ \Tw \cat C $ given by the deformed products
\begin{equation*}
μ_{\Tw \cat C_q}^k (α_k, …, α_1) = \sum_{n_0, …, n_k ≥ 0} μ_{\Add \cat C_q} (\underbrace{δ, …, δ}_{n_k}, α_k, …, α_1, \underbrace{δ, …, δ}_{n_0}).
\end{equation*}
\end{definition}

\begin{remark}
\label{rem:2B-ainfty-defo-liberaltwisted}
It is possible to define a variant $ \Tw'\cat C_q $ of the twisted completion of $ A_∞ $-deformations by allowing additional infinitesimal entries anywhere in the $ δ $-matrix. The objects of $ \Tw'\cat C_q $ shall be pairs
\begin{equation*}
(X, δ = δ_0 + δ'), \quad X ∈ \Add\cat C, \quad δ_0 ∈ \Hom_{\cat C}^1 (X, X), \quad δ' ∈ \mathfrak{m} \Hom_{\cat C}^1 (X, X).
\end{equation*}
Here we require only the leading part $ δ_0 $ to be upper triangular and satisfy the Maurer-Cartan equation with respect to $ μ_{\cat C} $. The infinitesimal part $ δ' $ can also lie below the diagonal. The category $ \Tw'\cat C_q $ is strictly speaking not a deformation of $ \Tw\cat C_q $ since it has more objects. For more insight, we refer to \papertwoA.
\end{remark}

\emph{Minimal models} of $ A_∞ $-deformations exist and we review this topic in \autoref{sec:2Bkadeishvili}.

\subsection{The deformed Kadeishvili theorem}
\label{sec:2Bkadeishvili}
In this section we recall the classical Kadeishvili construction and a special case of our deformed Kadeishvili construction from \papertwoA. The starting point is an $ A_∞ $-deformation $ \cat C_q $ of an $ A_∞ $-category $ \cat C $. Due to curvature and the failure of the differential to square to zero, it is initially unclear how to define the notion of minimal models and how to construct a minimal model for the given deformation $ \cat C_q $ explicitly. In the present section, we review our definition and construction from \papertwoA:

\begin{center}
\begin{tikzpicture}
\path (0, 0) node (A) {$ A_∞ $-deformation $ \cat C_q $} (8, 0) node[align=center] (B) {Minimal model $ \H\cat C_q $};
\path[draw, decorate, decoration={snake, amplitude=0.2em, post length=0.5em}, ->] ($ (A.east)!0.2!(B.west) $) to ($ (A.east)!0.8!(B.west) $);
\end{tikzpicture}
\end{center}
We do not attempt to review the deformed Kadeishvili construction in its full generality. Rather we focus on the special case where the curvature vanishes and an additional condition “$ D = 0 $” holds. In this case, the construction simplifies substantially. For the purposes of \autoref{sec:deformed} until \ref{sec:subdisk}, this special case suffices. For the purposes of \autoref{sec:sphere}, the special case is insufficient and we have to refer to \papertwoA\ for documentation of the general case.

We start by reviewing the definition of minimal models for $ A_∞ $-deformations from \papertwoA:

\begin{definition}
\label{def:2Bkadeishvili-minmodel}
Let $ \cat C $ be an $ A_∞ $-category, $ \cat C_q $ a deformation of $ \cat C $, and $ \H\cat C $ a minimal model for $ \cat C $. Then a \emph{minimal model} $ \H\cat C_q $ of $ \cat C_q $ is any $ A_∞ $-deformation of $ \H\cat C $ such that there exists a quasi-isomorphism $ F_q: \H\cat C_q → \cat C_q $.
\end{definition}

\begin{remark}
As we explain in \papertwoA, every $ A_∞ $-deformation has a minimal model. While a minimal model $ \H\cat C_q $ is a deformation of $ \H\cat C $, its differential and curvature need not vanish. Instead, $ \H\cat C_q $ carries an infinitesimal residue differential and curvature.
\end{remark}

Our second step is to recall the notion of homological splittings:

\begin{definition}
Let $ \cat C $ be an $ A_∞ $-category. Then a \emph{homological splitting} of $ \cat C $ consists of a direct sum decomposition
\begin{equation*}
\Hom_{\cat C} (X, Y) = H(X, Y) ⊕ I(X, Y) ⊕ R(X, Y), \quad ∀X, Y ∈ \cat C
\end{equation*}
for all its hom spaces, such that
\begin{equation*}
I(X, Y) = \Im(μ^1), \quad \Ker(μ^1) = H(X, Y) ⊕ I(X, Y), \quad ∀X, Y ∈ \cat C.
\end{equation*}
We frequently denote a homological splitting of $ \cat C $ simply by the letters $ H ⊕ I ⊕ R $, the dependence on $ X, Y ∈ \cat C $ understood implicitly.
\end{definition}

Given a category $ \cat C $, one obtains a homological splitting by choosing $ H $ as a space of cocycles that represents the cohomology of the hom complexes. One then chooses $ R $ as a complement to $ H $ in $ \Ker(μ^1) $. The notation $ I $ is simply a shorthand for the image of the differential. In terms of the direct sum decomposition $ \Hom_{\cat C} (X, Y) = H(X, Y) ⊕ I(X, Y) ⊕ R(X, Y) $, the differential reads
\begin{equation*}
μ^1 = \begin{pmatrix} 0 & 0 & 0 \\ 0 & 0 & * \\ 0 & 0 & 0 \end{pmatrix}.
\end{equation*}

\begin{remark}
Given a homological splitting $ H ⊕ I ⊕ R $, we write an element $ x ∈ \Hom_{\cat C} (X, Y) $ typically in tuple form as
\begin{equation*}
x = (h, μ^1 (r'), r), \quad \text{with } h ∈ H, ~ r' ∈ R, ~ r ∈ R.
\end{equation*}
\end{remark}

\begin{definition}
Let $ \cat C $ be an $ A_∞ $-category and $ H ⊕ I ⊕ R $ a homological splitting. The \emph{codifferential} is the map
\begin{align*}
h: \Hom_{\cat C} (X, Y) &→ R(X, Y), \\
(h, μ^1 (r'), r) &↦ r', \quad h ∈ H(X, Y), ~ r' ∈ R(X, Y), ~ r ∈ R(X, Y).
\end{align*}
The \emph{projection to cohomology} is the map
\begin{equation*}
π: \Hom_{\cat C} (X, Y) = H(X, Y) ⊕ I(X, Y) ⊕ R(X, Y) \twoheadrightarrow H(X, Y).
\end{equation*}
\end{definition}

At this stage we are ready to recall the classical Kadeishvili construction. Its intention is to construct a minimal model $ \H\cat C $ explicitly. The starting point for the construction is a homological splitting $ H ⊕ I ⊕ R $ of $ \cat C $. The result of the construction is an $ A_∞ $-structure on $ \{H(X, Y)\}_{X, Y ∈ \cat C} $ which can also be interpreted as an $ A_∞ $-structure on $ \{\H\Hom_{\cat C} (X, Y)\}_{X, Y ∈ \cat C} $, since $ H(X, Y) $ and $ \H\Hom_{\cat C} (X, Y) $ are isomorphic as graded vector spaces through the composition $ H(X, Y) \embeds \Ker(μ^1) \twoheadrightarrow \H\Hom_{\cat C} (X, Y) $. Specifically, the $ A_∞ $-structure on $ H(X, Y) $ is defined via trees. We fix terminology as follows:

\begin{definition}
\label{def:2Bkadeishvili-classical-treeshape}
A \emph{Kadeishvili tree shape} $ T $ is a rooted planar tree with $ n ≥ 2 $ leaves whose non-leaf nodes all have at least 2 children. A node in $ T $ is \emph{internal} if it is not a leaf and not the root. The number of internal nodes in $ T $ is denoted $ N_T $. We denote by $ \mathcal{T}_n $ the set of all Kadeishvili tree shapes with $ n $ leaves.

A \emph{Kadeishvili π-tree} $ (T, h_1, …, h_n) $ is a Kadeishvili tree shape $ T ∈ \mathcal{T}_n $ with $ n ≥ 2 $ leaves, together with a sequence $ h_1, …, h_n $ of cohomology elements $ h_i ∈ H(X_i, X_{i+1}) $. Decorate the leaves by $ h_1, …, h_n $ in sequence. Decorate every non-root node with the operation $ hμ $ and the root with the operation $ πμ $. Then the \emph{result} $ \Res(T, h_1, …, h_n) ∈ H(X_1, X_{n+1}) $ of the Kadeishvili π-tree is the result obtained by evaluating the tree from leaves to the root, according to the decorations.
\end{definition}

The construction of the $ A_∞ $-product $ μ_H $ on $ \H\cat C $ can be summarized as follows: Let $ h_1, …, h_n $ be cohomology elements with $ h_i ∈ H(X_i, X_{i+1}) $. Then their higher product is defined as
\begin{equation*}
μ_{\H\cat C} (h_n, …, h_1) = \sum_{T ∈ \mathcal{T}_n} (-1)^{N_T} \Res(T, h_1, …, h_n) ∈ H(X_1, X_{n+1}).
\end{equation*}

\begin{theorem}[{Kadeishvili, \cite[Chapter 6, 3.3.2]{Kontsevich-Soibelman}}]
Let $ \cat C $ be an $ A_∞ $-category. Then the products $ μ_{\H\cat C} $ define an $ A_∞ $-structure on $ \H\cat C $. Equipped with this structure, $ \H\cat C $ is a minimal model for $ \cat C $.
\end{theorem}

In the remainder of the section, we recall the deformed Kadeishvili construction from \papertwoA\ in a special case. The starting point is an $ A_∞ $-category $ \cat C $ with a deformation $ \cat C_q $ such that the curvature of $ \cat C_q $ vanishes and the deformed differential satisfies $ μ^1_q (H(X, Y)) ⊂ μ^1_q (B \htensor R(X, Y)) $. We termed this the “$ D = 0 $” case in \papertwoA\ and it turns out that in this case both the curvature $ μ^0_{\H\cat C_q} $ and the differential $ μ^1_{\H\cat C_q} $ of the minimal model vanish.

\begin{definition}
\label{def:2Bkadeishvili-deformed-counterpart}
Let $ \cat C $ be an $ A_∞ $-category and $ \cat C_q $ a deformation. Let $ H ⊕ I ⊕ R $ be a homological splitting of $ \cat C $. Assume $ \cat C_q $ is curvature-free with $ μ^1_q (B \htensor H(X, Y)) ⊂ μ^1_q (B \htensor R(X, Y)) $. Then the map $ μ^1_q $ is injective on $ B \htensor R(X, Y) $:
\begin{equation*}
μ^1_q \restr{B \htensor R(X, Y)}: B \htensor R(X, Y) \isoto μ^1_q (B \htensor R(X, Y)).
\end{equation*}
For $ h ∈ B \htensor H(X, Y) $ denote by $ ε_h ∈ B \htensor R(X, Y) $ the unique element such that $ μ^1_q (h - ε_h) = 0 $. The \emph{deformed counterpart} of $ h ∈ B \htensor H(X, Y) $ is the element $ h - ε_h $. Put
\begin{equation*}
H_q (X, Y) = \{h - ε_h \running h ∈ B \htensor H(X, Y)\}.
\end{equation*}
The \emph{deformed splitting} of $ \cat C_q $ is the decomposition of all hom spaces of $ \cat C_q $ as
\begin{equation*}
\Hom_{\cat C_q} (X, Y) = H_q (X, Y) ⊕ μ^1_q (B \htensor R(X, Y)) ⊕ (B \htensor R(X, Y)).
\end{equation*}
The correspondence between $ H_q (X, Y) $ and $ B \htensor H (X, Y) $ is denoted
\begin{align*}
φ: H_q (X, Y) &\verylongisoto B \htensor H(X, Y), \\
h - ε_h &\verylongmapsto h.
\end{align*}
The \emph{deformed codifferential} of $ \cat C_q $ is the $ B $-linear map
\begin{equation*}
h_q = (μ^1_q \restr{B \htensor R(X, Y)})^{-1}: μ^1_q (B \htensor R(X, Y)) \verylongto B \htensor R(X, Y).
\end{equation*}
The \emph{deformed projection} of $ \cat C_q $ is the $ B $-linear map
\begin{equation*}
π_q: \Hom_{\cat C_q} (X, Y) = H_q (X, Y) ⊕ μ^1_q (B \htensor R(X, Y)) ⊕ (B \htensor R(X, Y)) → H_q (X, Y).
\end{equation*}
\end{definition}

The meaning of the condition $ μ^1_q (B \htensor H(X, Y)) ⊂ μ^1_q (B \htensor R(X, Y)) $ becomes evident when we regard the shape of $ μ^1_q $ with respect to the deformed splitting:

\begin{lemma}[\papertwoA]
Let $ \cat C $ be an $ A_∞ $-category and $ \cat C_q $ a deformation. Let $ H ⊕ I ⊕ R $ be a homological splitting of $ \cat C $. Assume $ \cat C_q $ is curvature-free with $ μ^1_q (B \htensor H(X, Y)) ⊂ μ^1_q (B \htensor R(X, Y)) $. With respect to the deformed splitting of $ \cat C_q $, the differential $ μ^1_q $ takes the shape
\begin{equation*}
μ^1_q = \pmat{0 & 0 & 0 \\ 0 & 0 & * \\ 0 & 0 & 0}.
\end{equation*}
\end{lemma}

\begin{definition}
Let $ T $ be a Kadeishvili tree shape and $ h_1, …, h_k $ be compatible morphisms in $ \cat C_q $. Decorate the tree $ T $ by putting $ φ^{-1} (h_i) $ on the leaves, $ h_q μ_{\cat C_q} $ on every internal node and $ φ π_q μ_{\cat C_q} $ on the root. Define the $ \Res_q (T, h_1, …, h_k) $ as the result of the evaluation.
\end{definition}

We are now ready to describe the minimal model $ \H\cat C_q $. As a deformation of $ \H\cat C $, it carries the same objects as $ \cat C $ and its hom spaces are $ \Hom_{\H\cat C_q} (X, Y) = B \htensor \H\Hom_{\cat C} (X, Y) $. The deformed $ A_∞ $-structure takes the following shape:

\begin{theorem}[\papertwoA]
\label{th:2Bkadeishvili-th}
Let $ \cat C $ be an $ A_∞ $-category and $ \cat C_q $ a deformation. Let $ H ⊕ I ⊕ R $ be a homological splitting of $ \cat C $. Assume $ \cat C_q $ is curvature-free with $ μ^1_q (B \htensor H(X, Y)) ⊂ μ^1_q (B \htensor R(X, Y)) $. Then a minimal model $ \H\cat C_q $ is obtained by setting
\begin{equation*}
μ^1_{\H\cat C_q} = μ^0_{\H\cat C_q} = 0
\end{equation*}
and
\begin{equation*}
μ_{\H\cat C_q}^{k \geq 2} (h_k, …, h_1) = \sum_{T ∈ \mathcal{T}_k} (-1)^{N_T} \Res_q (T, h_1, …, h_k).
\end{equation*}
\end{theorem}

\section{Preliminaries on gentle algebras}
\label{sec:prelim-gtl}
In this section, we concisely recapitulate background on gentle algebras in order to bring the reader into touch with the relevant tools of this paper. We provide definitions of all preliminaries and explain alternative points of view on them. In particular, we will explain how every definition is used in the paper. We follow mostly \cite{Bocklandt} and \paperone.

\subsection{Punctured surfaces}
Punctured surfaces belong to the family of two-dimensional oriented manifolds, while at the same time facilitating singular behavior at the punctures.

\begin{definition}
\label{def:prelim-punctured-def}
A \emph{punctured surface} is a closed oriented surface $ S $ with a finite set of punctures $ M ⊂ S $. We assume that $ |M| ≥ 1 $, or $ |M| ≥ 3 $ if $ S $ is a sphere.
\end{definition}

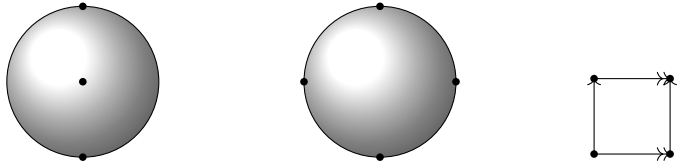
\begin{figure}
\center
\begin{subfigure}{0.2\linewidth}
\begin{tikzpicture}
\path[draw, shade, ball color=white] (0, 0) circle[radius=1];
\path[fill] (0, 0) circle[radius=0.05] (0, 1) circle[radius=0.05] (0, -1) circle[radius=0.05];
\end{tikzpicture}
\end{subfigure}
\begin{subfigure}{0.2\linewidth}
\center
\begin{tikzpicture}
\path[draw, shade, ball color=white] (0, 0) circle[radius=1];
\path[fill] (1, 0) circle[radius=0.05];
\path[fill] (-1, 0) circle[radius=0.05];
\path[fill] (0, 1) circle[radius=0.05];
\path[fill] (0, -1) circle[radius=0.05];
\end{tikzpicture}
\end{subfigure}
\begin{subfigure}{0.2\linewidth}
\center
\begin{tikzpicture}
\path[draw, ->>] (0, 0) to (1, 0);
\path[draw, ->] (1, 0) to (1, 1);
\path[draw, ->>] (0, 1) to (1, 1);
\path[draw, ->] (0, 0) to (0, 1);
\foreach \i in {0, 1} \foreach \j in {0, 1} \path[fill] (\i, \j) circle[radius=0.05];
\end{tikzpicture}
\end{subfigure}
\caption{The three-punctured sphere, the four-punctured sphere and the one-punctured torus}
\label{fig:prelim-punctured}
\end{figure}

A selection of popular punctured surfaces are depicted in \autoref{fig:prelim-punctured}. The condition $ |M| ≥ 1 $ and $ |M| ≥ 3 $ are merely cosmetic and will be explained in \autoref{sec:prelim-arc-systems}.

\begin{remark}
A punctured surface can alternatively be interpreted as a surface with $ S^1 $ boundaries: Let $ (S, M) $ be a punctured surface and regard one puncture $ q ∈ M $. The surface around $ q $ looks like a punctured disk. Now interpret the punctured disk as an infinitely long cylinder, glued to the rest of the surface. Cut off the cylinder at some distance. We obtain a surface with $ S^1 $ boundaries, which we interpret as markings. In other words, we have a marked surface with only $ S^1 $ boundaries in the sense of \cite{HKK}. For instance, cutting away disks around the punctures in the three-punctured sphere, we obtain the popular pair of pants surface.
\end{remark}

\subsection{Arc systems}
\label{sec:prelim-arc-systems}
In this section, we recall the notion of arcs and arc systems on punctured surfaces. We recall what it means for an arc system to be full, and explain how it cuts the surface into polygons. We fix some terminology regarding polygons, in particular the notion of a polygon's interior angles.

\begin{definition}
\label{def:prelim-arcsys-def}
Let $ (S, M) $ be a punctured surface. An \emph{arc} in $ S $ is a not necessarily closed curve $ γ: [0, 1] → S $ running from one puncture to another. An \emph{arc system} $ \cA $ on a punctured surface is a finite collection of arcs which meet only at the set $ M $ of punctures. Intersections and self-intersections are not allowed. The arc system satisfies the \emph{no monogons or digons} condition \emph{[NMD]} if
\begin{itemize}
\item No arc is a contractible loop in $ S \setminus M $.
\item No pair of distinct arcs is homotopic in $ S \setminus M $.
\end{itemize}
The arc system satisfies the \emph{no monogons or digons in the closed surface} condition \emph{[NMDC]} if
\begin{itemize}
\item No arc is a contractible loop in $ S $.
\item No pair of distinct arcs is homotopic in $ S $.
\end{itemize}
\end{definition}

\begin{figure}
\centering
\begin{subfigure}[b]{0.26\linewidth}
\centering
\begin{tikzpicture}
\path[shade, ball color=white] (0, 0) circle[radius=1];
\path[draw] (0, -1) arc(-90:90:1);
\path[draw] (0, 1) arc(90:270:1);
\path[fill] (0, 0) circle[radius=0.05] (0, 1) circle[radius=0.05] (0, -1) circle[radius=0.05];
\path[draw, thick, ->] (0, -1) -- (0, 0);
\path[draw, thick, ->] (0, 0) -- (0, 1);
\end{tikzpicture}
\caption{Three-punctured sphere}
\label{fig:prelim-arcsys-3p}
\end{subfigure}
\begin{subfigure}[b]{0.26\linewidth}
\centering
\begin{tikzpicture}
\path[draw, shade, ball color=white] (0, 0) circle[radius=1];
\path[draw, thick, ->] (1, 0) arc(0:90:1);
\path[draw, thick, ->] (1, 0) to[bend left=5] (-1, 0);
\path[fill] (1, 0) circle[radius=0.05];
\path[fill] (-1, 0) circle[radius=0.05];
\path[fill] (0, 1) circle[radius=0.05];
\path[fill] (0, -1) circle[radius=0.05];
\end{tikzpicture}
\caption{Four-punctured sphere}
\label{fig:prelim-arcsys-4p}
\end{subfigure}
\begin{subfigure}[b]{0.26\linewidth}
\centering
\begin{tikzpicture}
\path[draw, thick, ->] (0, 0) to node[below] {$ a $} (1, 0);
\path[draw, thick, ->] (1, 0) to node[right] {$ b $} (1, 1);
\path[draw, thick, ->] (0, 1) to node[above] {$ a $} (1, 1);
\path[draw, thick, ->] (0, 0) to node[left] {$ b $} (0, 1);
\foreach \i in {0, 1} \foreach \j in {0, 1} \path[fill] (\i, \j) circle[radius=0.05];
\end{tikzpicture}
\caption{One-punctured torus}
\label{fig:prelim-arcsys-1p-torus}
\end{subfigure}
\begin{subfigure}[b]{0.2\linewidth}
\centering
\begin{tikzpicture}
\path[draw, thick] (0, 0) to[out=60, in=0] (0, 2) to[out=180, in=120] (0, 0);
\path[fill] (0, 0) circle[radius=0.05];
\end{tikzpicture}
\caption{Not [NMD]}
\label{fig:prelim-arcsys-loop}
\end{subfigure}
\begin{subfigure}[b]{0.2\linewidth}
\centering
\begin{tikzpicture}
\path[draw, thick] (0, 0) to[out=60, in=300] (0, 2);
\path[draw, thick] (0, 0) to[out=120, in=240] (0, 2);
\path[fill] (0, 0) circle[radius=0.05] (0, 2) circle[radius=0.05];
\end{tikzpicture}
\caption{Not [NMD]}
\label{fig:prelim-arcsys-homotopic}
\end{subfigure}
\begin{subfigure}[b]{0.2\linewidth}
\centering
\begin{tikzpicture}
\path[draw, thick] (0, 0) -- (0, 1);
\path[draw, thick] (0, 0) to[out=60, in=0] (0, 2) to[out=180, in=120] (0, 0);
\path[fill] (0, 0) circle[radius=0.05] (0, 1) circle[radius=0.05];
\end{tikzpicture}
\caption{Not [NMDC]}
\label{fig:prelim-arcsys-loop-crossing}
\end{subfigure}
\begin{subfigure}[b]{0.26\linewidth}
\centering
\begin{tikzpicture}
\path[draw, thick] (0, 0) -- (0, 1);
\path[draw, thick] (0, 0) to[out=60, in=300] (0, 2);
\path[draw, thick] (0, 0) to[out=120, in=240] (0, 2);
\path[fill] (0, 0) circle[radius=0.05] (0, 1) circle[radius=0.05] (0, 2) circle[radius=0.05];
\end{tikzpicture}
\caption{Not [NMDC]}
\label{fig:prelim-arcsys-homotopic-crossing}
\end{subfigure}
\begin{subfigure}[b]{0.3\linewidth}
\centering
\begin{tikzpicture}
\path[draw, thick] (0, 0) -- (0, 1);
\path[draw, thick] (0, 0) -- (1, 1) -- (0, 2) -- (-1, 1) -- cycle;
\path[fill] (0, 0) circle[radius=0.05] (0, 1) circle[radius=0.05] (1, 1) circle[radius=0.05] (0, 2) circle[radius=0.05] (-1, 1) circle[radius=0.05];
\end{tikzpicture}
\caption{Allowed for [NMDC]}
\label{fig:prelim-arcsys-allowed}
\end{subfigure}
\caption{Arc systems and their properties}
\label{fig:prelim-arcs}
\end{figure}
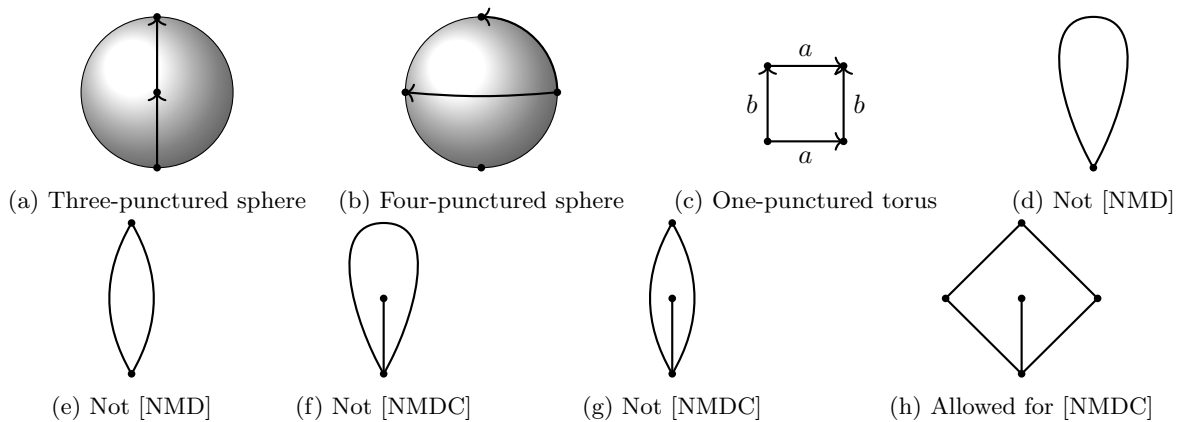

\begin{example}
In \autoref{fig:prelim-arcs}, we have depicted a few arbitrary arc systems on the three- and four-punctured sphere and one-punctured torus. The drawn three-punctured sphere has the north and south pole marked, as well as a point on the equator lying on the front half of the sphere. The four-punctured sphere has the north and south pole, as well as a point in the far east and far west marked.
\end{example}

\begin{example}
The arc system of \autoref{fig:prelim-arcsys-3p} consists of two half meridians lying in the frontal hemisphere. The arc system of \autoref{fig:prelim-arcsys-4p} consists of the frontal half of the equator and a northern half of a meridian. The one-punctured torus \autoref{fig:prelim-arcsys-1p-torus} is drawn as a gluing diagram. The arc system consists of the two standard generators of the torus.
\end{example}

\begin{remark}
The configurations banned by [NMD] are depicted in \autoref{fig:prelim-arcsys-loop} and \ref{fig:prelim-arcsys-homotopic}. The reason to ban these is that the definition of the $ A_∞ $-structure on the gentle algebras becomes a lot easier, avoiding a so-called monogon or digon rule. This makes checking the $ A_∞ $-axioms also more tractable. With the [NMDC] condition, we go a step further and ban also monogons and digons in the closed surface $ S $. More concretely, we ban loops which become contractible when the punctures are filled. Such a banned configuration is depicted in \autoref{fig:prelim-arcsys-loop-crossing}. Similarly, we ban pairs of homotopic arcs, the homotopy being allowed to cross punctures. Such a banned configuration is depicted in \autoref{fig:prelim-arcsys-homotopic-crossing}. The purpose of the [NMDC] condition is to avoid the monogon and digon rule also for the deformed gentle algebras.
\end{remark}

\begin{definition}
\label{def:prelim-arcsys-full}
An arc system is \emph{full} if it cuts the surface into contractible pieces. These pieces are the \emph{faces} or \emph{polygons} of the arc system.
\end{definition}

In other words, an arc system is full if its complement consists of a disjoint union of topological disks. We usually refer to these pieces as polygons to highlight that they are bounded by arcs of the arc system.

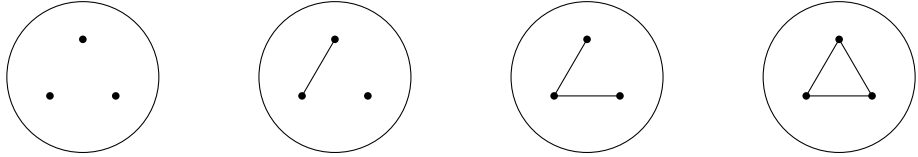
\begin{figure}
\centering
\begin{subfigure}[b]{0.2\linewidth}
\centering
\begin{tikzpicture}
\path[draw] (0, 0) circle[radius=1];
\path[fill] (0, 0) ++(90:0.5) coordinate (A) circle[radius=0.05];
\path[fill] (0, 0) ++(210:0.5) coordinate (B) circle[radius=0.05];
\path[fill] (0, 0) ++(330:0.5) coordinate (C) circle[radius=0.05];
\end{tikzpicture}
\end{subfigure}
\begin{subfigure}[b]{0.2\linewidth}
\centering
\begin{tikzpicture}
\path[draw] (0, 0) circle[radius=1];
\path[fill] (0, 0) ++(90:0.5) coordinate (A) circle[radius=0.05];
\path[fill] (0, 0) ++(210:0.5) coordinate (B) circle[radius=0.05];
\path[fill] (0, 0) ++(330:0.5) coordinate (C) circle[radius=0.05];
\path[draw] (A) -- (B);
\end{tikzpicture}
\end{subfigure}
\begin{subfigure}[b]{0.2\linewidth}
\centering
\begin{tikzpicture}
\path[draw] (0, 0) circle[radius=1];
\path[fill] (0, 0) ++(90:0.5) coordinate (A) circle[radius=0.05];
\path[fill] (0, 0) ++(210:0.5) coordinate (B) circle[radius=0.05];
\path[fill] (0, 0) ++(330:0.5) coordinate (C) circle[radius=0.05];
\path[draw] (A) -- (B) -- (C);
\end{tikzpicture}
\end{subfigure}
\begin{subfigure}[b]{0.2\linewidth}
\centering
\begin{tikzpicture}
\path[draw] (0, 0) circle[radius=1];
\path[fill] (0, 0) ++(90:0.5) coordinate (A) circle[radius=0.05];
\path[fill] (0, 0) ++(210:0.5) coordinate (B) circle[radius=0.05];
\path[fill] (0, 0) ++(330:0.5) coordinate (C) circle[radius=0.05];
\path[draw] (A) -- (B) -- (C) -- (A);
\end{tikzpicture}
\end{subfigure}
\caption{Arc systems on the three-punctured sphere}
\label{fig:prelim-3p-arc-systems}
\end{figure}

\begin{example}
Of the three arc systems presented in \autoref{fig:prelim-arcs}, only \ref{fig:prelim-arcsys-3p} and \ref{fig:prelim-arcsys-1p-torus} are full. Arc system \ref{fig:prelim-arcsys-4p} is not full, because the complement of the arcs is a topological disk with a puncture in its interior, the south pole. Removing any arc from \ref{fig:prelim-arcsys-3p} or \ref{fig:prelim-arcsys-1p-torus} also leads to a non-full arc system. Additional arcs may however be added to make or keep the arc system full. For example, in \ref{fig:prelim-arcsys-3p} one may add the full equator as arc and in \ref{fig:prelim-arcsys-1p-torus} one may add any diagonal as arc, but not both. All possible types of arc systems with [NMDC] on the three-punctured sphere are depicted in \autoref{fig:prelim-3p-arc-systems}. In these figures, the directions of the arcs is arbitrary. Only the arc systems in the third and fourth picture are full.
\end{example}

\begin{figure}
\centering
\begin{subfigure}[b]{0.3\linewidth}
\centering
\begin{tikzpicture}
\path[draw] (0, 0) -- ++(right:1) coordinate[pos=0.3] (1-start) coordinate[pos=0.7] (6-end) -- ++(60:1) coordinate[pos=0.3] (6-start) coordinate[pos=0.7] (5-end) -- ++(120:1) coordinate[pos=0.3] (5-start) coordinate[pos=0.7] (4-end) -- ++(left:1) coordinate[pos=0.3] (4-start) coordinate[pos=0.7] (3-end) -- ++(240:1) coordinate[pos=0.3] (3-start) coordinate[pos=0.7] (2-end) -- ++(300:1) coordinate[pos=0.3] (2-start) coordinate[pos=0.7] (1-end);
\foreach \i in {1, 2, 3, 4, 5, 6} {
\path[draw, ->, bend right=60] (\i-start) to (\i-end);
\path (\i-start) -- (\i-end) node[midway] {$ \i $};
};
\end{tikzpicture}
\caption{Order of angles}
\label{fig:prelim-clockwise}
\end{subfigure}
\begin{subfigure}[b]{0.3\linewidth}
\centering
\begin{tikzpicture}
\path[draw, ->] (0, 0) to[bend right] node[midway, right] {$ a $} (0, 1);
\path[draw, ->] (0, 0) to[bend left] node[midway, left] {$ a $} (0, 1);
\path[fill] (0, 0) circle[radius=0.05];
\path[fill] (0, 1) circle[radius=0.05];
\end{tikzpicture}
\caption{Twice the same arc}
\label{fig:prelim-arcsys-twice}
\end{subfigure}
\begin{subfigure}[b]{0.3\linewidth}
\centering
\begin{tikzpicture}
\path[draw, shade, ball color=white] (0, 0) circle[radius=1];
\path[draw, bend left, ->] (-0.3, 0) to (0.3, 0);
\path[fill] (-0.3, 0) circle[radius=0.05] (0.3, 0) circle[radius=0.05];
\end{tikzpicture}
\caption{Two-punctured sphere}
\label{fig:prelim-arcsys-onearc}
\end{subfigure}
\caption{Illustrations of arcs and polygons}
\end{figure}
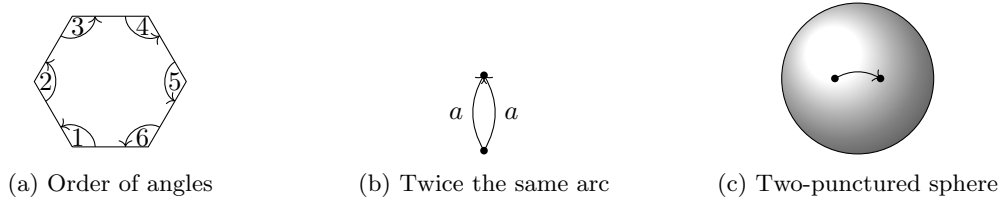

The reason we demand arc systems on spheres to have $ |M| ≥ 3 $ punctures becomes apparent: The [NMD] condition excludes the case of digons bounded by two different arcs, but we also desire to exclude the case where a digon is bounded by twice the same arc, depicted in \autoref{fig:prelim-arcsys-twice}. The only arc system with a digon bounded by twice the same arc is however the two-punctured sphere, depicted in \autoref{fig:prelim-arcsys-onearc}. This is the reason we require $ |M| ≥ 3 $.

\begin{definition}
The \emph{interior angles} of a polygon are the angles in the corners of the polygon. By \emph{angle}, we refer to the abstract entity (an interval starting at one arc and ending at the other, winding around their common endpoint) instead of the angle value.
\end{definition}

Since the punctured surface comes with an orientation, the interior angles of every polygon come with a natural cyclic (clockwise) order, see \autoref{fig:prelim-clockwise}. Working with arc systems often requires arguing with properties of the polygons and their angles. Some configurations of arcs and angles are allowed under the [NMDC] condition, others not.

\begin{remark}
In a [NMD] arc system, every polygon is bounded by a sequence of arcs with at least three interior angles in between. Indeed, zero angles bounding a polygon would mean the polygon is bounded by a single puncture. The punctured surface would necessarily be a one-punctured sphere, which we banned. A single angle bounding a polygon would mean that the polygon is bounded by a loop contractible in $ S \setminus M $, which we banned. Two angles bounding a polygon would mean they are equal, or distinct and homotopic in $ S \setminus M $. Both options are banned. In summary, every polygon in a [NMD] arc system is bounded by a sequence of arcs with at least three interior angles in between.
\end{remark}

\subsection{Dimers}
Dimer models, also referred to as brane tilings, originate in physicists' description of mirror symmetry. The idea is to describe arrangements of branes on the A-side of mirror symmetry in a surface graph. In a dimer model, adjacent nodes have opposite color. Dimer models can be seen as specific instances of punctured surfaces. A comprehensive reference is \cite{Bocklandt-dimers}.

\begin{definition}
\label{def:prelim-gtl-dimers-def}
A \emph{dimer} $ Q $ is a full arc system on a punctured surface such that
\begin{itemize}
\item every polygon is bounded by at least three arcs,
\item the arcs along the boundary of a polygon are all oriented in the same direction.
\end{itemize}
The letter $ Q $ also denotes the quiver, obtained from the arc system: Its vertex set $ Q_0 $ is the set of punctures and its arrow set $ Q_1 $ is the set of arcs. The underlying closed surface is denoted $ |Q| $.
\end{definition}

All polygons in a dimer are bounded either entirely clockwise or entirely anticlockwise. Neighboring polygons are bounded opposite: A polygon next to a clockwise polygon is anticlockwise, and a polygon next to an anticlockwise polygon is clockwise. The standard notation for a dimer is the letter $ Q $, minding the fact that the punctures together with the arcs can also be interpreted as a quiver embedded in a surface.

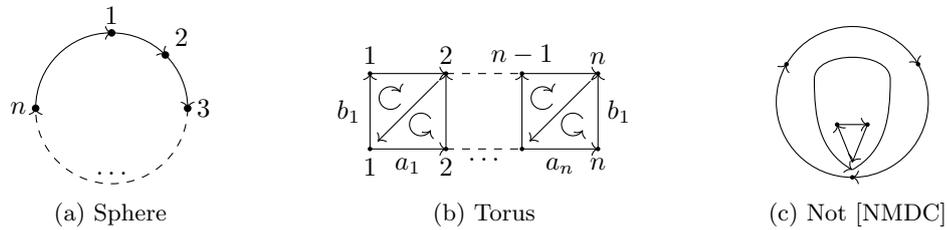
\begin{figure}
\centering
\begin{subfigure}{0.3\linewidth}
\centering
\begin{tikzpicture}
\path[draw, ->] (0, 1) node[above] {$ 1 $} coordinate (1) arc(90:45:1) coordinate (2) node[above right] {$ 2 $};
\path[draw, ->] (0, 0) ++(45:1) arc(45:0:1) coordinate (3) node[right] {$ 3 $};
\path[draw, dashed, ->] (1, 0) arc(0:-180:1) coordinate (n) node[left] {$ n $} node[midway, above] {$ … $};
\path[draw, ->] (-1, 0) arc(-180:-270:1);
\path[fill] (1) circle[radius=0.05] (2) circle[radius=0.05] (3) circle[radius=0.05] (n) circle[radius=0.05];
\end{tikzpicture}
\caption{Sphere}
\label{fig:prelim-std-dimers-sphere}
\end{subfigure}
\begin{subfigure}{0.3\linewidth}
\centering
\begin{tikzpicture}
\path[draw, ->] (0, 0) node[below] {$ 1 $} to node[midway, below] {$ a_1 $} (1, 0) node[below] {$ 2 $};
\path[draw, dashed] (1, 0) to node[midway, below] {$ … $} (2, 0);
\path[draw, ->] (2, 0) to node[midway, below] {$ a_n $} (3, 0) node[below] {$ n $};
\path[draw, ->] (3, 0) to node[midway, right] {$ b_1 $} (3, 1);
\path[draw, ->] (0, 0) to node[midway, left] {$ b_1 $} (0, 1) node[above] {$ 1 $};
\path[draw, ->] (0, 1) to (1, 1) node[above] {$ 2 $};
\path[draw, dashed] (1, 1) to (2, 1) node[above] {$ n-1 $};
\path[draw, ->] (2, 1) to (3, 1) node[above] {$ n $};
\foreach \i in {(0, 0), (0, 1), (1, 0), (1, 1), (2, 0), (2, 1), (3, 0), (3, 1)} \path[fill] \i circle[radius=0.03];
\path[draw, ->] (0.9, 0.9) to (0.1, 0.1); 
\path[draw, ->] (2.9, 0.9) to (2.1, 0.1); 
\path[draw, <-] (1, 1) to (1, 0);
\path[draw, <-] (2, 1) to (2, 0);
\path[draw, ->] (0.4, 0.6) arc(330:30:0.15);
\path[draw, ->] (0.8, 0.4) arc(30:330:0.15);
\path[draw, ->] (2.4, 0.6) arc(330:30:0.15);
\path[draw, ->] (2.8, 0.4) arc(30:330:0.15);
\end{tikzpicture}
\caption{Torus}
\label{fig:prelim-std-dimers-torus}
\end{subfigure}
\begin{subfigure}{0.3\linewidth}
\centering
\begin{tikzpicture}
\path[draw, ->] (0, -1) arc(270:150:1);
\path[draw, ->] (150:1) arc(150:30:1);
\path[draw, ->] (30:1) arc(30:-90:1);
\path[draw, ->] (0, -0.9) to[out=30, in=270]  ++(0.5, 1.2) to[out=90, in=90] ++(-1, 0) to[out=270, in=150] (0, -0.9);
\path[draw, ->] (0, -0.8) -- ++(-0.2, 0.5) coordinate (m1);
\path[draw, ->] (m1) -- ++(0.4, 0) coordinate (m2);
\path[draw, ->] (m2) -- ++(-0.2, -0.5);
\path[fill] (0, -1) circle[radius=0.03];
\path[fill] (150:1) circle[radius=0.03];
\path[fill] (30:1) circle[radius=0.03];
\path[fill] (m1) circle[radius=0.03];
\path[fill] (m2) circle[radius=0.03];
\end{tikzpicture}
\caption{Not [NMDC]}
\label{fig:prelim-std-dimers-NMDC}
\end{subfigure}
\caption{Standard dimer models and the [NMDC] condition}
\end{figure}

\begin{remark}
Every punctured surface has an arc system that is a dimer. Standard dimer models for the $ n $-punctured sphere ($ n ≥ 3 $) and $ n $-punctured torus ($ n ≥ 1 $) are depicted in \autoref{fig:prelim-std-dimers-sphere} and \ref{fig:prelim-std-dimers-torus}. A dimer automatically satisfies the [NMD] condition. There are however dimers which violate the [NMDC] condition, an example is depicted in \autoref{fig:prelim-std-dimers-NMDC}.
\end{remark}

\subsection{Gentle algebras}
\label{sec:prelim-gtl-gtl}
In this section, we recall gentle algebras associated with arc systems. We use a specific definition of gentle algebras, due to Bocklandt \cite{Bocklandt}. The reason they appear in this paper is that they form discrete models for Fukaya categories of punctured surfaces. In the present section, we describe only the algebra structure. The $ A_∞ $-structure will be added in \autoref{sec:prelim-gtl-infty}.

As their name suggests, gentle algebras are originally a type of finite-dimensional algebras. In \cite{Assem}, it was shown that so-called “unpunctured marked surface triangulations” naturally give rise to such gentle algebras. For readers familiar with Haiden-Katzarkov-Kontsevich's work \cite{HKK}, these are marked surfaces where all $ S^1 $ boundary components have at least one marking. The construction of gentle algebras from surfaces was subsequently carried over by Bocklandt \cite{Bocklandt} to the case of marked surfaces with full arc systems, as defined in \autoref{sec:prelim-arc-systems}. The definition is essential for this paper:

\begin{definition}
\label{def:prelim-gtl-gtl-def}
Let $ \cA $ be a full arc system on a punctured surface. Then the \emph{gentle algebra} (as ordinary algebra) $ \Gtl \cA = ℂ \rectified \cA / I $ is the quiver algebra with relations, where:
\begin{itemize}
\item The vertices of $ \rectified \cA $ are given by the arc midpoints of the arc system.
\item The arrows of $ \rectified \cA $ are given by the interior angles of the polygons.
\item The relations in $ I $ are given by all products of two consecutive interior angles of a polygon.
\end{itemize}
\end{definition}

\begin{figure}
\centering
\begin{subfigure}{0.9\linewidth}
\centering
\begin{tikzpicture}
\path[draw, thick] (1, 0) arc(0:90:1) coordinate[pos=0.3] (1-start) coordinate[pos=0.7] (2-end) node[midway, above right] {$ a $};
\path[draw, thick] (0, 1) arc(90:180:1) coordinate[pos=0.3] (2-start) coordinate[pos=0.7] (3-end) node[midway, above left] {$ d $};
\path[draw, thick] (-1, 0) arc(180:270:1) coordinate[pos=0.3] (3-start) coordinate[pos=0.7] (4-end) node[midway, below left] {$ c $};
\path[draw, thick] (0, -1) arc(-90:0:1) coordinate[pos=0.3] (4-start) coordinate[pos=0.7] (1-end) node[midway, below right] {$ b $};
\path[fill] (-1, 0) circle[radius=0.05] (0, 1) circle[radius=0.05] (0, -1) circle[radius=0.05] (1, 0) circle[radius=0.05];
\path[draw, ->, bend right=80, looseness=1.5] (1-start) to node[midway, right] {$ α $} (1-end);
\path[draw, ->, bend right=80, looseness=1.5] (2-start) to node[midway, above] {$ δ $} (2-end);
\path[draw, ->, bend right=80, looseness=1.5] (3-start) to node[midway, left] {$ γ $} (3-end);
\path[draw, ->, bend right=80, looseness=1.5] (4-start) to node[midway, below] {$ β $} (4-end);
\path[draw, ->, dashed, bend right=100, looseness=2.5] (1-end) to node[midway, left] {$ α' $} (1-start);
\path[draw, ->, dashed, bend right=100, looseness=2.5] (2-end) to node[midway, below] {$ δ' $} (2-start);
\path[draw, ->, dashed, bend right=100, looseness=2.5] (3-end) to node[midway, right] {$ γ' $} (3-start);
\path[draw, ->, dashed, bend right=100, looseness=2.5] (4-end) to node[midway, above] {$ β' $} (4-start);
\begin{scope}[shift={(3.25, 0)}]
\path (0, 1.8) node {$ \rectified \cA $};
\path[fill] (45:1) circle[radius=0.05] node[above right] {$ a $};
\path[fill] (135:1) circle[radius=0.05] node[above left] {$ b $};
\path[fill] (225:1) circle[radius=0.05] node[below left] {$ c $};
\path[fill] (315:1) circle[radius=0.05] node[below right] {$ d $};
\path[draw, ->, bend right] (50:1.1) to node[midway, above] {$ δ' $} (130:1.1);
\path[draw, ->, bend right] (130:0.9) to node[midway, above] {$ δ $} (50:0.9);
\path[draw, ->, bend right] (140:1.1) to node[midway, left] {$ γ' $} (220:1.1);
\path[draw, ->, bend right] (220:0.9) to node[midway, left] {$ γ $} (140:0.9);
\path[draw, ->, bend right] (230:1.1) to node[midway, below] {$ β' $} (310:1.1);
\path[draw, ->, bend right] (310:0.9) to node[midway, below] {$ β $} (230:0.9);
\path[draw, ->, bend right] (320:1.1) to node[midway, right] {$ α' $} (400:1.1);
\path[draw, ->, bend right] (400:0.9) to node[midway, right] {$ α $} (320:0.9);
\end{scope}
\begin{scope}[shift={(5.25, 0)}]
\path (0, 0) node[align=right] {Relations $ I $ \\ $ βα = 0 $ \\ $ γβ = 0 $ \\ $ δγ = 0 $ \\ $ αγ = 0 $ \\ $ α' β' = 0 $ \\ $ β' γ' = 0 $ \\ $ γ' δ' = 0 $ \\ $ δ' α' = 0 $};
\path (4.5, 0) node[align=left] {Basis for $ \Gtl \cA $ \\ $ \id_a (δδ')^i, \id_a (α'α)^i $ (full turns) \\ $ \id_b (αα')^i, \id_b (β'β)^i $ (full turns) \\ $ \id_c (ββ')^i, \id_c (γ'γ)^i $ (full turns) \\ $ \id_d (γγ')^i, \id_d (δ'δ)^i $ (full turns) \\ $ α (α' α)^i, ~ α' (α α')^i $ (full+1/2 turns) \\ $ β (β' β)^i, ~ β' (β β')^i $ (full+1/2 turns) \\ $ γ (γ' γ)^i, ~ γ' (γ γ')^i $ (full+1/2 turns) \\ $ δ (δ' δ)^i, ~ δ' (δ δ')^i $ (full+1/2 turns) \\ ($ i ≥ 0 $)};
\end{scope}
\end{tikzpicture}
\caption{Four-punctured sphere}
\label{fig:prelim-gtl-rectified-4p-sphere}
\end{subfigure}
\begin{subfigure}{0.9\linewidth}
\centering
\begin{tikzpicture}
\path[draw, thick, ->] (0, 0) to coordinate[pos=0.3] (1-start) coordinate[pos=0.7] (2-end) node[below] {$ a $} (2, 0);
\path[draw, thick, ->] (2, 0) to coordinate[pos=0.3] (2-start) coordinate[pos=0.7] (4-end) node[right] {$ b $} (2, 2);
\path[draw, thick, ->] (0, 2) to coordinate[pos=0.3] (3-end) coordinate[pos=0.7] (4-start) node[above] {$ a $} (2, 2);
\path[draw, thick, ->] (0, 0) to coordinate[pos=0.3] (1-end) coordinate[pos=0.7] (3-start) node[left] {$ b $} (0, 2);
\foreach \i in {0, 2} \foreach \j in {0, 2} \path[fill] (\i, \j) circle[radius=0.05];
\path[draw, bend right, ->] (1-start) to node[near end, below] {$ δ $} (1-end);
\path[draw, bend right, ->] (2-start) to node[near start, below] {$ γ $} (2-end);
\path[draw, bend right, ->] (3-start) to node[near start, above] {$ α $} (3-end);
\path[draw, bend right, ->] (4-start) to node[near end, above] {$ β $} (4-end);
\begin{scope}[shift={(3.5, 1)}]
\path (0, 1) node {$ \rectified \cA $};
\path[fill] (45:1) circle[radius=0.05] node[above right] {$ a $} (225:1) circle[radius=0.05] node[below left] {$ b $};
\path[draw, ->, bend right] (50:1.2) to node[midway, left] {$ δ $} (220:1.2);
\path[draw, ->, bend right] (50:0.9) to node[midway, below] {$ β $} (220:0.9);
\path[draw, ->, bend right] (230:1.2) to node[midway, right] {$ γ $} (400:1.2);
\path[draw, ->, bend right] (230:0.9) to node[midway, above] {$ α $} (400:0.9);
\end{scope}
\begin{scope}[shift={(6, 1)}]
\path (0, 0) node[align=right] {Relations $ I $ \\ $ βα = 0 $ \\ $ γβ = 0 $ \\ $ δγ = 0 $ \\ $ αδ = 0 $};
\path (5.5, 0) node[align=left] {Basis for $ \Gtl \cA $ \\ $ \id_a (αβγδ)^i, \id_a (γδαβ)^i $ (full turns) \\ $ \id_b (βγδα)^i, \id_b (δαβγ)^i $ (full turns) \\ $ α (βγδα)^i, β (γδαβ)^i $ (full+1/4 turns) \\ $ γ (δαβγ)^i, δ (αβγδ)^i $ (full+1/4 turns) \\ $ αβ (γδαβ)^i, βγ (δαβγ)^i $ (full+1/2 turns) \\ $ γδ (αβγδ)^i, δα (βγδα)^i $ (full+1/2 turns) \\ $ αβγ (δαβγ)^i, βγδ (αβγδ)^i $ (full+3/4 turns) \\ $ γδα (βγδα)^i, δαβ (γδαβ)^i $ (full+3/4 turns) \\ ($ i ≥ 0 $)};
\end{scope}
\end{tikzpicture}
\caption{One-punctured torus}
\label{fig:prelim-gtl-rectified-torus}
\end{subfigure}
\caption{Standard arc systems and their rectified quivers}
\label{fig:prelim-gtl-rectified}
\end{figure}
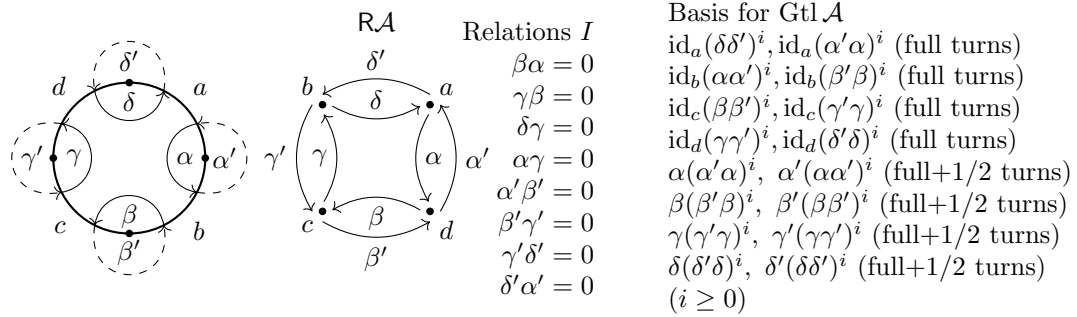
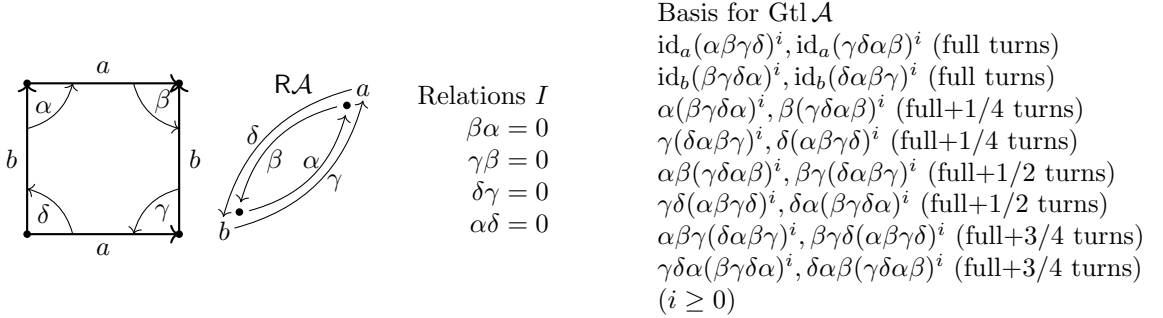

The quiver $ \rectified \cA $ has as many vertices as the arc system has arcs, as many arrows as the arc system has interior angles, and every polygon gives rise to as many relations as it has interior angles. The quiver $ \rectified \cA $ is called the rectified quiver in \cite{Bocklandt}. \autoref{fig:prelim-gtl-rectified} depicts some arc systems together with their rectified quivers. The left part of each graphic is the arc system itself, with arcs drawn thick. The interior angles are drawn as thin arrows; in the three- and four-punctured sphere, the dashed arrows mean the interior angles at the rear, invisible side of the sphere. The right part of each graphic depicts the rectified quiver together with its relations.

A vector space basis for the gentle algebra $ \Gtl \cA $ consists of all angles around punctures. The basis includes an identity $ \id_a $ for every arc $ a ∈ \cA $, which we may also view as an empty angle. The gentle algebra $ \Gtl \cA $ of an arc system is not finite-dimensional.

\begin{remark}
By nature, the algebra $ \Gtl \cA = ℂ \rectified \cA / I $ can be viewed as a $ ℂ $-linear category with objects being the arcs of the arc system. The hom spaces are spanned freely by the angles winding around punctures, starting at one arc and ending at another arc. In \cite{Bocklandt}, this interpretation of $ \Gtl \cA $ is also called the “gentle category”. We will consistently use the term gentle algebra instead, despite the slight inaccuracy.
\end{remark}

\begin{remark}
In the work of Haiden, Katzarkov and Kontsevich \cite{HKK}, gentle algebras were developed further under the name of “topological Fukaya categories“. This includes a generalization regarding the type of boundary allowed. The version of gentle algebra we defined above is called a surface with fully marked boundaries of $ S^1 $ type in \cite{HKK}. If one changes the boundary type to have at least one so-called “boundary arc” on every boundary component, the topological Fukaya category becomes finite-dimensional. The gentle algebras $ \Gtl \cA $ studied in the present paper are however infinite-dimensional by nature: One can keep winding around the punctures as often as one wants, obtaining morphisms of higher and higher length.
\end{remark}

It might be worthwhile comparing to the original definition due to \cite{Assem}: A finite-dimensional algebra presented as $ ℂQ / I $ is \emph{gentle} if
\begin{itemize}
\item At each vertex there start at most two arrows and there end at most two arrows.
\item The ideal $ I $ is generated by paths of length 2.
\item For every arrow $ β $, there is at most one arrow $ α $ such that $ αβ ∈ I $, and at most one arrow $ γ $ such that $ βγ ∈ I $.
\item For every arrow $ β $, there is at most one arrow $ α $ such that $ αβ \notin I $, and at most one arrow $ γ $ such that $ βγ \notin I $.
\end{itemize}
The same paper \cite{Assem} showed how surface triangulations naturally give rise to gentle algebras.

\subsection{The $ A_∞ $-structure on $ \Gtl \cA $}
\label{sec:prelim-gtl-infty}
In this section, we recall the $ A_∞ $-structure on $ \Gtl \cA $. The starting point is the description of $ \Gtl \cA $ as ordinary algebra in \autoref{sec:prelim-gtl-gtl}. The idea is to add $ A_∞ $-structure which captures the topology of the punctured surface. This $ A_∞ $-structure was introduced by Bocklandt \cite{Bocklandt} in order to define a discrete version of the wrapped Fukaya category. In the present section, we recall the $ A_∞ $-structure briefly and refer to \cite[Section 9.1]{Bocklandt-book} and \cite{HKK} and for more insight.

We start with a full arc system $ \cA $ with [NMD]. The first step towards the $ A_∞ $-structure is the grading. It is possible to put a $ ℤ $-grading on $ \Gtl \cA $ by viewing arcs and angles relative to a vector field on the surface, see \autoref{sec:uncurving-stringsbands}. The deformations of $ \Gtl \cA $ that we are interested in exist however only in the $ ℤ/2ℤ $-graded world. Consequentially, we define $ \Gtl \cA $ as a $ ℤ/2ℤ $-graded $ A_∞ $-category from the very beginning. The definition of degrees is depicted in \autoref{fig:prelim-gtl-infty-deg} and reads as follows:

\begin{definition}
The \emph{degree} $ |α| $ of an angle $ α: a → b $ is odd if one of the arcs $ a, b $ points towards the puncture, and one points away from the puncture. The degree of an angle is even if both arcs point away or both point towards the puncture.
\end{definition}

\begin{figure}
\centering
\begin{subfigure}{0.2\linewidth}
\centering
\begin{tikzpicture}
\path[draw, ->] (-1, 1) -- (0, 0) coordinate[midway] (2);
\path[draw, ->] (0, 0) -- (1, 1) coordinate[midway] (1);
\path[draw, ->, bend right] (1) to node[midway, below] {$ α $} (2);
\end{tikzpicture}
\caption{odd}
\end{subfigure}
\begin{subfigure}{0.2\linewidth}
\centering
\begin{tikzpicture}
\path[draw, <-] (-1, 1) -- (0, 0) coordinate[midway] (2);
\path[draw, <-] (0, 0) -- (1, 1) coordinate[midway] (1);
\path[draw, ->, bend right] (1) to node[midway, below] {$ α $} (2);
\end{tikzpicture}
\caption{odd}
\end{subfigure}
\begin{subfigure}{0.2\linewidth}
\centering
\begin{tikzpicture}
\path[draw, <-] (-1, 1) -- (0, 0) coordinate[midway] (2);
\path[draw, ->] (0, 0) -- (1, 1) coordinate[midway] (1);
\path[draw, ->, bend right] (1) to node[midway, below] {$ α $} (2);
\end{tikzpicture}
\caption{even}
\end{subfigure}
\begin{subfigure}{0.2\linewidth}
\centering
\begin{tikzpicture}
\path[draw, ->] (-1, 1) -- (0, 0) coordinate[midway] (2);
\path[draw, <-] (0, 0) -- (1, 1) coordinate[midway] (1);
\path[draw, ->, bend right] (1) to node[midway, below] {$ α $} (2);
\end{tikzpicture}
\caption{even}
\end{subfigure}
\caption{Degree of angle $ α $}
\label{fig:prelim-gtl-infty-deg}
\end{figure}

The second step in the definition of the $ A_∞ $-structure is the definition of the differential $ μ^1 $ and the product $ μ^2 $. The differential $ μ^1 $ is plainly set to zero. We keep the notation $ αβ $ for the concatenation of angles, and define the product $ μ^2 $ as its signed version:
\begin{equation*}
μ^1 ≔ 0, \quad μ^2 (α, β) ≔ (-1)^{|β|} αβ.
\end{equation*}

\begin{remark}
The reason we assume the [NMD] condition is that it simplifies the definition of $ μ^1 $ and $ μ^2 $. Indeed, can also define the $ A_∞ $-structure for arc systems without [NMD]. However, the definition of $ μ^1 $ and $ μ^2 $ then needs to be tweaked in order to capture the monogons and digons.
\end{remark}

\begin{figure}
\centering
\begin{subfigure}[b]{0.3\linewidth}
\centering
\begin{tikzpicture}
\path[draw] (0, 0) node[below] {3} -- ++(right:1) node[midway, above] {\tiny 2} node[below] {2} -- ++(72:1) node[midway, left] {\tiny 1} node[right] {1} -- ++(144:1) node[midway, below] {\tiny 5} node[above] {5} -- ++(216:1) node[midway, below] {\tiny 4} node[left] {4} -- ++(-72:1) node[midway, right] {\tiny 3};
\end{tikzpicture}
\caption{Standard polygon $ P_5 $}
\label{fig:prelim-gtl-infty-P5}
\end{subfigure}
\begin{subfigure}[b]{0.3\linewidth}
\centering
\begin{tikzpicture}
\path[draw] (0, 0) -- ++(left:1) coordinate[pos=0.3] (1-start) coordinate[pos=0.7] (2-end) -- ++(300:1) coordinate[pos=0.3] (2-start) coordinate[pos=0.7] (3-end) -- ++(right:1) coordinate[pos=0.3] (3-start) coordinate[pos=0.7] (4-end) -- ++(60:1) coordinate[pos=0.3] (4-start) coordinate[pos=0.7] (5-end) -- ++(left:1) coordinate[pos=0.3] (5-start) coordinate[pos=0.7] (1-end);
\path[draw] (0, 0) -- ++(240:1) (0, 0) -- ++(300:1);
\foreach \i in {2, 5} \path[draw, ->, bend right] (\i-start) to (\i-end);
\foreach \i in {3, 4} \path[draw, ->, bend right=60] (\i-start) to (\i-end);
\path[draw, ->, bend right=90, looseness=2] (1-start) to (1-end);
\end{tikzpicture}
\caption{A disk sequence}
\label{fig:prelim-gtl-infty-diskexample}
\end{subfigure}
\begin{subfigure}[b]{0.3\linewidth}
\centering
\begin{tikzpicture}
\path[draw] (0, 0) -- ++(left:1) coordinate (A) -- ++(60:1) coordinate (B) coordinate[pos=0.3] (3-end) coordinate[pos=0.7] (4-start) -- ++(right:1) coordinate[pos=0.3] (4-end) coordinate[pos=0.7] (5-start) coordinate (C) -- ++(300:1) coordinate[pos=0.3] (5-end) coordinate[pos=0.7] (6-start) coordinate (D) -- ++(240:1) coordinate[pos=0.3] (6-end) coordinate[pos=0.7] (7-start) coordinate (E) -- ++(left:1) coordinate[pos=0.4] (1-end) coordinate[pos=0.3] (7-end) coordinate[pos=0.7] (2-start) coordinate[pos=0.6] (8-start) coordinate (F) -- ++(120:1) coordinate[pos=0.3] (2-end) coordinate[pos=0.7] (3-start);
\path[draw] (0, 0) -- (F) coordinate[pos=0.6] (8-end) coordinate[pos=0.3] (9-start);
\path[draw] (0, 0) -- (E) coordinate[pos=0.6] (1-start) coordinate[pos=0.3] (9-end);
\foreach \i in {2, 3, 4, 5, 6, 7} \path[draw, ->, bend right=60] (\i-start) to (\i-end);
\path[draw] (A) -- (D) (B) -- (E) (C) -- (F);
\end{tikzpicture}
\caption{Not a disk sequence}
\label{fig:prelim-gtl-infty-no-disk}
\end{subfigure}
\caption{Illustration of discrete immersed disks}
\label{fig:prelim-gtl-infty-diskillustration}
\end{figure}
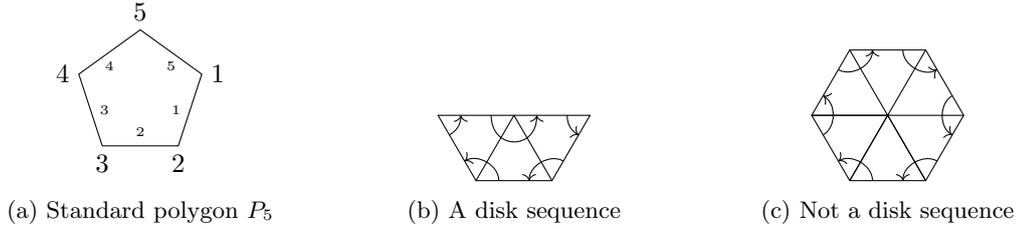

The third step is to define the higher products $ μ^{≥3} $ of $ \Gtl \cA $. They capture the topology of the arcs and angles. Roughly speaking, a higher product of a sequence of angles is nonzero if the sequence bounds a disk. Such a disk is given by an immersion of the standard polygon $ P_k $ into the surface $ S $, known as immersed disk. The domain of the immersion mapping is a standard polygon $ P_k $, depicted in \autoref{fig:prelim-gtl-infty-P5}. To distinguish this type of immersed disks from the type used for the Fukaya category, we shall refer to these disks as discrete immersed disks. The precise definition reads as follows:

\begin{definition}
Let $ \cA $ be a full arc system with [NMD]. A \emph{discrete immersed disk} in $ \cA $ consists of an oriented immersion $ D: P_k → S $ of a standard polygon $ P_k $ into the surface, such that
\begin{itemize}
\item The edges of the polygon are mapped to a sequence of arcs.
\item The immersion does not cover any punctures.
\end{itemize}
The immersion mapping $ D $ itself is only taken up to reparametrization. The sequence of \emph{interior angles} of $ D $ is the sequence of angles in $ \cA $ given as images of the interior angles of $ P_k $ under the map $ D $. An angle sequence $ α_1, …, α_k $ is a \emph{disk sequence} if it is the sequence of interior angles of some discrete immersed disk.
\end{definition}

To explain the definition in other words, the image of the interior of the polygon $ P_k $ consists only of polygon interiors of $ \cA $ and arcs between punctures, but not punctures themselves. The boundary of $ P_k $ is mapped to a sequence of arcs, and the corners inside $ P_k $ are mapped to an angle sequence in the arc system.

\begin{example}
The interior angles of a polygon, when written in clockwise order, are a disk sequence. In particular, if the arc system $ \cA $ has a triangle polygon in it, then there exists an disk sequence in $ \cA $ of just 3 angles. In every case, by the ban on loops and homotopic arcs, an disk sequence $ α_1, …, α_k $ consists of at least three angles, i.e.~$ k ≥ 3 $. In \autoref{fig:prelim-gtl-infty-P5}, we have depicted the schematic of a standard polygon. In \autoref{fig:prelim-gtl-infty-diskexample}, we have depicted a discrete immersed disk together with its sequence of interior angles. By definition, this sequence is a disk sequence. In \autoref{fig:prelim-gtl-infty-no-disk}, we have depicted an angle sequence which is not a disk sequence. The reason it is not a disk sequence is that there is a polygon immersion bounded by the drawn angles, but it covers the puncture at the center of the hexagon. Later on, we will however allow polygon immersions which cover punctures as part of the deformation $ \Gtl_q \cA $.
\end{example}

\begin{figure}
\centering
\begin{subfigure}{0.4\linewidth}
\centering
\begin{tikzpicture}
\path[draw] (0, 0) -- ++(left:1) -- ++(300:1) -- cycle coordinate[midway] (A1);
\path[draw] (0.1, 0) -- ++(240:1) coordinate[midway] (A2) -- ++(right:1) -- cycle coordinate[midway] (B1);
\path[draw] (0.2, 0) -- ++(300:1) coordinate[midway] (B2) -- ++(60:1) -- cycle;
\path[draw, ultra thick] ($ (A1)!-0.5!(A2) $) -- ($ (A2)!-0.5!(A1) $);
\path[draw, ultra thick] ($ (B1)!-0.5!(B2) $) -- ($ (B2)!-0.5!(B1) $);
\path (1.7, -0.3) node {\LARGE $ \rightsquigarrow $};
\begin{scope}[shift={(3, 0)}]
\path[draw] (0, 0) -- ++(left:1) coordinate[pos=0.3] (1-start) coordinate[pos=0.7] (2-end) -- ++(300:1) coordinate[pos=0.3] (2-start) coordinate[pos=0.7] (3-end) -- ++(right:1) coordinate[pos=0.3] (3-start) coordinate[pos=0.7] (4-end) -- ++(60:1) coordinate[pos=0.3] (4-start) coordinate[pos=0.7] (5-end) -- ++(left:1) coordinate[pos=0.3] (5-start) coordinate[pos=0.7] (1-end);
\path[draw] (0, 0) -- ++(240:1) (0, 0) -- ++(300:1);
\foreach \i in {2, 5} \path[draw, ->, bend right] (\i-start) to (\i-end);
\foreach \i in {3, 4} \path[draw, ->, bend right=60] (\i-start) to (\i-end);
\path[draw, ->, bend right=90, looseness=2] (1-start) to (1-end);
\end{scope}
\end{tikzpicture}
\end{subfigure}
\begin{subfigure}{0.6\linewidth}
\centering
\begin{tikzpicture}
\path[draw] (0, 0) -- ++(left:1) coordinate[midway] (A1) -- ++(300:1) -- cycle coordinate[midway] (D1);
\path[draw] (0.1, 0) -- ++(240:1) coordinate[midway] (D2) -- ++(right:1) -- cycle coordinate[midway] (E1);
\path[draw] (0.2, 0) -- ++(300:1) coordinate[midway] (E2) -- ++(60:1) -- cycle coordinate[midway] (F1);
\path[draw] (0, 0.1) -- ++(left:1) coordinate[midway] (A2) -- ++(60:1) -- cycle coordinate[midway] (B1);
\path[draw] (0.1, 0.1) -- ++(120:1) coordinate[midway] (B2) -- ++(right:1) -- cycle coordinate[midway] (C1);
\path[draw] (0.2, 0.1) -- ++(60:1) coordinate[midway] (C2) -- ++(300:1) -- cycle coordinate[midway] (G1);
\path[draw, dashed] (0.4, -0.1) -- ++(300:1) -- ++(60:1) -- cycle coordinate[midway] (G2);
\foreach \i in {A, B, C, D, E, G} \path[draw, ultra thick] ($ (\i 1)!-0.5!(\i 2) $) -- ($ (\i 2)!-0.5!(\i 1) $);
\path (2, 0) node {\LARGE $ \rightsquigarrow $};
\begin{scope}[shift={(3.5, 0)}]
\path[draw] (0, 0) -- ++(left:1) coordinate (A) -- ++(60:1) coordinate (B) coordinate[pos=0.3] (3-end) coordinate[pos=0.7] (4-start) -- ++(right:1) coordinate[pos=0.3] (4-end) coordinate[pos=0.7] (5-start) coordinate (C) -- ++(300:1) coordinate[pos=0.3] (5-end) coordinate[pos=0.7] (6-start) coordinate (D) -- ++(240:1) coordinate[pos=0.3] (6-end) coordinate[pos=0.7] (7-start) coordinate (E) -- ++(left:1) coordinate[pos=0.4] (1-end) coordinate[pos=0.3] (7-end) coordinate[pos=0.7] (2-start) coordinate[pos=0.6] (8-start) coordinate (F) -- ++(120:1) coordinate[pos=0.3] (2-end) coordinate[pos=0.7] (3-start);
\path[draw] (0, 0) -- (F) coordinate[pos=0.6] (8-end) coordinate[pos=0.3] (9-start);
\path[draw] (0, 0) -- (E) coordinate[pos=0.6] (1-start) coordinate[pos=0.3] (9-end);
\foreach \i in {2, 3, 4, 5, 6, 7} \path[draw, ->, bend right=60] (\i-start) to (\i-end);
\path[draw, ->, bend right] (8-start) to (8-end);
\path[draw, ->, bend right] (1-start) to (1-end);
\path[draw, ->] (9-start) arc(-120:90:0.3) arc(90:300:0.4);
\path[draw] (A) -- (D) (B) -- (E) (C) -- (F);
\path (1.5, 0) node {\LARGE $ = $};
\end{scope}
\begin{scope}[shift={(6.5, 0)}, xslant=0.2, yscale=0.8]
\path[draw] (0, 0) -- ++(300:1) coordinate[pos=0.3] (1-end) arc(-60:-240:1) coordinate (A);
\path[draw] (A) arc(120:-120:1.2) -- (0, 0) coordinate[pos=0.7] (1-start);
\path[draw, dashed, ->] (1-start) arc(-120:120:0.3) arc(120:300:0.4);
\path[fill] (0, 0) circle[radius=0.05];
\path[fill] (300:1) circle[radius=0.05];
\end{scope}
\end{tikzpicture}
\end{subfigure}
\caption{Stitching together polygons yields immersed disks}
\label{fig:prelim-gtl-infty-stitching}
\end{figure}
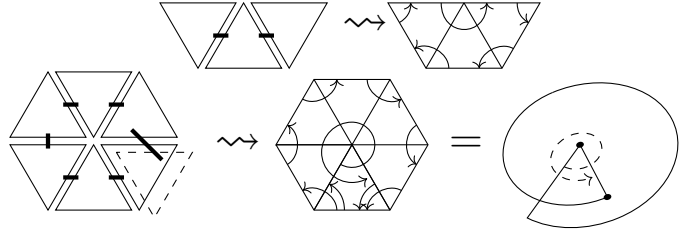

Disk sequences $ α_1, …, α_k $ can also be described combinatorically: They are either a polygon in $ \cA $, or stitched together from multiple polygons along arcs. \autoref{fig:prelim-gtl-infty-stitching} depicts two examples of stitching polygons together to form disk sequences. In every example, multiple triangles are stitched together to form a polygon. Thick connectors between two triangles indicate that these triangles are going to be stitched together along their shared edge. The first example is visually easy to grasp, since the three triangles are disjoint. In the second example, seven triangles are stitched together, with one triangle appearing twice. The result is a disk sequence of nine angles $ α_1, …, α_9 $, of which one is longer than a full turn. The sketch on the right of the “=” sign provides a visualization of this discrete immersed disk by thinking a third dimension into the picture. In that 3-dimensional sketch, the angle longer than a full turn is drawn dashed, the other angles are omitted and the outer boundary of the hexagon is depicted as a spiral instead of separate arcs.

\begin{remark}
In \autoref{sec:prelim-gtlq}, we change the terminology. From there on, discrete immersed disks are allowed to cover punctures.
\end{remark}

We are now ready to give the definition of the higher products $ μ^{≥3} $. Since these products are supposed to be multilinear, it suffices to define them on the basis of $ \Gtl \cA $ given by angles winding around punctures.

\begin{definition}
Let $ \cA $ be an full arc system with [NMD]. Then $ \Gtl \cA $ is the $ A_∞ $-category with objects being the arcs $ a ∈ \cA $, hom spaces spanned by angles, and $ A_∞ $-product $ μ $ defined by $ μ^1 = 0 $ and $ μ^2 (α, β) = (-1)^{|β|} αβ $. To define $ μ^{k≥3} $, let $ α_1, …, α_k $ be any disk sequence, let $ β $ be an angle composable with $ α_1 $, i.e.~$ β α_1 ≠ 0 $, and let $ γ $ be an angle post-composable with $ α_k $, i.e.~$ α_k γ ≠ 0 $. Then
\begin{equation*}
μ^k (β α_k, …, α_1) ≔ β, \quad μ^k (α_k, …, α_1 γ) ≔ (-1)^{|γ|} γ.
\end{equation*}
The higher products vanish on all angle sequences other than these.
\end{definition}

\begin{example}
Let us go through a few examples. Regard the one-punctured torus of \autoref{fig:prelim-gtl-rectified-torus}. The angles $ α $ and $ γ $ are odd, and $ β $ and $ δ $ are even. Angle degrees add up, for instance $ βγδ $ is odd and $ γδα $ is even. The product $ μ^2 $ has
\begin{equation*}
μ^2 (δ, γ) = - δγ, \quad μ^2 (αβγδα, βγδ) = - αβγδαβγδ = - (αβγδ)^2 \text{ and } μ^2 (β, α) = 0.
\end{equation*}
The higher product $ μ^3 $ vanishes because there are no triangles. From $ μ^4 $ onwards, we have higher products, for instance
\begin{equation*}
μ^4 (δ, γ, β, α) = \id_b \text{ and } μ^6 (γ, βγ, β, α, δα, δ) = \id_a.
\end{equation*}
A little less obvious is the $ A_∞ $-product
\begin{equation*}
μ^{12} (α, δ, γδ, γ, βγ, β, αβ, α, δα, δ, γ, βγδαβ) = \id_a,
\end{equation*}
which is “winds” one and a quarter times around the puncture, without covering the puncture itself though. The second disk sequence in \autoref{fig:prelim-gtl-infty-stitching} is very similar and also yields an identity.
\end{example}

After defining this structure, Bocklandt \cite{Bocklandt} proved that with this grading and products $ \Gtl \cA $ is indeed an $ A_∞ $-category. 

\begin{theorem}[\cite{Bocklandt}]
Let $ \cA $ be a full arc system with [NMD]. Then $ \Gtl \cA $ is an $ A_∞ $-category.
\end{theorem}

\begin{remark}
In \cite{Bocklandt}, the signs in the definition of the higher products $ μ^k $ on $ \Gtl \cA $ differ from the signs presented here. We follow the sign convention of \cite{HKK}.
\end{remark}

\begin{remark}
Every angle sequence $ α_1, …, α_k $ either bounds a unique discrete immersed disk or no disk at all. If it bounds a disk, then the products $ μ(α_k, …, α_1 γ) $ and $ μ(β α_k, …, α_1) $ are nonzero. If it bounds no disk, then the products vanish.
\end{remark}

\begin{remark}
Let us explain that that the degrees match. The $ A_∞ $-product $ μ^k $ is required to be of parity $ 2 - k $. If $ α_1, …, α_k $ is any disk sequence, then the total reduced degree
\begin{equation*}
‖α_1‖ + … + ‖α_k‖ ∈ ℤ/2ℤ
\end{equation*}
measures how often the boundary of the discrete immersed disk changes orientation when traversing it clockwise. Since a disk sequence traverses the boundary one full time, it ends up with the same orientation as it started. In other words, the total reduced degree of a disk sequence vanishes. This means that $ μ^k $ has the right parity.
\end{remark}

\begin{remark}
The interior angles of a discrete immersed disk are enumerated clockwise as $ α_1, …, α_k $, while the higher product consumes them only in the order $ α_k, …, α_1 $. This seemingly unusual order of the factors $ α_1, …, β α_k $ is due to the convention on $ A_∞ $-categories.
\end{remark}

\begin{remark}
It is not possible to write a given angle sequence $ γ_k, …, γ_1 $ as $ β α_k, …, α_1 $ or $ α_k, …, α_1 γ $ in two different ways. In fact, the immersion of the polygon is already determined by all angles but one, and an angle sequence of the form $ β α_k, …, α_1 $ cannot be written as $ α_k', …, α_1' γ $ with both $ α_1, …, α_k $ and $ α_1', …, α_k' $ being disk sequences. This is explained e.g.~in \cite{HKK}. We conclude that any angle sequence can be written in at most one way as $ β α_k, …, α_1 $ or $ α_k, …, α_1 γ $ with $ α_1, …, α_k $ a disk sequence. This makes the product $ μ^k (γ_k, …, γ_1) $ well-defined for every angle sequence $ γ_1, …, γ_k $.
\end{remark}

\subsection{The deformation $ \Gtl_q \cA $}
\label{sec:prelim-gtlq}
In this section we define the deformation $ \Gtl_q \cA $ of a gentle algebra $ \Gtl \cA $. It is a specific instance of the deformations constructed in \paperone\ and lies at the heart of the present paper. We give an explicit definition in order to provide a feel for this category. We provide a first glance concerning the use of $ \Gtl_q \cA $ in the later sections.

This paper is the second in a series of three, and the first paper \paperone\ was concerned with classifying the $ A_∞ $-deformations of $ \Gtl \cA $. A conclusion from that paper is that all deformations of $ \Gtl \cA $ up to gauge equivalence can be written down explicitly. In this paper, we select one of these deformations, which we call the deformed gentle algebra and denote by $ \Gtl_q \cA $. Other deformations of $ \Gtl \cA $ play no role anymore.

The deformation $ \Gtl_q \cA $ is very broad in the sense that it has a lot of deformation parameters, in fact one for each puncture. Any reader who wishes to work with the calculations of this paper can therefore freely set some of these deformation parameters to zero and still have an interesting deformation at hand. Conversely, the deformation $ \Gtl_q \cA $ is so broad that the reader who is interested in deformations not “covered” by $ \Gtl_q \cA $ can still derive qualitative expectations on the behavior of the other deformations.

\begin{remark}
Arguably, one would like to conduct the study of the present paper also for all the other deformations given in \paperone. The idea would be to use multiple parameters per puncture, so as to include deformations in $ \Gtl_q \cA $ that measure orbigons around punctures (see \paperone). One reason we restrict to the single deformation $ \Gtl_q \cA $ is that “orbigon deformations” are more difficult to handle than “disk deformations”. Another reason is that the relative Fukaya category $ \relFuk (S, M) $ also has only “disk deformations” as well, so a candidate for a small model of $ \relFuk (S, M) $ should only have “disk deformations” at all. This is why we only regard the deformation $ \Gtl_q \cA $, which has one parameter per puncture.
\end{remark}

The deformation base of $ \Gtl_q \cA $ is $ B = ℂ⟦M⟧ $. This is the commutative local ring of power series in $ |M| $ variables, one for each puncture. In fact, to capture the punctures covered by an immersion of a standard polygon, every puncture should have one deformation parameter. For this reason we use $ B = ℂ⟦M⟧ $. Every puncture $ q ∈ M $ gives rise to one deformation parameter, which is also denoted $ q $ and lies in the ring $ ℂ⟦M⟧ $ as one of the generators.

\begin{remark}
We use the letter “$ q $” as in three different meanings in this paper, depending on the context: First, the notation $ \Gtl_q \cA $ is fixed and the letter $ q $ does not have any meaning there. Second, whenever a specific puncture is considered, it is typically named $ q $. Third, whenever $ q $ is used multiplicatively in formulas, then it denotes the infinitesimal parameter $ q ∈ ℂ⟦M⟧ $. For example, if $ p, q ∈ M $ are punctures, then $ pq $ simply means the product $ pq ∈ ℂ⟦M⟧ $.
\end{remark}

\begin{figure}
\centering
\begin{subfigure}[b]{0.3\linewidth}
\centering
\begin{tikzpicture}[scale=0.35]
\path[draw, semithick] plot[smooth cycle, tension=1] coordinates {(4, 0) 
(5, 6) (5, 5) 
(2, 7.5) (2.5, 1) (10, 4) (7, 8)
(3, 5) (3.5, 6) 
(3, 1.5) (9, 5) (5, 10) (0, 8)
(4, 3) (3, 3) 
(8, 6) (2, 9) (0, 3)};
\end{tikzpicture}
\caption{Bennequin's curve}
\label{fig:prelim-gtlq-bennequin}
\end{subfigure}
\begin{subfigure}[b]{0.3\linewidth}
\centering
\begin{tikzpicture}
\path[draw, ->] (0, 0) node[left] {$ q $} -- (2, 0) node[right] {$ p $} node[midway, above] {$ a $};
\path[fill] (0, 0) circle[radius=0.05] (2, 0) circle[radius=0.05];
\path[draw, ->] (2, 0) ++(-170:0.5) arc(-170:170:0.5) node[midway, right] {$ ℓ_p $};
\path[draw, ->] (0, 0) ++(10:0.5) arc(10:350:0.5) node[midway, left] {$ ℓ_q $};
\path (1, -1) node {$ μ^0_q = q ℓ_q + p ℓ_p $};
\end{tikzpicture}
\caption{Assigning curvature to an arc}
\end{subfigure}
\begin{subfigure}[b]{0.3\linewidth}
\centering
\begin{tikzpicture}[scale=1]
\path[draw] (0, 0) -- ++(left:1) coordinate (A) -- ++(60:1) coordinate (B) coordinate[pos=0.3] (3-end) coordinate[pos=0.7] (4-start) -- ++(right:1) coordinate[pos=0.3] (4-end) coordinate[pos=0.7] (5-start) coordinate (C) -- ++(300:1) coordinate[pos=0.3] (5-end) coordinate[pos=0.7] (6-start) coordinate (D) -- ++(240:1) coordinate[pos=0.3] (6-end) coordinate[pos=0.7] (7-start) coordinate (E) -- ++(left:1) coordinate[pos=0.4] (1-end) coordinate[pos=0.3] (7-end) coordinate[pos=0.7] (2-start) coordinate[pos=0.6] (8-start) coordinate (F) -- ++(120:1) coordinate[pos=0.3] (2-end) coordinate[pos=0.7] (3-start);
\path[draw] (0, 0) -- (F) coordinate[pos=0.6] (8-end) coordinate[pos=0.3] (9-start);
\path[draw] (0, 0) -- (E) coordinate[pos=0.6] (1-start) coordinate[pos=0.3] (9-end);
\foreach \i in {2, 3, 4, 5, 6, 7} \path[draw, ->, bend right=60] (\i-start) to (\i-end);
\path[draw] (A) -- (D) (B) -- (E) (C) -- (F);
\path[fill] (0, 0) circle[radius=0.05];
\path (A) node[left] {$ α_1 $};
\path (B) node[left] {$ α_2 $};
\path (C) node[right] {$ α_3 $};
\path (D) node[right] {$ α_4 $};
\path (E) node[right] {$ α_5 $};
\path (F) node[left] {$ α_6 $};
\path (F) -- (A) node[midway, left] {$ a $};
\path (0, 0) node[above] {$ q $};
\path (0, -1.5) node {$ μ^6_q (α_6, …, α_1) = q \id_a $};
\end{tikzpicture}
\caption{Deformed product $ μ^6_q $}
\label{fig:prelim-gtl-infty-nodisk}
\end{subfigure}
\caption{Illustration of the deformation $ \Gtl_q \cA $}
\end{figure}
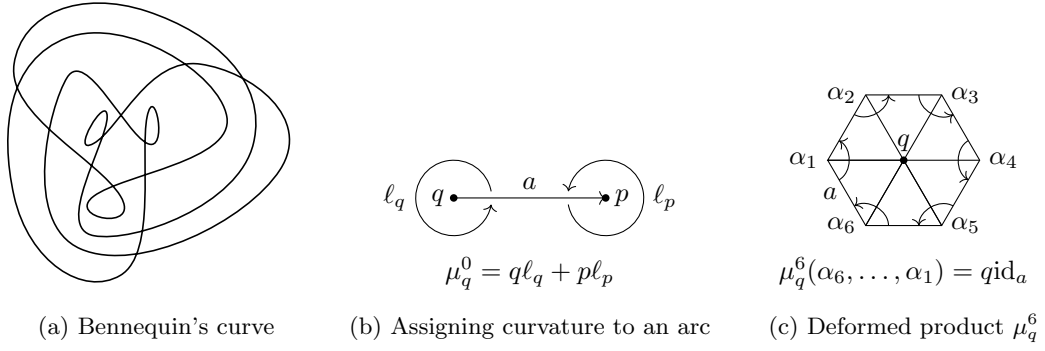

As a warm-up for the definition of $ \Gtl_q \cA $, recall that the angle sequence $ α_1, …, α_6 $ of \autoref{fig:prelim-gtl-infty-nodisk} is not a disk sequence. This means that $ μ^6 (α_6, …, α_1) = 0 $ in $ \Gtl \cA $. The deformation $ \Gtl_q \cA $ precisely changes this and similar higher products, while keeping the $ A_∞ $-relations intact. In short, the deformed higher products of $ \Gtl_q \cA $ precisely capture which and how often punctures are covered by an immersion of a standard polygon. From here on, we drop the requirement that a discrete immersed disk does not cover punctures:

\begin{definition}
Let $ \cA $ be a full arc system with [NMDC]. A \emph{discrete immersed disk} in $ \cA $ is an oriented immersion of a standard polygon $ P_k $ into $ S $ up to reparametrization such that the edges of the polygon are mapped to a sequence of arcs. A \emph{disk sequence} is an angle sequence together with a choice of discrete immersed disk of which it is the sequence of interior angles. We denote by $ M(α_1, …, α_k) $ the set of discrete immersed disks $ D $ with interior angles $ α_1, …, α_k $. For $ D ∈ M(α_1, …, α_k) $, we denote by $ q_D ∈ ℂ⟦M⟧ $ the product of the punctures covered by $ D $.
\end{definition}

\begin{remark}
\label{rem:prelim-gtlq-bennequin}
In contrast to \autoref{sec:prelim-gtl-infty}, for a given angle sequence $ α_1, …, α_k $ there might be multiple discrete immersed disks which have the same interior angle sequence $ α_1, …, α_k $. To see this, regard Bennequin's curve in \autoref{fig:prelim-gtlq-bennequin}. This smooth curve bounds five immersed disks which are not related by reparametrization. If we draw a fine enough grid in the surface and approximate the curve by arcs in the grid, then we obtain an angle sequence which bounds multiple distinct discrete immersed disks.
\end{remark}

The deformation $ \Gtl_q \cA $ has infinitesimal curvature, and there are three ways to describe the curvature: Each puncture $ q ∈ M $ contributes curvature $ qℓ_q ∈ (M) \htensor \Gtl \cA $ to $ \Gtl_q \cA $. Here $ (M) ⊂ ℂ⟦M⟧ $ denotes the maximal ideal of $ ℂ⟦M⟧ $ and $ ℓ_q $ denotes the sum of all full turns around $ q $, summed over all arc incidences at $ q $. In other words, the total curvature $ μ^0_q $ of $ \Gtl_q \cA $ is the sum over all puncture contributions:
\begin{equation*}
μ^0_q ≔ \sum_{q ∈ M} q ℓ_q.
\end{equation*}
In yet other words, we can describe the individual curvature of an arc $ a ∈ \cA $. It carries curvature $ μ^0_{q, a} $ given as the sum of the two turns around its endpoints, multiplied by the deformation parameters $ h(a), t(a) ∈ ℂ⟦M⟧ $ associated with the arc's endpoints.

In order to define the products on $ \Gtl_q \cA $, it suffices to describe them on basis elements of $ \Gtl \cA $. The continuous $ ℂ⟦M⟧ $-multilinear extension is then automatic, see \autoref{rem:2Bprelim-ainfty-continuityautomatic}. We are now ready to state the definition of $ \Gtl_q \cA $:

\begin{definition}
Let $ \cA $ be a full arc system with [NMDC]. The \emph{deformed gentle algebra} $ \Gtl_q \cA $ is the deformation of $ \Gtl \cA $ over $ ℂ⟦M⟧ $ given by:
\begin{itemize}
\item curvature $ μ^0_q = \sum_{q ∈ M} q ℓ_q $,
\item differential $ μ^1_q = 0 $ still vanishing,
\item product $ μ^2_q = μ^2 $ is not deformed,
\item higher products $ μ^{≥3}_q $ as follows: Let $ α_1, …, α_k $ be an angle sequence and let $ β, γ $ be angles such that $ βα_1 ≠ 0 $ and $ α_k γ ≠ 0 $. Then set
\begin{align*}
μ^k_q (β α_k, …, α_1) &= \sum_{D ∈ M(α_1, …, α_k)} q_D β, \\
μ^k_q (α_k, …, α_1 γ) &= \sum_{D ∈ M(α_1, …, α_k)} (-1)^{|γ|} q_D γ.
\end{align*}
\end{itemize}
\end{definition}

\begin{example}
For reference, let us go through a few example evaluations: The torus of \autoref{fig:prelim-gtl-rectified-torus} has one puncture, two arcs and four interior angles. The gentle algebra $ \Gtl \cA $ therefore has two objects and four generators of the morphism spaces. The figure also includes a list of basis elements for $ \Gtl \cA $. What is the deformation $ \Gtl_q \cA $ here? Since there is just one puncture $ q ∈ M = \{q\} $, the deformation base for $ \Gtl_q \cA $ is $ B = ℂ⟦q⟧ $. The arcs $ a $ and $ b $ get curvature
\begin{equation*}
μ^0_a = q αβγδ + q γδαβ \text{ and } μ^0_b = q βγδα + q δαβγ.
\end{equation*}
The product $ μ^2 $ remains non-deformed, for example we still have $ μ^2_q (δ, γ) = - δγ $ and $ μ^2_q (β, α) = 0 $ as in the non-deformed case, but we can now also insert parameters as in $ μ^2_q (q δ, q^2 γ) = - q^3 δγ $. The higher product $ μ^3 $ remains zero, because there are no triangles. The higher product $ μ^4 (δ, γ, β, α) = \id_b $ has still the non-deformed value, and $ μ^6 (γ, βγ, β, α, δα, δ) = \id_a $. Deformed products appear for example in
\begin{equation*}
μ^8_q (δ, γδ, γ, βγ, β, αβ, α, δα) = q \id_b,
\end{equation*}
the sequence is inscribed in a 2-by-2 rectangle covering the puncture once. More generally, we have the $ (m-1) $-by-$ (n-1) $ rectangles covering the puncture $ q $ precisely $ (m-1)(n-1) $ times:
\begin{equation*}
μ^{2m+2n}_q (δ, \underbrace{γδ, …, γδ}_{m-1}, γ, \underbrace{βγ, …, βγ}_{n-1}, β, \underbrace{αβ, …, αβ}_{m-1}, α, \underbrace{δα, …, δα}_{n-1}) = q^{(m-1)(n-1)}.
\end{equation*}
\end{example}

\begin{remark}
We view $ \Gtl \cA $ as a $ ℤ/2ℤ $-graded $ A_∞ $-category and $ \Gtl_q \cA $ as a $ ℤ/2ℤ $-graded deformation. Let us explain why the deformed products $ μ^k_q $ satisfy the requirement to be of degree $ 2 - k $. First, the curvature $ μ^0_{q, a} $ on every arc $ a ∈ \cA $ consists of a full turn and is automatically even. Second, regard e.g.~the deformed higher product
\begin{equation*}
μ^k_q (β α_k, …, α_1) = q_1 … q_m β.
\end{equation*}
Here the angles $ α_1, …, α_k $ are a disk sequence possibly covering the punctures $ q_1, …, q_m $. In fact, the disk sequence ends at the opposite side of the arc as it started, so the total reduced degree
\begin{equation*}
‖α_1‖ + … + ‖α_k‖
\end{equation*}
is still even, which means that $ β $ has the same parity as $ α_1, …, β α_k $ plus $ 2 - k $. This affirms that the products of $ \Gtl_q \cA $ have the right degree.
\end{remark}

\begin{remark}
In the $ ℤ $-graded world, the deformation $ \Gtl_q \cA $ does not exist. Indeed, the curvature $ μ^0_{q, a} $ of an arc $ a ∈ \cA $ is a full turn and its degree depends on the vector field used for the $ ℤ $-grading. It is still possible in the $ ℤ $-grading to define a deformation $ \Gtl_q \cA $ which only includes those punctures where the full turn has degree 2. This would then also give the deformed higher products the right degree. But such a deformation has far fewer deformation parameters and is less interesting than the $ ℤ/2ℤ $-graded version.
\end{remark}

In \paperone, we have already defined $ \Gtl_q \cA $ from a slightly more general starting point. In fact, the starting point of \paperone\ is a deformation base $ B $ and a deformation parameter $ r ∈ \mathfrak{m} Z(\Gtl \cA) $, where $ Z(\Gtl \cA) $ denotes the center of $ \Gtl \cA $ as an ordinary algebra. To the datum of $ r $, the construction in \paperone\ associates a deformed $ A_∞ $-structure $ \bbmu^r $. In order to obtain the specific case of $ μ_q $ from this construction, we have to put $ B = ℂ⟦M⟧ $ and $ r = \sum_{q ∈ M} q ℓ_q $. We have checked in \paperone\ that $ \Gtl_r \cA $ satisfies the curved $ A_∞ $-relations. In particular, this holds for the special case $ \Gtl_q \cA $:

\begin{theorem}[{\paperone}]
Let $ \cA $ be a full arc system with [NMDC]. Then $ \Gtl_q \cA $ is an $ A_∞ $-deformation of $ \Gtl \cA $.
\end{theorem}

\subsection{Zigzag paths}
\label{sec:prelim-gtl-zigzag}
In this section, we recall the notion of zigzag paths. These are combinatorical tools defined specifically for dimers. The idea is to follow the arrows of a dimer by alternatingly turning left and right. In the presentation of zigzag paths, we mostly follow \cite{Bocklandt} and \cite{Broomhead}.

Zigzag paths appear in this paper for two reasons: First, one uses them to define whether a dimer $ Q $ is zigzag consistent or not. Second, zigzag paths themselves can be realized as twisted complexes in $ \Tw\Gtl_q Q $, and the aim of this paper is to compute the minimal model of this category. Of course, it is not a coincidence that zigzag paths appear twice: The zigzag consistency of $ Q $ will help us perform the minimal model calculations of zigzag paths by ruling out lots of difficult cases. We recall zigzag paths as follows:

\begin{definition}
\label{def:prelim-gtl-zigzag-def}
Let $ Q $ be a dimer. A \emph{zigzag path} $ L $ is an infinite path $ … a_2 a_1 a_0 a_{-1} a_{-2} … $ of arcs in $ Q $ together with an alternating choice of “left” or “right” for every $ i ∈ ℕ $ such that
\begin{itemize}
\item $ a_{i+1} a_i $ lies in a clockwise polygon if $ i $ is assigned “right”,
\item $ a_{i+1} a_i $ lies in a counterclockwise polygon if $ i $ is assigned “left”.
\end{itemize}
We also say that $ L $ \emph{turns left} at $ a_i $ if $ i $ is assigned “left” and \emph{turns right} if $ a_i $ is assigned “right”. Two zigzag paths are identified if their paths including left/right indications differ only by integer shift.
\end{definition}

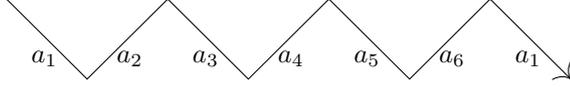
\begin{figure}
\centering
\begin{tikzpicture}
\path[draw, -{To[scale=2]}] (0, 0) -- ++(315:1.5) node[near end, left] {$ a_1 $} coordinate[midway] (alpha1-end) -- ++(45:1.5) node[near start, right] {$ a_2 $} coordinate[midway] (alpha1-start) -- ++(315:1.5) node[near end, left] {$ a_3 $} coordinate[midway] (alpha3-end) -- ++(45:1.5) node[near start, right] {$ a_4 $} coordinate[midway] (alpha3-start) -- ++(315:1.5) node[near end, left] {$ a_5 $} coordinate[midway] (alpha5-end) -- ++(45:1.5) node[near start, right] {$ a_6 $} coordinate[midway] (alpha5-start) -- ++(315:1.5) node[near end, left] {$ a_1 $} coordinate[midway] (additional);
\end{tikzpicture}
\caption{A zigzag path $ L $}
\label{fig:prelim-gtl-zigzag-sketch}
\end{figure}

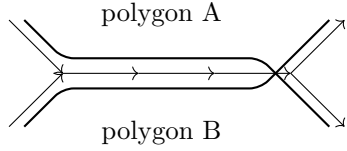
\begin{figure}
\centering
\begin{subfigure}{0.3\linewidth}
\centering
\begin{tikzpicture}
\path[draw, ->] (0, 0) -- ++(45:1) coordinate (1);
\path[draw, ->] (1) -- ++(right:1) coordinate (2);
\path[draw, ->] (2) -- ++(right:1) coordinate (3);
\path[draw, ->] (3) -- ++(right:1) coordinate (4);
\path[draw, ->] (4) -- ++(45:1) coordinate (5);
\path[draw, ->] (4) -- ++(315:1) coordinate (6);
\path[draw, <-] (1) -- ++(135:1) coordinate (7);
\path[draw, thick, rounded corners] (0.2, 0) -- ($ (1) + (0, -0.2) $) -- ($ (4) + (-0.4, -0.2) $) node[midway, below, shift={(0, -0.3)}] {polygon B} -- ($ (5) + (-0.2, 0) $);
\path[draw, thick, rounded corners] ($ (7) + (0.2, 0) $) -- ($ (1) + (0, 0.2) $) -- ($ (4) + (-0.4, 0.2) $) node[midway, above, shift={(0, 0.3)}] {polygon A} -- ($ (6) + (-0.2, 0) $);
\end{tikzpicture}
\end{subfigure}
\caption{Despite sharing multiple arcs, the two strands continue differently and do not finish a period.}
\label{fig:zigzag-period-bad}
\end{figure}

Since $ Q $ is finite, every zigzag path is periodic and has a shortest period $ i_0 ∈ ℕ $, which is defined as the smallest integer such that the zigzag path is invariant under shift by $ i_0 $. The period is not necessarily reached when an arc reappears in the zigzag path. The path may namely continue in a different way beyond that arc. In general, the period need not even be reached when a whole sequence of arcs reappears in the zigzag path. An example is depicted in \autoref{fig:zigzag-period-bad}.

\begin{definition}
The \emph{length} of a zigzag path is the shortest period $ i_0 ∈ ℕ $.
\end{definition}

Slightly simplified, a zigzag path $ L $ is a path in $ Q $ that turns alternatingly maximally right and maximally left in $ Q $. The typical shape of a zigzag path is drawn in \autoref{fig:prelim-gtl-zigzag-sketch}. If every puncture of $ Q $ has valence at least 4, then a path of length two cannot simultaneously lie in the boundary of a clockwise and a counterclockwise polygon. In this case, the left/right indication for zigzag paths is a superfluous part of the datum of a zigzag path. For other dimers $ Q $, the left/right indication is very important. An example is the $ n $-punctured sphere $ Q_n $ of \autoref{fig:prelim-std-dimers-sphere}. If $ n $ is odd, then $ Q_n $ has only one zigzag path, its length is $ 2n $. If $ n $ is even, then $ Q_n $ has two zigzag paths each of length $ n $. This way, we deviate slightly from the definition of \cite{Bocklandt}.

\subsection{Geometric consistency}
In this section, we recall what it means for a dimer to be geometrically consistent. This notion is important for the paper, because we will permanently work with a fixed geometrically consistent dimer, see \autoref{conv:alpha0-direction}. Geometric consistency is a specific instance of various consistency conditions which can be imposed on dimers. A summary can be found in \cite{Bocklandt-consistency}, which we also follow here. In this section, we recall universal covers and zigzag rays and then define geometric consistency.

As first step, we recall the universal cover of $ Q $. Recall that $ Q $ itself consists of a closed surface $ |Q| $ together with an arc system that satisfies the dimer property. Regard the universal cover $ |\tilde Q| $ of the closed surface $ |Q| $. We can construct an arc system on $ |\tilde Q| $ by lifting all punctures and arcs to the universal cover, in all possible ways. This gives an (infinite) arc system on $ |\tilde Q| $ which also has the dimer property. The universal cover together with its lifted arc system is simply denoted $ \tilde Q $.

\begin{figure}
\centering
\begin{subfigure}{0.35\linewidth}
\centering
\begin{tikzpicture}[scale=2]
\path[draw] (0, 0) -- ++(45:1) node[midway, left] {$ a_2^3 $} -- ++(315:0.7) node[midway, left] {$ a_1^3 $} coordinate (B);
\path[draw, -{To[width=0.2cm, length=0.3cm]}] (B) -- ++(up:1.2) node[midway, left] {$ a $} coordinate (A);
\path[draw] (A) -- ++(315:0.7) node[midway, right] {$ a_1^2 $} -- ++(45:1) node[midway, left] {$ a_2^2 $};
\path[draw] (A) -- ++(225:0.7) node[midway, left] {$ a_1^1 $} -- ++(135:1) node[midway, right] {$ a_2^1 $};
\path[draw] (A) -- ++(down:1.2) node[pos=0.6, right] {} -- ++(45:0.7) node[midway, right] {$ a_1^4 $} -- ++(315:1) node[midway, right] {$ a_2^4 $};
\end{tikzpicture}
\caption{Zigzag rays starting at $ a ∈ \tilde Q_1 $}
\label{fig:prelim-consistency-rays}
\end{subfigure}
\begin{subfigure}{0.35\linewidth}
\centering
\begin{tikzpicture}
\path[draw, ->] (0, 0) -- ++(45:1) coordinate (1);
\path[draw, ->] (1) -- ++(right:1) coordinate (2);
\path[draw, ->] (2) -- ++(right:1) coordinate (3);
\path[draw, ->] (3) -- ++(right:1) coordinate (4);
\path[draw, ->] (4) -- ++(45:1) coordinate (5);
\path[draw, ->] (4) -- ++(315:1) coordinate (6);
\path[draw, <-] (1) -- ++(135:1) coordinate (7);
\path ($ (1) + (0, -0.2) $) -- ($ (4) + (-0.4, -0.2) $) node[midway, below, shift={(0, -0.1)}] {polygon B};
\path ($ (1) + (0, 0.2) $) -- ($ (4) + (-0.4, 0.2) $) node[midway, above, shift={(0, 0.1)}] {polygon A};
\end{tikzpicture}
\caption{Not geometrically consistent}
\end{subfigure}
\begin{subfigure}{0.2\linewidth}
\centering
\begin{tikzpicture}
\begin{scope}[dashed, gray]
\path[draw, <-] (0.1, 0) -- (1.9, 0);
\path[draw, ->] (2, 0.1) -- (2, 1.9);
\path[draw, <-] (1.9, 2) -- (0.1, 2);
\path[draw, <-] (0, 0.1) -- (0, 1.9);
\path[draw, ->] (0.1, 0.1) -- (0.9, 0.9);
\path[draw, <-] (1.9, 0.1) -- (1.1, 0.9);
\path[draw, ->] (1.9, 1.9) -- (1.1, 1.1);
\path[draw, <-] (0.1, 1.9) -- (0.9, 1.1);
\end{scope}
\path[draw, thick, rounded corners] (-0.1, 0) -- (-0.1, 2.2) -- (1, 1.1) -- (2.1, 2.2) -- (2.1, 0);
\path[draw, thick, rounded corners] (-0.2, 2) -- (-0.2, -0.3) -- (1, 0.9) -- (2.2, -0.3) -- (2.2, 2);
\end{tikzpicture}
\caption{Not geometrically consistent}
\label{fig:prelim-consistency-pattern}
\end{subfigure}
\caption{On consistency}
\end{figure}
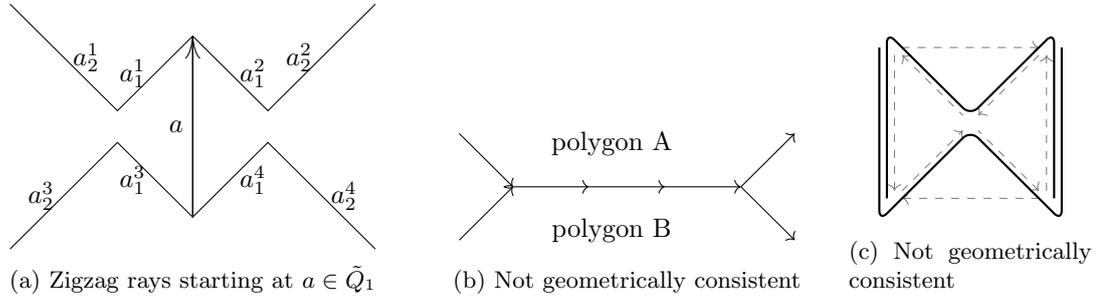

As second step, we recall the notion of zigzag rays. In contrast to zigzag paths, zigzag rays only run in one direction, starting from a given arc. Since we only need zigzag rays in the context of the universal cover, let us directly formulate their definition in $ \tilde Q $. The four zigzag rays starting at an arc $ a ∈ \tilde Q_1 $ are depicted in \autoref{fig:prelim-consistency-rays}.

\begin{definition}
Let $ Q $ be a dimer, $ \tilde Q $ its universal cover and $ a ∈ \tilde Q_1 $ an arc. Then the four zigzag rays starting at $ a $ are the sequences of arcs $ (a_i^1)_{i≥0} $, $ (a_i^2)_{i≥0} $, $ (a_i^3)_{i≥0} $ and $ (a_i^4)_{i≥0} $ in $ \tilde Q $ determined by $ a_0^1 = a_0^2 = a_0^3 = a_0^4 = a $ and the following properties:
\begin{itemize}
\item The sequences $ (a_i^1) $ and $ (a_i^2) $ satisfy $ h(a_i^{1/2}) = t(a_{i+1}^{1/2}) $.
\item The sequences $ (a_i^3) $ and $ (a_i^4) $ satisfy $ t(a_i^{3/4}) = h(a_{i+1}^{3/4}) $.
\item The path $ a_{i+1}^{1/2} a_i^{1/2} $ lies in the boundary of a counterclockwise polygon when $ i $ is odd/even, and clockwise when $ i $ is even/odd.
\item The path $ a_i^{3/4} a_{i+1}^{3/4} $ lies in the boundary of a counterclockwise polygon when $ i $ is odd/even, and clockwise when $ i $ is even/odd.
\end{itemize}
\end{definition}

A dimer is geometrically consistent if the zigzag rays starting with an arc $ a $ in the universal cover intersect nowhere, except at $ a $ itself. The precise definition reads as follows:

\begin{definition}
Let $ Q $ be a dimer. Then $ Q $ is \emph{geometrically consistent} if for every $ a ∈ \tilde Q_1 $ the four zigzag rays $ (a_i^1) $, $ (a_i^2) $, $ (a_i^3) $ and $ (a_i^4) $ satisfy the following property: Whenever $ a_i^k = a_j^l $, then $ i = j $ and $ k = l $, or $ i = j = 0 $.
\end{definition}

\begin{example}
A dimer $ Q $ on a sphere is never geometrically consistent, because $ \tilde Q = Q $ and therefore any zigzag rays in $ \tilde Q $ intersect after a while. There are plenty of geometrically consistent dimers on surfaces of genus $ g ≥ 1 $ though. For example, the $ n $-punctured torus dimer of \autoref{fig:prelim-std-dimers-torus} is geometrically consistent. Indeed, the universal cover of the torus is the real plane, and the torus dimer lifts to horizontally and vertically repeated copies of \autoref{fig:prelim-std-dimers-torus}. The zigzag rays then run away in different directions in the plane without ever coming closer to each other again. This geometry is a typical example of the toric zigzag fan, see for example \cite{Wong}.
\end{example}

\begin{remark}
If $ Q $ is geometrically consistent, then a zigzag path $ L $ on $ Q $ may return to an arc twice, however the segment of $ L $ between both occurrences is not allowed to be contractible. If it were contractible, then this segment would constitute a zigzag ray cutting itself (the case $ i = j $ and $ k = l $), contradicting geometric consistency.
\end{remark}

\begin{remark}
A geometrically consistent dimer automatically satisfies the [NMDC] condition.
\end{remark}

Geometric consistency is the strongest consistency condition one can require, apart from R-charge consistency. Indeed, geometric consistency is by definition a stronger version of so-called zigzag consistency, which in turn is known to be stronger than cancellation consistency, see \cite{Bocklandt-consistency}:

\begin{center}
\begin{tikzpicture}
\path (0, 0) node (A) {geometric consistency} (5, 0) node (B) {zigzag consistency} (10, 0) node (C) {cancellation consistency};
\path (A.east) -- (B.west) node[midway] {\Large $ \Longrightarrow $};
\path (B.east) -- (C.west) node[midway] {\Large $ \Longrightarrow $};
\path ($ (B) + (0, 1) $) node (D) {R-charge consistency};
\path[draw, double equal sign distance] (B.north) -- (D.south) node[midway, left] {$ g = 1 $};
\end{tikzpicture}
\end{center}

\subsection{Terminology for arcs and angles}
\label{sec:prelim-terminology}
In this section, we introduce technical terminology that we will use throughout the paper. This terminology is important to describe exactly what happens where in a dimer. It bears no mathematical creativity but is unavoidable for the sake of concise language.

\begin{figure}
\centering
\begin{subfigure}{0.2\linewidth}
\centering
\begin{tikzpicture}
\path[draw, ->] (-1, 1) -- (0, 0) coordinate[midway] (2);
\path[draw, ->] (0, 0) -- (1, 1) coordinate[midway] (1);
\path[draw, ->, bend right] (1) to node[midway, below] {$ α $} (2);
\end{tikzpicture}
\caption{tail/head}
\end{subfigure}
\begin{subfigure}{0.2\linewidth}
\centering
\begin{tikzpicture}
\path[draw, <-] (-1, 1) -- (0, 0) coordinate[midway] (2);
\path[draw, <-] (0, 0) -- (1, 1) coordinate[midway] (1);
\path[draw, ->, bend right] (1) to node[midway, below] {$ α $} (2);
\end{tikzpicture}
\caption{head/tail}
\end{subfigure}
\begin{subfigure}{0.2\linewidth}
\centering
\begin{tikzpicture}
\path[draw, <-] (-1, 1) -- (0, 0) coordinate[midway] (2);
\path[draw, ->] (0, 0) -- (1, 1) coordinate[midway] (1);
\path[draw, ->, bend right] (1) to node[midway, below] {$ α $} (2);
\end{tikzpicture}
\caption{tail/tail}
\end{subfigure}
\begin{subfigure}{0.2\linewidth}
\centering
\begin{tikzpicture}
\path[draw, ->] (-1, 1) -- (0, 0) coordinate[midway] (2);
\path[draw, <-] (0, 0) -- (1, 1) coordinate[midway] (1);
\path[draw, ->, bend right] (1) to node[midway, below] {$ α $} (2);
\end{tikzpicture}
\caption{head/head}
\end{subfigure}
\caption{The angle $ α $ starts/ends at …}
\label{fig:prelim-terminology-angle-headtail}
\end{figure}
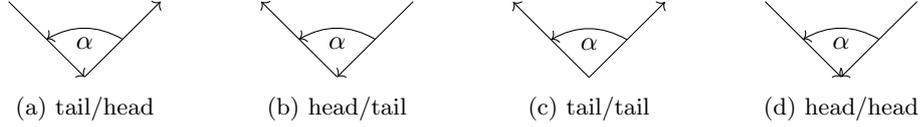

The first notion is that of an arc incidence. This is comparable to half-edges in a ribbon graph. Half-edges are not only a useful tool to describe graphs where one end of some edges is missing, but half-edges are also handy to describe incidences in a graph. Whenever we would like to sum over all edges incident at a given node, letting every loop contribute two (distinct) terms, the right entity to sum over is the set of half-edges incident at the node. Similarly in a dimer $ Q $, we would typically like to have a set of all incidences of arcs, where loops contribute both their “head part” and their “tail part”. With the terminology of head parts and tail parts, we can also talk about whether an angle starts at the head part or tail part of an arrow. This piece of terminology is depicted in \autoref{fig:prelim-terminology-angle-headtail}.

\begin{definition}
An \emph{arc incidence} at a puncture $ q ∈ Q_0 $ is either an incident \emph{head part} or an incident \emph{tail part} of some arc.
\end{definition}

For instance, a loop $ a ∈ Q_1 $ with $ h(a) = t(a) = q $ has two arc incidences at the puncture $ h(a) = t(a) $. The sample puncture in \autoref{fig:terminology-incidences} has four arc incidences. Correspondingly, there are four indecomposable angles around the puncture. In that figure, the loop is intended to be topologically nontrivial, indicated by the dots “$ … $”.

\begin{figure}
\centering
\begin{tikzpicture}
\path[draw, ->] (0.1, -0.1) to[out=300, in=270] (1.5, 0) to[out=90, in=60] (0.1, 0.1);
\path (0.75, 0) node {$ … $};
\path[fill] (0, 0) circle[radius=0.05];
\path[draw, ->] (-0.2, 1) -- (0, 0.1);
\path[draw, ->] (0, -0.1) -- (-0.2, -1);
\end{tikzpicture}
\caption{A puncture with four arc incidences}
\label{fig:terminology-incidences}
\end{figure}
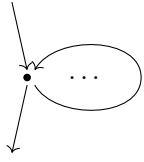

Let us now introduce some terminology for angles in $ Q $. Angles always have both an algebraic interpretation as basis morphisms for $ \Gtl Q $ and a geometric interpretation as winding around punctures in the surface $ Q $. We will therefore use double terminology from time to time: In algebraic contexts, we say an angle is an identity if it is the identity $ \id_a $ of some arc $ a ∈ Q_1 $. In geometric context, we call such an angle empty and all other angles non-empty. For instance, a typical usage in a geometric context would be to say that a certain angle $ α $ is non-empty and smaller than a full turn.

\begin{definition}
Let $ α $ be an angle in $ Q $. Then $ α $ is an \emph{empty} angle if it is the identity of some arc. Otherwise $ α $ is a \emph{non-empty} angle.
\end{definition}

Given an angle $ α $, we would like to distinguish whether it is composed of multiple smaller angles or not. By definition of the angles in $ Q $, the smallest units are the interior angles of polygons. This already gives us terminology for a geometric context: We can simply ask whether a given angle $ α $ is the interior angle of some polygon or not. We however also need terminology for the algebraic context. Examples are depicted in \autoref{fig:prelim-terminology-decomposable}. We fix terminology as follows:

\begin{definition}
An angle is \emph{decomposable} if it is the composition of two non-empty angles. An angle is \emph{indecomposable} if it is non-empty and not decomposable.
\end{definition}

\begin{remark}
An non-empty angle is indecomposable if it is an interior angle of some polygon, and indecomposable otherwise. We regard empty angles as neither decomposable nor indecomposable.
\end{remark}

\begin{figure}
\centering
\begin{subfigure}{0.3\linewidth}
\centering
\begin{tikzpicture}
\foreach \i in {0, 60, 120, 180, 240, 300} \path[draw] (0, 0) -- ++(\i:1.5);
\path[draw, ->] (60:0.9) arc(60:360:0.9);
\end{tikzpicture}
\end{subfigure}
\begin{subfigure}{0.3\linewidth}
\centering
\begin{tikzpicture}
\foreach \i in {0, 60, 120, 180, 240, 300} \path[draw] (0, 0) -- ++(\i:1.5);
\path[draw, ->] (60:0.9) arc(60:120:0.9);
\end{tikzpicture}
\end{subfigure}
\caption{A decomposable and an indecomposable angle}
\label{fig:prelim-terminology-decomposable}
\end{figure}
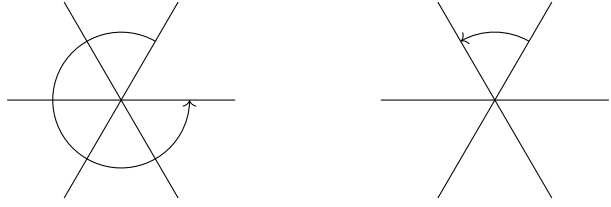

\begin{figure}
\centering
\begin{subfigure}{0.45\linewidth}
\centering
\begin{tikzpicture}[scale=2]
\path[draw] (0, 0) -- ++(45:1) node[midway, left] {1} -- ++(315:0.7) node[midway, left] {2} -- ++(up:1.2) node[midway, left] {3} coordinate (A) -- ++(315:0.7) node[midway, right] {4} -- ++(45:1) node[midway, left] {5};
\path[draw] (A) ++(0.1, 0.2) -- ++(215:0.7) node[midway, left] {41} -- ++(135:1) node[midway, right] {40};
\path[draw] (A) ++(0.1, 0.2) -- ++(down:1.5) node[pos=0.6, right] {42} -- ++(45:0.7) node[midway, right] {43} -- ++(315:1) node[midway, right] {44};
\end{tikzpicture}
\caption{Arc appearing twice in $ L $ with different index}
\label{fig:prelim-terminology-indexedarc}
\end{subfigure}
\begin{subfigure}{0.35\linewidth}
\begin{tikzpicture}
\path[draw] (0, 0) -- ++(45:1) -- ++(315:1) coordinate (A) -- ++(45:1) -- ++(315:1) -- ++(45:1) -- ++(315:1) coordinate (B) -- ++(45:1);
\path[draw, decoration={brace, amplitude=8pt, mirror}, decorate] ($ (A) + (0, -0.2) $) -- ($ (B) + (0, -0.2) $) node[midway, below, shift={(0, -0.3)}] {segment};
\end{tikzpicture}
\caption{An indexed segment}
\label{fig:prelim-terminology-indexedsegment}
\end{subfigure}
\caption{Terminology in a zigzag path}
\end{figure}
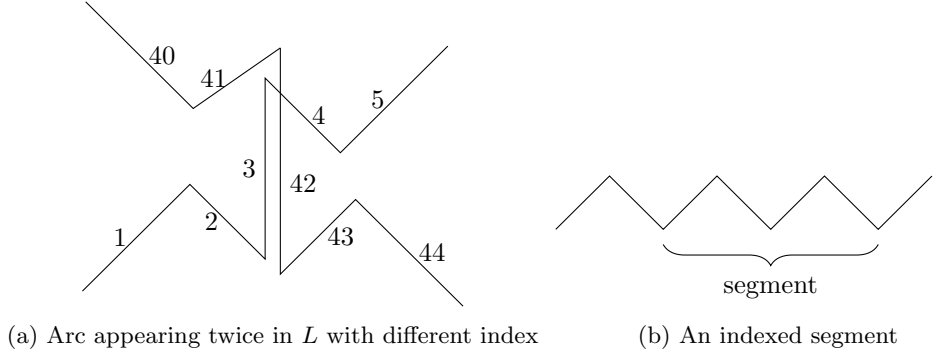

Let us introduce terminology for locations on zigzag paths. Loosely speaking, we want to define an “indexed arc” as an arc $ a ∈ Q_1 $ lying on $ L $, but remember whether $ L $ turns left or right after $ a $. For example, let $ 2k $ be the length of $ L $, then $ L $ has precisely $ 2k $ indexed arcs. \autoref{fig:prelim-terminology-indexedarc} features a visual explanation: Some arc $ a_3 = a_{42} $ appears twice while traversing $ L $, one time at index $ 3 $ and one time at index $ 42 $. The arc itself is the same in $ Q $, but different as indexed arcs of $ L $. This amount of precision gives rise to further names for relative positions on $ L $. For instance, we can regard indexed segments, depicted in \autoref{fig:prelim-terminology-indexedsegment}. We fix terminology as follows:

\begin{definition}
Let $ L $ be a zigzag path, given by an infinite path $ … a_1 a_0 a_{-1} … $ together with left/right indications.
\begin{itemize}
\item An \emph{indexed arc} on $ L $ is a tuple $ (a_i, i) $ consisting of one of the arcs on $ L $ together with its index modulo the period length of the zigzag path.
\item The \emph{next arc} after $ (a_i, i) $ is the indexed arc $ (a_{i+1}, i+1) $.
\item The \emph{previous arc} before $ (a_i, i) $ is the indexed arc $ (a_{i-1}, i-1) $.
\item Two \emph{consecutive} indexed arcs are two indexed arcs on $ L $ which can be written in the form $ (a_{i+1}, i+1) $ and $ (a_i, i) $ or the other way around.
\item An \emph{indexed segment} of length $ k $ on $ L $ is the datum of a tuple $ (a_i, …, a_{i+k-1}, i) $ of arcs on $ L $, remembering the index $ i $ modulo the period length of $ L $.
\end{itemize}
\end{definition}

\begin{figure}
\centering
\begin{subfigure}{0.3\linewidth}
\centering
\begin{tikzpicture}
\path[draw] (0:1) -- (60:1) coordinate[pos=0.3] (1-start) coordinate[pos=0.7] (2-end) -- (120:1) coordinate[pos=0.3] (2-start) coordinate[pos=0.7] (3-end) -- (180:1) coordinate[pos=0.3] (3-start) coordinate[pos=0.7] (4-end) -- (240:1) coordinate[pos=0.3] (4-start) coordinate[pos=0.7] (5-end) -- (300:1) coordinate[pos=0.3] (5-start) coordinate[pos=0.7] (6-end) -- (0:1) coordinate[pos=0.3] (6-start) coordinate[pos=0.7] (1-end);
\path[draw] (-1, 0) -- (1, 0);
\path[draw] (120:1) -- (300:1);
\path[draw] (60:1) -- (240:1);
\path[draw, thick] (240:1) -- (300:1);
\path[draw, thick] (300:1) -- (0:1);
\path[draw] (300:1) -- ++(right:1);
\path[draw] (300:1) -- ++(300:1);
\path[draw] (300:1) -- ++(240:1) coordinate[pos=0.3] (out-end);
\path[draw, ->, bend right=60] (1-start) to (1-end);
\path[draw, ->, bend right=60] (2-start) to (2-end);
\path[draw, ->, bend right=60] (3-start) to (3-end);
\path[draw, ->, bend right=60] (4-start) to (4-end);
\path[draw, ->, bend right=60] (5-start) to node[at start, below left] {1} (5-end);
\path[draw, ->, bend right=90, looseness=1.5] (6-start) to node[at start, right] {6} (out-end);
\path[fill] (300:1) circle[radius=0.05];
\begin{scope}[shift={(0, -2.5)}]
\path node[align=center] {slots inside: 2 \\ slots outside: 4};
\end{scope}
\end{tikzpicture}
\caption{A final-out disk}
\end{subfigure}
\begin{subfigure}{0.3\linewidth}
\centering
\begin{tikzpicture}
\path[draw] (0:1) -- (60:1) coordinate[pos=0.3] (1-start) coordinate[pos=0.7] (2-end) -- (120:1) coordinate[pos=0.3] (2-start) coordinate[pos=0.7] (3-end) -- (180:1) coordinate[pos=0.3] (3-start) coordinate[pos=0.7] (4-end) -- (240:1) coordinate[pos=0.3] (4-start) coordinate[pos=0.7] (5-end) -- (300:1) coordinate[pos=0.3] (5-start) coordinate[pos=0.7] (6-end) -- (0:1) coordinate[pos=0.3] (6-start) coordinate[pos=0.7] (1-end);
\path[draw] (-1, 0) -- (1, 0);
\path[draw] (120:1) -- (300:1);
\path[draw] (60:1) -- (240:1);
\path[draw, thick] (240:1) -- (300:1);
\path[draw, ->, bend right=60] (1-start) to (1-end);
\path[draw, ->, bend right=60] (2-start) to (2-end);
\path[draw, ->, bend right=60] (3-start) to (3-end);
\path[draw, ->, bend right=60] (4-start) to (4-end);
\path[draw, ->, bend right=60] (5-start) to node[at start, below left] {1} (5-end);
\path[draw, ->, bend right=60] (6-start) to node[at end, below right] {6} (6-end);
\begin{scope}[shift={(0, -2.5)}]
\end{scope}
\end{tikzpicture}
\caption{An all-in disk}
\end{subfigure}
\begin{subfigure}{0.3\linewidth}
\centering
\begin{tikzpicture}
\path[draw] (0:1) -- (60:1) coordinate[pos=0.3] (1-start) coordinate[pos=0.7] (2-end) -- (120:1) coordinate[pos=0.3] (2-start) coordinate[pos=0.7] (3-end) -- (180:1) coordinate[pos=0.3] (3-start) coordinate[pos=0.7] (4-end) -- (240:1) coordinate[pos=0.3] (4-start) coordinate[pos=0.7] (5-end) -- (300:1) coordinate[pos=0.3] (5-start) coordinate[pos=0.7] (6-end) -- (0:1) coordinate[pos=0.3] (6-start) coordinate[pos=0.7] (1-end);
\path[draw] (-1, 0) -- (1, 0);
\path[draw] (120:1) -- (300:1);
\path[draw] (60:1) -- (240:1);
\path[draw, thick] (180:1) -- (240:1);
\path[draw, thick] (240:1) -- (300:1);
\path[draw] (240:1) -- ++(left:1);
\path[draw] (240:1) -- ++(240:1);
\path[draw] (240:1) -- ++(300:1) coordinate[pos=0.3] (out-start);
\path[draw, ->, bend right=60] (1-start) to (1-end);
\path[draw, ->, bend right=60] (2-start) to (2-end);
\path[draw, ->, bend right=60] (3-start) to (3-end);
\path[draw, ->, bend right=60] (4-start) to (4-end);
\path[draw, ->, bend right=60] (6-start) to node[at end, below right] {6} (6-end);
\path[draw, ->, bend right=90, looseness=1.5] (out-start) to node[at end, left] {1} (5-end);
\path[fill] (240:1) circle[radius=0.05];
\begin{scope}[shift={(0, -2.5)}]
\path node[align=center] {slots inside: 2 \\ slots outside: 4};
\end{scope}
\end{tikzpicture}
\caption{A first-out disk}
\end{subfigure}
\caption{Illustration of final-out, all-in and first-out disks}
\label{fig:prelim-terminology-firstfinal}
\end{figure}
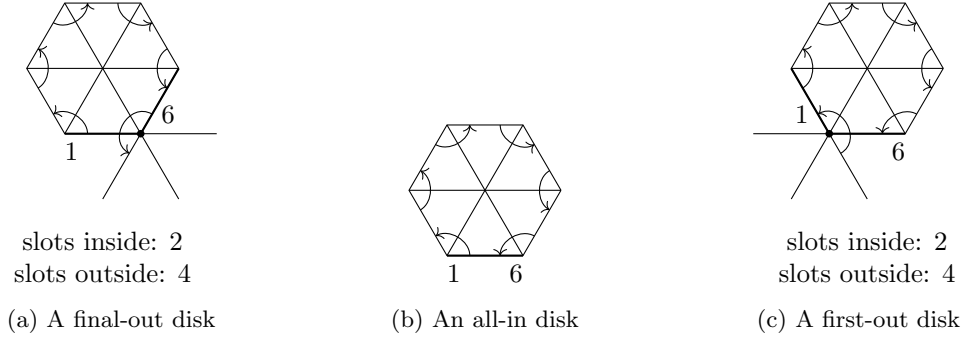

Let us introduce terminology for angle sequences. We already have the very fortunate notion of disk sequences available, but in order to analyze products in $ \Tw\Gtl Q $ we need flexible terminology to distinguish between the two rules that define higher products in $ \Gtl Q $. Let $ α_1, …, α_k $ be a disk sequence. Recall that the discrete immersed disk contained in the data of the disk sequence contributes to the product $ μ^k_{\Gtl_q Q} (α_k, …, α_1) $. Now if $ β $ is an angle such that $ β α_k ≠ 0 $ and $ γ $ is an angle such that $ α_1 γ ≠ 0 $, then the discrete immersed disk also contributes to the products $ μ^k_{\Gtl_q Q} (β α_k, …, α_1) $ and $ μ^k_{\Gtl_q Q} (α_k, …, α_1 γ) $. We want to call these contributions final-out and first-out, respectively. This terminology is depicted in \autoref{fig:prelim-terminology-firstfinal}. A more formal definition reads as follows:

\begin{definition}
\label{def:prelim-terminology-firstfinal}
Let $ Q $ be a dimer. Let $ α_1, …, α_k $ be a disk sequence in $ Q $ with discrete immersed disk $ D $. Let $ β, γ $ be non-empty angles such that $ β α_k ≠ 0 $ and $ α_1 γ ≠ 0 $. Then:
\begin{itemize}
\item The sequence $ α_1, …, β α_k $ together with D is a \emph{final-out} disk. We call $ β $ the outside morphism and $ α_k $ the inside morphism. We call $ t(α_1) $ the first arc and $ t(α_k) $ the final arc.
\item The sequence $ α_1 γ, …, α_k $ together with $ D $ is a \emph{first-out} disk. We call $ γ $ the outside morphism and $ α_1 $ the inside morphism. We call $ h(α_1) $ the first arc and $ h(α_k) $ the final arc.
\item The sequence $ α_1, …, α_k $ together with $ D $ is an \emph{all-in} disk. We call the arc $ t(α_1) = h(α_k) $ the first, equivalently final arc.
\end{itemize}
We may call an angle sequence together with a discrete immersed disk a \emph{some-out disk} if it is first-out or final-out. In the case of a some-out disk, the first and final arc share an endpoint, the \emph{concluding puncture} of the disk. In the case of an all-in disk, the first and final arc coincide, which is the \emph{concluding arc} of the disk.
\end{definition}

Loosely speaking, all contributions to $ μ_{\Gtl_q Q}^{≥3} $ come from first-out, final-out or all-in disks. Some-out means first-out or final-out. For a some-out disk, the concluding puncture is the one around which the first or final angle winds and it is very important. The first and final arcs are those arcs that neighbor the concluding puncture. For an all-in disk, the first and final arcs are the same and this single arc is very important. Whenever we refer to first-out, final-out or all-in disks, we typicall pass the datum of the discrete immersed disk implicitly. In \autoref{fig:prelim-terminology-firstfinal}, the first and final arc are drawn thick and the concluding puncture is marked with a dot. We may use wording like “towards the concluding puncture” when referring to the behavior of a sequence of arcs, viewed in the direction of the concluding puncture.

Last but not least, we shall give some means to measure how large the inside angle is by counting the “slots” inside and outside the disk. The terminology is included in \autoref{fig:prelim-terminology-firstfinal}. We formalize this as follows:

\begin{definition}
Let $ Q $ be a dimer. Let $ α_1, …, α_k $ be a some-out disk and $ γ $ be its inside morphism. Write $ ℓ $ for one full turn around the concluding puncture, starting at the final arc of the disk. Write $ γ = γ' ℓ^n $ for some $ n $ such that $ γ' $ is strictly smaller than one full turn. Take the complementary angle $ (γ')^c $ such that $ (γ')^c γ' = ℓ $.
\begin{itemize}
\item The number of \emph{slots inside} the disk is the number of indecomposable angles that $ γ' $ consists of.
\item The number of \emph{slots outside} the disk is the number of indecomposable angles that $ (γ')^c $ consists of.
\end{itemize}
\end{definition}

\section{Preliminaries of Fukaya categories}
\label{sec:fukaya}
In this section, we recall basics of Fukaya categories. One after another, we recall the construction of the Fukaya pre-category, Fukaya category, relative Fukaya pre-category and relative Fukaya category. The core aim of the paper is to define the category $ \DefZigzagCat $ and interpret its minimal model $ \H\DefZigzagCat $ as a part of the relative Fukaya category. The present section aims to facilitate this understanding by preparing the view from the side of Fukaya categories. We have therefore included a dedicated description of the subcategory of the relative Fukaya category given by so-called zigzag curves in \autoref{sec:prelim-fukaya-zigzag}. Our main references are \cite{Abouzaid, CHL}. We comment on results of Efimov, Sheridan and Perutz.

\subsection{The exact Fukaya pre-category}
\label{sec:prelim-fukaya-pre}
In this section, we review exact Fukaya pre-categories. They are not an immediate necessity for this paper, since we only work with the discrete model $ \Gtl_q Q $. The main result however ties $ \H\Tw\Gtl_q Q $ to the relative Fukaya category, so we will benefit from a review. We follow a combination of the highly recommendable sources \cite{Abouzaid}, \cite{Efimov} and \cite[Chapter 6]{Bocklandt-book}.

In symplectic geometry, one aims at defining a fully-fledged $ A_∞ $-category $ \Fuk X $ from a symplectic manifold $ X $. A Fukaya category is supposed to have closed Lagrangians as objects and intersection points as basis elements for the hom spaces. The products $ μ^{≥2} $ are supposed to be formed from immersed disks between Lagrangians. For Lagrangians lying in general position, this construction works well. It is however not clear what the endomorphism space of a single Lagrangian $ L $ should be. We would expect it to be a finite-dimensional vector space, and it should be equal for all small Hamiltonian deformations of $ L $. This makes the full set of hom spaces and $ A_∞ $ products of a Fukaya category very hard to define.

The difficulty in defining a fully-fledged Fukaya category $ \Fuk X $ has led to the introduction of pre-categories as partial remedy: Products need not be defined on all sequences of morphisms, only on a choice of transversal sequences.

\begin{definition}
\label{def:prelim-fukaya-precat}
An $ A_∞ $-pre-category $ \cat C $ consists of the following data:
\begin{enumerate}
\item a set of objects $ \Ob \cat C $,
\item for every $ N ≥ 1 $ a set $ (\cat C)^N_{\transversal} ⊂ (\Ob \cat C)^N $ of \emph{transversal sequences}, with $ (\cat C)^1_{\transversal} = \Ob \cat C $,
\item for every $ (X, Y) ∈ (\cat C)^2_{\transversal} $ a graded hom space $ \Hom(X, Y) $,
\item for every transversal sequence $ X_1, …, X_{N+1} $ with $ N ≥ 1 $ a degree $ 2-N $ product map
\begin{equation*}
μ^N: \Hom(X_N, X_{N+1}) \tensor … \tensor \Hom(X_1, X_2) → \Hom(X_1, X_{N+1}),
\end{equation*}
\end{enumerate}
such that each subsequence $ (X_{i_1}, …, X_{i_l}) $ with $ 1 ≤ i_1 < … < i_l ≤ n $ of a transversal sequence $ X_1, …, X_n $ is transversal as well, and the $ A_∞ $-relation holds for $ \Hom(X_N, X_{N+1}) \tensor … \tensor \Hom(X_1, X_2) $ whenever $ X_1, …, X_{N+1} $ is a transversal sequence.
\end{definition}

\begin{remark}
Staring at the definition seems to imply that the condition on transversal sequences is arbitrarily weak: Setting $ (\cat C)^N_{\transversal} = ø $ for all $ N ≥ 2 $ is possible, and yields a completely vacuous notion of pre-category. The point of \cite{Efimov} is that if one strengthens the conditions suitably, then giving a pre-$ A_∞ $-category is the same as giving a full $ A_∞ $-category. We will comment on this later on. In particular, we will ensure that our definition of transversal sequences is such that it satisfies the condition in \cite{Efimov}.
\end{remark}

Abouzaid's exposition \cite{Abouzaid} exhibits the Fukaya pre-category of a surface with boundary. In particular, we get from his paper a direct construction of the Fukaya pre-category of a punctured surface, by interpreting the punctures as boundary circles. We deviate from Abouzaid's definition by only including exact Lagrangians in the category. This makes it possible to dispose of the Novikov ring and work over $ ℂ $ instead. We are now ready for the first definitions.

Our aim here is to write down the definition of the exact Fukaya pre-category of a punctured surface, such that it is a pre-category in the sense of \autoref{def:prelim-fukaya-precat}. Before we give the definition, we have to recall several concepts from \cite{Abouzaid}: teardrops, spin structures, unobstructed curves, exact curves, transversal sequences of unobstructed curves, degrees of intersection points, immersed disks between unobstructed curves, and the Abouzaid sign rule. We recall these terminologies one by one.

\begin{definition}
A \emph{teardrop} of a curve $ X: S^1 → S \setminus M $ is an immersion of the monogon $ P_1 $ into $ S \setminus M $ which is bounded by a segment of $ X $, such that the corner coming from $ P_1 $ is convex.
\end{definition}

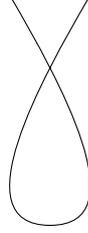
\begin{figure}
\centering
\begin{subfigure}{0.3\linewidth}
\centering
\begin{tikzpicture}
\path[draw] (0, 0) to[out=300, in=0] (0.5, -3) to[out=180, in=240] (1, 0);
\end{tikzpicture}
\caption{Teardrop}
\label{fig:prelim-fukaya-teardrop}
\end{subfigure}
\caption{Terminology for unobstructed curves}
\label{fig:prelim-fukaya-curveterminology}
\end{figure}

A teardrop is depicted in \autoref{fig:prelim-fukaya-teardrop}. In contrast, an interval winding around a puncture does not constitute a teardrop. A curve has a teardrop if and only if it contains an interval that is contractible in $ S \setminus M $.

\begin{definition}
A \emph{spin structure} on a curve consists of putting an arbitrary number of “stars” on distinct points of the curve. We also call these stars the \emph{\# signs} on a curve.
\end{definition}

We regard the stars as $ \# $ signs when we think of them as a negative sign $ -1 $. The number of $ \# $ signs on a curve is arbitrary, but the resulting isomorphism class of the curve will in fact only depend on the parity of this number. In other words, zero or one $ \# $ sign suffice in practice.

\begin{definition}
An \emph{unobstructed curve} in $ S \setminus M $ is a smooth closed immersed curve $ X: S^1 → S \setminus M $ with a choice of spin structure, such that $ X $ is not contractible and does not bound a teardrop.
\end{definition}

Let us now recall the notion of exact curves, a subset of the unobstructed curves. Exact curves serve as objects of the Fukaya (pre-)category. To introduce the notion, put an exact symplectic form $ ω = dθ $ on $ S \setminus M $. The 1-form $ θ $ is then also referred to as the Liouville form.

\begin{definition}
An unobstructed curve $ X: S^1 → S \setminus M $ is \emph{exact} if $ X^* θ $ is an exact 1-form, in other words if
\begin{equation*}
∃f: S^1 → ℝ: \quad X^* θ = df, \quad \text{or} \quad ∫_X θ = 0.
\end{equation*}
\end{definition}

These two conditions are equivalent because $ ∫_X θ = ∫_{S^1} X^* θ $. The latter integral vanishes if and only if $ X^* θ ∈ Ω^1 (S^1) $ has a primitive $ f $.

\begin{definition}
\label{def:prelim-fukaya-pre-transversal}
A sequence $ (X_1, …, X_N) $ of unobstructed curves is \emph{transversal} if
\begin{itemize}
\item For $ i < j $ the curves $ X_i $ and $ X_j $ have only transversal intersection points.
\item For $ i < j < k $ the curves $ X_i, X_j, X_k $ have no triple intersection: $ X_i ∩ X_j ∩ X_k = ø $.
\end{itemize}
\end{definition}

According to the definition, an unobstructed curve is allowed to intersect itself, just as a self-intersection of one unobstructed curve $ X_i $ is allowed to further intersect with a second unobstructed curve $ X_j $.

\begin{figure}
\centering
\begin{subfigure}{0.3\linewidth}
\begin{tikzpicture}
\path[draw] (0, 0) -- (2, 0) node[right] {$ X_1 $};
\path[draw] (1, 0) -- ++(135:1) node[left] {$ X_2 $} (1, 0) -- ++(315:1);
\path[fill] (1, 0) circle[radius=0.05];
\path (0.8, -0.2) node {$ p $};
\end{tikzpicture}
\caption{Transversal intersection}
\end{subfigure}
\begin{subfigure}{0.3\linewidth}
\begin{tikzpicture}
\path[draw] (0, 0) -- (3, 0) node[near start, above] {$ X_1 $};
\path[draw] (1, 1) to[bend left] coordinate[midway] (A) (1, -1);
\path[draw] (2, 1) to[bend left] coordinate[midway] (B) (2, -1);
\path[fill] (A) circle[radius=0.05];
\path[fill] (B) circle[radius=0.05];
\path (1, 1) node[left] {$ X_2 $};
\path (2, 1) node[right] {$ X_3 $};
\end{tikzpicture}
\caption{Transversal sequence}
\end{subfigure}
\caption{Transversality}
\label{fig:prelim-fukaya-transversal}
\end{figure}
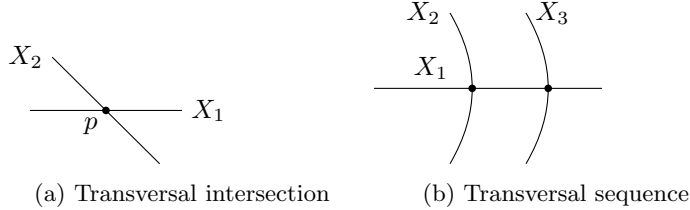

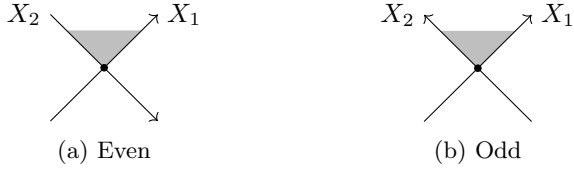
\begin{figure}
\centering
\begin{subfigure}{0.3\linewidth}
\centering
\begin{tikzpicture}
\path[fill, lightgray] (0, 0) -- (45:0.7) -- (135:0.7) -- cycle;
\path[draw] (0, 0) -- ++(135:1) node[left] {$ X_2 $} (0, 0) -- ++(225:1);
\path[draw, ->] (0, 0) -- ++(315:1);
\path[draw, ->] (0, 0) -- ++(45:1) node[right] {$ X_1 $};
\path[fill] (0, 0) circle[radius=0.05];
\end{tikzpicture}
\caption{Even}
\label{fig:prelim-fukaya-intersectioneven}
\end{subfigure}
\begin{subfigure}{0.3\linewidth}
\centering
\begin{tikzpicture}
\path[fill, lightgray] (0, 0) -- (45:0.7) -- (135:0.7) -- cycle;
\path[draw, ->] (0, 0) -- ++(135:1) node[left] {$ X_2 $};
\path[draw] (0, 0) -- ++(225:1) (0, 0) -- ++(315:1);
\path[draw, ->] (0, 0) -- ++(45:1) node[right] {$ X_1 $};
\path[fill] (0, 0) circle[radius=0.05];
\end{tikzpicture}
\caption{Odd}
\label{fig:prelim-fukaya-intersectionodd}
\end{subfigure}
\caption{Intersection degree}
\label{fig:prelim-fukaya-intersectiondeg}
\end{figure}

Next, let us recall the degree assigned to an intersection point $ p ∈ X_1 ∩ X_2 $. The idea is that the intersection $ p ∈ X_1 ∩ X_2 $ serve as generators of $ \Hom(X_1, X_2) $, so we have to assign a degree. Since the surface $ S $ and the curves $ X_1, X_2 $ are oriented, we can distinguish the direction of $ X_1 $ and $ X_2 $ relative to each other at $ p $. The degree we assign is depicted in \autoref{fig:prelim-fukaya-intersectiondeg}. In that figure, the shaded area has no meaning in this definition, but indicates for the convenience of the reader how we are going to use such intersection points as corners of immersed disks. Note that $ p $ can be interpreted both as element of $ \Hom(X_1, X_2) $ and $ \Hom(X_2, X_1) $. In fact, it has opposite parity in both hom spaces.

\begin{definition}
Let $ p $ be a transversal intersection point of $ X_1 $ and $ X_2 $. Then $ p $ as morphism from $ X_1 $ to $ X_2 $ is denoted $ p: X_1 → X_2 $. The morphism $ p: X_1 → X_2 $ is odd if a neighborhood of $ p ∈ S $ can be identified in an oriented way with a neighborhood of the origin in $ ℝ^2 $, mapping $ X_1 $ to the oriented $ x $-axis and $ X_2 $ to the oriented $ y $-axis. Otherwise $ p $ is even.
\end{definition}

Let us recall the notion of smooth immersed disks between unobstructed curves. Despite their name, the disks have corners and are therefore actually polygons. We stick to the classical terminology however. Recall that $ P_{N+1} $ denotes the standard oriented polygon in $ ℝ^2 $, with indexed $ N+1 $ clocNwise indexed edges and $ N+1 $ indexed corners. The $ i $-th corner lies between the $ i $-th and $ (i+1) $-th edge.

\begin{definition}
Let $ X_1, …, X_{N+1} $ be a transversal sequence of $ N+1 ≥ 2 $ unobstructed curves. Let $ p_1, …, p_N $ be a sequence of intersection points $ p_i: X_i → X_{i+1} $ and let $ p ∈ X_1 → X_{N+1} $ be another intersection point. A \emph{smooth immersed disk} with \emph{inputs} $ p_1, …, p_N $ and \emph{output} $ p $ consists of an orientation-preserving polygon immersion $ D: P_{N+1} → S \setminus M $ up to reparametrization, such that
\begin{itemize}
\item the corners of $ D $ are all convex,
\item the $ i $-th edge of $ P_{N+1} $ lands on $ X_i $ for $ 1 ≤ i ≤ N+1 $,
\item the $ i $-th corner of $ P_{N+1} $ lands on $ p_i $.
\end{itemize}
We denote by $ M(p_1, …, p_N, p) $ the set of smooth immersed disks with inputs $ p_1, …, p_N $ and output $ p $.
\end{definition}

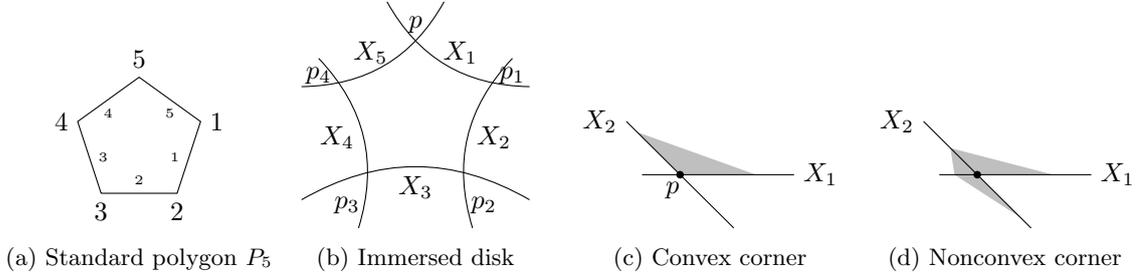
\begin{figure}
\centering
\begin{subfigure}{0.22\linewidth}
\centering
\begin{tikzpicture}
\path[draw] (0, 0) node[below] {3} -- ++(right:1) node[midway, above] {\tiny 2} node[below] {2} -- ++(72:1) node[midway, left] {\tiny 1} node[right] {1} -- ++(144:1) node[midway, below] {\tiny 5} node[above] {5} -- ++(216:1) node[midway, below] {\tiny 4} node[left] {4} -- ++(-72:1) node[midway, right] {\tiny 3};
\end{tikzpicture}
\caption{Standard polygon $ P_5 $}
\end{subfigure}
\begin{subfigure}{0.22\linewidth}
\centering
\begin{tikzpicture}[scale=0.75]
\path[draw] (0, 0) to[bend left] node[midway, below] {$ X_3 $} (4, 0);
\path[draw] (0, 2) to[bend right] node[pos=0.4, above] {$ X_5 $} (2.5, 3.5);
\path[draw] (4, 2) to[bend left] node[pos=0.4, above] {$ X_1 $} (1.5, 3.5);
\path[draw] (0.3, 2.5) to[bend left] node[pos=0.5, left] {$ X_4 $} (1, -0.5);
\path[draw] (3.7, 2.5) to[bend right] node[pos=0.5, right] {$ X_2 $} (3, -0.5);
\path (0.3, 2.2) node {$ p_4 $};
\path (3.7, 2.2) node {$ p_1 $};
\path (0.8, -0.1) node {$ p_3 $};
\path (3.2, -0.1) node {$ p_2 $};
\path (2, 3.1) node {$ p $};
\end{tikzpicture}
\caption{Immersed disk}
\end{subfigure}
\begin{subfigure}{0.25\linewidth}
\centering
\begin{tikzpicture}
\path[fill, lightgray] (0.5, 0) ++(135:0.8) -- (0.5, 0) -- ++(1, 0) -- cycle;
\path[draw] (0, 0) -- ++(right:2) node[right] {$ X_1 $};
\path[draw] (0.5, 0) -- ++(135:1) node[left] {$ X_2 $} (0.5, 0) -- ++(315:1);
\path[fill] (0.5, 0) circle[radius=0.05];
\path (0.4, -0.2) node {$ p $};
\end{tikzpicture}
\caption{Convex corner}
\end{subfigure}
\begin{subfigure}{0.22\linewidth}
\centering
\begin{tikzpicture}
\path[fill, lightgray] (1.5, 0) -- ($ (0.5, 0) + (135:0.5) $) -- (0.2, 0) -- ($ (0.5, 0) + (315:0.8) $) -- (0.5, 0) -- cycle;
\path[draw] (0, 0) -- ++(right:2) node[right] {$ X_1 $};
\path[draw] (0.5, 0) -- ++(135:1) node[left] {$ X_2 $} (0.5, 0) -- ++(315:1);
\path[fill] (0.5, 0) circle[radius=0.05];
\end{tikzpicture}
\caption{Nonconvex corner}
\end{subfigure}
\caption{Immersed disks}
\label{fig:prelim-fukaya-disk}
\end{figure}

By convexity of the corners, we mean that the image of any interior angle of the polygon is strictly smaller than half a full turn. Here, an interior angle of the polygon is interpreted as a very small curve near any corner of $ P_{N+1} $, and a half turn is the natural half turn given by the two sides of the tangent line of $ X_i $ at $ p_i ∈ S $. All this is depicted in \autoref{fig:prelim-fukaya-disk}.

\begin{remark}
Regarding the numbering of the disk inputs, we deviate from the Fukaya-theoretic literature. More precisely, the standard convention \cite{Abouzaid} is to number the disk inputs in counterclockwise order. Instead, we number the disk inputs in clockwise order. The difference is necessary in order to match with the convention for gentle algebras \cite{Bocklandt}.
\end{remark}

The orientation of a curve $ X_i $ involved in a smooth immersed disk $ D $ need not agree with the orientation of $ ∂P_{N+1} $. We can give the boundary $ ∂P_{N+1} $ the orientation pointing in clockwise direction and distinguish whether $ X_i $ agrees with this orientation or not:

\begin{definition}
Let $ D ∈ M(p_1, …, p_N, p) $ be a smooth immersed disk, with $ p_i ∈ X_i ∩ X_{i+1} $. Then:
\begin{itemize}
\item $ X_i $ is \emph{oriented clockwise with $ D $} if the orientation of $ X_i $ agrees with clockwise orientation of $ ∂P_{N+1} $,
\item $ X_i $ is \emph{oriented counterclockwise with $ D $} if the orientation of $ X_i $ is opposite to the clockwise orientation of $ ∂P_{N+1} $.
\end{itemize}
\end{definition}

The differences in orientation give rise to what we call the Abouzaid sign of the disk. This sign is taken from \cite{Abouzaid} and is the surface world incarnation of a sign rule in higher dimensions.

\begin{definition}
\label{def:prelim-fukaya-pre-Abouzaid}
Let $ D ∈ M(p_1, …, p_N, p) $ be a smooth immersed disk with inputs $ p_i: X_i → X_{i+1} $ and output $ p: X_1 → X_{N+1} $. Then the \emph{Abouzaid sign} $ \Abouzaid(D) ∈ ℤ/2ℤ $ is the number of indices $ i $ such that $ p_i $ is odd and $ X_{i+1} $ is oriented counterclockwise with $ D $, plus one if $ p $ is odd and $ X_{N+1} $ is oriented counterclockwise with $ D $, plus the number of $ \# $ signs from the spin structure on the boundary of the disk.
\end{definition}

With all devices ready, we can recall the construction of the Fukaya pre-category.

\begin{definition}
Let $ (S, M) $ be a punctured surface, with exact symplectic form $ ω = dθ $. Then the exact \emph{Fukaya pre-category} $ \Fuk^{\pre} (S, M) $ is defined as follows:
\begin{itemize}
\item The objects are the exact unobstructed curves in $ S \setminus M  $ with chosen spin structures.
\item The set $ \cat C^N_{\transversal} $ consists of the tranversal sequences according to \autoref{def:prelim-fukaya-pre-transversal}.
\item For transversal $ X, Y $, the hom space $ \Hom(X, Y) $ is freely spanned over $ ℂ $ by the intersection points $ p ∈ X ∩ Y $:
\begin{equation*}
\Hom(X, Y) = \bigoplus_{p ∈ X ∩ Y} ℂp.
\end{equation*}
\item For any sequence $ p_1, …, p_N $ of intersection points $ p_i: X_i → X_{i+1} $ and $ p: X_1 → X_{N+1} $, the higher product is defined as
\begin{equation*}
μ^N (p_N, …, p_1) = \sum_{p ∈ X_1 ∩ X_{N+1}} \sum_{D ∈ M(p_1, …, p_N, p)} (-1)^{\Abouzaid(D)} p.
\end{equation*}
\end{itemize}
\end{definition}

\begin{figure}
\centering
\begin{tikzpicture}
\path[draw] (0, 0) -- ++(4, 0) -- ++(up:1) -- ++(left:4) -- ++(down:1);
\path[fill, lightgray] (0.3, 0.3) -- (1.3, 0.3) -- (1.3, 0.7) -- (0.3, 0.7) -- cycle;
\path[fill] (0.3, 0.3) circle[radius=0.05];
\path[fill] (1.3, 0.3) circle[radius=0.05];
\path[fill] (1.3, 0.7) circle[radius=0.05];
\path[fill] (0.3, 0.7) circle[radius=0.05];
\foreach \i in {0, 1, 2, 3, 4} \path[draw] (\i, 0) -- (\i, 1);
\path[draw] (0.15, 0.3) -- (3.75, 0.3);
\path[draw] (0.15, 0.7) -- (3.75, 0.7);
\foreach \i in {0.3, 1.3, 2.3, 3.3} \path[draw] (\i, 0.15) -- (\i, 0.85);
\begin{scope}[shift={(4.5, 0)}]
\path[draw] (0, 0) -- ++(4, 0) -- ++(up:1) -- ++(left:4) -- ++(down:1);
\path[fill, lightgray] (0.3, 0.3) -- (2.3, 0.3) -- (2.3, 0.7) -- (0.3, 0.7) -- cycle;
\path[fill] (0.3, 0.3) circle[radius=0.05];
\path[fill] (2.3, 0.3) circle[radius=0.05];
\path[fill] (2.3, 0.7) circle[radius=0.05];
\path[fill] (0.3, 0.7) circle[radius=0.05];
\foreach \i in {0, 1, 2, 3, 4} \path[draw] (\i, 0) -- (\i, 1);
\path[draw] (0.15, 0.3) -- (3.75, 0.3);
\path[draw] (0.15, 0.7) -- (3.75, 0.7);
\foreach \i in {0.3, 1.3, 2.3, 3.3} \path[draw] (\i, 0.15) -- (\i, 0.85);
\end{scope}
\begin{scope}[shift={(9, 0)}]
\path[draw] (0, 0) -- ++(4, 0) -- ++(up:1) -- ++(left:4) -- ++(down:1);
\path[fill, lightgray] (0.3, 0.3) -- (3.3, 0.3) -- (3.3, 0.7) -- (0.3, 0.7) -- cycle;
\path[fill] (0.3, 0.3) circle[radius=0.05];
\path[fill] (3.3, 0.3) circle[radius=0.05];
\path[fill] (3.3, 0.7) circle[radius=0.05];
\path[fill] (0.3, 0.7) circle[radius=0.05];
\foreach \i in {0, 1, 2, 3, 4} \path[draw] (\i, 0) -- (\i, 1);
\path[draw] (0.15, 0.3) -- (3.75, 0.3);
\path[draw] (0.15, 0.7) -- (3.75, 0.7);
\foreach \i in {0.3, 1.3, 2.3, 3.3} \path[draw] (\i, 0.15) -- (\i, 0.85);
\end{scope}
\end{tikzpicture}
\caption{Non-exact curves on a torus have endless disks}
\end{figure}
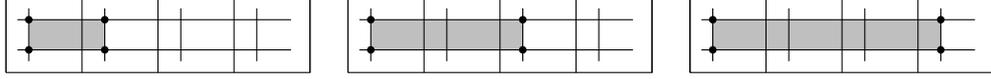

\begin{theorem}[{\cite{Abouzaid}}]
$ \preFuk (S, M) $ is an $ A_∞ $-pre-category.
\end{theorem}

\begin{remark}
Spin structures determine the signs in the higher products of the Fukaya category. A spin structure on a Lagrangian $ X $ can however also be seen as a special case of local system on $ X $: One bakes the spin structure into the local system on $ X $. Upon passing to a version of Fukaya category where each Lagrangian comes with a local system of any dimension assigned, the Fukaya category roughly becomes closed under taking cones. In fact, taking a cone amounts to adding up the local systems.
\end{remark}

\begin{remark}
There are two reasons we only include exact Lagrangians in the Fukaya pre-category. First, we do not need non-exact curves in this paper at all, since the zigzag curves are already exact. The second reason is due to the Novikov field. Indeed, including non-exact curves allows the set $ M(p_1, …, p_N, p) $ of immersed disks to be infinite which requires the technical insertion of the Novikov field. The higher product coming from an immersed disk $ D $ then gets multiplied by the formal power $ t^{ω(D)} $, where $ ω(D) $ denotes the symplectic area of $ D $. This renders the Fukaya pre-category an $ A_∞ $-pre-category over the Novikov field. Since the aim of this paper is to compare the relative Fukaya category to the gentle algebra defined over $ ℂ $, we have decided to avoid the Novikov field early on. It is an interesting question what a discrete model would look like for non-exact parts of the relative Fukaya category that can only be defined with the Novikov field. The discrete model would then need to be defined over the Novikov field instead of $ ℂ $ and its higher products would need to be defined upon a notion of (algebraic) symplectic area.
\end{remark}

\subsection{The exact Fukaya category}
\label{sec:prelim-fukaya-cat}
In this section, we recall the notion of Fukaya category. We explain how one passes from the pre-category of the previous section to an actual category. In particular, we intend to make the reader acquainted with the endomorphism spaces in the Fukaya category and how to extend the $ A_∞ $-products from transversal sequences to all sequences of morphisms.

The characteristic property of the Fukaya category is that its transversal part is precisely the Fukaya pre-category. By transversal parts we mean the following:

\begin{definition}
Let $ \cat C $ be an $ A_∞ $-pre-category and $ \cat D $ an $ A_∞ $-category. Assume $ \Ob \cat C = \Ob \cat D $. Then the \emph{transversal part} of $ \cat D $ (with respect to $ \cat C $) is the $ A_∞ $-pre-category $ (\cat D)_{\transversal} $ defined by
\begin{align*}
\big((\cat D)_{\transversal}\big)^N_{\transversal} &= (\cat C)^N_{\transversal}, \\
μ^N_{(\cat D)_{\transversal}} &= μ^N_{\cat D} \restr{(\cat C)^{N+1}_{\transversal}}.
\end{align*}
The category $ \cat D $ \emph{agrees on the transversal part} with $ \cat C $ if $ \cat C = (\cat D)_{\transversal} $.
\end{definition}

We provide an ad-hoc definition of the Fukaya category as follows:

\begin{definition}
\label{def:prelim-fukaya-cat-def}
Let $ (S, M) $ be a punctured surface. Then the exact \emph{Fukaya category} $ \Fuk (S, M) $ is any $ A_∞ $-category which agrees with $ \preFuk (S, M) $ on the transversal part.
\end{definition}

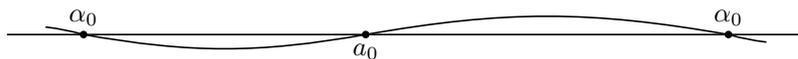
\begin{figure}
\centering
\begin{tikzpicture}
\path[->] (0, 0) -- ++(315:1.5) coordinate[midway] (stop-1) -- ++(45:1.5) coordinate[midway] (stop-2) -- ++(315:1.5) -- ++(45:1.5) -- ++(315:1.5) coordinate[midway] (stop-5) -- ++(45:1.5) -- ++(315:1.5) -- ++(45:1.5) -- ++(315:1.5) coordinate[midway] (stop-9) -- ++(45:1.5) coordinate[midway] (stop-10);
\path[draw, semithick] ($ (stop-1)!-0.05!(stop-10) $) -- ($ (stop-1)!1.05!(stop-10) $);
\path ($ (stop-9)!0.5!(stop-10) $) coordinate (coid-right);
\path ($ (stop-1)!0.5!(stop-2) $) coordinate (coid-left);
\path[draw, semithick] ($ (coid-right) + (0.5, -0.1) $) to[bend left=4] (coid-right) to[bend right=10] (stop-5) to[bend left=10] (coid-left) to[bend right=4] ++(-0.5, 0.1);
\path[fill] (coid-left) circle[radius=0.05] node[above] {$ α_0 $};
\path[fill] (stop-5) circle[radius=0.05] node[below] {$ a_0 $};
\path[fill] (coid-right) circle[radius=0.05] node[above] {$ α_0 $};
\end{tikzpicture}
\caption{A Lagrangian and its Hamiltonian deformation}
\label{fig:prelim-fukaya-cat-Ham}
\end{figure}

Explicit construction of the Fukaya category exist. The standard reference is Seidel's work \cite{Seidel}. The idea is to apply Hamiltonian deformations to make nontransversal pairs of Lagrangians transversal. Most importantly, hom spaces are then also defined for non-transversal pairs. For example, the endomorphism space $ \End(L) $ of a Lagrangian $ L $ contains an \emph{identity} and a \emph{co-identity} element, which we may in the Fukaya category context denote $ \id $ and $ \id^* $. The philosophy is that Hamiltonian deformation of $ L $ yields a transversal version of $ L $, intersecting precisely twice with $ L $. This is depicted in \autoref{fig:prelim-fukaya-cat-Ham}.

\begin{remark}
There exist approaches of constructing the Fukaya category from the Fukaya pre-category by purely categorical methods. In \cite{Efimov}, Efimov proved the conjecture attributed to Kontsevich-Soibelman that every $ A_∞ $-pre-category is quasi-equivalent to an $ A_∞ $-category as $ A_∞ $-pre-categories. The quasi-equivalence relation for $ A_∞ $-pre-categories is defined in \cite[Definition 2.18/2.19]{Efimov}. For ordinary $ A_∞ $-categories, this notion coincides with the ordinary notion of quasi-equivalence.

Efimov's theorem implies there is an $ A_∞ $-category quasi-equivalent to $ \preFuk (S, M) $. This $ A_∞ $-category in almost but not quite (a model for) the Fukaya pre-category in the sense of \autoref{def:prelim-fukaya-cat-def}, since it may have larger hom spaces than the Fukaya pre-category even on the transversal sequences.
\end{remark}

The higher products of non-transversal sequences become very difficult to grasp, since multiple Hamiltonian deformations may need to be performed on the same Lagrangian in order to make all intersections transversal. This results in ambiguities, resolved by providing additional deformation data. In summary, the higher products cannot be determined by simply staring at them. In contrast, our paper provides an explicit realization also of these higher products on non-transversal sequence, at least on zigzag paths. For more on exact Fukaya categories, we refer to \cite{Bocklandt-book} and \cite{Seidel}.

\subsection{The relative exact Fukaya pre-category}
\label{sec:prelim-fukaya-relpre}
In this section, we recall the relative exact Fukaya pre-category for punctured surfaces. The starting point of the relative exact Fukaya pre-category is the exact Fukaya pre-category. The idea is to deform the products by allowing the disk to cover punctures. The resulting object is what we will call an $ A_∞ $-pre-category deformation.

The history of the subject can be traced back fairly accurately: In \cite{Seidel-relative},  Seidel introduced the idea of deforming the Fukaya category by working relative to a divisor. Twenty years later, the relative Fukaya category was finally constructed in \cite{Sheridan-Perutz}. Its versality as a deformation of the ordinary Fukaya category was investigated in \cite{Sheridan}. Lekili, Perutz and Polishchuk \cite{Lekili-Perutz, Lekili-Polishchuk} proved deformed mirror symmetry for the $ n $-punctured torus, apparently the first use of the relative Fukaya category in mirror symmetry.

\begin{definition}
Let $ \cat C $ be an $ A_∞ $-pre-category and $ B $ a deformation base, e.g.~$ B = ℂ⟦q⟧ $. Then an \emph{$ A_∞ $-pre-category deformation} $ \cat C_q $ of $ \cat C $ is an (infinitesimally curved) and $ B $-linear $ A_∞ $-pre-category structure on $ B \htensor \cat C $. More precisely, this means that $ \cat C_q $ has:
\begin{itemize}
\item the same objects as $ \cat C $,
\item hom spaces $ \Hom_{\cat C_q} (X, Y) = B \htensor \Hom_{\cat C} (X, Y) $ for $ X, Y ∈ \cat C^2_{\transversal} $,
\item $ B $-multilinear products of degree $ 2 - N $
\begin{equation*}
μ_q^N: \Hom_{\cat C_q} (X_N, X_{N+1}) \tensor … \tensor \Hom_{\cat C_q} (X_1, X_2) → \Hom_{\cat C_q} (X_1, X_{N+1}), ~ N ≥ 1
\end{equation*}
for all transversal sequences $ X_1, …, X_{N+1} $,
\item curvature of degree $ 2 $ for every object $ X ∈ \cat C $ with $ (X, X) ∈ \cat C^2_{\transversal} $:
\begin{equation*}
μ_{q, X}^0 ∈ \mathfrak{m} \Hom_{\cat C_q} (X, X),
\end{equation*}
\end{itemize}
such that $ μ_q $ reduces to $ μ $ once the maximal ideal $ \mathfrak{m} ⊂ B $ is divided out, and $ μ_q $ satisfies the curved $ A_∞ $ relations on transversal sequences.
\end{definition}

With this definition in mind, we would like to define the relative version of the Fukaya pre-category. The idea is to define the higher products by counting smooth immersed disks, which are now also allowed to cover punctures. By abuse of terminology, we shall keep using the terminology of “smooth immersed disks” even for those smooth immersed disks which cover punctures:

\begin{definition}
\label{def:prelim-fukaya-relpre-disk}
Let $ X_1, …, X_{N+1} $ be a transversal sequence of $ N+1 ≥ 2 $ unobstructed curves. Let $ p_1, …, p_N $ be a sequence of intersection points $ p_i: X_i → X_{i+1} $ and let $ p ∈ X_1 → X_{N+1} $ be another intersection point. A \emph{smooth immersed disk} with \emph{inputs} $ p_1, …, p_N $ and \emph{output} $ p $ consists of an orientation-preserving polygon immersion $ D: P_{N+1} → S $ up to reparametrization, such that
\begin{itemize}
\item the corners of $ D $ are all convex,
\item the $ i $-th edge of $ P_{N+1} $ lands on $ X_i $ for $ 1 ≤ i ≤ N+1 $,
\item the $ i $-th corner of $ P_{N+1} $ lands on $ p_i $.
\end{itemize}
We denote by $ M_q (p_1, …, p_N, p) $ the set of smooth immersed disks with inputs $ p_1, …, p_N $ and output $ p $.
\end{definition}

The deformation base of the relative exact Fukaya pre-category is $ B = ℂ⟦M⟧ $. This is the power series ring with one variable for each puncture $ m ∈ M $. Correspondingly, every puncture $ q ∈ M $ given an element $ q ∈ ℂ⟦M⟧ $. Multiple punctures $ q_1, …, q_s ∈ M $ can be multiplied to form the element $ q_1 … q_s ∈ ℂ⟦M⟧ $.

\begin{definition}
\label{def:prelim-fukaya-relpre-Abouzaid}
Let $ D ∈ M_q (p_1, …, p_N, p) $. Then the \emph{Abouzaid sign} $ \Abouzaid(D) $ is defined precisely as in \autoref{def:prelim-fukaya-pre-Abouzaid}. The \emph{$ q $-parameter} $ \punctures(D) ∈ ℂ⟦M⟧ $ is defined as the product of all the punctures reached by the interior of $ P_{N+1} $ under $ D $, counting multiplicities.
\end{definition}

The parameter $ \punctures(D) $ is very similar to the deformation parameter in the higher products of $ \Gtl_q Q $ in \autoref{sec:prelim-gtlq}.

\begin{definition}
\label{def:prelim-fukaya-relpre-def}
The \emph{relative exact Fukaya pre-category} $ \relpreFuk (S, M) $ is the $ A_∞ $-pre-category deformation of $ \preFuk (S, M) $ over $ B = ℂ⟦M⟧ $ given by the deformed $ A_∞ $-products on transversal sequences
\begin{equation*}
μ^N_q (p_N, …, p_1) = \sum_{p: X_1 → X_{N+1}} \sum_{D ∈ M_q (p_1, …, p_N, p)} (-1)^{\Abouzaid(D)} \punctures(D) p.
\end{equation*}
\end{definition}

Checking that $ \relpreFuk (S, M) $ is really an $ A_∞ $-pre-category deformation involves two parts. The first part consists of checking that the higher products $ μ^N_q $ are well-defined. Indeed, exactness guarantees that the set of disks not covering any puncture is finite, but the case of disks covering some punctures is non-trivial. A likely successful procedure is as follows: For each monomial $ q = q_1 … q_s ∈ ℂ⟦M⟧ $, use exactness of the curves to bound the maximum size of disks that cover precisely the punctures $ q_1, …, q_s $. A standard argument is then that the Gromov compactness theorem implies the number of disks is finite. See e.g.~ \cite[Section 6.2.3]{Bocklandt}.

The second check consists of evaluating the $ A_∞ $-relations on the transversal sequences. This boils down to re-doing the work of Abouzaid \cite{Abouzaid}, now allowing the disks to cover punctures. The procedure should be straightforward and conclude that the $ A_∞ $-relations still hold. In total, this renders $ \relpreFuk (S, M) $ an $ A_∞ $-pre-category deformation of $ \preFuk (S, M) $.

\subsection{The relative exact Fukaya category}
\label{sec:prelim-fukaya-rel}
In this section, we recall the notion of relative exact Fukaya category. The starting point is the relative exact Fukaya pre-category which we sketched in \autoref{sec:prelim-fukaya-relpre}. It is an $ A_∞ $-pre-category deformation of $ \preFuk (S, M) $. In \autoref{sec:prelim-fukaya-cat}, we have seen that the category $ \preFuk (S, M) $ has a lift to an actual $ A_∞ $-category $ \Fuk (S, M) $. The question we discuss in this section is what a lift of the deformation $ \relpreFuk (S, M) $ to a deformation $ \relFuk (S, M) $ of $ \Fuk (S, M) $ should look like. Our desired lifting procedure is best captured graphically as follows:

\begin{center}
\begin{tikzpicture}
\path (4, 2) node (B) {$ \Fuk $} (0, 0) node (C) {$ \relpreFuk $} (4, 0) node (D) {$ \preFuk $};
\path[draw] (C.east) to node[midway, above] {defo} (D.west);
\path[draw, {Hooks[right, scale=2]}->] (D.north) to (B.south);
\path (5.5, 1) node {\LARGE $ \rightsquigarrow $};
\begin{scope}[shift={(7, 0)}]
\path (0, 2) node (A) {$ \relFuk $} (4, 2) node (B) {$ \Fuk $} (0, 0) node (C) {$ \relpreFuk $} (4, 0) node (D) {$ \preFuk $};
\path[draw, {Hooks[right, scale=2]}->] (C.north) to (A.south);
\path[draw] (C.east) to node[midway, above] {defo} (D.west);
\path[draw, {Hooks[right, scale=2]}->] (D.north) to (B.south);
\path[draw] (A.east) to node[midway, above] {defo} (B.west);
\end{scope}
\end{tikzpicture}
\end{center}

\begin{definition}
\label{def:prelim-fukaya-agree}
Let $ \cat C $ be an $ A_∞ $-pre-category and $ \cat C_q $ an $ A_∞ $-pre-category deformation. Let $ \cat D $ be an $ A_∞ $-category and $ \cat D_q $ an $ A_∞ $-deformation. Assume $ \Ob \cat C = \Ob \cat D $. Then the \emph{transversal part} of $ \cat D_q $ (with respect to $ \cat C_q $) is the $ A_∞ $-pre-category deformation $ (\cat D_q)_{\transversal} $ of $ (\cat D)_{\transversal} $ defined by
\begin{equation*}
μ^N_{(\cat D_q)_{\transversal}} = μ^N_{\cat D_q} \restr{(\cat C)^{N+1}_{\transversal}}.
\end{equation*}
\end{definition}

We may say that the deformation $ \cat D_q $ \emph{agrees on the transversal part} with $ \cat C_q $ if $ \cat C_q = (\cat D_q)_{\transversal} $. For sake of explicitness in \autoref{sec:subdisk-main}, we provide the following terminology:

\begin{definition}
\label{def:prelim-fukaya-defiso}
Let $ \cat C $ and $ \cat D $ be $ A_∞ $-pre-categories and $ \cat C_q $ and $ \cat D_q $ be $ A_∞ $-pre-category deformations. Then a \emph{strict isomorphism} $ F_q: \cat C_q → \cat D_q $ of $ A_∞ $-pre-category deformations consists of
\begin{itemize}
\item a bijection $ F_q: \Ob\cat C → \Ob\cat D $ such that
\begin{equation*}
∀N ≥ 1: \quad (\cat D)^N_{\transversal} = \{(F_q X_1, …, F_q X_N) \running (X_1, …, X_N) ∈ (\cat C)^N_{\transversal}\},
\end{equation*}
\item for every $ X, Y ∈ (\cat C)^2_{\transversal} $ a $ B $-linear isomorphism $ F^1: \Hom_{\cat C_q} (X, Y) → \Hom_{\cat D_q} (F_q X, F_q Y) $ of degree $ 0 $ such that
\begin{align*}
& ∀N ≥ 1, \quad (X_1, …, X_{N+1}) ∈ (\cat C)^{N+1}_{\transversal}, \quad a_i ∈ \Hom_{\cat C} (X_i, X_{i+1}): \\
& \quad F_q^1 (μ_{\cat C_q} (a_N, …, a_1)) = μ_{\cat D_q} (F_q^1 (a_N), …, F_q^1 (a_1)).
\end{align*}
\end{itemize}
\end{definition}

We provide an ad-hoc definition of the relative Fukaya category as follows:

\begin{definition}
\label{def:prelim-fukaya-rel-def}
Let $ (S, M) $ be a punctured surface. Then the \emph{relative Fukaya category} $ \relFuk (S, M) $ is any $ A_∞ $-deformation of $ \Fuk (S, M) $ such that $ \relFuk (S, M)_{\transversal} = \relpreFuk (S, M) $.
\end{definition}

\begin{remark}
Sheridan and Perutz \cite{Sheridan-Perutz} provide explicit constructions of relative Fukaya categories.
\end{remark}

\begin{remark}
All concrete models of the Fukaya category give rise to a priori different notions of relative Fukaya category. It should be possible to show that these are in fact all isomorphic.
\end{remark}

\subsection{Zigzag paths as Lagrangians}
\label{sec:prelim-fukaya-zigzag}
In the present section, we show how to interpret zigzag paths as objects in Fukaya categories. The starting point is a dimer $ Q $. It gives rise to a collection of zigzag paths and we show how to turn them into curves in $ |Q| $ which we call “zigzag curves”. We recall how to make these curves objects of the relative Fukaya pre-category and look at their hom spaces and higher products. The material can also be found in \cite[Chapter 10]{CHL}.

The Fukaya category of a dimer is simply defined as the Fukaya category of its underlying punctured surface. More precisely, let $ Q $ be a dimer. Then $ Q $ includes the datum of a punctured surface $ (|Q|, Q_0) $ and we write
\begin{equation*}
\begin{aligned}
\preFuk Q &= \preFuk (|Q|, Q_0), && \Fuk Q = \Fuk (|Q|, Q_0), \\
\relpreFuk Q &= \relpreFuk (|Q|, Q_0), && \relFuk Q = \Fuk (|Q|, Q_0).
\end{aligned}
\end{equation*}

The first step in this section is to turn zigzag paths into curves. Let $ L $ be a zigzag path in $ Q $. Then $ L $ consists by definition of a path $ … a_1 a_0 a_{-1} … $ of composable arcs in $ Q $, together with left/right indications. The idea to produce a curve $ \smooth L ⊂ |Q| \setminus Q_0 $ from $ L $ is to stitch together the arcs $ a_i $ in sequence, minding the left/right indication. The precise definition reads as follows:

\begin{definition}
\label{def:prelim-fukaya-zigzag-curve}
Let $ Q $ be a dimer and $ L $ a zigzag path in $ Q $. Then the associated \emph{zigzag curve} $ \smooth L: S^1 → S \setminus M $ is defined by connecting the midpoints of the arcs $ …, a_{-1}, a_0, a_1, … $ by means of the \emph{angle cutting procedure}:
\begin{itemize}
\item If $ L $ turns left at index $ i $, then the midpoint of $ a_i $ is connected to the midpoint of $ a_{i+1} $ by turning clockwise around the puncture $ h(a_i) = t(a_{i+1}) $.
\item If $ L $ turns right at index $ i $, then the midpoint of $ a_i $ is connected to the midpoint of $ a_{i+1} $ by turning counterclockwise around the puncture $ h(a_i) = t(a_{i+1}) $.
\end{itemize}
The connecting arc between the midpoints of $ a_i $ and $ a_{i+1} $ is to be chosen close enough to the puncture that it does not intersect with the zigzag curves associated with the other zigzag paths.
\end{definition}

\begin{figure}
\centering
\begin{tikzpicture}
\begin{scope}[scale=0.5]
\path[draw, -{To[scale=2]}] (0, 0) -- ++(315:1.5) 
coordinate[midway] (alpha1-end) -- ++(45:1.5) 
coordinate[midway] (alpha1-start) -- ++(315:1.5) 
coordinate[midway] (alpha3-end) -- ++(45:1.5) 
coordinate[midway] (alpha3-start) -- ++(315:1.5) 
coordinate[midway] (alpha5-end) -- ++(45:1.5) 
coordinate[midway] (alpha5-start) -- ++(315:1.5) 
coordinate[midway] (additional);
\path (4, -2) node {zigzag path};
\end{scope}
\path (6, -0.5) node {\Large $ \longleftrightarrow $};
\begin{scope}[shift={(8, 0)}, scale=0.5]
\path[draw, gray!50, -{To[scale=2]}] (0, 0) -- ++(315:1.5) 
coordinate[midway] (alpha1-end) -- ++(45:1.5) 
coordinate[midway] (alpha1-start) -- ++(315:1.5) 
coordinate[midway] (alpha3-end) -- ++(45:1.5) 
coordinate[midway] (alpha3-start) -- ++(315:1.5) 
coordinate[midway] (alpha5-end) -- ++(45:1.5) 
coordinate[midway] (alpha5-start) -- ++(315:1.5) 
coordinate[midway] (additional);
\path[draw, very thick, rounded corners] ($ (alpha1-end)!-0.5!(alpha1-start) + (0, 0.2) $) -- (alpha1-end) -- ($ (alpha1-start)!0.5!(alpha1-end) + (0, -0.1) $) -- (alpha1-start) -- ($ (alpha1-start)!0.5!(alpha3-end) + (0, 0.1) $) -- (alpha3-end) -- ($ (alpha3-start)!0.5!(alpha3-end) + (0, -0.1) $) -- (alpha3-start) -- ($ (alpha3-start)!0.5!(alpha5-end) + (0, 0.1) $) -- (alpha5-end) -- ($ (alpha5-start)!0.5!(alpha5-end) + (0, -0.1) $) -- (alpha5-start) -- ($ (alpha5-start)!0.5!(additional) + (0, 0.1) $) -- (additional) -- ($ (additional)!-0.5!(alpha5-start) + (0, -0.2) $);
\path (4, -2) node {zigzag curve};
\end{scope}
\end{tikzpicture}
\caption{Zigzag path and zigzag curve}
\label{fig:prelim-fukaya-zigzag-correspondence}
\end{figure}
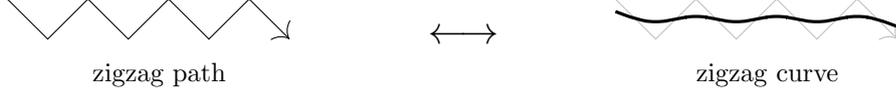

In \autoref{fig:prelim-fukaya-zigzag-correspondence}, we have depicted a zigzag path together with its associated zigzag curve. By definition of $ \smooth L $, we have the freedom to deform $ \smooth L $ a bit near the punctures. We quote the following lemma from \cite{CHL} which claims that a small deformation can be chosen in such a way that that $ \smooth L $ becomes an exact Lagrangian. The condition is that the Liouville form $ α $ near the punctures takes the shape $ α = dθ/r $ where $ r $ is the distance from the puncture and $ θ $ the polar angle.

\begin{lemma}[{\cite[Lemma 10.6]{CHL}}]
Let $ Q $ be a dimer. Pick a symplectic form $ ω = dα $ on $ |Q| \setminus Q_0 $ such that $ α = dθ/r $ near the punctures. Then the curve $ \smooth L: S^1 → |Q| \setminus Q_0 $ can be constructed in such a way that it is exact with respect to $ ω $.
\end{lemma}

Zigzag curves in a general dimer are not contractible and do not bound teardrops in $ |Q| \setminus Q_0 $. If $ Q $ is geometrically consistent, then even in the closed surface $ |Q| $ a zigzag curve is not contractible and does not bound a teardrop. Upon specification of spin structures and a symplectic form, the zigzag curves $ \smooth L $ define object of the exact Fukaya category $ \Fuk Q $. In analogy to the “category of zigzag paths” constructed combinatorially via $ \Gtl_q Q $ in \autoref{sec:deformed}, we may give the subcategory of $ \relFuk Q $ consisting of these objects a name:

\begin{definition}
\label{def:prelim-fukaya-zigzag-restr}
Let $ Q $ be a geometrically consistent dimer. Pick a choice of spin structure for every zigzag curve in $ Q $. Then we denote by
\begin{equation*}
\preFuk Q \restr{\Ob\ZigzagCat}, \quad \Fuk Q \restr{\Ob\ZigzagCat}, \quad \relpreFuk Q \restr{\Ob\ZigzagCat}, \quad \relFuk Q \restr{\Ob\ZigzagCat}
\end{equation*}
the subcategories of the (relative) Fukaya (pre-)categories given by the zigzag curves with chosen spin structure.
\end{definition}

An intersection between two zigzag paths shall be defined as a shared arc of the two zigzag paths. More precisely, we use the following terminology:

\begin{definition}
\label{def:prelim-fukaya-discreteintersection}
Let $ Q $ be a dimer and $ L_1, L_2 $ be two zigzag paths. An \emph{(indexed) intersection} between $ L_1 $ and $ L_2 $ is a pair $ (a_i, b_j) $ such that
\begin{itemize}
\item $ a_i $ is an indexed arc of $ L_1 $,
\item $ b_j $ is an indexed arc of $ L_2 $,
\item $ a_i = b_j $ as arcs in $ Q $,
\item $ L_1 $ turns left at $ a_i $ and $ L_2 $ turns right at $ b_j $, or the other way around.
\end{itemize}
\end{definition}

\begin{remark}
The cautious wording of \autoref{def:prelim-fukaya-discreteintersection} is necessary in order to make a transversal self-intersection count double.
\end{remark}

\begin{figure}
\centering
\begin{subfigure}{0.3\linewidth}
\centering
\begin{tikzpicture}
\path[draw, ->] (0, 0) -- ++(30:2) coordinate[pos=0.5] (target-tail) -- ++(up:1.5) coordinate[pos=0.5] (target-middle) coordinate (5-head);
\path[draw] (5-head) -- ++(30:2) coordinate[pos=0.5] (target-head) node[near end, below] {$ \smooth L_2 $};
\path[draw, ->] (3.1, 0) -- ++(150:2) coordinate[pos=0.5] (source-tail) -- ++(up:1.5) coordinate[pos=0.5] (source-middle) coordinate (2-head);
\path[draw] (2-head) -- ++(150:2) node[near end, below] {$ \smooth L_1 $} coordinate[pos=0.5] (source-head);
\path (source-middle) -- (target-middle) coordinate[pos=0.5] (middle);
\path[draw, thick, rounded corners, ->] ($ (middle)!2!(target-tail) $) -- (middle) coordinate[pos=0.7] (target2)-- ($ (middle)!2!(target-head) $) coordinate[pos=0.2] (target);
\path[draw, thick, rounded corners, ->] ($ (middle)!2!(source-tail) $) -- (middle) coordinate[pos=0.7] (source) -- ($ (middle)!2!(source-head) $) coordinate[pos=0.2] (source2);
\path[fill] (middle) circle(0.1) node[left, shift={(left:0.1)}] {$ p $};
\path[draw, ->, bend right=60] (source) to (target);
\path[draw, ->, bend right=60] (source2) to (target2);
\end{tikzpicture}
\caption{$ p: \smooth L_1 → \smooth L_2 $ is odd}
\end{subfigure}
\begin{subfigure}{0.3\linewidth}
\centering
\begin{tikzpicture}
\path[draw, ->] (0, 0) -- ++(30:2) coordinate[pos=0.5] (target-tail) -- ++(up:1.5) coordinate[pos=0.5] (target-middle) coordinate (5-head);
\path[draw] (5-head) -- ++(30:2) coordinate[pos=0.5] (target-head) node[near end, below] {$ \smooth L_1 $};
\path[draw, ->] (3.1, 0) -- ++(150:2) coordinate[pos=0.5] (source-tail) -- ++(up:1.5) coordinate[pos=0.5] (source-middle) coordinate (2-head);
\path[draw] (2-head) -- ++(150:2) node[near end, below] {$ \smooth L_2 $} coordinate[pos=0.5] (source-head);
\path (source-middle) -- (target-middle) coordinate[pos=0.5] (middle);
\path[draw, thick, rounded corners, ->] ($ (middle)!2!(target-tail) $) -- (middle) coordinate[pos=0.3] (target2) -- ($ (middle)!2!(target-head) $) coordinate[pos=0.7] (target);
\path[draw, thick, rounded corners, ->] ($ (middle)!2!(source-tail) $) -- (middle) coordinate[pos=0.3] (source2) -- ($ (middle)!2!(source-head) $) coordinate[pos=0.7] (source);
\path[fill] (middle) circle(0.1) node[left, shift={(left:0.1)}] {$ p $};
\path[draw, ->, bend right=40] (target) to (source);
\path[draw, ->, bend right=40] (target2) to (source2);
\end{tikzpicture}
\caption{$ p: \smooth L_2 → \smooth L_1 $ is even}
\end{subfigure}
\caption{Shared arcs between $ L_1 $ and $ L_2 $ correspond to intersections $ p ∈ \smooth L_1 ∩ \smooth L_2 $}
\label{fig:prelim-fukaya-zigzag-intersection}
\end{figure}

In \autoref{th:prelim-fukaya-zigzag-Hom}, we explain that intersections between two distinct zigzag curves $ \smooth L_1 $ and $ \smooth L_2 $ are precisely the same as (indexed) intersections between $ L_1 $ and $ L_2 $. This is depicted in \autoref{fig:prelim-fukaya-zigzag-intersection}. In particular, the hom spaces of $ \preFuk Q \restr{\Ob\ZigzagCat} $ can be identified with spans of (indexed) intersections of zigzag paths. To make this work also in case $ L_1 = L_2 $, we have to use a model for the Fukaya category as discussed in \autoref{sec:prelim-fukaya-cat}, in which also the endomorphism spaces are spanned by transversal intersections, plus an identity and a co-identity morphism.

\begin{lemma}
\label{th:prelim-fukaya-zigzag-Hom}
Let $ Q $ be a dimer and $ L_1, L_2 $ be two zigzag paths. Then transversal intersections points of $ \smooth L_1 $ and $ \smooth L_2 $, counting transversal self-intersections double if $ L_1 = L_2 $, are in one-to-one correspondence with indexed intersections between $ L_1 $ and $ L_2 $. Therefore:
\begin{equation}
\label{eq:prelim-fukaya-zigzag-Hom}
\Hom_{\Fuk Q} (\smooth L_1, \smooth L_2) = \vspan_ℂ \{\text{intersections } (a, b) \text{ of } L_1, L_2\} \quad [⊕ ℂ\id ⊕ ℂ\coid \text{ if } L_1 = L_2].
\end{equation}
\end{lemma}

\begin{proof}
The one-to-one correspondence is a simple inspection. It is worth noting that a transversal self-intersection gives rise to two intersection points between $ \smooth L_1 $ and $ \smooth L_2 $, according to the double counting, and two indexed intersections between $ L_1 $ and $ L_2 $.

The description of the hom space \eqref{eq:prelim-fukaya-zigzag-Hom} for $ L_1 ≠ L_2 $ follows from the definition of $ \preFuk Q $ and the requirement that $ (\Fuk Q)_{\transversal} = \preFuk Q $. For $ L_1 = L_2 $, \eqref{eq:prelim-fukaya-zigzag-Hom} follows from our choice for $ \Fuk Q $, which describes endomorphism spaces as spans of transversal intersection points plus identity and co-identity. This finishes the proof.
\end{proof}

The zigzag curves automatically become objects in the relative Fukaya pre-category and the relative Fukaya category, by virtue of \autoref{def:prelim-fukaya-cat-def} and \ref{def:prelim-fukaya-rel-def}. If we choose $ \Fuk Q $ such that the endomorphism spaces are spanned by transversal intersections plus identity and co-identity, then we know the hom spaces of $ \relFuk Q $:

\begin{lemma}
\label{th:prelim-fukaya-zigzag-relHom}
Let $ Q $ be a dimer. Then $ \relFuk Q $ is an $ A_∞ $-deformation of $ \Fuk Q $. Let $ L_1, L_2 $ be two zigzag paths. Then
\begin{equation*}
\Hom_{\relFuk Q} (\smooth L_1, \smooth L_2) = B \htensor \vspan \{\text{intersections } (a, b) \text{ of } L_1, L_2\} \quad [⊕ B\id ⊕ B\coid \text{ if } L_1 = L_2].
\end{equation*}
\end{lemma}

\begin{proof}
By virtue of \autoref{def:prelim-fukaya-rel-def}, $ \relFuk Q $ is an $ A_∞ $-deformation of $ \Fuk Q $. In particular, its hom spaces are given by
\begin{equation*}
\Hom_{\relFuk Q} (\smooth L_1, \smooth L_2) = B \htensor \Hom_{\Fuk Q} (\smooth L_1, \smooth L_2).
\end{equation*}
Using the combinatorical description of $ \Hom_{\Fuk Q} (L_1, L_2) $ from \autoref{th:prelim-fukaya-zigzag-Hom} finishes the proof.
\end{proof}

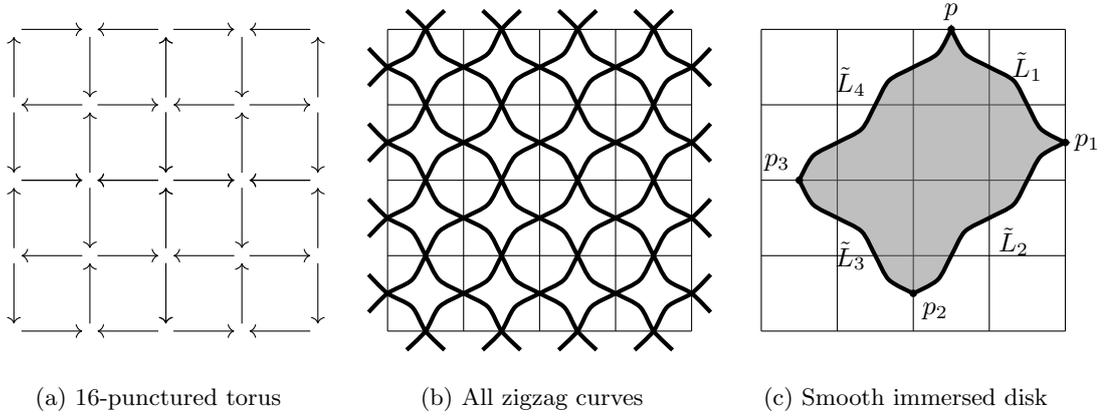
\begin{figure}
\centering
\begin{subfigure}{0.3\linewidth}
\centering
\begin{tikzpicture}
\path[use as bounding box] (-0.5, -0.5) -- (4.5, 4.5);
\newcommand{\arrowbetween}[2]{($ (#1)!0.1!(#2) $) -- ($ (#2)!0.1!(#1) $)}
\newcommand{\thetorus}{%
\path[draw, ->] \arrowbetween{0, 0}{1, 0};
\path[draw, ->] \arrowbetween{2, 0}{1, 0};
\path[draw, ->] \arrowbetween{0, 2}{1, 2};
\path[draw, ->] \arrowbetween{2, 2}{1, 2};
\path[draw, ->] \arrowbetween{1, 1}{0, 1};
\path[draw, ->] \arrowbetween{1, 1}{2, 1};
\path[draw, ->] \arrowbetween{0, 1}{0, 0};
\path[draw, ->] \arrowbetween{0, 1}{0, 2};
\path[draw, ->] \arrowbetween{2, 1}{2, 0};
\path[draw, ->] \arrowbetween{2, 1}{2, 2};
\path[draw, ->] \arrowbetween{1, 0}{1, 1};
\path[draw, ->] \arrowbetween{1, 2}{1, 1};}
\begin{scope}[shift={(2, 0)}] \thetorus \end{scope}
\begin{scope}[shift={(2, 2)}] \thetorus \end{scope}
\begin{scope}[shift={(0, 0)}] \thetorus \end{scope}
\begin{scope}[shift={(0, 2)}] \thetorus \end{scope}
\end{tikzpicture}
\caption{16-punctured torus}
\label{fig:prelim-fukaya-zigzag-products1}
\end{subfigure}
\begin{subfigure}{0.3\linewidth}
\centering
\begin{tikzpicture}
\path[use as bounding box] (-0.5, -0.5) -- (4.5, 4.5);
\path[draw] (0, 0) -- ++(right:1) coordinate[midway] (1-0) -- ++(right:1) coordinate[midway] (3-0) -- ++(right:1) coordinate[midway] (5-0) -- ++(right:1) coordinate[midway] (7-0);
\path[draw] (0, 1) -- ++(right:1) coordinate[midway] (1-2) -- ++(right:1) coordinate[midway] (3-2) -- ++(right:1) coordinate[midway] (5-2) -- ++(right:1) coordinate[midway] (7-2);
\path[draw] (0, 2) -- ++(right:1) coordinate[midway] (1-4) -- ++(right:1) coordinate[midway] (3-4) -- ++(right:1) coordinate[midway] (5-4) -- ++(right:1) coordinate[midway] (7-4);
\path[draw] (0, 3) -- ++(right:1) coordinate[midway] (1-6) -- ++(right:1) coordinate[midway] (3-6) -- ++(right:1) coordinate[midway] (5-6) -- ++(right:1) coordinate[midway] (7-6);
\path[draw] (0, 4) -- ++(right:1) coordinate[midway] (1-8) -- ++(right:1) coordinate[midway] (3-8) -- ++(right:1) coordinate[midway] (5-8) -- ++(right:1) coordinate[midway] (7-8);
\path[draw] (0, 0) -- ++(up:1) coordinate[midway] (0-1) -- ++(up:1) coordinate[midway] (0-3) -- ++(up:1) coordinate[midway] (0-5) -- ++(up:1) coordinate[midway] (0-7);
\path[draw] (1, 0) -- ++(up:1) coordinate[midway] (2-1) -- ++(up:1) coordinate[midway] (2-3) -- ++(up:1) coordinate[midway] (2-5) -- ++(up:1) coordinate[midway] (2-7);
\path[draw] (2, 0) -- ++(up:1) coordinate[midway] (4-1) -- ++(up:1) coordinate[midway] (4-3) -- ++(up:1) coordinate[midway] (4-5) -- ++(up:1) coordinate[midway] (4-7);
\path[draw] (3, 0) -- ++(up:1) coordinate[midway] (6-1) -- ++(up:1) coordinate[midway] (6-3) -- ++(up:1) coordinate[midway] (6-5) -- ++(up:1) coordinate[midway] (6-7);
\path[draw] (4, 0) -- ++(up:1) coordinate[midway] (8-1) -- ++(up:1) coordinate[midway] (8-3) -- ++(up:1) coordinate[midway] (8-5) -- ++(up:1) coordinate[midway] (8-7);
\path[draw, ultra thick] plot[smooth] coordinates{($ (0-1)!-0.5!(1-2) $) (0-1) ($ (0-1)!0.5!(1-2) + (315:0.1) $) (1-2) ($ (1-2)!0.5!(2-3) + (135:0.1) $) (2-3) ($ (2-3)!0.5!(3-4) + (315:0.1) $) (3-4) ($ (3-4)!0.5!(4-5) + (135:0.1) $) (4-5) ($ (4-5)!0.5!(5-6) + (315:0.1) $) (5-6) ($ (5-6)!0.5!(6-7) + (135:0.1) $) (6-7) ($ (6-7)!0.5!(7-8) + (315:0.1) $) (7-8) ($ (7-8)!-0.5!(6-7) $)};
\path[draw, ultra thick] plot[smooth] coordinates{($ (0-3)!-0.5!(1-4) $) (0-3) ($ (0-3)!0.5!(1-4) + (315:0.1) $) (1-4) ($ (1-4)!0.5!(2-5) + (135:0.1) $) (2-5) ($ (2-5)!0.5!(3-6) + (315:0.1) $) (3-6) ($ (3-6)!0.5!(4-7) + (135:0.1) $) (4-7) ($ (4-7)!0.5!(5-8) + (315:0.1) $) (5-8) ($ (5-8)!-0.5!(4-7) $)};
\path[draw, ultra thick] plot[smooth] coordinates{($ (0-5)!-0.5!(1-6) $) (0-5) ($ (0-5)!0.5!(1-6) + (315:0.1) $) (1-6) ($ (1-6)!0.5!(2-7) + (135:0.1) $) (2-7) ($ (2-7)!0.5!(3-8) + (315:0.1) $) (3-8) ($ (3-8)!-0.5!(2-7) $)};
\path[draw, ultra thick] plot[smooth] coordinates{($ (0-7)!-0.5!(1-8) $) (0-7) ($ (0-7)!0.5!(1-8) + (315:0.1) $) (1-8) ($ (1-8)!-0.5!(0-7) $)};
\path[draw, ultra thick] plot[smooth] coordinates{($ (1-0)!-0.5!(2-1) $) (1-0) ($ (1-0)!0.5!(2-1) + (135:0.1) $) (2-1) ($ (2-1)!0.5!(3-2) + (315:0.1) $) (3-2) ($ (3-2)!0.5!(4-3) + (135:0.1) $) (4-3) ($ (4-3)!0.5!(5-4) + (315:0.1) $) (5-4) ($ (5-4)!0.5!(6-5) + (135:0.1) $) (6-5) ($ (6-5)!0.5!(7-6) + (315:0.1) $) (7-6) ($ (7-6)!0.5!(8-7) + (135:0.1) $) (8-7) ($ (8-7)!-0.5!(7-6) $)};
\path[draw, ultra thick] plot[smooth] coordinates{($ (3-0)!-0.5!(4-1) $) (3-0) ($ (3-0)!0.5!(4-1) + (135:0.1) $) (4-1) ($ (4-1)!0.5!(5-2) + (315:0.1) $) (5-2) ($ (5-2)!0.5!(6-3) + (135:0.1) $) (6-3) ($ (6-3)!0.5!(7-4) + (315:0.1) $) (7-4) ($ (7-4)!0.5!(8-5) + (135:0.1) $) (8-5) ($ (8-5)!-0.5!(7-4) $)};
\path[draw, ultra thick] plot[smooth] coordinates{($ (5-0)!-0.5!(6-1) $) (5-0) ($ (5-0)!0.5!(6-1) + (135:0.1) $) (6-1) ($ (6-1)!0.5!(7-2) + (315:0.1) $) (7-2) ($ (7-2)!0.5!(8-3) + (135:0.1) $) (8-3) ($ (8-3)!-0.5!(7-2) $)};
\path[draw, ultra thick] plot[smooth] coordinates{($ (7-0)!-0.5!(8-1) $) (7-0) ($ (7-0)!0.5!(8-1) + (135:0.1) $) (8-1) ($ (8-1)!-0.5!(7-0) $)};
\path[draw, ultra thick] plot[smooth] coordinates{($ (0-1)!-0.5!(1-0) $) (0-1) ($ (0-1)!0.5!(1-0) + (45:0.1) $) (1-0) ($ (1-0)!-0.5!(0-1) $)};
\path[draw, ultra thick] plot[smooth] coordinates{($ (0-3)!-0.5!(1-2) $) (0-3) ($ (0-3)!0.5!(1-2) + (45:0.1) $) (1-2) ($ (1-2)!0.5!(2-1) + (225:0.1) $) (2-1) ($ (2-1)!0.5!(3-0) + (45:0.1) $) (3-0) ($ (3-0)!-0.5!(2-1) $)};
\path[draw, ultra thick] plot[smooth] coordinates{($ (0-5)!-0.5!(1-4) $) (0-5) ($ (0-5)!0.5!(1-4) + (45:0.1) $) (1-4) ($ (1-4)!0.5!(2-3) + (225:0.1) $) (2-3) ($ (2-3)!0.5!(3-2) + (45:0.1) $) (3-2) ($ (3-2)!0.5!(4-1) + (225:0.1) $) (4-1) ($ (4-1)!0.5!(5-0) + (45:0.1) $) (5-0) ($ (5-0)!-0.5!(4-1) $)};
\path[draw, ultra thick] plot[smooth] coordinates{($ (0-7)!-0.5!(1-6) $) (0-7) ($ (0-7)!0.5!(1-6) + (45:0.1) $) (1-6) ($ (1-6)!0.5!(2-5) + (225:0.1) $) (2-5) ($ (2-5)!0.5!(3-4) + (45:0.1) $) (3-4) ($ (3-4)!0.5!(4-3) + (225:0.1) $) (4-3) ($ (4-3)!0.5!(5-2) + (45:0.1) $) (5-2) ($ (5-2)!0.5!(6-1) + (225:0.1) $) (6-1) ($ (6-1)!0.5!(7-0) + (45:0.1) $) (7-0) ($ (7-0)!-0.5!(6-1) $)};
\path[draw, ultra thick] plot[smooth] coordinates{($ (1-8)!-0.5!(2-7) $) (1-8) ($ (1-8)!0.5!(2-7) + (225:0.1) $) (2-7) ($ (2-7)!0.5!(3-6) + (45:0.1) $) (3-6) ($ (3-6)!0.5!(4-5) + (225:0.1) $) (4-5) ($ (4-5)!0.5!(5-4) + (45:0.1) $) (5-4) ($ (5-4)!0.5!(6-3) + (225:0.1) $) (6-3) ($ (6-3)!0.5!(7-2) + (45:0.1) $) (7-2) ($ (7-2)!0.5!(8-1) + (225:0.1) $) (8-1) ($ (8-1)!-0.5!(7-2) $)};
\path[draw, ultra thick] plot[smooth] coordinates{($ (3-8)!-0.5!(4-7) $) (3-8) ($ (3-8)!0.5!(4-7) + (225:0.1) $) (4-7) ($ (4-7)!0.5!(5-6) + (45:0.1) $) (5-6) ($ (5-6)!0.5!(6-5) + (225:0.1) $) (6-5) ($ (6-5)!0.5!(7-4) + (45:0.1) $) (7-4) ($ (7-4)!0.5!(8-3) + (225:0.1) $) (8-3) ($ (8-3)!-0.5!(7-4) $)};
\path[draw, ultra thick] plot[smooth] coordinates{($ (5-8)!-0.5!(6-7) $) (5-8) ($ (5-8)!0.5!(6-7) + (225:0.1) $) (6-7) ($ (6-7)!0.5!(7-6) + (45:0.1) $) (7-6) ($ (7-6)!0.5!(8-5) + (225:0.1) $) (8-5) ($ (8-5)!-0.5!(7-6) $)};
\path[draw, ultra thick] plot[smooth] coordinates{($ (7-8)!-0.5!(8-7) $) (7-8) ($ (7-8)!0.5!(8-7) + (225:0.1) $) (8-7) ($ (8-7)!-0.5!(7-8) $)};
\end{tikzpicture}
\caption{All zigzag curves}
\label{fig:prelim-fukaya-zigzag-products2}
\end{subfigure}
\begin{subfigure}{0.3\linewidth}
\centering
\begin{tikzpicture}
\path[use as bounding box] (-0.5, -0.5) -- (4.5, 4.5);
\path[draw] (0, 0) -- ++(right:1) coordinate[midway] (1-0) -- ++(right:1) coordinate[midway] (3-0) -- ++(right:1) coordinate[midway] (5-0) -- ++(right:1) coordinate[midway] (7-0);
\path[draw] (0, 1) -- ++(right:1) coordinate[midway] (1-2) -- ++(right:1) coordinate[midway] (3-2) -- ++(right:1) coordinate[midway] (5-2) -- ++(right:1) coordinate[midway] (7-2);
\path[draw] (0, 2) -- ++(right:1) coordinate[midway] (1-4) -- ++(right:1) coordinate[midway] (3-4) -- ++(right:1) coordinate[midway] (5-4) -- ++(right:1) coordinate[midway] (7-4);
\path[draw] (0, 3) -- ++(right:1) coordinate[midway] (1-6) -- ++(right:1) coordinate[midway] (3-6) -- ++(right:1) coordinate[midway] (5-6) -- ++(right:1) coordinate[midway] (7-6);
\path[draw] (0, 4) -- ++(right:1) coordinate[midway] (1-8) -- ++(right:1) coordinate[midway] (3-8) -- ++(right:1) coordinate[midway] (5-8) -- ++(right:1) coordinate[midway] (7-8);
\path[draw] (0, 0) -- ++(up:1) coordinate[midway] (0-1) -- ++(up:1) coordinate[midway] (0-3) -- ++(up:1) coordinate[midway] (0-5) -- ++(up:1) coordinate[midway] (0-7);
\path[draw] (1, 0) -- ++(up:1) coordinate[midway] (2-1) -- ++(up:1) coordinate[midway] (2-3) -- ++(up:1) coordinate[midway] (2-5) -- ++(up:1) coordinate[midway] (2-7);
\path[draw] (2, 0) -- ++(up:1) coordinate[midway] (4-1) -- ++(up:1) coordinate[midway] (4-3) -- ++(up:1) coordinate[midway] (4-5) -- ++(up:1) coordinate[midway] (4-7);
\path[draw] (3, 0) -- ++(up:1) coordinate[midway] (6-1) -- ++(up:1) coordinate[midway] (6-3) -- ++(up:1) coordinate[midway] (6-5) -- ++(up:1) coordinate[midway] (6-7);
\path[draw] (4, 0) -- ++(up:1) coordinate[midway] (8-1) -- ++(up:1) coordinate[midway] (8-3) -- ++(up:1) coordinate[midway] (8-5) -- ++(up:1) coordinate[midway] (8-7);
\path[fill=gray, fill opacity=0.5, draw, ultra thick] plot[smooth] coordinates{(4-1) ($ (4-1)!0.5!(5-2) + (315:0.1) $) (5-2) ($ (5-2)!0.5!(6-3) + (135:0.1) $) (6-3) ($ (6-3)!0.5!(7-4) + (315:0.1) $) (7-4) ($ (7-4)!0.5!(8-5) + (135:0.1) $) (8-5)} -- plot[smooth] coordinates{(8-5) ($ (7-6)!0.5!(8-5) + (225:0.1) $) (7-6) ($ (6-7)!0.5!(7-6) + (45:0.1) $) (6-7) ($ (5-8)!0.5!(6-7) + (225:0.1) $) (5-8)} -- plot[smooth] coordinates{(5-8) ($ (4-7)!0.5!(5-8) + (315:0.1) $) (4-7) ($ (3-6)!0.5!(4-7) + (135:0.1) $) (3-6) ($ (2-5)!0.5!(3-6) + (315:0.1) $) (2-5) ($ (1-4)!0.5!(2-5) + (135:0.1) $) (1-4)} -- plot[smooth] coordinates{(1-4) ($ (1-4)!0.5!(2-3) + (225:0.1) $) (2-3) ($ (2-3)!0.5!(3-2) + (45:0.1) $) (3-2) ($ (3-2)!0.5!(4-1) + (225:0.1) $) (4-1)};
\path[fill] (4-1) circle[radius=0.05] node[below right] {$ p_2 $};
\path[fill] (8-5) circle[radius=0.05] node[right] {$ p_1 $};
\path[fill] (1-4) circle[radius=0.05] node[above left] {$ p_3 $};
\path[fill] (5-8) circle[radius=0.05] node[above] {$ p $};
\path (7-6) node[above, shift={(up:0.2)}] {$ \smooth L_1 $};
\path (6-3) node[below right] {$ \smooth L_2 $};
\path (3-2) node[left] {$ \smooth L_3 $};
\path (3-6) node[above left] {$ \smooth L_4 $};
\end{tikzpicture}
\caption{Smooth immersed disk}
\label{fig:prelim-fukaya-zigzag-products3}
\end{subfigure}
\caption{Illustration of zigzag curves}
\end{figure}

Among zigzag curves, it is easy to describe the transversal sequences:

\begin{lemma}
\label{th:prelim-fukaya-zigzag-transversal}
Let $ Q $ be a dimer and $ L_1, …, L_{N+1} $ be a sequence of $ N+1 ≥ 1 $ zigzag paths in $ Q $. Then the sequence of zigzag curves $ (\smooth L_1, …, \smooth L_{N+1}) $ is transversal if and only if the zigzag paths $ L_i $ are pairwise distinct.
\end{lemma}

\begin{proof}
Assume $ (\smooth L_1, …, \smooth L_{N+1}) $ is a transversal sequence. By definition, $ \smooth L_i $ and $ \smooth L_j $ for $ i < j $ only have transversal intersection points. Then in particular $ \smooth L_i ≠ \smooth L_j $, hence $ L_i ≠ L_j $ as zigzag paths. Conversely, assume all zigzag paths are pairwise distinct. Then for $ i < j $ the zigzag curves $ \smooth L_i $ and $ \smooth L_j $ have only transversal intersection points. Moreover, for $ i < j < k $ the zigzag curves $ \smooth L_i $, $ \smooth L_j $ and $ \smooth L_k $ have no common intersection point at all, since an intersection point of zigzag curves is shared between at most two zigzag curves. This shows that $ (\smooth L_1, …, \smooth L_{N+1}) $ is a transversal sequence according to \autoref{def:prelim-fukaya-pre-transversal} and finishes the proof.
\end{proof}

We can compute some $ A_∞ $-products in the category of zigzag curves $ \relFuk Q \restr{\Ob\ZigzagCat} ⊂ \relFuk Q $. Let $ p_1, …, p_N $ be intersection points with $ p_i: \smooth L_i → \smooth L_{i+1} $. If $ (\smooth L_1, …, \smooth L_{N+1}) $ is a transversal sequence, then the product $ μ_{\relFuk Q} (p_N, …, p_1) $ agrees with the product of the relative Fukaya pre-category, which is by definition enumerated by smooth immersed disks with inputs $ p_1, …, p_N $ and arbitrary output $ p: \smooth L_1 → \smooth L_{N+1} $.

\begin{example}
In \autoref{fig:prelim-fukaya-zigzag-products1}, we have depicted a 16-punctured torus dimer. There are 8 zigzag paths and zigzag curves, depicted in \autoref{fig:prelim-fukaya-zigzag-products2}. A sample smooth immersed disk bounded by four zigzag curves is depicted in \autoref{fig:prelim-fukaya-zigzag-products3}. In the figure, we have denoted the four involved zigzag curves by $ \smooth L_1 $, $ \smooth L_2 $, $ \smooth L_3 $, $ \smooth L_4 $. The disk has
\begin{equation*}
\text{inputs} \quad p_1: \smooth L_1 → \smooth L_2, \quad p_2: \smooth L_2 → \smooth L_3, \quad p_3: \smooth L_3 → \smooth L_4, \qquad \text{output} \quad p: \smooth L_1 → \smooth L_4.
\end{equation*}
The inputs $ p_1 $ is odd, and the inputs $ p_2, p_3 $ are even. The output $ p $ is even. It covers six punctures which we denote by $ q_1, …, q_6 $. The contribution of this smooth immersed disk to $ μ^3 (p_3, p_2, p_1) $ is then
\begin{equation*}
± q_1 q_2 q_3 q_4 q_5 q_6 ~ p.
\end{equation*}
\end{example}

\begin{remark}
When $ p_1, …, p_N $ are morphisms $ p_i ∈ \Hom_{\Fuk Q} (\smooth L_i, \smooth L_{i+1}) $ and the sequence $ (\smooth L_1, …, \smooth L_{N+1}) $ is not transversal, then the product $ μ_{\relFuk Q} (p_N, …, p_1) $ is unpredictable. In the present paper, we define a category $ \H\DefZigzagCat ⊂ \HTw\Gtl_q Q $ which has the property that its transversal part agrees with $ \relpreFuk \restr{\Ob\ZigzagCat} $. The products of $ \H\DefZigzagCat $ are explicitly constructed in \autoref{sec:subdisk}. They provide a candidate for describing the products among non-transversal sequences in $ \relFuk Q \restr{\Ob\ZigzagCat} $.

It seems likely that $ \H\DefZigzagCat $ is (gauge equivalent to) $ \relFuk Q \restr{\Ob\ZigzagCat} $ (for any model of $ \relFuk Q $). Our main result is however no guarantee for this, since $ \relFuk Q $ is defined as a lift of the entire relative pre-category and taking subcategories and lifting pre-categories to categories need not commute: Every subcategory of a lift is a lift of the subcategory, but not the other way around.
\end{remark}

Combining \autoref{def:prelim-fukaya-relpre-def}, \autoref{def:prelim-fukaya-rel-def}, \autoref{th:prelim-fukaya-zigzag-Hom} and \autoref{th:prelim-fukaya-zigzag-transversal}, we summarize our findings as follows:

\begin{corollary}
\label{th:prelim-fukaya-zigzag-restr}
Let $ Q $ be a dimer. Then the $ A_∞ $-pre-category $ \preFuk Q \restr{\Ob\ZigzagCat} $ and its $ A_∞ $-pre-category deformation $ \relpreFuk Q \restr{\Ob\ZigzagCat} $ are described as follows:
\begin{itemize}
\item The objects are the zigzag curves $ \smooth L $ for all zigzag paths $ L $, with chosen spin structure.
\item The set of transversal sequences is
\begin{align*}
(\preFuk Q \restr{\Ob\ZigzagCat})^N_{\transversal} = \{(\smooth L_1, …, \smooth L_N) \running ∀i<j: ~ L_i ≠ L_j\}.
\end{align*}
\item For $ \smooth L_1, \smooth L_2 $ with $ L_1 ≠ L_2 $, the hom space is
\begin{equation*}
\Hom_{\preFuk Q \restr{\Ob\ZigzagCat}} (\smooth L_1, \smooth L_2) = \vspan\{p ∈ \smooth L_1 ∩ \smooth L_2\},
\end{equation*}
\item For $ (\smooth L_1, …, \smooth L_{N+1}) ∈ (\preFuk Q \restr{\Ob\ZigzagCat})^{N+1}_{\transversal} $ and $ p_i ∈ \smooth L_i ∩ \smooth L_{i+1} $ we have
\begin{equation*}
μ^N_{\preFuk Q \restr{\Ob\ZigzagCat}} (p_N, …, p_1) = \sum_{p ∈ \smooth L_1 ∩ \smooth L_{N+1}} \sum_{D ∈ M (p_1, …, p_N, p)} (-1)^{\Abouzaid(D)} p.
\end{equation*}
\item For $ (\smooth L_1, …, \smooth L_{N+1}) ∈ (\preFuk Q \restr{\Ob\ZigzagCat})^{N+1}_{\transversal} $ and $ p_i ∈ \smooth L_i ∩ \smooth L_{i+1} $ we have
\begin{equation*}
μ^N_{\relpreFuk Q \restr{\Ob\ZigzagCat}} (p_N, …, p_1) = \sum_{p ∈ \smooth L_1 ∩ \smooth L_{N+1}} \sum_{D ∈ M_q (p_1, …, p_N, p)} (-1)^{\Abouzaid(D)} \punctures(D) p.
\end{equation*}
\end{itemize}
\end{corollary}

\section{Uncurving of strings and bands}
\label{sec:uncurving}
In this section, we show how to remove curvature of band objects in derived deformed gentle algebras. The aim is to find the uncurvable objects of $ \Tw\Gtl_q \cA $ according to the uncurving theory presented in \papertwoA. We start by recalling how objects of $ \Tw\Gtl \cA $ can be classified geometrically as strings and bands, due to \cite{HKK}. All objects of $ \Tw\Gtl \cA $ can naturally be interpreted as objects in $ \Tw\Gtl_q \cA $. However, they do not satisfy the Maurer-Cartan equation of $ \Tw\Gtl_q \cA $ itself and become curved objects. In the present section, we introduce a method to reduce curvature of these curved objects, which we call the “complementary angle trick”.

\begin{center}
\begin{tikzpicture}
\path (0, 0) node[align=center] (A) {\textbf{Band object} \\ $ (X, δ) ∈ \Tw\Gtl \cA $} (6, 0) node[align=center] (B) {\textbf{Curved object} \\ $ (X, δ) ∈ \Tw\Gtl_q \cA $} (12, 0) node[align=center] (C) {\textbf{Uncurved object} \\ $ (X, δ_q) ∈ \Tw'\Gtl_q \cA $};
\path[draw, ->] ($ (A.east)!0.2!(B.west) $) -- ($ (A.east)!0.8!(B.west) $) node[midway, above] {upon deformation};
\path[draw, ->] ($ (B.east)!0.2!(C.west) $) -- ($ (B.east)!0.8!(C.west) $) node[midway, above] {uncurving};
\end{tikzpicture}
\end{center}

We show that the complementary angle trick succeeds in removing curvature from band objects which satisfy a technical condition. Subsequently, we show how to drop this condition. While not an ideal term, we will constantly refer to removing curvature as “uncurving”. The starting point is a full arc system with [NMDC]. We find that uncurvability differs between band objects and string objects and depends also on the topology of the objects. Our findings can be summarized as follows:
\begin{itemize}
\item Uncurvability of band objects: Our main criterion states that a band object can be uncurved in case the underlying curve, regarded as a curve in $ S $, is not contractible and does not include a teardrop.
\item Uncurvability of string objects: The general rule is that string objects cannot be uncurved. There are exceptions, for example a string where both ends touch each other may be uncurvable if the deformation parameter lies in $ \mathfrak{m}^2 $.
\end{itemize}

\subsection{Strings and bands}
\label{sec:uncurving-stringsbands}
In this section, we recall the classification of objects in $ \HTw\Gtl \cA $. This classification is due to Haiden, Katzarkov and Kontsevich and categorizes the objects into two classes, the so-called string and band objects. Roughly speaking, a string object is a non-closed curve running between two punctures of $ \cA $ and a band object is a closed curve that avoids the punctures of $ \cA $. In the present section, we recall the precise classification and how to realize string and band objects explicitly as twisted complexes.

Originally, the gentle algebra $ \Gtl \cA $ was introduced in \cite{Bocklandt} to provide a combinatorial description of the wrapped Fukaya category of $ S \setminus M $. It was shown in that paper's appendix that $ \Gtl \cA $ indeed embeds into $ \wFuk(S \setminus M) $. The question arose whether this embedding is essentially surjective, upon passing to the derived category $ \HTw\Gtl \cA $. If one puts suitable restrictions on the geometry of the objects allowed in $ \wFuk(S \setminus M) $ and works with a $ ℤ $-grading, the embedding is indeed essentially surjective. The fact that an arc system suffices to generate the wrapped Fukaya category has apparently been folklore for longer, and was affirmed by Haiden, Katzarkov and Kontsevich \cite{HKK}.

Until now, we have defined $ \Gtl \cA $ as a $ ℤ/2ℤ $-graded $ A_∞ $-category. In order to state and discuss the classification of objects, we need to recall $ ℤ $-gradings $ \Gtl^ℤ \cA $ on $ \Gtl \cA $. The procedure to upgrade $ \Gtl \cA $ to a $ ℤ $-graded $ A_∞ $-category $ \Gtl^ℤ \cA $ is as follows:
Choose a line field on the surface $ S $, with singularities allowed at the punctures. Choose a grading of the arcs $ a ∈ \cA $ relative to the line field. Define the degree of an angle $ α: a → b $ as the rounding off of the amount it turns relative to the grading of $ a $ and $ b $.

A precise description of this procedure, including definitions of the notions of line field, arc grading, degree of an angle can be found in \cite{HKK}. With this in mind, we are ready to recall the classification of objects.

\begin{theorem}[\cite{HKK}]
Let $ \cA $ be a full arc system with [NMD]. Up to isomorphism, the objects of $ \HTw\Gtl^ℤ \cA $ can be classified as direct sums of the following:
\begin{itemize}
\item \emph{String objects:} graded curves in $ S $, starting and ending at two punctures, but otherwise avoiding punctures and not bounding a teardrop in $ S \setminus M $. The curve is to be considered up to homotopy, keeping endpoints fixed.
\item \emph{Band objects:} closed graded curves in $ S $ with indecomposable local system, avoiding punctures and not bounding a teardrop in $ S \setminus M $. The curve is to be considered up to homotopy and the local system up to isomorphism.
\end{itemize}
\end{theorem}

\begin{figure}
\centering
\begin{subfigure}{0.3\linewidth}
\centering
\begin{tikzpicture}
\path[draw] (0, 0) arc(-90:270:1.5 and 1);
\path[fill] (0, 0) circle[radius=0.05] (0, 2) circle[radius=0.05];
\path[fill] (-0.5, 1) coordinate (A) circle[radius=0.05] (0.5, 1) coordinate (B) circle[radius=0.05];
\path[draw] (0, 2) to[out=270, in=0] ($ (A) + (0, -0.3) $) to[out=180, in=270] ($ (A) + (-0.4, 0) $) to[out=90, in=180] ($ (A) + (0, 0.3) $) to[out=0, in=90] (0, 0);
\end{tikzpicture}
\end{subfigure}
\begin{subfigure}{0.3\linewidth}
\centering
\begin{tikzpicture}
\path[draw] (0, 0) arc(-90:270:1.5 and 1);
\path[fill] (0, 0) circle[radius=0.05] (0, 2) circle[radius=0.05];
\path[fill] (-0.5, 1) coordinate (A) circle[radius=0.05] (0.5, 1) coordinate (B) circle[radius=0.05];
\path[draw] ($ (A) + (0, 0.3) $) to[out=0, in=180] ($ (B) + (0, -0.3) $) to[out=0, in=270] ($ (B) + (0.4, 0) $) to[out=90, in=0] ($ (B) + (0, 0.3) $) to[out=180, in=0] ($ (A) + (0, -0.3) $) to[out=180, in=270] ($ (A) + (-0.4, 0) $) to[out=90, in=180] cycle;
\end{tikzpicture}
\end{subfigure}
\caption{A string object and a band object on the 4-punctured sphere}
\label{fig:strings-bands-examples}
\end{figure}
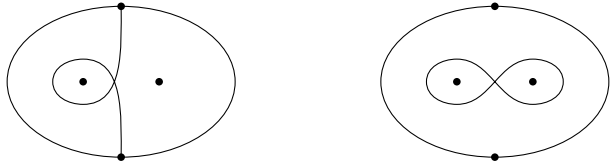

An example string and band object are depicted in \autoref{fig:strings-bands-examples}. Actual representatives of string and band objects as twisted complexes in $ \HTw\Gtl^ℤ \cA $ can be provided by an explicit construction, which we now describe in some detail. The input datum of both constructions is a graded curve of one of the two types above.

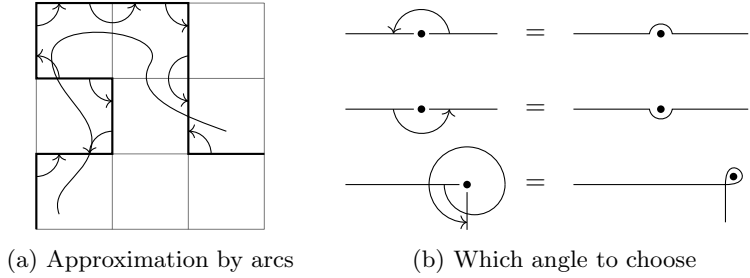
\begin{figure}
\centering
\begin{subfigure}{0.25\linewidth}
\centering
\begin{tikzpicture}
\path[draw, gray] (0, 0) grid (3, 3);
\path[draw] (0.3, 0.2) to[out=110, in=250] (0.7, 1) to[out=70, in=270] (0.2, 2.3) to[out=90, in=60] (1.5, 2.2) to[out=240, in=160] (2.5, 1.3);
\path[draw, thick] (0, 0) -- (0, 1) coordinate[pos=0.7] (1) -- (1, 1) coordinate[pos=0.3] (2) coordinate[pos=0.7] (3) -- (1, 2) coordinate[pos=0.3] (4) coordinate[pos=0.7] (5) -- (0, 2) coordinate[pos=0.3] (6) coordinate[pos=0.7] (7) -- (0,  3) coordinate[pos=0.3] (8) coordinate[pos=0.7] (9) -- (2, 3) coordinate[pos=0.15] (10) coordinate[pos=0.85] (11) coordinate[pos=0.35] (15) coordinate[pos=0.65] (16) -- (2, 1) coordinate[pos=0.15] (12)  coordinate[pos=0.85] (13) coordinate[pos=0.35] (17) coordinate[pos=0.65] (18) -- (3, 1) coordinate[pos=0.3] (14);
\path[draw, ->, bend right=45] (1) to (2);
\path[draw, ->, bend right=45] (4) to (3);
\path[draw, ->, bend right=45] (6) to (5);
\path[draw, ->, bend right=45] (7) to (8);
\path[draw, ->, bend right=45] (9) to (10);
\path[draw, ->, bend right=45] (11) to (12);
\path[draw, ->, bend right=45] (14) to (13);
\path[draw, ->, bend right=80, looseness=1.5] (15) to (16);
\path[draw, ->, bend right=80, looseness=1.5] (17) to (18);
\end{tikzpicture}
\caption{Approximation by arcs}
\end{subfigure}
\begin{subfigure}{0.4\linewidth}
\centering
\begin{tikzpicture}
\path[draw] (-3, 0) -- (-2.1, 0) coordinate[pos=0.7] (1) (-1.9, 0) -- (-1, 0) coordinate[pos=0.3] (2);
\path[draw, ->, bend right=80, looseness=1.5] (2) to (1);
\path[draw] (0, 0) -- (1, 0) to[bend left=80, looseness=1.5] (1.3, 0) -- (2.3, 0);
\path[draw] (-3, -1) -- (-2.1, -1) coordinate[pos=0.7] (3) (-1.9, -1) -- (-1, -1) coordinate[pos=0.3] (4);
\path[draw, ->, bend right=80, looseness=1.5] (3) to (4);
\path[draw] (0, -1) -- (1, -1) to[bend right=80, looseness=1.5] (1.3, -1) -- (2.3, -1);
\path[draw] (-3, -2) -- (-1.5, -2) (-1.4, -2.1) -- ++(down:0.5);
\path (-1.4, -2) coordinate (C);
\path[draw] (C) ++(180:0.3) arc(180:360:0.4);
\path[draw, ->] (C) ++(0:0.5) arc(0:270:0.5);
\path[draw] (0, -2) -- ++(right:2) coordinate (F) .. controls ++(0.5, 0) and ++(0, 0.5) .. ++(0, 0) -- ++(down:0.5);
\path (-0.5, 0) node {=} (-0.5, -1) node {=} (-0.5, -2) node {=};
\path[fill] (-2, 0) circle[radius=0.05] (-2, -1) circle[radius=0.05] (-1.4, -2) circle[radius=0.05];
\path[fill] (1.15, 0) circle[radius=0.05] (1.15, -1) circle[radius=0.05] (F) ++(45:0.15) circle[radius=0.05];
\end{tikzpicture}
\caption{Which angle to choose}
\label{fig:strings-bands-angles}
\end{subfigure}
\caption{Smooth versus discrete}
\label{fig:strings-bands-approximation}
\end{figure}

A string object can be constructed by approximating the curve by a sequence of arcs $ a_1, …, a_k $, and adding up these arcs to form a twisted complex
\begin{equation*}
(a_1 [s_1] ⊕ … ⊕ a_k [s_k], δ) ∈ \Tw\Gtl^ℤ \cA.
\end{equation*}
The shifts $ s_i ∈ ℤ $ are defined such that the degree of $ a_i [s_i] $ with respect to the line field equals the inherent grading of the curve. Let us explain how $ δ $ is found: At every endpoint between two consecutive arcs $ a_i, a_{i+1} $ of the sequence, determine whether the curve runs to the left or to the right of the puncture. The curve either follows an angle $ a_i → a_{i+1} $ or an angle $ a_{i+1} → a_i $. Insert this angle into the $ δ $ matrix. Possibly, the order of the arcs as summands of the twisted complex needs to be reordered to make $ δ $ upper triangular. In summary, the $ δ $ matrix indicates how the arcs $ a_1, …, a_k $ are stitched together. The procedure is depicted in \autoref{fig:strings-bands-approximation}.

\begin{remark}
\label{th:strings-bands-properties}
Let us record the following properties: The twisted differential $ δ $ automatically becomes homogeneous of degree $ 1 $ due to the chosen shifts $ s_i $. Identity angles $ α_i $ are not allowed, and can in fact be avoided when choosing the approximation by arcs. Two consecutive angles $ α_i, α_{i+1} $ are never composable because they turn around the opposite ends of their common arc. One might wonder about composability in case this arc is a loop. Indeed the two consecutive angles may be composable in $ \Gtl \cA $, but are not composable in the order enforced by the ordering on the angles in the twisted complex.
\end{remark}

Band objects can be constructed in a fashion similar to strings. In contrast to strings, their ends are however also stitched together by an angle, and the local system is manifested in the $ δ $ matrix. Let us explain these steps in more detail.
\begin{itemize}
\item Insert not only the angles between two consecutive arcs $ a_i $ and $ a_{i+1} $, but also the angle between the last arc $ a_k $ and the first arc $ a_1 $ into the $ δ $ matrix. Reorder the arcs to make $ δ $ upper triangular. In contrast to the case of strings, such a reordering need not exist. This happens if all arcs are connected entirely cyclically. Abort in this case, and approximate the curve by a different sequence of arcs. It is shown below that this is possible, without inserting identities into $ δ $.
\item The shifts are chosen such that $ δ $ has degree $ 1 $. In contrast to strings which have an entry less in their $ δ $ matrix, it requires a check that this can be done in a consist manner. After walking around the curve one full cycle, do we end up with the same degree shift as we started with? The answer is yes, and the reason is that the curve was required to be graded with respect to the line field.
\item If the local system is of dimension $ d > 1 $, duplicate all arcs in the twisted complex so that each arc appears $ d $ times. Also duplicate the angles so each angle appears $ d $ times, running between the $ i $-th copy of some arc and the $ i $-th copy of the next (without running from one copy to the other).
\item If the local system is non-trivial, represent it as a matrix $ M = (m_{pq}) ∈ ℂ^{d × d} $ and insert it into the $ δ $ matrix as follows: Choose two consecutive arcs $ a_i, a_{i+1} $ in the representation of the curve by arcs. Then change the $ δ $ entry running from the $ q $-th copy of $ a_i $ to the $ p $-th copy of $ a_{i+1} $ to $ m_{pq} $. For instance in case $ d = 2 $, the part of $ δ $ matrix between $ a_i $ and $ a_{i+1} $ shall look like
\begin{equation*}
\left(\begin{array}{c|c} \begin{array}{cc} & \end{array} & \begin{array}{cc} m_{11} α_i & m_{12} α_i \\ m_{21} α_i & m_{22} α_i \end{array} \\\hline \begin{array}{c} \\ \end{array} & \end{array}\right).
\end{equation*}
This overwrites the default identity matrix at those entries written to $ δ $ in step (2). In case the angle between $ a_i $ an $ a_{i+1} $ runs from $ a_{i+1} $ to $ a_i $ instead, use the inverse of $ M $ instead of $ M $ itself. It does not matter which pair of consecutive arcs is chosen. In fact, $ M $ could be arbitrarily factorized into matrices, one for each pair of consecutive arcs, and the values could be written to $ δ $ per pair. It does not matter which factorization we choose: The isomorphism class of the resulting object in $ \HTw\Gtl^ℤ \cA $ only depends on the product of the factors.
\end{itemize}

\begin{lemma}
\label{th:strings-bands-cyclicity}
Let $ \cA $ be a full arc system with [NMD]. Then every string and band has an approximation by arcs where no angle is an identity and all arcs can be ordered such that $ δ $ is upper triangular.
\end{lemma}

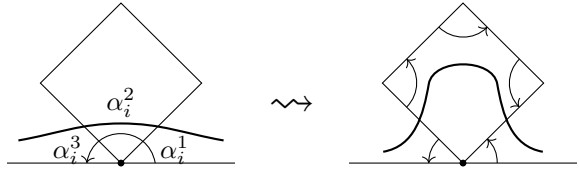
\begin{figure}
\centering
\begin{tikzpicture}[scale=1.5]
\path[draw] (0, 0) -- ++(right:1) coordinate[pos=0.7] (0-end) coordinate (q) -- ++(right:1) coordinate[pos=0.3] (0-start);
\path[draw] (q) -- ++(45:1) -- ++(135:1) -- ++(225:1) -- cycle;
\path[draw, ->, bend right=80, looseness=1.5] (0-start) to node[very near end, left] {$ α_i^3 $} node[midway, above, shift={(0, 0.1)}] {$ α_i^2 $} node[very near start, right] {$ α_i^1 $} (0-end);
\path[draw, thick] (1.9, 0.2) to[out=170, in=0] (1, 0.35) to[out=180, in=10] (0.1, 0.2);
\path (2.5, 0.5) node {\LARGE $ \rightsquigarrow $};
\path[fill] (q) circle[radius=0.03];
\begin{scope}[shift={(3, 0)}]
\path[draw] (0, 0) -- ++(right:1) coordinate[pos=0.7] (1-end) coordinate (q2) -- ++(right:1) coordinate[pos=0.3] (2-start);
\path[draw] (q2) -- ++(45:1) coordinate[pos=0.3] (2-end) coordinate[pos=0.7] (3-end) -- ++(135:1) coordinate[pos=0.3] (3-start) coordinate[pos=0.7] (4-end) -- ++(225:1) coordinate[pos=0.3] (4-start) coordinate[pos=0.7] (5-end) -- cycle coordinate[pos=0.3] (5-start) coordinate[pos=0.7] (1-start);
\foreach \i in {1, 2} \path[draw, ->, bend right=22] (\i-start) to (\i-end);
\foreach \i in {3, 4, 5} \path[draw, ->, bend right=45] (\i-start) to (\i-end);
\path[draw, thick] ($ (q2) + (0.7, 0.1) $) to[out=170, in=280] ($ (q2) + (0.3, 0.7) $) to[out=100, in=80] ($ (q2) + (-0.3, 0.7) $) to[out=260, in=10] ($ (q2) + (-0.7, 0.1) $);
\path[fill] (q2) circle[radius=0.03];
\end{scope}
\end{tikzpicture}
\caption{Removing cyclicity}
\label{fig:strings-bands-angle-3}
\end{figure}

\begin{proof}
First, choose some arbitrary approximation of the curve such that no connecting angle is the identity. This is always possible. The rest of the proof consists of tweaking this approximation such that $ δ $ becomes upper triangular.

For strings, there is always an ordering of the arcs in which $ δ $ is upper triangular, and we are done. For bands however, such an ordering need not exist. That is, the arcs might be connected cyclically. The remaining task in this proof is to break the cyclicity in the band case by tweaking the arc collection.

We may assume that one of the angles $ α_i $ consists of at least three indecomposable components. Otherwise, choose a different approximation by arcs where one angle winds a little more around some puncture.

Regard such an angle $ α_i $ that consists of at least three indecomposable components, and split it into a product $ α_i = α_i^3 α_i^2 α_i^1 $ of three non-empty angles such that $ α_i^2 $ is indecomposable. In particular, $ α_i^2 $ is an interior angle of some polygon. We now modify the angle sequence $ α_1, …, α_k $ by flipping $ α_i $ over to the opposite side of this polygon, see \autoref{fig:strings-bands-angle-3}. This tweak yields a non-cyclic approximation where all angles are still non-empty.
\end{proof}

\begin{figure}
\centering
\begin{tikzpicture}
\path (0, 0) node {$ ℤ $-graded $ \Gtl^ℤ \cA $};
\path (0, -1) node {\LARGE \rotatebox{270}{$ \rightsquigarrow $}};
\path (0, -2) node[align=center] {$ ℤ $-graded $ \HTw\Gtl^ℤ \cA $ \\ classification: strings and bands};
\path (10, 0) node {$ ℤ/2ℤ $-graded $ \Gtl \cA $};
\path (10, -1) node {\LARGE \rotatebox{270}{$ \rightsquigarrow $}};
\path (10, -2) node[align=center] {$ ℤ/2ℤ $-graded $ \HTw\Gtl \cA $ \\ classification unknown};
\path[draw, ->] (3, -2) to[bend left=0] node[midway, above, align=center] {Strings and bands can also be \\ constructed $ ℤ/2ℤ $-graded} (7, -2);
\end{tikzpicture}
\caption{Overview on objects of $ ℤ $-vs. $ ℤ/2ℤ $-graded gentle algebras}
\label{fig:strings-bands-overview}
\end{figure}

Let us discuss this classification in the context of the grading question. What are the objects up to isomorphism of $ \HTw\Gtl \cA $? Here $ \Gtl \cA $ and its twisted completion are taken as usual with $ ℤ/2ℤ $-grading. The answer is, there are both more and less objects. For example, twisted complexes differing only by even shifts are unequal in the $ ℤ $-grading, but are identified in the $ ℤ/2ℤ $-grading. On the other hand, some twisted complexes that can be made with respect to the $ ℤ/2ℤ $-grading cannot be constructed in the $ ℤ $-grading. One might say, the $ ℤ/2ℤ $-graded $ \HTw\Gtl \cA $ category contains all objects of $ \HTw\Gtl^ℤ \cA $ for every possible $ ℤ $-grading, plus additional objects that cannot be obtained from a $ ℤ $-graded version, modulo identifying objects differing by shifts. This is depicted in the overview \autoref{fig:strings-bands-overview}. We are however not aware of a concise classification of the objects of $ \HTw\Gtl \cA $.

A broad class of objects in $ \HTw\Gtl \cA $ can however be constructed by forming $ ℤ/2ℤ $-graded strings and bands, corresponding to curves on $ S $ without grading requirements. Let us define what we mean by this:
\begin{itemize}
\item \emph{$ ℤ/2ℤ $-graded string objects:} Stitch arcs together as in the $ ℤ $-graded case, and shift arcs with $ s_i ∈ ℤ/2ℤ $ such that the degree of $ δ $ is odd.
\item \emph{$ ℤ/2ℤ $-graded band objects:} Stitch arcs together as in the $ ℤ $-graded case, shift arcs with $ s_i ∈ ℤ/2ℤ $ such that the degree of $ δ $ is odd. The shifts are consistent: After cycling around the curve once, we end up with the same shift because the curve is orientable.
\end{itemize}

\subsection{Complementary angle trick}
\label{sec:uncurving-trick}
In this section, we introduce a method to uncurve objects of $ \Tw\Gtl_q \cA $. The starting point is the classification of objects in $ \Tw\Gtl \cA $ recalled in \autoref{sec:uncurving-stringsbands}. The idea to uncurve these objects is to infinitesimally deform their $ δ $-matrix by inserting infinitesimal multiples of the complements of the angles already present in the $ δ $-matrix. We therefore call this method the “complementary angle trick”. In the present section, we show that this trick successfully uncurves $ ℤ/2ℤ $-graded band objects under certain conditions.

We will start the setup in a slightly more general approach: We take $ \cA $ to denote a full arc system with [NMDC] and we take the category $ \Gtl_r \cA $ to be the associated deformed gentle algebra constructed in \paperone. This deformed gentle algebra depends on a parameter
\begin{equation*}
r = r_0 1 + \sum_{\substack{q ∈ M \\ n ≥ 1}} r_{q, n} ℓ_q^n ∈ \mathfrak{m} Z(\Gtl \cA).
\end{equation*}
Here $ B $ is a chosen deformation base with maximal ideal $ \mathfrak{m} $. We say that $ r $ is \emph{without 1-component} if $ r_0 = 0 $. For simplicity, we write $ μ_q $ for the product of $ \Gtl_r \cA $.

\begin{example}
The category $ \Gtl_r \cA $ may simply be the standard deformation $ \Gtl_r \cA = \Gtl_q \cA $ over $ B = ℂ⟦M⟧ $, detailed in \autoref{sec:prelim-gtlq}. It is determined by the specific parameter
\begin{equation*}
r = \sum_{q ∈ M} q ℓ_q ∈ (M) Z(\Gtl \cA).
\end{equation*}
This parameter is without 1-component.
\end{example}

Recall the notion of uncurvability: Let $ \cat C $ be an $ A_∞ $-category and $ \cat C_q $ a deformation of $ \cat C $. Then an object $ X ∈ \cat C_q $ is uncurvable if there is an odd $ S ∈ \mathfrak{m} \End_{\cat C} (X) $ such that $ μ^0_{q, X} + μ^1_q (S) + μ^2_q (S, S) + … = 0 $. For twisted complexes $ X = (\sum X_i [s_i], δ) $, this means to find an infinitesimal deformation $ δ_q $ of $ δ $ such that the curvature of the twisted complex becomes zero:
\begin{equation*}
μ_{q, X}^0 + μ_{\Add\Gtl_r \cA}^1 (δ_q) + μ^2_{\Add\Gtl_r \cA} (δ_q, δ_q) + … = 0.
\end{equation*}
The infinitesimal part of $ δ_q $ is allowed to lie anywhere in the matrix, not restricted to the upper-triangular part.

Let us now describe our “complementary angle trick”. Regard a band object $ X = (a_1 [s_1] ⊕ …, δ) $. To simplify the discussion, we assume its local system is one-dimensional with transition value simply equal to $ 1 $. In the twisted complex $ X $, every arc then only appears once (apart from arcs appearing multiple times in the approximation), and every angle $ α_i $ appearing in $ δ $ appears with a coefficient of $ +1 $.

The curvature $ μ^0_{q, X} $ of $ X $ consists by definition of the sum of the curvatures $ μ^0_{a_i} $ of the constituent arcs of $ X $. Since we are regarding a standard deformation $ \Gtl_r \cA $ of $ \Gtl \cA $, we know this curvature explicitly: Every arc carries an infinitesimal amount of turns around both of its endpoints as curvature. Since $ X $ is a band object, both endpoints of every arc $ a_i $ are connected to the predecessor or successor arc $ a_{i-1} $ and $ a_i $ by angles $ α_i $.

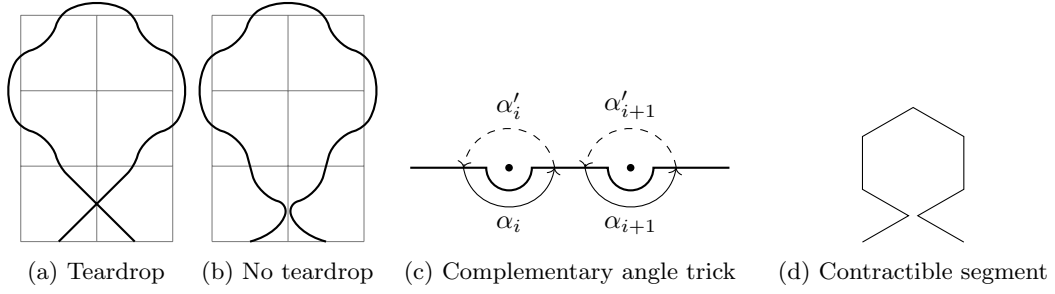
\begin{figure}
\centering
\begin{subfigure}[b]{0.15\linewidth}
\centering
\begin{tikzpicture}
\path[draw, gray] (0, 0) grid (2, 3);
\path[draw, thick, rounded corners] (1.5, 0) to (1, 0.5) to (0.5, 1) to[bend right] (0, 1.5) to[bend left] (0, 2.5) to[bend right] (0.5, 3) to[bend left] (1.5, 3) to[bend right] (2, 2.5) to[bend left] (2, 1.5) to[bend right] (1.5, 1) to (1, 0.5) to (0.5, 0);
\end{tikzpicture}
\caption{Teardrop}
\label{fig:uncurving-trick-teardrop}
\end{subfigure}
\begin{subfigure}[b]{0.15\linewidth}
\centering
\begin{tikzpicture}
\path[draw, gray] (0, 0) grid (2, 3);
\path[draw, thick, rounded corners] (0.5, 0) to[bend right] (1, 0.5) to[bend left] (0.5, 1) to[bend right] (0, 1.5) to[bend left] (0, 2.5) to[bend right] (0.5, 3) to[bend left] (1.5, 3) to[bend right] (2, 2.5) to[bend left] (2, 1.5) to[bend right] (1.5, 1) to[bend left] (1, 0.5) to[bend right] (1.5, 0);
\end{tikzpicture}
\caption{No teardrop}
\label{fig:uncurving-trick-noteardrop}
\end{subfigure}
\begin{subfigure}[b]{0.3\linewidth}
\centering
\begin{tikzpicture}
\path[draw, thick] (0, 0) -- ++(right:1) coordinate[pos=0.7] (1) arc(180:360:0.3) -- ++(right:1) coordinate[pos=0.3] (2) coordinate[pos=0.7] (3) arc(180:360:0.3) -- ++(right:1) coordinate[pos=0.3] (4);
\path[fill] (1.3, 0) circle[radius=0.05];
\path[fill] (2.9, 0) circle[radius=0.05];
\path[draw, bend right=80, looseness=1.5, ->] (1) to node[midway, below] {$ α_i $} (2);
\path[draw, bend right=80, looseness=1.5, ->] (3) to node[midway, below] {$ α_{i+1} $} (4);
\path[draw, bend right=80, looseness=1.5, dashed, ->] (2) to node[midway, above] {$ α_i' $} (1);
\path[draw, bend right=80, looseness=1.5, dashed, ->] (4) to node[midway, above] {$ α_{i+1}' $} (3);
\end{tikzpicture}
\caption{Complementary angle trick}
\label{fig:uncurving-trick-trick}
\end{subfigure}
\begin{subfigure}[b]{0.25\linewidth}
\centering
\begin{tikzpicture}[scale=0.7]
\path[draw] (0, 0) -- ++(30:1) -- ++(150:1) -- ++(up:1) -- ++(30:1.1) -- ++(330:1.1) -- ++(down:1) -- ++(210:1) -- ++(330:1);
\end{tikzpicture}
\caption{Contractible segment}
\label{fig:uncurving-trick-contractible}
\end{subfigure}
\caption{Illustration of the complementary angle tick and its technicalities}
\end{figure}

The trick to uncurving is to add the complements $ α_i' $ of these angles to the $ δ $ matrix, depicted in \autoref{fig:uncurving-trick-trick}. Generically denote by $ ℓ $ a full turn around any puncture. We denote by $ r ℓ^{-1} α_i' $ the element of $ B \htensor \Gtl \cA $ obtained from $ r $ by extracting the part that winds around the same puncture as $ α_i' $, shortening all angles by one full turn, and multiplying by $ α_i' $. Naturally, the element $ r ℓ^{-1} α_i' $ can be interpreted as an odd morphism lying in $ \Hom_{\Add\Gtl_r \cA} (X, X) $.

\begin{definition}
Let $ \cA $ be a full arc system with [NMDC], $ B $ a deformation base and $ r ∈ \mathfrak{m} Z(\Gtl \cA) $ a parameter without 1-component. Regard a band object $ X = (a_1 [s_1] ⊕ …, δ) $ with trivial 1-dimensional local system. Assume its $ δ $-angles $ α_1, …, α_k $ are all shorter than a full turn and not identities. Then the \emph{complementary angle trick} associates to $ X $ the twisted complex $ (a_1 [s_1] ⊕ …, δ_q) ∈ \Tw'\Gtl_q \cA $ with $ δ_q $ given by
\begin{equation*}
δ_q ≔ δ + δ' = \sum_{i = 1}^k α_i + r ℓ^{-1} α_i'.
\end{equation*}
\end{definition}

\begin{example}
Regard the deformation $ \Gtl_r \cA $ with parameter $ r = q ℓ_p $ over $ B = ℂ⟦q⟧ $, with $ ℓ_p $ denoting the central element consisting of single turns around the puncture $ p $. Then the deformation entry added to the $ δ $ matrix is simply $ r ℓ^{-1} α_i' = q α_i' $ for every $ α_i $ which winds around $ p $.
\end{example}

\begin{remark}
The deformation of $ δ $ to $ δ_q $ happens precisely on the opposite side of the diagonal of the matrix: If $ α_i $ is in the $ δ $ matrix as angle from $ a_i $ to $ a_{i+1} $, then the angle $ r ℓ^{-1} α_i' $ is inserted within the $ δ $ matrix as morphism from $ a_{i+1} $ to $ α_i $. If $ α_i $ was an angle from $ a_{i+1} $ to $ a_i $, then $ r ℓ^{-1} α_i' $ is inserted as morphism from $ a_i $ to $ a_{i+1} $. This way we obtain a matrix $ δ_q $ with an infinitesimal lower-triangular part.
\end{remark}

\begin{remark}
The complementary angle trick also works when the local system is higher-dimensional and non-trivial. Recall that the $ δ $ matrix encodes the transition matrix $ M $ of a higher-dimensional local system by carrying its entries $ m_{pq} $ in front of the $ α_i $ morphism from the $ q $-th copy of $ a_i $ to the $ p $-th copy of $ a_{i+1} $. In case the angle runs in the opposite direction, the $ δ $ matrix encodes $ M $ by carrying the entries $ m^{pq} ≔ (M^{-1})_{pq} $ of the inverse matrix.

In order to uncurve this band object, we include the inverse matrix $ M^{-1} $ in the uncurving deformation of $ δ $, or $ M $ in case the angle runs in opposite direction. For instance if $ d = 2 $ and $ α_i $ runs from $ a_i $ to $ a_{i+1} $, the part of the $ δ_q $ matrix between $ a_i $ and $ a_{i+1} $ shall read
\begin{equation*}
\text{part of } δ_q = \left(\begin{array}{c|c} \begin{array}{cc} & \end{array} & \begin{array}{cc} m_{11} α_i & m_{12} α_i \\ m_{21} α_i & m_{22} α_i \end{array} \\\hline \begin{array}{cc} m^{11} r ℓ^{-1} α_i' & m^{12} r ℓ^{-1} α_i' \\ m^{21} r ℓ^{-1} α_i' & m^{22} r ℓ^{-1} α_i' \end{array} & \end{array}\right).
\end{equation*}
\end{remark}

We are now ready to check that the complementary angle trick succeeds in uncurving band objects. The trick however comes with strict conditions. We use the following terminology:

\begin{definition}
An \emph{indexed arc} on $ X $ is one of the arcs $ a_i $ of $ X $, remembering the index $ i $. An \emph{segment} of indexed arcs on $ X $ is a sequence of consecutive arcs $ a_i, a_{i+1}, …, a_{i+j} $, remembering the indices. An indexed segment is \emph{contractible in $ S $} if it returns to the same puncture as it started from and the loop defined this way is contractible in the closed surface $ S $.
\end{definition}

\begin{remark}
A band object which bounds a teardrop in $ S $ is depicted in \autoref{fig:uncurving-trick-teardrop}. Meanwhile, the shape depicted in \autoref{fig:uncurving-trick-noteardrop} does not constitute a teardrop. A contractible segment of indexed arcs is depicted in \autoref{fig:uncurving-trick-contractible}. The existence of a contractible segment of arcs on $ X $ does not imply that the underlying curve of $ X $ has a teardrop.
\end{remark}

As we shall see in \autoref{th:uncurving-trick-works}, the complementary angle trick succeeds in uncurving $ X ∈ \Tw\Gtl_r \cA $ when we assume the following three conditions on $ X $:

\begin{itemize}
\item The underlying curve of $ X $, regarded as a curve in the closed surface $ S $, is not contractible and does not bound a teardrop.
\item All angles $ α_i $ in the $ δ $ matrix are non-identities and strictly smaller than a full turn.
\item No segment of indexed arcs of $ X $ is contractible.
\end{itemize}

In \autoref{sec:trick}, we explain how to abandon the condition that segments of indexed arcs are not contractible. Without the condition, one has to add further angles to $ δ_q $ for every location where $ X $ comes close to itself, other than only the complementary angles $ α_i' $. It is interesting to note that geometrically these additional angles can be interpreted as “complementary to segments of indexed arcs” of $ X $.

\begin{lemma}
\label{th:uncurving-trick-works}
Let $ \cA $ be a full arc system with [NMDC]. Regard a standard deformation $ \Gtl_r \cA $ of $ \Gtl \cA $ by some $ r ∈ \mathfrak{m} Z(\Gtl \cA) $ without $ 1 $-component. Let $ X ∈ \Tw\Gtl \cA $ be a $ ℤ/2ℤ $-graded band object whose underlying curve in $ S $ is not contractible and does not bound a teardrop. Assume that all angles in $ X $ are non-identities and shorter than full turns. Then $ X ∈ \Tw\Gtl_r \cA $ is uncurvable.
\end{lemma}

\begin{proof}
Without loss of generality, we assume that $ X $ has one-dimensional local system and all transition values in $ δ $ are $ 1 $. The case with contractible segments is dealt with in \autoref{sec:trick}. We shall therefore assume that $ X $ has no contractible segments of indexed arcs.

Now let us prove that the complementary angle trick successfully uncurves $ X $. This entails checking that the curvature of $ (⊕ a_i [s_i], δ_q) ∈ \Tw'\Gtl_r \cA $ vanishes. Explicitly, the curvature is
\begin{equation*}
\sum_{k ≥ 0} μ^k_{\Add\Gtl_r \cA} (δ_q, …, δ_k).
\end{equation*}
The summand at $ k = 0 $ is the curvature $ μ_{X, q}^0 $, explicitly the sum of the curvatures of the individual arcs. Note that $ μ_q^1 $ vanishes due to [NMDC]. In the first step of the proof, we show that $ μ_{X, q}^0 $ precisely cancels $ μ_{\Add\Gtl_r \cA}^2 (δ_q, δ_q) $. In the second step of the proof, we show that all the higher terms $ μ_q^{≥3} (δ_q, …) $ vanish.

Let us analyze $ μ_q^2 (δ_q, δ_q) $. Since $ δ_q = δ + δ' $ is the sum of the original $ δ $ and the modification $ δ' $ due to the complementary angle trick, we need to check the original part $ μ_q^2 (δ, δ) $ and the new components $ μ_q^2 (δ, δ') $, $ μ_q (δ', δ) $ and $ μ_q^2 (δ', δ') $. Recall also that the product $ μ_q^2 $ is not deformed: It is merely the $ B $-linear extension of the original $ μ^2 $ by assumption of [NMDC]. As observed in \autoref{th:strings-bands-properties}, we have $ μ_q^2 (δ, δ) = 0 $. Similarly, there are no products of complementary angles possible, so $ μ_q^2 (δ', δ') = 0 $.

Meanwhile, we have
\begin{equation*}
μ^2_{\Add\Gtl_r \cA} (δ, δ') = \sum_{i=1}^k - α_i r ℓ^{-1} α_i' = - r ∈ \mathfrak{m} Z(\Gtl \cA),
\end{equation*}
where $ r $ on the right-hand side is interpreted as linear combination of powers of full turns starting at those arc ends of the arc approximation where the angles $ α_i $ enter. The minus sign comes from the sign convention for $ \Add\Gtl_r \cA $. Similarly,
\begin{equation*}
μ^2_{\Add\Gtl_r \cA} (δ, δ') = \sum_{i=1}^k - r ℓ^{-1} α_i' α_i = - r ∈ \mathfrak{m} Z(\Gtl \cA),
\end{equation*}
where $ r $ on the right-hand side is interpreted as linear combination of powers of full turns starting at those arc ends of the arc approximation where the angles $ α_i $ leave.

We have used that all angles in the twisted differential of $ X $ are non-identities and shorter than full turns: Identities would give extra terms in the products $ μ^2 (δ, δ') $ and $ μ^2 (δ', δ) $. Moreover, taking complementary angles is only possible if all angles in $ X $ are at most full turns. A precise full turn would in turn give an identity in $ δ' $, hence an undesired contribution to e.g.~$ μ^2 (δ', δ) $. In short, we assumed just the right condition so that nothing but the right terms appears in $ μ^2 (δ_q, δ_q) $. We conclude
\begin{equation*}
μ_{X, q}^0 + μ_{\Add\Gtl_r \cA}^2 (δ_q, δ_q) = 0.
\end{equation*}
For the second step of the proof, we show that $ μ_q^{k≥3} (δ_q, …, δ_q) = 0 $. Assume $ D $ is an orbigon contributing to this product. Then the boundary of $ D $ is a contractible indexed segment of $ X $, in contradiction to the assumption that there are no contractible indexed segments. This shows that $ μ_q^{k≥3} (δ_q, …, δ_q) = 0 $. Finally, we conclude that the curvature of the deformed twisted complex $ X_q = (⊕ a_i [s_i], δ_q) $ vanishes. This finishes the proof.
\end{proof}

\subsection{The uncurvable objects}
\label{sec:uncurving-gtl}
In this section, we show that most band objects in $ \Tw\Gtl_q \cA $ are uncurvable. The starting point is the classification of band objects recalled in \autoref{sec:uncurving-stringsbands} and the complementary angle trick defined in \autoref{sec:uncurving-trick}. The goal is to show that a band object is uncurvable if its underlying curve in $ S $ is topologically nontrivial and does not bound a teardrop. We have already shown in \autoref{th:uncurving-trick-works} that the complementary angle trick succeeds in uncurving these objects under the technical condition that all angles in the $ δ $-matrix of $ X $ are non-identities and shorter than full turns. In the present section, we show how to abandon this technical condition.

Our starting point is again a full arc system $ \cA $ with [NMDC] and a deformation $ \Gtl_r \cA $ with $ r ∈ \mathfrak{m} Z(\Gtl \cA) $ a deformation parameter without 1-component. It is our wish to apply the complementary angle trick to every band object whose underlying curve in $ S $ is topologically nontrivial and does not bound a teardrop. Combining \autoref{rem:2Bprelim-ainfty-defo-uncurvingqi} and \autoref{th:uncurving-trick-works}, we would be done if every such band object has a twisted complex representation where all angles are non-identities and shorter than full turns. This is however not the case:

\begin{remark}
\label{rem:uncurving-representation-norep}
There are arc systems in which some bands fail to have representatives which satisfy the requirements of \autoref{th:uncurving-trick-works}. For example, regard the 1-punctured torus with two arcs depicted in \autoref{fig:uncurving-trick-norep}. Its horizontal, or vertical, band cannot be represented as a twisted complex with all angles non-identities and shorter than a full turn. As soon as we divide the 4-gon into two triangles, the band suddenly has a desired representation.
\end{remark}

There are however arc systems which guarantee the existence of representatives suitable for \autoref{th:uncurving-trick-works} for every band object whose underlying curve in $ S $ is topologically nontrivial and does not bound a teardrop. Following \autoref{rem:uncurving-representation-norep}, the idea is to simply require that the arc system contains only of triangles:

\begin{figure}
\centering
\begin{subfigure}{0.3\linewidth}
\centering
\begin{tikzpicture}[scale=1.2]
\path[draw, ->] (0, 0) -- (1, 0) node[midway, below] {$ a $} coordinate[pos=0.3] (alpha-start) coordinate[pos=0.7] (beta-end);
\path[draw, ->] (0, 1) -- (1, 1) node[midway, above] {$ a $} coordinate[pos=0.3] (delta-end) coordinate[pos=0.7] (gamma-start);
\path[draw, ->] (0, 0) -- (0, 1) node[pos=0.55, left] {$ b $} coordinate[pos=0.3] (alpha-end) coordinate[pos=0.7] (delta-start);
\path[draw, ->] (1, 0) -- (1, 1) node[pos=0.55, right] {$ b $} coordinate[pos=0.3] (beta-start) coordinate[pos=0.7] (gamma-end);
\path[draw, ->, bend right=45] (alpha-start) to node[near end, below] {\small $ α $} (alpha-end);
\path[draw, ->, bend right=45] (beta-start) to node[near end, right] {\small $ β $} (beta-end);
\path[draw, ->, bend right=45] (gamma-start) to node[near start, right] {\small $ γ $} (gamma-end);
\path[draw, ->, bend right=45] (delta-start) to node[near end, left] {\small $ δ $} (delta-end);
\path[draw, thick] (-0.5, 0.4) -- (1.5, 0.4);
\end{tikzpicture}
\end{subfigure}
\begin{subfigure}{0.3\linewidth}
\centering
\begin{tikzpicture}[scale=1.2]
\path[draw, ->] (0, 0) -- (1, 0) node[midway, below] {$ a $};
\path[draw, ->] (0, 1) -- (1, 1) node[midway, above] {$ a $};
\path[draw, ->] (0, 0) -- (0, 1) node[midway, left] {$ b $} coordinate[pos=0.3] (1-end);
\path[draw, ->] (1, 0) -- (1, 1) node[midway, right] {$ b $} coordinate[pos=0.7] (2-end);
\path[draw] (0, 0) -- (1, 1) node[midway, above] {$ c $} coordinate[pos=0.3] (1-start) coordinate[pos=0.7] (2-start);
\path[draw, thick] (0, 0.5) to[bend left] (0.5, 0.5) to[bend right] (1, 0.5);
\path[draw, ->, bend right=45] (1-start) to (1-end);
\path[draw, ->, bend right=45] (2-start) to (2-end);
\end{tikzpicture}
\end{subfigure}
\caption{The horizontal band in the 1-punctured torus. This band only has a suitable twisted complex representation when the 4-gon is divided into triangles.}
\label{fig:uncurving-trick-norep}
\end{figure}
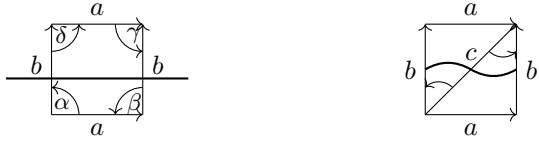

\begin{lemma}
\label{th:uncurving-representation-cyclicity}
Let $ \cA $ be a full arc system with [NMDC] and assume all polygons in $ \cA $ are triangles. Let $ X ∈ \HTw\Gtl \cA $ be a string or band object whose underlying curve in $ S $ is nontrivial and does not bound a teardrop. Then $ X $ has a twisted complex representation in which all angles $ α_i $ are non-identities and shorter than a full turn.
\end{lemma}

\begin{figure}
\centering
\begin{subfigure}{0.49\linewidth}
\centering
\begin{tikzpicture}[scale=1.5]
\path[draw] (0, 0) -- ++(right:1) coordinate[pos=0.7] (0-end) coordinate (q) -- ++(right:1) coordinate[pos=0.3] (0-start);
\path[draw] (q) -- ++(60:1) -- ++(left:1) -- cycle;
\path[draw, ->, bend right=80, looseness=1.5] (0-start) to (0-end);
\path[draw, thick] (1.7, 0.2) to[out=170, in=0] (1, 0.4) to[out=180, in=10] (0.3, 0.2);
\path (2.5, 0.5) node {\LARGE $ \rightsquigarrow $};
\path[fill] (q) circle[radius=0.03];
\begin{scope}[shift={(2.5, 0)}]
\path[draw] (0, 0) -- ++(right:1) coordinate[pos=0.7] (1-end) coordinate (q2) -- ++(right:1) coordinate[pos=0.3] (2-start);
\path[draw] (q2) -- ++(60:1) coordinate[pos=0.3] (2-end) coordinate[pos=0.7] (3-end) -- ++(left:1) coordinate[pos=0.3] (3-start) coordinate[pos=0.7] (4-end) -- cycle coordinate[pos=0.3] (4-start) coordinate[pos=0.7] (1-start);
\foreach \i in {1, 2, 3, 4} \path[draw, ->, bend right=30] (\i-start) to (\i-end);
\path[draw, thick] ($ (q2) + (0.7, 0.1) $) to[out=170, in=280] ($ (q2) + (0.3, 0.5) $) to[out=100, in=80] ($ (q2) + (-0.3, 0.5) $) to[out=260, in=10] ($ (q2) + (-0.7, 0.1) $);
\path[fill] (q2) circle[radius=0.03];
\end{scope}
\end{tikzpicture}
\caption*{Case (a)}
\end{subfigure}
\begin{subfigure}{0.49\linewidth}
\centering
\begin{tikzpicture}[scale=1.5]
\path[draw] (0, 0) -- ++(right:0.5) coordinate[pos=0.6] (1-end) coordinate (q) -- ++(right:1) coordinate[pos=0.2] (1-start) coordinate[pos=0.8] (2-end) coordinate (q2) -- ++(right:0.5) coordinate[pos=0.4] (2-start);
\path[draw, ->, bend right=80, looseness=1.5] (1-start) to (1-end);
\path[draw, ->, bend right=80, looseness=1.5] (2-start) to (2-end);
\path[draw] (q) -- ++(60:1) -- ++(300:1);
\path[fill] (q) circle[radius=0.03] (q2) circle[radius=0.03];
\path[draw, thick] (2, 0.1) to[out=120, in=0] (1, 0.1) to[out=180, in=60] (0, 0.1);
\path (2.5, 0.5) node {\LARGE $ \rightsquigarrow $};
\begin{scope}[shift={(2.5, 0)}]
\path[draw] (0, 0) -- ++(right:0.5) coordinate[pos=0.6] (1-end) coordinate (q) -- ++(right:1) coordinate (q2) -- ++(right:0.5) coordinate[pos=0.4] (2-start);
\path[draw] (q) -- ++(60:1) coordinate[pos=0.2] (1-start) coordinate[pos=0.8] (3-start) -- ++(300:1) coordinate[pos=0.8] (2-end) coordinate[pos=0.2] (3-end);
\path[draw, ->, bend right=60] (1-start) to (1-end);
\path[draw, ->, bend right=60] (2-start) to (2-end);
\path[draw, ->, bend right=30] (3-start) to (3-end);
\path[fill] (q) circle[radius=0.03] (q2) circle[radius=0.03];
\path[draw, thick] (2, 0.1) to[out=150, in=0] (1, 0.5) to[out=180, in=30] (0, 0.1);
\end{scope}
\end{tikzpicture}
\caption*{Case (c)}
\end{subfigure}
\begin{subfigure}{0.24\linewidth}
\centering
\begin{tikzpicture}[scale=1.5]
\path[draw] (0, 0) -- ++(right:1) coordinate[pos=0.3] (1-start) coordinate[pos=0.7] (2-end) -- ++(120:1);
\path[draw, <-] (1, 0) -- ++(120:1) coordinate (top) node[near end, right] {$ a $};
\path[draw, <-] (top) -- ++(240:1) node[near start, left] {$ a $};
\path[draw, ->] (0, 0) ++(0:0.2) arc(0:120:0.2) node[below] {$ α_1 $};
\path[draw, ->] (1, 0) ++(60:0.2) arc(60:180:0.2) node[at start, below right] {$ α_1 $};
\path[fill] (0.7, 0.4) circle[radius=0.02] (0.3, 0.4) circle[radius=0.02];
\path[draw, thick] (1.2, 0.2) to[out=120, in=0] (0.5, 0.1) to[out=180, in=60] (-0.2, 0.2);
\path (1.5, 0.5) node {\LARGE $ \lightning $};
\end{tikzpicture}
\caption*{Case (b)}
\end{subfigure}
\begin{subfigure}{0.24\linewidth}
\centering
\begin{tikzpicture}[scale=1.5]
\path[draw] (0, 0) -- ++(60:1) coordinate[pos=0.3] (1-end) coordinate[pos=0.7] (2-start) -- ++(300:1) -- ++(left:1) coordinate[pos=0.3] (2-end) coordinate[pos=0.7] (1-start) ;
\path[draw, ->, bend right=30] (1-start) to node[midway, right] {$ α_1 $} (1-end);
\path[draw, ->] (2-start) arc(240:360:0.3) node[near end, right] {$ α_2 $};
\path[draw, <-] (2-end) arc(180:60:0.3) node[at end, below] {$ α_2 $};
\path[draw, ->] (1, 0) -- (60:1);
\path (1.5, 0.5) node {\LARGE $ \lightning $};
\end{tikzpicture}
\caption*{Case (d)}
\end{subfigure}
\begin{subfigure}{0.49\linewidth}
\centering
\begin{tikzpicture}[scale=1.5]
\path[draw] (0, 0) -- ++(30:1) coordinate (top) coordinate[pos=0.2] (3-end) coordinate[pos=0.8] (1-start) -- ++(right:0.5) coordinate[pos=0.4] (1-end);
\path[draw] (0, 0) -- ++(330:1) coordinate[pos=0.8] (2-end) coordinate[pos=0.2] (3-start) -- ++(right:0.5) coordinate[pos=0.4] (2-start);
\path[draw] (top) -- ++(down:1);
\path[draw, ->, bend right=45] (1-start) to (1-end);
\path[draw, ->, bend right=45] (2-start) to (2-end);
\path[draw, ->, bend right=30] (3-start) to (3-end);
\path[draw, thick] (1.5, -0.3) to[out=160, in=270] (0.4, 0) to[out=90, in=200] (1.5, 0.3);
\path (1.5, 0) node {\LARGE $ \rightsquigarrow $};
\begin{scope}[shift={(2, 0)}]
\path[draw] (0, 0) -- ++(30:1) coordinate (top) -- ++(right:0.5) coordinate[pos=0.4] (1-end);
\path[draw] (0, 0) -- ++(330:1) -- ++(right:0.5) coordinate[pos=0.4] (2-start);
\path[draw] (top) -- ++(down:1) coordinate[pos=0.2] (1-start) coordinate[pos=0.8] (2-end);
\path[draw, ->, bend right=45] (1-start) to (1-end);
\path[draw, ->, bend right=45] (2-start) to (2-end);
\path[draw, thick] (1.5, -0.3) to[out=160, in=270] (1, 0) to[out=90, in=200] (1.5, 0.3);
\end{scope}
\end{tikzpicture}
\caption*{Case (e)}
\end{subfigure}
\caption{Removing cyclicity in the proof of \autoref{th:uncurving-representation-cyclicity}}
\label{fig:uncurving-representation-cyclicity}
\end{figure}
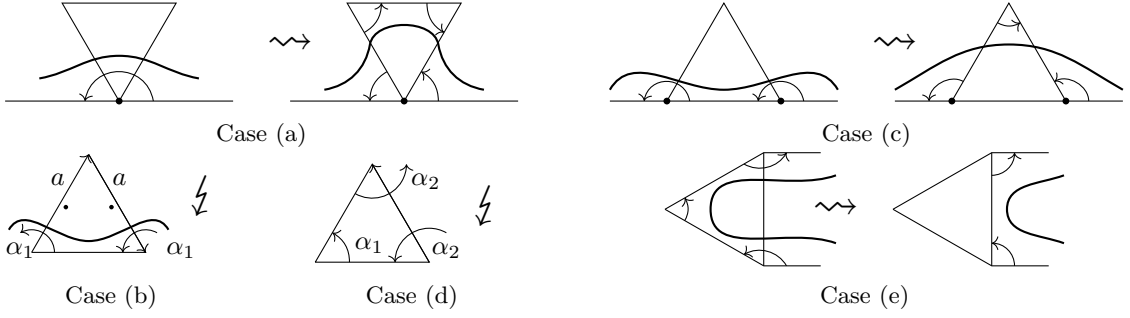

\begin{proof}
We proceed as in the proof of \autoref{th:strings-bands-cyclicity}, and go a little further. Let us repeat the steps: Choose an initial approximation of $ X $ by arcs $ a_i $ such all angles $ α_i $ are non-identities and strictly shorter than a full turn. The rest of the proof is concerned with tweaking the approximation so as to make $ δ $ upper triangular. In case of a string object, we are done.
Let us inspect the given arc collection $ α_1, …, α_k $ with its connecting angles $ α_1, …, α_k $. By assumption, all angles $ α_i $ go from arc $ a_i $ to $ a_{i+1} $, or all the other way around. Let us assume the former is the case: that $ α_i $ runs from $ a_i $ to $ a_{i+1} $. Moreover, any pair of consecutive arcs $ α_i, α_{i+1} $ runs at opposite ends of the arc $ a_{i+1} $, and all are non-identities.

Our strategy to break cyclicity is to pull some consecutive arcs with their angles following the interior of a polygon to the opposite side of the polygon. Let us make this concrete and distinguish the following hierarchy of cases: (a) there is an angle $ α_i $ with at least three indecomposable components, (b) $ k = 1 $, (c) there are two consecutive decomposable angles, (d) $ k = 2 $, (e) one angle $ α_i $ is indecomposable. By this hierarchy of cases, we mean that case (b) shall include that (a) does not hold; (c) shall include that (a) and (b) do not hold, etc. Samples for all cases are depicted in \autoref{fig:uncurving-representation-cyclicity}.

Regard case (a). Then we can flip a part of $ α_i $ to the other side of a triangle.

Regard case (b). If $ α_1 $ has just one indecomposable component, we have an immediate contradiction. If $ α_1 $ has two indecomposable components, then the arc in the middle of $ α_1 $ appears twice in the triangle, with equal orientation. This is also a contradiction, since the triangle is then not embedded anymore: In \autoref{fig:uncurving-representation-cyclicity}, the two dots inside the triangle would need to be equal, rendering the triangle non-embedded.

Regard case (c). Regard two consecutive decomposable angles $ α_i $ and $ α_{i+1} $. Then the last indecomposable part of $ α_i $ and the first indecomposable part of $ α_{i+1} $ are interior angles of a triangle. We can now flip these parts of the angles to the opposite side of the triangle. As depicted in the figure, this suffices to break cyclicity.

Regard case (d). It is impossible that both $ α_1 $ and $ α_2 $ are indecomposable. Indeed, this would mean that the curve partially winds around the interior of a triangle. Since we are excluding case (c), we can assume that $ α_1 $ is indecomposable and $ α_2 $ consists of two indecomposable parts. Then $ α_2 $ crosses both the head and the tail side of some arc, which is impossible.

Regard case (e). Let $ α_i $ be the indecomposable angle and regard $ α_{i-1} $ as well as $ α_{i+1} $. Since the curve does not bound a teardrop, the angle $ α_{i-1} $ enters the polygon from outside and $ α_{i+1} $ leaves the triangle. In other words, both are longer than the preceding and succeeding interior angles of the triangle. The tweak we apply to the angle sequence is to cut away the parts of $ α_{i-1} $ and $ α_{i+1} $ lying inside the triangle and deleting $ α_i $. A shorter sequence remains to be dealt with.

In every step, the angles that are already present become only shorter. Moreover, all angles that are inserted new in a step are interior angles of a polygon by choice. Since all polygons in $ \cA $ are triangles and all arcs in $ \cA $ are non-contractible in $ S $, all interior angles are shorter than a full turn. In total, we end up with an approximation where all angles are non-identities and shorter than full turns, as well as $ δ $ being upper triangular.
\end{proof}

\begin{theorem}
\label{th:uncurving-th}
Let $ \cA $ be a full arc system with [NMDC]. Let $ r ∈ \mathfrak{m} Z(\Gtl \cA) $ be a deformation parameter without $ 1 $-component. Then all $ ℤ/2ℤ $-graded band objects whose underlying curves in $ S $ are topologically nontrivial and do not bound a teardrop are uncurvable.
\end{theorem}

\begin{proof}
The proof consists of two steps. The first is to observe that we have already proven the case when all polygons in $ \cA $ are triangles. The second is to extend to the general case.

For the first step, let us assume that all polygons in $ \cA $ are triangles. Then by \autoref{th:uncurving-representation-cyclicity}, the band object $ X $ has a twisted complex representation where all angles are non-identities and shorter than full turns. By \autoref{th:uncurving-trick-works}, the complementary angle trick now successfully removes the curvature of $ X $.

For the second step, cut all polygons that are not triangles yet into pieces by adding arcs. Let us denote the resulting marked surface with arc system by $ \cA' $. We have an embedding $ \Gtl \cA ⊂ \Gtl \cA' $, and correspondingly
\begin{equation*}
i: \Tw\Gtl \cA ⊂ \Tw\Gtl \cA'.
\end{equation*}
It is well-known that this map is actually a quasi-equivalence, see \cite{Bocklandt}. The reason is that all additional arcs of $ \cA' $ can be built up to quasi-isomorphism as twisted complexes of arcs in $ \cA $. Correspondingly, we also have a quasi-isomorphism
\begin{equation*}
π: \Tw\Gtl \cA' → \Tw\Gtl \cA
\end{equation*}
that sends an object $ i(X) $ to some object $ π(i(X)) $ quasi-isomorphic to $ X $. Now choose $ X $ to be a band object as in the hypothesis, i.e.~topologically nontrivial in $ S $ and not bounding a teardrop. Then $ i(X) ∈ \Tw\Gtl \cA' $ is uncurvable by the first step of the proof. According to \autoref{rem:2Bprelim-ainfty-defo-uncurvingqi}, also $ π(i(X)) $ is uncurvable. Since $ X $ itself is quasi-isomorphic to $ π(i(X)) $, it is uncurvable as well by \autoref{rem:2Bprelim-ainfty-defo-uncurvingqi}. This finishes the proof.
\end{proof}

\begin{remark}
Let us get back to the 1-punctured torus with two arcs of \autoref{rem:uncurving-representation-norep}. We have seen that the horizontal band has no twisted complex representation that is suitable for the uncurving trick. According to \autoref{th:uncurving-th}, it can be uncurved nevertheless and by \autoref{rem:2Bprelim-ainfty-defo-uncurvingqi} this must in fact be possible for every chosen twisted complex representation.

Denote the angles in the 1-punctured torus by $ α $, $ β $, $ γ $, $ δ $ as in \autoref{fig:uncurving-trick-norep}. Pick the twisted complex representation
\begin{equation*}
X = \left(b ⊕ a ⊕ b[1], \pmat{0 & γ & \id \\ 0 & 0 & δ \\ 0 & 0 & 0}\right).
\end{equation*}
By experimenting, we have found the uncurved twisted complex
\begin{equation*}
X_q = \left(b ⊕ a ⊕ b[1], \pmat{0 & γ & \id \\ - q βαδ & q βα & δ \\ - q αδγβ & q γβα & 0}\right).
\end{equation*}
In other words, apart from the expected complementary angles $ q βαδ $, $ q γβα $ and $ q αδγβ $, we also have to insert $ q βα $. This difficulty is the reason we restricted the complementary angle trick to the case where all polygons are triangles.
\end{remark}

\begin{remark}
There are bands objects which are uncurvable but do not fall under the requirements of \autoref{th:uncurving-th}. Deriving finer criteria is however increasingly difficult. For instance, uncurvability of band objects representing curves with a teardrop depends on whether the deformation parameter includes $ ℓ^s $ for low $ s $.

For string objects, the situation is complicated as well. Most string objects cannot be uncurved, but there are exceptions. The underlying curve of a string object can typically be interpreted as an arc candidate for some full arc system if it has no self-intersections. In this case, it is not uncurvable for general deformation parameter $ r $. A string object whose underlying curve is a loop may however be uncurvable if $ r ∈ \mathfrak{m}^2 Z(\Gtl \cA) $.
\end{remark}

\section{The category of zigzag paths}
\label{sec:splitting}
In this section, we define the category $ \ZigzagCat $ and construct for it an explicit homological splitting. This category $ \ZigzagCat $ is a new, discrete analog of the smooth zigzag category studied by Cho, Hong and Lau \cite{CHL}. We follow their idea of including all zigzag curves with chosen spin structure into a category, except that we realize the curves as twisted complexes over the gentle algebra instead of objects in the Fukaya category.

After defining this category $ \ZigzagCat $, the second step in this section is to analyze the morphisms between the objects. We introduce terminology to handle locations where two zigzag paths come close to each other: situations of type A, B, C and D. We show that every morphism between two zigzag paths can be written as a linear combination of angles, and that each angle comes from a unique A, B, C or D type situation.

Finally, we use this classification of morphisms to provide an explicit homological splitting for the category $ \ZigzagCat $. This entails identifying a cohomology space $ H $ for every hom space, and finding complementary spaces $ R $ for every pair of zigzag paths. The splitting of individual morphisms into $ H $, $ I $ and $ R $ parts is collected in \autoref{tab:coh_splitting-splitting-verification}, which will be referred to throughout the paper.

\subsection{Category of zigzag paths}
\label{sec:prelim-zigzagcat}
In this section, we define the category $ \ZigzagCat $ of zigzag paths. The idea is to turn zigzag paths in $ Q $ into twisted complexes in $ \Tw\Gtl Q $. The construction of the twisted complexes requires the additional input data of a spin structure for every zigzag path. We fix terminology and notation for these spin structures. The result of the construction is the subcategory $ \ZigzagCat ⊂ \Tw\Gtl Q $.

\begin{center}
\begin{tikzpicture}
\path (0, 0) node (A) {Zigzag paths in $ Q $} (8, 0) node[align=center] (B) {Category of zigzag paths $ \ZigzagCat ⊂ \Tw\Gtl Q $};
\path[draw, decorate, decoration={snake, amplitude=0.2em, post length=0.5em}, ->] ($ (A.east)!0.2!(B.west) $) to ($ (A.east)!0.8!(B.west) $);
\end{tikzpicture}
\end{center}

Generally, arcs as objects in $ \Gtl Q $ can be stitched together to form twisted complexes which model curves in the Fukaya category. This is a well-known method, explicit in \cite[Section 9.2]{Bocklandt-book} and implicit in \cite{HKK}. We have detailed it in \autoref{sec:uncurving-stringsbands} and shall now apply it specifically to zigzag paths. Recall that a zigzag path has necessarily even length, because it alternates between turning left and right.

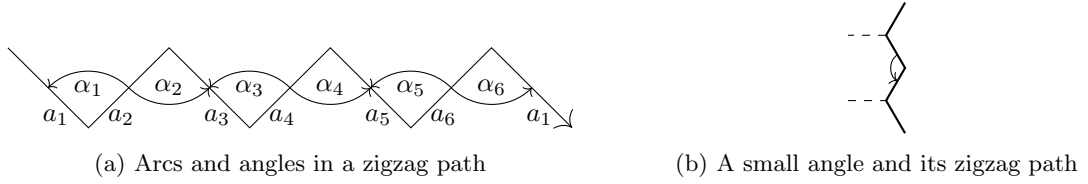
\begin{figure}
\centering
\begin{subfigure}{0.55\linewidth}
\centering
\begin{tikzpicture}
\path[draw, -{To[scale=2]}] (0, 0) -- ++(315:1.5) node[very near end, left] {$ a_1 $} coordinate[midway] (alpha1-end) -- ++(45:1.5) node[very near start, right] {$ a_2 $} coordinate[midway] (alpha1-start) -- ++(315:1.5) node[very near end, left] {$ a_3 $} coordinate[midway] (alpha3-end) -- ++(45:1.5) node[very near start, right] {$ a_4 $} coordinate[midway] (alpha3-start) -- ++(315:1.5) node[very near end, left] {$ a_5 $} coordinate[midway] (alpha5-end) -- ++(45:1.5) node[very near start, right] {$ a_6 $} coordinate[midway] (alpha5-start) -- ++(315:1.5) node[very near end, left] {$ a_1 $} coordinate[midway] (additional);
\path[draw, ->, bend right=45] (alpha1-start) to node[midway, below] {$ α_1 $} (alpha1-end);
\path[draw, ->, bend right=45] (alpha1-start) to node[midway, above] {$ α_2 $} (alpha3-end);
\path[draw, ->, bend right=45] (alpha3-start) to node[midway, below] {$ α_3 $} (alpha3-end);
\path[draw, ->, bend right=45] (alpha3-start) to node[midway, above] {$ α_4 $} (alpha5-end);
\path[draw, ->, bend right=45] (alpha5-start) to node[midway, below] {$ α_5 $} (alpha5-end);
\path[draw, ->, bend right=45] (alpha5-start) to node[midway, above] {$ α_6 $} (additional);
\end{tikzpicture}
\caption{Arcs and angles in a zigzag path}
\label{fig:zigzag-path}
\end{subfigure}
\begin{subfigure}{0.4\linewidth}
\centering
\begin{tikzpicture}[scale=0.5]
\path[draw, dashed] (0, 0) -- (1, 0);
\path[draw, thick] ($ (1, 0) + (300:1) $) -- (1, 0) -- ++(60:1) coordinate[pos=0.6] (1-end) -- ++(120:1) coordinate[pos=0.4] (1-start) coordinate (A) -- ++(60:1);
\path[draw, dashed] (A) -- ++(left:1);
\path[draw, ->, bend right=60] (1-start) to (1-end);
\end{tikzpicture}
\caption{A small angle and its zigzag path}
\label{fig:prelim-zigzagcat-small-determines}
\end{subfigure}
\caption{Illustration of zigzag paths and small angles}
\end{figure}

Given a zigzag path $ L $, we sometimes need to refer to specific angles surrounding $ L $. We set up this terminology as follows: Assume the consecutive arcs of the zigzag path are $ a_1, a_2, … $. Let us assume $ Q $ is geometrically consistent, so that every puncture in $ Q $ has valency at least four. Thus, the mere data of the arc sequence $ a_1, a_2, … $ already determines for every $ i ∈ ℕ $ whether $ L $ turns left or right at (the head or tail of) $ a_i $. For those $ i ∈ ℕ $ where $ L $ turns left at $ a_i $, denote by $ α_i $ the angle that winds around the common puncture $ h(a_i) = t(a_{i+1}) $, is shorter than a full turn, starts at the head of $ a_i $ and ends at the tail of $ a_{i+1} $. For those $ i ∈ ℕ $ where $ L $ turns right at $ a_i $, denote by $ α_i $ the angle that winds around the common puncture $ h(a_i) = t(a_{i+1}) $, is shorter than a full turn, starts at the tail of $ a_{i+1} $ and ending at the head of $ a_i $. The notation is depicted in \autoref{fig:zigzag-path}. We also call $ α_i $ the \emph{small angle} in $ L $ between $ a_i $ and $ a_{i+1} $.

\begin{remark}
\label{th:zigzagcat-arc-distinction}
Every (indexed) arc $ a_i $ either has precisely two small angles leaving it and no small angle ending at it, or two small angles ending at it and no small angles leaving it. In fact, the arcs of those types alternate along $ L $. The reader can easily convince himself of this fact by regarding \autoref{fig:zigzag-path}.
\end{remark}

Turning a zigzag path into a twisted complex requires the datum of a spin structure. For our purposes, this simply entails choosing a sign $ (-1)^{\#α_i} $ for every small angle $ α_i $ in the zigzag path. Writing the sign additively, we fix the terminology as follows:

\begin{definition}
A \emph{spin structure} on a zigzag path $ L $ is a choice of signs $ \# α_i ∈ ℤ/2ℤ $ for each of its small angles $ α_i $.
\end{definition}

The notation $ \# α_i $ makes sense: If the small angle between $ a_i $ and $ a_{i+1} $ is equal to the small angle between $ a_j $ and $ a_{j+1} $, then $ i $ and $ j $ differ precisely by a period of $ L $. In short, every small angle of $ L $ appears only once as small angle of $ L $ and therefore the notation $ \# α_i $ makes sense for the scope of a single zigzag path.

Giving a spin structure simultaneously for all zigzag paths of $ Q $ is equivalent to giving a sign $ \# α ∈ ℤ/2ℤ $ for all indecomposable angles $ α $ in $ Q $. Indeed, if two zigzag paths $ L_1 $ and $ L_2 $ share a small angle, then $ L_1 $ and $ L_2 $ are actually equal. This fact is depicted in \autoref{fig:prelim-zigzagcat-small-determines}. In that figure, the dashed arcs indicate that the drawn angle shall be an interior angle of a polygon. In summary, giving a collection of spin structures for all zigzag paths is equivalent to choosing a sign $ \# α ∈ ℤ/2ℤ $ for every indecomposable angle $ α $ in $ Q $.

\begin{definition}
\label{def:zigzag-complex}
Let $ Q $ be a geometrically consistent dimer and $ L $ a zigzag path with spin structure. Write the arcs of $ L $ as $ a_1, …, a_{2k} $, chosen such that $ L $ turns left from $ a_1 $ to $ a_2 $. Let $ α_1, α_2, … $ be the small angles in $ L $ between $ a_1 $ and $ a_2 $, etc. Then the twisted complex associated with $ L $ is given by
\begin{equation*}
L = (a_1 \oplus a_3 \oplus … \oplus a_k \oplus a_2 \oplus … \oplus a_{2k}, δ)
\end{equation*}
with twisting differential
\begin{equation*}
δ = \left[
\begin{array}{c|c}
0 & \begin{array}{ccccc}
(-1)^{\#α_1} α_1 & 0 & … & 0 & (-1)^{\#α_{2k}} α_{2k} \\
(-1)^{\#α_2} α_2 & (-1)^{\#α_3} α_3 & … & 0 & 0 \\
0 & (-1)^{\#α_4} α_4 & … & 0 & 0 \\
… & … & … & … & … \\
0 & 0 & … & (-1)^{\#α_{2k-3}} α_{2k-3} & 0 \\
0 & 0 & … & (-1)^{\#α_{2k-2}} α_{2k-2} & (-1)^{\#α_{2k-1}} α_{2k-1}
\end{array} \\\hline
0 & 0
\end{array}\right].
\end{equation*}
\end{definition}

In \autoref{def:zigzag-complex}, we have introduced abuse of notation twice: Using the letter “$ L $” both for a zigzag path and its twisted complex, and simply calling both uses a “zigzag path”. The intention behind this abuse is to switch seamlessly between between both uses. A typical sentence in this paper will be: “Regard the endomorphisms of some zigzag path $ L $.” In that sentence, it is clear that $ L $ shall be a zigzag path and we regard the endomorphisms of its associated twisted complex.

The twisted complex defined in \autoref{def:zigzag-complex} is indeed a well-defined object of $ \Tw\Gtl Q $, i.e.~$ δ $ satisfies the Maurer-Cartan equation:

\begin{lemma}
Let $ L $ be a zigzag path with spin structure. Then its twisting differential $ δ $ satisfies the Maurer-Cartan equation, so that $ L $ indeed lies in $ \Tw\Gtl Q $.
\end{lemma}

\begin{proof}
The Maurer-Cartan equation reads
\begin{equation*}
μ^1_{\Add\Gtl Q} (δ) + μ^2_{\Add\Gtl Q} (δ, δ) + μ^3_{\Add\Gtl Q} (δ, δ, δ) + … = 0.
\end{equation*}
Proving the Maurer-Cartan equation therefore boils down to showing that for any sequence of compatible angles $ α_1, …, α_k $ appearing in the $ δ $ matrix of $ L $ we have
\begin{equation*}
μ^k_{\Gtl Q} (α_k, …, α_1) = 0.
\end{equation*}
Let us check all such terms. Since the differential $ μ^1_{\Gtl Q} $ vanishes, we can assume $ k ≥ 2 $. Assume $ a_1 $ starts at arc $ a_1 $ and ends at $ a_2 $. In order to have any nonzero contribution $ μ^k_{\Gtl Q} (…, α_1) $, there would need to be a small angle starting at $ a_2 $. According to \autoref{th:zigzagcat-arc-distinction}, every arc however admits either only incoming or only outgoing small angles. Since $ α_1 $ is already an incoming angle for $ a_2 $, we conclude that $ a_2 $ has no outgoing angles. Therefore the product $ μ^k_{\Gtl Q} (…, α_1) $ vanishes. The reader can convince themself of this visually by drawing a zigzag path together with its small angles and trying to draw an immersed disk bounded solely by small angles. We conclude that $ δ $ satisfies the Maurer-Cartan equation.
\end{proof}

\begin{definition}
\label{def:splitting-zigzagcat-def}
Choose a spin structure for all zigzag paths $ L_1, …, L_N $ in $ Q $. Then the \emph{category of zigzag paths} $ \ZigzagCat ⊂ \Tw\Gtl Q $ is the (full) subcategory $ \ZigzagCat = \{L_1, …, L_N\} $ of $ \Tw\Gtl Q $ consisting of all zigzag paths with their single chosen spin structure.
\end{definition}

Every zigzag path only appears once in $ \ZigzagCat $. We do not allow the same zigzag path multiple times in $ \ZigzagCat $ with different spin structures, since this is not the goal of this paper and it would make calculations more complicated. Since $ \ZigzagCat $ depends on $ Q $ and the choice of spin structures, denoting this category by the letter $ \ZigzagCat $ denotes a slight abuse of notation.

\subsection{ABCD situations}
\label{sec:splitting-situations}
In this section, we provide a basis for the hom spaces between zigzag paths. We depart from two zigzag paths $ L_1 $ and $ L_2 $ and analyze their hom space in the category $ \ZigzagCat $. The first step is to introduce a notion of elementary morphisms from $ L_1 $ to $ L_2 $. This way, every morphism from $ L_1 $ to $ L_2 $ can be written as a linear combination of elementary morphisms. We then classify elementary morphisms according to the geometry of $ L_1 $ and $ L_2 $ in the surroundings of the morphisms. This gives rise to a classification into four types A, B, C, D. Elementary morphisms associated with these four types provide a basis for $ \Hom_{\ZigzagCat} (L_1, L_2) $. Ultimately, the A, B, C, D types will accompany us during our entire journey to the minimal model.

Regard two zigzag paths $ L_1 $ and $ L_2 $. It is our aim to provide a basis for $ \Hom(L_1, L_2) $. Recall from \autoref{sec:prelim-terminology} that an indexed arc of a zigzag path consists of an arc lying on the zigzag path, together with the datum of whether the zigzag path turns left or right at the head (equivalently tail) of the arc. Recall from \autoref{def:zigzag-complex} that every zigzag path comes with an associated twisted complex over $ \Gtl Q $. Therefore any morphism $ ε ∈ \Hom_{\ZigzagCat} (L_1, L_2) $ can be uniquely written as a linear combination of angles from indexed arcs of $ L_1 $ to indexed arcs of $ L_2 $. In other words, every angle from an indexed arc of $ L_1 $ to an indexed arc of $ L_2 $ gives rise to a morphism $ ε ∈ \Hom_{\ZigzagCat} (L_1, L_2) $, and all morphisms can be obtained as sums of such “elementary morphisms”. Let us make this precise:

\begin{definition}
An \emph{elementary morphism} $ ε: L_1 → L_2 $ is an angle from an indexed arc of $ L_1 $ to an indexed arc in $ L_2 $, interpreted as morphism between twisted complexes.
\end{definition}

We now move on to defining A, B, C and D situations. As a preparation, regard two consecutive arcs $ a, b $ in an elementary polygon of $ Q $. This pair defines a zigzag path $ L $, namely the one starting with $ a $ and turning then to $ b $. This way, $ a $ and $ b $ become indexed arcs of $ L $. Indeed, the zigzag path never traverses $ a $ followed by $ b $ again, except after cycling once fully through $ L $.

\begin{definition}
Let $ Q $ be a geometrically consistent dimer and $ L_1 $ and $ L_2 $ be zigzag paths.

An \emph{A situation} from $ L_1 $ to $ L_2 $ consists of two consecutive indexed arcs of $ L_1 $ and two consecutive indexed arcs of $ L_2 $ such that both midpoints (head of the first arc, equally tail of the second arc) are equal and all four arcs provide distinct incidences at this point. As in \autoref{fig:coh_splitting-situation-A}, the two arcs of $ L_1 $ are cyclically denoted 1, 2. The two arcs of $ L_2 $ are cyclically denoted 3, 4. The angles involved are denoted $ α $, $ β $, $ γ $ and $ β' $. An \emph{elementary angle of this A situation} is a product of such angles running from 1/2 to 3/4. Concretely, these are the angles in the A section of \autoref{tab:coh_splitting-splitting-verification}.

A \emph{B (resp.~C) situation} from $ L_1 $ to $ L_2 $ consists of an indexed arc 2 of $ L_1 $ and an indexed arc 5 of $ L_2 $ such that $ 2=5 $ as arcs in $ Q $, $ L_1 $ turns left (resp.~right) at the head and tail of 2, and $ L_2 $ turns right (resp.~left) at the head and tail of 2. As in \autoref{fig:coh_splitting-situation-B} (resp.~\ref{fig:coh_splitting-situation-C}), the neighboring indexed arcs are denoted 1, 3, 4, 6. An \emph{elementary angle of this B (resp.~C) situation} is a composition of angles in the figure, including $ \id_{2→5} $, that runs from 1/2/3 to 4/5/6. Concretely, these are the angles in the B (resp.~C) section of \autoref{tab:coh_splitting-splitting-verification}.

If $ L_1 = L_2 $, then a \emph{D situation} from $ L_1 = L_2 $ to itself consists of a single indexed arc of $ L_1 = L_2 $. If $ L_1 = L_2 $ turns right at the head of the arc, this arc is denoted 1 and the next indexed arc is denoted 2. If $ L_1 = L_2 $ turns left, the arc is denoted 2 and the next indexed arc is denoted 1. The angles are named as in \autoref{fig:coh_splitting-situation-D}. An \emph{elementary angle of this D situation} is a composition of angles in the figure that runs from the first arc to itself or to the second, from the second to the first, or from the second to itself by at least one full turn.

An elementary morphism may be annotated with its type of situation to enhance clarity: $ β $ (A), $ α_3 $ (B), $ \id $ (C), $ \id $ (D), etc.
\end{definition}

\begin{figure}
\begin{subfigure}{0.2\linewidth}
\centering
\begin{tikzpicture}
\path[draw] (0, 0) -- ++(45:2) coordinate[pos=0.7] (alpha-end) node[left=1cm] {$ L_1 $} node[near start, below] {2} -- ++(135:2) coordinate[pos=0.3] (alpha-start) node[near end, above] {1};
\path[draw] (3, 0) -- ++(135:2) coordinate[pos=0.7] (gamma-start) node[right=1cm] {$ L_2 $} node[near start, below] {3} -- ++(45:2) coordinate[pos=0.3] (gamma-end) node[near end, above] {4};
\path[draw, ->, bend right=45] (alpha-start) to node[midway, right] {$ α $} (alpha-end);
\path[draw, ->, bend right=45] (gamma-start) to node[midway, left] {$ γ $} (gamma-end);
\path[draw, ->, bend right=45] (gamma-end) to node[midway, below] {$ β' $} (alpha-start);
\path[draw, ->, bend right=45] (alpha-end) to node[midway, above] {$ β $} (gamma-start);
\end{tikzpicture}
\caption{Situation A}
\label{fig:coh_splitting-situation-A}
\end{subfigure}
\hfill
\begin{subfigure}{0.24\linewidth}
\centering
\begin{tikzpicture}
\path[draw] (0, 0) -- ++(30:2) coordinate[midway] (beta-start) node[near start, below] {4} coordinate (stop-1);
\path[draw, ->] (stop-1) -- ++(up:1.5) coordinate[pos=0.3] (alpha4-end) coordinate[pos=0.7] (alpha2-start) node[midway, right] {5} coordinate[pos=1] (stop-2);
\path[draw, -{To[scale=2]}] (stop-2) -- ++(30:2) node[near end, above] {$ L_2 $} node[near end, below] {6} coordinate[pos=0.3] (alpha2-end);
\path[draw] (3.1, 0) -- ++(150:2) coordinate[midway] (beta-end) node[near start, below] {3} coordinate (stop-3);
\path[draw, ->] (stop-3) -- ++(up:1.5) coordinate[pos=0.3] (alpha3-start) coordinate[pos=0.7] (alpha1-end) node[midway, left] {2} coordinate[pos=1] (stop-4);
\path[draw, -{To[scale=2]}] (stop-4) -- ++(150:2) node[near end, above] {$ L_1 $} node[near end, below] {1} coordinate[pos=0.3] (alpha1-start);
\path[draw, ->, bend right=60] (beta-start) to node[midway, above] {$ β $} (beta-end);
\path[draw, ->, bend right=60] (beta-end) to node[midway, left] {$ α_4 $} (alpha4-end);
\path[draw, ->, bend right=60] (alpha3-start) to node[midway, right] {$ α_3 $} (beta-start);
\path[draw, ->, bend right=60] (alpha2-start) to node[midway, left] {$ α_2 $} (alpha2-end);
\path[draw, ->, bend right=60] (alpha1-start) to node[midway, right] {$ α_1 $} (alpha1-end);
\path[draw, ->, bend right=70] (alpha2-end) to node[midway, below] {$ β' $} (alpha1-start);
\end{tikzpicture}
\caption{Situation B}
\label{fig:coh_splitting-situation-B}
\end{subfigure}
\hfill
\begin{subfigure}{0.24\linewidth}
\centering
\begin{tikzpicture}
\path[draw] (0, 0) -- ++(30:2) coordinate[midway] (beta-start) node[near start, below] {3} coordinate (stop-1);
\path[draw, ->] (stop-1) -- ++(up:1.5) coordinate[pos=0.3] (alpha4-end) coordinate[pos=0.7] (alpha2-start) node[midway, right] {2} coordinate (stop-2);
\path[draw, -{To[scale=2]}] (stop-2) -- ++(30:2) node[near end, above] {$ L_1 $} node[near end, below] {1} coordinate[pos=0.3] (alpha2-end);
\path[draw] (3.1, 0) -- ++(150:2) coordinate[midway] (beta-end) node[near start, below] {4} coordinate (stop-3);
\path[draw, ->] (stop-3) -- ++(up:1.5) coordinate[pos=0.3] (alpha3-start) coordinate[pos=0.7] (alpha1-end) node[midway, left] {5} coordinate (stop-4);
\path[draw, -{To[scale=2]}] (stop-4) -- ++(150:2) node[near end, above] {$ L_2 $} node[near end, below] {6} coordinate[pos=0.3] (alpha1-start);
\path[draw, ->, bend right=60] (beta-start) to node[midway, above] {$ β $} (beta-end);
\path[draw, ->, bend right=60] (beta-end) to node[midway, left] {$ α_4 $} (alpha4-end);
\path[draw, ->, bend right=60] (alpha3-start) to node[midway, right] {$ α_3 $} (beta-start);
\path[draw, ->, bend right=60] (alpha2-start) to node[midway, left] {$ α_2 $} (alpha2-end);
\path[draw, ->, bend right=60] (alpha1-start) to node[midway, right] {$ α_1 $} (alpha1-end);
\path[draw, ->, bend right=70] (alpha2-end) to node[midway, below] {$ β' $} (alpha1-start);
\end{tikzpicture}
\caption{Situation C}
\label{fig:coh_splitting-situation-C}
\end{subfigure}
\hfill
\begin{subfigure}{0.25\linewidth}
\centering
\begin{tikzpicture}
\path[draw] (0, 0) -- (1, 2) node[midway, left] {1} coordinate[pos=0.7] (alpha-start) -- (2, 0) coordinate[pos=0.3] (alpha-end) node[midway, right] {2};
\path[draw, ->, bend right=30] (alpha-start) to node[midway, above] {$ α $} (alpha-end);
\path[draw, ->] (alpha-end) .. controls +(0.5, 0.4) and +(0.5, 0) .. (1, 2.5) node[below] {$ α' $};
\path[draw, ->] (1, 2.5) .. controls +(-0.5, 0) and +(-0.5, 0.4) .. (alpha-start);
\end{tikzpicture}
\caption{Situation D}
\label{fig:coh_splitting-situation-D}
\end{subfigure}
\caption{All elementary morphisms $ ε: L_1 → L_2 $ are contained in one of these situations.}
\label{fig:situation-ABCD}
\end{figure}
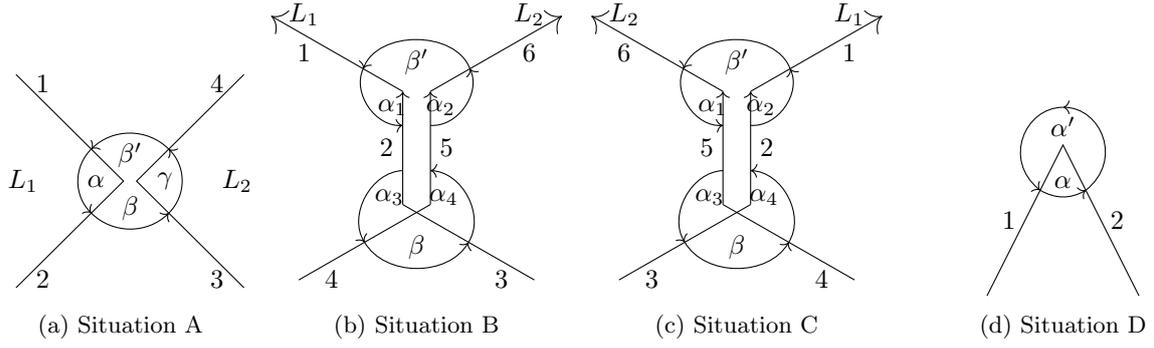
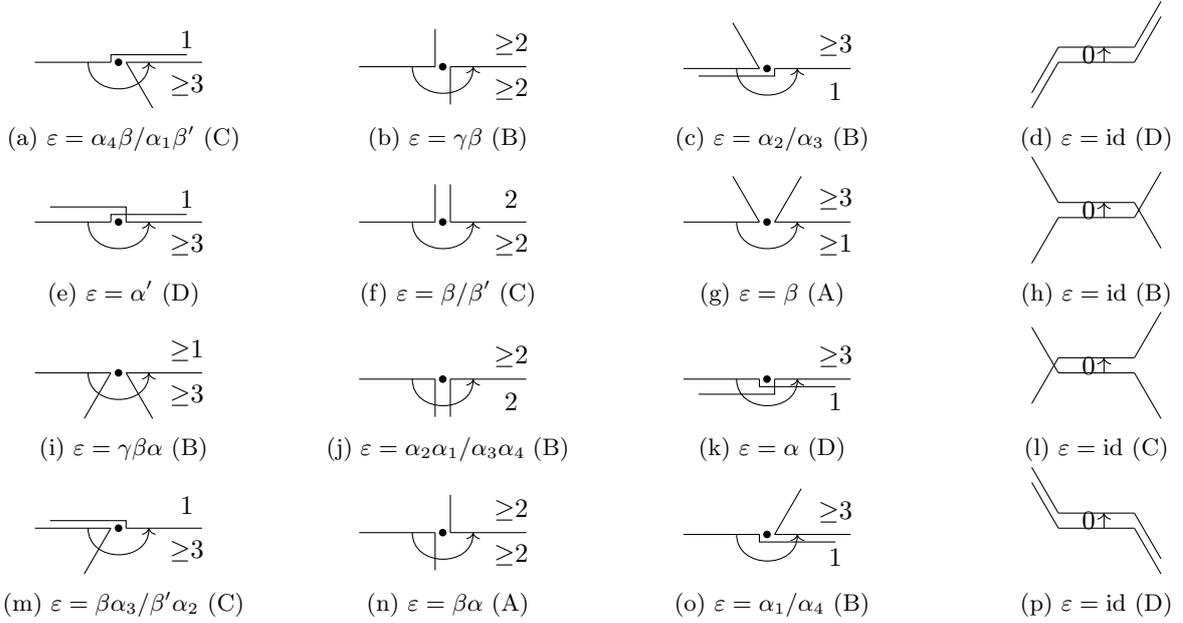
\begin{figure}
\begin{subfigure}{0.2\linewidth}
\centering
\begin{tikzpicture}
\path[draw] (0, 0) -- ++(right:1) coordinate[pos=0.7] (eps-start) coordinate (a) ++(right:0.1) coordinate (mid) ++(right:0.1) coordinate (b) -- ++(right:1) coordinate[pos=0.3] (eps-end);
\path[draw, ->, bend right=90, looseness=1.5] (eps-start) to (eps-end);
\path[fill] (mid) circle[radius=0.05];
\path[draw] (a) -- ++(up:0.1) -- ++(right:1);
\path[draw] (b) -- ++(300:0.7);
\path (2, 0.3) node {1};
\path (2, -0.3) node {≥3};
\end{tikzpicture}
\caption{$ ε = α_4 β / α_1 β' $ (C)}
\end{subfigure}
\hfill
\begin{subfigure}{0.2\linewidth}
\centering
\begin{tikzpicture}
\path[draw] (0, 0) -- ++(right:1) coordinate[pos=0.7] (eps-start) coordinate (a) ++(right:0.1) coordinate (mid) ++(right:0.1) coordinate (b) -- ++(right:1) coordinate[pos=0.3] (eps-end);
\path[draw, ->, bend right=90, looseness=1.5] (eps-start) to (eps-end);
\path[fill] (mid) circle[radius=0.05];
\path[draw] (a) -- ++(up:0.5);
\path[draw] (b) -- ++(down:0.5);
\path (2, 0.3) node {≥2};
\path (2, -0.3) node {≥2};
\end{tikzpicture}
\caption{$ ε = γβ $ (B)}
\end{subfigure}
\hfill
\begin{subfigure}{0.2\linewidth}
\centering
\begin{tikzpicture}
\path[draw] (0, 0) -- ++(right:1) coordinate[pos=0.7] (eps-start) coordinate (a) ++(right:0.1) coordinate (mid) ++(right:0.1) coordinate (b) -- ++(right:1) coordinate[pos=0.3] (eps-end);
\path[draw, ->, bend right=90, looseness=1.5] (eps-start) to (eps-end);
\path[fill] (mid) circle[radius=0.05];
\path[draw] (a) -- ++(120:0.7);
\path[draw] (b) -- ++(down:0.1) -- ++(left:1);
\path (2, 0.3) node {≥3};
\path (2, -0.3) node {1};
\end{tikzpicture}
\caption{$ ε = α_2 / α_3 $ (B)}
\end{subfigure}
\hfill
\begin{subfigure}{0.2\linewidth}
\centering
\begin{tikzpicture}
\path[draw] (0, 0) -- ++(right:1) coordinate[pos=0.6] (eps-start) coordinate (a) ++(up:0.2) coordinate (b) -- ++(left:1) coordinate[pos=0.4] (eps-end) coordinate (end);
\path[draw, ->] (eps-start) to (eps-end);
\path[draw] (0, 0) -- ++(240:0.7);
\path[draw] (end) -- ++(240:0.7);
\path[draw] (a) -- ++(60:0.7);
\path[draw] (b) -- ++(60:0.7);
\path (0.4, 0.1) node {0};
\end{tikzpicture}
\caption{$ ε = \id $ (D)}
\end{subfigure}
\hfill
\begin{subfigure}{0.2\linewidth}
\centering
\begin{tikzpicture}
\path[draw] (0, 0) -- ++(right:1) coordinate[pos=0.7] (eps-start) coordinate (a) ++(right:0.1) coordinate (mid) ++(right:0.1) coordinate (b) -- ++(right:1) coordinate[pos=0.3] (eps-end);
\path[draw, ->, bend right=90, looseness=1.5] (eps-start) to (eps-end);
\path[fill] (mid) circle[radius=0.05];
\path[draw] (a) -- ++(up:0.1) -- ++(right:1);
\path[draw] (b) -- ++(up:0.2) -- ++(left:1);
\path (2, 0.3) node {1};
\path (2, -0.3) node {≥3};
\end{tikzpicture}
\caption{$ ε = α' $ (D)}
\end{subfigure}
\hfill
\begin{subfigure}{0.2\linewidth}
\centering
\begin{tikzpicture}
\path[draw] (0, 0) -- ++(right:1) coordinate[pos=0.7] (eps-start) coordinate (a) ++(right:0.1) coordinate (mid) ++(right:0.1) coordinate (b) -- ++(right:1) coordinate[pos=0.3] (eps-end);
\path[draw, ->, bend right=90, looseness=1.5] (eps-start) to (eps-end);
\path[fill] (mid) circle[radius=0.05];
\path[draw] (a) -- ++(up:0.5);
\path[draw] (b) -- ++(up:0.5);
\path (2, 0.3) node {2};
\path (2, -0.3) node {≥2};
\end{tikzpicture}
\caption{$ ε = β/β' $ (C)}
\end{subfigure}
\hfill
\begin{subfigure}{0.2\linewidth}
\centering
\begin{tikzpicture}
\path[draw] (0, 0) -- ++(right:1) coordinate[pos=0.7] (eps-start) coordinate (a) ++(right:0.1) coordinate (mid) ++(right:0.1) coordinate (b) -- ++(right:1) coordinate[pos=0.3] (eps-end);
\path[draw, ->, bend right=90, looseness=1.5] (eps-start) to (eps-end);
\path[fill] (mid) circle[radius=0.05];
\path[draw] (a) -- ++(120:0.7);
\path[draw] (b) -- ++(60:0.7);
\path (2, 0.3) node {≥3};
\path (2, -0.3) node {≥1};
\end{tikzpicture}
\caption{$ ε = β $ (A)}
\end{subfigure}
\hfill
\begin{subfigure}{0.2\linewidth}
\centering
\begin{tikzpicture}
\path[draw] (0, 0) -- ++(right:1) coordinate[pos=0.6] (eps-start) coordinate (a) ++(up:0.2) coordinate (b) -- ++(left:1) coordinate[pos=0.4] (eps-end) coordinate (end);
\path[draw, ->] (eps-start) to (eps-end);
\path[draw] (0, 0) -- ++(240:0.7);
\path[draw] (end) -- ++(120:0.7);
\path[draw] (a) -- ++(60:0.7);
\path[draw] (b) -- ++(300:0.7);
\path (0.4, 0.1) node {0};
\end{tikzpicture}
\caption{$ ε = \id $ (B)}
\end{subfigure}
\hfill
\begin{subfigure}{0.2\linewidth}
\centering
\begin{tikzpicture}
\path[draw] (0, 0) -- ++(right:1) coordinate[pos=0.7] (eps-start) coordinate (a) ++(right:0.1) coordinate (mid) ++(right:0.1) coordinate (b) -- ++(right:1) coordinate[pos=0.3] (eps-end);
\path[draw, ->, bend right=90, looseness=1.5] (eps-start) to (eps-end);
\path[fill] (mid) circle[radius=0.05];
\path[draw] (a) -- ++(240:0.7);
\path[draw] (b) -- ++(300:0.7);
\path (2, 0.3) node {≥1};
\path (2, -0.3) node {≥3};
\end{tikzpicture}
\caption{$ ε = γβα $ (B)}
\end{subfigure}
\hfill
\begin{subfigure}{0.2\linewidth}
\centering
\begin{tikzpicture}
\path[draw] (0, 0) -- ++(right:1) coordinate[pos=0.7] (eps-start) coordinate (a) ++(right:0.1) coordinate (mid) ++(right:0.1) coordinate (b) -- ++(right:1) coordinate[pos=0.3] (eps-end);
\path[draw, ->, bend right=90, looseness=1.5] (eps-start) to (eps-end);
\path[fill] (mid) circle[radius=0.05];
\path[draw] (a) -- ++(down:0.5);
\path[draw] (b) -- ++(down:0.5);
\path (2, 0.3) node {≥2};
\path (2, -0.3) node {2};
\end{tikzpicture}
\caption{$ ε = α_2 α_1 / α_3 α_4 $ (B)}
\end{subfigure}
\hfill
\begin{subfigure}{0.2\linewidth}
\centering
\begin{tikzpicture}
\path[draw] (0, 0) -- ++(right:1) coordinate[pos=0.7] (eps-start) coordinate (a) ++(right:0.1) coordinate (mid) ++(right:0.1) coordinate (b) -- ++(right:1) coordinate[pos=0.3] (eps-end);
\path[draw, ->, bend right=90, looseness=1.5] (eps-start) to (eps-end);
\path[fill] (mid) circle[radius=0.05];
\path[draw] (a) -- ++(down:0.1) -- ++(right:1);
\path[draw] (b) -- ++(down:0.2) -- ++(left:1);
\path (2, 0.3) node {≥3};
\path (2, -0.3) node {1};
\end{tikzpicture}
\caption{$ ε = α $ (D)}
\end{subfigure}
\hfill
\begin{subfigure}{0.2\linewidth}
\centering
\begin{tikzpicture}
\path[draw] (0, 0) -- ++(right:1) coordinate[pos=0.6] (eps-start) coordinate (a) ++(up:0.2) coordinate (b) -- ++(left:1) coordinate[pos=0.4] (eps-end) coordinate (end);
\path[draw, ->] (eps-start) to (eps-end);
\path[draw] (0, 0) -- ++(120:0.7);
\path[draw] (end) -- ++(240:0.7);
\path[draw] (a) -- ++(300:0.7);
\path[draw] (b) -- ++(60:0.7);
\path (0.4, 0.1) node {0};
\end{tikzpicture}
\caption{$ ε = \id $ (C)}
\end{subfigure}
\hfill
\begin{subfigure}{0.2\linewidth}
\centering
\begin{tikzpicture}
\path[draw] (0, 0) -- ++(right:1) coordinate[pos=0.7] (eps-start) coordinate (a) ++(right:0.1) coordinate (mid) ++(right:0.1) coordinate (b) -- ++(right:1) coordinate[pos=0.3] (eps-end);
\path[draw, ->, bend right=90, looseness=1.5] (eps-start) to (eps-end);
\path[fill] (mid) circle[radius=0.05];
\path[draw] (a) -- ++(240:0.7);
\path[draw] (b) -- ++(up:0.1) -- ++(left:1);
\path (2, 0.3) node {1};
\path (2, -0.3) node {≥3};
\end{tikzpicture}
\caption{$ ε = β α_3 / β' α_2 $ (C)}
\end{subfigure}
\hfill
\begin{subfigure}{0.2\linewidth}
\centering
\begin{tikzpicture}
\path[draw] (0, 0) -- ++(right:1) coordinate[pos=0.7] (eps-start) coordinate (a) ++(right:0.1) coordinate (mid) ++(right:0.1) coordinate (b) -- ++(right:1) coordinate[pos=0.3] (eps-end);
\path[draw, ->, bend right=90, looseness=1.5] (eps-start) to (eps-end);
\path[fill] (mid) circle[radius=0.05];
\path[draw] (a) -- ++(down:0.5);
\path[draw] (b) -- ++(up:0.5);
\path (2, 0.3) node {≥2};
\path (2, -0.3) node {≥2};
\end{tikzpicture}
\caption{$ ε = βα $ (A)}
\end{subfigure}
\hfill
\begin{subfigure}{0.2\linewidth}
\centering
\begin{tikzpicture}
\path[draw] (0, 0) -- ++(right:1) coordinate[pos=0.7] (eps-start) coordinate (a) ++(right:0.1) coordinate (mid) ++(right:0.1) coordinate (b) -- ++(right:1) coordinate[pos=0.3] (eps-end);
\path[draw, ->, bend right=90, looseness=1.5] (eps-start) to (eps-end);
\path[fill] (mid) circle[radius=0.05];
\path[draw] (a) -- ++(down:0.1) -- ++(right:1);
\path[draw] (b) -- ++(60:0.7);
\path (2, 0.3) node {≥3};
\path (2, -0.3) node {1};
\end{tikzpicture}
\caption{$ ε = α_1 / α_4 $ (B)}
\end{subfigure}
\hfill
\begin{subfigure}{0.2\linewidth}
\centering
\begin{tikzpicture}
\path[draw] (0, 0) -- ++(right:1) coordinate[pos=0.6] (eps-start) coordinate (a) ++(up:0.2) coordinate (b) -- ++(left:1) coordinate[pos=0.4] (eps-end) coordinate (end);
\path[draw, ->] (eps-start) to (eps-end);
\path[draw] (0, 0) -- ++(120:0.7);
\path[draw] (end) -- ++(120:0.7);
\path[draw] (a) -- ++(300:0.7);
\path[draw] (b) -- ++(300:0.7);
\path (0.4, 0.1) node {0};
\end{tikzpicture}
\caption{$ ε = \id $ (D)}
\end{subfigure}
\caption{Examining case-by-case which situation $ ε $ belongs to. Here $ ε: L_1 → L_2 $ is an elementary morphism strictly shorter than a full turn. Each of the cases corresponds to the option whether $ L_1 $ and $ L_2 $ turns left or right at $ ε $ and how many indecomposable angles $ ε $ consists of. In the figures, the latter number is indicated below the target arc. The number above the target arc indicates the number of indecomposable angles in the complement of $ ε $. Due to zigzag consistency, the sum of the two numbers is at least 4.}
\label{fig:coh_splitting-situation-distinction}
\end{figure}

\begin{remark}
Any arcs involved in a situation are allowed to be equal. The distinction and ordering of arcs only concerns the local behavior of the arc ends at the common punctures.
\end{remark}

We show that the A, B, C and D situations classify elementary morphisms uniquely. One may think that, for example, an elementary morphism $ α_3 α_4 $ of a B situation equals an elementary morphism $ β' α_2 $ of a C situation, or even $ α_2 α_1 $ of the same B situation. We now explain that this is not the case.

\begin{proposition}
\label{th:coh_splitting-situations-exhaustive}
Let $ Q $ be a geometrically consistent dimer and let $ L_1, L_2 $ be zigzag paths in $ Q $. Then any elementary morphism $ ε: L_1 → L_2 $ is an elementary morphism of precisely one A, B, C or D situation. Moreover, it is assigned only once as an elementary morphism of that situation.
\end{proposition}

\begin{proof}
We show that any elementary morphism $ ε: L_1 → L_2 $ appears in an A, B, C or D situation. In case $ ε $ is strictly shorter than a full turn, a case-by-case study is performed in \autoref{fig:coh_splitting-situation-distinction}. Now if $ ε $ is a full turn or longer, write $ ε = ε' ℓ^n $, where $ ℓ^n $ denotes a number of full turns. Then $ ε' $ itself is shorter than a full turn and hence is an elementary morphism of an A, B, C or D situation. Then by definition also $ ε $ is an elementary morphism of the same situation. Note this also applies in case $ ε' $ is $ \id $ (B), $ \id $ (C) or $ \id $ (D), since the turns $ (α_1 β' α_2)^i $ (B), $ (α_4 β α_3)^i $ (B), $ (α_1 β' α_2)^i $ (C), $ (α_4 β α_3)^i $ (C), $ α α' $ (D) and $ α' α $ (D) are also elementary morphisms of B, C and D situations.

For uniqueness, realize that for any elementary morphism, as distinguished in \autoref{fig:coh_splitting-situation-distinction}, we can read off the entire situation and which of its elementary morphisms it concerns by inspecting which arc ends coincide, where the zigzag path turns at $ ε $ and in which directions the arrows point.
\end{proof}

Classifying elementary morphisms into A, B, C, D situations is extremely handy. During this paper, we often need to indicate generic morphisms of these four types. For instance, we may say that a certain morphism is an “$ α_3 $ morphism”. By this we mean that it is an $ α_3 $ morphism of one certain B situation.


\subsection{Homological splitting}
\label{sec:splitting-splitting}
In this section, we introduce a homological splitting for the category $ \ZigzagCat $ of zigzag paths. The starting point is the category $ \ZigzagCat $ together with the description of its hom spaces according to \autoref{sec:splitting-situations}. The first step in this section is to fix once and for all the requirements on $ Q $ and additional data we assume for the rest of the paper. The second step to define the splitting $ H ⊕ I ⊕ R $ by giving an explicit basis in terms of the A, B, C, D situations in $ Q $. The idea is to reflect the geometry of the associated zigzag curves as far as possible. It is worth memorizing some basis elements, for instance $ β $ (A) morphisms always belong to $ R $.

Let us fix once and for all the requirements on $ Q $ and the description of additional data we require. Apart from assuming that $ Q $ is a fixed geometrically consistent dimer with spin structures chosen for every of its zigzag paths, we also require the data of what we call idenity and co-identity locations. The idea is that the choice of homological splitting in this paper is not entirely canonical, but depends on these two kinds of choices. The datum of an identity location on a zigzag path $ L $ entails the choice of an indexed arc $ a_0 $ on $ L $. The midpoint of $ a_0 $ is to be thought of as location of the identity intersection point between the associated zigzag curve $ \smooth L $ and its Hamiltonian deformation. The datum of a co-identity location on a zigzag path $ L $ entails the choice of a small angle $ α_0 $ on $ L $. The midpoint of this angle is to be thought of as the location of the co-identity intersection point between the associated zigzag curve $ \smooth L $ and its Hamiltonian deformation. For the visual meaning of identity and co-identity locations, we refer to \autoref{fig:prelim-fukaya-cat-Ham} and \ref{fig:subdisk-coidentity}. With these use cases in mind, we can state our convention as follows:

\begin{convention}
\label{conv:alpha0-direction}
$ Q $ is a geometrically consistent dimer equipped with a choice of spin structure, identity location $ a_0 $ and co-identity location $ α_0 $ on every zigzag path. The category $ \ZigzagCat $ contains every zigzag path once, with the chosen spin structure. The co-identity $ α_0 $ shall be chosen to lie in a counterclockwise polygon.
\end{convention}

We assume \autoref{conv:alpha0-direction} throughout \autoref{sec:deformed}, \ref{sec:resultcomp} and \ref{sec:subdisk} without further notice. Only in the statement of the main theorems will we mention again that the convention is assumed. In \autoref{sec:classification}, the convention is assumed as well, while \autoref{sec:sphere} specifically deals with the sphere case where we abandon the convention. The requirement that $ α_0 $ shall lie in a counterclockwise polygon is required to make certain calculations work. For more details, we refer to the discussion in \autoref{sec:whydoes}.

We are now ready to define the homological splitting. For $ H $, the idea is to reflect the intersection geometry of the associated zigzag curves as far as possible. The basis elements for $ H $ defined in this section will be used throughout the paper and we refer to them as cohomology basis elements. For $ R $, we have to make slightly arbitrary choices of basis elements.

\begin{definition}
Let $ L_1 $ and $ L_2 $ be two zigzag paths. Denote by $ H ⊂ \Hom_{\ZigzagCat} (L_1, L_2) $ the subspace spanned by the \emph{cohomology basis morphisms}

\begin{itemize}
\setlength\itemsep{0.1cm}
\item for every B situation, $ (-1)^{\# α_3 + 1} α_3 + (-1)^{\# α_4} α_4 $,
\item for every C situation, $ \id_{2→5} $,
\item if $ L_1 = L_2 $, then $ (-1)^{\#α_0 + 1} α_0 $ and $ \sum_a \id_a $.
\end{itemize}

Denote by $ R ⊂ \Hom_{\ZigzagCat} (L_1, L_2) $ the subspace spanned by the following elements, collected from all situations from $ L_1 $ to $ L_2 $:

\begin{itemize}
\setlength\itemsep{0.2cm}
\item $ γβ (αβ'γβ)^i $ (A), $ β (αβ'γβ)^i $ (A),
\item $ (α_4 β α_3)^i, ~i≥1 $ (B), $ (α_1 β' α_2)^i, ~i≥1 $ (B), $ \id_{2→5} $ (B), $ α_3 (α_4 β α_3)^i $ (B), $ α_1 (β' α_2 α_1)^i $ (B),
\item $ β' (α_2 α_1 β')^i $ (C), $ β (α_3 α_4 β)^i $ (C), $ α_1 β' (α_2 α_1 β')^i $ (C), $ β α_3 (α_4 β α_3)^i $ (C),
\item $ \id_a, ~a≠a_0 $ (D), $ α' (αα')^i $ (D), $ (α'α)^i, ~i≥1 $ (D).
\end{itemize}

Denote by $ I $ the image of the twisted differential $ μ^1_{\ZigzagCat}: \Hom_{\ZigzagCat} (L_1, L_2) → \Hom_{\ZigzagCat} (L_1, L_2) $. Then $ R $, $ I $ and $ H $ form the \emph{(standard) splitting} of $ \Hom_{\ZigzagCat} (L_1, L_2) $.
\end{definition}

\begin{figure}
\centering
\begin{subfigure}{0.3\linewidth}
\centering
\begin{tikzpicture}
\path[draw, ->] (0, 0) -- ++(30:2) coordinate[pos=0.5] (target-tail) -- ++(up:1.5) coordinate[pos=0.5] (target-middle) coordinate (5-head);
\path[draw] (5-head) -- ++(30:2) coordinate[pos=0.5] (target-head) node[near end, below] {$ L_2 $};
\path[draw, ->] (3.1, 0) -- ++(150:2) coordinate[pos=0.5] (source-tail) -- ++(up:1.5) coordinate[pos=0.5] (source-middle) coordinate (2-head);
\path[draw] (2-head) -- ++(150:2) node[near end, below] {$ L_1 $} coordinate[pos=0.5] (source-head);
\path (source-middle) -- (target-middle) coordinate[pos=0.5] (middle);
\path[draw, thick, rounded corners, ->] ($ (middle)!2!(target-tail) $) -- (middle) coordinate[pos=0.7] (target2)-- ($ (middle)!2!(target-head) $) coordinate[pos=0.2] (target);
\path[draw, thick, rounded corners, ->] ($ (middle)!2!(source-tail) $) -- (middle) coordinate[pos=0.7] (source) -- ($ (middle)!2!(source-head) $) coordinate[pos=0.2] (source2);
\path[fill] (middle) circle(0.1);
\path[draw, ->, bend right=60] (source) to (target);
\path[draw, ->, bend right=60] (source2) to (target2);
\end{tikzpicture}
\caption{Situation B intersection}
\end{subfigure}
\begin{subfigure}{0.3\linewidth}
\centering
\begin{tikzpicture}
\path[draw, ->] (0, 0) -- ++(30:2) coordinate[pos=0.5] (target-tail) -- ++(up:1.5) coordinate[pos=0.5] (target-middle) coordinate (5-head);
\path[draw] (5-head) -- ++(30:2) coordinate[pos=0.5] (target-head) node[near end, below] {$ L_1 $};
\path[draw, ->] (3.1, 0) -- ++(150:2) coordinate[pos=0.5] (source-tail) -- ++(up:1.5) coordinate[pos=0.5] (source-middle) coordinate (2-head);
\path[draw] (2-head) -- ++(150:2) node[near end, below] {$ L_2 $} coordinate[pos=0.5] (source-head);
\path (source-middle) -- (target-middle) coordinate[pos=0.5] (middle);
\path[draw, thick, rounded corners, ->] ($ (middle)!2!(target-tail) $) -- (middle) coordinate[pos=0.3] (target2) -- ($ (middle)!2!(target-head) $) coordinate[pos=0.7] (target);
\path[draw, thick, rounded corners, ->] ($ (middle)!2!(source-tail) $) -- (middle) coordinate[pos=0.3] (source2) -- ($ (middle)!2!(source-head) $) coordinate[pos=0.7] (source);
\path[fill] (middle) circle(0.1);
\path[draw, ->, bend right=40] (target) to (source);
\path[draw, ->, bend right=40] (target2) to (source2);
\end{tikzpicture}
\caption{Situation C intersection}
\end{subfigure}
\caption{The cohomology elements $ (-1)^{\#α_3} α_3 + (-1)^{\#α_4 + 1} α_4 $ (B) and $ \id $ (C) are the odd resp.~even generators of Floer cohomology.}
\label{fig:splitting-splitting-BCinterpretation}
\end{figure}
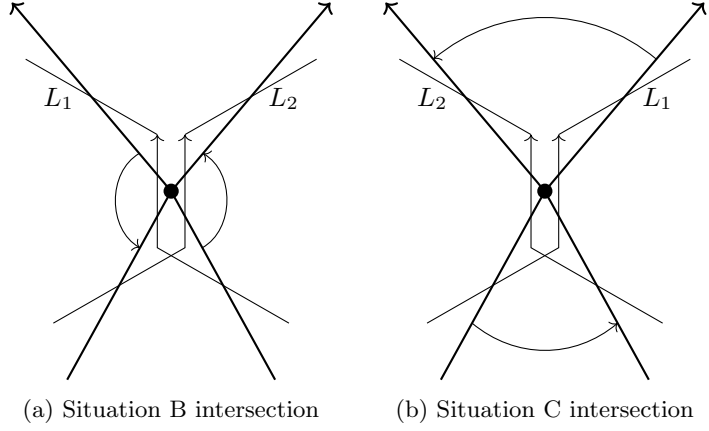

The cohomology basis elements are to be interpreted as follows: The odd element $ (-1)^{\#α_3 + 1} α_3 + (-1)^{\#α_4} α_4 $ corresponds to the odd intersection between $ \smooth L_1 $ and $ \smooth L_2 $ at the midpoint of the arc $ 2=5 $. In contrast, the even element $ \id_{2→5} $ corresponds to the even intersection between $ \smooth L_1 $ and $ \smooth L_2 $ at the midpoint of the arc $ 2=5 $. The odd element $ (-1)^{\#α_0 + 1} α_0 $ corresponds to the co-identity element of the zigzag curve $ \smooth L_1 $. The even element $ \sum_a \id_a $ corresponds to the identity element of the zigzag curve $ \smooth L_1 $. Two of the geometric interpretations are depicted in \autoref{fig:splitting-splitting-BCinterpretation}.

The splitting is not that hard to find and has at least been anticipated in \cite[Section 9.2]{Bocklandt-book}. For more information we refer to the discussion in \autoref{sec:literature}. Our next step is to show that the standard splitting indeed provides a homological splitting of $ \ZigzagCat $.

\begin{table}
\centering
\begin{tabular}{@{}llp{0.75\linewidth}@{}}
Situation A & 1→3 & $ βα (β'γβα)^i = μ^1_{\ZigzagCat} ((-1)^{\#α + 1} β (αβ'γβ)^i) + (-1)^{\#α + \#γ + ∥β∥} γβ(αβ'γβ)^i $ without triangle degeneration \\
& 1→3 & $ βα = μ^1_{\ZigzagCat} ((-1)^{\#α + 1} β) + (-1)^{\#α + \#γ + ∥β∥} γβ + (-1)^{\#α + \#γ_1 + \#γ_2} \id_{L_1 → L_2} $ in case of triangle degeneration \\
& 1→4 & $ γβα (β'γβα)^i = μ^1_{\ZigzagCat} ((-1)^{\#α + 1} γβ (αβ'γβ)^i) $ \\
& 2→3 & $ β (αβ'γβ)^i ∈ R $ \\
& 2→4 & $ γβ (αβ'γβ)^i ∈ R $ \\\hline
%
Situation B & 1→5 & $ α_1 (β' α_2 α_1)^i ∈ R $ \\
& 1→6 & $ α_2 α_1 (β' α_2 α_1)^i = μ^1_{\ZigzagCat} ((-1)^{\#α_2 + 1} α_1 (β' α_2 α_1)^i) $ \\
& 2→4 & $ α_3 (α_4 β α_3)^i ∈ R $ \\
& 2→5 & $ (α_1 β' α_2)^i ∈ R, ~ i ≥ 1 $ \\
& 2→5 & $ (α_4 β α_3)^i ∈ R, ~ i ≥ 1 $ \\
& 2→5 & $ \id_{2→5} ∈ R $ \\
& 2→6 & $ α_2 (α_1 β' α_2)^i = μ^1_{\ZigzagCat} ((-1)^{\#α_2} (α_1 β' α_2)^i) + (-1)^{\#α_1 + \#α_2} α_1 (β' α_2 α_1)^i, ~ i ≥ 1 $ \\
& 2→6 & $ α_2 = (-1)^{\#α_2} ((-1)^{\#α_3 + 1} α_3 + (-1)^{\#α_4} α_4) + μ^1_{\ZigzagCat} ((-1)^{\#α_2} \id_{2→5}) + (-1)^{\#α_1 + \#α_2} α_1 $ \\
& 3→4 & $ α_3 α_4 (β α_3 α_4)^i = μ^1_{\ZigzagCat} ((-1)^{\#α_4 + 1} α_3 (α_4 β α_3)^i) $ \\
& 3→5 & $ α_4 (β α_3 α_4)^i = μ^1_{\ZigzagCat} ((-1)^{\#α_4 + 1} (α_4 β α_3)^i) - α_3 (α_4 β α_3)^i, ~ i ≥ 1 $, \\
& 3→5 & $ α_4 = (-1)^{\#α_4} ((-1)^{\#α_3 + 1} α_3 + (-1)^{\#α_4} α_4) + (-1)^{\#α_3 + \#α_4} α_3 $ \\\hline
%
Situation C & 1→5 & $ α_1 β' (α_2 α_1 β')^i ∈ R $ \\
& 1→6 & $ β' (α_2 α_1 β')^i ∈ R $ \\
& 2→4 & $ β α_3 (α_4 β α_3)^i ∈ R $ \\
& 2→5 & $ (α_1 β' α_2)^i = μ^1_{\ZigzagCat} ((-1)^{\#α_2 + 1} α_1 β' (α_2 α_1 β')^{i-1}), ~ i ≥ 1$ \\
& 2→5 & $ (α_4 β α_3)^i = μ^1_{\ZigzagCat} ((-1)^{\#α_4 + 1} β α_3 (α_4 β α_3)^{i-1}), ~ i ≥ 1 $ \\
& 2→5 & $ \id_{2→5} ∈ H $, \\
& 2→6 & $ β' α_2 (α_1 β' α_2)^i = μ^1_{\ZigzagCat} ((-1)^{\#α_2 + 1} β' (α_2 α_1 β')^i) + (-1)^{\#α_1 + \#α_2} α_1 β' (α_2 α_1 β')^i $ \\
& 3→4 & $ β (α_3 α_4 β)^i ∈ R $ \\
& 3→5 & $ α_4 β (α_3 α_4 β)^i = μ^1_{\ZigzagCat} ((-1)^{\#α_4} β (α_3 α_4 β)^i) + (-1)^{\#α_3 + α_4} β α_3 (α_4 β α_3)^i $ \\\hline
%
Situation D & 1→1 & $ \id_{1→1} ∈ R $ if $ 1 ≠ a_0 $ \\
& 1→1 & $ \id_{1→1} = \sum \id_a - \sum_{a ≠ a_0} \id_a $ if $ 1 = a_0 $ \\
& 1→1 & $ (α' α)^i ∈ R $, $ i ≥ 1 $ \\
& 1→2 & $ α = α_0 ∈ H $ if $ α = α_0 $ \\
& 1→2 & $ α = (-1)^{\#α + \#α_0 + 1} α_0 + μ^1_{\ZigzagCat} ((-1)^{\#α + 1} (\id_{a_1} + … + \id_{a_k})) $ if $ α ≠ α_0 $ in case 1 with $ k $ odd \\
& 1→2 & $ α = (-1)^{\#α + \#α_0} α_0 + μ^1_{\ZigzagCat} ((-1)^{\#α + 1} (\id_{a_1} + … + \id_{a_k})) $ if $ α ≠ α_0 $ in case 1 with $ k $ even \\
& 1→2 & $ α = (-1)^{\#α + \#α_0 + 1} α_0 + μ^1_{\ZigzagCat} ((-1)^{\#α} (\id_{a_1} + … + \id_{a_k})) $ if $ α ≠ α_0 $ in case 2 with $ k $ odd \\
& 1→2 & $ α = (-1)^{\#α + \#α_0} α_0 + μ^1_{\ZigzagCat} ((-1)^{\#α + 1} (\id_{a_1} + … + \id_{a_k})) $ if $ α ≠ α_0 $ in case 2 with $ k $ odd \\
& 1→2 & $ α (α' α)^i = μ^1_{\ZigzagCat} ((-1)^{\#α} (α' α)^i) $, $ i ≥ 1 $ \\
& 2→1 & $ α' (α α')^i ∈ R $ \\
& 2→2 & $ \id_{2→2} ∈ R $ if $ 2 ≠ a_0 $ \\
& 2→2 & $ \id_{2→2} = \sum \id_a - \sum_{a ≠ a_0} \id_a $ if $ 2 = a_0 $ \\
& 2→2 & $ (α α')^i = μ^1_{\ZigzagCat} ((-1)^{\#α + 1} α' (α α')^{i-1}) - (α' α)^i $, $ i ≥ 1 $
\end{tabular}
\caption{Verification of the homological splitting}
\label{tab:coh_splitting-splitting-verification}
\end{table}

\begin{lemma}
\label{th:coh_splitting-splitting-exhaustive}
Let $ L_1 $ and $ L_2 $ be zigzag paths. Then $ \Hom_{\ZigzagCat} (L_1, L_2) = H + I + R $.
\end{lemma}

\begin{proof}
Any morphism $ L_1 → L_2 $ is a sum of elementary morphisms. By \autoref{th:coh_splitting-situations-exhaustive}, any elementary morphisms belongs to an A, B, C or D situation. Given this case distinction, \autoref{tab:coh_splitting-splitting-verification} shows that any such morphism can be written as a sum of elements in $ R $, $ I $ and $ H $.

It is worth commenting on the fact that the equations in \autoref{tab:coh_splitting-splitting-verification} actually hold true. This is due to diligent evaluation of the twisted differential
\begin{equation*}
μ^1_{\ZigzagCat} (ε) = μ^1 (ε) + μ^2 (δ, ε) + μ^2 (ε, δ) + ….
\end{equation*}
Let us examine the possible $ μ^{≥3} $ terms appear. For any elementary morphism $ ε: L_1 → L_2 $, at most one $ δ $ can be inserted upfront and at the back. Hence only $ μ^3 $ disks can appear. Inspecting the direction of $ δ $ morphisms, we conclude that the only possible disks appear in the case of $ μ^1_{\ZigzagCat} (β) $ of situation A and $ μ^1_{\ZigzagCat} (β) $ and $ μ^1_{\ZigzagCat} (β') $ of situation C.

First, let us examine the case of $ μ^1_{\ZigzagCat} (β) $ in situation A. Denote the corresponding $ δ $-angles by $ γ_1 $ and $ γ_2 $, and the next arcs of $ L_1 $ and $ L_2 $ by 5 and 6. Then traversing 6, $ γ_1 $, 2, $ β $, 3, $ γ_2 $, 5 must bound a disk and hence 5 and 6 are equal arcs and the disk is a simple polygon. Then indeed $ μ^1_{\ZigzagCat} (β) $ includes an $ \id_{5→6} $ morphism. See \autoref{fig:coh_splitting-triangle}.

Next, we rule out the possibility that $ μ^1_{\ZigzagCat} (β) $ or $ μ^1_{\ZigzagCat} (β') $ of situation C has a disk contribution. Carry out the same analysis as in situation A and find that $ β $ or $ β' $ is an elementary polygon angle, respectively. Together with the fact that the neighboring angles $ α_3 $ and $ α_4 $ are also elementary polygon angles comprising a full turn around puncture $ t(2) $, this contradicts the fact that $ Q $ is a dimer. We conclude that the differential $ μ^1_{\ZigzagCat} (β) $ in situation A remains as the only one that may include non-obvious terms.

Let us explain the meaning of case 1 and 2 in \autoref{tab:coh_splitting-splitting-verification}. This case distinction appears when we try to write an $ α $ angle with $ α ≠ α_0 $ in terms of $ H $, $ I $ and $ R $. Obviously, precisely one of the following two is the case:
\begin{itemize}
\item[1.] The segment of $ L_1 $ starting with the target 2 of $ α $ and continuing in the direction of $ α $ first hits the source or target of $ α_0 $ before hitting $ a_0 $.
\item[2.] The segment of $ L_1 $ starting with the source 1 of $ α $ and continuing in the opposite direction of $ α $ first hits the source or target of $ α_0 $ before hitting $ a_0 $.
\end{itemize}
The two cases are depicted in \autoref{fig:coh_splitting-identity-location}.

In case 1, put $ a_1 = 2 $, the indexed target of $ α $. Denote by $ a_1, …, a_k $ the segment of $ L $ until $ a_k $ is either the indexed source or target of $ α_0 $, whichever comes first. In case 2, put $ a_1 = 1 $, the indexed source of $ α $. Denote by $ a_1, …, a_k $ the segment of $ L $ until $ a_k $ is either the indexed source or target of $ α_0 $, whichever comes first.

In both cases, let $ α_i $ be the angle from $ a_{i-1} $ to $ a_i $ if $ i $ is odd, or from $ a_i $ to $ a_{i-1} $ if $ i $ is even. In particular, put $ α_1 = α $ and $ α_{k+1} = α_0 $.

In case 1, we have
\begin{alignat*}{2}
μ^1_{\ZigzagCat} (\id_{a_i}) &= (-1)^{\#α_i + 1} α_i + (-1)^{\#α_{i+1} + 1} α_{i+1}, &&\text{ if } i \text{ odd}, \\
μ^1_{\ZigzagCat} (\id_{a_i}) &= (-1)^{\#α_i} α_i + (-1)^{\#α_{i+1}} α_{i+1}, &&\text{ if } i \text{ even}.
\end{alignat*}
Adding these up, we obtain
\begin{alignat*}{2}
μ^1_{\ZigzagCat} (\id_{a_1} + … + \id_{a_k}) &= (-1)^{\#α + 1} α + (-1)^{\#α_0 + 1} α_0, &&\text{ if } k \text{ odd}, \\
μ^1_{\ZigzagCat} (\id_{a_1} + … + \id_{a_k}) &= (-1)^{\#α + 1} α + (-1)^{\#α_0} α_0, &&\text{ if } k \text{ even}.
\end{alignat*}
In case 2, we similarly get
\begin{alignat*}{2}
μ^1_{\ZigzagCat} (\id_{a_1} + … + \id_{a_k}) &= (-1)^{\#α} α + (-1)^{\#α_0} α_0, &&\text{ if } k \text{ odd}, \\
μ^1_{\ZigzagCat} (\id_{a_1} + … + \id_{a_k}) &= (-1)^{\#α} α + (-1)^{\#α_0 + 1} α_0, &&\text{ if } k \text{ even}.
\end{alignat*}
This precisely verifies the equations in \autoref{tab:coh_splitting-splitting-verification} concerning $ α $ (D).
\end{proof}

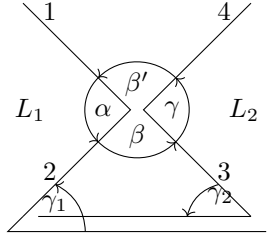
\begin{figure}
\centering
\begin{tikzpicture}
\path[draw] (0, 0) -- ++(45:2) coordinate[pos=0.3] (gamma1-end) coordinate[pos=0.7] (alpha-end) node[left=1cm] {$ L_1 $} node[near start, above] {2} -- ++(135:2) coordinate[pos=0.3] (alpha-start) node[near end, above] {1};
\path[draw] (3, 0) -- ++(135:2) coordinate[pos=0.3] (gamma2-start) coordinate[pos=0.7] (gamma-start) node[right=1cm] {$ L_2 $} node[near start, above] {3} -- ++(45:2) coordinate[pos=0.3] (gamma-end) node[near end, above] {4};
\path[draw, ->, bend right=45] (alpha-start) to node[midway, right] {$ α $} (alpha-end);
\path[draw, ->, bend right=45] (gamma-start) to node[midway, left] {$ γ $} (gamma-end);
\path[draw, ->, bend right=45] (gamma-end) to node[midway, below] {$ β' $} (alpha-start);
\path[draw, ->, bend right=45] (alpha-end) to node[midway, above] {$ β $} (gamma-start);
\path[draw] (0, 0) -- (-0.2, -0.2) -- (3.2, -0.2) coordinate[pos=0.3] (gamma1-start);
\path[draw] (3, 0) -- (0.2, 0) coordinate[pos=0.3] (gamma2-end);
\path[draw, ->, bend right=30] (gamma1-start) to node[left] {$ γ_1 $} (gamma1-end);
\path[draw, ->, bend right=30] (gamma2-start) to node[right] {$ γ_2 $} (gamma2-end);
\end{tikzpicture}
\caption{Triangle degeneration in situation A}
\label{fig:coh_splitting-triangle}
\end{figure}
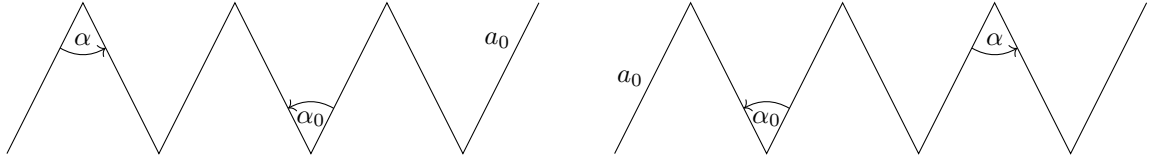
\begin{figure}
\begin{subfigure}{0.49\linewidth}
\centering
\begin{tikzpicture}
\path[draw] (0, 0) -- (1, 2) coordinate[pos=0.7] (alpha-start) -- (2, 0) coordinate[pos=0.3] (alpha-end) -- (3, 2) -- (4, 0) coordinate[pos=0.7] (alpha0-end) -- (5, 2) coordinate[pos=0.3] (alpha0-start) -- (6, 0) -- (7, 2) node[near end, left] {$ a_0 $};
\path[draw, ->, bend right=30] (alpha-start) to node[midway, above] {$ α $} (alpha-end);
\path[draw, ->, bend right=30] (alpha0-start) to node[midway, below] {$ α_0 $} (alpha0-end);
\end{tikzpicture}
\caption{Case 1: When following $ L $ in the direction of $ α $, the co-identity $ α_0 $ is hit before $ a_0 $.}
\end{subfigure}
\begin{subfigure}{0.49\linewidth}
\centering
\begin{tikzpicture}
\path[draw] (0, 0) -- (1, 2) node[midway, left] {$ a_0 $} -- (2, 0) coordinate[pos=0.7] (alpha0-end) -- (3, 2) coordinate[pos=0.3] (alpha0-start) -- (4, 0) -- (5, 2) coordinate[pos=0.7] (alpha-start) -- (6, 0) coordinate[pos=0.3] (alpha-end) -- (7, 2);
\path[draw, ->, bend right=30] (alpha-start) to node[midway, above] {$ α $} (alpha-end);
\path[draw, ->, bend right=30] (alpha0-start) to node[midway, below] {$ α_0 $} (alpha0-end);
\end{tikzpicture}
\caption{Case 2: When following $ L $ in the direction against $ α $, the co-identity $ α_0 $ is hit before $ a_0 $.}
\end{subfigure}
\caption{Position of an angle $ α $ between $ α_0 $ and $ a_0 $.}
\label{fig:coh_splitting-identity-location}
\end{figure}

\begin{proposition}
Let $ L_1 $ and $ L_2 $ be zigzag paths. Then the standard splitting $ H $, $ I $, $ R $ defines a homological splitting of $ \Hom_{\ZigzagCat} (L_1, L_2) $.
\end{proposition}

\begin{proof}
We need to check that $ I $ is the image of $ μ^1_{\ZigzagCat} $, that $ H + I $ is the kernel of $ μ^1_{\ZigzagCat} $ and that $ \Hom_{\ZigzagCat} (L_1, L_2) = H ⊕ I ⊕ R $. First, note that $ I $ is the image of $ μ^1_{\ZigzagCat} $ by definition. Next let us check that $ μ^1_{\ZigzagCat} (H) = 0 $. In our situational formalism, this is a simple calculation:
\begin{align*}
μ^1_{\ZigzagCat} ((-1)^{\#α_3 + 1} α_3 + (-1)^{\#α_4} α_4) &= μ^2 ((-1)^{\#α_3 + 1} α_3, (-1)^{\#α_4} α_4) + μ^2 ((-1)^{\#α_3} α_3, (-1)^{\#α_4} α_4) + μ^{≥3} (…) \\
&= (-1)^{\#α_3 + \#α_4} α_3 α_4 + (-1)^{\#α_3 + \#α_4 + 1} α_3 α_4 = 0, \\
μ^1_{\ZigzagCat} (\id \text{ (C)}) &= μ^{≥3} (…, \id, …) = 0, \\
μ^1_{\ZigzagCat} (\sum \id_a) &= μ^2 (\sum \id_a, δ) + μ^2 (δ, \sum \id_a) + μ^{≥3} (…, \id, …) = -δ + δ = 0, \\
μ^1_{\ZigzagCat} ((-1)^{\#α_0 + 1} α_0) &= (-1)^{\#α_0 + 1} μ^{≥3} (…, α_0, …) = 0.
\end{align*}
Indeed, both $ \id $ (C) and $ α_0 $ (D) have no $ δ $ morphisms that can be multiplied upfront or at the back, and correspondingly also produce no $ μ^{≥3} $ disks. Since $ I ⊂ \Ker(μ^1_{\ZigzagCat}) $, we conclude $ H + I ⊂ \Ker(μ^1_{\ZigzagCat}) $.

Next, recall from \autoref{th:coh_splitting-splitting-exhaustive} that we have $ \Hom_{\ZigzagCat} (L_1, L_2) = H + I + R $. Let us verify by dimension counting that this sum is direct. For $ N ∈ ℕ $, denote by $ E_N ⊂ E := \Hom_{\ZigzagCat} (L_1, L_2) $ the subspace spanned by elementary angles of at most $ N $ full turns. This is a finite dimensional space.

\autoref{tab:coh_splitting-splitting-verification} implies that any element of $ E_N $ can be written as a sum of elements of $ H ∩ E_N $, $ I ∩ E_N $ and $ R ∩ E_N $ if $ N ≥ 1 $. The sum of these three spaces is direct, since the sum of their dimensions is less than or equal to the dimension of $ E_N $:
\newlength{\dimHwidth}
\settowidth{\dimHwidth}{$ \dim (H ∩ E_N) $}
\begin{align*}
\underset{+}{\dim (H ∩ E_N)} &= \#B + \#C ~[+ 2 \text{ if } L_1 = L_2], \\
\underset{+}{\dim (I ∩ E_N)} &≤ 2N \# A + (4N - 1) \#B + 4N \#C ~[+ 2N \#D - 1 \text{ if } L_1 = L_2], \\
\underset{\rotatebox[origin=c]{270}{≤}}{\dim (R ∩ E_N)} &= 2N \# A + (4N + 1) \#B + 4N \#C ~[+ (2N + 1) \#D - 1 \text{ if } L_1 = L_2], \\
\makebox[\dimHwidth]{$ \dim E_N $} &= 4N \#A + (8N + 1) \#B + (8N + 1) \#C ~[+ (4N + 1) \#D \text{ if } L_1 = L_2].
\end{align*}
Here $ \#A $, $ \#B $, $ \#C $ and $ \#D $ denote the number of situations of type A, B, C and D that appear from $ L_1 $ to $ L_2 $.

The filtration of $ E $ by $ E_N $ is exhaustive. If the zero element of $ E $ is a nontrivial sum of elements of $ H $, $ I $ and $ R $, then by picking the maximum number of turns $ N $ involved, we obtain a contradiction to $ E_N $ being the direct sum of $ H ∩ E_N $, $ I ∩ E_N $ and $ R ∩ E_N $. We conclude that $ E = H ⊕ I ⊕ R $.

Finally, let us argue that $ \Ker(μ^1_{\ZigzagCat}) = H ⊕ I $. Indeed, the equivalent statement that $ μ^1_{\ZigzagCat}: R → I $ is injective can be checked by an elementary situational calculation. Note that if $ ε $ is from a certain situation, then $ μ^1_{\ZigzagCat} (ε) $ is from the same situation, apart from triangle degenerations and situation D morphisms.

Alternatively, note the cohomology of the complex $ \Hom_{\ZigzagCat} (L_1, L_2) $ is precisely Floer cohomology: The category $ \HTw\Gtl Q $ is nothing else than the wrapped Fukaya category of $ Q $. Now Floer cohomology has a basis element for every intersection of $ L_1 $ and $ L_2 $, plus an identity and a co-identity if $ L_1 = L_2 $. We conclude
\begin{equation*}
\dim (\Ker(μ^1_{\ZigzagCat})/I) = \#B + \#C ~[+ 2 \text{ if } L_1 = L_2] = \dim H.
\end{equation*}
This implies $ \Ker(μ^1_{\ZigzagCat}) = H ⊕ I $ and finishes the proof.
\end{proof}

\section{The deformed category of zigzag paths}
\label{sec:deformed}
In this section, we define and study the deformed category $ \DefZigzagCat $ of zigzag paths. This category is the deformed version of $ \ZigzagCat $. Its deformation comes from the deformation $ \Gtl_q Q $ of the gentle algebra. Already in the definition of $ \DefZigzagCat $, we apply the complementary angle trick from \autoref{sec:uncurving-gtl} in order to remove the curvature. We then analyze the differential $ μ^1_{\DefZigzagCat} $, in fact investigate how $ μ^1_{\DefZigzagCat} $ interacts with the homological splitting of $ \ZigzagCat $ from \autoref{sec:splitting-splitting}. We show that $ \DefZigzagCat $ falls under the regime of the simplified deformed Kadeishvili construction from \autoref{sec:2Bkadeishvili}. We provide explicit formulas for the deformed codifferential $ h_q $ and deformed projection $ π_q $ in terms of “tails” of the morphisms involved.

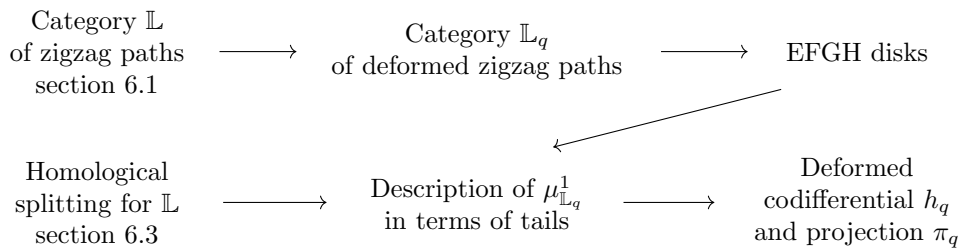
\begin{center}
\begin{tikzpicture}
\path (0, 0) node[align=center] (A) {Category $ \ZigzagCat $\\ of zigzag paths \\ \autoref{sec:prelim-zigzagcat}};
\path (0, -2) node[align=center] (B) {Homological \\ splitting for $ \ZigzagCat $ \\ \autoref{sec:splitting-splitting}};
\path (5, 0) node[align=center] (C) {Category $ \DefZigzagCat $ \\ of deformed zigzag paths};
\path (10, 0) node[align=center] (D) {EFGH disks};
\path (5, -2) node[align=center] (E) {Description of $ μ^1_{\DefZigzagCat} $ \\ in terms of tails};
\path (10, -2) node[align=center] (F) {Deformed \\ codifferential $ h_q $ \\ and projection $ π_q $};
%
\newcommand{\nodeconnect}[2]{\path[draw, ->] ($ (#1.east)!0.2!(#2.west) $) -- ($ (#1.east)!0.8!(#2.west) $);}
\nodeconnect AC
\nodeconnect CD
\nodeconnect BE
\nodeconnect EF
\path[draw, ->] ($ (D.south)!0.2!(E.north) $) -- ($ (D.south)!0.8!(E.north) $);
\end{tikzpicture}
\end{center}

These “tails” arise from applying the general Kadeishvili construction to the specific case of $ \DefZigzagCat $. According to the general deformed Kadeishvili theorem, we need to find for every cohomology basis element $ h ∈ H $ a certain deformed version $ φ^{-1} (h) $ such that $ φ^{-1} (h) $ and $ h $ only differ by $ R $ terms. In case of $ \DefZigzagCat $, we can explicitly describe $ φ^{-1} (h) $ for each of the cohomology basis elements in $ H $. The explicit description forces us to define and make us of what we call the tails of morphisms. Tails work as follows: Let $ ε: L_1 → L_2 $ be an elementary morphism. Look at all locations where $ L_1 $ and $ L_2 $ come close to each other and bound a discrete immersed disk together with $ ε $. It turns out that these locations of closeness have a hierarchical structure, which we organize in a tree. This tree is the tail of $ ε $ and carries the data of many secondary A, B, C or D situations which lie far away from $ ε $. We use the angles contained in this tail to construct explicitly the deformed cohomology basis elements $ φ^{-1} $. In other words, tails are the right tool to convert the rather algebraic requirement of the deformed Kadeishvili theorem into a geometric interpretation in the specific case of $ \DefZigzagCat $.

\begin{remark}
\label{rem:deformed-mushorthand}
From this section onwards, we typically write $ μ $ or $ μ_q $ for the product of $ \Add\Gtl_q Q $:
\begin{equation*}
μ_q ≔ μ ≔ μ_{\Add\Gtl_q Q}.
\end{equation*}
The reason for this notation is that we frequently expand products of the twisted completion $ \Tw\Gtl_q Q $ in terms of the products of the additive completion. This shorthand is meant to facilitate this expansion and mixes well with writing $ μ $ or $ μ_q $ for the product of $ \Gtl_q Q $. The product of $ \Gtl Q $ is irrelevant and is never meant. We keep the shorthand $ μ_q $ until \autoref{sec:subdisk}.
\end{remark}

\subsection{Deformed zigzag paths}
\label{sec:deformed-uncurving}
In this section, we define the category $ \DefZigzagCat $ of deformed zigzag paths. The starting point is the non-deformed category $ \ZigzagCat $ of zigzag paths. Taking the same twisted complexes gives a subcategory of $ \Tw\Gtl_q Q $. The aim of the entire paper is to compute a minimal model for this subcategory. According to the deformed Kadeishvili theorem, the first step is to gauge away as much curvature as possible. The aim of the present section is to conduct this uncurving procedure, and to define $ \DefZigzagCat $ to be the resulting uncurved category. The essential tool for uncurving is the complementary angle trick of \autoref{sec:uncurving-gtl}.

\begin{center}
\begin{tikzpicture}
\path (0, 0) node[align=center] (A) {Zigzag paths \\ $ \ZigzagCat ⊂ \Tw\Gtl Q $};
\path (0, -1.5) node[align=center] (D) {Deformation $ \Gtl_q Q $};
\path (5, -0.75) node[align=center] (B) {Zigzag paths \\ in $ \Tw\Gtl_q Q $};
\path (11, -0.75) node[align=center] (C) {Deformed zigzag paths \\ $ \DefZigzagCat ⊂ \Tw'\Gtl_q Q $};
\path[draw, ->] ($ (A.east)!0.2!(B.west) $) -- ($ (A.east)!0.8!(B.west) $);
\path[draw, ->] ($ (D.east)!0.2!(B.west) $) -- ($ (D.east)!0.8!(B.west) $);
\path[draw, ->] ($ (B.east)!0.2!(C.west) $) -- ($ (B.east)!0.8!(C.west) $) node[midway, above] {uncurving};
\end{tikzpicture}
\end{center}

Recall that the zigzag paths are objects in $ \Tw\Gtl Q $. They can also be interpreted as objects in $ \Tw\Gtl_q Q $ by definition of the deformed twisted completion. As an object of $ \Tw\Gtl_q Q $, every zigzag path has curvature. Our approach is to uncurve every zigzag path by means of the complementary angle trick of \autoref{sec:uncurving-gtl}. For every zigzag path, the trick gives rise to a twisted complex with also infinitesimal below-diagonal entries. By definition, this is an element of the category $ \Tw'\Gtl_q Q $ which we recalled in \autoref{rem:2B-ainfty-defo-liberaltwisted}. We shall call this object a deformed zigzag path because its $ δ $-matrix has been deformed. The precise definition for deformed zigzag paths reads as follows:

\begin{definition}
Let $ L $ be a zigzag path of $ Q $, with twisted complex
\begin{equation*}
L = \left(a_1 ⊕ a_3 ⊕ … ⊕ a_k ⊕ a_2 ⊕ … ⊕ a_{2k}, δ\right)
\end{equation*}
as in \autoref{def:zigzag-complex}. Then the corresponding \emph{deformed zigzag path} is the following object of $ \Tw'\Gtl_q $, also denoted $ L $:
\begin{align*}
L &= \left(a_1 ⊕ a_3 ⊕ … ⊕ a_k ⊕ a_2 ⊕ … ⊕ a_{2k}, δ\right), \\
δ &= \left[
\begin{array}{c|c}
0 & \text{ditto} \\\hline
\begin{array}{ccccc}
(-1)^{\# α_1} q_1 α_1' & (-1)^{\#α_2} q_2 α_2' & 0 & … & 0 \\
0 & (-1)^{\#α_3} q_3 α_3' & (-1)^{\#α_4} q_4 α_4' & … & 0 \\
… & … & … & … & … \\
0 & 0 & 0 & … & (-1)^{\#α_{2k-2}} q_{2k-2} α_{2k-2}' \\
(-1)^{\#α_{2k}} q_{2k} α_{2k}' & 0 & 0 & … & (-1)^{\#α_{2k-1}} q_{2k-1} α_{2k-1}'
\end{array} & 0
\end{array}\right].
\end{align*}
Here “ditto” denotes the same matrix entries as in \autoref{def:zigzag-complex}. The letter $ q_i $ denotes the puncture around which $ α_i $ winds.
\end{definition}

The deformed zigzag paths are objects of $ \Tw'\Gtl_q Q $. The entries of their $ δ $-matrix are angles $ α_i $ and $ α_i' $, which we call the \emph{inner} respectively \emph{outer $ δ $-angles} of $ L $. The main interest of the present paper is in the subcategory of all deformed zigzag paths, each with their associated choice of spin structure. We give the category consisting of these objects a name:

\begin{definition}
\label{def:deformed-zigzagcat-def}
The \emph{category of deformed zigzag paths} is the full subcategory $ \DefZigzagCat ⊂ \Tw'\Gtl_q Q $ consisting of the deformed zigzag paths.
\end{definition}

The category $ \DefZigzagCat $ is a deformation of $ \ZigzagCat $. This follows directly from the definition of the liberal twisted completion construction $ \Tw'\cat C_q $, see \papertwoA. It is completely acceptable that we have inserted below-diagonal entries into the twisted complex. In what follows, we explain why $ \DefZigzagCat $ is indeed curvature-free. It is in fact a consequence of the more general property of the complementary angle trick \autoref{th:uncurving-trick-works}, but we provide here a direct proof as well. The direct proof builds on the fact that segments of zigzag paths in geometrically consistent dimers cannot bound discrete immersed disks. In \autoref{fig:disks-zigzag-boundary} we have depicted a fictitious discrete immersed disk bounded by a zigzag path and we shall prove that this situation is indeed impossible:

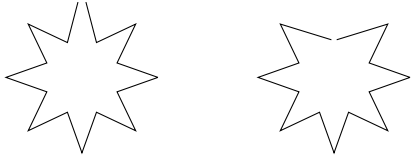
\begin{figure}
\centering
\begin{subfigure}{0.2\linewidth}
\centering
\begin{tikzpicture}
\path[draw] (0:1) -- (22.5:0.5) -- (45:1) -- (67.5:0.5) -- (87:1) (93:1) -- (112.5:0.5) -- (135:1) -- (157.5:0.5) -- (180:1) -- (-157.5:0.5) -- (-135:1) -- (-112.5:0.5) -- (-90:1) -- (-67.5:0.5) -- (-45:1) -- (-22.5:0.5) -- (0:1);
\end{tikzpicture}
\end{subfigure}
\begin{subfigure}{0.2\linewidth}
\centering
\begin{tikzpicture}
\path[draw] (0:1) -- (22.5:0.5) -- (45:1) -- (86:0.5) (94:0.5) -- (135:1) -- (157.5:0.5) -- (180:1) -- (-157.5:0.5) -- (-135:1) -- (-112.5:0.5) -- (-90:1) -- (-67.5:0.5) -- (-45:1) -- (-22.5:0.5) -- (0:1);
\end{tikzpicture}
\end{subfigure}
\caption{Fictitious discrete immersed disks bounded by segments of zigzag paths}
\label{fig:disks-zigzag-boundary}
\end{figure}

\begin{lemma}
\label{th:disks-zigzag-boundary}
Let $ Q $ be a zigzag consistent dimer and $ L $ a zigzag path. Then a segment of $ L $ cannot bound a discrete immersed disk.
\end{lemma}

\begin{proof}
Assume towards contradiction that $ D: P_k → |Q| $ is a discrete immersed disk bounded by a segment of $ L $. Then at least one of the consecutive interior angles of the discrete immersed disk must be an outer $ δ $-angle $ α' $ and hence consist of at least two indecomposable angles. The rest of the argument is a standard contradiction from geometric consistency theory: Construct a zigzag path $ L' $ at $ α' $ that points inside the discrete immersed disk. Follow $ L' $ until it leaves the discrete immersed disk. To be precise, “leaving” refers not to the leaving an area in $ |Q| $, but to touching the boundary $ ∂P_k $ in the domain of the underlying polygon immersion $ D: P_k → |Q| $. Either way, the final arc $ a $ before definitely leaving the discrete immersed disk lies on the boundary of the disk due to the zigzag nature of $ L $. The $ L $ and $ L' $ segments from $ α' $ until $ a $ are homotopic, since both lie in the discrete immersed disk. We obtain a contradiction with geometric consistency. This shows that $ L $ cannot bound a discrete immersed disk.
\end{proof}

\begin{lemma}
\label{th:deformed-uncurving-curvaturefree}
Every deformed zigzag path $ L $ is curvature-free.
\end{lemma}

\begin{proof}
The curvature of $ L $ as object in $ \Tw\Gtl_q Q $ is
\begin{equation*}
μ^0_{\Add\Gtl_q Q, L} + μ^1_{\Add\Gtl_q Q} (δ) + μ^2_{\Add\Gtl_q Q} (δ, δ) + …
\end{equation*}
It is our task to prove that this curvature vanishes. A first observation is that the two terms $ μ_{\Add\Gtl_q Q, L}^0 $ and $ μ_{\Add\Gtl_q Q}^2 (δ, δ) $ precisely cancel each other and we have $ μ_{\Add\Gtl_q Q}^1 = 0 $. By \autoref{th:disks-zigzag-boundary}, a segment of $ L $ cannot bound a discrete immersed disk and we conclude that $ μ^{≥3}_{\Add\Gtl_q Q} (δ, …) = 0 $. Adding up all 4 terms, we see that the curvature of the deformed zigzag path vanishes. This finishes the proof.
\end{proof}

\subsection{EFGH disks}
\label{sec:deformed-differential}
In this section, we develop elementary understanding of the differential $ μ^1_{\DefZigzagCat} $. The starting point is the description of the category $ \DefZigzagCat $ by explicit twisted complexes. The differential $ μ^1_{\DefZigzagCat} $ does not vanish, but counts those discrete immersed disks where apart from one single angle all interior angles stem from the $ δ $-matrix of the twisted complex. In this section, we classify these disks into four types which we call E, F, G and H disks. The goal is to be able to say:

\begin{center}
\begin{tikzpicture}
\path (0, 0) node (A) {Terms occurring in $ μ^1_{\DefZigzagCat} (ε) $} (8, 0) node (B) {E, F, G, H disks of $ ε $.};
\path[draw, <->] ($ (A.east)!0.2!(B.west) $) -- ($ (A.east)!0.8!(B.west) $);
\end{tikzpicture}
\end{center}

We start with a basic analysis of $ μ^1_{\DefZigzagCat} $. Recall we write $ μ_q $ for $ μ_{\Add\Gtl_q Q} $. Let $ ε: L_1 → L_2 $ be an elementary morphism. Then we have
\begin{equation*}
μ^1_{\DefZigzagCat} (ε) = \sum_{k, l ≥ 0} μ_q (\underbrace{δ, …, δ}_k, ε, \underbrace{δ, …, δ}_l).
\end{equation*}
Here each $ δ $ insertion stands for the $ δ $-matrix of $ L_1 $ or $ L_2 $, depending on whether $ δ $ stands right or left of $ ε $. The individual summands $ μ^1_q (δ, …, ε, …, δ) $ can again be expanded by writing the $ δ $-matrix as the sum of its entries, the $ δ $-angles. The elegant way to capture all the terms arising this way is as follows: We say a disk \emph{made of} $ μ_q (δ, …, ε, …, δ) $ is a final-out, first-out or all-in disk with angle sequence consisting of $ ε $, preceded and succeeded by an arbitrary number of $ δ $-angles. This way, we have enumerated all terms contributing to $ μ^1_{\DefZigzagCat} (ε) $.

\begin{figure}
\begin{subfigure}{0.23\linewidth}
\centering
\begin{tikzpicture}
\path[draw] (0, 0) -- ++(45:2) coordinate[pos=0.7] (alphat-end) coordinate[pos=0.6] (add-3) -- ++(135:2) coordinate[pos=0.3] (alphat-start) coordinate[pos=0.4] (add-2) coordinate[pos=0.7] (eps1-start) -- ++(45:2) coordinate[pos=0.3] (eps1-end) coordinate[pos=0.7] (alpha-end) node[near end, left] {$ L_1 $};
\path[draw] (3, 0) -- ++(135:2) coordinate[pos=0.7] (gammat-start) coordinate[pos=0.6] (add-4) -- ++(45:2) coordinate[pos=0.3] (gammat-end) coordinate[pos=0.4] (add-1) coordinate[pos=0.7] (eps2-end) -- ++(135:2) coordinate[pos=0.3] (eps2-start) coordinate[pos=0.7] (gamma-start) node[near end, right] {$ L_2 $};
\path[draw, bend right=45] (alphat-end) to node[midway, above] {$ β $} (gammat-start);
\path[draw, bend right=45] (gammat-start) to node[midway, left] {$ γ $} (gammat-end);
\path[draw, ->, bend right=45] (gammat-end) to node[midway, below] {$ β' $} (alphat-start);
\path[draw, bend right=45] (add-1) to (add-2);
\path[draw, bend right=45] (add-2) to (add-3);
\path[draw, ->, bend right=45] (add-3) to (add-4);
\path[draw, ->, bend right=45] (eps1-start) to (eps1-end);
\path[draw, ->, bend right=45] (eps2-start) to (eps2-end);
\path[draw, ->, bend right=45] (alpha-end) to node[midway, above] {$ ε $} (gamma-start);
\end{tikzpicture}
\caption{Type E disk}
\label{fig:disks-type-E}
\end{subfigure}
\begin{subfigure}{0.23\linewidth}
\centering
\begin{tikzpicture}
\path[draw] (-0.7, -0.7) -- (0, 0) coordinate[pos=0.5] (alpha-end) -- ++(135:1.5) coordinate[pos=0.4] (alpha-start) -- ++(45:1) coordinate (checkpoint) -- ++(135:1) -- ++(45:1.5) coordinate[pos=0.7] (eps-start);
\path[draw] (0.2, 1) -- (0.1, 0) coordinate[pos=0.4] (gamma-end) -- ++(45:1.5) coordinate[pos=0.4] (gamma-start) -- ++(right:0.5) -- ($ (checkpoint) + (0.8, 0) $) -- ++(45:1) -- ++(135:1.5) coordinate[pos=0.7] (eps-end);
\path[draw, ->, bend right=45] (eps-start) to node[midway, above] {$ ε $} (eps-end);
\path[draw, bend right=90, looseness=1.5] (alpha-end) to node[midway, above] {$ β $} (gamma-start);
\path[draw, bend right=15] (gamma-start) to node[pos=0.6, below] {$ γ $} (gamma-end);
\path[draw, ->, bend right=20] (gamma-end) to node[midway, below] {$ β' $} (alpha-start);
\end{tikzpicture}
\caption{Type F disk}
\label{fig:disks-type-F}
\end{subfigure}
\begin{subfigure}{0.23\linewidth}
\centering
\begin{tikzpicture}
\path[draw] (0, 0) -- ++(90:2) coordinate[pos=0.7] (alphat-end) -- ++(135:2) coordinate[pos=0.3] (alphat-start) coordinate[pos=0.2] (betat-short-end) coordinate[pos=0.7] (eps1-start) -- ++(45:2) coordinate[pos=0.3] (eps1-end) coordinate[pos=0.7] (alpha-end) node[near end, left] {$ L_1 $};
\path[draw] (0.2, 0) -- ++(90:2) coordinate[pos=0.8] (gammat-start) -- ++(45:2) coordinate[pos=0.3] (gammat-end) coordinate[pos=0.2] (gammat-short-end) coordinate[pos=0.7] (eps2-end) -- ++(135:2) coordinate[pos=0.3] (eps2-start) coordinate[pos=0.7] (gamma-start) node[near end, right] {$ L_2 $};
\path[draw, ->, bend right=70] (alphat-start) to node[midway, right] {$ α_1 $} (alphat-end);
\path[draw, bend right=70, looseness=1.5] (gammat-start) to node[midway, left, shift={(0.1, 0)}] {$ α_2 $} (gammat-short-end);
\path[draw, bend right=45] (gammat-end) to node[midway, below, shift={(0, -0.2)}] {$ β $} (alphat-start);
\path[draw, ->, bend right=45] (gammat-short-end) to (betat-short-end);
\path[draw, ->, bend right=45] (eps1-start) to (eps1-end);
\path[draw, ->, bend right=45] (eps2-start) to (eps2-end);
\path[draw, ->, bend right=45] (alpha-end) to node[midway, above] {$ ε $} (gamma-start);
\end{tikzpicture}
\caption{Type G disk}
\label{fig:disks-type-G}
\end{subfigure}
\begin{subfigure}{0.23\linewidth}
\centering
\begin{tikzpicture}
\path[draw] (0, 0) -- ++(225:1.5) coordinate[pos=0.4] (eps-start) -- ++(315:1) -- ++(225:1) coordinate (L1-left) coordinate[pos=0.3] (alpha3-end) -- ++(right:2) coordinate[pos=0.4] (alpha4-end) coordinate[pos=0.75] (alpha2-start) coordinate (L1-right) -- ++(250:1) coordinate[pos=0.5] (alpha2-end);
\path[draw] (0.1, 0) -- ++(315:1.5) coordinate[pos=0.4] (eps-end) -- ++(225:1) -- ++(315:0.85) coordinate (L2-right) coordinate[pos=0.3] (alpha1-start) -- ++(left:1.9) coordinate[pos=0.35] (alpha1-end) coordinate[pos=0.7] (alpha3-start) coordinate (L2-left) -- ++(290:1) coordinate[pos=0.5] (alpha4-start);
\path[draw, ->] (L1-left) to (L1-right);
\path[draw, ->] (L2-left) to (L2-right);
\path[draw, ->, bend right=45] (eps-start) to node[midway, below] {$ ε $} (eps-end);
\path[draw, bend right=20] (alpha3-start) to node[at end, below] {$ α_3 $} (alpha3-end);
\path[draw, ->, bend right=110, looseness=4] (alpha3-end) to node[midway, right] {$ β $} (alpha4-start);
\path[draw, bend right=110, looseness=3] (alpha2-end) to node[midway, left] {$ β' $} (alpha1-start);
\path[draw, ->, bend right=25] (alpha1-start) to node[at start, below] {$ α_1 $} (alpha1-end);
\end{tikzpicture}
\caption{Type H disk}
\label{fig:disks-type-H}
\end{subfigure}
\caption{Illustration of type E, F, G, H disks. The type F disk is depicted in case $ L_1 $ turns outside and $ L_2 $ turns inside the discrete immersed disk at the concluding puncture. For the type G disk, the naming of $ α_1 $ and $ α_2 $ is in case the arc outside the disk points upwards (type G1). For the type H disk, the naming of $ α_1 $ and $ α_3 $ are in case the concluding arc points to the right.}
\label{fig:disks-type-EFGH}
\end{figure}
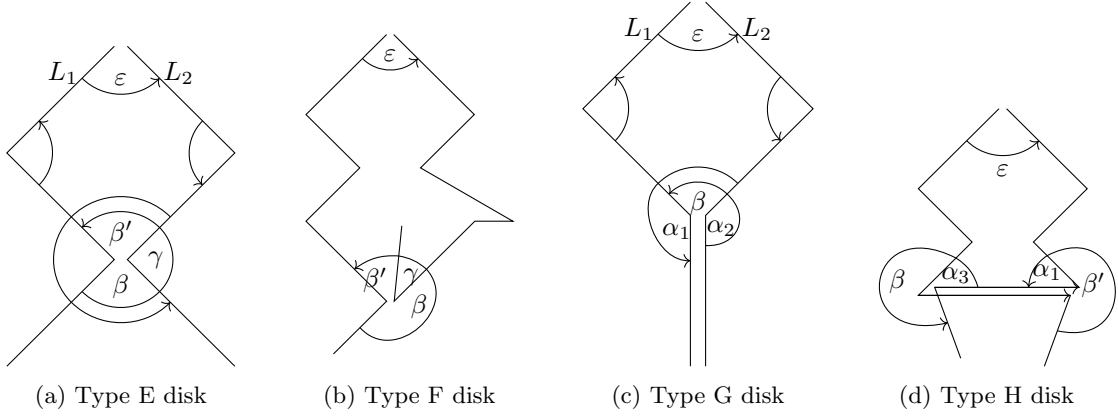

In \autoref{def:deformed-differential-EFGH}, we categorize the disks made of $ μ_q (δ, …, ε, …, δ) $ into four types. As a starting point, every disk made of $ μ_q (δ, …, ε, …, δ) $ is by definition a discrete immersed disk, together with possibly an outside morphism $ β $ or $ γ $ in the terminology of \autoref{def:prelim-terminology-firstfinal}. In the categorization and its proof, we make heavy use of the slots, concluding puncture and concluding arc terminology introduced in \autoref{sec:prelim-terminology}. For instance, we may say that $ L_1 $ “turns right towards the concluding puncture of the disk”. As an example, in \autoref{fig:disks-type-E} the zigzag path $ L_1 $ turns right towards the concluding puncture and $ L_2 $ turns left towards the concluding puncture. The categorization is obtained by a case distinction based on the behavior of $ L_1 $ and $ L_2 $ towards the concluding puncture.

\begin{definition}
\label{def:deformed-differential-EFGH}
Let $ ε: L_1 → L_2 $ be an elementary morphism. Then a disk that can be made of $ μ_q (δ, …, ε, …, δ) $ is of
\begin{itemize}
\item \emph{type E} if it is some-out, $ L_1 $ turns right towards the concluding puncture, $ L_2 $ turn left towards the concluding puncture, and there are at least 3 slots outside the disk,
\item \emph{type F} if it is some-out, there are at least two slots both inside and outside the disk, and $ L_1 $ turns left and $ L_2 $ right towards the concluding puncture, or the other way around,
\item \emph{type G} if it is some-out, $ L_1 $ turns right towards the concluding puncture and $ L_2 $ turns left towards to concluding puncture and there are 2 slots outside the disk,
\item \emph{type H} if it is all-in; or if it is first-out, $ L_2 $ turns right towards the concluding puncture and there is only 1 slot inside the disk; or if it is final-out, $ L_1 $ turns left towards the concluding puncture and there is only 1 slot inside the disk.
\end{itemize}
A disk of type G is of \emph{type G1} if the first, shared, arc of $ L_1 $ and $ L_2 $ at the concluding puncture, outside the disk, is oriented towards the concluding puncture, and of \emph{type G2} if the arc is oriented away from the puncture.
\end{definition}

The terminology is depicted in \autoref{fig:disks-type-EFGH}. We now show that the types E, F, G, H indeed provide an exhaustive classification of disks that can be made of $ μ_q (δ, …, ε, …, δ) $. During the proof, we will frequently use slots terminology from \autoref{sec:prelim-terminology}. We also show that some of the disk types come in pairs or triples. For instance, type E disks come in pairs. By this, we mean that every type E disk comes naturally with a distinct partner also of type E.

\begin{lemma}
Let $ ε: L_1 → L_2 $ be an elementary morphism. Then type E, F, G, H provide an exhaustive classification of disks that can be made from $ μ_q (δ, …, ε, …, δ) $. Type E disks come naturally in pairs, type F disks come alone, type G disks come in pairs and type H disks come in triples.
\end{lemma}

\begin{proof}
We first prove that every disk that can be made of $ μ_q (δ, …, ε, …, δ) $ falls under one of the four types and then comment on the pairs and triples.

Let $ D $ be a disk that can be made of $ μ_q (δ, …, ε, …, δ) $. Then $ D $ is either first-out, final-out or all-in. We shall in all three cases that $ D $ falls under our classification. The simplest way to understand the proof is by trying to recognize the properties we derive about $ D $ in \autoref{fig:disks-type-EFGH}.

Assume $ D $ is first-out. Then the first morphism must be a $ δ $ insertion and not $ ε $, otherwise $ L_2 $ would bound a discrete immersed disk, in conflict with geometric consistency. The $ δ $-insertion necessarily concerns an outer $ δ $-angle, as opposed to an inner $ δ $-angle, and we conclude that $ L_1 $ turns right at the concluding puncture. There are at least two slots outside of $ D $ at the concluding puncture. If $ L_2 $ turns left at the concluding puncture, then $ D $ is of type E or G, depending on the number of slots outside $ D $. If $ L_2 $ turns right at the concluding puncture, then $ D $ is of type F or H, depending on the number of slots inside $ D $.

A similar classification holds if $ D $ is final-out. Finally assume $ D $ is all-in. Then the concluding arc belongs to both $ L_1 $ and $ L_2 $. Including the concluding arc, the $ L_1 $ and $ L_2 $ segments bounding the disk are at least 2 arcs long, since otherwise the $ L_2 $ or $ L_1 $ segment would bound a discrete immersed disk in conflict with geometric consistency. Let us analyze how $ L_1 $ and $ L_2 $ continue beyond the concluding arc of $ D $, away from their segments that bound $ D $. Imagine that $ L_1 $ turns left at the head (and tail) of the concluding arc. Then $ L_1 $ enters the interior of $ D $, in conflict with geometric consistency, or lands on the arc of $ L_2 $ before the concluding arc, which would render $ L_1 = L_2 $ and mean that $ L_1 $ bounds a discrete immersed disk. We conclude that $ L_1 $ turns right and $ L_2 $ turns left at the concluding arc. In particular, there are at least two slots on the outside of the disk. This constitutes a type H disk. We have finished the first part of the proof.

For the second part of the proof, let us comment on the pairs and triples. We shall here restrict to the case of type E disks, since the other cases are similar. To show that type E disks come in pairs, the idea is to simply match two type E disks with each other by swapping the $ δ $ insertions: Regard a first-out type E disk. Then its first angle is a $ δ $ insertion. Remove this $ δ $ insertion and instead append a $ δ $ insertion as final angle. The result is a final-out disk, the desired partner disk. The partner disk is also of type E. The first angle of the first-out disk and the final angle of its final-out partner disk are depicted in \autoref{fig:disks-type-E} as well. Both disks have the same underlying discrete immersed disk, up to cyclically permuting the inputs and the output by one. The other types F, G, H are similar, and the relevant first/final angles are depicted in \autoref{fig:disks-type-EFGH} as well. This finishes the proof.
\end{proof}

\begin{figure}
\centering
\begin{subfigure}{0.3\linewidth}
\centering
\begin{tikzpicture}
\path[draw] (0, 0) -- ++(down:0.5) -- ++(210:1) -- ++(330:0.5) -- ++(210:0.5) -- ++(330:0.9) -- ++(210:0.7) coordinate[midway] (eps-start);
\path[draw] (0.1, 0) -- ++(down:0.5) -- ++(330:1) -- ++(210:0.5) -- ++(330:0.5) -- ++(210:0.9) -- ++(330:0.7) coordinate[midway] (eps-end);
\path[draw, ->, bend right=60] (eps-start) to node[midway, above] {$ ε $} (eps-end);
\end{tikzpicture}
\end{subfigure}
or
\begin{subfigure}{0.3\linewidth}
\centering
\begin{tikzpicture}
\path[draw] (0, 0) -- ++(240:2) coordinate[midway] (eps-start) node[midway, left] {$ L_1 $};
\path[draw] (0, -0.2) -- ++(300:1.8) coordinate[pos=0.45] (eps-end) node[near end, left] {$ L_2 $};
\path[draw] (0, 0) -- ++(300:2);
\path[draw, ->, bend right=30] (eps-start) to node[midway, above] {$ ε $} (eps-end);
\end{tikzpicture}
\end{subfigure}
or
\begin{subfigure}{0.3\linewidth}
\centering
\begin{tikzpicture}
\path[draw] (0, 0) -- ++(240:2);
\path[draw] (0, -0.2) -- ++(240:1.8) coordinate[pos=0.45] (eps-start) node[near end, right] {$ L_1 $};
\path[draw] (0, 0) -- ++(300:2) coordinate[midway] (eps-end) node[midway, right] {$ L_2 $};
\path[draw, ->, bend right=30] (eps-start) to node[midway, above] {$ ε $} (eps-end);
\end{tikzpicture}
\end{subfigure}
\caption{If $ L_1 $ and $ L_2 $ intersect above or at $ ε $, then $ ε $ has only type E disks.}
\label{fig:disks-intersection-above-eps}
\end{figure}
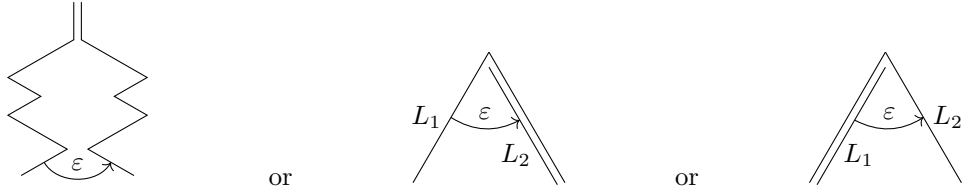

\subsection{Deformed differential}
\label{sec:deformed-tails}
In this section, we investigate the precise shape of the differential $ μ^1_{\DefZigzagCat} $. The starting point is the description of possible output of $ μ^1_{\DefZigzagCat} $ in terms of E, F, G, H disks, according to \autoref{sec:deformed-differential}. For the purposes of the deformed Kadeishvili theorem, this description would not be sufficient. We therefore trace the shape of $ μ^1_{\DefZigzagCat} $ even further. As announced, this leads to the data structure which we call the tail of an elementary morphism.

Let $ ε: L_1 → L_2 $ be an elementary morphism. We have seen that the disks that can be made of $ μ_q (δ, …, ε, …, δ) $ are of type E, F, G, H and come in groups, which we call the \emph{disk shapes} of $ ε $. In order to describe $ μ^1_{\DefZigzagCat} $, we need to capture the situation near the concluding puncture or arc of the disk shapes. For instance for a type E disk shape $ D $, we shall assign to $ D $ the A situation at the concluding puncture of $ D $. An A situation is already given by specifying the its angles $ α $, $ β $, $ γ $ and $ β' $ as in \autoref{fig:situation-ABCD}. In generality, we fix the following notation:

\begin{definition}
Let $ ε: L_1 → L_2 $ be an elementary morphism.
\begin{itemize}
\item For each type E disk shape $ D $ of $ ε $, denote by $ (α^D, β^D, γ^D, β'{}^D) $ the A situation at the concluding puncture. Let $ s^D ∈ ℤ $ be the sum of the total $ \# $ signs of all $ δ $ insertions along the disk, including both $ \#α^D $ and $ \#γ^D $.
\item For each type F disk shape $ D $ of $ ε $, denote by $ (α^D, β^D, γ^D, β'^D) $ the A situation at the concluding puncture. Let $ s^D ∈ ℤ $ be the sum of the total $ \# $ signs of all $ δ $ insertions along the disk, including $ \#α_D $ if $ D $ is first-out and including $ \#γ_D $ if $ D $ is final-out.
\item For each type G disk shape $ D $ of $ ε $, denote by $ (α_1^D, α_2^D, α_3^D, α_4^D, β^D, β'{}^D) $ the B situation at the concluding puncture. Let $ s^D ∈ ℤ $ be the sum of the total $ \# $ signs of all $ δ $ insertions along the disk, including both $ \#α_1 $ and $ \#α_2 $ or $ \#α_3 $ and $ \#α_4 $, depending on the orientation of the arc $ 2^D = 5^D $.
\item For each type H disk shape $ D $ of $ ε $, denote by $ (α_1^D, α_2^D, α_3^D, α_4^D, β^D, β'{}^D) $ the C situation at the concluding puncture. Let $ s^D ∈ ℤ $ be the sum of the total $ \# $ signs of all $ δ $ insertions along the disk, including $ \# α_1 $ and $ \#α_3 $ or $ \#α_2 $ and $ \#α_4 $, depending on the orientation of the concluding arc.
\end{itemize}
In all cases, the element $ q^D ∈ ℂ⟦Q_0⟧ $ is the product of all punctures covered by the discrete immersed disk, including those punctures on the boundary whose $ δ $ insertion is an outer $ δ $ angle, and including the concluding puncture or both endpoints of the concluding arc.
\end{definition}

In \autoref{def:deformed-tails-tail}, we define tails of elementary morphisms. The motivation is as follows: To apply the deformed Kadeishvili theorem to $ \DefZigzagCat $, we need to provide the deformed counterparts of the cohomology basis elements from \autoref{sec:splitting-splitting}. For instance, for a given situation B cohomology basis element $ h = (-1)^{\#α_3 + 1} α_3 + (-1)^{\#α_4} α_4 ∈ H $ we need to find a deformed counterpart $ φ^{-1} (h) $ such that $ φ^{-1} (h) $ and $ h $ only differ by infinitesimal $ R $ terms and $ μ^1_{\DefZigzagCat} (φ^{-1} (h)) = 0 $. The first step is to note that $ μ^1_{\DefZigzagCat} (h) $ is described, among others, by the type E disks which can be made of its two components $ α_3 $ and $ α_4 $. Every type E disk shape $ D $ of $ α_3 $ gives a contribution to $ μ^1_{\DefZigzagCat} (α_3) $ of
\begin{equation}
\label{eq:deformed-tails-alpha3contrib}
(-1)^{s^D + \#γ^D + |γ^D β^D|} q^D γ^D β^D + (-1)^{s^D + \#α^D} q^D β^D α^D.
\end{equation}
Similarly, every type E disk of $ α_4 $ gives a contribution to $ μ^1_{\DefZigzagCat} (α_4) $. Counting these contributions together, we see that $ μ^1_{\DefZigzagCat} (h) $ already contains two terms for every type E disk shape of $ α_3 $ plus two terms for every type E disk shape of $ α_4 $. We see that $ μ^1_{\DefZigzagCat} (h) $ is far from zero. The deformed counterpart $ φ^{-1} (h) $ is given by adding $ R $ terms to $ h $ such that $ μ^1_{\DefZigzagCat} $ eventually becomes zero. We observe that adding a multiple of $ β^D $ for every disk shape $ D $ of $ α_3 $ or $ α_4 $ does the trick in that it kills the two terms in \eqref{eq:deformed-tails-alpha3contrib}. However, $ μ^1_{\DefZigzagCat} (β^D) $ does not equate only to the two terms in \eqref{eq:deformed-tails-alpha3contrib}, but also to terms coming from the E, F, G, H disks which can be made of $ β^D $ itself. In turn, we have to kill these terms by adding yet more $ R $ terms, and every time we add $ R $ terms we obtain new $ R $ terms which we kill again. This gives rise to a recursive terms killing process which we can fortunately organize in a hierarchical structure, the tail of $ α_3 $ and $ α_4 $. The precise definition for general elementary morphisms reads as follows:

\begin{definition}
\label{def:deformed-tails-tail}
Let $ ε: L_1 → L_2 $ be an elementary morphism. Then its \emph{tail} is the tree $ T $ defined as follows. Insert $ ε $ as root. For each disk shape $ D $ of $ ε $, attach $ D $ as a child, annotated additionally with the type of $ D $. Continue inductively: For each leaf $ D ∈ T $ of type E, attach all disk shapes of $ β^D $ as children, annotated with their types.

Let $ D ∈ T $ be a node of type E. Denote by $ D_0 = ε, …, D_n = D $ be the sequence of nodes from the root till $ D $. Set
\begin{equation*}
S^D = \sum_{i = 1}^n s^{D_i}, \quad Q^D = \prod_{i = 1}^n q^{D_i}.
\end{equation*}
The morphism $ ε $ is \emph{E-preserving} if its tail $ T $ only consists of type E disks, apart from $ ε $ itself.
\end{definition}

The schematic of tails is depicted in \autoref{fig:deformed-tails-example-large}. Roughly speaking, a tail collects sequences of type E disks where every $ β^D $ morphism serves as $ ε $ for the next item in the sequence. The tail also collects type F, G or H disks but does not trace them any further. A sample elementary morphism together with its tail is depicted in \autoref{fig:deformed-tails-example}. In drawing the elementary morphism and its disks, we have neglected the zigzag nature of the zigzag paths. In drawing the tail, we have only depicted the tree structure and the type indication on all nodes and ignored the discrete immersed disk and situation data.

The typical tail is best imagined as a linear chain of type E disks with possibly an F, G or H disk at the end. Theoretically, nonlinear tails exist, but they require an angle sequence which bounds more than a single discrete immersed disks. Such angle sequences exist, but are very large, see also \autoref{rem:prelim-gtlq-bennequin}.

\begin{figure}
\centering
\begin{subfigure}{0.3\linewidth}
\centering
\begin{tikzpicture}
\path node (A) {$ ε $}
node[below of=A] (B) {E} edge (A)
node[below of=B] (C) {E} edge (B)
node[below of=C] (D) {F} edge (C)
node[left of=D] {E} edge (C) node[right of=D] {H} edge (C);
\end{tikzpicture}
\caption{Tail in general}
\label{fig:deformed-tails-example-large}
\end{subfigure}
\begin{subfigure}{0.4\linewidth}
\centering
\begin{tikzpicture}
\path[draw] (0, 0) -- ++(330:1) -- ++(210:1) coordinate[pos=0.6] (1-start)-- ++(down:1) -- ++(330:1) -- ++(210:1) coordinate[pos=0.6] (2-start) -- ++(down:1) -- ++(330:1) -- ++(60:0.7);
\path[draw] (2, 0) -- ++(210:1) -- ++(330:1) coordinate[pos=0.6] (1-end) -- ++(down:1) -- ++(210:1) -- ++(330:1) coordinate[pos=0.6] (2-end) -- ++(down:1) -- ++(210:1) -- ++(330:1) coordinate[pos=0.6] (3-end);
\path[draw, ->, bend right=60] (1-start) to node[pos=0.5, above] {$ ε $} (1-end);
\path[draw, ->, bend right=60] (2-start) to (2-end);
\begin{scope}[shift={(3.5, -1.5)}]
\path node (A) {$ ε $}
node[below of=A] (B) {E} edge (A)
node[below of=B] {G} edge (B);
\end{scope}
\end{tikzpicture}
\caption{An elementary morphism and its tail}
\label{fig:deformed-tails-example-small}
\end{subfigure}
\caption{Illustration of tails}
\label{fig:deformed-tails-example}
\end{figure}

Depending on the further knowledge of an elementary morphism $ ε $, we can say a more about the structure of its tail. In fact, every morphism from a B or C situation is E-preserving by virtue of geometric consistency. The following lemma makes this precise. Its premise is depicted in \autoref{fig:disks-intersection-above-eps}.

\begin{lemma}
Let $ ε: L_1 → L_2 $ be an elementary morphism. Suppose that above $ ε $ the zigzag paths $ L_1 $ and $ L_2 $ intersect and their segments from $ ε $ until the intersection are homotopic. Or suppose that at $ ε $, the zigzag path $ L_1 $ turns to the the target arc of $ ε $ or $ L_2 $ turns to the source arc of $ ε $. Then $ ε $ is E-preserving. In particular, this applies if $ ε $ is a morphism from a B or C situation.
\end{lemma}

\begin{proof}
The first observation is that $ ε $ cannot have disk shapes of type G or H, since these would create a digon with the intersection above $ ε $. We argue that $ ε $ can also not have type F disk shapes. Indeed, the ray of $ L_1 $ or $ L_2 $ that turns into the interior of a discrete immersed disk would leave the disk at some point, creating a contractible self-intersection of $ L_1 $ or $ L_2 $ or a digon with the intersection above $ ε $. We conclude that $ ε $ has only disks of type E. The same argument can now be applied inductively to all children of $ ε $. Ultimately, the entire tail of $ ε $ consists only of type E disk shapes and we conclude $ ε $ is E-preserving.

Elementary morphisms of a B or C situation automatically satisfy the premises of the lemma, simply because the involved zigzag paths intersect at the arc which we called $ 2=5 $ in \autoref{fig:situation-ABCD}. This finishes the proof.
\end{proof}

As we will see, tails are indeed the right tool to describe the deformed counterparts of cohomology basis elements. In \autoref{th:tail-sum}, we prepare for this by explicitly decomposing $ μ^1_{\DefZigzagCat} $ with respect to the decomposition
\begin{equation*}
\Hom_{\DefZigzagCat} (L_1, L_2) = (ℂ⟦Q_0⟧ \htensor H) ⊕ μ^1_{\DefZigzagCat} (ℂ⟦Q_0⟧ \htensor R) ⊕ (ℂ⟦Q_0⟧ \htensor R).
\end{equation*}
Here $ H $ and $ R $ refer to the standard splitting for $ \ZigzagCat $ defined in \autoref{sec:splitting-splitting}.

\begin{lemma}
\label{th:tail-sum}
Let $ ε: L_1 → L_2 $ be an elementary morphism. Then modulo $ R $ we have
\begin{equation}
\begin{aligned}
\label{eq:tail-sum}
μ^1_{\DefZigzagCat} (ε) &= μ^2 (ε, δ) + μ^2 (δ, ε) \\
&~ + μ^1_{\DefZigzagCat} \left(\sum_{\substack{D ∈ T \\ \text{of type E} \\ D ≠ ε}} (-1)^{S^D + 1} Q^D β^D + \sum_{\substack{D ∈ T \\ \text{of type G1}}} (-1)^{S^D + 1} Q^D \id_{2^D → 5^D}\right) \\
&~ + \sum_{\substack{D ∈ T \\ \text{of type G1}}} (-1)^{S^D + 1} Q^D \big((-1)^{\#α_3^D + 1} α_3^D + (-1)^{\#α_4^D} α_4^D\big) \\
&~ + \sum_{\substack{D ∈ T \\ \text{of type G2}}} (-1)^{S^D} Q^D \big((-1)^{\#α_3^D + 1} α_3^D + (-1)^{\#α_4^D} α_4^D\big) \\
&~ + \sum_{\substack{D ∈ T \\ \text{of type H}}} (-1)^{S^D} Q^D \id_{2^D → 5^D}.
\end{aligned}
\end{equation}
\end{lemma}

\begin{proof}
Let us evaluate the right-hand side. We have
\begin{align*}
μ^1_{\DefZigzagCat} &\left(\sum_{\substack{D ∈ T \\ \text{of type G1}}} (-1)^{S^D + 1} Q^D \id_{2^D → 5^D}\right)\\
&= \sum_{\substack{D ∈ T \\ \text{of type G1}}} (-1)^{S^D + 1} Q^D \big((-1)^{\#α_1^D + 1} α_1^D + (-1)^{\#α_2^D} α_2^D + (-1)^{\#α_3^D} α_3^D + (-1)^{\#α_4^D + 1} α_4^D\big).
\end{align*}
Further for $ D ∈ T \setminus \{ε\} $ of type E, we have modulo $ R $ that
\begin{align*}
μ^1_{\DefZigzagCat} ((-1)^{S^D + 1} Q^D β^D) &= (-1)^{S^D + \#γ^D + |β^D| + 1} Q^D γ^D β^D + (-1)^{S^D + \#α^D} Q^D β^D α^D \\
& + \sum_{\substack{E ∈ C_T (D) \\ \text{of type E}}} (-1)^{S^D + 1} Q^D \big((-1)^{s^E + \#γ^E + |β^E| + 1} q^E γ^E β^E + (-1)^{s^E + \#α^E} q^E β^E α^E\big) \\
& + \sum_{\substack{E ∈ C_T (D) \\ \text{of type G1}}} (-1)^{S^D + 1} Q^D \big((-1)^{s^E + \#α_2^E + 1} q^E α_2^E + (-1)^{s^E + \#α_1^E} q^E α_1^E\big) \\
& + \sum_{\substack{E ∈ C_T (D) \\ \text{of type G2}}} (-1)^{S^D + 1} Q^D \big((-1)^{s^E + \#α_3^E + 1} q^E α_3^E + (-1)^{s^E + \#α_4^E} q^E α_4^E\big) \\
& + \sum_{\substack{E ∈ C_T (D) \\ \text{of type H}}} (-1)^{S^D + 1} Q^D (-1)^{s^E} q^E \id_{2^E → 5^E}.
\end{align*}
Here, we have stripped off type F disks and the two first- and final-out type H disks. Both yield multiples of $ β $ (A), $ β $ (C) and $ β' $ (C), which lie in $ R $. Let us now add up the right-hand side of \eqref{eq:tail-sum}. This becomes a telescopic sum: The E, G2 and H terms cancel pairwise and the G2 terms cancel in triples. Only terms coming directly from the root remain. Modulo $ R $, the right-hand side of \eqref{eq:tail-sum} now reads
\begin{align*}
& μ_q^2 (ε, δ) + μ_q^2 (δ, ε) \\
& ~ + \sum_{\substack{D ∈ C_T (ε) \\ \text{of type E}}} (-1)^{S^D + \#γ^D + |β^D| + 1} Q^D γ^D β^D + (-1)^{S^D + \#α^D} Q^D β^D α^D \\
& ~ + \sum_{\substack{D ∈ C_T (ε) \\ \text{of type G1}}} (-1)^{S^D + \#α_2^D + 1} Q^D α_2^D + (-1)^{S^D + \#α_1^D} Q^D α_1^D \\
& ~ + \sum_{\substack{D ∈ C_T (ε) \\ \text{of type G2}}} (-1)^{S^D + \#α_3^D + 1} Q^D α_3^D + (-1)^{S^D + \#α_4^D} Q^D α_4^D \\
& ~ + \sum_{\substack{D ∈ C_T (ε) \\ \text{of type H}}} (-1)^{S^D} Q^D \id_{2^D → 5^D}.
\end{align*}
Modulo $ R $, this is precisely $ μ^1_{\DefZigzagCat} (ε) $. Indeed, the terms missing for $ μ^1_{\DefZigzagCat} (ε) $ are type F disks and the two first- and final-out type H disks, which again lie in $ R $. This finishes the proof.
\end{proof}

\subsection{Deformed cohomology basis elements}
\label{sec:deformed-splitting}
In this section, we compute the deformed cohomology basis elements of $ \DefZigzagCat $. The starting point is the homological splitting $ H ⊕ I ⊕ R $ for $ \ZigzagCat $ from \autoref{sec:splitting-splitting}. This splitting itself is not a homological splitting for the deformed category $ \DefZigzagCat $. Rather, we show that $ \DefZigzagCat $ together with $ H ⊕ I ⊕ R $ falls under the “$ D = 0 $” case of our deformed Kadeishvili theorem studied in \autoref{sec:2Bkadeishvili}. Accordingly, the category $ \DefZigzagCat $ comes with an associated homological splitting, including a list of deformed counterparts $ φ^{-1} (h) $ of the cohomology basis elements $ h ∈ H $. In the present section, we compute all these deformed cohomology basis elements.

\begin{proposition}
\label{th:deformed-cohomology-basis}
For the category $ \DefZigzagCat ⊂ \Tw\Gtl_q $ we have
\begin{equation*}
μ^1_{\DefZigzagCat} (H) ⊂ μ^1_{\DefZigzagCat} (ℂ⟦Q_0⟧ \htensor R).
\end{equation*}
Hence the deformed Kadeishvili construction of \autoref{sec:2Bkadeishvili} applies to $ \DefZigzagCat $. It produces a deformed homological splitting
\begin{equation*}
\DefZigzagCat = H_q ⊕ μ^1_{\DefZigzagCat} (ℂ⟦Q_0⟧ \htensor R) ⊕ (ℂ⟦Q_0⟧ \htensor R)
\end{equation*}
and a minimal model $ \H\DefZigzagCat $. For each cohomology basis element $ h ∈ H $, we obtain a deformed counterpart $ φ^{-1} (h) $, explicitly given as follows:
\begin{itemize}
\item The deformed counterpart of $ (-1)^{\#α_3 + 1} α_3 + (-1)^{\#α_4} α_4 $ (B) is
\begin{equation}
\label{eq:deformed-cohomology-basis-alpha34}
\begin{aligned}
& (-1)^{\#α_3 + 1} α_3 &+& \sum_{D ∈ T(α_3) \setminus \{α_3\}} (-1)^{\#α_3 + S^D + 1} Q^D β^D \\
+ & (-1)^{\#α_4} α_4 &+& \sum_{D ∈ T(α_4) \setminus \{α_4\}} (-1)^{\#α_4 + S^D} Q^D β^D.
\end{aligned}
\end{equation}
\item The deformed counterpart of $ \id_{2→5} $ (C) is
\begin{equation}
\begin{aligned}
\label{eq:deformed-cohomology-basis-idC}
\id_{2→5} & + (-1)^{\#α_1 + \#α_2} q_1 \left(β' + \sum_{D ∈ T(β') \setminus \{β'\}} (-1)^{S^D} Q^D β^D\right) \\
& + (-1)^{\#α_3 + \#α_4} q_2 \left(β + \sum_{D ∈ T(β) \setminus \{β\}} (-1)^{S^D} Q^D β^D\right).
\end{aligned}
\end{equation}
\item The deformed counterpart of $ \sum_a \id_a $ (D) is
\begin{equation*}
\sum \id_a.
\end{equation*}
\item The deformed counterpart of $ (-1)^{\#α_0 + 1} α_0 $ (D) is
\begin{equation*}
(-1)^{\#α_0 + 1} α_0 + (-1)^{\#α_0} q α_0'.
\end{equation*}
\end{itemize}
In the case of $ \id_{2→5} $, the punctures $ q_1, q_2 ∈ Q_0 $ are the head and tail of arc 2. In the case of $ (-1)^{\#α_0 + 1} α_0 $, the puncture $ q $ is the one around which $ α_0 $ turns. The codifferential is denoted $ h_q $ and the projection onto $ H_q $ is denoted $ π_q $.
\end{proposition}

\begin{proof}
All the $ β $ (A) and $ α_0' $ (D) morphisms added in the claimed deformed cohomology basis elements lie in $ ℂ⟦Q_0⟧ \htensor R $. In order to show $ μ^1_{\DefZigzagCat} (H) ⊂ μ^1_{\DefZigzagCat} (ℂ⟦Q_0⟧ \htensor R) $, it therefore suffices to show that $ μ^1_{\DefZigzagCat} $ vanishes on all the four types of claimed deformed cohomology basis elements. We check this for all four types individually.

For the situation B type the vanishing amounts to applying \autoref{th:tail-sum} to $ α_3 $ and $ α_4 $, and adding the results. Since $ α_3 $ and $ α_4 $ are E-preserving, the complicated sums over type G and H disks vanish.

For the situation C type, note that $ β $ (C) and $ β' $ (C) are E-preserving. Therefore two applications of \autoref{th:tail-sum} and a direct computation give the following results, whose sum renders the desired differential indeed zero:
\begin{align*}
μ^1_{\DefZigzagCat} \left(β' + \sum_{D ∈ T(β') \setminus \{β'\}} (-1)^{S^D} Q^D β^D\right) &= (-1)^{\#α_1} α_1 β' + (-1)^{\#α_2 + 1} β' α_2, \\
μ^1_{\DefZigzagCat} \left(β + \sum_{D ∈ T(β) \setminus \{β\}} (-1)^{S^D} Q^D β^D\right) &= (-1)^{\#α_4} α_4 β + (-1)^{\#α_3 + 1} β α_3, \\
μ^1_{\DefZigzagCat} (\id_{2→5}) &= (-1)^{\#α_2 + 1} q_1 α_1 β' + (-1)^{\#α_1} q_1 β' α_2 \\
& \quad + (-1)^{\#α_3 + 1} q_2 α_4 β + (-1)^{\#α_4} q_2 β α_3,
\end{align*}
For the situation D identity $ \sum_a \id_a $ (D), note that no disk sequences can be made with an identity. The ordinary product of $ \id_a $ with a neighboring $ α_i' $ each appears twice in $ μ^1_{\DefZigzagCat} $ and they cancel each other:
\begin{equation*}
μ^1_{\DefZigzagCat} \left(\sum \id_a\right) = 0,
\end{equation*}
For the situation D co-identity $ (-1)^{\#α_0 + 1} α_0 $ (D), note that no disks can be made of $ μ_q (δ, …, α_0', …, δ) $ due to consistency. We get
\begin{equation*}
μ_{\DefZigzagCat}^1 ((-1)^{\#α_0 + 1} α_0 + (-1)^{\#α_0} q α_0') = - (- q α_0 α_0' - q α_0' α_0 + q α_0 α_0' + q α_0' α_0) = 0.
\end{equation*}
This finishes the proof.
\end{proof}

\subsection{Deformed codifferential and projection}
\label{sec:deformed-codif}
In this section, we compute part of the deformed codifferential and deformed projection for $ \DefZigzagCat $. The starting point is the homological splitting $ H ⊕ I ⊕ R $ for $ \ZigzagCat $. In \autoref{sec:deformed-splitting}, we have verified that $ \DefZigzagCat $ satisfies the requirements of the deformed Kadeishvili construction of \autoref{sec:2Bkadeishvili} so that we obtain a deformed homological splitting $ H_q ⊕ μ^1_{\DefZigzagCat} (ℂ⟦Q_0⟧ \htensor R) ⊕ (ℂ⟦Q_0⟧ \htensor R) $. According to \autoref{def:2Bkadeishvili-deformed-counterpart}, there is an associated deformed codifferential $ h_q $ given by projecting morphisms onto $ μ^1_{\DefZigzagCat} (ℂ⟦Q_0⟧ \htensor R) $ and finding the $ R $ preimage element under $ μ^1_{\DefZigzagCat} $. In the present section, we examine this procedure for the morphisms $ βα $ and $ β $ of A situations and indicate how one proceeds for other types of morphisms.

Let us recall how the (non-deformed) codifferential for $ \ZigzagCat $ works. Regard an A situation, given by angles $ (α, β, γ, β') $. Then $ β $ lies in $ R $, while $ βα $ lies in $ I+R $. According to \autoref{tab:coh_splitting-splitting-verification}, we have $ h(βα) = (-1)^{\#α + 1} β $. For the deformed codifferential, we however have to add terms. It is namely not true that $ μ^1_{\DefZigzagCat} ((-1)^{\#α + 1} β) = βα + (-1)^{\#α + \#γ + ∥β∥ + 1} γβ $. Rather, the differential $ μ^1_{\DefZigzagCat} (β) $ is computed by the formula \eqref{eq:tail-sum}. The formula implies we have to subtract terms from the expression $ (-1)^{\#α + 1} β $ in order to make its $ μ^1_{\DefZigzagCat} $ image up to $ R $-terms equal to $ βα $. In the following proposition, we compute these terms.

\begin{proposition}
\label{th:deformed-codif-ba}
Let $ (β, α, γ, β') $ denote an A situation $ L_1 → L_2 $. Denote by $ T(β) $ the tail of $ β $. Then we have
\begin{align*}
h_q (βα) &= (-1)^{\#α + 1} β + \subsum{D ∈ T(β) \setminus \{β\} \\ \text{of type E}} (-1)^{S^D + \#α + 1} Q^D β^D + \subsum{D ∈ T(β) \\ \text{of type G1}} (-1)^{S^D + \#α + 1} Q^D \id_{2^D → 5^D}, \\
φ π_q (βα) &= \subsum{D ∈ T(β) \\ \text{of type G1}} (-1)^{S^D + \#α + 1} Q^D \big((-1)^{\#α_3^D + 1} α_3^D + (-1)^{\#α_4^D} α_4^D\big) \\
& ~ + \subsum{D ∈ T(β) \\ \text{of type G2}} (-1)^{S^D + \#α} Q^D \big((-1)^{\#α_3^D + 1} α_3^D + (-1)^{\#α_4^D} α_4^D\big) \\
& ~ + \subsum{D ∈ T(β) \\ \text{of type H}} (-1)^{S^D + \#α} Q^D \id_{2^D → 5^D}.
\end{align*}
\end{proposition}

\begin{proof}
Apply $ h_q $ and $ φ π_q $ on both sides of \autoref{th:tail-sum} with $ ε = β $.
\end{proof}

\begin{remark}
The morphisms of type $ βα $ (A) are the most important morphisms to which we would like to apply the deformed codifferential $ h_q $. There are only very few cases where we need to apply $ h_q $ to other morphisms. In fact, $ h_q $ vanishes by definition on all deformed cohomology basis elements and all elements of $ R $. The only interesting elementary morphisms which do not lie in $ H_q ⊕ (ℂ⟦Q_0⟧ \htensor R) $ are $ α_4 β $ (C), $ β' α_2 $ (C) and $ α_3 α_4 $ (B). For these three types of morphisms, one obtains formulas for their $ h_q $ and $ π_q $ values by applying \eqref{eq:tail-sum} to $ ε = β $ (C) or $ ε = β' $ (C) or $ ε = α_3 $ (B). The results are formulas very similar to \autoref{th:deformed-codif-ba}. In contrast to $ β $ (A), these three morphisms have the benefit of being E-preserving. Therefore all complicated G and H terms on the right-hand side of \eqref{eq:tail-sum} do not even appear and only $ μ^2 (δ, ε) + μ^2 (ε, δ) $ and the sum over type E nodes remain.
\end{remark}

\section{Result components of Kadeishvili trees}
\label{sec:resultcomp}
In this section, we develop a first glance at the minimal model $ \H\DefZigzagCat $. The starting point is the knowledge of $ \DefZigzagCat $ established in \autoref{sec:deformed} and the deformed Kadeishvili construction recalled in \autoref{sec:2Bkadeishvili}. According to the deformed Kadeishvili construction, the deformed $ A_∞ $-structure on the minimal model $ \H\DefZigzagCat $ is determined by Kadeishvili trees. The outstanding task is therefore to enumerate and analyze all results from all possible Kadeishvili trees.

In \autoref{sec:resultcomp-trees}, we explain which products are to be computed and set up notation. In \autoref{sec:resultcomp-possible}, we list possible types of morphisms resulting from Kadeishvili trees. In \autoref{sec:resultcomp-definition}, we introduce a notion of “result components” which allows us to systematically track terms arising from Kadeishvili trees. In \autoref{sec:resultcomp-classification}, we conclude the section with a semi-explicit, inductive characterization of result components:

\begin{center}
\begin{tikzpicture}
\path (0, 0) node[align=center] (A) {\textbf{Individual output terms} \\ of $ μ_{\H\DefZigzagCat} (h_k, …, h_1) $} (8, 0) node[align=center] (B) {\textbf{Result components} \\ classified by \autoref{tab:resultcomp-classification}};
\path[draw, <->] ($ (A.east)!0.2!(B.west) $) -- ($ (A.east)!0.8!(B.west) $);
\end{tikzpicture}
\end{center}

\subsection{Kadeishvili trees}
\label{sec:resultcomp-trees}
This section explains how our minimal model theorem applies to $ \DefZigzagCat $ specifically. We explain which trees need to be investigated, and which not. We also set up specific terminology. Note that we keep writing $ μ ≔ μ_q ≔ μ_{\Add \Gtl_q Q} $, see \autoref{rem:deformed-mushorthand}.

The deformed Kadeishvili construction instructs us to start with the hom spaces. If $ L_1 $ and $ L_2 $ are zigzag paths, then $ \Hom_{\H\DefZigzagCat} (L_1, L_2) = ℂ⟦Q_0⟧ \htensor H $, where $ H $ denotes the cohomology in the standard splitting of $ \Tw\Gtl Q $. The higher products on $ \H\DefZigzagCat $ are obtained as outputs of Kadeishvili trees, with $ φ^{-1}: ℂ⟦Q_0⟧ \htensor H → H_q $ applied at all leaves, $ h_q μ_{\DefZigzagCat} $ applied at all non-leaf nodes, and $ φ π_q μ_{\DefZigzagCat} $ applied at the root.

We do not need to calculate all trees. In fact, we observe directly from the Kadeishvili construction that $ μ_{\H\DefZigzagCat} $ is strictly unital, with the same unit morphisms as $ \DefZigzagCat $. More precisely, we already know that for every zigzag path $ L $ and compatible morphism $ h $ we have
\begin{align*}
μ^1_{\H\DefZigzagCat} &= 0, & μ^2_{\H\DefZigzagCat} (h, \id_L) &= h \\
μ^{≥3}_{\H\DefZigzagCat} (…, \id_L, …) &= 0, & μ^2_{\H\DefZigzagCat} (\id_L, h) &= (-1)^{|h|} h.
\end{align*}
It therefore suffices to regard Kadeishvili trees whose inputs are all non-identity cohomology basis elements. At this point, it is clever to have terminology available to study not only entire Kadeishvili trees, but also their subtrees. We fix the following terminology:

\begin{definition}
An \emph{h-tree} is an ordered tree $ T $ where each non-leaf node has at least two children, together with non-identity deformed basis cohomology elements $ h_1, …, h_N $ on the leaves from right to left, with $ h_i: L_i → L_{i+1} $.

A \emph{π-tree} is an ordered tree $ T $ with at least three nodes where each non-leaf node has at least two children, together with non-identity deformed basis cohomology elements $ h_1, …, h_N $ on the leaves from right to left, with $ h_i: L_i → L_{i+1} $.
\end{definition}

Both h-trees and π-trees can be evaluated. They have \emph{results} or \emph{outputs}. When evaluating an h-tree, we put $ h_q μ $ on every non-leaf node. When evaluating a π-tree, we put $ h_q μ $ on every non-leaf non-root node, and $ φ π_q μ $ on the root.

\begin{remark}
\label{rem:resultcomp-trees-productsremark}
Often we will make statements about “products”. The datum of a “product” shall then typically include all of its inputs. For example “a product $ μ^{≥3} (…) $” refers to a choice of arity $ k ≥ 3 $, a collection of compatible morphisms $ a_1, …, a_k $ and the result of the product itself.
\end{remark}

\begin{remark}
Since we abbreviate $ μ = μ_q = μ_{\Add\Gtl_q Q} $, a product $ μ^2 (a, b) $ stands simply for the product of angles and does not include discrete immersed disk terms like $ μ^3 (δ, a, b) $ stemming from twisted completion. Similarly, a product $ μ^{≥3} $ always stands for a single discrete immersed disk, and following the rule explained in \autoref{rem:resultcomp-trees-productsremark} it includes the datum of input morphisms some of which may be $ δ $-morphisms. Note that $ δ $-morphisms are always spelt out as $ α $ (D) or $ α' $ (D).
\end{remark}

\begin{remark}
We occasionally group $ β $ (C) and $ β' $ (C) as $ β/β' $ (C) due to their similar nature. We ignore signs and deformation parameters $ q ∈ ℂ⟦Q_0⟧ $ in this section. For example, we may say that a product like $ μ^{≥3} (α, α', α_3, …) $ is $ \id $ (D), meaning that it is equal to some arc identity, possibly multiplied by a sign and deformation parameters. We abbreviate a situation B cohomology basis element $ (-1)^{\#α_3 + 1} α_3 + (-1)^{\#α_4} α_4 $ simply as $ α_3 + α_4 $.
\end{remark}

\subsection{Possible tree output}
\label{sec:resultcomp-possible}
In this section, we analyze the possible types of output of Kadeishvili trees. The starting point is the observation that the result of a Kadeishvili tree is a linear combination of elementary morphism, but not every elementary morphism can actually appear. In the present section, we compose a tight list of possible elementary morphisms that can result from Kadeishvili trees.

\begin{table}
\centering
\renewcommand{\arraystretch}{1.1}
\begin{tabular}{cr|cccccc}
$ m_2 \backslash m_1 $ & & $ β $(A) & $ α_4 $(B) & $ α_3 $(B) & $ β/β' $(C) & $ α_0 $(D) & $ α_0' $(D) \\\hline
$ β $        & $  μ^2 = $ & imp & imp & $ β $(A)/$ β' $(C) & imp & $ βα $(A) & imp \\
             & $  h_q = $ &   0 & 0 & 0 & 0 & $ β $(A)+E/G & 0 \\
             & $ φπ_q = $ &   0 & 0 & 0 & 0 & G/H & 0 \\\hline
$ α_4 $      & $  μ^2 = $ & $ β' $(C)/$ β $(A) & imp & imp & $ α' $(D) & imp & imp \\
             & $  h_q = $ & 0 & 0 & 0 & 0 & 0 & 0 \\
             & $ φπ_q = $ & 0 & 0 & 0 & 0 & 0 & 0 \\\hline
$ α_3 $      & $  μ^2 = $ & imp & imp & imp & imp & imp & imp \\
             & $  h_q = $ & 0 & 0 & 0 & 0 & 0 & 0 \\
             & $ φπ_q = $ & 0 & 0 & 0 & 0 & 0 & 0 \\\hline
$ β $/$ β' $ & $  μ^2 = $ & imp & imp & $ α' $ & imp & $ β α_3 $/imp & imp \\
             & $  h_q = $ & 0 & 0 & 0 & 0 & 0 & 0 \\
             & $ φπ_q = $ & 0 & 0 & 0 & 0 & 0 & 0 \\\hline
$ α_0 $      & $  μ^2 = $ & $ γβ $ & $ α_3 α_4 $ & imp & imp/$ α_1 β' $ & imp & $ α_0 α_0' $ \\
             & $  h_q = $ & 0 & $ α_3 $+E & 0 & 0 & 0 & $ α_0' $ \\
             & $ φπ_q = $ & 0 & 0 & 0 & 0 & 0 & 0 \\\hline
$ α_0' $     & $  μ^2 = $ & imp & imp & $ α_4 β α_3 $ & imp & $ α_0' α_0 $ & imp \\
             & $  h_q = $ & 0 & 0 & 0 & 0 & 0 & 0 \\
             & $ φπ_q = $ & 0 & 0 & 0 & 0 & 0 & 0 \\\hline
$ \id $(B)   & $  μ^2 = $ & imp & imp & $ β $(A) & imp & $ α_1 $ & imp \\
             & $  h_q = $ & 0 & 0 & 0 & 0 & 0 & 0 \\
             & $ φπ_q = $ & 0 & 0 & 0 & 0 & 0 & 0 \\\hline
$ \id $(C)   & $  μ^2 = $ & $ γβ $/$ α_2 $/$ α_3 $ & $ α $(D) & imp & $ γβ $ & imp & $ α_4 β $ \\
             & $  h_q = $ & $ \id $(B) & $ \id $(D) & 0 & 0 & 0 & $ β $(C)+E \\
             & $ φπ_q = $ & 0/$ α_3 + α_4 $ & $ α_0 $ & 0 & 0 & 0 & 0 \\\hline
$ \id $(D)   & $  μ^2 = $ & $ β $(A) & $ α_4 $(B) & $ α_3 $(B) & $ β $/$ β' $(C) & $ α_0 $(D) & $ α_0' $(D) \\
             & $  h_q = $ & 0 & 0 & 0 & 0 & 0 & 0 \\
             & $ φπ_q = $ & 0 & $ α_3 + α_4 $ & 0 & 0 & $ α_0 $(D) & 0
\end{tabular}
\caption{Multiplication scheme}
\label{tab:components-multiplication}
\end{table}

\begin{table}
\ContinuedFloat
\centering
\renewcommand{\arraystretch}{1.1}
\begin{tabular}{cr|ccc}
$ m_2 \backslash m_1 $ & & $ \id $ (B) & $ \id $ (C) & $ \id $ (D) \\\hline
$ β $ (A) & $  μ^2 = $ & imp & $ βα $/$ α_1 $/$ α_4 $ & $ β $ (A) \\
          & $  h_q = $ & 0 & $ β $(A)+E/G & 0 \\
          & $ φπ_q = $ & 0 & G/H/$ α_3 + α_4 $ & 0 \\\hline
$ α_4 $   & $  μ^2 = $ & $ β $ (A) & imp & $ α_4 $ \\
          & $  h_q = $ & 0 & 0 & 0 \\
          & $ φπ_q = $ & 0 & 0 & $ α_3 + α_4 $ \\\hline
$ α_3 $   & $  μ^2 = $ & imp & $ α $ (D) & $ α_3 $ \\
          & $  h_q = $ & 0 & $ \id $ (D) & 0 \\
          & $ φπ_q = $ & 0 & $ α_0 $ & 0 \\\hline
$ β $/$ β' $ (C) & $  μ^2 = $ & imp & $ βα $ & $ β $/$ β' $ (C) \\
          & $  h_q = $ & 0 & $ β $(A)+E/G & 0 \\
          & $ φπ_q = $ & 0 & G/H & 0 \\\hline
$ α_0 $   & $  μ^2 = $ & $ α_3 $ & imp & $ α_0 $ \\
          & $  h_q = $ & 0 & 0 & 0 \\
          & $ φπ_q = $ & 0 & 0 & $ α_0 $ \\\hline
$ α_0' $  & $  μ^2 = $ & imp & $ β' α_2 $ & $ α_0' $ \\
          & $  h_q = $ & 0 & $ β' $(C)+E & 0 \\
          & $ φπ_q = $ & 0 & 0 & 0 \\\hline
$ \id $(B) & $ μ^2 = $ & imp & $ \id $ (D) & $ \id $ (B) \\
          & $  h_q = $ & 0 & 0 & 0 \\
          & $ φπ_q = $ & 0 & $ 1 $ if $ a=a_0 $ & 0 \\\hline
$ \id $(C) & $ μ^2 = $ & $ \id $ (D) & imp & $ \id $ (C) \\
          & $  h_q = $ & 0 & 0 & 0 \\
          & $ φπ_q = $ & $ 1 $ if $ a=a_0 $ & 0 & $ \id $ (C) \\\hline
$ \id $(D) & $ μ^2 = $ & $ \id $ (B) & $ \id $ (C) & $ \id $ (D) \\
          & $  h_q = $ & 0 & 0 & 0 \\
          & $ φπ_q = $ & 0 & $ \id $ (C) & $ 1 $ if $ a=a_0 $
\end{tabular}
\caption{Multiplication scheme (continued)}
\end{table}

By nature, a π-tree can only have a linear combination of cohomology basis elements as output. Similarly, an h-tree can only have a linear combination of $ R $ basis morphims as output. Every node of an h- or π-tree carries an evaluation result itself. To further narrow down on the possible output of the tree, we have to investigate what happens at every node in the tree. For instance, we claim that angle length cannot grow to infinity as we go from leaves to root. In fact, we claim there is a list $ S $ of morphisms which is stable under evaluations, in the sense that any $ h_q μ $ applied to a sequence of morphisms from $ S $ yields a morphism from $ S $ again. The explicit list reads as follows:

\begin{center}
$ S ~~ ≔ ~~ $ $ β $ (A), $ \id $ (B), $ α_3 $ (B), $ α_4 $ (B), $ \id $ (C), $ β $ (C), $ β' $ (C), $ \id $ (D), $ α_0 $ (D) and $ α_0' $ (D).
\end{center}

We claim that this list $ S $ is preserved under evaluations $ h_q μ $. Before we prove this, let us prepare reasoning. For all three cases of $ h_q μ^2 $, first-out $ h_q μ^{≥3} $ and final-out $ h_q μ^{≥3} $ evaluation, we have set up product schemes which indicate the type of output from in principle any kind of evaluations with arbitrary inputs from $ S $.

These product schemes are found in \autoref{tab:components-multiplication}, \autoref{fig:components-first-out} and \ref{fig:components-final-out}. They are generally structured by the three keys $ μ $, $ h_q $ and $ φπ_q $. The schemes should universally be read as follows: A product of morphisms of given types may yield only the types of morphisms indicated in the $ μ $ row. Application of the codifferential yields the morphism indicated in the $ h_q $ row. Of course, vanishing products are also possible. For later use, the possible results of $ φ π_q μ $ have been collected in the $ φ π_q $ row. With these product schemes in mind, we are ready to prove that the list $ S $ is preserved:

\begin{lemma}
Let $ T $ be an h-tree. Then its output contains only $ β $ (A), $ \id $ (B), $ α_3 $ (B), $ α_4 $ (B), $ \id $ (C), $ β $ (C), $ β' $ (C), $ \id $ (D), $ α_0 $ (D) and $ α_0' $ (D) terms.
\end{lemma}

\begin{figure}
\centering
\begin{subfigure}{0.24\linewidth}
\centering
\begin{tikzpicture}
\path[draw] (0, 0) coordinate (out) -- (0, 1) coordinate[pos=0.4] (m-start);
\path[draw] (0.1, -0.1) coordinate (start) -- ++(315:1) coordinate[pos=0.4] (m-end);
\path[draw] (-0.1, -0.1) coordinate (end) -- ++(225:1);
\path[draw, ->, bend right=110, looseness=3] (m-start) to (m-end);
\path[draw, decoration={border, amplitude=25pt, angle=90, segment length=4pt, mirror}, color=gray, decorate] ($ (start) + (-0.1, -0.1) $) -- ++(315:1);
\path[draw] (out) -- (start);
\path[draw] (end) -- ++(135:1);
\end{tikzpicture}
\caption{$ δ $ insertion \\ $ μ = γβ $ \\ $ h_q = 0 $ \\ $ φ π_q = 0 $}
\end{subfigure}
\begin{subfigure}{0.24\linewidth}
\centering
\begin{tikzpicture}
\path[draw] (0, 0) coordinate (out) -- (0, 1) coordinate[pos=0.4] (m-start);
\path[draw] (0.1, -0.1) coordinate (start) -- ++(315:1) coordinate[pos=0.4] (m-end);
\path[draw] (-0.1, -0.1) coordinate (end) -- ++(225:1);
\path[draw, ->, bend right=110, looseness=3] (m-start) to (m-end);
\path[draw, decoration={border, amplitude=25pt, angle=90, segment length=4pt, mirror}, color=gray, decorate] ($ (start) + (-0.1, -0.1) $) -- ++(315:1);
\path[draw] (out) -- (start);
\path[draw] (end) -- ++(90:1);
\end{tikzpicture}
\caption{$ δ $ insertion \\ $ μ = α_2 \text{ or } α_3 $ (B) \\ $ h_q = \id $ (B) or $ 0 $ \\ $ φ π_q = α_3 + α_4 $}
\end{subfigure}
\begin{subfigure}{0.24\linewidth}
\centering
\begin{tikzpicture}
\path[draw] (0, 0) coordinate (out) -- (0, 1) coordinate[pos=0.4] (m-start);
\path[draw] (0.1, -0.1) coordinate (start) -- ++(315:1) coordinate[pos=0.4] (m-end);
\path[draw] (-0.1, -0.1) coordinate (end) -- ++(225:1);
\path[draw, ->, bend right=110, looseness=3] (m-start) to (m-end);
\path[draw, decoration={border, amplitude=25pt, angle=90, segment length=4pt, mirror}, color=gray, decorate] ($ (start) + (-0.1, -0.1) $) -- ++(315:1);
\path[draw] (out) -- (start);
\path[draw] (end) -- ++(270:1);
\end{tikzpicture}
\caption{$ δ $ insertion \\ $ μ = β $ (A) \\ $ h_q = 0 $ \\ $ φ π_q = 0 $}
\end{subfigure}
\begin{subfigure}{0.24\linewidth}
\centering
\begin{tikzpicture}
\path[draw] (0, 0) coordinate (out) -- (0, 1) coordinate[pos=0.4] (m-start);
\path[draw] (0.1, -0.1) coordinate (start) -- ++(315:1) coordinate[pos=0.4] (m-end);
\path[draw] (-0.1, -0.1) coordinate (end) -- ++(225:1);
\path[draw, ->, bend right=110, looseness=3] (m-start) to (m-end);
\path[draw, decoration={border, amplitude=25pt, angle=90, segment length=4pt, mirror}, color=gray, decorate] ($ (start) + (-0.1, -0.1) $) -- ++(315:1);
\path[draw] (out) -- (start);
\path[draw] (end) -- ++(315:1);
\end{tikzpicture}
\caption{$ δ $ insertion \\ $ μ = β \text{ or } β' $ (C) \\ $ h_q = 0 $ \\ $ φ π_q = 0 $}
\end{subfigure}
\begin{subfigure}{0.24\linewidth}
\centering
\begin{tikzpicture}
\path[draw] (0, 0) coordinate (out) -- (0, 1) coordinate[pos=0.4] (m-start);
\path[draw] (0.1, -0.1) coordinate (start) -- ++(315:1) coordinate[pos=0.4] (m-end);
\path[draw] (-0.1, -0.1) coordinate (end) -- ++(225:1);
\path[draw, ->, bend right=110, looseness=3] (m-start) to (m-end);
\path[draw, decoration={border, amplitude=25pt, angle=90, segment length=4pt, mirror}, color=gray, decorate] ($ (start) + (-0.1, -0.1) $) -- ++(315:1);
\path[draw] (out) -- ++(30:1);
\path[draw] (start) -- ++(30:1);
\path[draw] (end) -- ++(150:1);
\end{tikzpicture}
\caption{$ β \text{ or } β' $ (C) \\ $ μ = γβ $ \\ $ h_q = 0 $ \\ $ φ π_q = 0 $}
\end{subfigure}
\begin{subfigure}{0.24\linewidth}
\centering
\begin{tikzpicture}
\path[draw] (0, 0) coordinate (out) -- (0, 1) coordinate[pos=0.4] (m-start);
\path[draw] (0.1, -0.1) coordinate (start) -- ++(315:1) coordinate[pos=0.4] (m-end);
\path[draw] (-0.1, -0.1) coordinate (end) -- ++(225:1);
\path[draw, ->, bend right=110, looseness=3] (m-start) to (m-end);
\path[draw, decoration={border, amplitude=25pt, angle=90, segment length=4pt, mirror}, color=gray, decorate] ($ (start) + (-0.1, -0.1) $) -- ++(315:1);
\path[draw] (out) -- ++(30:1);
\path[draw] (start) -- ++(30:1);
\path[draw] (end) -- ++(90:1);
\end{tikzpicture}
\caption{$ β \text{ or } β' $ (C) \\ $ μ = α_2 \text{ or } α_3 $ \\ $ h_q = \id $ (B) or $ 0 $ \\ $ φ π_q = α_3 + α_4 $}
\end{subfigure}
\begin{subfigure}{0.24\linewidth}
\centering
\begin{tikzpicture}
\path[draw] (0, 0) coordinate (out) -- (0, 1) coordinate[pos=0.4] (m-start);
\path[draw] (0.1, -0.1) coordinate (start) -- ++(315:1) coordinate[pos=0.4] (m-end);
\path[draw] (-0.1, -0.1) coordinate (end) -- ++(225:1);
\path[draw, ->, bend right=110, looseness=3] (m-start) to (m-end);
\path[draw, decoration={border, amplitude=25pt, angle=90, segment length=4pt, mirror}, color=gray, decorate] ($ (start) + (-0.1, -0.1) $) -- ++(315:1);
\path[draw] (out) -- ++(30:1);
\path[draw] (start) -- ++(30:1);
\path[draw] (end) -- ++(270:1);
\end{tikzpicture}
\caption{$ β \text{ or } β' $ (C) \\ $ μ = β $ (A) \\ $ h_q = 0 $ \\ $ φ π_q = 0 $}
\end{subfigure}
\begin{subfigure}{0.24\linewidth}
\centering
\begin{tikzpicture}
\path[draw] (0, 0) coordinate (out) -- (0, 1) coordinate[pos=0.4] (m-start);
\path[draw] (0.1, -0.1) coordinate (start) -- ++(315:1) coordinate[pos=0.4] (m-end);
\path[draw] (-0.1, -0.1) coordinate (end) -- ++(225:1);
\path[draw, ->, bend right=110, looseness=3] (m-start) to (m-end);
\path[draw, decoration={border, amplitude=25pt, angle=90, segment length=4pt, mirror}, color=gray, decorate] ($ (start) + (-0.1, -0.1) $) -- ++(315:1);
\path[draw] (out) -- ++(30:1);
\path[draw] (start) -- ++(30:1);
\path[draw] (end) -- ++(315:1);
\end{tikzpicture}
\caption{$ β \text{ or } β' $ (C) \\ $ μ = β $ (A) \\ $ h_q = 0 $ \\ $ φ π_q = 0 $}
\end{subfigure}
\begin{subfigure}{0.24\linewidth}
\centering
\begin{tikzpicture}
\path[draw] (0, 0) coordinate (out) -- (0, 1) coordinate[pos=0.4] (m-start);
\path[draw] (0.1, -0.1) coordinate (start) -- ++(315:1) coordinate[pos=0.4] (m-end);
\path[draw] (-0.1, -0.1) coordinate (end) -- ++(225:1);
\path[draw, ->, bend right=110, looseness=3] (m-start) to (m-end);
\path[draw, decoration={border, amplitude=25pt, angle=90, segment length=4pt, mirror}, color=gray, decorate] ($ (start) + (-0.1, -0.1) $) -- ++(315:1);
\path[draw] (out) -- ++(60:1);
\path[draw] (start) -- ++(345:1);
\path[draw] (end) -- ++(150:1);
\end{tikzpicture}
\caption{$ β $ (A) \\ $ μ = γβ $ \\ $ h_q = 0 $ \\ $ φ π_q = 0 $}
\end{subfigure}
\begin{subfigure}{0.24\linewidth}
\centering
\begin{tikzpicture}
\path[draw] (0, 0) coordinate (out) -- (0, 1) coordinate[pos=0.4] (m-start);
\path[draw] (0.1, -0.1) coordinate (start) -- ++(315:1) coordinate[pos=0.4] (m-end);
\path[draw] (-0.1, -0.1) coordinate (end) -- ++(225:1);
\path[draw, ->, bend right=110, looseness=3] (m-start) to (m-end);
\path[draw, decoration={border, amplitude=25pt, angle=90, segment length=4pt, mirror}, color=gray, decorate] ($ (start) + (-0.1, -0.1) $) -- ++(315:1);
\path[draw] (out) -- ++(60:1);
\path[draw] (start) -- ++(345:1);
\path[draw] (end) -- ++(90:1);
\end{tikzpicture}
\caption{$ β $ (A) \\ $ μ = α_2 \text{ or } α_3 $ \\ $ h_q = \id $ (B) or $ 0 $ \\ $ φ π_q = α_3 + α_4 $}
\end{subfigure}
\begin{subfigure}{0.24\linewidth}
\centering
\begin{tikzpicture}
\path[draw] (0, 0) coordinate (out) -- (0, 1) coordinate[pos=0.4] (m-start);
\path[draw] (0.1, -0.1) coordinate (start) -- ++(315:1) coordinate[pos=0.4] (m-end);
\path[draw] (-0.1, -0.1) coordinate (end) -- ++(225:1);
\path[draw, ->, bend right=110, looseness=3] (m-start) to (m-end);
\path[draw, decoration={border, amplitude=25pt, angle=90, segment length=4pt, mirror}, color=gray, decorate] ($ (start) + (-0.1, -0.1) $) -- ++(315:1);
\path[draw] (out) -- ++(60:1);
\path[draw] (start) -- ++(345:1);
\path[draw] (end) -- ++(270:1);
\end{tikzpicture}
\caption{$ β $ (A) \\ $ μ = β $ (A) \\ $ h_q = 0 $ \\ $ φ π_q = 0 $}
\end{subfigure}
\begin{subfigure}{0.24\linewidth}
\centering
\begin{tikzpicture}
\path[draw] (0, 0) coordinate (out) -- (0, 1) coordinate[pos=0.4] (m-start);
\path[draw] (0.1, -0.1) coordinate (start) -- ++(315:1) coordinate[pos=0.4] (m-end);
\path[draw] (-0.1, -0.1) coordinate (end) -- ++(225:1);
\path[draw, ->, bend right=110, looseness=3] (m-start) to (m-end);
\path[draw, decoration={border, amplitude=25pt, angle=90, segment length=4pt, mirror}, color=gray, decorate] ($ (start) + (-0.1, -0.1) $) -- ++(315:1);
\path[draw] (out) -- ++(60:1);
\path[draw] (start) -- ++(345:1);
\path[draw] (end) -- ++(315:1);
\end{tikzpicture}
\caption{$ β $ (A) \\ $ μ = β $ (A) \\ $ h_q = 0 $ \\ $ φ π_q = 0 $}
\end{subfigure}
\caption{Possible result components of first-out disks in Kadeishvili trees}
\label{fig:components-first-out}
\end{figure}
\begin{figure}
\centering
\begin{subfigure}{0.24\linewidth}
\centering
\begin{tikzpicture}
\path[draw] (0, 0) coordinate (out) -- (0, 1) coordinate[pos=0.4] (m-end);
\path[draw] (0.1, -0.1) coordinate (start) -- ++(315:1);
\path[draw] (-0.1, -0.1) coordinate (end) -- ++(225:1) coordinate[pos=0.4] (m-start);
\path[draw, ->, bend right=110, looseness=3] (m-start) to (m-end);
\path[draw, decoration={border, amplitude=25pt, angle=90, segment length=4pt, mirror}, color=gray, decorate] ($ (start) + (-0.1, -0.1) $) -- ++(315:1);
\path[draw] (end) -- (out);
\path[draw] (start) -- ++(30:1);
\end{tikzpicture}
\caption{$ δ $ insertion \\ $ μ = βα $ \\ $ h_q = β $ (A) + E tail \\ $ φ π_q = $ G/H}
\end{subfigure}
\begin{subfigure}{0.24\linewidth}
\centering
\begin{tikzpicture}
\path[draw] (0, 0) coordinate (out) -- (0, 1) coordinate[pos=0.4] (m-end);
\path[draw] (0.1, -0.1) coordinate (start) -- ++(315:1);
\path[draw] (-0.1, -0.1) coordinate (end) -- ++(225:1) coordinate[pos=0.4] (m-start);
\path[draw, ->, bend right=110, looseness=3] (m-start) to (m-end);
\path[draw, decoration={border, amplitude=25pt, angle=90, segment length=4pt, mirror}, color=gray, decorate] ($ (start) + (-0.1, -0.1) $) -- ++(315:1);
\path[draw] (end) -- (out);
\path[draw] (start) -- ++(up:1);
\end{tikzpicture}
\caption{$ δ $ insertion \\ $ μ = α_1 $ or $ α_4 $ \\ $ h_q = 0 $ \\ $ φ π_q = 0 $ or $ α_3 + α_4 $}
\end{subfigure}
\begin{subfigure}{0.24\linewidth}
\centering
\begin{tikzpicture}
\path[draw] (0, 0) coordinate (out) -- (0, 1) coordinate[pos=0.4] (m-end);
\path[draw] (0.1, -0.1) coordinate (start) -- ++(315:1);
\path[draw] (-0.1, -0.1) coordinate (end) -- ++(225:1) coordinate[pos=0.4] (m-start);
\path[draw, ->, bend right=110, looseness=3] (m-start) to (m-end);
\path[draw, decoration={border, amplitude=25pt, angle=90, segment length=4pt, mirror}, color=gray, decorate] ($ (start) + (-0.1, -0.1) $) -- ++(315:1);
\path[draw] (end) -- (out);
\path[draw] (start) -- ++(down:1);
\end{tikzpicture}
\caption{$ δ $ insertion \\ $ μ = β $ (A) \\ $ h_q = 0 $ \\ $ φ π_q = 0 $}
\end{subfigure}
\begin{subfigure}{0.24\linewidth}
\centering
\begin{tikzpicture}
\path[draw] (0, 0) coordinate (out) -- (0, 1) coordinate[pos=0.4] (m-end);
\path[draw] (0.1, -0.1) coordinate (start) -- ++(315:1);
\path[draw] (-0.1, -0.1) coordinate (end) -- ++(225:1) coordinate[pos=0.4] (m-start);
\path[draw, ->, bend right=110, looseness=3] (m-start) to (m-end);
\path[draw, decoration={border, amplitude=25pt, angle=90, segment length=4pt, mirror}, color=gray, decorate] ($ (start) + (-0.1, -0.1) $) -- ++(315:1);
\path[draw] (end) -- (out);
\path[draw] (start) -- ++(225:1);
\end{tikzpicture}
\caption{$ δ $ insertion \\ $ μ = β/β' $ (C) \\ $ h_q = 0 $ \\ $ φ π_q = 0 $}
\end{subfigure}
\begin{subfigure}{0.24\linewidth}
\centering
\begin{tikzpicture}
\path[draw] (0, 0) coordinate (out) -- (0, 1) coordinate[pos=0.4] (m-end);
\path[draw] (0.1, -0.1) coordinate (start) -- ++(315:1);
\path[draw] (-0.1, -0.1) coordinate (end) -- ++(225:1) coordinate[pos=0.4] (m-start);
\path[draw, ->, bend right=110, looseness=3] (m-start) to (m-end);
\path[draw, decoration={border, amplitude=25pt, angle=90, segment length=4pt, mirror}, color=gray, decorate] ($ (start) + (-0.1, -0.1) $) -- ++(315:1);
\path[draw] (end) -- ++(150:1);
\path[draw] (out) -- ++(150:1);
\path[draw] (start) -- ++(30:1);
\end{tikzpicture}
\caption{$ β/β' $ (C) \\ $ μ = βα $ \\ $ h_q = β $ (A) + E tail \\ $ φ π_q = $ G/H}
\end{subfigure}
\begin{subfigure}{0.24\linewidth}
\centering
\begin{tikzpicture}
\path[draw] (0, 0) coordinate (out) -- (0, 1) coordinate[pos=0.4] (m-end);
\path[draw] (0.1, -0.1) coordinate (start) -- ++(315:1);
\path[draw] (-0.1, -0.1) coordinate (end) -- ++(225:1) coordinate[pos=0.4] (m-start);
\path[draw, ->, bend right=110, looseness=3] (m-start) to (m-end);
\path[draw, decoration={border, amplitude=25pt, angle=90, segment length=4pt, mirror}, color=gray, decorate] ($ (start) + (-0.1, -0.1) $) -- ++(315:1);
\path[draw] (end) -- ++(150:1);
\path[draw] (out) -- ++(150:1);
\path[draw] (start) -- ++(up:1);
\end{tikzpicture}
\caption{$ β/β' $ (C) \\ $ μ = α_1 $ or $ α_4 $ \\ $ h_q = 0 $ \\ $ φ π_q = 0 $ or $ α_3 + α_4 $}
\end{subfigure}
\begin{subfigure}{0.24\linewidth}
\centering
\begin{tikzpicture}
\path[draw] (0, 0) coordinate (out) -- (0, 1) coordinate[pos=0.4] (m-end);
\path[draw] (0.1, -0.1) coordinate (start) -- ++(315:1);
\path[draw] (-0.1, -0.1) coordinate (end) -- ++(225:1) coordinate[pos=0.4] (m-start);
\path[draw, ->, bend right=110, looseness=3] (m-start) to (m-end);
\path[draw, decoration={border, amplitude=25pt, angle=90, segment length=4pt, mirror}, color=gray, decorate] ($ (start) + (-0.1, -0.1) $) -- ++(315:1);
\path[draw] (end) -- ++(150:1);
\path[draw] (out) -- ++(150:1);
\path[draw] (start) -- ++(down:1);
\end{tikzpicture}
\caption{$ β/β' $ (C) \\ $ μ = β $ (A) \\ $ h_q = 0 $ \\ $ φ π_q = 0 $}
\end{subfigure}
\begin{subfigure}{0.24\linewidth}
\centering
\begin{tikzpicture}
\path[draw] (0, 0) coordinate (out) -- (0, 1) coordinate[pos=0.4] (m-end);
\path[draw] (0.1, -0.1) coordinate (start) -- ++(315:1);
\path[draw] (-0.1, -0.1) coordinate (end) -- ++(225:1) coordinate[pos=0.4] (m-start);
\path[draw, ->, bend right=110, looseness=3] (m-start) to (m-end);
\path[draw, decoration={border, amplitude=25pt, angle=90, segment length=4pt, mirror}, color=gray, decorate] ($ (start) + (-0.1, -0.1) $) -- ++(315:1);
\path[draw] (end) -- ++(150:1);
\path[draw] (out) -- ++(150:1);
\path[draw] (start) -- ++(225:1);
\end{tikzpicture}
\caption{$ β/β' $ (C) \\ $ μ = β $ (A) \\ $ h_q = 0 $ \\ $ φ π_q = 0 $}
\end{subfigure}
\begin{subfigure}{0.24\linewidth}
\centering
\begin{tikzpicture}
\path[draw] (0, 0) coordinate (out) -- (0, 1) coordinate[pos=0.4] (m-end);
\path[draw] (0.1, -0.1) coordinate (start) -- ++(315:1);
\path[draw] (-0.1, -0.1) coordinate (end) -- ++(225:1) coordinate[pos=0.4] (m-start);
\path[draw, ->, bend right=110, looseness=3] (m-start) to (m-end);
\path[draw, decoration={border, amplitude=25pt, angle=90, segment length=4pt, mirror}, color=gray, decorate] ($ (start) + (-0.1, -0.1) $) -- ++(315:1);
\path[draw] (end) -- ++(195:1);
\path[draw] (out) -- ++(120:1);
\path[draw] (start) -- ++(30:1);
\end{tikzpicture}
\caption{$ β $ (A) \\ $ μ = βα $ \\ $ h_q = β $ (A) + E tail \\ $ φ π_q = $ G/H}
\end{subfigure}
\begin{subfigure}{0.24\linewidth}
\centering
\begin{tikzpicture}
\path[draw] (0, 0) coordinate (out) -- (0, 1) coordinate[pos=0.4] (m-end);
\path[draw] (0.1, -0.1) coordinate (start) -- ++(315:1);
\path[draw] (-0.1, -0.1) coordinate (end) -- ++(225:1) coordinate[pos=0.4] (m-start);
\path[draw, ->, bend right=110, looseness=3] (m-start) to (m-end);
\path[draw, decoration={border, amplitude=25pt, angle=90, segment length=4pt, mirror}, color=gray, decorate] ($ (start) + (-0.1, -0.1) $) -- ++(315:1);
\path[draw] (end) -- ++(195:1);
\path[draw] (out) -- ++(120:1);
\path[draw] (start) -- ++(up:1);
\end{tikzpicture}
\caption{$ β $ (A) \\ $ μ = α_1 $ or $ α_4 $ \\ $ h_q = 0 $ \\ $ φ π_q = 0 $ or $ α_3 + α_4 $}
\end{subfigure}
\begin{subfigure}{0.24\linewidth}
\centering
\begin{tikzpicture}
\path[draw] (0, 0) coordinate (out) -- (0, 1) coordinate[pos=0.4] (m-end);
\path[draw] (0.1, -0.1) coordinate (start) -- ++(315:1);
\path[draw] (-0.1, -0.1) coordinate (end) -- ++(225:1) coordinate[pos=0.4] (m-start);
\path[draw, ->, bend right=110, looseness=3] (m-start) to (m-end);
\path[draw, decoration={border, amplitude=25pt, angle=90, segment length=4pt, mirror}, color=gray, decorate] ($ (start) + (-0.1, -0.1) $) -- ++(315:1);
\path[draw] (end) -- ++(195:1);
\path[draw] (out) -- ++(120:1);
\path[draw] (start) -- ++(down:1);
\end{tikzpicture}
\caption{$ β $ (A) \\ $ μ = β $ (A) \\ $ h_q = 0 $ \\ $ φ π_q = 0 $}
\end{subfigure}
\begin{subfigure}{0.24\linewidth}
\centering
\begin{tikzpicture}
\path[draw] (0, 0) coordinate (out) -- (0, 1) coordinate[pos=0.4] (m-end);
\path[draw] (0.1, -0.1) coordinate (start) -- ++(315:1);
\path[draw] (-0.1, -0.1) coordinate (end) -- ++(225:1) coordinate[pos=0.4] (m-start);
\path[draw, ->, bend right=110, looseness=3] (m-start) to (m-end);
\path[draw, decoration={border, amplitude=25pt, angle=90, segment length=4pt, mirror}, color=gray, decorate] ($ (start) + (-0.1, -0.1) $) -- ++(315:1);
\path[draw] (end) -- ++(195:1);
\path[draw] (out) -- ++(120:1);
\path[draw] (start) -- ++(225:1);
\end{tikzpicture}
\caption{$ β $ (A) \\ $ μ = β $ (A) \\ $ h_q = 0 $ \\ $ φ π_q = 0 $}
\end{subfigure}
\caption{Possible result components of final-out disks in Kadeishvili trees}
\label{fig:components-final-out}
\end{figure}

\begin{proof}
Let $ S $ be the set of all elementary morphisms of these types. We prove the claim by induction on the tree size. For an h-tree with just one node, leaf and root at the same time, there is nothing to show, since the result of this tree is the deformed basis cohomology element itself which only contains terms from $ S $.

Now for any arbitrarily large tree, output components are of the form
\begin{equation}
\label{eq:trees-root-shape}
h_q (μ^{≥3} (m_k, …, m_1)) \text{ and } h_q (μ^2 (m_2, m_1)),
\end{equation}
where by induction hypothesis each $ m_i $ is from $ S $, or may in case of $ μ^{≥3} $ also be a $ δ $ insertion. We will now check all possible terms that can occur in \eqref{eq:trees-root-shape}.

In the case of a discrete immersed disk $ μ^{≥3} (m_k, …, m_1) $, assume the disk is all-in. Then the $ μ^{≥3} $ result is an identity from situation B, C or D. Its image under $ h_q $ vanishes in all three cases.

Assume now the disk is first-out. Then the first morphism is by assumption one in $ S $. In particular, it is strictly smaller than one full turn. Since the disk is first-out, the first morphism $ m_1 $ consists of at least two indecomposable angles. We can now compute the $ μ^{≥3} $ result on a case-by-case basis, distinguishing after the type of $ m_1 $. The results are shown in \autoref{fig:components-first-out}. The figures omit the case of $ α_0' $ which is similar to that of an outer $ δ $ insertion. Hatching indicates the interior of the disk. Similarly, the results for final-out disks are shown in \autoref{fig:components-final-out}.

We conclude that for any discrete immersed disk, the output $ h_q (μ^{≥3} (m_k, …, m_1)) $ consists only of terms lying in $ S $. Next, we check the terms occurring in a simple composition $ μ^2 (m_2, m_1)) $. Table \ref{tab:components-multiplication} contains the results of such multiplications, “imp” denoting an impossible combination, hence vanishing product. The table also lists their images under $ h_q $, abbreviating tail terms $ β $ (A) as +E and tail terms $ \id $ (B) as +G. We conclude that as long as factors lie in $ S $, their image under $ h_q μ^2 $ also has components only in $ S $.
\end{proof}

In the product tables, the expression +G/H appears under the $ φ π_q $ key. We have used this abbreviation +G/H to denote the terms of $ π_q (βα) $ according to \autoref{th:deformed-codif-ba}. In other words, +G/H simply denotes tail terms of the form $ α_3 + α_4 $ (B) and $ \id $ (C). 

\subsection{Result components}
\label{sec:resultcomp-definition}
In this section, we introduce our notion of “result components”. The reason for this notion is that any evaluation $ h_q μ $ in a Kadeishvili tree may in principle yield a large number of terms. Any single of these terms may yield multiple terms again upon the next evaluation in the tree. The notion of result components serves to get grip on these terms. After the definition, we provide some terminology and a few examples.

Before we state the precise definition, let us illustrate the idea: Regard an evaluation task like computing $ (3x + 5y)(2x + 3y) $ or even $ (3x + 5y)(2x + 3y)(x - y) $. These evaluations can be represented by the trees

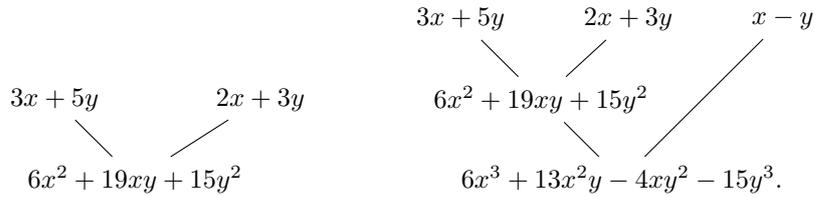
\begin{center}
\begin{tikzpicture}
\path[node distance=1.5cm] node (A) {$ 3x + 5y $} node[right of=A, anchor=west, shift={(0.5, 0)}] (B) {$ 2x + 3y $}
node[below right of=A] {$ 6x^2 + 19xy + 15y^2 $} edge (A) edge (B);
\end{tikzpicture}
\hspace{1cm}
\begin{tikzpicture}
\path[node distance=1.5cm] node (A) {$ 3x + 5y $} node[right of=A, anchor=west] (B) {$ 2x + 3y $} node[right of=B, anchor=west] (C) {$ x - y $}
node[below right of=A] (D) {$ 6x^2 + 19xy + 15y^2 $} edge (A) edge (B)
node[below right of=D] {$ 6x^3 + 13x^2 y - 4 x y^2 - 15 y^3 $.} edge (D) edge (C);
\end{tikzpicture}
\end{center}

The result expression $ 6x^2 + 19xy + 15y^2 $ is concise, but does not include information on how it was derived from the individual factors. The idea behind result components is to retain this information instead. For example, the left tree should have four result components:

\begin{center}
\begin{tikzpicture}
\path node (A) {$ 6x^2 $ derived from multiplying $ 3x $ and $ 2x $,} (A.east) node[anchor=west, shift={(1, 0)}] (B) {$ 9xy $ derived from multiplying $ 3x $ and $ 3y $,}
(A.west) node[below, anchor=north west, shift={(down:0.2)}] {$ 10xy $ derived from multiplying $ 5y $ and $ 2x $,}
(B.west) node[below, anchor=north west, shift={(down:0.2)}] {$ 15y^2 $ derived from multiplying $ 5y $ and $ 3y $.};
\end{tikzpicture}
\end{center}
In other words, the left tree has four distinct result components, even though the result can be abbreviated to only three terms. The tree on the right has just four result terms, while there are eight distinct result components. For example, one of these eight result components consists of the choice of $ 3x $ on the leftmost leaf, $ 2x $ on the middle leaf and $ -y $ on the rightmost leaf.

Let us prepare for result components of Kadeishvili trees: In contrast to the simple multiplication trees above, the leaves of a Kadeishvili tree are labeled by deformed cohomology basis elements. These elements consist of a finite or infinite amount of \emph{additive components}. For example, let $ α_3 + α_4 $ denote a certain cohomology basis element and assume both $ α_3 $ and $ α_4 $ have a tail each consisting of one type E disk with $ β $ morphisms denoted $ β_1 $ and $ β_2 $ respectively. Then the deformed cohomology basis element reads $ α_3 + α_4 + β_1 + β_2 $ and is defined to have four distinct additive components, even though technically it may be possible that $ β_1 = β_2 $.

Similarly, node evaluations in Kadeisvhili trees may yield a large amount of \emph{additive components}. For example, an evaluation $ h_q (βα) $ may yield an expression like $ β + β_1 + β_2 + \id \text{(B)} $ according to \autoref{th:deformed-codif-ba}. This evaluation is defined to have four distinct additive components. In other words, an additive component always refers to one of the main terms or a choice of one of the tail terms. We are finally ready to define result components of Kadeishvili trees:

\begin{definition}
\label{def:resultcomp-def}
The \emph{restriction} of an h-tree or π-tree $ (T, h_1, …, h_N) $ at a non-root node $ P ∈ T $ is the h-tree from $ P $ up to all leaves, together with the corresponding subset of $ (h_1, …, h_N) $. A \emph{result component} of an h-tree or π-tree is defined inductively as follows:
\begin{itemize}
\item A result component of an h-tree with only one node consists of an additive component appearing in the corresponding $ h_1 $.
\item A result component of an h-tree with at least three nodes consists of result components $ r_1, …, r_k $ of the restrictions at the ordered children of the root, together with choices $ n_0, …, n_k ≥ 0 $ of $ δ $ insertions, and an additive component appearing in
\begin{equation*}
h_q μ_q (\underbrace{δ, …, δ}_{n_k}, r_k, …, r_1, \underbrace{δ, …, δ}_{n_0}).
\end{equation*}
\item A result component of a π-tree consists of result components $ r_1, …, r_k $ of the restrictions at the ordered children of the root, together with choices $ n_0, …, n_k ≥ 0 $ of $ δ $ insertions, and an additive component appearing in
\begin{equation*}
φ π_q μ_q (\underbrace{δ, …, δ}_{n_k}, r_k, …, r_1, \underbrace{δ, …, δ}_{n_0}).
\end{equation*}
\item For π-trees, the result components $ (-1)^{\#α_3} α_3 $ and $ (-1)^{\#α_4 + 1} α_4 $ shall be grouped together as one result component. Also, the result components $ \id $ (D) shall be grouped together as one result component.
\item Additive components arising from different tail nodes in the evaluation of $ h_q $ or $ φ π_q $ shall be kept distinct as result components.
\end{itemize}
A \emph{direct morphism} is a result component of a one-node h-tree. A result component \emph{derives} from the result components $ r_1, …, r_k $, and from all the result components they derive from themselves. A direct morphism derives from nothing. A \emph{tail result component} is one that comes from a tail additive component of the final $ h_q $ or $ φ π_q $ evaluation. Tail additive components of direct morphisms are also counted as tail result components. Any other result component is a \emph{main result component}. The class of all result components of π-trees is denoted $ \PiTr $.
\end{definition}

\begin{example}
Regard an evaluation $ h_q (βα) $. Its result consists of a tower of $ β $ (A) morphisms. They need not be distinct as morphisms, but shall be treated as distinct result components. The $ β $ (A) with the lowest $ q $ power is the main result component, all others are tail result components.
\end{example}

\begin{figure}
\centering
\begin{tikzpicture}
\path node (A) {$ α_4 $} node[right of=A] (B) {$ β $(A)} node[right of=B] (C) {$ α_4 $} node[right of=C] (D) {$ α_4 $}
node[below right of=B] {$ φ π_q μ^9 = \id $ (D)} edge (A) edge (B) edge (C) edge (D);
\end{tikzpicture}
\hspace{1cm}
\begin{tikzpicture}[scale=0.6, rotate=-90]
\path[draw, <-] (0, 0) -- ++(60:1) coordinate (left-down) -- ++(right:2) -- ++(135:1.2) -- ++(45:1) -- ++(132:1.3) coordinate (right);
\path[draw] (left-down) ++(up:0.2) -- ++(right:2) -- ++(300:1);
\path[draw] (left-down) ++(up:0.2) -- ++(45:1) -- ++(135:1) -- ($ (right) + (left:0.1) $) coordinate (left);
\path[draw, ->] (left) -- ++(left:1.3) coordinate[pos=0.3] (beta-end) -- ++(45:1) -- ++(135:1.5) coordinate[pos=0.7] (1-start) coordinate (1) -- ++(45:1) -- ++(135:1);
\path[draw, ->] (1) ++(right:0.2) -- ++(45:1) coordinate[pos=0.3] (1-end) coordinate[pos=0.5] (2-start) coordinate (2) -- ++(right:1.5) -- ++(45:1);
\path[draw, ->] (2) ++(down:0.2) -- ++(right:1.5) coordinate[midway] (switch) coordinate[pos=0.3] (2-end) coordinate[pos=0.7] (3-start) coordinate (3) -- ++(315:1);
\path[draw] (3) ++(down:0.2) -- ++(315:1) coordinate[pos=0.3] (3-end) -- ++(225:1) -- ++(315:1) -- (right) coordinate[pos=0.7] (beta-start);
\path[draw, ->, bend right=80, looseness=1.5] (beta-start) to node[midway, right] {\small $ β $} (beta-end);
\path[draw, ->, bend right=60] (1-start) to node[pos=0, left] {\small $ α_4 $} (1-end);
\path[draw, ->, bend right=60] (2-start) to node[midway, left] {\small $ α_4 $} (2-end);
\path[draw, ->, bend right=60] (3-start) to node[midway, left] {\small $ α_4 $} (3-end);
\path (switch) node[left] {\small $ L $};
\path[draw] (switch) ++(up:0.3) -- ++(down:0.4);
\path[draw] (switch) ++(0.05, 0.3) -- ++(down:0.4);
\end{tikzpicture}
\caption{A π-tree with a concrete result component}
\label{fig:resultcomp-definition-example}
\end{figure}
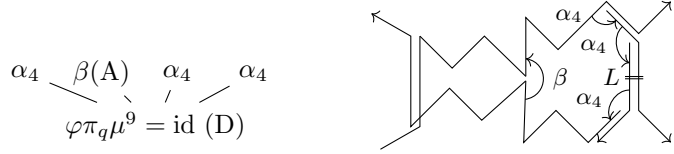

\begin{example}
A sample π-tree with a concrete result component is depicted in \autoref{fig:resultcomp-definition-example}. The inputs of this π-tree are four $ β_3 + β_4 $ morphisms, corresponding to the four intersection points between the zigzag curves. One of the zigzag paths is denoted $ L $. The first input morphism departs from $ L $ and the fourth ends on $ L $. The angles depicted are the main $ α_4 $ components of the first, third and fourth input, as well as the first tail component of the second input. The double stroke on the rightmost arrow indicated the separation between the first and the fourth morphisms of the sequence. In case the rightmost arc is the identity location of the first/final zigzag path, then the output of the π-tree is the identity. Otherwise the output vanishes. This identity result component is derived from the main components of the first, the first tail component of the second, and the main components of the third and fourth input morphisms, with $ n_1 = 3 $ many $ δ $'s after the first morphism and $ n_2 = 2 $ many $ δ $'s after the second morphism. This example illustrates a nontrivial result component of a π-tree and shows how tail components lead to results. In contrast, no single result component of this π-tree derives from the main components of all four input morphisms.
\end{example}

\subsection{Classification of result components}
\label{sec:resultcomp-classification}
In this section, we provide a semi-explicit, inductive characterization of result components of Kadeishvili trees. To understand what this means, recall from \autoref{sec:resultcomp-possible} that only certain types of morphisms can appear as result components. For each of these types, we will describe all possible in which it is derived from simpler result components. This description is recursive, and has to remain so until we match result components with pieces of smooth immersed disks later on.

\begin{table}
\centering
\begin{tabular}{@{}c@{\hspace{0.05\linewidth}}c@{\hspace{0.05\linewidth}}p{0.5\linewidth}@{}}
\textbf{Type} & \textbf{Cue} & \textbf{Possible ways of derivation} \\\hline
$ α_0 $ & h & direct \\
$ \id $ (C) & h & direct \\
$ α_4 $ & h & direct \\
$ α_3 $ & h & direct or \autoref{fig:subdisk-alpha3} \\
$ α_0' $ & h & direct or \autoref{fig:subdisk-alpha0p} \\
$ β/β' $ (C) & h & direct or \autoref{fig:subdisk-betaC} or \ref{fig:subdisk-betaCp} \\
$ β $ (A) & direct & tail of some $ α_3 $, $ α_4 $, $ β $ (C) or $ β' $ (C) \\
$ β $ (A) & main $ h_q μ^{≥3} $ & final-out disk, with final morphism an outer $ δ $ insertion, $ β $ (A), $ β/β' $ (C) or $ α_0' $ \\
$ β $ (A) & tail $ h_q μ $ & comes with corresponding main result component $ β $ (A), $ β/β' $ (C) or $ α_3 $, example see \autoref{fig:subdisk-betaA-tail} \\
$ β $ (A) & main $ h_q μ^2 $ & \autoref{fig:subdisk-betaA-mult-main} \\
$ \id $ (D) & h & \autoref{fig:subdisk-idD-1} or \ref{fig:subdisk-idD-2} \\
$ \id $ (B) & h & \autoref{fig:components-first-out} or \ref{fig:subdisk-idB} \\
$ α_3 + α_4 $ & main $ φ π_q μ^{≥3} $ & \autoref{fig:components-first-out} or \ref{fig:components-final-out} \\
$ α_3 + α_4 $ & main $ φ π_q μ^2 $ & \autoref{fig:subdisk-alpha34-mult-main} or \ref{fig:subdisk-degenerate} \\
$ α_3 + α_4 $ & tail $ φ π_q μ $ & tail of a certain $ φ π_q (βα) $, with $ βα $ itself being a $ μ^2 $ product or a disk of \autoref{fig:components-final-out} \\
$ \id $ (C) & main $ φ π_q μ^2 $ & \autoref{fig:subdisk-degenerate} \\
$ \id $ (C) & main $ φ π_q μ^{≥3} $ & all-in disk of type H, whose inner morphisms may be $ δ $ insertions, $ β $ (A), $ α_3 $ (B), $ α_4 $ (B), $ β/β' $ (C), $ α_0 $ (D), $ α_0' $ (D), example see \autoref{fig:subdisk-idC-disk-main} \\
$ \id $ (C) & tail $ φ π_q μ $ & tail of a certain $ φ π_q (βα) $ evaluation \\
$ \id $ (D) & π & \autoref{fig:subdisk-idD} \\
$ α_0 $ & π & \autoref{fig:subdisk-alpha0}
\end{tabular}
\caption{Classification of result components}
\label{tab:resultcomp-classification}
\end{table}

\begin{figure}
\begin{subfigure}{0.5\linewidth}
\centering
\begin{tikzpicture}
\path node (alpha0) {$ α_0 $} node[right of=alpha0] (alpha4) {$ α_4 $}
node[below right of=alpha0] {$ h_q μ^2 = α_3 $} edge (alpha0) edge (alpha4);
\end{tikzpicture}
\begin{tikzpicture}
\path[draw] (-0.2, 0) -- ++(right:1.2) coordinate[pos=0.3] (switch-target) coordinate[pos=0.6] (alpha0-end) coordinate[pos=0.3] (L2-bot) -- ++(up:1.5) coordinate[pos=0.4] (alpha4-end) coordinate[pos=0.5] (L1-mid) -- ++(right:1);
\path[draw] (1.9, -0.1) -- ++(left:1) coordinate[pos=0.5] (alpha4-start) -- ++(up:1.5) coordinate[pos=0.4] (alpha0-start) coordinate[pos=0.5] (L2-mid) -- ++(left:1);
\path[draw, fill] ($ 0.5*(L1-mid) + 0.5*(L2-mid) $) circle (1pt) -- ($ 0.125*(L1-mid) + 0.125*(L2-mid) + 0.75*(L2-bot) $) coordinate[pos=0.66] (switch);
\path[draw, ->, bend right=45] (alpha0-start) to node[near start, below] {$ α_0 $} (alpha0-end);
\path[draw, ->, bend right=45] (alpha4-start) to node[near end, below] {$ α_4 $} (alpha4-end);
\path[draw, fill] (switch) circle(1pt);
\path[draw] (switch) to[out=220, in=90] (switch-target);
\end{tikzpicture}
\caption{Tree with $ α_3 $ as result component}
\label{fig:subdisk-alpha3}
\end{subfigure}
\begin{subfigure}{0.49\linewidth}
\centering
\begin{tikzpicture}
\path node (in1) {$ α_0 $} node[right of=in1] (dots) {…} node[right of=dots] (in2) {$ α_0 $} node[right of=in2] (in3) {$ α_0' $}
node[below left of=in3] (m1) {} edge (in2) edge (in3)
node[below left of=m1] {$ h_q μ^2 = α_0' $} edge (in1) edge[densely dotted, thick] (m1);
\path[draw, decoration=brace, decorate] (in1.north) to node[midway, above, shift={(0, 0.2)}] {$ ≥ 0 $} (in2.north);
\end{tikzpicture}
\begin{tikzpicture}
\path[draw] (0, 0) -- ++(60:1.5) coordinate[pos=0.4] (L1-target) coordinate[pos=0.5] (L2-target) coordinate[pos=0.6] (L3-target) coordinate[pos=0.62] (L4-target) coordinate[pos=0.7] (alpha0p-end) -- ++(300:1.5) coordinate[pos=0.6] (L1-source) coordinate[pos=0.3] (alpha0p-start);
\path[draw] (L1-source) to coordinate[pos=0.5] (switch1) ($ 0.5*(L1-source) + 0.5*(L1-target) $);
\path[draw] (switch1) to[out=150, in=0] coordinate[pos=0.66] (switch2) ($ 0.5*(switch1) + 0.5*(L2-target) $);
\path[draw] (switch2) to[out=150, in=0] coordinate[pos=0.66] (switch3) ($ 0.25*(switch2) + 0.75*(L3-target) $);
\path[draw] (switch3) to[out=150, in=10] (L4-target);
\path[fill] (switch1) circle (1pt);
\path[fill] (switch2) circle (1pt);
\path[fill] (switch3) circle (1pt);
\path[draw] (switch1) to[out=330, in=180] ++(0.15, -0.1);
\path[draw] (switch2) to[out=320, in=180] ++(0.25, -0.25);
\path[draw] (switch3) to[out=320, in=180] ++(0.25, -0.25);
\path[draw, ->, bend right=150, looseness=12] (alpha0p-start) to node[midway, below] {$ α_0' $} (alpha0p-end);
\end{tikzpicture}
\caption{Trees with $ α_0' $ as result component}
\label{fig:subdisk-alpha0p}
\end{subfigure}
\begin{subfigure}{0.49\linewidth}
\centering
\begin{tikzpicture}
\path node (mid1) {$ \id $ (C)} node[right of=mid1] (rin1) {$ α_0 $} node[right of=rin1] (dots2) {…} node[right of=dots2] (rin2) {$ α_0 $} node[right of=rin2] (right) {$ α_0' $}
node[below left of=right] (m1) {$ h_q μ^2 $} edge (rin2) edge (right)
node[below left of=m1, shift={(-0.2, -0.2)}] (m2) {$ h_q μ^2 $} edge (rin1) edge[densely dotted, thick] (m1)
node[below left of=m2] (m3) {$ h_q μ^2 = β $ (C)} edge (mid1) edge (m2);
\path[draw, decoration=brace, decorate] (rin1.north) to node[midway, above, shift={(0, 0.2)}] {$ ≥ 0 $} (rin2.north);
\begin{scope}[shift={(5, -2.3)}]
\path[draw] (0, 0) -- ++(up:1.5) coordinate[pos=0.55] (L4-source) -- ++(230:1.5) coordinate[pos=0.3] (beta-end) coordinate[pos=0.5] (L4-target);
\path[draw] (0.1, 0) -- ++(up:1.5) coordinate[pos=0.35] (L1-target) coordinate[pos=0.5] (L2-target) coordinate[pos=0.55] (L3-target) -- ++(310:1.5) coordinate[pos=0.3] (beta-start) coordinate[pos=0.5] (L1-source);
\path[draw] (L1-source) -- ($ 0.75*(L1-target) + 0.25*(L1-source) $) coordinate[pos=0.44] (switch1);
\path[draw] (switch1) to[out=150, in=25] coordinate[pos=0.66] (switch2) ($ 0.25*(switch1) + 0.75*(L2-target) $);
\path[draw] (switch2) to[out=150, in=15] ($ 0.5*(L4-source) + 0.5*(L3-target) $) coordinate (mid);
\path[draw] (mid) -- (L4-target);
\path[fill] (switch1) circle (1pt);
\path[fill] (switch2) circle (1pt);
\path[fill] (mid) circle (1pt);
\path[draw, ->, bend right=135, looseness=6] (beta-start) to node[midway, below] {$ β $ (C)} (beta-end);
\path[draw] (switch1) to[out=330, in=180] ++(0.15, 0);
\path[draw] (switch2) to[out=320, in=180] ++(0.25, -0.1);
\end{scope}
\end{tikzpicture}
\caption{Trees with $ β $ (C) as result component.}
\label{fig:subdisk-betaC}
\end{subfigure}
\begin{subfigure}{0.49\linewidth}
\centering
\begin{tikzpicture}
\path node (in1) {$ α_0 $} node[right of=in1] (dots1) {…} node[right of=dots1] (in2) {$ α_0 $} node[right of=in2] (mid1) {$ α_0' $} node[right of=mid1] (mid2) {$ \id $ (C)}
node[below right of=in2] (m1) {$ h_q μ^2 $} edge (in2) edge (mid1)
node[below left of=m1, shift={(0, -0.2)}] (m2) {$ h_q μ^2 $} edge[densely dotted, thick] (m1) edge (in1)
node[below right of=m2] (m5) {$ h_q μ^2 = β' $ (C)} edge (m2) edge (mid2);
\path[draw, decoration={brace}, decorate] (in1.north) to node[midway, above, shift={(0, 0.2)}] {$ ≥ 0 $} (in2.north);
\begin{scope}[shift={(5, -2.3)}]
\path[draw] (0, 0) -- ++(up:1.5) coordinate[pos=0.55] (L2-source) -- ++(230:1.5) coordinate[pos=0.3] (beta-end) coordinate[pos=0.45] (L4-target) coordinate[pos=0.5] (L3-target) coordinate[pos=0.6] (L2-target);
\path[draw] (0.1, 0) -- ++(up:1.5) coordinate[pos=0.55] (L1-target) -- ++(310:1.5) coordinate[pos=0.3] (beta-start) coordinate[pos=0.5] (L1-source);
\path[draw] ($ 0.5*(L1-target) + 0.5*(L2-source) $) coordinate (mid) -- ($ 0.75*(L2-target) + 0.25*(L2-source) $) coordinate[pos=0.44] (switch1);
\path[draw] (switch1) to[out=150, in=25] coordinate[pos=0.66] (switch2) ($ 0.25*(switch1) + 0.75*(L3-target) $);
\path[draw] (switch2) to[out=150, in=0] (L4-target);
\path[draw] (L1-source) -- (mid);
\path[fill] (switch1) circle (1pt);
\path[fill] (switch2) circle (1pt);
\path[fill] (mid) circle (1pt);
\path[draw, ->, bend right=135, looseness=6] (beta-start) to node[midway, below] {$ β' $ (C)} (beta-end);
\path[draw] (switch1) to[out=330, in=180] ++(0.15, -0.1);
\path[draw] (switch2) to[out=320, in=180] ++(0.25, -0.25);
\end{scope}
\end{tikzpicture}
\caption{Trees with $ β' $ (C) as result component}
\label{fig:subdisk-betaCp}
\end{subfigure}
\caption{Classification of $ α_3 $, $ α_0' $, $ β $ (C) and $ β' $ (C) result components}
\end{figure}
\begin{figure}
\begin{subfigure}{0.99\linewidth}
\begin{tikzpicture}
\path node (idC) {$ \id $ (C)} node[right of=idC] (in1) {$ α_0 $} node[right of=in1] (dots) {…} node[right of=dots] (in2) {$ α_0 $} node[right of=in2] (in3) {$ α_0' $} node[right of=in3] (idC2) {$ \id $ (C)} node[right of=idC2] (dots2) {…} node[right of=dots2] (extra-alpha0) {$ α_0 $} node[right of=extra-alpha0] (dots3) {…} node[right of=dots3] (extra-idC) {$ \id $ (C)} node[right of=extra-idC] {…}
node[below right of=in2] (m1) {$ α_0' $} edge (in2) edge (in3)
node[below left of=m1] (m2) {$ α_0' $} edge[densely dotted, thick] (m1) edge (in1)
node[below left of=m2] (m3) {$ β $ (C)} edge (idC) edge (m2)
node[below right of=m3] (m4) {$ β $ (A)} edge (m3) edge (idC2)
node[below right of=m4] (m5) {$ β $ (A)} edge[densely dotted, thick] (m4) edge (extra-alpha0)
node[below right of=m5] (m6) {$ h_q μ^2 = β $ (A)} edge[densely dotted, thick] (m5) edge (extra-idC);
\path[draw, decoration=brace, decorate] (in1.north) to node[midway, above, shift={(0, 0.1)}] {$ ≥ 0 $} (in2.north);
\path[draw, ultra thick, gray!50] ($ (m4.south east) + (0, -0.15) $) to[bend left=30] ($ (idC.north west) + (-0.3, 0) $) to[bend left=30] (idC2.north east) to[bend left=20] ($ (m4.south east) + (0, -0.15) $);
\path[draw] ($ (extra-idC)!0.5!(m6) + (2, 0) $) -- ++(30:2) coordinate[pos=0.5] (4-target) -- ++(240:2) coordinate[pos=0.5] (4-source) coordinate (4-end);
\path[draw] (4-end) ++(330:0.25) -- ++(60:2) coordinate[pos=0.5] (3-target) -- ++(300:2) coordinate[pos=0.5] (3-source) coordinate (3-end);
\path[draw] (3-end) ++(30:0.25) -- ++(120:2) coordinate[pos=0.5] (2-target) -- ++(330:2) coordinate[pos=0.5] (2-source) coordinate (2-end);
\path[bend left] ($ 0.5*(2-target) + 0.5*(3-source) $) to coordinate[pos=0.33] (3-1) coordinate[pos=0.66] (3-2) coordinate[pos=0] (mid23) coordinate[pos=1] (mid34) ($ 0.5*(3-target) + 0.5*(4-source) $);
\path[draw, semithick] (mid23) -- (3-1);
\path[draw, semithick, bend right] (3-1) to (3-2);
\path[draw, semithick, bend right] (3-2) to (mid34);
\path[fill] (mid23) circle(1.5pt) (3-1) circle(1.5pt) (3-2) circle(1.5pt) (mid34) circle(1.5pt);
\path[draw, semithick] (2-source) -- (mid23);
\path[draw, semithick] (mid34) -- (4-target);
\end{tikzpicture}
\end{subfigure}
\begin{subfigure}{0.99\linewidth}
\begin{tikzpicture}
\path node (in1) {$ α_0 $} node[right of=in1] (dots) {…} node[right of=dots] (in2) {$ α_0 $} node[right of=in2] (in3) {$ α_0' $} node[right of=in3] (idC) {$ \id $ (C)} node[right of=idC] (idC2) {$ \id $ (C)} node[right of=idC2] (dots2) {…} node[right of=dots2] (extra-alpha0) {$ α_0 $} node[right of=extra-alpha0] (dots3) {…} node[right of=dots3] (extra-idC) {$ \id $ (C)} node[right of=extra-idC] {…}
node[below right of=in2] (m1) {$ α_0' $} edge (in2) edge (in3)
node[below left of=m1] (m2) {$ α_0' $} edge[densely dotted, thick] (m1) edge (in1)
node[below right of=m2] (m3) {$ β' $ (C)} edge (m2) edge (idC)
node[below right of=m3] (m4) {$ β $ (A)} edge (m3) edge (idC2)
node[below right of=m4] (m5) {$ β $ (A)} edge[densely dotted, thick] (m4) edge (extra-alpha0)
node[below right of=m5] (m6) {$ h_q μ^2 = β $ (A)} edge[densely dotted, thick] (m5) edge (extra-idC);
\path[draw, ultra thick, gray!50] ($ (m4.south east) + (0, -0.15) $) to[bend left=30] ($ (in1.north west) + (-0.3, 0) $) to[bend left=30] (idC2.north east) to[bend left=20] ($ (m4.south east) + (0, -0.15) $);
\path[draw, decoration=brace, decorate] (in1.north) to node[midway, above, shift={(0, 0.1)}] {$ ≥ 0 $} (in2.north);
\path[draw] ($ (extra-idC)!0.5!(m6) + (2, 0) $) -- ++(30:2) coordinate[pos=0.5] (3-target) -- ++(down:2) coordinate[pos=0.5] (3-source) coordinate (3-end);
\path[draw] (3-end) ++(0.25, 0) -- ++(up:2) coordinate[pos=0.5] (2-target) -- ++(300:2) coordinate[pos=0.5] (2-source) coordinate (2-end);
\path[draw] (2-end) ++(30:0.25) -- ++(120:2) coordinate[pos=0.5] (1-target) -- ++(330:2) coordinate[pos=0.5] (1-source);
\path[bend left] ($ 0.5*(2-target) + 0.5*(3-source) $) to coordinate[pos=0.33] (3-1) coordinate[pos=0.66] (3-2) coordinate[pos=0] (mid23) (3-target);
\path[draw, semithick] (1-source) -- ($ 0.5*(1-target) + 0.5*(2-source) $) coordinate (mid12);
\path[draw, semithick] (mid12) -- (mid23);
\path[draw, semithick] (mid23) to (3-1);
\path[draw, semithick, bend right] (3-1) to (3-2);
\path[draw, semithick, bend right] (3-2) to (3-target);
\path[fill] (mid12) circle(1.5pt) (mid23) circle(1.5pt) (3-1) circle(1.5pt) (3-2) circle(1.5pt);
\end{tikzpicture}
\end{subfigure}
\begin{subfigure}{0.6\linewidth}
\begin{tikzpicture}
\path node (topdots1) {…} node[right of=topdots1] (topdots2) {…}
node[below right of=topdots1] (in1) {$ β $ (A)} edge (topdots1.west) edge ($ 0.5*(topdots1.east) + 0.5*(topdots2.west) $) edge (topdots2.east)
node[right of=in1] (dots1) {…} node[right of=dots1] (idC) {$ \id $ (C)} node[right of=idC] (dots2) {…} node[right of=dots2] (alpha0) {$ α_0 $} node[right of=alpha0] (dots3) {…}
node[below right of=in1] (m1) {$ β $ (A)} edge (in1) edge (idC)
node[below right of=m1] (m2) {$ β $ (A)} edge[densely dotted, thick] (m1)
node[below right of=m2] {$ h_q μ^2 = β $ (A)} edge (m2) edge (alpha0);
\path[draw, decoration=brace, decorate] ($ (dots1.north) + (0, 0.2) $) to node[midway, above, shift={(0, 0.1)}] {$ ≥ 1 $} ($ (dots3.north) + (0, 0.2) $);
\end{tikzpicture}
\end{subfigure}
\begin{subfigure}{0.39\linewidth}
\begin{tikzpicture}
\path[draw] (0, 0) -- ++(300:2) coordinate[pos=0.5] (4-target) coordinate (4-tip) -- ++(150:2) coordinate[pos=0.5] (4-source) coordinate (4-end);
\path[draw] (4-tip) ++(0.1, -0.4) -- ++(240:2) -- ++(60:2) coordinate[pos=0.5] (3-target) -- ++(300:2) coordinate[pos=0.5] (3-source) coordinate (3-end);
\path[draw] (3-end) ++(30:0.25) -- ++(120:2) coordinate[pos=0.5] (2-target) -- ++(330:2) coordinate[pos=0.5] (2-source) coordinate (2-end);
\path[draw] (2-end) ++(60:0.25) -- ++(150:2) coordinate[pos=0.5] (1-target) -- ++(30:2) coordinate[pos=0.5] (1-source) coordinate (1-end);
\path[bend left] (1-source) to coordinate[pos=0.33] (1-1) coordinate[pos=0.66] (1-2) coordinate[pos=1] (mid12) ($ 0.5*(1-target) + 0.5*(2-source) $);
\path[bend left] ($ 0.5*(2-target) + 0.5*(3-source) $) to coordinate[pos=0.33] (3-1) coordinate[pos=0.66] (3-2) coordinate[pos=0] (mid23) coordinate[pos=1] (mid34) (3-target);
\path[draw, semithick] (1-source) -- (1-1);
\path[draw, semithick, bend right] (1-1) to (1-2);
\path[draw, semithick, bend right] (1-2) to (mid12);
\path[fill] (1-1) circle(1.5pt) (1-2) circle(1.5pt) (mid12) circle(1.5pt);
\path[draw, semithick] (mid23) -- (3-1);
\path[draw, semithick, bend right] (3-1) to (3-2);
\path[draw, semithick, bend right] (3-2) to (mid34);
\path[fill] (mid23) circle(1.5pt) (3-1) circle(1.5pt) (3-2) circle(1.5pt);
\path[draw, semithick] (mid12) -- (mid23);
\path[draw, semithick] (3-target) -- ++(-0.5, 0) node[left] {…};
\path[draw, semithick] (4-target) -- ($ (4-target)!1.75!(4-source) $) node[below] {…};
\end{tikzpicture}
\end{subfigure}
\caption{Trees with $ β $ (A) as result component. In the first and second tree, the framed part is essential and further $ \id $ (C) and $ α_0 $ inputs are optional. In the third tree, at least one $ \id $ (C) or $ α_0 $ is required and further ones are optional. The $ β $ (A) on the left is supposed to be a direct, $ h_q μ^{≥3} $ or tail $ h_q μ^2 $ result component.}
\label{fig:subdisk-betaA-mult-main}
\end{figure}

\begin{figure}
\centering
\begin{tikzpicture}
\path node (topdots) {…}
node[below of=topdots] (betaA) {$ β $ (A)} edge (topdots.west) edge (topdots.east) node[right of=betaA] (idC) {$ \id $ (C)}
node[below right of=betaA] {$ φ π_q μ^2 = α_3 + α_4 $} edge (idC) edge (betaA);
\end{tikzpicture}
\begin{tikzpicture}
\path[draw, <-] (0, 0) to coordinate[pos=0.4] (L1-left) coordinate[pos=0.5] (idC-end) ++(330:1.1) to coordinate[pos=0.6] (L1-right) coordinate[at end] (L1-target) ++(30:1.1);
\path[draw, <-] (-0.1, -0.1) to coordinate[pos=0.4] (L3-top) coordinate[pos=0.8] (beta-end) coordinate[pos=0.6] (idC-start) ++(330:1) to coordinate[midway] (L3-bottom) ++(210:1);
\path[draw] (L1-target) ++(0.1, -0.1) -- ++(210:1) coordinate[pos=0.8] (beta-start) coordinate[pos=0.4] (L2-top) -- ++(330:1) coordinate[midway] (L2-bottom);
\path ($ (L3-top)!0.5!(L1-left) $) coordinate (mid-left);
\path ($ (L2-top)!0.5!(L1-right) $) coordinate (mid-right);
\path[draw, rounded corners, semithick] ($ (L3-bottom) + (0.1, -0.3) $) node[below] {…} to[bend left=10] (mid-left) to[bend left=45] (mid-right) to[bend left=10] ($ (L2-bottom) + (-0.1, -0.3) $) node[below] {…};
\path[fill] (mid-left) circle[radius=0.05] node[above, shift={(0, 0.2)}] {out};
\path[fill] (mid-right) circle[radius=0.05];
\path[draw, bend right=75, ->] (beta-start) to (beta-end);
\path[draw, ->] (idC-start) to (idC-end);
\end{tikzpicture}
\begin{tikzpicture}
\path node (topdots) {…}
node[below of=topdots] (betaA) {$ β $ (A)} edge (topdots.west) edge (topdots.east) node[left of=betaA] (idC) {$ \id $ (C)}
node[below right of=idC] {$ φ π_q μ^2 = α_3 + α_4 $} edge (idC) edge (betaA);
\end{tikzpicture}
\begin{tikzpicture}
\path[draw] (0, 0) to coordinate[pos=0.4] (L1-left) coordinate[at end] (L1-mid) coordinate[pos=0.5] (idC-end) ++(330:1.1);
\path[draw, <-] (L1-mid) to coordinate[pos=0.6] (L1-right) coordinate[at end] (L1-target) ++(30:1.1);
\path[draw] (-0.1, -0.1) to coordinate[pos=0.4] (L3-top) coordinate[pos=0.6] (idC-start) coordinate[pos=0.8] (beta-end) ++(330:1) to coordinate[midway] (L3-bottom) ++(210:1);
\path[draw, ->] (L1-target) ++(0.1, -0.1) to coordinate[pos=0.4] (L2-top) coordinate[at end] (L2-mid) coordinate[pos=0.8] (beta-start) ++(210:1);
\path[draw] (L2-mid) -- ++(330:1) coordinate[midway] (L2-bottom);
\path ($ (L3-top)!0.5!(L1-left) $) coordinate (mid-left);
\path ($ (L2-top)!0.5!(L1-right) $) coordinate (mid-right);
\path[draw, rounded corners, semithick] ($ (L3-bottom) + (0.1, -0.3) $) node[below] {…} to[bend left=10] (mid-left) to[bend left=45] (mid-right) to[bend left=10] ($ (L2-bottom) + (-0.1, -0.3) $) node[below] {…};
\path[fill] (mid-left) circle[radius=0.05];
\path[fill] (mid-right) circle[radius=0.05] node[above, shift={(0, 0.2)}] {out};
\path[draw, ->, bend right=75] (beta-start) to (beta-end);
\path[draw, ->] (idC-start) to (idC-end);
\end{tikzpicture}
\caption{Trees with $ α_3 + α_4 $ as main result component of $ π_q μ^2 $}
\label{fig:subdisk-alpha34-mult-main}
\end{figure}
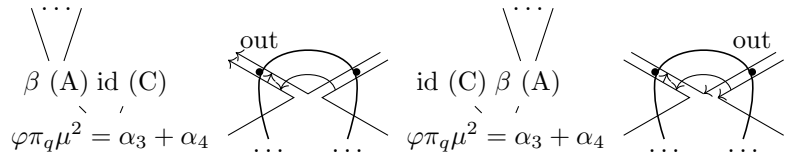
\begin{figure}
\centering
\begin{subfigure}{0.24\linewidth}
\begin{tikzpicture}
\path node (1) {$ \id $ (C)} node[right of=1] (2) {$ β $ (A)}
node[below right of=1] {$ h_q μ^2 = \id $ (B) (main)} edge (1) edge (2);
\end{tikzpicture}
\end{subfigure}
\begin{subfigure}{0.24\linewidth}
\begin{tikzpicture}
\path node (1) {$ β $ (A)} node[right of=1] (2) {$ \id $ (C)}
node[below right of=1] {$ h_q μ^2 = \id $ (B) (tail)} edge (1) edge (2);
\end{tikzpicture}
\end{subfigure}
\begin{subfigure}{0.24\linewidth}
\begin{tikzpicture}
\path node (1) {$ β/β' $} node[right of=1] (2) {$ \id $ (C)}
node[below right of=1] {$ h_q μ^2 = \id $ (B) (tail)} edge (1) edge (2);
\end{tikzpicture}
\end{subfigure}
\begin{subfigure}{0.24\linewidth}
\begin{tikzpicture}
\path node (1) {$ β $ (A)} node[right of=1] (2) {$ α_0 $}
node[below right of=1] {$ h_q μ^2 = \id $ (B) (tail)} edge (1) edge (2);
\end{tikzpicture}
\end{subfigure}
\caption{Trees with $ \id $ (B) as result component. No subdisk is assigned, but the trees are used for trees of $ \id $ (D).}
\label{fig:subdisk-idB}
\end{figure}
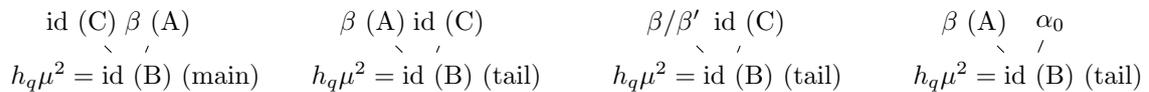
\input{tree_figures/subdisk-idD.tex}
\input{tree_figures/subdisk-degenerate.tex}
\begin{figure}
\begin{subfigure}{\linewidth}
\centering
\begin{tikzpicture}
\path node (A) {$ α_3 $} node[right of=A] (B) {$ \id $ (C)} node[right of=B] (C) {$ α_0 $}
node[below right of=A] (D) {$ \id $ (D)} edge (A) edge (B)
node[below right of=D] {$ φ π_q μ^2 = α_0 $} edge (D) edge (C);
\end{tikzpicture}
\begin{tikzpicture}
\path node (A) {$ \id $ (C)} node[right of=A] (B) {$ α_4 $} node[right of=B] (C) {$ α_0 $}
node[below right of=A] (D) {$ \id $ (D)} edge (A) edge (B)
node[below right of=D] {$ φ π_q μ^2 = α_0 $} edge (D) edge (C);
\end{tikzpicture}
%
\begin{tikzpicture}[scale=1.5]
\path[draw, <-] (0, 0) -- ++(60:1) coordinate[midway] (1) -- ++(300:1) coordinate[midway] (2) -- ++(60:1) coordinate[pos=0.1] (alpha0-start) coordinate[midway] (3) -- ++(300:1) coordinate[midway] (4) coordinate[pos=0.9] (alpha0-end) -- ++(60:1) coordinate[midway] (5) -- ++(300:1) coordinate[midway] (6) -- ++(60:1) coordinate[midway] (7) -- ++(300:1) coordinate[midway] (8);
\path[draw, semithick] ($ (1)!-0.05!(8) $) -- ($ (8)!-0.05!(1) $);
\path ($ (3)!0.5!(4) $) coordinate (alpha0);
\path ($ (3)!0.5!(4) + (-0.3, 0.2) $) coordinate (out) node[above] {out};
\path (6) ++(0.15, -0.1) coordinate (in1);
\path (6) ++(0.05, -0.2) coordinate (in2);
\path (in1) ++(right:0.5) coordinate (upright-cont);
\path (in2) ++(right:0.5) coordinate (lowright-cont);
\path[draw, fill=gray, fill opacity=0.5, semithick] (out) to[bend left] (alpha0) to[bend right] ++(0.2, -0.1) to[bend right=0] (in1) to (in2) to ($ (4) + (down:0.2) $) to[bend left] (out);
\path[draw, semithick] (in1) -- (upright-cont) (in2) -- (lowright-cont);
\path[draw, semithick] ($ (2) + (up:0.2) $) to[bend left=5] (out);
\path[draw, semithick] (out) to[bend right] ($ (2) + (up:0.4) $);
\path[fill] \foreach \i in {(in1), (in2), (alpha0), (out)} {\i circle[radius=0.04]};
\path[draw, ->, bend right=30] (alpha0-start) to node[midway, above, shift={(0, -0.1)}] {$ α_0 $} (alpha0-end);
\path (8) node[below right] {$ a_0 $};
\end{tikzpicture}
\caption{The inputs of these two trees consists $ α_3 $ (B) and $ \id $ (C) lying on the same arc, and the co-identity. The evaluation of $ h_q $ in the tree gives a sum of identities $ \id $ (D) ranging over all arcs lying between the arc and the co-identity. The subdisk is depicted for the second tree. It is a wedge lying between the zigzag curve and its Hamiltonian deformation. Since $ α_0 $ is the first input, the wedge lies on the side of $ α_0 $ where $ α_0 $ points to. In the case of the figure it lies to the right of $ α_0 $.}
\end{subfigure}

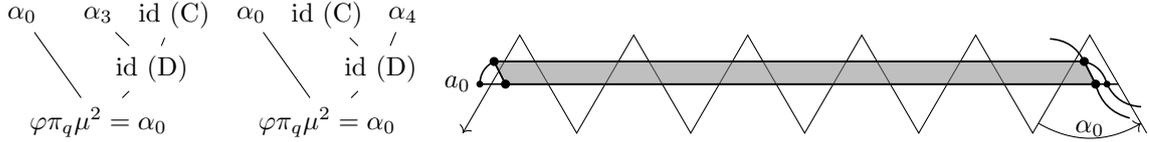
\begin{subfigure}{\linewidth}
\centering
\begin{tikzpicture}
\path node (A) {$ α_0 $} node[right of=A] (B) {$ α_3 $} node[right of=B] (C) {$ \id $ (C)}
node[below right of=B] (D) {$ \id $ (D)} edge (B) edge (C)
node[below left of=D] {$ φ π_q μ^2 = α_0 $} edge (D) edge (A);
\end{tikzpicture}
\begin{tikzpicture}
\path node (A) {$ α_0 $} node[right of=A] (B) {$ \id $ (C)} node[right of=B] (C) {$ α_4 $}
node[below right of=B] (D) {$ \id $ (D)} edge (B) edge (C)
node[below left of=D] {$ φ π_q μ^2 = α_0 $} edge (D) edge (A);
\end{tikzpicture}
%
\begin{tikzpicture}[scale=1.5]
\path[draw, <-] (0, 0) -- ++(60:1) coordinate[pos=0.5] (left-end) -- ++(300:1) coordinate[pos=0.5] (left-low) -- ++(60:1) -- ++(300:1) -- ++(60:1) -- ++(300:1) -- ++(60:1) -- ++(300:1) -- ++(60:1) -- ++(300:1) -- ++(60:1) coordinate[pos=0.1] (alpha0-start) coordinate[pos=0.5] (right-low) -- ++(300:1) coordinate[pos=0.3] coordinate[pos=0.5] (right-end) coordinate[pos=0.9] (alpha0-end);
\path ($ (right-low)!0.8!(right-end) $) coordinate (switch-right);
\path[draw, semithick] (left-end) ++(-0.1, 0) -- (right-end);
\path[fill] (switch-right) circle[radius=0.03];
\path[fill] (left-end) ++(-0.1, 0) circle[radius=0.03] node[left] {$ a_0 $};
\path[draw, semithick] ($ (right-low)!0.6!(right-end) $) coordinate (right-corner-low) -- ++(-0.1, 0.2) coordinate (right-corner-high) to[bend right] ++(-0.3, 0.2);
\path[draw, semithick] (right-corner-low) to[bend right] ++(0.3, -0.3);
\path[draw, semithick] ($ (left-low)!0.75!(left-end) $) coordinate (left-corner-low) -- ++(-0.1, 0.2) coordinate (left-corner-high);
\path[draw, semithick] (left-end) ++(-0.1, 0) to[bend left] (left-corner-high) to (right-corner-high) to[bend left] (switch-right) to[bend right] ++(0.3, -0.2);
\path[draw, semithick, fill=gray, fill opacity=0.5] (left-corner-high) to (right-corner-high) to (right-corner-low) to (left-corner-low) to cycle;
\path[fill] (right-corner-low) circle[radius=0.04];
\path[fill] (right-corner-high) circle[radius=0.04];
\path[fill] (left-corner-low) circle[radius=0.04];
\path[fill] (left-corner-high) circle[radius=0.04];
\path[draw, ->, bend right=30] (alpha0-start) to node[midway, above, shift={(0, -0.1)}] {$ α_0 $} (alpha0-end);
\end{tikzpicture}
\caption{The subdisk of the first tree is depicted. Since $ α_0 $ is the last input, the wedge lies on the opposite side of where $ α_0 $ points to. In the case of the figure it lies to the left of $ α_0 $.}
\end{subfigure}

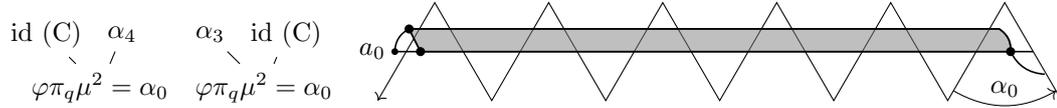
\begin{subfigure}{\linewidth}
\centering
\begin{tikzpicture}
\path node (A) {$ \id $ (C)} node[right of=A] (B) {$ α_4 $}
node[below right of=A] {$ φ π_q μ^2 = α_0 $} edge (A) edge (B);
\end{tikzpicture}
\begin{tikzpicture}
\path node (A) {$ α_3 $} node[right of=A] (B) {$ \id $ (C)}
node[below right of=A] {$ φ π_q μ^2 = α_0 $} edge (A) edge (B);
\end{tikzpicture}
%
\begin{tikzpicture}[scale=1.5]
\path[draw, <-] (0, 0) -- ++(60:1) coordinate[pos=0.5] (left-end) -- ++(300:1) coordinate[pos=0.5] (left-low) -- ++(60:1) -- ++(300:1) -- ++(60:1) -- ++(300:1) -- ++(60:1) -- ++(300:1) -- ++(60:1) -- ++(300:1) -- ++(60:1) coordinate[pos=0.1] (alpha0-start) coordinate[pos=0.5] (right-low) -- ++(300:1) coordinate[pos=0.3] coordinate[pos=0.5] (right-end) coordinate[pos=0.9] (alpha0-end);
\path ($ (right-low)!0.6!(right-end) $) coordinate (switch-right);
\path[draw, semithick] (left-end) ++(-0.1, 0) -- (right-end);
\path[fill] (switch-right) circle[radius=0.03];
\path[fill] (left-end) ++(-0.1, 0) circle[radius=0.03] node[left] {$ a_0 $};
\path[draw, semithick] ($ (left-low)!0.75!(left-end) $) coordinate (left-corner-low) -- ++(-0.1, 0.2) coordinate (left-corner-high);
\path[draw, semithick] (left-end) ++(-0.1, 0) to[bend left] (left-corner-high) to (right-corner-high) to[bend left] (switch-right) to[bend right] ++(0.3, -0.2);
\path[draw, semithick, fill=gray, fill opacity=0.5] (left-corner-high) to (right-corner-high) to[bend left] (switch-right) to (left-corner-low) to cycle;
\path[fill] (right-corner-low) circle[radius=0.04];
\path[fill] (left-corner-low) circle[radius=0.04];
\path[fill] (left-corner-high) circle[radius=0.04];
\path[draw, ->, bend right=30] (alpha0-start) to node[midway, above] {$ α_0 $} (alpha0-end);
\end{tikzpicture}
\caption{These two trees have only two inputs and yield the co-identity directly. The subdisk is depicted for the second tree and features a sample case where the two inputs and the output lie maximally far apart, namely at the identity and at the co-identity.}
\end{subfigure}

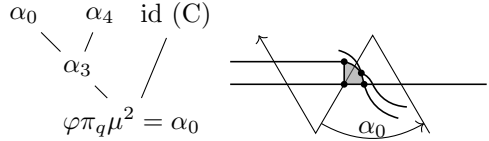
\begin{subfigure}{\linewidth}
\centering
\begin{tikzpicture}
\path node (A) {$ α_0 $} node[right of=A] (B) {$ α_4 $} node[right of=B] (C) {$ \id $ (C)}
node[below right of=A] (D) {$ α_3 $} edge (A) edge (B)
node[below right of=D] {$ φ π_q μ^2 = α_0 $} edge (D) edge (C);
\end{tikzpicture}
%
\begin{tikzpicture}[scale=1.5]
\path[draw, <-] (0, 0) -- ++(300:1) coordinate[midway] (1) -- ++(60:1) coordinate[midway] (2) coordinate[pos=0.1] (alpha0-start) -- ++(300:1) coordinate[midway] (3) coordinate[pos=0.9] (alpha0-end);
\path ($ (2)!0.5!(3) $) coordinate (alpha0);
\path ($ (2) + (0, 0.2) $) coordinate (in2);
\path ($ (alpha0) + (-0.1, 0.1) $) coordinate (in3);
\path ($ (2)!0.7!(alpha0) $) coordinate (out);
\path[draw, semithick] (2) to ++(left:1);
\path[draw, semithick] (in2) to ++(left:1);
\path[draw, semithick] (alpha0) to ++(right:1);
\path[draw, semithick] (alpha0) to[bend right] ++(0.3, -0.2);
\path[draw, semithick] (in3) to[bend right] ++(-0.2, 0.2);
\path[draw, semithick] (out) to[bend right] ++(0.3, -0.3);
\path[draw, semithick] (in3) to[bend left=20] (alpha0);
\path[draw, semithick] (out) -- (alpha0);
\path[draw, semithick, fill=gray, fill opacity=0.5] (in2) to[bend left=20] (in3) to (out) to (2) to cycle;
\path[fill] \foreach \i in {(2), (in2), (in3), (out)} {\i circle[radius=0.03]};
\path[draw, ->, bend right=30] (alpha0-start) to node[midway, above, shift={(0, -0.1)}] {$ α_0 $} (alpha0-end);
\end{tikzpicture}
\caption{This tree has the special feature that all its angles neighbor one single arc. The subdisk is tiny and concentrated around the midpoint of this arc.}
\end{subfigure}
\caption{Trees with $ α_0 $ as result component}
\label{fig:subdisk-alpha0}
\end{figure}

For example, regard a result component $ α_3 $ of an h-tree. We are interested in how this $ α_3 $ can possibly have been derived. A glance at the multiplication and disk tables \ref{tab:components-multiplication} and \ref{fig:components-first-out}, \ref{fig:components-final-out} reveals that it must be a product $ h_q μ^2 (α_0, α_4) $. In turn, we are interested in how $ α_0 $ and $ α_4 $ could have been derived. Another glance at the multiplication and disk tables reveals that both are necessarily direct. We conclude that any result component $ α_3 $ of an h-tree is necessarily the result component of the tree $ h_q μ^2 (α_0, α_3 + α_4) $ with leaves $ α_0 $ and $ α_3 + α_4 $.

We have conducted this investigation for all types of morphisms, resulting in the classification of \autoref{tab:resultcomp-classification}. Let us explain here how to read this table: The first and second column specify a type of result component. More precisely, the first column fixes the type of morphism. The second column sets further conditions on the type of result component. For example, the second column may indicate that only result components of h-trees shall be considered, or only tail result components of $ φ π_q μ^{≥3} $. The third column then provides a list of ways in which a result component of the specified type can be derived.

For example, we have seen before that an $ α_3 $ result component of an h-tree is necessarily direct or derived from the tree in \autoref{fig:subdisk-alpha3}. This is reflected in the fourth row of the classification table. As another example, we read off from the table that an $ \id $ (C) main result component of $ φ π_q μ^2 $ is necessarily a result component of one of the four $ \id $ (C) trees in \autoref{fig:subdisk-degenerate}.

As a final example, our classification of $ β $ (A) tail result components of h-trees is relatively implicit: These result components come from an evaluation $ h_q (βα) $, $ h_q (α_4 β) $, $ h_q (β' α_2) $ or $ h_q (α_3 α_4) $. These evaluations produce a $ β $ (A), $ β $ (C), $ β' $ (C) or $ α_3 $ (B) main result component. The tail $ β $ (A) result component then sits at the tail of these morphisms. In other words, whenever we encounter a tail $ β $ (A) result component of an h-tree, we will make reference to its associated main result components for further inspection.

\begin{remark}
Let us comment on two specific cases of \autoref{tab:resultcomp-classification}: Both $ \id $ (C) and $ α_3 + α_4 $ tail result component of $ φ π_q μ $ necessarily come from a certain $ φ π_q (βα) $ evaluation. The result component $ βα $ itself produces also main and tail result components from $ h_q (βα) $. We will use this observation later as a tool to abbreviate the construction of subdisks.
\end{remark}

\begin{remark}
The figures referenced in \autoref{tab:resultcomp-classification} show more than only trees: They depict trees and subdisks side-by-side. At present, the subdisks may simply be ignored and only trees count. We have chosen this way of presentation to facilitate retrospection during the reading of \autoref{sec:subdisk}.
\end{remark}

\begin{lemma}
\label{th:subdisk-assignment}
The result components classification of \autoref{tab:resultcomp-classification} is complete: The named types of result components can only be derived in the given way.
\end{lemma}

\begin{proof}
The checks are detailed in \autoref{sec:classification-shape}.
\end{proof}

\section{From trees to disks}
\label{sec:subdisk}
In this section, we show how to transform a result component of a π-tree into a kind of smooth immersed disk. Simply speaking, we draw all intersection points and connect them in a way dictated by the result component. The result is a matching between result components and certain types of smooth immersed disks. It leads to our main theorem which is a precise characterization of the minimal model $ \H\DefZigzagCat $ in terms of smooth immersed disks.

\begin{center}
\begin{tikzpicture}
\path (0, 0) node[align=center] (A) {\textbf{Result components} \\ with inputs $ h_1, …, h_N $} (10, 0) node[align=center] (B) {\textbf{CR, ID, DS, DW disks} \\ with inputs $ h_1, …, h_N $};
\path[draw, <->] ($ (A.east)!0.2!(B.west) $) -- ($ (A.east)!0.8!(B.west) $) node[midway, above] {Subdisk mapping $ \Subdisk $};
\end{tikzpicture}
\end{center}

Our bijection between result components of π-trees and smooth immersed disks is denoted $ \Subdisk $. The domain of the mapping $ \Subdisk $ is the set $ \PiTr $ of result components of π-trees. The precise terminology is that $ \Subdisk $ sends a result component $ r ∈ \PiTr $ to its associated “subdisk” $ \Subdisk(r) $. We define the mapping $ \Subdisk $ inductively over tree size, by also defining subdisks for result components of h-trees: If a result component $ r $ derives from certain result components $ r_1, …, r_N $ closer to the leaves, then the subdisk of $ r $ is defined by gluing the subdisks of $ r_1, …, r_N $.

In \autoref{sec:subdisk-protocol}, we introduce a protocol which lays down how subdisks may be glued together along handles. This way, subdisks of π-trees are closed disks, while subdisks of h-trees rather look like half a disk. In \autoref{sec:subdisk-shapeless}, we introduce a precise container class $ \SLd $ meant to capture the shape of subdisks of π-trees. In \autoref{sec:subdisk-subdisk}, we show explicitly how to draw the subdisk associated with a result component. In \autoref{sec:subdisk-types}, we define the subdisk mapping $ \Subdisk: \PiTr → \SLd $ and classify its image. The class of disks reached by $ \Subdisk $ decomposes into four visually distinguished types: CR, ID, DS and DW. In \autoref{sec:subdisk-minmodel}, we finish our computational journey and state our precise description of the minimal model $ \H\DefZigzagCat $. In \autoref{sec:subdisk-main}, we state our main theorem. The proofs of intermediate classification results and sign computations have been placed in \autoref{sec:classification}. We continue using the shorthand $ μ ≔ μ_q ≔ μ_{\Add\Gtl_q Q} $, see \autoref{rem:deformed-mushorthand}.

\subsection{The subdisk protocol}
\label{sec:subdisk-protocol}
In this section, we introduce our protocol for subdisk handles. The purpose of this protocol is to give an accurate description of the handles with which we will glue subdisks together. Recall from \autoref{sec:prelim-fukaya-zigzag} that every zigzag path $ L $ comes with an associated zigzag curve $ \smooth L $. The subdisk of a result component of an h-tree should consist of a sequence of intersection points and segments of the zigzag curves involved, filled with half a disk. While the intersection points and zigzag segments can be located anywhere on the dimer, both endpoints of the sequence should be located near the value of the result component itself. For a given type of result component, we wish that the endpoints follow a predictable pattern to facilitate gluing of subdisks. The protocol presented here is meant to define this local pattern, although we will give no precise definition what kind of object a subdisk of an h-tree is from a global view.

\begin{remark}
The subdisk protocol enjoys the following characteristics:
\begin{itemize}
\item The protocol applies to every type of morphism that can appear as result component of h-trees, namely $ α_3 $, $ α_4 $, $ β/β' $ (C), $ α_0 $, $ α_0' $ and $ β $ (A).
\item For each morphism $ ε: L_1 → L_2 $ of these types, the protocol defines a germ (small interval) of the zigzag curves $ \tilde L_1 $ and $ \tilde L_2 $.
\item Every germ comes with a \emph{handle}. A handle is an indication which endpoint is its gluable outside, and an indication which surface side is regarded as disk inside and which as disk outside.
\end{itemize}
\end{remark}

\begin{figure}
\centering
\begin{subfigure}{0.2\linewidth}
\begin{tikzpicture}
\path[draw] (0, 0) -- ++(315:1) coordinate[pos=0.5] (low);
\path[draw] (0, 0) -- ++(right:1) coordinate[pos=0.4] (alpha3-end) coordinate[pos=0.5] (L2-bot) -- ++(up:1.5) coordinate[pos=0.5] (L1-mid) -- ++(right:1);
\path[draw] (1.9, -0.1) -- ++(left:1) -- ++(up:1.5) coordinate[pos=0.4] (alpha3-start) coordinate[pos=0.5] (L2-mid) -- ++(left:1) coordinate[pos=0.5] (L1-top);
\path[draw, very thick, -{Circle[scale=0.4]}] ($ 0.5*(L1-mid) + 0.5*(L2-mid) $) coordinate (mid) -- ++(135:0.3);
\path[draw, gray, decorate, decoration={border, amplitude=5pt, angle=90, segment length=0.75pt}] (mid) to ++(135:0.3);
\path[draw, very thick, -{Circle[scale=0.4]}] ($ (mid)!0.8!(L2-bot) $) -- (L2-bot);
\path[draw, gray, decorate, decoration={border, amplitude=5pt, angle=90, segment length=0.75pt, mirror}] ($ (mid)!0.8!(L2-bot) $) to (L2-bot);
\path[bend right=45] (alpha3-start) to node[near start, below] {$ α_3 $} (alpha3-end);
\path[draw, densely dotted] (mid) -- (L1-top);
\path[draw, densely dotted] (L2-bot) -- (low);
\end{tikzpicture}
\end{subfigure}
\begin{subfigure}{0.2\linewidth}
\begin{tikzpicture}
\path[draw] (0, 0) -- ++(right:1) -- ++(up:1.5) coordinate[pos=0.4] (alpha4-end) coordinate[pos=0.5] (L1-mid) -- ++(right:1);
\path[draw] (1.9, -0.1) -- ++(left:1) coordinate[pos=0.5] (alpha4-start) coordinate[pos=0.4] (L1-bot) -- ++(up:1.5) coordinate[pos=0.5] (L2-mid) -- ++(left:1);
\path[draw, very thick, -{Circle[scale=0.4]}] ($ 0.5*(L1-mid) + 0.5*(L2-mid) $) coordinate (mid) -- ++(45:0.3);
\path[draw, gray, decorate, decoration={border, amplitude=5pt, angle=90, segment length=0.75pt, mirror}] (mid) to ++(45:0.3);
\path[draw, very thick, -{Circle[scale=0.4]}] ($ (mid)!0.8!(L1-bot) $) -- (L1-bot);
\path[draw, gray, decorate, decoration={border, amplitude=5pt, angle=90, segment length=0.75pt}] ($ (mid)!0.8!(L1-bot) $) to (L1-bot);
\path (L1-bot) node[above left] {$ α_4 $};
\end{tikzpicture}
\end{subfigure}
\begin{subfigure}{0.2\linewidth}
\begin{tikzpicture}
\path[draw] (0, 0) -- ++(up:1) -- ++(left:1) coordinate[pos=0.5] (L2-mid) coordinate[pos=0.4] (beta-end);
\path[draw] (0.2, 0) -- ++(up:1) -- ++(right:1) coordinate[pos=0.5] (L1-mid) coordinate[pos=0.4] (beta-start);
\path[draw, very thick, {Circle[scale=0.4]}-] (L1-mid) -- ++(225:0.3);
\path[draw, gray, decorate, decoration={border, amplitude=5pt, angle=90, segment length=0.75pt, mirror}] (L1-mid) to ++(225:0.3);
\path[draw, very thick, {Circle[scale=0.4]}-] (L2-mid) -- ++(315:0.3);
\path[draw, gray, decorate, decoration={border, amplitude=5pt, angle=90, segment length=0.75pt}] (L2-mid) to ++(315:0.3);
\path[bend right=90, looseness=2] (beta-start) to node[midway, below] {$ β/β' $} (beta-end);
\end{tikzpicture}
\end{subfigure}
\begin{subfigure}{0.2\linewidth}
\begin{tikzpicture}
\path[draw] (0, -0.2) -- ++(120:2) coordinate[pos=0.8] (alpha0-end) coordinate[pos=0.5] (L2-mid);
\path[draw] (0, -0.2) -- ++(60:2) coordinate[pos=0.7] (alpha0-start) coordinate[pos=0.5] (L1-mid);
\path[draw, very thick, {Circle[scale=0.4]}-] (L1-mid) -- ++(left:0.3);
\path[draw, gray, decorate, decoration={border, amplitude=5pt, angle=90, segment length=0.75pt, mirror}] (L1-mid) to ++(left:0.3);
\path[draw, very thick, {Circle[scale=0.4]}-] (L2-mid) -- ++(right:0.3);
\path[draw, gray, decorate, decoration={border, amplitude=5pt, angle=90, segment length=0.75pt}] (L2-mid) to ++(right:0.3);
\path[bend right=20] (alpha0-start) to node[midway, below] {$ α_0 $} (alpha0-end);
\end{tikzpicture}
\end{subfigure}
\begin{subfigure}{0.2\linewidth}
\begin{tikzpicture}
\path[draw] (0, -0.2) -- ++(120:2) coordinate[pos=0.5] (L1-mid) coordinate[pos=0.2] (alpha0p-start);
\path[draw] (0, -0.2) -- ++(60:2) coordinate[pos=0.5] (L2-mid) coordinate[pos=0.15] (alpha0p-end);
\path[draw, very thick, {Circle[scale=0.4]}-] (L1-mid) -- ++(right:0.3);
\path[draw, gray, decorate, decoration={border, amplitude=5pt, angle=90, segment length=0.75pt, mirror}] (L1-mid) to ++(right:0.3);
\path[draw, very thick, {Circle[scale=0.4]}-] (L2-mid) -- ++(left:0.3);
\path[draw, gray, decorate, decoration={border, amplitude=5pt, angle=90, segment length=0.75pt}] (L2-mid) to ++(left:0.3);
\path[bend right=135, looseness=12] (alpha0p-start) to node[midway, above] {$ α_0' $} (alpha0p-end);
\end{tikzpicture}
\end{subfigure}
%
%
\begin{subfigure}{0.2\linewidth}
\centering
\begin{tikzpicture}
\path[draw] (0, 0) -- ++(20:1.5) coordinate[pos=0.5] (end) -- ++(160:1.5) coordinate[pos=0.5] (beta-end);
\path[draw] (2.9, 0) -- ++(160:1.5) coordinate[pos=0.5] (start) -- ++(20:1.5) coordinate[pos=0.5] (beta-start);
\path[draw, very thick, {Circle[scale=0.4]}-] (beta-start) -- ++(down:0.3);
\path[draw, very thick, {Circle[scale=0.4]}-] (end) -- ++(down:0.3);
\path[draw, gray, decorate, decoration={border, amplitude=5pt, angle=90, segment length=0.75pt, mirror}] (beta-start) to ++(down:0.3);
\path[draw, gray, decorate, decoration={border, amplitude=5pt, angle=90, segment length=0.75pt}] (end) to ++(down:0.3);
\path[bend right=65] (beta-start) to node[midway, below] {$ β $} (beta-end);
\end{tikzpicture}
\caption*{short} 
\end{subfigure}
\begin{subfigure}{0.2\linewidth}
\centering
\begin{tikzpicture}
\path[draw] (0, 0) -- ++(20:1.5) coordinate[pos=0.5] (end) -- ++(160:1.5) coordinate[pos=0.5] (beta-end);
\path[draw] (2.9, 0) -- ++(160:1.5) coordinate[pos=0.5] (start) -- ++(20:1.5) coordinate[pos=0.5] (beta-start);
\path[draw, very thick, {Circle[scale=0.4]}-] (beta-start) -- ++(down:0.3);
\path[draw, gray, decorate, decoration={border, amplitude=5pt, angle=90, segment length=0.75pt, mirror}] (beta-start) to ++(down:0.3);
\path[draw, very thick, {Circle[scale=0.4]}-] (beta-end) -- ++(down:0.3);
\path[draw, gray, decorate, decoration={border, amplitude=5pt, angle=90, segment length=0.75pt}] (beta-end) to ++(down:0.3);
\path[bend right=65] (beta-start) to node[midway, below] {$ β $} (beta-end);
\end{tikzpicture}
\caption*{long}
\end{subfigure}
\caption{Subdisk protocol}
\label{fig:subdisk-protocol}
\end{figure}
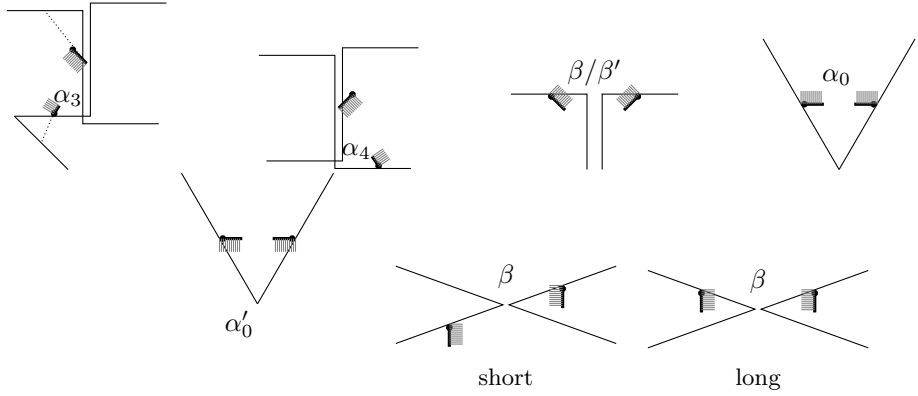

With these characteristics in mind, the protocol is defined in \autoref{fig:subdisk-protocol}. The germ intervals are drawn thickly, the gluable endpoints are drawn by dots and the disk inside is drawn hatched. Only $ β $ (A) comes in two variants: a short and a long version. We use the short version for tail components of $ h_q μ^2 $ and all components of $ h_q μ^{≥3} $, and the long version for direct morphisms and inputs of $ h_q μ^{≥3} $. We will explain the reason of this distinction in \autoref{rem:subdisk-subdisk-shortlong}.

\begin{remark}
\label{rem:subdisk-subdisk-shortlong}
The distinction between short and long version of $ β $ (A) protocol is due to a general phenomenon of subdisks. Namely, we will glue subdisks together by prolonging and subsequently connecting their handles. Sometimes, morphisms lie so close to each other that their handles connect without need for prolongation. In fact, if we drew every $ β $ (A) result component as the long version, it would strictly speaking not be possible to draw the right subdisks in some cases.

The best example is \autoref{fig:subdisk-betaA-mult-main}: The co-identities can be drawn one after another next to the connector of the $ β $ (A) input. With a long version, we would have to shorten the connector before drawing the co-identities. In convening a short and long version, we prioritized the convention that handles can always be prolonged and never shortened.
\end{remark}

\input{tree_figures/subdisk-connecting.tex}

\begin{remark}
The most easily imaginable result components come from discrete immersed disks $ h_q μ^{≥3} $. In these disks, some angles may be result components, while some are $ δ $-morphisms. By the subdisk protocol, the result components have handles assigned. These handles do not immediately connect to each other. Instead, the result components lie apart by as many arcs as the number $ n_i $ of $ δ $-morphisms between them. To facilitate smooth connections, we need to connect the handles by means of the \emph{angle cutting} procedure laid out in \autoref{def:prelim-fukaya-zigzag-curve}.
\end{remark}

We are now ready to use the protocol for the first time:

\begin{lemma}
\label{th:subdisk-delta-connect}
Let $ r_1, …, r_N $ be a sequence of result components $ r_i: L_i → L_{i+1} $. Assume the values of these result components are the consecutive angles of a discrete immersed disk when complemented with $ δ $-morphisms. Then their subdisk handles and the cuttings of the $ δ $-angles connect smoothly. Here, all short $ β $ (A) handles shall be extended to long ones first.
\end{lemma}

\begin{proof}
The sequence of angles of the discrete immersed disk is a mix of result components and $ δ $-morphisms. To check that everything connects smoothly, it suffices to check two neighbors at a time. These may be either two $ δ $-morphisms, two result components, or one $ δ $-morphism to the left or right of a result component. The first case of two $ δ $-morphisms is trivial.

The second case of two result components is checked in \autoref{fig:subdisk-connecting}. In this figure, all possible pairs of consecutive disk angles are checked for smoothness.

The third case of one result component and one $ δ $-morphism is an automatic feature of the subdisk protocol. The example case of $ α_3 $ is depicted by dotted lines in \autoref{fig:subdisk-protocol}.
\end{proof}

\subsection{Shapeless disks}
\label{sec:subdisk-shapeless}
In this section, we introduce shapeless disks as a container type for subdisks of π-trees. Recall that we intend to define a mapping between result components of π-trees and certain types of smooth immersed disks. In the present section, we define a suitable codomain for this map. Our solution is a broad container format, which we call shapeless disks. A shapeless disk consists of intersection points of zigzag curves with curve segments in between, filled by a disk immersion up to reparametrization. The specialty of shapeless disks is that intersection points may occur multiple times, with zero distance between each other. The set $ \SLd $ of all shapeless disks will serve as codomain of the map $ \Subdisk: \PiTr → \SLd $.

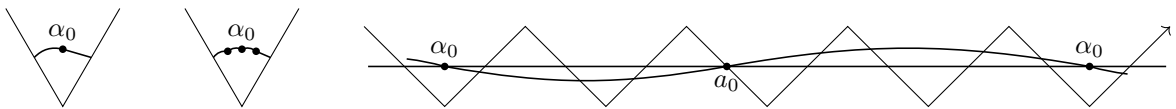
\begin{figure}
\centering
\begin{subfigure}{0.14\linewidth}
\begin{tikzpicture}
\path[draw] (0, 0) -- ++(300:1.5) coordinate[midway] (end) -- ++(60:1.5) coordinate[midway] (start);
\path[bend right] (start) to coordinate[midway] (switch) (end);
\path[draw, semithick] (start) to (switch) to[bend right] (end);
\path[fill] (switch) circle[radius=0.05] node[above] {$ α_0 $};
\end{tikzpicture}
\end{subfigure}
\begin{subfigure}{0.14\linewidth}
\begin{tikzpicture}
\path[draw] (0, 0) -- ++(300:1.5) coordinate[midway] (end) -- ++(60:1.5) coordinate[midway] (start);
\path[bend right] (start) to coordinate[pos=0.25] (switch1) coordinate[pos=0.5] (switch2) coordinate[pos=0.75] (switch3) (end);
\path[draw, semithick] (start) to (switch1) to[bend right] (switch2) to[bend right] (switch3) to[bend right] (end);
\path[fill] (switch1) circle[radius=0.05];
\path[fill] (switch2) circle[radius=0.05] node[above] {$ α_0 $};
\path[fill] (switch3) circle[radius=0.05];
\end{tikzpicture}
\end{subfigure}
\begin{subfigure}{0.7\linewidth}
\begin{tikzpicture}
\path[draw, ->] (0, 0) -- ++(315:1.5) coordinate[midway] (stop-1) -- ++(45:1.5) coordinate[midway] (stop-2) -- ++(315:1.5) -- ++(45:1.5) -- ++(315:1.5) coordinate[midway] (stop-5) -- ++(45:1.5) -- ++(315:1.5) -- ++(45:1.5) -- ++(315:1.5) coordinate[midway] (stop-9) -- ++(45:1.5) coordinate[midway] (stop-10);
\path[draw, semithick] ($ (stop-1)!-0.05!(stop-10) $) -- ($ (stop-1)!1.05!(stop-10) $);
\path ($ (stop-9)!0.5!(stop-10) $) coordinate (coid-right);
\path ($ (stop-1)!0.5!(stop-2) $) coordinate (coid-left);
\path[draw, semithick] ($ (coid-right) + (0.5, -0.1) $) to[bend left=4] (coid-right) to[bend right=10] (stop-5) to[bend left=10] (coid-left) to[bend right=4] ++(-0.5, 0.1);
\path[fill] (coid-left) circle[radius=0.05] node[above] {$ α_0 $};
\path[fill] (stop-5) circle[radius=0.05] node[below] {$ a_0 $};
\path[fill] (coid-right) circle[radius=0.05] node[above] {$ α_0 $};
\end{tikzpicture}
\end{subfigure}
\caption{A co-identity of $ L $ in a subdisk is drawn as a switch from $ L $ to its Hamiltonian deformation.}
\label{fig:subdisk-coidentity}
\end{figure}

\begin{remark}
In the definition of SL disks, it is essential that all cohomology basis elements are understood as intersection points between the associated zigzag curves. This correspondence is defined in \autoref{sec:prelim-fukaya-zigzag}. In particular, an identity $ \id_L $ is viewed as the even intersection point of $ \smooth L $ and its Hamiltonian deformation $ \smooth L' $, located at the midpoint of the identity location arc $ a_0 $ of $ L $. The co-identity $ α_0 $ is viewed as the odd intersection point between $ \smooth L $ and $ \smooth L' $, located at the midpoint of the chosen co-identity angle $ α_0 $. The Hamiltonian deformation $ \smooth L' $ goes right of $ \smooth L $ at $ α_0 $ and left of $ \smooth L $ at $ a_0 $, see \autoref{fig:subdisk-coidentity}.
\end{remark}

We are now ready to define shapeless disks.

\begin{definition}
\label{def:subdisk-shapeless-def}
Let $ N ≥ 0 $ and let $ L_1, …, L_{N+1} $ be a sequence of zigzag paths. Let $ h_i: L_i → L_{i+1} $ be cohomology basis elements. An \emph{SL disk} (shapeless disk) consists of
\begin{itemize}
\item an output cohomology basis element $ t: L_{N+1} → L_1 $,
\item a possibly empty $ \smooth L_1 $ segment from $ t $ to $ h_1 $,
\item for every $ i = 1, …, N $ a possibly empty $ \smooth L_{i+1} $ segment from $ h_i $ to $ h_{i+1} $,
\item a possibly empty $ \smooth L_{N+1} $ segment from $ h_N $ to $ t $,
\item an oriented polygon immersion $ D: P_{N+1} → |Q| $ up to reparametrization,
\end{itemize}
such that $ D $ has convex corners and traces the segments of $ \smooth L_1, …, \smooth L_{N+1} $ one after another. More precisely, $ D $ shall map the $ i $-th corner to $ h_i $, the $ N+1 $-th corner to $ t $, the edge between $ i $-th and $ i+1 $-th corner to the $ \smooth L_{i+1} $ segment and the edge between $ N+1 $-th and 1st corner to the $ \smooth L_1 $ segment, lying on the right side of this chain of segments. The mapping $ D $ need not be an immersion on the boundary. The disk may have infinitesimally small area. The class of SL disks is denoted $ \SLd $.
\end{definition}

\begin{figure}
\centering
\begin{subfigure}[b]{0.3\linewidth}
\centering
\begin{tikzpicture}
\path[draw, semithick] (0, 0) to[out=0, in=315, looseness=2] coordinate[pos=0.6] (node) (1.5, 1.5) to[out=135, in=90, looseness=2] (0, 0);
\path[fill] (0, 0) circle[radius=0.05] node[below] {$ t $};
\path (node) node[right] {$ \smooth L_1 $};
\path (0.75, 0.75) node {$ D $};
\end{tikzpicture}
\caption{$ N = 0 $}
\label{fig:subdisk-shapeless-SL0}
\end{subfigure}
\begin{subfigure}[b]{0.3\linewidth}
\centering
\begin{tikzpicture}
\path[draw, semithick] (0, 0) to[bend left=30] coordinate[midway] (node1) (3, 0) to[bend left=30] coordinate[midway] (node2) (0, 0);
\path[fill] (0, 0) circle[radius=0.05] node[below] {$ t $} (3, 0) circle[radius=0.05] node[below] {$ h_1 $};
\path (node1) node[above] {$ \smooth L_1 $};
\path (node2) node[below] {$ \smooth L_2 $};
\end{tikzpicture}
\caption{$ N = 1 $}
\label{fig:subdisk-shapeless-SL1}
\end{subfigure}
\begin{subfigure}[b]{0.3\linewidth}
\centering
\begin{tikzpicture}
\path[draw, semithick] (0, 0) coordinate (h2) to[bend left=5] ++(2, 0) coordinate (h1) to[bend left=8] ++(1, 0) coordinate (t) to[bend left=5] ++(-1, 2) coordinate (h6) to[bend left=5] ++(-1, 0) coordinate (h5) to[bend left=5] ++(-1, -1.6) coordinate (h4) to ++(down:0.2) coordinate (h3) to ++(down:0.2) coordinate (h2);
\path[fill] \foreach \i in {h1, h2, h3, h4, h5, h6, t} {(\i) circle[radius=0.05]};
\path (h1) node[below] {$ h_1 $};
\path (h2) node[below left] {$ h_2 $};
\path (h3) node[right] {$ h_3 $};
\path (h4) node[above left] {$ h_4 $};
\path (h5) node[above left] {$ h_5 $};
\path (h6) node[above right] {$ h_6 $};
\path (t) node[below right] {$ t $};
\path (t) -- (h1) node[midway, below] {\tiny $ \smooth L_1 $};
\path (h1) -- (h2) node[midway, below] {\tiny $ \smooth L_2 $};
\path (h2) -- (h3) node[midway, left] {\tiny $ \smooth L_3 $};
\path (h3) -- (h4) node[midway, left] {\tiny $ \smooth L_4 $};
\path (h4) -- (h5) node[midway, left] {\tiny $ \smooth L_5 $};
\path (h5) -- (h6) node[midway, above] {\tiny $ \smooth L_6 $};
\path (h6) -- (t) node[midway, right] {\tiny $ \smooth L_7 $};
\path (1.5, 0.9) node {$ D $};
\end{tikzpicture}
\caption{$ N = 6 $}
\label{fig:subdisk-shapeless-SLN}
\end{subfigure}
\caption{These figures illustrate SL disks with a given number of $ N = 0, 1, 6 $ inputs. In the SL disk with six inputs, the two zigzag paths $ L_1 $ and $ L_2 $ are supposed to be equal and the input $ h_1 $ is supposed to be the co-identity of $ L_1 $. The three inputs $ h_2 $, $ h_3 $, $ h_4 $ lie infinitesimally close to each other. The way we have portrayed them is meant to imply $ L_2 = L_4 $ and $ L_3 = L_5 $. The morphisms $ h_2 $, $ h_3 $, $ h_4 $ change back and forth from $ L_2 $ to $ L_5 $. Allowing this distinctive behavior is the reason for our definition of SL disks.}
\label{fig:subdisk-shapeless-SL}
\end{figure}
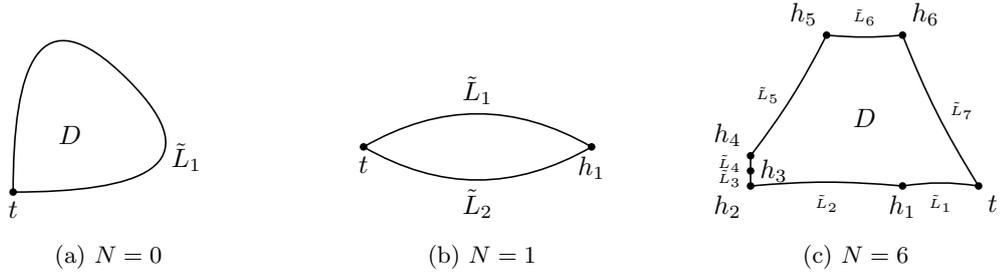

\begin{remark}
The notion of SL disks is depicted in \autoref{fig:subdisk-shapeless-SL}. We may refer to an empty segment also as a segment of \emph{infinitesimally small length} and say that the two endpoints of the segment are \emph{infinitesimally close}.
\end{remark}

\begin{remark}
The definition of an SL disk entails the option of infinitesimally small area and empty zigzag curve segments. It is impossible to draw these accurately, so we have opted to visually inflate every infinitesimally small area and empty segments and draw them as substantial area and short but visible segments in all drawings. While the definition of SL disks does technically not involve any Hamiltonian deformations, we always draw co-identity and identity as switches from $ \smooth L $ to Hamiltonian deformation, see \autoref{fig:subdisk-coidentity}. We draw stacked co-identities as repeated switches from $ \smooth L $ to $ \smooth L' $ to $ \smooth L'' $ etc.~with infinitesimally small distance in between, in line with the Fukaya-theoretic viewpoint.
\end{remark}

\begin{remark}
An SL disk is in principle allowed to have as few as zero or one inputs. An SL disk without inputs is a monogon, an SL disk with a single input is a digon. Under the present assumption that $ Q $ is geometrically consistent, an SL disk automatically has a minimum of two inputs. There are a few exceptions: The monogon with infinitesimally small area, located at an arbitrary intersection point, constitutes an SL disk without inputs. The digon bounded by two infinitesimally small segments of two intersecting zigzag curves, located at a single intersection, constitutes an SL disk with a single input. The digon bounded by a zigzag path and itself with input the identity and output the co-identity constitutes an SL disk with a single input. For geometrically consistent $ Q $, all SL disks with less than two inputs have infinitesimal area. They are an artifact of the definition and will not be used.
\end{remark}

\subsection{Constructing subdisks}
\label{sec:subdisk-subdisk}
In this section, we define subdisks for most Kadeishvili trees. As announced, the procedure is an inductive drawing construction, taking into account the way a given result component was derived. The reader has encountered many subdisk drawings already, spread out over figures from \autoref{sec:resultcomp}. Here we will explain these drawings and add more.

\begin{definition}
Regard a $ β $ (A), $ α_3 $ (B), $ β/β' $ (C) or $ α_0' $ (D) result component of an h-tree or an $ α_3 + α_4 $ (B), $ \id $ (C), $ \sum \id_a $ (D) or $ α_0 $ (D) result component of a π-tree. Then its \emph{subdisk} is defined inductively by the catalog presented in the rest of this section.
\end{definition}

\begin{figure}
\begin{subfigure}{0.2\linewidth}
\begin{tikzpicture}
\path[draw] (0, 0) -- ++(right:1) coordinate[pos=0.4] (alpha3-end) coordinate[pos=0.5] (L2-bot) -- ++(up:1.5) coordinate[pos=0.5] (L1-mid) -- ++(right:1);
\path[draw] (1.9, -0.1) -- ++(left:1) -- ++(up:1.5) coordinate[pos=0.4] (alpha3-start) coordinate[pos=0.5] (L2-mid) -- ++(left:1) coordinate[pos=0.5] (L1-top);
\path[draw, fill, semithick] ($ 0.5*(L1-mid) + 0.5*(L2-mid) $) circle (0.8pt) -- (L2-bot);
\path[draw, ->, bend right=45] (alpha3-start) to node[near start, below] {$ α_3 $} (alpha3-end);
\end{tikzpicture}
\end{subfigure}
\begin{subfigure}{0.2\linewidth}
\begin{tikzpicture}
\path[draw] (0, 0) -- ++(right:1) -- ++(up:1.5) coordinate[pos=0.4] (alpha4-end) coordinate[pos=0.5] (L1-mid) -- ++(right:1);
\path[draw] (1.9, -0.1) -- ++(left:1) coordinate[pos=0.5] (alpha4-start) coordinate[pos=0.4] (L1-bot) -- ++(up:1.5) coordinate[pos=0.5] (L2-mid) -- ++(left:1);
\path[draw, fill, semithick] ($ 0.5*(L1-mid) + 0.5*(L2-mid) $) circle (0.8pt) -- (L1-bot);
\path[draw, ->, bend right=45] (alpha4-start) to node[near start, left] {$ α_4 $} (alpha4-end);
\end{tikzpicture}
\end{subfigure}
\begin{subfigure}{0.2\linewidth}
\begin{tikzpicture}
\path[draw] (0, 0) -- ++(up:1) -- ++(left:1) coordinate[pos=0.5] (L2-mid) coordinate[pos=0.4] (beta-end);
\path[draw] (0.2, 0) -- ++(up:1) -- ++(right:1) coordinate[pos=0.5] (L1-mid) coordinate[pos=0.4] (beta-start);
\path[draw, fill, semithick] (L1-mid) -- (0.1, 0.5) circle (1pt) -- (L2-mid);
\path[draw, ->, bend right=90, looseness=2] (beta-start) to node[midway, below] {$ β/β' $} (beta-end);
\end{tikzpicture}
\end{subfigure}
\begin{subfigure}{0.2\linewidth}
\begin{tikzpicture}
\path[draw] (0, 0) -- ++(120:1.5) coordinate[pos=0.8] (alpha0-end) coordinate[pos=0.6] (L2-mid);
\path[draw] (0, 0) -- ++(60:1.5);
\path[draw] (0, -0.2) -- ++(120:2);
\path[draw] (0, -0.2) -- ++(60:2) coordinate[pos=0.7] (alpha0-start) coordinate[pos=0.5] (L1-mid);
\path[draw, semithick] (L1-mid) to ++(left:0.5) [fill] circle (0.8pt) coordinate (mid);
\path[draw, semithick] (mid) to[bend right] (L2-mid);
\path[draw, ->, bend right=20] (alpha0-start) to node[midway, below] {$ α_0 $} (alpha0-end);
\end{tikzpicture}
\end{subfigure}
\begin{subfigure}{0.2\linewidth}
\begin{tikzpicture}
\path[draw] (0, 0) -- ++(120:1.5);
\path[draw] (0, 0) -- ++(60:1.5) coordinate[pos=0.4] (L2-mid) coordinate[pos=0.15] (alpha0p-end);
\path[draw] (0, -0.2) -- ++(120:2) coordinate[pos=0.5] (L1-mid) coordinate[pos=0.2] (alpha0p-start);
\path[draw] (0, -0.2) -- ++(60:2);
\path[draw, semithick] (L1-mid) to ++(right:0.7) to ++(left:0.2) [fill] circle (0.8pt) coordinate (mid);
\path[draw, semithick] (mid) to[bend right] (L2-mid);
\path[draw, ->, bend right=135, looseness=12] (alpha0p-start) to node[midway, above] {$ α_0' $} (alpha0p-end);
\end{tikzpicture}
\end{subfigure}
\begin{subfigure}{0.2\linewidth}
\begin{tikzpicture}
\path[draw] (0, 0) -- ++(20:1.5) coordinate[pos=0.5] (end) -- ++(160:1.5) coordinate[pos=0.5] (beta-end);
\path[draw] (2.9, 0) -- ++(160:1.5) coordinate[pos=0.5] (start) -- ++(20:1.5) coordinate[pos=0.5] (beta-start);
\path[draw, semithick] (start) -- ++(down:0.4) node[below] {…};
\path[draw, semithick] (end) -- ++(down:0.4) node[below] {…};
\path[draw, ->, bend right=65] (beta-start) to node[midway, below] {$ β $} (beta-end);
\end{tikzpicture}
\end{subfigure}
\caption{Subdisks of direct morphisms}
\label{fig:subdisk-direct}
\end{figure}
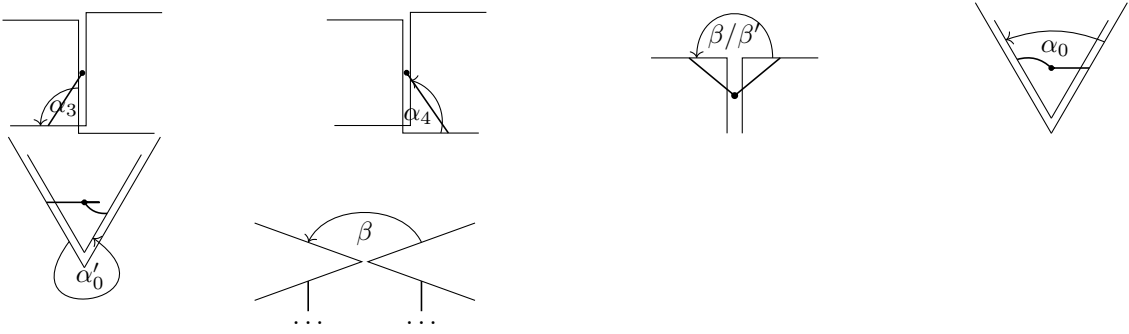

\paragraph*{Any $ α_3 $, $ α_4 $, $ β/β' $ (C), $ α_0 $, $ α_0' $ result component.} Depending on whether direct or not, their subdisks are given in \autoref{fig:subdisk-direct}, \ref{fig:subdisk-alpha3}, \ref{fig:subdisk-betaC}, \ref{fig:subdisk-betaCp} and \ref{fig:subdisk-alpha0p}.

\paragraph*{Direct $ β $ (A).} Note it is necessarily part of a tail of some morphism $ ε $, which is either $ α_3 $, $ α_4 $ or $ β/β' $ (C). The subdisk of $ β $ (A) is now obtained by taking the subdisk of $ ε $ and connecting it all the way around the tail disks by cutting the $ δ $ angles, continuing up until the given $ β $ (A) component. Finish with the short subdisk version of $ β $ (A).

\paragraph*{A $ β $ (A) main result component of $ h_q μ^2 $.} Its subdisk is shown in \autoref{fig:subdisk-betaA-mult-main}.

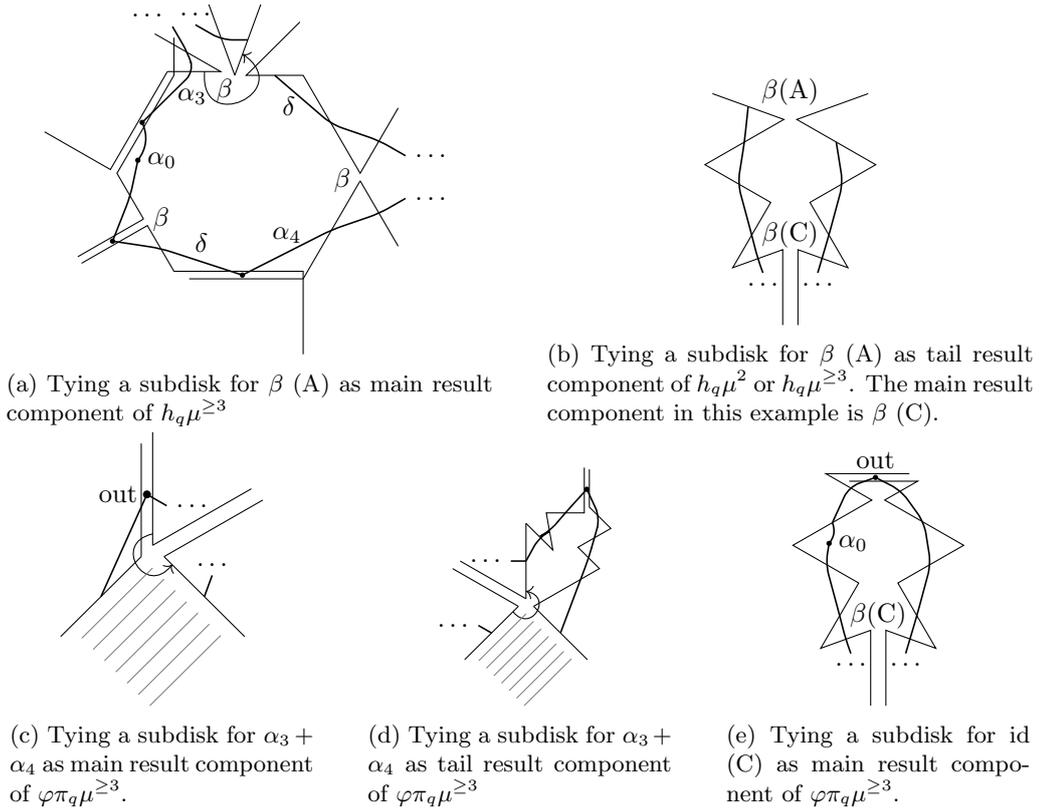
\begin{figure}
\centering
\begin{subfigure}{0.4\linewidth}
\centering
\begin{tikzpicture}
\path[draw] (0, 0) -- ++(225:1) coordinate[midway] (c1) coordinate (start) -- ++(right:0.75) coordinate[midway] (c2) -- ++(300:1.5) coordinate (betaA-top) coordinate[midway] (c3) -- ++(60:1) coordinate[midway] (c4);
\path[draw] (betaA-top) ++(0, -0.1) node[left] {$ β $} -- ++(300:1) coordinate[midway] (c5);
\path[draw] (betaA-top) ++(0, -0.1) -- ++(240:1.5) coordinate (alpha4-center) coordinate[midway] (c6) -- ++(left:1.5) coordinate[midway] (c7);
\path[draw] (alpha4-center) -- ++(down:1);
\path[draw] (alpha4-center) -- ++(up:0.1) -- ++(left:1.7) coordinate[midway] (c8) -- ++(120:0.7) coordinate (beta-L1) coordinate[midway] (c9) -- ++(210:1) coordinate[midway] (c10);
\path[draw] (beta-L1) ++(120:0.1) node[right] {$ β $} -- ++(210:1) coordinate[midway] (c11);
\path[draw] (beta-L1) ++(120:0.1) -- ++(120:0.7) coordinate[midway] (c12) -- ++(60:1.5) coordinate (alpha3-center) coordinate[midway] (c13) -- ++(90:0.5);
\path[draw] (alpha3-center) ++(150:0.1) -- ++(240:1.5) coordinate[midway] (c14) -- ++(150:1);
\path[draw] (alpha3-center) ++(150:0.1) -- ++(right:0.7) coordinate[midway] (c15) coordinate[pos=0.7] (beta-start) -- ++(150:1) coordinate[midway] (c16);
\path[draw] (start) ++(-0.15, 0) -- ++(70:1) coordinate[midway] (c17) coordinate[pos=0.3] (beta-end);
\path[draw] (start) ++(-0.15, 0) -- ++(110:1) coordinate[midway] (c18);
\path[draw, semithick, rounded corners] (c2) -- (c3) node[near start, below] {$ δ $} -- (c4) -- ++(330:0.4) node[right] {…};
\path[draw, semithick, rounded corners] (c5) ++(30:0.4) node[right] {…} -- (c5) -- (c6) -- ($ (c7)!0.5!(c8) $) coordinate (p1) node[midway, above] {$ α_4 $};
\path[draw, semithick, rounded corners] (p1) -- (c9) node[midway, above] {$ δ $} -- ($ (c10)!0.5!(c11) $) coordinate (p2);
\path[draw, semithick, rounded corners] (p2) -- (c12) -- ($ (c12)!0.5!(c13) $) coordinate (switch) node[right] {$ α_0 $};
\path[draw, semithick, bend right] (switch) to coordinate[pos=1] (p3) ($ (c13)!0.5!(c14) $);
\path[draw, semithick, rounded corners] (p3) -- (c15) node[midway, right] {$ α_3 $} -- (c16) -- ++(120:0.4) node[above left] {…};
\path[draw, semithick, rounded corners] (c17) -- (c18) -- ++(150:0.4) node[above] {…};
\path[fill] (p1) circle (1pt);
\path[fill] (p2) circle (1pt);
\path[fill] (p3) circle (1pt);
\path[fill] (switch) circle (1pt);
\path[draw, ->, bend right=120, looseness=4] (beta-start) to node[very near start, right] {$ β $} (beta-end);
\end{tikzpicture}
\caption{Tying a subdisk for $ β $ (A) as main result component of $ h_q μ^{≥3} $}
\label{fig:subdisk-betaA-disk-main}
\end{subfigure}
\hspace{0.03\linewidth}
\begin{subfigure}{0.4\linewidth}
\centering
\begin{tikzpicture}
\path[draw] (0, 0) -- ++(up:1) -- ++(200:0.7) coordinate[pos=0.4] (betaC-end) coordinate[midway] (c1) -- ++(60:1) coordinate[midway] (c2) -- ++(150:1) coordinate[midway] (c3) -- ++(30:1.2) coordinate[midway] (c4) -- ++(160:1) coordinate[pos=0.4] (betaA-end) coordinate[midway] (c5);
\path[draw] (0.2, 0) -- ++(up:1) -- ++(340:0.7) coordinate[midway] (c6) coordinate[pos=0.4] (betaC-start) -- ++(120:1) coordinate[midway] (c7) -- ++(30:1) coordinate[midway] (c8) -- ++(150:1.2) coordinate[midway] (c9) -- ++(20:1) coordinate[pos=0.4] (betaA-start) coordinate[midway] (c10);
\path[bend right=110, looseness=3] (betaC-start) to node[midway, below] {$ β $(C)} (betaC-end);
\path[bend right=80, looseness=2] (betaA-start) to node[midway, below] {$ β $(A)} (betaA-end);
\path[draw, semithick, rounded corners] (c1) ++(290:0.2) node[below] {…} -- (c1) -- (c2) -- (c3) -- (c4) -- (c5);
\path[draw, semithick, rounded corners] (c6) ++(250:0.2) node[below] {…} -- (c7) -- (c8) -- (c9);
\end{tikzpicture}
\caption{Tying a subdisk for $ β $ (A) as tail result component of $ h_q μ^2 $ or $ h_q μ^{≥3} $. The main result component in this example is $ β $ (C).}
\label{fig:subdisk-betaA-tail}
\end{subfigure}
\hspace{0.03\linewidth}
\begin{subfigure}{0.25\linewidth}
\centering
\begin{tikzpicture}[scale=1.5]
\path[draw] (0, 0) coordinate (out) -- (0, 1) coordinate[pos=0.5] (c1) coordinate[pos=0.1] (m-start);
\path[draw] (0.1, -0.1) coordinate (start) -- ++(315:1) coordinate[pos=0.5] (c2) coordinate[pos=0.1] (m-end);
\path[draw] (-0.1, -0.1) coordinate (end) -- ++(225:1) coordinate[pos=0.5] (c3);
\path[draw, ->, bend right=110, looseness=3] (m-start) to (m-end);
\path[draw, decoration={border, amplitude=37.5pt, angle=90, segment length=6pt, mirror}, color=gray, decorate] ($ (start) + (-0.1, -0.1) $) -- ++(315:1);
\path[draw] (out) -- ++(30:1);
\path[draw] (start) -- ++(30:1);
\path[draw] (end) -- ++(90:1) coordinate[pos=0.5] (c4);
\path[draw, semithick, rounded corners] (c3) -- ($ (c1)!0.5!(c4) $) coordinate (mid);
\path[draw, semithick] (mid) -- ++(330:0.2) node[right] {…};
\path[draw, semithick] (c2) -- ++(70:0.2) node[above] {…};
\path[fill] (mid) circle (1pt) node[left] {out};
\end{tikzpicture}
\caption{Tying a subdisk for $ α_3 + α_4 $ as main result component of $ φ π_q μ^{≥3} $.}
\label{fig:subdisk-alpha34-disk-main}
\end{subfigure}
\hspace{0.03\linewidth}
\begin{subfigure}{0.25\linewidth}
\centering
\begin{tikzpicture}
\path[draw] (0, 0) coordinate (out) -- (0, 1) coordinate[pos=0.1] (m-end) coordinate[midway] (c2);
\path[draw] (0.1, -0.1) coordinate (start) -- ++(315:1) coordinate[midway] (c11);
\path[draw] (-0.1, -0.1) coordinate (end) -- ++(225:1) coordinate[pos=0.1] (m-start) coordinate[midway] (c1);
\path[draw, ->, bend right=110, looseness=3] (m-start) to (m-end);
\path[draw, decoration={border, amplitude=25pt, angle=90, segment length=4pt, mirror}, color=gray, decorate] ($ (start) + (-0.1, -0.1) $) -- ++(315:1);
\path[draw] (end) -- ++(150:1);
\path[draw] (out) -- ++(150:1);
\path[draw] (start) -- ++(30:1) coordinate[midway] (c10) -- ++(135:0.4) coordinate[midway] (c10) -- ++(30:0.5) coordinate[midway] (c9) -- ++(135:0.4) coordinate[midway] (c8) -- ++(up:0.5) coordinate[midway] (c7);
\path[draw] (0, 1) -- ++(315:0.5) coordinate[midway] (c3) -- ++(100:0.5) coordinate[midway] (c4) -- ++(right:0.5) coordinate[midway] (c5) -- ++(up:0.6) coordinate[midway] (c6);
\path[draw, semithick] (c1) -- ++(150:0.2) node[left] {…};
\path[draw, semithick] (c2) -- ++(left:0.2) node[left] {…};
\path[draw, semithick, rounded corners] (c2) -- (c3) -- (c4) -- (c5) -- ($ (c6)!0.5!(c7) $) coordinate (mid);
\path[draw, semithick, rounded corners] (mid) -- (c8) -- (c9) -- (c10) -- (c11);
\path[fill] (mid) circle (1pt);
\end{tikzpicture}
\caption{Tying a subdisk for $ α_3 + α_4 $ as tail result component of $ φ π_q μ^{≥3} $}
\label{fig:subdisk-alpha34-disk-tail}
\end{subfigure}
\hspace{0.03\linewidth}
\begin{subfigure}{0.25\linewidth}
\centering
\begin{tikzpicture}
\path[draw] (0, 0) -- ++(up:1) -- ++(200:0.7) coordinate[pos=0.4] (betaC-end) coordinate[midway] (c1) -- ++(60:1) coordinate[midway] (c2) -- ++(150:1) coordinate[midway] (c3) -- ++(30:1.2) coordinate[midway] (c4) -- ++(150:0.7) coordinate[midway] (e1) -- ++(right:1.1) coordinate[pos=0.4] (betaA-end) coordinate[midway] (c5);
\path[draw] (0.2, 0) -- ++(up:1) -- ++(340:0.7) coordinate[midway] (c6) coordinate[pos=0.4] (betaC-start) -- ++(120:1) coordinate[midway] (c7) -- ++(30:1) coordinate[midway] (c8) -- ++(150:1.2) coordinate[midway] (c9) -- ++(30:0.5) coordinate[midway] (e2) -- ++(left:0.9) coordinate[pos=0.4] (betaA-start) coordinate[midway] (c10);
\path[draw, semithick, rounded corners] (c1) ++(290:0.2) node[below] {…} -- (c1) -- (c2) -- (c3) -- ($ (c3)!0.5!(c4) $) coordinate (switch);
\path[draw, semithick, bend right] (switch) to (c4);
\path[draw, semithick, rounded corners] (c4) -- (e1) -- ($ (c5)!0.5!(c10) $) coordinate (mid);
\path[draw, semithick, rounded corners] (c6) ++(250:0.2) node[below] {…} -- (c7) -- (c8) -- (c9) -- (e2) -- (mid);
\path[fill] (mid) circle (1pt) node[above] {out};
\path[fill] (switch) circle (1pt) node[right] {$ α_0 $};
\path[bend right=110, looseness=3] (betaC-start) to node[midway, below] {$ β $(C)} (betaC-end);
\end{tikzpicture}
\caption{Tying a subdisk for $ \id $ (C) as main result component of $ φ π_q μ^{≥3} $.}
\label{fig:subdisk-idC-disk-main}
\end{subfigure}
\caption{Further examples of how to tie subdisks}
\label{fig:subdisk-examples}
\end{figure}

\paragraph*{A $ β $ (A) main result component of $ h_q μ^{≥3} $.} The given discrete immersed disk is necessarily final-out, with final morphism an outer $ δ $ insertion, $ β $ (A) or $ β/β' $ (C). The result components that may be used in this higher product are $ δ $ insertions, $ β $ (A), $ β/β' $ (C), $ α_3 $ (B), $ α_4 $ (B), $ α_0 $ (D), $ α_0' $ (D). All of them have subdisk handles assigned. Close all $ β $ (A) handles. Connect the handles of all morphisms around the discrete immersed disk in clockwise order, following the $ δ $ insertions. Note that this produces a smooth curve according to \autoref{th:subdisk-delta-connect}. Finish with the short version of $ β $ (A). An example is shown in \autoref{fig:subdisk-betaA-disk-main}.

\paragraph*{A $ β $ (A) tail result component of $ h_q μ^2 $ or $ h_q μ^{≥3} $.} The corresponding main result component is a $ β $ (A), $ β/β' $ (C) or $ α_3 $. Now the subdisk of the $ β $ (A) tail result component is obtained by taking the subdisk of the main result component, closing it if it is a $ β $ (A), and connecting it all the way around the tail disk by cutting the $ δ $ angles, continuing up until the given $ β $ (A) tail result component. Finish with the short subdisk version. An example is shown in \autoref{fig:subdisk-betaA-tail}.

\paragraph*{An $ α_3 + α_4 $ main result component of $ φ π_q μ^{≥3} $.} The disk is then one of \autoref{fig:components-first-out} or \ref{fig:components-final-out}. In all cases, connect the handles all around the disk as in the $ β $ (A) case. If the disk is first-out, cut the $ δ $ angle at the beginning of the disk. If the disk is final-out, cut the $ δ $ angle at the end of the disk. Finally, close the disk with an output mark. An example is shown in \autoref{fig:subdisk-alpha34-disk-main}.

\paragraph*{An $ α_3 + α_4 $ main result component of $ φ π_q μ^2 $.} The entire tree is then one of those in \autoref{fig:subdisk-alpha34-mult-main}, where the subdisks are also depicted.

\paragraph*{An $ α_3 + α_4 $ tail result component of $ φ π_q μ^{≥3} $ or $ φ π_q μ^2 $.} It comes from a type G disk in a certain $ φ π_q (βα) $ evaluation of a product $ μ^2 $ or one of the disks $ μ^{≥3} $ of \autoref{fig:components-final-out}. Note that this very same $ β $ (A) appears as main result component of the $ h_q (βα) $ evaluation and we have already assigned a subdisk with short $ β $ (A) version to it. Now obtain the subdisk of $ α_3 + α_4 $ from the subdisk of $ β $ (A) by closing the subdisk and connecting it all the way up around the disk by cutting the $ δ $ angles, and finally finishing with an output mark at the 2/5 arc of the G situation. An example for $ φ π_q μ^{≥3} $ is shown in \autoref{fig:subdisk-alpha34-disk-tail}.

\paragraph*{An $ \id $ (C) main result component of $ φ π_q μ^2 $.} Its subdisk is depicted in \autoref{fig:subdisk-degenerate}.

\paragraph*{An $ \id $ (C) main result component of $ φ π_q μ^{≥3} $.} The disk is then all-in and of type H. Its inner morphisms may be $ δ $ insertions, $ β $ (A), $ α_3 $ (B), $ α_4 $ (B), $ β/β' $ (C), $ α_0 $ (D), $ α_0' $ (D). Connect them all and finish with an output mark on the concluding 2/5 arc of the disk. An example is shown in \autoref{fig:subdisk-idC-disk-main}.

\paragraph*{An $ \id $ (C) tail result component of $ φ π_q μ^2 $ or $ φ π_q μ^{≥3} $.} It comes from a type H disk in a certain $ φ π_q (βα) $ evaluation. Note that this $ β $ (A) already appears as a main result component and has a subdisk assigned. Now obtain the subdisk of $ \id $ (C) from closing the subdisk of $ β $ (A) and connecting it all the way up until the 2/5 concluding arc of the type H disk. Finish with an output mark.

\paragraph*{An $ \id $ (D) result component.} Its subdisk is depicted in \autoref{fig:subdisk-idD}.

\paragraph*{An $ α_0 $ result component of $ φ π_q μ^2 $ or $ φ π_q μ^{≥3} $.} Its subdisk is depicted in \autoref{fig:subdisk-alpha0}.

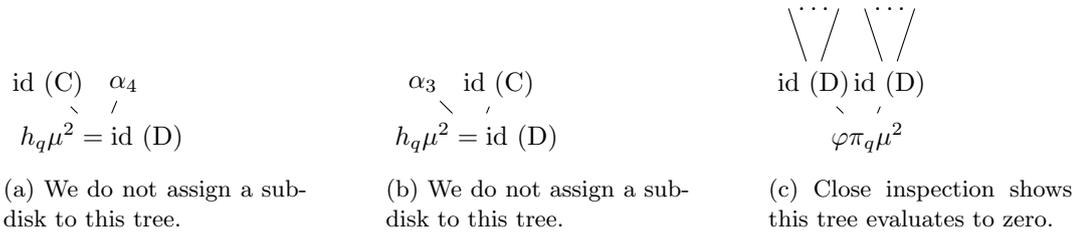
\begin{figure}
\centering
\begin{subfigure}{0.25\linewidth}
\begin{tikzpicture}
\path node (in1) {$ \id $ (C)} node[right of=in1] (in2) {$ α_4 $}
node[below right of=in1] {$ h_q μ^2 = \id $ (D)} edge (in1) edge (in2);
\end{tikzpicture}
\caption{We do not assign a subdisk to this tree.}
\label{fig:subdisk-idD-1}
\end{subfigure}
\hspace{0.05\linewidth}
\begin{subfigure}{0.25\linewidth}
\begin{tikzpicture}
\path node (in1) {$ α_3 $} node[right of=in1] (in2) {$ \id $ (C)}
node[below right of=in1] {$ h_q μ^2 = \id $ (D)} edge (in1) edge (in2);
\end{tikzpicture}
\caption{We do not assign a subdisk to this tree.}
\label{fig:subdisk-idD-2}
\end{subfigure}
\hspace{0.05\linewidth}
\begin{subfigure}{0.25\linewidth}
\begin{tikzpicture}
\path node (topdots1) {…} node[right of=topdots1] (topdots2) {…}
node[below of=topdots] (in1) {$ \id $ (D)} edge (topdots1.west) edge (topdots1.east)
node[right of=in1] (in2) {$ \id $ (D)} edge (topdots2.west) edge (topdots2.east)
node[below right of=in1] {$ φ π_q μ^2 $} edge (in1) edge (in2);
\end{tikzpicture}
\caption{Close inspection shows this tree evaluates to zero.}
\label{fig:subdisk-idD-double}
\end{subfigure}
\caption{Miscellaneous trees.}
\end{figure}

\begin{figure}
\centering
\begin{subfigure}{0.4\linewidth}
\begin{tikzpicture}
\path node (in1) {$ β/β' $} node[right of=in1] (in2) {$ \id $(C)} node[right of=in2] (in3) {$ \id $(C)}
node[below right of=in1] (m1) {$ βα $} edge (in1) edge (in2)
node[below right of=m1] {$ h_q μ^2 $} edge (m1) edge (in3);
\end{tikzpicture}
\begin{tikzpicture}
\path[draw] (0, 0) -- ++(down:0.7) coordinate[pos=0.5] (idC-start) coordinate[pos=0.7] (idB-end) coordinate[midway] (L3-top) -- ++(left:0.7) coordinate[midway] (1) -- ++(310:0.5) coordinate[midway] (2) -- ++(210:0.5) coordinate[midway] (3) -- ++(330:0.5) coordinate[midway] (4) -- ++(210:0.5) coordinate[midway] (5) -- ++(right:0.8) coordinate[midway] (8) -- ++(270:0.5);
\path[draw] (0.2, 0) -- ++(down:0.7) coordinate[pos=0.5] (idC-end) coordinate[pos=0.7] (idB-start) coordinate[pos=0.2] (L1-top) coordinate[midway] (r1) -- ++(right:0.7) coordinate[midway] (r2) -- ++(220:0.5) coordinate[midway] (r3) -- ++(330:0.5) coordinate[midway] (r4) -- ++(210:0.9) coordinate[midway] (r5) -- ++(right:1) coordinate[midway] (L1-bottom) ++(down:0.1) -- ++(left:1) coordinate[midway] (L2-top) -- ++(270:0.5);
\path[draw, ->] (0, 0) to ++(down:0.7);
\path[draw, ->] (0.2, 0) to ++(down:0.7);
\path[draw, semithick, rounded corners] ($ (8) + (0, -0.1) $) node[below] {…} -- (5) -- (4) -- (3) -- (2) -- (1) to[bend right=10] ($ (L3-top) + (-0.1, 0) $) coordinate (out);
\path[draw, semithick] (out) -- ($ (L1-top) + (0.1, 0) $) coordinate (idC);
\path[draw, semithick, rounded corners] (idC) -- (r2) -- (r3) -- (r4) -- (r5) -- ($ (L1-bottom)!0.5! (L2-top) $) coordinate (idC2);
\path[draw, semithick] (idC2) -- ++(-0.1, -0.2) node[below] {…};
\path[fill] (out) circle[radius=0.05] node[left] {out};
\path[fill] (idC) circle[radius=0.05] node[right] {$ \id $ (C)};
\path[fill] (idC2) circle[radius=0.05];
\path[draw, ->] (idC-start) to (idC-end);
\path[draw, ->] (idB-start) to (idB-end);
\end{tikzpicture}
\end{subfigure}
\begin{subfigure}{0.4\linewidth}
\centering
\begin{tikzpicture}
\path node (in1) {$ β $(A)} node[right of=in1] (in2) {$ α_0 $} node[right of=in2] (in3) {$ \id $(C)}
node[below right of=in1] (m1) {$ βα $} edge (in1) edge (in2)
node[below right of=m1] {$ h_q μ^2 $} edge (m1) edge (in3);
\end{tikzpicture}
\begin{tikzpicture}
\path[draw] (0, 0) -- ++(down:0.7) coordinate[pos=0.5] (idC-start) coordinate[pos=0.7] (idB-end) coordinate[midway] (L3-top) -- ++(left:0.7) coordinate[midway] (1) -- ++(310:0.5) coordinate[midway] (2) -- ++(210:0.5) coordinate[midway] (3) -- ++(330:0.5) coordinate[midway] (4) -- ++(210:0.5) coordinate[midway] (5) -- ++(right:0.8) coordinate[midway] (8) -- ++(240:0.5);
\path[draw] (0.2, 0) -- ++(down:0.7) coordinate[pos=0.5] (idC-end) coordinate[pos=0.7] (idB-start) coordinate[pos=0.2] (L1-top) coordinate[midway] (r1) -- ++(right:0.7) coordinate[midway] (r2) -- ++(220:0.5) coordinate[midway] (r3) -- ++(330:0.5) coordinate[midway] (r4) -- ++(210:0.9) coordinate[midway] (r5) -- ++(right:1) coordinate[midway] (L1-bottom) ++(down:0.1) -- ++(left:1) coordinate[midway] (L2-top) -- ++(300:0.5);
\path[draw, ->] (0, 0) to ++(down:0.7);
\path[draw, ->] (0.2, 0) to ++(down:0.7);
\path[draw, semithick, rounded corners] ($ (8) + (0, -0.1) $) node[below] {…} -- (5) -- (4) -- (3) -- (2) -- (1) to[bend right=10] ($ (L3-top) + (-0.1, 0) $) coordinate (out);
\path[draw, semithick] (out) -- ($ (L1-top) + (0.1, 0) $) coordinate (idC);
\path[draw, semithick, rounded corners] (idC) -- (r2) -- (r3) -- (r4) -- (r5) -- ($ (L2-top) + (-0.05, -0.1) $) coordinate (switch);
\path[draw, semithick] (switch) to[bend right=30] ++(0, -0.2) node[below] {…};
\path[fill] (out) circle[radius=0.05] node[left] {out};
\path[fill] (idC) circle[radius=0.05] node[right] {$ \id $ (C)};
\path[fill] (switch) circle[radius=0.05];
\path[draw, ->] (idC-start) to (idC-end);
\path[draw, ->] (idB-start) to (idB-end);
\end{tikzpicture}
\end{subfigure}
\caption{In \autoref{fig:subdisk-idD-idC-first}, the $ \id $ (B) component can impossibly come from $ μ^2 (β/β' \text{(C)}, \id \text{(C)}) $ or $ μ^2 (β \text{(A)}, α_0) $, because the arrow directions along the disk mismatch resp.~because the arrow direction of $ α_0 $ contradicts \autoref{conv:alpha0-direction}. The four resulting trees have no $ \id $ (D) result components. Of the four trees, the two with $ \id $ (C) right after the output mark are depicted here.}
\label{fig:subdisk-idD-impossible}
\end{figure}
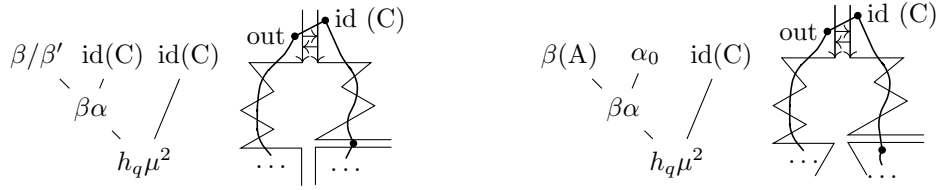

\begin{remark}
We have associated subdisks to all result components of all π-trees. Because they are difficult to draw consistently, we do not assign subdisks to $ \id $ (B) result components and $ \id $ (D) result components of h-trees.
\end{remark}

The reader who has read the catalog of subdisk definitions may feel unsure what these subdisks actually are. To ease his pain, we remind him that subdisks are specific collections of data defined in \autoref{sec:subdisk-protocol} and \autoref{sec:subdisk-shapeless}. Let us explain and record that the subdisks defined in the catalog actually satisfy these conditions:

\begin{lemma}
\label{th:subdisk-subdisk-Dconstruction}
Subdisks are well-defined. Subdisks of h-trees respect the subdisk protocol. Subdisks of π-trees are SL disks, providing a map $ \Subdisk: \PiTr → \SLd $.
\end{lemma}

\begin{proof}
This is easy and follows from induction on tree size. We shall not check all cases, but explain the line of argument. The base case of induction are the subdisks of the direct morphisms, which are depicted in \autoref{fig:subdisk-direct} and indeed respect the subdisk protocol.

As induction hypothesis, assume the subdisk of any result component of an h-tree with less than $ N $ inputs already respects the subdisk protocol. Regard a result component $ r $ of an h- or π-tree with $ N $ inputs. Assume $ r $ is derived from result components $ r_1, …, r_k $. Each of these result components $ r_i $ has less than $ N $ inputs and hence their subdisks respects the subdisk protocol.

According to the catalog, the subdisk of $ r $ is constructed by gluing or extending the subdisks of $ r_1, …, r_k $. At these points, the catalog typically invokes \autoref{th:subdisk-delta-connect}. This invokation is indeed possible since $ r_1, …, r_k $ all respect the subdisk protocol. The final step of the catalog entry is to finish the drawing somewhere near $ r $ itself. This step is indicated in individual pictures, from which it is evident that the finish respects the subdisk protocol respectively is an SL disk. This completes the induction.
\end{proof}

\subsection{The four types of disks}
\label{sec:subdisk-types}
In this section, we exhibit the image of $ \Subdisk: \PiTr → \SLd $. More precisely, we group result components of π-trees into four different types, according to the shape of their subdisk. These four types of result components go by the name CR, ID, DS and DW result components. We will also define four types of shapeless disks, which are meant to coincide with the image of these types of result components under $ \Subdisk $:

\begin{center}
\begin{tabular}{ccc}
\textbf{Geometry} & \textbf{Result component} & \textbf{Shapeless disk} \\\hline
Degenerate strip & DS result component & DS disk \\
Degenerate wedge & DW result component & DW disk \\
Identity degenerate & ID result component & ID disk \\
Co-identity rule & CR result component & CR disk
\end{tabular}
\end{center}

Recall from \autoref{sec:subdisk-subdisk} that every result component of a π-tree comes with a subdisk assigned. All of these subdisks are SL disks, but some are more special than others. For example, the subdisks depicted in \autoref{fig:subdisk-degenerate} are all degenerate: There are two zigzag curve segments with infinitesimally small length. In contrast, all segments in the subdisk in \autoref{fig:subdisk-idC-disk-main} are nonempty. We exploit these differences in subdisks to define four classes of result components:

\begin{definition}
\label{def:subdisk-types-rescomp}
A result component $ r ∈ \PiTr $ is a
\begin{itemize}
\item \emph{DS result component} if it is the result component of one of the 8 trees of \autoref{fig:subdisk-degenerate}.
\item \emph{DW result component} if it is the result component of one of the 7 trees of \autoref{fig:subdisk-alpha0}.
\item \emph{ID result component} if it is a result component of one of the trees in \autoref{fig:subdisk-idD-idC-first} or \ref{fig:subdisk-idD-idC-last}, or a result component of \autoref{fig:subdisk-idD-all-in} where the first angle of the discrete immersed disk is an $ α_3 $ or the final angle of the disk is an $ α_4 $.
\item \emph{CR result component} otherwise.
\end{itemize}
The classes of DS, DW, ID and CR result components are denoted $ \DSr, \DWr, \IDr, \CRr ⊂ \PiTr $ respectively.
\end{definition}

\begin{remark}
We have chosen the acronyms to reflect the amount of degeneracy allowed in the subdisks: Subdisks of DS result components are “degenerate strips”. Subdisks of DW result components are “degenerate wedges”. Subdisks of ID result components are “identity degenerate”, having an identity output and one of the inputs lying infinitesimally close to it. CR are mostly regular and satisfy the “co-identity rule”.
\end{remark}

The remainder of this section is devoted to defining the notions of CR, ID, DS and DW disks. These four classes are subsets of $ \SLd $ and meant to be explicitly constructible: For every imaginable SL disk, the reader should be able to determine whether it concerns a CR, ID, DS or DW disk or non of those. Ultimately, we will prove that these very explicit classes of disks are precisely the images of $ \CRr $, $ \IDr $, $ \DSr $ and $ \DWr $ under $ \Subdisk $.

\begin{definition}
\label{def:subdisk-types-CRdisk}
A \emph{CR disk} is an SL disk all of whose segments are nonempty, with the exception that multiple stacked co-identity inputs connected by infinitesimally short $ \smooth L_i $ segments are allowed, as long as their zigzag curve is oriented clockwise with the disk. The class of CR disks is denoted $ \CRd $.
\end{definition}

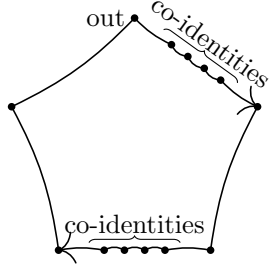
\begin{figure}
\centering
\begin{tikzpicture}
\path[draw, semithick, {To[scale=2]}-] (0, 0) coordinate(h3) to[bend left=10] ++(0.6, 0) coordinate(m1) to[bend left=30] ++(0.26666, 0) coordinate(m2) to[bend left=30] ++(0.26666, 0) coordinate (m3) to[bend left=30] ++(0.266666, 0) coordinate (m4) to[bend left=10] ++(0.6, 0) coordinate(h2) to[bend left=10] ++(72:2) coordinate(h1);
\path[draw, semithick, {To[scale=2]}-] (h1) to[bend left=10] ++(144:0.6) coordinate(m5) to[bend left=30] ++(144:0.26666) coordinate (m6) to[bend left=30] ++(144:0.26666) coordinate(m7) to[bend left=30] ++(144:0.266666) coordinate (m8) to[bend left=10] ++(144:0.6) coordinate(t) to[bend left=10] ++(216:2) coordinate(h4) to[bend left=10] (0, 0);
\path[fill] \foreach \i in {h1, h2, h3, h4, m1, m2, m3, m4, m5, m6, m7, m8, t} {(\i) circle[radius=0.05]};
\path (t) node[left] {out};
\path[draw, decorate, decoration={brace}] ($ (m1) + (-0.2, 0.1) $) to node[midway, above] {co-identities} ($ (m4) + (0.2, 0.1) $);
\path[draw, decorate, decoration={brace}] ($ (m8) + (144:0.2) + (54:0.1) $) to node[midway, sloped, above]  {co-identities} ($ (m5) + (144:-0.2) + (54:0.1) $);
\end{tikzpicture}
\caption{This picture depicts the schematic of CR disks. The specific CR disk depicted here has twelve inputs, of which four are of type B or C and eight are co-identities. The eight co-identities come in two stacks, each consisting of four co-identities lying infinitesimally close to each other.}
\label{fig:subdisk-types-CRdisk}
\end{figure}

\begin{remark}
The behavior of CR disks is depicted in \autoref{fig:subdisk-types-CRdisk}. We remark that in any CR disk, whenever a co-identity appears in the angle cut just before or after an intersection of type B, it appears only once due to arrow directions.
\end{remark}

\begin{definition}
\label{def:subdisk-types-IDdisk}
An \emph{ID disk} is an SL disk satisfying the following conditions:
\begin{itemize}
\item The output is the identity of a zigzag path,
\item Precisely one input, the \emph{degenerate input}, is infinitesimally close to the output,
\item The degenerate input is of type B or C,
\item The disk becomes CR upon excision of the output and substitution of the output mark by the degenerate input,
\item In case of a degenerate B input, it precedes respectively succeeds the output mark if $ L_1 $ is oriented clockwise respectively counterclockwise with the disk,
\item In case of a degenerate C input, the source zigzag path of the degenerate input is counterclockwise and the target zigzag path is clockwise.
\end{itemize}
The class of ID disks is denoted $ \IDd $.
\end{definition}

\begin{figure}
\centering
\begin{subfigure}{0.3\linewidth}
\centering
\begin{tikzpicture}
\path (2, 0) -- ++(72:2) coordinate (A) -- ++(144:2) coordinate (B);
\path[draw, semithick, {To[scale=2]}-] (0, 0) coordinate (h3) to[bend left=5] ++(0.6, 0) coordinate (m1) to[bend left=30] ++(0.26666, 0) coordinate (m2) to[bend left=30] ++(0.2666666, 0) coordinate (m3) to[bend left=30] ++(0.26666, 0) coordinate (m4) to[bend left=5] ++(0.6, 0) coordinate (h2) to[bend left=5] ++(72:1.8) coordinate (h1);
\path[draw, semithick, {To[scale=2]}-] (h1) to[bend left=10] ($ (A)!0.9!(B) $) coordinate (t) to (B) coordinate (h5);
\path[draw, semithick, {To[scale=2]}-] (h5) to ++(216:2) coordinate (h4) to (0, 0);
\path[fill] \foreach \i in {h1, h2, h3, h4, h5, t, m1, m2, m3, m4} {(\i) circle[radius=0.05]};
\path (t) node[right] {out};
\path[draw, decorate, decoration={brace}] ($ (m1) + (-0.2, 0.1) $) to node[midway, above] {\small co-identities} ($ (m4) + (0.2, 0.1) $);
\end{tikzpicture}
\caption{Clockwise B input}
\label{fig:subdisk-types-IDdiskCB}
\end{subfigure}
\begin{subfigure}{0.3\linewidth}
\centering
\begin{tikzpicture}
\path (2, 0) -- ++(72:2) -- ++(144:2) coordinate (A) -- ++(216:2) coordinate (B);
\path[draw, semithick, {To[scale=2]}-] (0, 0) coordinate (h3) to[bend left=5] ++(0.6, 0) coordinate (m1) to[bend left=30] ++(0.26666, 0) coordinate (m2) to[bend left=30] ++(0.2666666, 0) coordinate (m3) to[bend left=30] ++(0.26666, 0) coordinate (m4) to[bend left=5] ++(0.6, 0) coordinate (h2) to[bend left=5] ++(72:1.8) coordinate (h1);
\path[draw, semithick, -{To[scale=2]}] (h1) to[bend left=10] (A) coordinate (h0);
\path[draw, semithick, -{To[scale=2]}] (h0) to ($ (A)!0.1!(B) $) coordinate (t) to[bend left=5] (B) coordinate (h5);
\path[draw, semithick] (h5) to (0, 0);
\path[draw, semithick] (h4) to (0, 0);
\path[fill] \foreach \i in {h0, h1, h2, h3, h5, t, m1, m2, m3, m4} {(\i) circle[radius=0.05]};
\path (t) node[left] {out};
\path[draw, decorate, decoration={brace}] ($ (m1) + (-0.2, 0.1) $) to node[midway, above] {\small co-identities} ($ (m4) + (0.2, 0.1) $);
\end{tikzpicture}
\caption{Counterclockwise B input}
\label{fig:subdisk-types-IDdiskCCB}
\end{subfigure}
\begin{subfigure}{0.3\linewidth}
\centering
\begin{tikzpicture}
\path (2, 0) -- ++(72:2) coordinate (A) -- ++(144:2) coordinate (B);
\path[draw, semithick, {To[scale=2]}-] (0, 0) coordinate (h3) to[bend left=5] ++(0.6, 0) coordinate (m1) to[bend left=30] ++(0.26666, 0) coordinate (m2) to[bend left=30] ++(0.2666666, 0) coordinate (m3) to[bend left=30] ++(0.26666, 0) coordinate (m4) to[bend left=5] ++(0.6, 0) coordinate (h2) to[bend left=5] ++(72:1.8) coordinate (h1);
\path[draw, semithick, {To[scale=2]}-] (h1) to[bend left=10] ($ (A)!0.9!(B) $) coordinate (t) to (B) coordinate (h5);
\path[draw, semithick, -{To[scale=2]}] (h5) to ++(216:2) coordinate (h4);
\path[draw, semithick] (h4) to (0, 0);
\path[fill] \foreach \i in {h1, h2, h3, h4, h5, t, m1, m2, m3, m4} {(\i) circle[radius=0.05]};
\path (t) node[right] {out};
\path[draw, decorate, decoration={brace}] ($ (m1) + (-0.2, 0.1) $) to node[midway, above] {\small co-identities} ($ (m4) + (0.2, 0.1) $);
\end{tikzpicture}
\caption{C input}
\label{fig:subdisk-types-IDdiskC}
\end{subfigure}
\caption{This picture depicts the schematic of ID disks, categorized according to whether the degenerate input is of type B or C. Each of the specific ID disks depicted here has nine inputs, of which five are of type B or C and four consist of a stack of co-identities. The degenerate input is the one at the top corner. For the case of degenerate B input, we have depicted both the clockwise and the counterclockwise case. For the case of degenerate C input, we have depicted only the case where the degenerate input precedes the output mark.}
\label{fig:subdisk-types-IDdisk}
\end{figure}
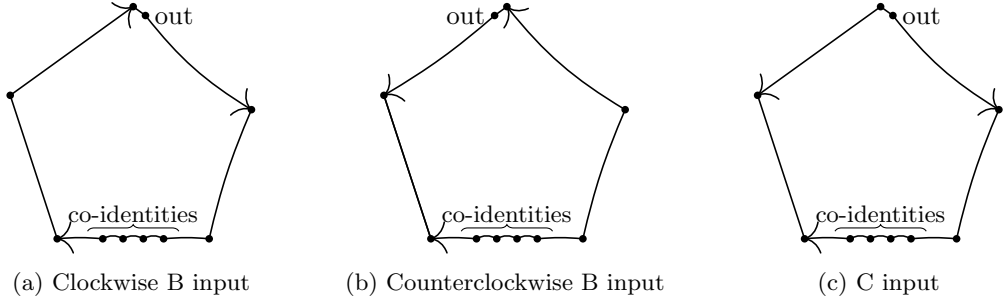

\begin{remark}
The two conditions of \autoref{def:subdisk-types-IDdisk} specific to the B and C case can be formulated in more relaxed terms. In case of a degenerate B input, the three zigzag paths given by the source and target of the degenerate input and the output all have the same orientation. We can make the precedence of degenerate input and output therefore dependent on any of the three, instead of $ L_1 $. In case of a degenerate C input, the requirement regarding orientations equivalently requires that the source and target zigzag path of the degenerate input are always oriented “towards” the disk, instead of “away from” the disk. This is visually depicted in \autoref{fig:subdisk-types-IDdisk}.
\end{remark}

\begin{definition}
\label{def:subdisk-types-DSdisk}
A \emph{DS disk} is an immersed strip fitting into one of the two digons bounded by a zigzag curve $ \smooth L $ and its Hamiltonian deformation $ \smooth L' $. More precisely, the strip is a 4-gon bounded by $ \smooth L $, $ \smooth L' $ and two (indexed) arcs $ a $ and $ b $ lying on $ L $. The arc $ b $ is the one lying closer to the co-identity. Two inputs lie on the midpoint of the arc $ a $ and one input on the midpoint of the arc $ b $. The output mark lies on the midpoint of the arc $ b $. There are corner cases in which additional conditions apply:
\begin{itemize}
\item If $ a = a_0 $, then $ L $ is oriented away from the co-identity.
\item If $ a = b $, then either (a) $ L $ is oriented away from the co-identity and turns left at $ a = b $, and the input at the $ b $ side is odd/final, or (b) $ L $ is oriented towards the co-identity and turns right at $ a = b $, or (c) $ L $ is oriented towards the co-identity and turns left at $ a = b $, and the $ b $ side input is odd/final.
\item If $ a = b = a_0 $, then both conditions must be met: $ L $ is oriented away from the co-identity, turns left at $ a = b $ and the input at the $ b $ side is odd/final.
\end{itemize}
The class of DS disks is denoted $ \DSd $.
\end{definition}

The behavior of DS disks is best observed in \autoref{fig:subdisk-degenerate}. In the definition of DS disks, we have used terminology that $ \smooth L $ may be oriented towards or away from the co-identity. Indeed, a strip lies in the digon between identity and co-identity. This brings a distinction whether $ \smooth L $ is oriented towards the co-identity and away from the identity, or away from the co-identity and towards the identity. Of course, the whole definition with its corner cases is designed to capture precisely the result components of \autoref{fig:subdisk-degenerate}.

\begin{definition}
\label{def:subdisk-types-DWdisk}
A \emph{DW disk} is one of the following:
\begin{itemize}
\item A 3-gon sitting between a zigzag curve $ \smooth L $ and its Hamiltonian deformation, bounded on one side by an arc $ a $ of $ L $ and on the other side by the co-identity. The output mark is placed at the co-identity. It is allowed that $ a = a_0 $ if $ L $ is oriented away from the co-identity.
\item A 4-gon, obtained from the first option by inserting an additional co-identity input infinitesimally preceding the output. The condition is that $ L $ is oriented away from the co-identity and that $ a ≠ t(α_0) $.
\item A 4-gon, obtained from the first option by inserting an additional co-identity input infinitesimally succeeding the output. The condition is that $ L $ is oriented towards the co-identity.
\end{itemize}
The set of DW disks is denoted by $ \DWd $.
\end{definition}

The behavior of DW disks is best observed in \autoref{fig:subdisk-alpha0}. The definition distinguishes three types of DW disks. To be more precise with the conditions, observe that a DW disk the second type is allowed to have $ a = a_0 $ while a DW disk of the third type is required to have $ a ≠ a_0 $. In the second type, the arc $ a $ is supposed to be not the tail arc $ t(α_0) $ of the co-identity angle $ α_0 $. In the third type, the assertion $ a ≠ t(α_0) $ holds automatically, since the co-identity angle $ α_0 $ is located in a counterclockwise polygon of $ Q $ and $ L $ is supposed to be oriented towards the co-identity. All DW disks have infinitesimal area, but precisely two nonempty zigzag segments. In fact, the distance between the midpoint of the arc $ a $ and the midpoint of the co-identity angle $ α_0 $ is at least half an angle in size.

We have constructed the definitions of CR, ID, DS and DW disks such that the subdisk of a CR result component is a CR disk, and so on:

\begin{lemma}
\label{th:subdisk-types-shape}
The subdisk of a CR, ID, DS or DW result component is a CR, ID, DS or DW disk, respectively.
\end{lemma}

\begin{proof}
The inspection is performed in \autoref{sec:classification-subdiskshape}.
\end{proof}

In fact, we will prove and discuss later that our definition of CR, ID, DS and DW disks is also sharp: Every CR, ID, DS and DW disk is actually reached as a subdisk of some result component.

\begin{figure}
\centering
\begin{tikzpicture}
\path[draw] (0, 0) -- ++(300:1) coordinate[midway] (med)  -- ++(60:1) -- ++(300:1) -- ++(60:1) -- ++(300:1) -- ++(60:1) -- ++(300:1) -- ++(60:1);
\path[draw, semithick] (med) -- ++(240:0.12) coordinate (down);
\path[draw, semithick] (med) -- ++(60:0.12) coordinate (up);
\path[draw, semithick, fill=gray, fill opacity=0.5] (down) -- (up) -- ++(right:3.8) -- ++(240:0.24) -- cycle;
\path[fill] (up) circle[radius=0.05] node[above] {$ h_2 $};
\path[fill] (down) circle[radius=0.05] node[left] {$ h_1 $};
\end{tikzpicture}
\caption{This picture explains that one expects only very few non-transversal disks among zigzag curves in the Fukaya category. We depict a smooth immersed disk bounded by zigzag curves and assume it has two transversal intersections $ h_1 $, $ h_2 $ on the boundary which lie infinitesimally close to each other. Due to the zigzag nature and transversality, $ h_1 $ and $ h_2 $ must be the even and odd intersection points located at the midpoint of one single arc of $ Q $. This means that the source zigzag path of $ h_1 $ is the target zigzag path of $ h_2 $ and the entire disk is then very thin.}
\label{fig:subdisk-types-thin}
\end{figure}
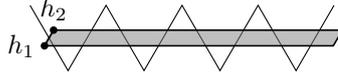

\begin{remark}
It is very pleasant that most subdisks are rather regular in the sense that their zigzag curve segments are non-empty. The only irregularities are found in stacked co-identities of CR and ID disks, the degenerate output of ID disks, and the two irregular types of DS and DW disks. Viewed geometrically, this is not really a surprise: In \autoref{fig:subdisk-types-thin} we argue that smooth immersed disks with non-transversal intersections lying infinitesimally close to each other are very thin. The conclusion is that within the Fukaya category, one expects only very few irregular disks between zigzag curves. The DS and DW disks provide the exact representation-theoretic witness of this phenomenon.
\end{remark}

\subsection{The minimal model}
\label{sec:subdisk-minmodel}
In this section, we tie together the computation of $ \H\DefZigzagCat $: Its $ A_∞ $-structure is defined in terms of π-trees. All the result components of a π-tree can again be matched with CR, ID, DS and DW disks:

\begin{center}
\begin{tikzpicture}
\path (-6, 0) node[align=center] {$ \H\DefZigzagCat $ \\ Minimal model};
\path[draw, {To[scale=1.5]}-{To[scale=1.5]}] (-4, -0.2) -- (-2, -0.2) node[midway, above] {Kadeishvili};
\path (0, 0) node[align=center] {$ \PiTr $ \\ Result components};
\path[draw, {To[scale=1.5]}-{To[scale=1.5]}] (2, -0.2) -- (4, -0.2) node[midway, above] {$ \Subdisk $};
\path (6, 0) node[align=center] {CR, ID, DS, DW \\ Immersed disks};
\end{tikzpicture}
\end{center}

This correspondence allows us to express the minimal model $ \H\DefZigzagCat $ in terms of disks. The present section is meant to spell out the details and provide intuition.

Our first step is to get more grip on the subdisk mapping $ \Subdisk: \PiTr → \SLd $. We have already seen in \autoref{sec:subdisk-types} that $ \Subdisk $ sends CR result components to CR disks, and so on. In the following lemma, we affirm that all CR, ID, DS and DW disks are actually reached by $ \Subdisk $.

\begin{lemma}
\label{th:subdisk-minmodel-bijection}
The classes of CR, ID, DS and DW result components are disjoint, as are the classes of CR, ID, DS and DW disks. The subdisk mapping $ \Subdisk $ bijectively sends each of the four result component classes to its disk counterpart:

\begin{center}
\begin{tikzpicture}
\path (0, 0) node (A) {$ \CRr $} (1, 0) node {$ \disjoint $} (2, 0) node (B) {$ \IDr $} (3, 0) node {$ \disjoint $} (4, 0) node (C) {$ \DSr $} (5, 0) node {$ \disjoint $} (6, 0) node (D) {$ \DWr $};
\path (0, -1) node (A') {$ \CRd $} (1, -1) node {$ \disjoint $} (2, -1) node (B') {$ \IDd $} (3, -1) node {$ \disjoint $} (4, -1) node (C') {$ \DSd $} (5, -1) node {$ \disjoint $} (6, -1) node (D') {$ \DWd $};
\foreach \i in {A, B, C, D} {\path[draw, ->] (\i.south) -- (\i'.north) node[midway, left] {$ \Subdisk $} node[midway, right] {\rotatebox{90}{\large $ \sim $}};};
\path (-2, 0) node (1) {$ \PiTr $} ($ (1.east)!0.5!(A.west) $) node {$ = $};
\path (8, -1) node (2) {$ \SLd $} ($ (2.west)!0.5!(D'.east) $) node {$ ⊂ $};
\end{tikzpicture}
\end{center}
\end{lemma}

\begin{proof}
The inclusions are checked and an explicit inverse map is constructed in \autoref{sec:classification}.
\end{proof}

At this point, we can already describe the minimal model $ \H\DefZigzagCat $ in a rigged way by means of disks: Let $ r $ be a result component of a π-tree. Then its subdisk $ \Subdisk(r) $ comes with a designated output mark, in particular we can read off its output morphism $ \disktarget(\Subdisk(r)) $. In fact, the value of $ r $ is equal to $ \disktarget(\Subdisk(r)) $, at least when sign and $ q $-parameters are stripped off. For example, a subdisk of an $ α_3 + α_4 $ result component has output mark at this very same $ α_3 + α_4 $ morphism, by construction. Even though we currently have to recover signs and $ q $-parameters from result component instead of disk, this enables us to largely describe the product in terms of disks:
\begin{align}
\label{eq:subdisk-minmodel-rigged}
\begin{split}
μ^N (h_N, …, h_1) &= \sum_{\substack{r ∈ \PiTr \\ r \text{ has inputs } h_1, …, h_N}} r \\
&= \sum_{\substack{D ∈ \CRd \disjoint \IDd \disjoint \DSd \disjoint \DWd \\ D \text{ has inputs } h_1, …, h_N}} \Subdisk^{-1} (D) \\
&= \sum_{\substack{D ∈ \CRd \disjoint \IDd \disjoint \DSd \disjoint \DWd \\ D \text{ has inputs } h_1, …, h_N}} (-1)^{\sgn(\Subdisk^{-1} (D))} \qparam(\Subdisk^{-1} (D)) \disktarget(D).
\end{split}
\end{align}
Here $ \disktarget(D) $ denotes the output morphism of the disk $ D $, and $ \sgn(r) $ and $ \qparam(r) $ temporarily denote the sign and $ q $-parameter of a result component. As announced, the sign $ (-1)^{\sgn(\Subdisk^{-1} (D))} $ and $ q $-parameter $ \qparam(\Subdisk^{-1} (D)) $ in this formula are only recovered from the result component $ \Subdisk^{-1} (D) $ instead of $ D $ itself.

Our next step is to write signs and $ q $-parameters in terms of $ D $ instead of recovering them from $ \Subdisk^{-1} (D) $. In fact, \autoref{th:subdisk-minmodel-signs} will show that the sign is precisely the Abouzaid sign of $ D $, which the reader may recall from the context of Fukaya categories in \autoref{sec:prelim-fukaya-pre}. Moreover, the $ q $-parameter is precisely the product of all punctures covered by $ D $, counted with multiplicities. Before we make these statements, let us fix the Abouzaid sign terminology in our present context:

\begin{definition}
\label{def:subdisk-minmodel-Abouzaid}
Let $ D $ be an SL disk. Then its \emph{Abouzaid sign} $ \Abouzaid(D) ∈ ℤ/2ℤ $ is the sum of all $ \# $ signs around $ D $, plus the number of odd inputs $ h_i: L_i → L_{i+1} $ where $ L_{i+1} $ is oriented counterclockwise with $ D $, plus one if its output $ t: L_1 → L_{N+1} $ is odd and $ L_{N+1} $ is oriented counterclockwise. The \emph{$ \bm{q} $-parameter} $ \punctures(D) ∈ ℂ⟦Q_0⟧ $ of $ D $ is defined as the product of all punctures covered by $ D $, counted with multiplicities.
\end{definition}

\begin{lemma}
\label{th:subdisk-minmodel-signs}
Let $ r ∈ \PiTr $ be the result component of a π-tree. Then its sign is equal to the Abouzaid sign of its subdisk $ \Subdisk(r) $ and its $ q $-parameter $ ∈ ℂ⟦Q_0⟧ $ is equal to the product of all punctures covered by $ \Subdisk(r) $, counted with multiplicities:

\begin{center}
\begin{tikzpicture}
\path (0, 0.1) node {\textbf{In terms of} $ \mathbf{r} $} (7, 0.1) node {\textbf{In terms of} $ \mathbf{D} $};
\path (0, -0.5) node {Sign $ \sgn(r) ∈ ℤ/2ℤ $} (7, -0.5) node {Sign $ \Abouzaid(D) ∈ ℤ/2ℤ $};
\path (0, -1) node {$ q $-parameter $ \qparam(r) ∈ ℂ⟦Q_0⟧ $} (7, -1) node {$ q $-parameter $ \punctures(D) ∈ ℂ⟦Q_0⟧ $};
\path (3.5, -0.5) node {\LARGE $ = $};
\path (3.5, -1) node {\LARGE $ = $};
\end{tikzpicture}
\end{center}
\end{lemma}

\begin{proof}
Both checks can be performed in an inductive fashion. The signs checks are detailed in \autoref{sec:classification-signs}. The checks for $ q $-parameters are easier and left to the reader.
\end{proof}

With the help of this lemma, we are ready to translate the rigged formula \eqref{eq:subdisk-minmodel-rigged} into the soothing description of the minimal model purely in terms of disks. For simplicity we denote the identity element of a zigzag path $ L $ by $ \id_L = \sum_a \id_a $.

\begin{theorem}
\label{th:subdisk-minmodel-th}
Let $ Q $ be a geometrically consistent dimer. Regard the category $ \DefZigzagCat ⊂ \Tw'\Gtl_q Q $ of deformed zigzag paths according to \autoref{conv:alpha0-direction}. Then the $ A_∞ $-structure of the minimal model $ \H\DefZigzagCat $ is described as follows:
\begin{itemize}
\item The curvature and differential vanish:
\begin{equation*}
μ^0_{\H\DefZigzagCat} = μ^1_{\H\DefZigzagCat} = 0.
\end{equation*}
\item The minimal model is unital: For every cohomology basis element $ h: L_1 → L_2 $ we have
\begin{align*}
μ^{≥3}_{\H\DefZigzagCat} \left(…, \id_{L_1}, …\right) &= 0, \\
μ^2_{\H\DefZigzagCat} \left(h, \id_{L_1}\right) = (-1)^{|h|} μ^2_{\H\DefZigzagCat} \left(\id_{L_2}, h\right) &= h.
\end{align*}
\item The products are given by CR, ID, DS and DW disks: Let $ N ≥ 2 $ and let $ h_1, …, h_N $ be a sequence of non-identity cohomology basis morphisms with $ h_i: L_i → L_{i+1} $. Then their product is given by
\begin{equation*}
μ^N (h_N, …, h_1) = \sum_{\substack{D ∈ \CRd \disjoint \IDd \disjoint \DSd \disjoint \DWd \\ D \text{ has inputs } h_1, …, h_N}} (-1)^{\Abouzaid(D)} \punctures(D) \disktarget(D).
\end{equation*}
\end{itemize}
\end{theorem}

\begin{proof}
This is a summary of our journey. As we have observed earlier, the minimal model $ \H\DefZigzagCat $ has vanishing differential and curvature. It is also unital with the same identities as $ \ZigzagCat $. When $ h_1, …, h_N $ are cohomology basis elements, the rigged formula \eqref{eq:subdisk-minmodel-rigged} gives
\begin{align*}
μ^N (h_N, …, h_1) &= \sum_{\substack{D ∈ \CRd \disjoint \IDd \disjoint \DSd \disjoint \DWd \\ D \text{ has inputs } h_1, …, h_N}} (-1)^{\sgn(\Subdisk^{-1} (D))} \qparam(\Subdisk^{-1} (D)) \disktarget(D) \\
&= \sum_{\substack{D ∈ \CRd \disjoint \IDd \disjoint \DSd \disjoint \DWd \\ D \text{ has inputs } h_1, …, h_N}} (-1)^{\Abouzaid(D)} \punctures(D) \disktarget(D).
\end{align*}
In the second row, we have inserted \autoref{th:subdisk-minmodel-signs}. This finishes the proof.
\end{proof}

\subsection{Main result}
\label{sec:subdisk-main}
In this section, we present our main result. It ties together the “discrete relative Fukaya category” $ \H\DefZigzagCat $ and the “smooth relative Fukaya category” $ \relFuk Q $:

\begin{center}
\begin{tikzpicture}
\path (0, 0) node[align=center] (A) {\textbf{Discrete relative} \\ $ (\H\DefZigzagCat)_{\transversal} $} (8, 0) node[align=center] (B) {\textbf{Smooth relative} \\ $ \relpreFuk Q \restr{\Ob\ZigzagCat} $};
\path[draw, <->] ($ (A.east)!0.2!(B.west) $) -- ($ (A.east)!0.8!(B.west) $);
\end{tikzpicture}
\end{center}

The starting point on the discrete side is the explicit description of the minimal model $ \H\DefZigzagCat $ due to \autoref{th:subdisk-minmodel-th}. The starting point on the smooth side is the explicit description of the subcategory $ \relpreFuk Q \restr{\Ob\ZigzagCat} $ from \autoref{th:prelim-fukaya-zigzag-restr}. The main result entails a strict isomorphism between the transversal part $ (\H\DefZigzagCat)_{\transversal} $ on one side and $ \relpreFuk Q \restr{\Ob\ZigzagCat} $ on the other side. In what follows, we recall a few specific properties of $ \relpreFuk \restr{\Ob\ZigzagCat} $ and a few similarities with $ \H\DefZigzagCat $.

\begin{remark}
In \autoref{sec:fukaya} we have elaborated on the construction of Fukaya categories. More specifically, we have defined the categories $ \preFuk Q $, $ \Fuk Q $, $ \relpreFuk Q $, $ \relFuk Q $ and their subcategories given by zigzag curves. Most importanty, recall from \autoref{def:prelim-fukaya-relpre-def} that $ \relpreFuk Q $ denotes the relative Fukaya pre-category of $ Q $. In \autoref{sec:prelim-fukaya-zigzag}, we have provided an extensive elaboration on how zigzag paths can be interpreted as objects in $ \relpreFuk Q $. In particular, a zigzag path $ L ∈ \H\DefZigzagCat $ corresponds to a zigzag curve $ \smooth L ∈ \relpreFuk Q $. Recall from \autoref{def:prelim-fukaya-zigzag-restr} that $ \relpreFuk Q \restr{\Ob\ZigzagCat} $ is the $ A_∞ $-pre-category defined as the subcategory of $ \relpreFuk Q $ given by zigzag curves $ \smooth L $, together with the spin structure dictated by $ L ∈ \ZigzagCat $. As we have seen in \autoref{th:prelim-fukaya-zigzag-transversal}, a sequence of zigzag curves $ (\smooth L_1, …, \smooth L_{N+1}) $ is transversal in $ \relpreFuk Q \restr{\Ob\ZigzagCat} $ if and only if the zigzag paths $ L_i $ are pairwise distinct. We have described the subcategory $ \relpreFuk Q \restr{\Ob\ZigzagCat} $ more explicitly in \autoref{th:prelim-fukaya-zigzag-restr}.
\end{remark}

\begin{remark}
In \autoref{th:prelim-fukaya-zigzag-Hom}, we have identified basis elements for the hom spaces $ \Hom_{\Fuk Q} (\smooth L_1, \smooth L_2) $ with intersection points of $ \smooth L_1 $ and $ \smooth L_2 $. In case $ \smooth L_1 = \smooth L_2 $, the intersection points only refer to the transversal self-intersections, plus the identity and co-identity self-intersections. In \autoref{sec:splitting-splitting}, we have seen that basis elements for $ \Hom_{\H\ZigzagCat} (L_1, L_2) $ are identified with intersection points between $ \smooth L_1 $ and $ \smooth L_2 $ as well. In case $ L_1 = L_2 $, the intersection points only refer to the transversal self-intersections, plus identity and co-identity self-intersections:

\begin{center}
\begin{tabular}{@{}ccc@{}}
\textbf{Category of zigzag paths} & \textbf{Geometry} & \textbf{Fukaya category} \\\hline
Zigzag path & Zigzag curve & Zigzag curve \\
$ L $ & $ \smooth L $ & $ \smooth L $ \\\hline
Cohomology basis element & Intersection point & Basis element \\
$ h: L_1 → L_2 $ & $ p ∈ \smooth L_1 ∩ \smooth L_2 $ & $ p: \smooth L_1 → \smooth L_2 $
\end{tabular}
\end{center}
\end{remark}

\begin{remark}
\label{th:subdisk-main-homspaces}
As laid out in \autoref{th:prelim-fukaya-zigzag-relHom}, the relative Fukaya category $ \relFuk Q $ is a deformation of $ \Fuk Q $. As such, its hom spaces are the $ B $-enlargement of the hom spaces of $ \Fuk Q $, see also \autoref{th:prelim-fukaya-zigzag-relHom}:
\begin{equation*}
\Hom_{\relFuk Q} (\smooth L_1, \smooth L_2) = B \htensor \Hom_{\Fuk Q} (\smooth L_1, \smooth L_2).
\end{equation*}
Similarly, $ \H\DefZigzagCat $ is a deformation of $ \H\ZigzagCat $ by construction. Its hom spaces are the $ B $-enlargement of the hom spaces of $ \H\ZigzagCat $:
\begin{equation*}
\Hom_{\H\DefZigzagCat} (L_1, L_2) = B \htensor \Hom_{\H\ZigzagCat} (L_1, L_2).
\end{equation*}
The identification of the basis elements of $ \Hom_{\Fuk Q} (\smooth L_1, \smooth L_2) $ and $ \Hom_{\H\DefZigzagCat} (L_1, L_2) $ provides an explicit $ B $-linear identification of the hom spaces $ \Hom_{\relFuk Q} (\smooth L_1, \smooth L_2) $ and $ \Hom_{\H\DefZigzagCat} (L_1, L_2) $.
\end{remark}

In \autoref{th:subdisk-main-onlyCR}, we examine CR, DS, ID and DS disks in the case that the sequence of input zigzag paths is transversal. The notable outcome is that only CR disks remain, which can in turn be interpreted directly as smooth immersed disks. This establishes the desired link between the minimal model $ \H\DefZigzagCat $ and $ \relpreFuk Q \restr{\Ob\ZigzagCat} $ which we will expand in \autoref{th:subdisk-main-th}.

\begin{lemma}
\label{th:subdisk-main-onlyCR}
Let $ L_1, …, L_{N+1} $ be a sequence of zigzag paths in $ Q $ such that $ \smooth L_1, …, \smooth L_{N+1} $ is a transversal sequence. Let $ h_i: L_i → L_{i+1} $ for $ 1 ≤ i ≤ N $ and $ h: L_1 → L_{N+1} $ be cohomology basis elements in $ \H\ZigzagCat $. Denote by $ p_i: \smooth L_i → \smooth L_{i+1} $ and $ p: \smooth L_1 → \smooth L_{N+1} $ the corresponding basis elements in $ \Fuk Q \restr{\Ob\ZigzagCat} $. Then:
\begin{enumerate}
\item There are no ID, DS and DW disks with inputs $ h_1, …, h_N $.
\item There is a bijection
\begin{align*}
Φ: \left\{\substack{\text{CR disks} \\ \text{with inputs } h_1, …, h_N \\ \text{and output } h}\right\} \verylongisoto \left\{\substack{\text{Smooth immersed disks} \\ \text{with inputs } p_1, …, p_N \\ \text{and output } p}\right\}.
\end{align*}
\item The Abouzaid signs agree: $ \Abouzaid(D) = \Abouzaid(Φ(D)) $.
\item The $ q $-paramaters agree: $ \punctures(D) = \punctures(Φ(D)) $.
\end{enumerate}
\end{lemma}

\begin{proof}
We explain the four statements one after another. For the first statement, let $ D $ be an ID, DS or DW disk with inputs $ h_1, …, h_N $. Then necessarily at least two of the zigzag curves $ L_1, …, L_{N+1} $ are equal. Therefore $ (\smooth L_1, …, \smooth L_{N+1}) $ is not a transversal sequence, in contradiction with the assumption. This shows that there are no ID, DS or DW disks with input $ h_1, …, h_N $. In other words, there can only be CR disks with inputs $ h_1, …, h_N $ among the four types of disks.

For the second statement, pick a CR disk with inputs $ h_1, …, h_N $ and output $ h $. Since all zigzag paths $ L_1, …, L_{N+1} $ are pairwise distinct, the sequence $ h_1, …, h_N $ does not contain any co-identities. Therefore all zigzag curve segments involved in the CR disk $ D $ are non-empty. This way $ D $ immediately constitutes a smooth immersed disk in the sense of \autoref{def:prelim-fukaya-relpre-disk}. We denote this smooth immersed disk by $ Φ(D) $. The smooth immersed disk $ Φ(D) $ has inputs $ p_1, …, p_N $ and output $ p $, precisely as desired. This sets up the desired mapping $ Φ $. The map $ Φ $ is clearly injective, since a CR disk contains as much information about the polygon immersion $ D: P_{N+1} → |Q| $ as does a smooth immersed disk. For instance, the two notions of CR disks and smooth immersed disks both identify immersions related by reparametrization. The map $ Φ $ is also surjective, since a smooth immersed disk with inputs $ p_1, …, p_N $ and output $ p $ can immediately be interpreted as a CR disk. This shows that $ Φ $ is a bijection.

For the third statement, let $ D $ be a CR disk with inputs $ h_1, …, h_N $ and output $ h $. According to \autoref{def:subdisk-minmodel-Abouzaid}, the Abouzaid sign $ \Abouzaid(D) ∈ ℤ/2ℤ $ is the sum of all $ \# $ signs on the boundary of $ D $, plus the number of odd inputs $ h_i $ where $ L_{i+1} $ is oriented counterclockwise with $ D $, plus one if the output $ h: L_1 → L_{N+1} $ is odd and $ L_{N+1} $ is oriented counterclockwise. This is exactly the same as the definition of the Abouzaid sign of $ Φ(D) $, see \autoref{def:prelim-fukaya-relpre-Abouzaid} and \ref{def:prelim-fukaya-pre-Abouzaid}. This shows $ \Abouzaid(D) = \Abouzaid(Φ(D)) $.

For the fourth statement, let $ D $ be a CR disk with inputs $ h_1, …, h_N $ and output $ h $. According to \autoref{def:subdisk-minmodel-Abouzaid}, the $ q $-parameter $ \punctures(D) ∈ ℂ⟦Q_0⟧ $ is the product of all punctures covered by $ D $, counting punctures multiple times if they are covered multiple times. This is exactly the same as the definition of the $ q $-parameter of $ Φ(D) $, see \autoref{def:prelim-fukaya-relpre-Abouzaid}. This finishes the proof.
\end{proof}

Our main theorem shows that the transversal part of $ \H\DefZigzagCat $ agrees with the subcategory $ \relpreFuk Q \restr{\Ob\ZigzagCat} $ of the relative Fukaya pre-category. For sake of logical independence, we repeat the setup here: The starting point is a geometrically consistent dimer $ Q $. We assume \autoref{conv:alpha0-direction}. We denote by $ \DefZigzagCat ⊂ \Tw'\Gtl_q Q $ the category of deformed zigzag paths according to \autoref{def:deformed-zigzagcat-def}. We denote by $ \H\DefZigzagCat $ the minimal model of $ \DefZigzagCat $, described explicitly in \autoref{th:subdisk-minmodel-th}. We denote by $ \relpreFuk Q \restr{\Ob\ZigzagCat} $ the subcategory of the relative Fukaya pre-category of $ Q $, described explicitly in \autoref{th:prelim-fukaya-zigzag-restr}. We denote by $ (\H\DefZigzagCat)_{\transversal} $ the transversal part of $ \H\DefZigzagCat $ with respect to $ \relpreFuk Q \restr{\Ob\ZigzagCat} $, as defined in \autoref{def:prelim-fukaya-agree}. The notion of strict isomorphism of deformed $ A_∞ $-pre-categories is provided in \autoref{def:prelim-fukaya-defiso}. Under this terminology, we state our main theorem as follows:

\begin{theorem}
\label{th:subdisk-main-th}
Let $ Q $ be a geometrically consistent dimer and assume \autoref{conv:alpha0-direction}. Then there is a strict isomorphism of deformed $ A_∞ $-pre-categories
\begin{equation*}
F_q: (\H\DefZigzagCat)_{\transversal} \verylongisoto \relpreFuk Q \restr{\Ob\ZigzagCat}.
\end{equation*}
The functor $ F_q $ sends a zigzag path $ L ∈ \H\DefZigzagCat $ to the associated zigzag curve $ \smooth L $ and a cohomology basis element $ h: L_1 → L_2 $ to the associated intersection point $ p: \smooth L_1 → \smooth L_2 $.
\end{theorem}

\begin{proof}
This follows directly from \autoref{th:subdisk-minmodel-th}, but we state the details. The starting point is the description of the $ A_∞ $-deformation $ \H\DefZigzagCat $ of $ \H\ZigzagCat $ from \autoref{th:subdisk-minmodel-th} and the description of the $ A_∞ $-pre-category deformation $ \relpreFuk Q \restr{\Ob\ZigzagCat} $ of $ \preFuk Q \restr{\Ob\ZigzagCat} $ from \autoref{th:prelim-fukaya-zigzag-restr}. We have detailed the definition of the transversal part $ (\H\DefZigzagCat)_{\transversal} $ with respect to $ \preFuk Q \restr{\Ob\ZigzagCat} $ in \autoref{def:prelim-fukaya-agree}.

To construct the functor $ F_q $ according to \autoref{def:prelim-fukaya-defiso}, we have to execute four steps: (1) to set up a bijection between the objects of $ \H\ZigzagCat $ and $ \preFuk Q \restr{\Ob\ZigzagCat} $, (2) to show that the transversal sequences of $ (\H\ZigzagCat)_{\transversal} $ and $ \preFuk Q \restr{\Ob\ZigzagCat} $ agree under the bijection on objects, (3) to set up a $ ℂ⟦Q_0⟧ $-linear degree $ 0 $ isomorphism between the hom spaces of $ (\H\DefZigzagCat)_{\transversal} $ and $ \relpreFuk Q \restr{\Ob\ZigzagCat} $, (4) to show that the higher products of $ (\H\DefZigzagCat)_{\transversal} $ and $ \relpreFuk Q \restr{\Ob\ZigzagCat} $ agree under the identification of hom spaces.

For step (1), the bijection between objects of $ \H\ZigzagCat $ and $ \preFuk Q \restr{\Ob\ZigzagCat} $ consists simply of mapping a zigzag path $ L ∈ \H\ZigzagCat $ to its associated zigzag curve $ \smooth L ∈ \preFuk Q \restr{\Ob\ZigzagCat} $:
\begin{align*}
F_q: \Ob(\H\ZigzagCat) &\verylongisoto \Ob(\preFuk Q \restr{\Ob\ZigzagCat}) \\
L &\verylongmapsto \smooth L.
\end{align*}

For step (2), we have to explain that the transversal sequences of $ (\H\ZigzagCat)_{\transversal} $ are precisely the transversal sequences of $ \preFuk Q \restr{\Ob\ZigzagCat} $ under the identification of $ L $ with $ \smooth L $. In fact, this is immediate from the definition of $ (\H\DefZigzagCat)_{\transversal} $ as transversal part of $ \H\DefZigzagCat $ with respect to $ \relpreFuk Q \restr{\Ob\ZigzagCat} $ under the identification of $ L $ with $ \smooth L $. Explicitly, a sequence $ (L_1, …, L_N) $ in $ (\H\ZigzagCat)_{\transversal} $ is by definition transversal if and only if $ (\smooth L_1, …, \smooth L_N) $ is transversal.

For step (3), we have to set up a $ ℂ⟦Q_0⟧ $-linear identification between the hom spaces of $ (\H\DefZigzagCat)_{\transversal} $ and $ \relpreFuk Q \restr{\Ob\ZigzagCat} $. In order to define this identification, let $ (L_1, L_2) ∈ (\H\ZigzagCat)^2_{\transversal} $. Then we have $ L_1 ≠ L_2 $. We now set up the identification by sending a cohomology basis element $ h: L_1 → L_2 $ to its associated intersection point $ p: \smooth L_1 → \smooth L_2 $, as described in \autoref{th:subdisk-main-homspaces}:
\begin{align*}
F_q^1: \Hom_{(\H\DefZigzagCat)_{\transversal}} (L_1, L_2) &\verylongisoto \Hom_{\relpreFuk Q \restr{\Ob\ZigzagCat}} (\smooth L_1, \smooth L_2) \\
h & \verylongmapsto p.
\end{align*}
For step (4), we have to show that $ F_q $ preserves the products of $ (\H\DefZigzagCat)_{\transversal} $ and $ \relpreFuk Q \restr{\Ob\ZigzagCat} $. Pick $ N ≥ 1 $ and let $ (L_1, …, L_{N+1}) ∈ (\H\ZigzagCat)^{N+1}_{\transversal} $. Pick basis elements $ h_i ∈ \Hom_{\H\ZigzagCat} (L_i, L_{i+1}) $ and let $ p_i ∈ \Hom_{\preFuk Q \restr{\Ob\ZigzagCat}} (\smooth L_i, \smooth L_{i+1}) $ be the associated intersection points. We get
\begin{align*}
F_q^1 (μ_{(\H\DefZigzagCat)_{\transversal}} (h_N, …, h_1)) &= F_q^1 \bigg(\sum_{\substack{D ∈ \CRd \disjoint \IDd \disjoint \DSd \disjoint \DWd \\ D \text{ has inputs } h_1, …, h_N}} (-1)^{\Abouzaid(D)} \punctures(D) \disktarget(D)\bigg) \\
&= \sum_{\substack{D ∈ \CRd \\ D \text{ has inputs } h_1, …, h_N}} (-1)^{\Abouzaid(D)} \punctures(D) F_q^1 (\disktarget(D)) \\
&= \sum_{h: L_1 → L_{N+1}} \sum_{\substack{D ∈ \CRd \\ D \text{ has inputs } h_1, …, h_N \\ \text{and output } h}} (-1)^{\Abouzaid(D)} \punctures(D) F_q^1 (h) \\
&= \sum_{p ∈ \smooth L_1 ∩ \smooth L_{N+1}} \sum_{D ∈ M_q (p_1, …, p_N, p)} (-1)^{\Abouzaid(Φ^{-1} (D))} \punctures(Φ^{-1} (D)) p \\
&= \sum_{p ∈ \smooth L_1 ∩ \smooth L_{N+1}} \sum_{D ∈ M_q (p_1, …, p_N, p)} (-1)^{\Abouzaid(D)} \punctures(D) p \\
&= μ_{\relpreFuk Q \restr{\Ob\ZigzagCat}} (p_N, …, p_1).
\end{align*}
In the first row, we have inserted the description of $ μ_{\H\DefZigzagCat} $ from \autoref{th:subdisk-minmodel-th}. In the second row, we have pulled $ F_q^1 $ into the sum and used that there are no ID, DS and DW disks with inputs $ h_1, …, h_N $ according to \autoref{th:subdisk-main-onlyCR}. In the third row, we have turned the sum into a double sum ranging over the possible output basis elements $ h $. The notation $ h: L_1 → L_{N+1} $ used is a slight abuse: The sum is intendend to run over the basis elements $ h ∈ \Hom_{\H\ZigzagCat} (L_1, L_{N+1}) $. In the fourth row, we have re-enumerated the summands as smooth immersed disks instead of CR disks. This enumeration uses the bijection $ Φ $ set up in \autoref{th:subdisk-main-onlyCR}. In the fifth row, we have used that $ Φ^{-1} (D) $ and $ D $ have the same Abouzaid sign and $ q $-parameter, according to \autoref{th:subdisk-main-onlyCR}. In the sixth row, we have inserted the definition of the products of $ \relpreFuk Q \restr{\Ob\ZigzagCat} $, according to \autoref{th:prelim-fukaya-zigzag-restr}.

This proves step (4) and finishes the construction of the strict isomorphism $ F_q $ between $ (\H\DefZigzagCat)_{\transversal} $ and $ \relpreFuk Q \restr{\Ob\ZigzagCat} $. By construction, $ F_q $ sends $ L $ to $ \smooth L $ and a cohomology basis element $ h: L_1 → L_2 $ to the associated intersection point $ p: \smooth L_1 → \smooth L_2 $. This finishes the proof.
\end{proof}

\begin{remark}
The main result shows that $ \HTw\Gtl_q Q $ is a candidate for a relative wrapped Fukaya category. Also, the subcategory of $ \HTw\Gtl_q Q $ given by $ ℤ/2ℤ $-graded band objects is a candidate model for $ \relFuk Q $ in the sense of \autoref{def:prelim-fukaya-rel-def}.

It seems likely that $ \H\DefZigzagCat $ is (gauge equivalent to) the subcategory of zigzag paths in (any model for) $ \relFuk Q $. Our main result is however no guarantee for this, since taking subcategories and lifting pre-categories to categories need not commute: Every subcategory of a lift is a lift of the subcategory, but not the other way around. For further discussion we refer to \autoref{sec:whyshould-candidate}.
\end{remark}

\appendix

\section{Examples}
\label{sec:2examples}
We provide here an example of a CR and an ID disk together with their matching result components. The aim is to demonstrate in practice how one finds the preimage of a given CR or ID disk under the subdisk mapping. The examples illustrate the strong geometric aspect of the subdisk construction and its inverse construction. On the other hand, the examples demonstrate the sheer amount of case distinctions and precision work required for reconstructing the result component from a given CR or ID disk.

\subsection{ID disk}
We present here an example pair of an ID disk and its matching result component. We depart from the view of the ID disk and construct by inspection the corresponding Kadeishvili tree and result component. The ID disk we present is very small. The reader can use this example to get a feeling how the presence of the degenerate input and the small size of the disk are translated into the Kadeishvili tree.

The example ID disk is presented in \autoref{fig:examples-ID-hexagon}. It is situated at a puncture with six incident arcs and six incident polygons. Every pair of two neighboring incident arcs makes for a zigzag path, and the smoothed zigzag curves $ \smooth L_1, …, \smooth L_6 $ have six intersections around the puncture. These intersections alone bound a hexagon.

Recall that every zigzag curve is supposed to have locations assigned of identity and co-identity morphisms. In our example, the identity morphism of the sixth curve is supposed to lie on the curve's second arc when reading clockwise around the puncture. The disk we present makes use of this feature in order to be a disk of ID type with six inputs and an identity output.

Let us explain precisely the data of this disk: The disk's inputs are the six intersection points $ h_1, …, h_6 $ of $ \smooth L_1, …, \smooth L_6 $ around the puncture, and the output is the identity of $ \smooth L_6 $. The identity output lies infinitesimally close to the first input. More precisely, the first input is the degenerate input of this ID disk and succeeds the output mark. The arcs are oriented so that $ \smooth L_1 $ is oriented clockwise and $ \smooth L_6 $ is counterclockwise with the disk. This turns the disk into an ID disk according to definition.

\begingroup
\newcommand{\hexagondraw}{\path[draw] (0, 0) -- ++(up:1) coordinate[midway] (A) coordinate[pos=0.7] (Ai) -- ++(330:1) coordinate[midway] (B) coordinate[pos=0.3] (Bi) ++(60:0.2) -- ++(150:1) coordinate[midway] (C) coordinate[pos=0.7] (Ci) -- ++(30:1) coordinate[midway] (D) coordinate[pos=0.3] (Di) ++(120:0.2) -- ++(210:1) coordinate[midway] (E) coordinate[pos=0.7] (Ei) -- ++(up:1) coordinate[midway] (F) coordinate[pos=0.3] (Fi) ++(left:0.2) -- ++(down:1) coordinate[midway] (G) coordinate[pos=0.7] (Gi) -- ++(150:1) coordinate[midway] (H) coordinate[pos=0.3] (Hi) ++(240:0.2) -- ++(330:1) coordinate[midway] (I) coordinate[pos=0.7] (Ii) -- ++(210:1) coordinate[midway] (J) coordinate[pos=0.3] (Ji) ++(300:0.2) -- ++(30:1) coordinate[midway] (K) coordinate[pos=0.7] (Ki) -- ++(down:1) coordinate[midway] (L) coordinate[pos=0.3] (Li);
\path ($ (F)!0.5!(G) $) coordinate (1);
\path ($ (D)!0.5!(E) $) coordinate (2);
\path ($ (B)!0.5!(C) $) coordinate (3);
\path ($ (L)!0.5!(A) $) coordinate (4);
\path ($ (J)!0.5!(K) $) coordinate (5);
\path ($ (H)!0.5!(I) $) coordinate (6);}

\begin{figure}
\centering
\begin{tikzpicture}[scale=1.5]
\hexagondraw
\path (1) -- (6) coordinate[pos=0.4] (out);
\path (out) node[below] {\tiny out};
\foreach \i in {1, 2, 3, 4, 5, 6, out} \path[fill] (\i) circle[radius=0.03];
\path[draw, semithick, -{To[scale=1.5]}] (1) -- (2);
\path[draw, semithick] (2) -- (3);
\path[draw, semithick] (3) -- (4);
\path[draw, semithick] (4) -- (5);
\path[draw, semithick] (5) -- (6);
\path[draw, semithick, {To[scale=1.5]}-] (6) -- (out) to[bend right] (1);
\path (1) -- (2) node[midway, above] {\small $ \smooth L_1 $};
\path (2) -- (3) node[midway, right] {\small $ \smooth L_2 $};
\path (3) -- (4) node[midway, below] {\small $ \smooth L_3 $};
\path (4) -- (5) node[midway, below] {\small $ \smooth L_4 $};
\path (5) -- (6) node[midway, left] {\small $ \smooth L_5 $};
\path (6) -- (out) node[midway, above] {\small $ \smooth L_6 $};
\end{tikzpicture}
\caption{The immersed disk: a hexagon}
\label{fig:examples-ID-hexagon}
\end{figure}
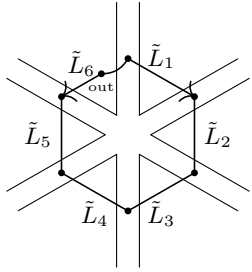

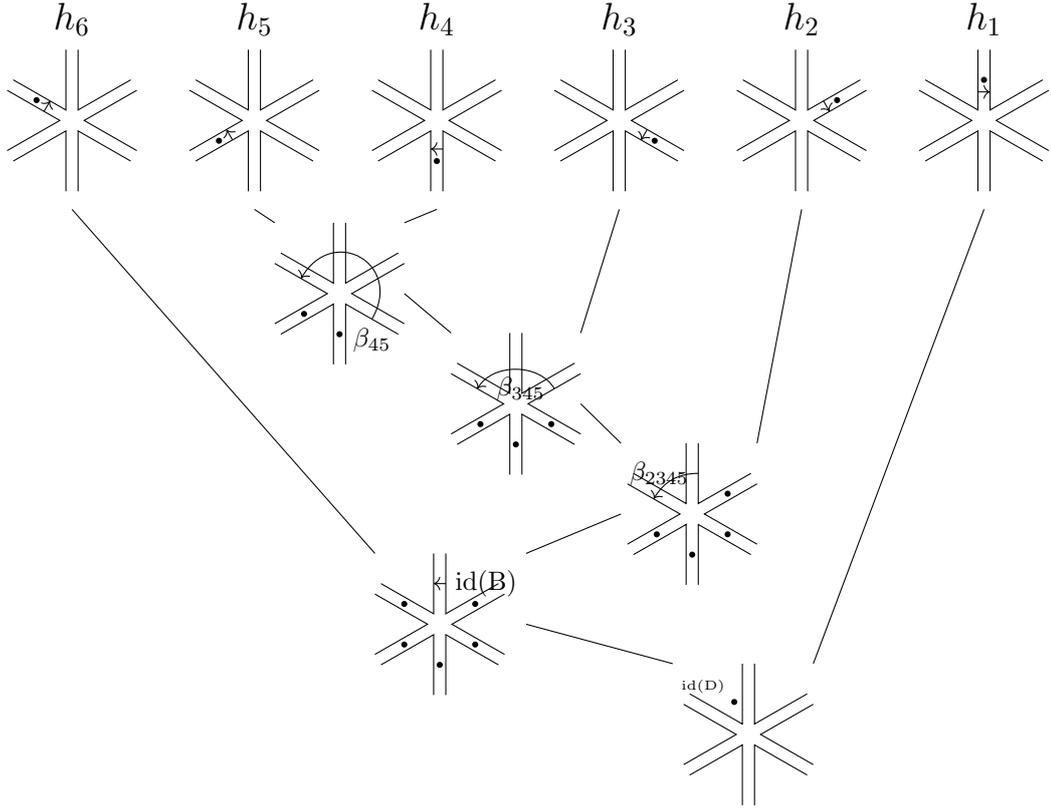
\begin{figure}
\centering
\begin{tikzpicture}[scale=0.8]
\begin{scope}[shift={(0, 8)}, local bounding box=i6]
\hexagondraw
\foreach \i in {6} \path[fill] (\i) circle[radius=0.05];
\path[draw, ->] (Ii) to (Hi);
\path (1) ++(up:1) node {\Large $ h_6 $};
\end{scope}
%
\begin{scope}[shift={(3, 8)}, local bounding box=i5]
\hexagondraw
\foreach \i in {5} \path[fill] (\i) circle[radius=0.05];
\path[draw, ->] (Ki) to (Ji);
\path (1) ++(up:1) node {\Large $ h_5 $};
\end{scope}
%
\begin{scope}[shift={(6, 8)}, local bounding box=i4]
\hexagondraw
\foreach \i in {4} \path[fill] (\i) circle[radius=0.05];
\path[draw, ->] (Ai) to (Li);
\path (1) ++(up:1) node {\Large $ h_4 $};
\end{scope}
%
\begin{scope}[shift={(9, 8)}, local bounding box=i3]
\hexagondraw
\foreach \i in {3} \path[fill] (\i) circle[radius=0.05];
\path[draw, ->] (Ci) to (Bi);
\path (1) ++(up:1) node {\Large $ h_3 $};
\end{scope}
%
\begin{scope}[shift={(12, 8)}, local bounding box=i2]
\hexagondraw
\foreach \i in {2} \path[fill] (\i) circle[radius=0.05];
\path[draw, ->] (Ei) to (Di);
\path (1) ++(up:1) node {\Large $ h_2 $};
\end{scope}
%
\begin{scope}[shift={(15, 8)}, local bounding box=i1]
\hexagondraw
\foreach \i in {1} \path[fill] (\i) circle[radius=0.05];
\path[draw, ->] (Gi) to (Fi);
\path (1) ++(up:1) node {\Large $ h_1 $};
\end{scope}
%
\begin{scope}[shift={($ (i4)!0.5!(i5) + (0, -4.5) $)}, local bounding box=b45]
\hexagondraw
\path (1) -- (6) coordinate[pos=0.4] (out);
\foreach \i in {4, 5} \path[fill] (\i) circle[radius=0.05];
\path[draw, ->, bend right=90, looseness=2] (B) to node[at start, below] {$ β_{45} $} (I);
\end{scope}
%
\begin{scope}[shift={($ (b45) + (3, -3) $)}, local bounding box=b345]
\hexagondraw
\path (1) -- (6) coordinate[pos=0.4] (out);
\foreach \i in {3, 4, 5} \path[fill] (\i) circle[radius=0.05];
\path[draw, ->, bend right=60] (D) to node[at start, left] {$ β_{345} $} (I);
\end{scope}
%
\begin{scope}[shift={($ (b345) + (3, -3) $)}, local bounding box=b2345]
\hexagondraw
\path (1) -- (6) coordinate[pos=0.4] (out);
\foreach \i in {2, 3, 4, 5} \path[fill] (\i) circle[radius=0.05];
\path[draw, ->, bend right=30] (F) to node[at start, left] {$ β_{2345} $} (I);
\end{scope}
%
\begin{scope}[shift={($ (b2345) + (-4, -3) $)}, local bounding box=id23456]
\hexagondraw
\foreach \i in {2, 3, 4, 5, 6} \path[fill] (\i) circle[radius=0.05];
\path[draw, ->] (F) -- (G) node[at start, right] {$ \id $(B)};
\end{scope}
%
\begin{scope}[shift={($ (id23456) + (5, -3) $)}, local bounding box=idD]
\hexagondraw
\path (1) -- (6) coordinate[pos=0.4] (out);
\foreach \i in {out} \path[fill] (\i) circle[radius=0.05];
\path (out) node[above left] {\tiny $ \id $(D)};
\end{scope}
\path[draw] (i5.south) ++(down:0.3) -- (b45.north west);
\path[draw] (i4.south) ++(down:0.3) -- (b45.north east);
\path[draw] (b45.east) -- (b345.north west);
\path[draw] (i3.south) ++(down:0.3) -- (b345.north east);
\path[draw] (b345.east) -- (b2345.north west);
\path[draw] (i2.south) ++(down:0.3) -- (b2345.north east);
\path[draw] (b2345.west) -- (id23456.north east);
\path[draw] (i6.south) ++(down:0.3) -- (id23456.north west);
\path[draw] (id23456.east) -- (idD.north west);
\path[draw] (i1.south) ++(down:0.3) -- (idD.north east);
\end{tikzpicture}
\caption{The Kadeishvili tree, with result components depicted at each node}
\label{fig:examples-ID-rescomp}
\end{figure}
\endgroup

To interpret the meaning of this disk, we have to understand its features: First, the disk is very small. The difficulty for us lies instead in the fact that the disk has no direct discrete analog: There is no space between the arcs to form a higher product of $ \Gtl_q Q $. The disk is however an ID disk and as such has a nonzero value in our explicit model $ \H\DefZigzagCat $. Second, the disk is unpredictable in the relative Fukaya category. Indeed, the sequence $ \smooth L_1, …, \smooth L_6, \smooth L_1 $ is not transversal. The disk even has an input lying infinitesimally close to the output. All this means the higher product in the relative Fukaya is unpredictable. The smallness and unpredictability features are nicely illustrated by the fact that the simplest possible Kadeishvili tree $ φπ_q μ(h_6, …, h_1) $ simply gives a zero result.

There exists is a single Kadeishvili π-tree that yields a nonzero output for the input sequence $ h_1, …, h_6 $. In fact, there is only a single result component, depicted in \autoref{fig:examples-ID-rescomp}. We explain here the precise data of this result component and why it is the only result component.

\begin{description}
\item[Step 1: Combining $ h_5 $ and $ h_4 $:] The deformed cohomology basis element $ h_5 $ is the sum of an arc identity and a $ β' $ (C) morphism. The deformed cohomology basis element $ h_4 $ is the sum of an arc identity and a $ β $ (C) morphism. Combining both into the product $ μ^2_q (h_5, h_4) $ gives the sum of two terms, the first of which comes from the arc identity of $ h_5 $ and the $ β $ (C) morphism of $ h_4 $, the second comes from the $ β' $ (C) morphism of $ h_5 $ and the arc identity of $ h_4 $. The product morphisms are both situation A morphisms, the first is of type $ γβ $ (A) and the second of type $ βα $ (A). Applying the codifferential $ h_q $ to the first gives zero, since $ γβ $ (A) lies in the $ R $-part of $ \Hom_{\DefZigzagCat} (L_3, L_5) $. Applying of the codifferential $ h_q $ to the second product gives a certain angle $ β_{45} $, depicted explicitly in the figure. In summary, only the second term survives and is used for the result component.
\item[Step 2: Combining with $ h_3 $:] The deformed cohomology basis element $ h_3 $ consists of an arc identity and a $ β' $ (C) morphism. Combining with $ β_{45} $ gives the product $ μ^2_q (β_{45}, \id \text{(C)}) $ which is a situation A morphism of type $ βα $ (A). Applying the codifferential $ h_q $ gives a certain angle $ β_{345} ∈ \Hom_{\DefZigzagCat} (L_2, L_5) $, depicted explicitly in the figure.
\item[Step 3: Combining with $ h_2 $:] Analogous to the previous step, gives angle $ β_{2345} ∈ \Hom_{\DefZigzagCat} (L_1, L_5) $.
\item[Step 4: Combining with $ h_6 $:] The deformed cohomology basis element $ h_6 $ consists of an arc identity and a $ β $ (C) morphism. Combining with $ β_{2345} $ gives the product $ μ^2_q (\id \text{(C)}, β_{2345}) $ which is a situation B morphism of type $ α_2 $ (B). Application of the codifferential $ h_q $ gives the identity angle $ \id \text{(B)} ∈ \Hom_{\DefZigzagCat} (L_1, L_6) $.
\item[Step 5: Combining with $ h_1 $:] The deformed cohomology basis element $ h_1 $ consists of an arc identity $ \id $ (C) and a $ β' $ (C) morphism. Combining with $ \id $ (B) gives the product $ μ^2_q (\id\text{(B)}, \id\text{(C)}) $ which is the arc identity on the second arc of $ L_6 $ read in clockwise direction. By assumption, this second arc is the identity location of $ L_6 $. Application of the projection $ π_q $ gives the identity morphism $ \id_{L_6} ∈ \Hom_{\DefZigzagCat} (L_6, L_6) $. This is the final result component of the Kadeishvili π-tree.
\item[$ q $-parameters:] In the past five steps, we have ignored signs and $ q $-parameters in the result component. The $ q $-parameter in fact consists of the single puncture located at the center of the hexagon. This parameter enters the result component in Step 1 by means of the $ β' $ (C) morphism of $ h_5 $.
\end{description}

The result component described above and depicted in \autoref{fig:examples-ID-rescomp} is the only result component of the sequence $ h_1, …, h_6 $, at least as far as displayed in the figure. For the specific π-tree, the single result component is the only result component. Indeed we have exhausted in every step all possible products $ μ^2_q $ or $ μ^{≥3}_q $ and all terms in their codifferential, apart from possible tail terms which lie far away and are not visible in the figure. Other π-trees with the same input sequence $ h_1, …, h_6 $ do not yield any result components. Indeed, one might for instance try to combine $ β_{2345} $ with $ h_1 $ before combining with $ h_6 $. However, the product $ μ^2_q (β_{2345}, \id\text{(C)}) $ gives a situation B morphism $ α_1 $ (B) whose codifferential vanishes. This explains how the result component depicted in \autoref{fig:examples-ID-rescomp} is the only result component and illustrates the delicate nature of matching disks with result components.

\input{examples/cr.tex}

\section{Uncurving of band objects}
\label{sec:trick}
In this section, we prove \autoref{th:uncurving-trick-works}, which concerns uncurvability of band objects. It is our task to explicitly construct a deformed twisted complex such that a given band object becomes curvature-free. More precisely, the starting point consists of:
\begin{itemize}
\item a punctured surface $ (S, M) $,
\item a full arc system $ \cA $ with [NMDC],
\item a standard deformation $ \Gtl_r \cA $ where $ r ∈ \mathfrak{m} Z(\Gtl \cA) $ is without 1-component,
\item a band object organized in a twisted complex $ X = (⊕ a_i [s_i], δ) $,
\end{itemize}
such that the underlying curve of $ X $ in the closed surface $ S $ is not contractible and does not bound a teardrop, and every connecting angle $ α_i $ in $ δ $ is longer than an identity and shorter than a full turn. It is our aim to construct a deformation $ δ_q ∈ \Hom^1_{\Add\Gtl_r \cA} (X, X) $ such that $ X_q = (⊕ a_i [s_i], δ_q) $ has vanishing curvature.

The present section is logically independent of the computation of the minimal model $ \H\DefZigzagCat $. However, it builds directly on its methods in three ways: First, we introduce a notion of situations to characterize types of angles between arcs of $ X $, similar to the notion of situations for morphisms between zigzag paths. Second, we build a protocol which describes what kind of angles we may encounter while gathering the complementary angles, similar to the E, F, G, H disks or the subdisk protocol for zigzag paths. Third, we build up the list of complementary angles in a recursive way, a bit similar to the way we work with tails of morphims or subdisks for zigzag paths.

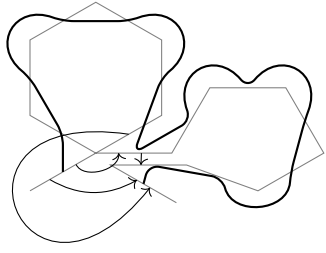
\begin{figure}
\centering
\begin{tikzpicture}
\path[use as bounding box] (-0.5, -1) -- (4, 3);
\path[draw, gray] (0, 0) \foreach \i\j in {1/30, 2/150, 3/90, 4/30, 5/330, 6/270, 6A/210, 7/0, 8/60, 9/0, 10/300, 11/210, 12/160, 13/180, 14/330} { -- ++(\j:1) coordinate[pos=0.3] (\i-start) coordinate[midway] (\i-mid) coordinate[pos=0.7] (\i-end)};
\path[draw, thick, rounded corners] (1-mid) to (2-mid) to (3-mid) to[bend left=80, looseness=2] (4-mid) to (5-mid) to[bend left=80, looseness=2] (6-mid) to (7-mid) to (8-mid) to[bend left=80, looseness=2] (9-mid) to[bend left=80, looseness=2] (10-mid) to (11-mid) to[bend left=80, looseness=1.5] (12-mid) to (13-mid) to (14-mid);
\path[draw, ->, bend right] (1-start) to ($ (14-mid)!0.5!(14-start) $);
\path[draw, ->, bend right=80] ($ (1-end) $) to (7-start);
\path[draw, ->, bend right=120, looseness=9] (6A-mid) to ($ (14-mid)!0.5!(14-end) $);
\path[draw, ->] ($ (7-mid)!0.5!(7-end) $) to ($ (13-mid)!0.5!(13-end) $);
\end{tikzpicture}
\caption{An example which requires four additional complementary angles}
\label{fig:trick-intro-example}
\end{figure}

The idea to construct $ δ_q $ is to add not only complementary angles for all connecting angles $ α_i $ of $ X $, but also complementary angles at locations where the strands of $ X $ come close to each other. An example for $ X $ is depicted in \autoref{fig:trick-intro-example}. Apart from the complementary angles for the connecting angles $ α_i $, this example requires four additional angles to be added to $ δ_q $. In terms of the situational formalism of \autoref{sec:trick-situations}, the four angles are of type $ \id $ (A'), $ β $ (A'), $ β $ (A) and $ β (A) $.

In \autoref{sec:trick-situations}, we inspect the possible configurations of arcs and angles around $ X $ and introduce notions of situations. In \autoref{sec:trick-protocol}, we define a temporary type of possible angles we may add in order to produce $ δ_q $, and investigate their possible products. In \autoref{sec:trick-flower}, we construct $ δ_q $ in a recursive way and show that the curvature of $ X_q = (⊕ a_i [s_i], δ_q) $ vanishes.

\subsection{Situations}
\label{sec:trick-situations}
In this section, we examine the possible configurations of arcs and angles along $ X $. We capture these types configurations in the notion of A, A', A'', ID and D situations. The main difference with the notion of situations for zigzag paths is that two different indexed arcs which are equal as arcs of $ \cA $ need not determine an intersection of the underlying curve. For this reason, we obtain a slightly larger amount of different situations.


\begin{definition}
Regard the band object $ X = (⊕ a_i [s_i], δ) $. An \emph{indexed arc} on $ X $ is a choice of arc $ a_i $, remembering the index $ i $. If $ a_i $ is an indexed arc of $ X $, then the \emph{strand} of $ X $ at $ a_i $ refers to the portion of $ X $ given by the neighboring indexed arcs $ …, a_{i-1}, a_i, a_{i+1}, … $. An \emph{elementary morphism} $ ε: a_i → a_j $ on $ X $ is a single angle between two indexed arcs $ a_i $, $ a_j $ on $ X $, interpreted as $ ε ∈ \Hom_{\Add\Gtl \cA} (X, X) $. The \emph{source strand} or \emph{target strand} of an elementary morphism $ ε: a_i → a_j $ is the strand of $ X $ at $ a_i $ or $ a_j $, respectively.

Two angles in $ \cA $ which wind around a common puncture have an \emph{overlap} if they contain a shared indecomposable angle. Otherwise, the two angles are \emph{disjoint}. For every indexed arc $ a_i $, the object $ X $ contains a distinction whether $ X $ \emph{turns left} or \emph{right} towards a given endpoint of $ a_i $.
\end{definition}

Recall that an \emph{arc incidence} at a puncture is slightly different from an arc incident at a puncture: An arc incidence includes the datum whether it concerns the head or tail of the arc. A loop incident at a puncture gives rise to two arc incidences. We generically denote a full turn by $ ℓ $. We may call the $ α_i $ angles of $ X $ also the \emph{turning angles}.

\begin{figure}
\centering
\begin{subfigure}{0.2\linewidth}
\centering
\begin{tikzpicture}
\path[draw] (0, 0) node[below] {…} -- ++(120:1) coordinate[pos=0.7] (1-end) -- ++(60:1) coordinate[pos=0.3] (1-start) coordinate[pos=0.7] (2-end) -- ++(20:1) coordinate[pos=0.3] (2-start) coordinate[pos=0.7] (3-start) -- ++(135:1) coordinate[pos=0.3] (3-end) coordinate[pos=0.7] (4-start)  -- ++(20:1) coordinate[pos=0.3] (4-end) coordinate[pos=0.7] (5-end) -- ++(315:1) coordinate[pos=0.3] (5-start) coordinate[pos=0.7] (6-end)  -- ++(60:1) coordinate[pos=0.3] (6-start) node[right] {…};
\path[draw, thick, bend right=100, looseness=2] (1-start) to (1-end);
\path[draw, thick, bend right=100, looseness=2] (2-start) to (2-end);
\path[draw, thick, bend right=150, looseness=8] (3-start) to (3-end);
\path[draw, thick, bend right=20] (4-start) to (4-end);
\path[draw, thick, bend right=100, looseness=2] (5-start) to (5-end);
\path[draw, thick, bend right=30] (6-start) to (6-end);
\end{tikzpicture}
\caption{We draw the angles $ α_i $ thickly.}
\end{subfigure}
\begin{subfigure}{0.2\linewidth}
\centering
\begin{tikzpicture}
\path[draw] (0, 0) -- ++(30:1.5) coordinate[midway] (1-end) coordinate[pos=0.3] (a-end) coordinate (m) -- ++(150:1.5) coordinate[midway] (1-start) coordinate[pos=0.7] (a-start);
\path[draw] (m)++(right:0.2) -- ++(30:1.5) coordinate[midway] (2-end) coordinate[pos=0.7] (g-end) (m)++(right:0.2) -- ++(330:1.5) coordinate[midway] (2-start) coordinate[pos=0.7] (g-start);
\path[draw, thick, bend right] (1-start) to (1-end);
\path[draw, thick, bend right] (2-start) to (2-end);
\path[draw, ->, bend right] (a-start) to node[midway, left] {$ α $} (a-end);
\path[draw, ->, bend right] (g-start) to node[midway, right] {$ γ $} (g-end);
\path[draw, ->, bend right=60] (a-end) to node[midway, below] {$ β $} (g-start);
\end{tikzpicture}
\caption{A situation}
\label{fig:trick-situation-A}
\end{subfigure}
\begin{subfigure}{0.2\linewidth}
\centering
\begin{tikzpicture}
\path[draw] (0, 0) -- ++(30:1.5) coordinate[midway] (1-end) coordinate[pos=0.3] (a-end) coordinate (m) -- ++(90:1.5) coordinate[midway] (1-start) coordinate[pos=0.7] (a-start);
\path[draw] (m)++(right:0.2) -- ++(90:1.5) coordinate[midway] (2-end) coordinate[pos=0.7] (g-end) (m)++(right:0.2) -- ++(330:1.5) coordinate[midway] (2-start) coordinate[pos=0.7] (g-start);
\path[draw, thick, bend right] (1-start) to (1-end);
\path[draw, thick, bend right] (2-start) to (2-end);
\path[draw, ->, bend right] (a-start) to node[midway, left] {$ α $} (a-end);
\path[draw, ->, bend right] (g-start) to node[midway, right] {$ γ $} (g-end);
\path[draw, ->, bend right=60] (a-end) to node[midway, below] {$ β $} (g-start);
\end{tikzpicture}
\caption{A' situation}
\label{fig:trick-situation-A'}
\end{subfigure}
\begin{subfigure}{0.2\linewidth}
\centering
\begin{tikzpicture}
\path[draw] (0, 0) -- ++(330:1.5) coordinate[midway] (1-end) coordinate[pos=0.3] (a-end) coordinate (m) -- ++(down:1.5) coordinate[midway] (1-start) coordinate[pos=0.7] (a-start);
\path[draw] (m)++(right:0.2) -- ++(30:1.5) coordinate[midway] (2-end) coordinate[pos=0.7] (g-end) (m)++(right:0.2) -- ++(down:1.5) coordinate[midway] (2-start) coordinate[pos=0.7] (g-start);
\path[draw, thick, bend right] (1-end) to (1-start);
\path[draw, thick, bend right] (2-start) to (2-end);
\path[draw, ->, bend right] (a-end) to node[midway, left] {$ α $} (a-start);
\path[draw, ->, bend right] (g-start) to node[midway, right] {$ γ $} (g-end);
\path[draw, ->, bend right=60] (g-end) to node[midway, below] {$ β' $} (a-end);
\end{tikzpicture}
\caption{A'' situation}
\label{fig:trick-situation-A''}
\end{subfigure}
\begin{subfigure}{0.2\linewidth}
\centering
\begin{tikzpicture}
\path[draw] (0, 0) -- ++(up:1) coordinate[pos=0.8] (m1) node[pos=0.2, left] {$ a_i $} coordinate[midway] (1-start) coordinate[pos=0.6] (2-end);
\path[draw] (0.2, 0) -- ++(up:1) coordinate[pos=0.8] (m2) node[pos=0.2, right] {$ a_j $} coordinate[midway] (1-end) coordinate[pos=0.6] (3-start);
\path[draw, ->] (1-start) to (1-end);
\path[fill] (0.1, 1.1) circle[radius=0.05];
\path[draw] (0, 1.2) -- ++(up:1) coordinate[pos=0.2] (m3) coordinate[pos=0.4] (2-start);
\path[draw] (0.2, 1.2) -- ++(up:1) coordinate[pos=0.2] (m4) coordinate[pos=0.4] (3-end);
\path[draw, thick, bend right=90, looseness=1.5] (m3) to (m1);
\path[draw, thick, bend right=90, looseness=1.5] (m2) to (m4);
\path[draw, ->, bend right=90, looseness=1.5] (2-start) to node[midway, left] {$ α $}(2-end);
\path[draw, ->, bend right=90, looseness=1.5] (3-start) to node[midway, right] {$ γ $} (3-end);
\end{tikzpicture}
\caption{ID situation}
\label{fig:trick-situation-ID}
\end{subfigure}
\begin{subfigure}{0.2\linewidth}
\centering
\begin{tikzpicture}
\path[draw] (0, 0) -- ++(60:1.5) coordinate[pos=0.6] (1-start) coordinate[pos=0.3] (a-start) -- ++(300:1.5) coordinate[pos=0.4] (1-end) coordinate[pos=0.7] (a-end);
\path[draw, thick, bend right] (1-start) to (1-end);
\path[draw, ->, bend right] (a-start) to node[midway, below] {$ α $} (a-end);
\path[draw, ->] (a-end) to[out=30, in=0] ($ (a-end)!0.5!(a-start) + (up:1.5) $) node[above] {$ α' $} to[out=180, in=150] (a-start);
\end{tikzpicture}
\caption{D situation}
\label{fig:trick-situation-D}
\end{subfigure}
\caption{Illustration of A, A', A'', ID and D situations}
\end{figure}
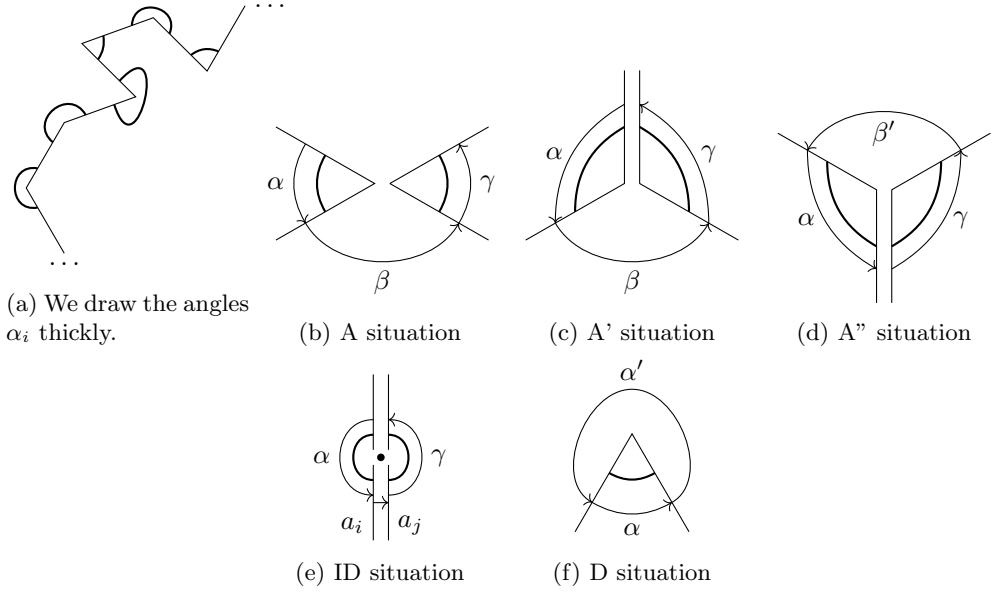

\begin{definition}
We define the following types of situations on $ X $:
\begin{itemize}
\item A \emph{type A situation} consists of a puncture together with incidences of two indexed arc $ a_i $ and $ a_j $ such that (a) $ X $ turns left at $ a_i $ towards the puncture, and (b) $ X $ turns right at $ a_j $ towards the puncture, and (c) the turning angles do not overlap. The associated angles of the situation are denoted $ α $, $ β $, $ γ $ and $ β' $, as in \autoref{fig:trick-situation-A}.
\item A \emph{type A' situation} consists of a puncture $ q $ together with incidences of two distinct indexed arcs $ a_i $ and $ a_j $ which are equal as arc of $ \cA $, such that the strand of $ a_i $ turns right at $ a_i $ towards $ q $ and the strand of $ a_j $ turns left towards $ q $. The associated angles are denoted $ α $, $ β $, $ γ $ and $ \id $, as in \autoref{fig:trick-situation-A'}.
\item A \emph{type A'' situation} consists of a puncture $ q $ together with incidences of two distinct indexed arcs $ a_i $ and $ a_j $ which are equal as arcs of $ \cA $, such that the strand of $ a_i $ turns left at $ a_i $ towards $ q $ and the strand of $ a_j $ turns right towards $ q $. The associated angles are denoted $ α $, $ \id $, $ γ $, $ β' $, as in \autoref{fig:trick-situation-A''}.
\item A \emph{type ID situation} consists of a puncture $ q $ together with incidences of two distinct indexed arcs $ a_i $ and $ a_j $ which are equal as arcs of $ \cA $, such that (a) the strand of $ a_i $ turns left at $ a_i $ towards $ q $, and (b) the strand of $ a_j $ turns right towards $ q $, and (c) the turning angles together form a full turn. The associated angles are denoted $ α $, $ \id $ and $ γ $, as in \autoref{fig:trick-situation-ID}.
\item A \emph{type D situation} consists of a pair of consecutive arcs of $ X $. It gives rise to $ α $ and $ α' $ angles: A type $ α $ angle of $ X $ is one of the $ α_i $ angles of $ X $. A type $ α' $ angle is one of the $ α_i' $ angles of $ X $, see \autoref{fig:trick-situation-D}.
\end{itemize}
\end{definition}

\begin{remark}
When an angle comes from a situation of a certain type, we typically indicate the type in brackets for clarity. For instance, a certain angle may qualify as a “$ β $ (A) angle”. When working with situations, we may from time to time indicate the situations by its associated angles, for instance referring to a “type A situation $ (α_R, β_R, γ_r, β'_R) $”.
\end{remark}

The difference with zigzag paths is that the arcs of $ X $ are not oriented in the same direction, and the $ α $ (D) angles are not necessarily indecomposable. Another important difference is that a zigzag segment is never contractible, while a segment of $ X $ may be contractible. In our uncurving construction, we form $ δ_q $ by inserting angles whenever a contractible segment lies “above” the angle. In order to make this precise, we use the following terminology:

\begin{definition}
\label{def:trick-situations-above}
For an A/A'/A''/ID situation, \emph{above} refers to tracing $ X $ in the opposite direction of $ α $ and in the direction of $ γ $.
\end{definition}

The direction “above” in the situation figures \autoref{fig:trick-situation-A} till \ref{fig:trick-situation-ID} is the natural upwards direction on paper.

\subsection{Uncurving protocol}
\label{sec:trick-protocol}
In this section, we examine possible products which can be made from angles of A, A', A'', ID and D situations. The core tool is a notion of angles with balloons. The datum of balloons makes it possible to safely examine possible products $ μ^2 $ and $ μ^{≥3}_r $ between angles. All of the angles we later insert into $ δ_q $ are of this type, but come with additional data which is irrelevant and not accessible at the current stage of the construction. This way, angles with balloons serve as a “protocol” which greatly facilitates the construction of $ δ_q $.

The aim of the examination is to draw maximally strong conclusions on the configurations of arcs and angles from the fact that the underlying curve of $ X $ is not contractible and does not bound a teardrop. The notion of angles with balloons is the cheapest way to incorporate this property of $ X $ into individual angles. Angles with balloons are elementary angles whose $ X $ strand lying above the angle is contractible:


\begin{definition}
An \emph{angle with balloon} on $ X $ is a $ βℓ^m $ (A/A') or $ \id ℓ^m $ (ID/A'') angle whose $ X $ segment above $ β $ or $ \id $ is contractible. We also count $ α $ (D) and $ α'ℓ^m $ (D) as angles with balloons.
\end{definition}

\begin{figure}
\centering
\begin{subfigure}[b]{0.24\linewidth}
\centering
\begin{tikzpicture}
\path[draw] (0, 0) -- ++(30:1.5) coordinate[midway] (1-end) coordinate[pos=0.3] (a-end) coordinate (m) -- ++(150:1.5) coordinate (stop-left) coordinate[midway] (1-start) coordinate[pos=0.7] (a-start);
\path[draw] (m)++(right:0.2) -- ++(30:1.5) coordinate[midway] (2-end) coordinate[pos=0.7] (g-end)coordinate (stop-right) (m)++(right:0.2) -- ++(330:1.5) coordinate[midway] (2-start) coordinate[pos=0.7] (g-start);
\path[draw, thick, bend right] (1-start) to (1-end);
\path[draw, thick, bend right] (2-start) to (2-end);
\path[draw, ->, bend right=60] (a-end) to node[midway, below] {$ β $} (g-start);
\path[draw, dashed] (stop-right) to[out=70, in=0] ($ (stop-left)!0.5!(stop-right) + (up:2) $) to[out=180, in=110] (stop-left);
\path (stop-left) -- (stop-right) node[midway, above, align=center] {contractible \\ segment};
\end{tikzpicture}
\caption{$ β $ (A) with balloon}
\end{subfigure}
\begin{subfigure}[b]{0.24\linewidth}
\centering
\begin{tikzpicture}
\path[draw] (0, 0) -- ++(30:1.5) coordinate[midway] (1-end) coordinate[pos=0.3] (a-end) coordinate (m) -- ++(90:1.5) coordinate[midway] (1-start) coordinate[pos=0.7] (a-start) coordinate (A);
\path[draw] (m)++(right:0.2) -- ++(90:1.5) coordinate[midway] (2-end) coordinate[pos=0.7] (g-end) coordinate(B) (m)++(right:0.2) -- ++(330:1.5) coordinate[midway] (2-start) coordinate[pos=0.7] (g-start);
\path[draw, thick, bend right] (1-start) to (1-end);
\path[draw, thick, bend right] (2-start) to (2-end);
\path[draw, ->, bend right=60] (a-end) to node[midway, below] {$ β $} (g-start);
\path[draw, dashed] (B) to[out=30, in=270] ++(30:1.5) to[out=90, in=0] ($ (A)!0.5!(B) + (up:1.2) $) node[below, align=center] {contractible \\ segment} to[out=180, in=90] ($ (A)+(150:1.5) $) to[out=270, in=150] (A);
\end{tikzpicture}
\caption{$ β $ (A') with balloon}
\end{subfigure}
\begin{subfigure}[b]{0.24\linewidth}
\centering
\begin{tikzpicture}
\path[draw] (0, 0) -- ++(330:1.5) coordinate[midway] (1-end) coordinate[pos=0.3] (a-end) coordinate (m) -- ++(down:1.5) coordinate[midway] (1-start) coordinate[pos=0.7] (a-start);
\path[draw] (m)++(right:0.2) -- ++(30:1.5) coordinate (B) coordinate[midway] (2-end) coordinate[pos=0.7] (g-end) (m)++(right:0.2) -- ++(down:1.5) coordinate[midway] (2-start) coordinate[pos=0.7] (g-start);
\path[draw, thick, bend right] (1-end) to (1-start);
\path[draw, thick, bend right] (2-start) to (2-end);
\path[draw, ->] (a-start) to (g-start);
\path[draw, bend right=120, looseness=2, dashed] (B) to (0, 0);
\path (B) -- (0, 0) node[midway, above, align=center] {contractible \\ segment};
\end{tikzpicture}
\caption{$ \id $ (A'') with balloon}
\end{subfigure}
\begin{subfigure}[b]{0.24\linewidth}
\centering
\begin{tikzpicture}
\path[draw] (0, 0) -- ++(up:1) coordinate[pos=0.8] (m1) node[pos=0.2, left] {$ a_i $} coordinate[midway] (1-start) coordinate[pos=0.6] (2-end);
\path[draw] (0.2, 0) -- ++(up:1) coordinate[pos=0.8] (m2) node[pos=0.2, right] {$ a_j $} coordinate[midway] (1-end) coordinate[pos=0.6] (3-start);
\path[draw, ->] (1-start) to (1-end);
\path[fill] (0.1, 1.1) circle[radius=0.05];
\path[draw] (0, 1.2) -- ++(up:1) coordinate[pos=0.2] (m3) coordinate[pos=0.4] (2-start) coordinate(A);
\path[draw] (0.2, 1.2) -- ++(up:1) coordinate[pos=0.2] (m4) coordinate[pos=0.4] (3-end) coordinate(B);
\path[draw, thick, bend right=90, looseness=1.5] (m3) to (m1);
\path[draw, thick, bend right=90, looseness=1.5] (m2) to (m4);
\path[draw, dashed] (B) to[out=30, in=270] ++(30:1.5) to[out=90, in=0] ($ (A)!0.5!(B) + (up:1.2) $) node[below, align=center] {contractible \\ segment} to[out=180, in=90] ($ (A)+(150:1.5) $) to[out=270, in=150] (A);
\end{tikzpicture}
\caption{$ \id $ (ID) with balloon}
\end{subfigure}
\caption{Illustration of angles with balloons}
\label{fig:trick-protocol-balloons}
\end{figure}
\begin{figure}
\centering
\begin{subfigure}{0.3\linewidth}
\centering
\begin{tikzpicture}
\path[draw] (0, 0)  coordinate(F-start) -- ++(150:1) coordinate(A-end) coordinate[pos=0.3] (6-end) coordinate[pos=0.7] (1-start) ++(up:0.2) coordinate(A-start) -- ++(up:1)  coordinate(B-end) coordinate[pos=0.3] (1-end) coordinate[pos=0.7] (2-start)  ++(up:0.2)  coordinate(B-start) -- ++(30:1)  coordinate(C-end)  coordinate[pos=0.3] (2-end) coordinate[pos=0.7] (3-start)  ++(right:0.2)  coordinate(C-start) -- ++(330:1)   coordinate(D-end) coordinate[pos=0.3] (3-end) coordinate[pos=0.7] (4-start)  ++(down:0.2)  coordinate(D-start) -- ++(down:1)  coordinate(E-end) coordinate[pos=0.3] (4-end) coordinate[pos=0.7] (5-start)  ++(down:0.2)  coordinate(E-start) -- ++(210:1)  coordinate(F-end) coordinate[pos=0.3] (5-end) coordinate[pos=0.7] (6-start);
\path[draw, ->, bend right=60] (1-start) to (1-end);
\path[draw, ->, bend right=60] (2-start) to (2-end);
\path[draw, ->, bend right=60] (3-start) to (3-end);
\path[draw, ->, bend right=60] (4-start) to (4-end);
\path[draw, ->, bend right=60] (5-start) to (5-end);
\path[draw, ->, bend right=60] (6-start) to (6-end);
\path[draw, dashed] (A-start) to[out=135, in=90] ++(left:1) to[out=270, in=225] (A-end);
\path[draw, dashed] (B-start) to[out=135, in=90] ++(left:1) to[out=270, in=225] (B-end);
\path[draw, dashed] (C-start) to[out=30, in=0] ++(-0.1, 0.5) to[out=180, in=150] (C-end);
\path[draw, dashed] (D-end) to[out=45, in=90] ++(right:1) to[out=270, in=315] (D-start);
\path[draw, dashed] (E-end) to[out=45, in=90] ++(right:1) to[out=270, in=315] (E-start);
\path[draw, dashed] (F-start) to[out=240, in=180] ++(0.1, -0.5) to[out=0, in=330] (F-end);
\end{tikzpicture}
\caption{Fictitious all-in orbigon}
\label{fig:trick-protocol-allin}
\end{subfigure}
\hspace{0.05\linewidth}
\begin{subfigure}{0.3\linewidth}
\centering
\begin{tikzpicture}
\path[draw] (0, 0)  coordinate(F-start) -- ++(150:1) coordinate[pos=0.4] (m1) coordinate(A-end) coordinate[pos=0.3] (6-end) coordinate[pos=0.7] (1-start) ++(up:0.2) coordinate(A-start) -- ++(up:1)  coordinate(B-end) coordinate[pos=0.3] (1-end) coordinate[pos=0.7] (2-start)  ++(up:0.2)  coordinate(B-start) -- ++(30:1)  coordinate(C-end)  coordinate[pos=0.3] (2-end) coordinate[pos=0.7] (3-start)  ++(right:0.2)  coordinate(C-start) -- ++(330:1)   coordinate(D-end) coordinate[pos=0.3] (3-end) coordinate[pos=0.7] (4-start)  ++(down:0.2)  coordinate(D-start) -- ++(down:1)  coordinate(E-end) coordinate[pos=0.3] (4-end) coordinate[pos=0.7] (5-start)  ++(down:0.2)  coordinate(E-start) -- ++(210:1)  coordinate(F-end) coordinate[pos=0.3] (5-end) coordinate[pos=0.7] (6-start) coordinate[pos=0.6] (m2);
\path[draw, ->, bend right=60] (1-start) to (1-end);
\path[draw, ->, bend right=60] (2-start) to (2-end);
\path[draw, ->, bend right=60] (3-start) to (3-end);
\path[draw, ->, bend right=60] (4-start) to (4-end);
\path[draw, ->, bend right=60] (5-start) to (5-end);
\path[draw, dashed] (A-start) to[out=135, in=90] ++(left:1) to[out=270, in=225] (A-end);
\path[draw, dashed] (B-start) to[out=135, in=90] ++(left:1) to[out=270, in=225] (B-end);
\path[draw, dashed] (C-start) to[out=30, in=0] ++(-0.1, 0.5) to[out=180, in=150] (C-end);
\path[draw, dashed] (D-end) to[out=45, in=90] ++(right:1) to[out=270, in=315] (D-start);
\path[draw, dashed] (E-end) to[out=45, in=90] ++(right:1) to[out=270, in=315] (E-start);
\path[draw] (F-start) -- ++(210:1) coordinate[pos=0.4] (m3);
\path[draw] (F-end) -- ++(330:1) coordinate[pos=0.4] (m4);
\path[draw, thick, bend right] (m1) to (m3);
\path[draw, thick, bend right] (m4) to (m2);
\end{tikzpicture}
\caption{Correct some-out orbigon}
\label{fig:trick-protocol-someout}
\end{subfigure}
\caption{Illustration of orbigons made of angles with balloons}
\end{figure}

Angles with balloons are depicted in \autoref{fig:trick-protocol-balloons}. The balloons facilitate a lot of tricks and desirable properties. With the above/below terminology from \autoref{def:trick-situations-above}, we can say that the balloon always lies above the angle. Moreover, an $ \id $ angle with balloon directly determines the turning directions of its source and target strands:

\begin{lemma}
\label{th:trick-protocol-idstrands}
Let $ \id $ (ID/A'') be an identity angle with balloon. Then above, its source strand turns left and its target strand turns right. Similarly below, its source strand turns right and its target strand turns left.
\end{lemma}

\begin{proof}
In both ID and A'' situations, the turning on the above side is predetermined. For the below side, note that turning in the opposite direction would immediately constitute a teardrop in the underlying curve of $ X $, contradicting the assumption that there are no teardrops. This finishes the proof.
\end{proof}

We now examine products and higher products of angles with balloons. When $ h_1, …, h_k $ are elementary angles, then an (additive) \emph{contribution} to a product $ μ^{k≥2}_r (h_k, …, h_1) $ is simply the product itself in case $ k = 2 $, or an orbigon contributing to the product in case $ k ≥ 3 $. As a first step, we can show that when a higher product of angles with balloons are taken, no contribution to the higher product is all-in:

\begin{lemma}
\label{th:trick-protocol-outside}
Let $ h_1, …, h_k $ be a sequence of $ k ≥ 3 $ angles with balloons. Then any contribution to $ μ^k_r (h_k, …, h_1) $ is first-out or final-out. At the concluding puncture, the first and final strands of $ X $ both turn outside the disk and their turning angles are disjoint.
\end{lemma}

\begin{proof}
Regard a given contribution, denoted $ D $. It is our task to show that $ D $ is not all-in. We have depicted a fictitious all-in contribution in \autoref{fig:trick-protocol-allin}. With this figure, the reason that $ D $ cannot be all-in is immediate: All interior angles of the orbigon are angles with balloons. Correspondingly, the $ k $-many $ X $ segments above the interior angles are all contractible. Since the orbigon itself is contractible, we conclude that the curve underlying $ X $ is contractible, in contradiction with our assumption.

We conclude that every contribution is first-out or final-out. It is our task to show that the first and final strands of the orbigon both outside the disk at the concluding puncture, instead of inside. This is a standard argument concerning contractibility. We have depicted the desired situation in \autoref{fig:trick-protocol-someout}. Towards a contradiction, assume that one of the strands turns inside the disk. View the orbigon as a bracketed discrete immersed disk. Denote the polygon immersion of the discrete disk by $ D: P_k → S $. Then at the concluding puncture, the strand of $ X $ turns inside the interior of $ P_k $. Tracing the strand further, at some point it necessarily intersects itself or leaves $ P_k $. This constitutes a teardrop for the underlying curve of $ X $, in contradiction with the assumption that $ X $ has no teardrop. This shows that both strands turns outwards.

Finally, let us explain why the two turning angles outside the disk are disjoint. Indeed, if the two turning angles outside the disk are not disjoint, the orbigon with its interior angles constitutes a teardrop for the underlying curve of $ X $, in contradiction with the assumption. This finishes the proof.
\end{proof}

According to \autoref{th:trick-protocol-outside}, the situation around the concluding puncture of an orbigon gives rise to an A/A'' situation, which we call the \emph{concluding situation} of the contribution. Thanks to this characterization, we can inspect the contributions to $ μ^{≥2}_r (h_k, …, h_1) $ in more detail, depending on the type of the angles $ h_1, …, h_k $. We have expressed this inspection in \autoref{th:trick-protocol-products}. There, we have detailed only the statement of the first item, because the others statements are analogous and can be interpreted from the figures.

\input{trick/fig_protocol.tex}

\begin{lemma}
\label{th:trick-protocol-products}
Let $ h_1, …, h_k $ be a sequence of angles with balloons. Let $ D $ be a contribution to the product $ μ^{≥2}_r (h_k, …, h_1) $. Then $ D $ can impossibly concern a product of the type $ μ^2 (β ℓ^m, β ℓ^n) $, $ μ^2 (β ℓ^m, α' ℓ^n) $, $ μ^2 (α' ℓ^m, β ℓ^n) $ or $ μ^2 (\id ℓ^m, \id ℓ^n) $. Instead, it falls under one of the following cases, and we make the following refining statements:
\begin{itemize}
\item If $ D $ concerns a product $ μ^{≥3}_r (β ℓ^m \text{(A)}, …) $ with $ β ℓ^m $ final-out: Let $ (α_1, β_1, γ_1) $ be the angles associated with the situation of $ β $. The concluding situation $ (α_2, β_2, γ_2) $ of $ D $ is necessarily of type A. The angles $ β_1 α_1 $ and $ α_2 $ are disjoint. Depending on the size of the angle between $ β_1 $ and $ α_2 $, we obtain a well-defined A situation $ (α_R, β_R, γ_R) $ or an A'' situation $ (α_R, γ_R, β'_R) $ (\autoref{fig:trick-protocol-betaA_finalout}).
\item If $ D $ concerns a product $ μ^{≥3}_r (β ℓ^m \text{(A')}, …) $ with $ β ℓ^m $ final-out: The concluding situation of $ D $ is of type A and we obtain a well-defined A situation $ (α_R, β_R, γ_R) $ or an A'' situation $ (α_R, γ_R, β'_R) $ (\autoref{fig:trick-protocol-betaA'_finalout}).
\item If $ D $ concerns a product $ μ^{≥3}_r (…, β ℓ^m \text{(A)}) $ with $ β ℓ^m $ first-out: The concluding situation of $ D $ is of type A and we obtain a well-defined A situation $ (α_R, β_R, γ_R) $ or A'' situation $ (α_R, γ_R, β'_R) $ (\autoref{fig:trick-protocol-betaA_firstout}).
\item If $ D $ concerns a product $ μ^{≥3}_r (…, β ℓ^m \text{(A')}) $ with $ β ℓ^m $ first-out: The concluding situation of $ D $ is of type A and we obtain a well-defined A situation $ (α_R, β_R, γ_R) $ or A'' situation $ (α_R, γ_R, β'_R) $ (\autoref{fig:trick-protocol-betaA'_firstout}).
\item If $ D $ concerns a product $ μ^{≥3}_r (α' ℓ^m, …) $ with $ α' ℓ^m $ (D) final-out: We obtain a well-defined A or A'' situation (\autoref{fig:trick-protocol-alpha'_finalout}).
\item If $ D $ concerns a product $ μ^{≥3}_r (…, α' ℓ^m) $ with $ α' ℓ^m $ (D) final-out: We obtain a well-defined A or A'' situation (\autoref{fig:trick-protocol-alpha'_firstout}).
\item If $ D $ concerns a product $ μ^2 (\id ℓ^m \text{(ID/A'')}, βℓ^n \text{(A)}) $: We obtain a well-defined A situation or A'' situation (\autoref{fig:trick-protocol-betaA_id}).
\item If $ D $ concerns a product $ μ^2 (\id ℓ^m \text{(ID/A'')}, βℓ^n \text{(A')}) $: We obtain a well-defined A situation or A'' situation (\autoref{fig:trick-protocol-betaA'_id}).
\item If $ D $ concerns a product $ μ^2 (β ℓ^m \text{(A)}, \id ℓ^n \text{(ID/A'')}) $: We obtain a well-defined A situation or A'' situation (\autoref{fig:trick-protocol-id_betaA}).
\item If $ D $ concerns a product $ μ^2 (β ℓ^m \text{(A')}, \id ℓ^n \text{(ID/A'')}) $: We obtain a well-defined A situation or A'' situation (\autoref{fig:trick-protocol-id_betaA'}).
\item $ D $ may also concern a product of one of the types
\begin{align*}
& μ^2 (α, α'ℓ^m) , \quad μ^2 (α'ℓ^m, α) , \quad μ^2 (α, \id ℓ^m \text{(ID/A'')}) ,\quad   μ^2 (\id ℓ^m \text{(ID/A'')}, α), \\
& μ^2 (α' ℓ^m, \id ℓ^n \text{(ID/A'')}) , \quad  μ^2 (\id ℓ^m \text{(ID/A'')}, α' ℓ^n) , \quad μ^2 (β ℓ^m \text{(A/A')}, α) \quad \text{or} \quad μ^2 (α, β ℓ^m \text{(A/A')}) .
\end{align*}
\end{itemize}
In all cases, the $ β $ angle of the resulting A or A'' situation comes again with a balloon.
\end{lemma}

\begin{proof}
Our first task is to show that there are no products of the type $ μ^2 (β ℓ^m, β ℓ^n) $. But this is obvious from the definition and exactly the same as in the case of zigzag paths: The two $ β $ angles would need to wind around the same puncture $ q $. In order to have the product $ μ^2 (β ℓ^m, β ℓ^n) $, $ X $ necessarily turns right at $ q $ when viewed from $ β ℓ^m $, but needs to turn left when viewed from $ β ℓ^n $. This shows there is no product $ μ^2 (β ℓ^m, β ℓ^n) $. Similarly, there are no products $ μ^2 (α' ℓ^m, β ℓ^n) $ and $ μ^2 (β ℓ^m, α' ℓ^n) $. It is also easy to see that there is no product of type $ μ^2 (\id ℓ^m, \id ℓ^n) $. Indeed, from \autoref{th:trick-protocol-idstrands} we obtain two contradicting statements regarding the turning of the target strand of the first, equivalently source strand of the second identity.

Let us now filter out a few possible contributions to higher products. Thanks to \autoref{th:trick-protocol-outside}, we already know that there are no all-in contributions. Furthermore, the two strands both turn outside the disk, so there are no options to form a $ μ^{≥3}_r (α, …) $ or $ μ^{≥3}_r (…, α) $ product.

We now dedicate ourselves to working through the list of viable products and verifying their properties. We highly recommend taking the figures as a visual aid for the arguments.

We start by regarding a result $ D $ of a product $ μ^{≥3}_r (β ℓ^m \text{(A/A')}, …) $ with $ β ℓ^m $ final-out. The strand of the first arc of the orbigon turns right at the concluding puncture. Its turning angle is clearly at most as long as the remaining part of $ β $, because otherwise the orbigon $ D $ and the balloon of $ β $ would constitute a teardrop. We obtain a well-defined child: If the turning angle is strictly shorter, we obtain an A situation. If the turning angle is equally long, we obtain instead an A'' situation. In both cases, the resulting $ β_R $ or $ \id_R $ comes with a balloon, obtained as the joining of the balloon of $ β $ and the orbigon $ D $. The case of $ μ^{≥3}_r (…, β) $ with $ β $ first-out is very similar to the case of $ μ^{≥3}_r (β, …) $.

Regard now a result $ D $ of a product $ μ^{≥3}_r (α' ℓ^m, …) $ with $ α' ℓ^m $ final-out or $ μ^{≥3}_r (…, α' ℓ^m) $ with $ α' ℓ^m $ first-out. The either case, we obtain an A or A'' situation at the concluding puncture by \autoref{th:trick-protocol-outside}. In case of final-out $ α' ℓ^m $, the result $ D $ is of the form $ β_R α_R ℓ^m $ (A) or $ α_R ℓ^m $ (A''). In case of first-out $ α' ℓ^m $, the result $ D $ is of the form $ γ_R β_R ℓ^m $ (A) or $ γ_R ℓ^m $ (A''). In either case, the $ β_R $ or $ \id_R $ angle of the resulting situation comes with a balloon.

Now regard a product $ μ^2 (β ℓ^m \text{(A/A')}, \id ℓ^n) $. By definition, the identity comes with a balloon. The $ β $ angle does not wind around the above end of $ \id $, but around the below end of $ \id $, because the target strand of $ \id $ turns right at the above puncture, rendering a composition with $ β $ impossible. Next, we note that the target strand of the identity turns right at the below side of the source strand and the turning angle is at most $ β $, for otherwise the combination of the balloons of the identity and $ β $ would constitute a teardrop. As a result, space remains for defining a type A situation $ (α_R, β_R, γ_R) $ or type A'' situation $ (α_R, \id_R, γ_R, β'_R) $. Its $ β_R $ or $ \id_R $ angle comes with a balloon again, namely the combination of the balloons of $ β $ and $ \id $. The case of $ μ^2 (\id ℓ^m, β ℓ^n\text{(A/A')}) $ is similar.

For the products $ μ^2 (α, α' ℓ^m) $ and $ μ^2 (α' ℓ^m, α) $, we are not supposed to define anything. Note that it is clear that for this product to exist, $ α $ and $ α' $ need to be precisely complementary angles. For a product $ μ^2 (β ℓ^m \text{(A/A')}, α) $ or $ μ^2 (α, β ℓ^m \text{(A/A')}) $, it is interesting to note that $ α $ must be the $ α $ or $ γ $ angle from the same situation as $ β $. This finishes the proof.
\end{proof}

\subsection{Flowers}
\label{sec:trick-flower}
In this section, we construct the deformed twisted differential $ δ_q $ which uncurves $ X $. We start by introducing flowers, a combinatorical gadget that recursively keeps track of all contractible segments of $ X $. Such a flower includes by construction all the terms we need to insert into $ δ_q $. More precisely, we define $ δ_q $ as the sum over the values of all flowers of $ X $ and show that $ X_q = (⊕ a_i [s_i], δ_q) $ is indeed curvature-free.

\begin{center}
\begin{tikzpicture}
\path (0, 0) node (A) {Flowers of $ X $} (6, 0) node (B) {Terms for $ δ_q $};
\path[draw, ->] ($ (A.east)!0.2!(B.west) $) -- ($ (A.east)!0.8!(B.west) $);
\end{tikzpicture}
\end{center}

The core idea of our construction is best explained as follows. We start with $ δ_q $ containing only the angles $ α_i $ and their complements $ r ℓ^{-1} α_i' $. This already makes $ μ^0_X $ and $ μ^2 (δ_q, δ_q) $ cancel, however we get a potentially unlimited amount of orbigon contributions from $ μ^{k≥3}_r (δ_q, …, δ_q) $. For each of these terms, we need to insert an additional term into $ δ_q $ in order to make it cancel out. The new terms inserted into $ δ_q $ however can give rise to further disturbing terms in $ μ^{k≥3}_r (δ_q, …, δ_q) $ and we need to iteratively repeat this process.

In every step of the process, we should remember the entire history of how a given term was formed, much like the notion of tails or result components for zigzag paths. The tool of flowers which we define here systematically keeps track of the appearing terms. Since a term typically appears recursively for every orbigon that can be formed from already existing terms, the orbigons get stitched together much like a flower. All the terms we insert into $ δ_q $ come from orbigons, in particular we can ensure inductively that they come naturally with balloons. This way \autoref{th:trick-protocol-products} applies and facilitates the friction-free definition of flowers.

\begin{figure}
\centering
\begin{subfigure}[b]{0.25\linewidth}
\centering
\begin{tikzpicture}[yscale=0.5, xscale=0.8]
\path[draw] (0, 0) -- ++(up:1.5) coordinate[midway] (m1) coordinate (stop-left);
\path[draw] (0.2, 0) -- ++(up:1.5) coordinate[midway] (m2) coordinate (stop-right);
\path[draw, ->] (m1) -- (m2);
\path[draw, dashed, looseness=2] (stop-right) to[out=45, in=0] ($ (stop-right)!0.5!(stop-left) + (up:2) $) to[out=180, in=135] (stop-left);
\path[fill] (0.1, -0.1) circle[radius=0.05];
\foreach \i in {1.7, 3.4, 5.1} {
\begin{scope}[shift={(down:\i)}]
\path[fill] (0.1, -0.1) circle[radius=0.05];
\path[draw] (0, 0) -- ++(up:1.5) coordinate[midway] (m1) coordinate (stop-left) coordinate[pos=0.7] (arc-left);
\path[draw] (0.2, 0) -- ++(up:1.5) coordinate[midway] (m2) coordinate (stop-right) coordinate[pos=0.7] (arc-right);
\path[draw, ->] (m1) -- (m2);
\path[draw, thick, bend right] (arc-right) arc(-90:90:0.5);
\path[draw, thick, bend left] (arc-left) arc(270:90:0.5);
\end{scope}}
\path (0.8, -4.5) node {$ \id_S $};
\end{tikzpicture}
\caption{An ID stem flower}
\label{fig:trick-flower-IDstem}
\end{subfigure}
\begin{subfigure}[b]{0.25\linewidth}
\centering
\begin{tikzpicture}
\path[draw] (0, 0) -- ++(right:1) coordinate[midway] (beta-start) -- ++(100:1.5) coordinate (9) ++(up:0.2) coordinate (10) -- ++(120:1) coordinate (11) ++(up:0.2) coordinate (12) -- ++(45:1) coordinate  (13) ++(right:0.2) coordinate (14) -- ++(315:1) coordinate (15) ++(down:0.2) coordinate (16) -- ++(240:1) coordinate (17) ++(down:0.2) coordinate (18) -- ++(260:1.5) -- ++(right:1) coordinate[midway] (beta-end);
\foreach \i\j\k in {9/10/170, 11/12/120, 13/14/90, 15/16/60, 17/18/10} {%
\path[draw, dashed] (\j) to[bend right=80, looseness=3] ($ (\i)!0.5!(\j) + (\k:0.4) $) to[bend right=80, looseness=2] (\i);};
\path[draw, bend right=80, looseness=1.5, ->] (beta-start) to node[midway, above] {$ β_R $} (beta-end);
\end{tikzpicture}
\caption{An orbigon flower}
\end{subfigure}
\begin{subfigure}[b]{0.4\linewidth}
\centering
\begin{tikzpicture}
\path[draw] (0, 0) -- ++(60:1) coordinate[midway] (beta-start) -- ++(210:1) coordinate (1) ++(left:0.2) coordinate (2) -- ++(150:0.8) coordinate (3) ++(150:0.2) coordinate (4) -- ++(120:0.5) coordinate (5) ++(up:0.2) coordinate (6) -- ++(30:1) coordinate (7) ++(right:0.2) coordinate (8) -- ++(right:1) -- ++(300:1) -- ++(100:1.5) coordinate (9) ++(up:0.2) coordinate (10) -- ++(120:1) coordinate (11) ++(up:0.2) coordinate (12) -- ++(45:1) coordinate  (13) ++(right:0.2) coordinate (14) -- ++(315:1) coordinate (15) ++(down:0.2) coordinate (16) -- ++(240:1) coordinate (17) ++(down:0.2) coordinate (18) -- ++(260:1.5) -- ++(60:1) -- ++(0:1) coordinate (19) ++(right:0.2) coordinate (20) -- ++(330:1) coordinate (21) ++(down:0.2) coordinate (22) -- ++(240:0.5) coordinate (23) ++(210:0.2) coordinate (24) -- ++(210:0.8) coordinate (25) ++(left:0.2) coordinate (26) -- ++(150:1) -- ++(300:1) coordinate[midway] (beta-end);
\foreach \i\j\k in {1/2/240, 3/4/240, 5/6/180, 7/8/120, 9/10/170, 11/12/120, 13/14/90, 15/16/60, 17/18/10, 19/20/60, 21/22/0, 23/24/300, 25/26/300} {%
\path[draw, dashed] (\j) to[bend right=80, looseness=3] ($ (\i)!0.5!(\j) + (\k:0.4) $) to[bend right=80, looseness=2] (\i);};
\path[draw, ->, bend right] (beta-start) to node[midway, above] {$ β_R $} (beta-end);
\end{tikzpicture}
\caption{A compound flower}
\end{subfigure}
\caption{Illustration of flowers}
\end{figure}
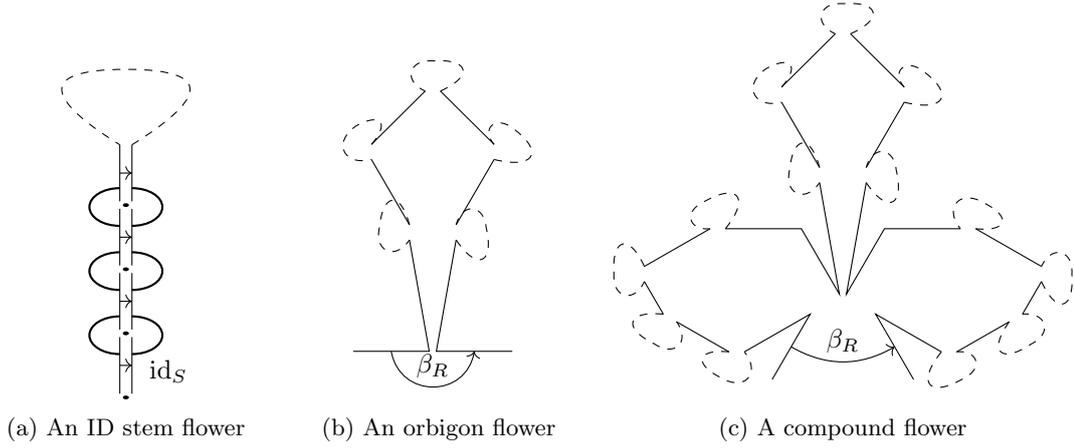

Not only orbigons need to be taken into account. It is possible that at the concluding puncture of an orbigon, the two strands which were separated by the orbigon now come together and keep traveling in parallel for a while, as in \autoref{fig:trick-flower-IDstem}. On such occasions, we have to insert a whole sequence of identities into $ δ_q $. After a while, the two strands may separate again and wildly continue forming orbigons.

The rule of thumb could be memorized as follows:
\begin{itemize}
\item $ \id $ (A'') needs to be inserted when the strands come together.
\item $ \id $ (ID) needs to be inserted when the strands keep running together.
\item $ β $ (A') needs to be inserted when the strands separate. 
\item $ β $ (A) needs to be inserted when the strands come together and immediately separate again.
\end{itemize}

The construction of flowers is so technical because of the large amount of complexity observed while constructing $ δ_q $. Explicitly, it concerns the following complications: First, segments of $ X $ can bound multiple orbigons. This makes that there is no linear way of enumerating the terms we need to add to $ δ_q $. Instead, the terms will “cross-pollinate” each other. Second, for every orbigon that can be made of terms already present in $ δ_q $, we typically add a new $ β $ (A) angle to $ δ_q $. This $ β $ angle already creates two new products $ μ^2 (β, α) $ and $ μ^2 (γ, β) $. Cancelling them entails working with four terms in total. Third, the identities need separate creation and cancellation procedures.

\begin{remark}
We write
\begin{equation*}
r = \sum_{\substack{m ≥ 1 \\ q ∈ M}} r_{q, m} ℓ_q^m ∈ \mathfrak{m} Z(\Gtl \cA).
\end{equation*}
Recall that an orbigon $ D $ comes with a deformation parameter $ r_D $. The element $ r_D $ lies in the deformation base $ B $ and is defined as the product of all $ r_{q, m} $ ranging over all orbifold points of $ D $.
\end{remark}

\begin{definition}
A \emph{flower} comes with the datum of an \emph{concluding situation}, which is an A, A', A'' or ID situation. A flower always comes with the datum of a \emph{value}, which is a $ B $-multiple of $ \id ℓ^m $ (A''), $ \id ℓ^m $ (ID), $ βℓ^m $ (A) or $ βℓ^m $ (A') for some $ m ≥ 0 $, depending on the type of the concluding situation.

We distinguish five types of flowers, namely $ α $ flowers, $ α' $ flowers, orbigon flowers, ID stem flowers and compound flowers. The orbigon flowers remember which flowers their orbigon is made of and the compound flowers each remember which flowers they form a compound of. The complete recursive definition of flowers is given in the catalog below. Wherever a flower $ F $ is operated on in a formula, its value is meant.
\end{definition}

\begin{description}
\item[$ α $ flowers:] Any $ α $ (D) angle of $ X $ determines a flower $ F $.
\begin{itemize}
\item The concluding situation of $ F $ is the type D situation determined by $ α $.
\item The value of $ F $ is $ α $.
\end{itemize}
\item[$ α' $ flowers:] Any complementary angle $ α'ℓ^{m-1} $ (D) with $ m ≥ 1 $ determines a flower $ F $.
\begin{itemize}
\item The concluding situation of $ F $ is the type D situation determined by $ α' $.
\item The value of $ F $ is $ r_{q, m} α' ℓ^{m-1} $.
\end{itemize}
\item[Orbigon flowers:] Let $ F_1, …, F_k $ be a sequence of flowers, together with the datum of a type A situation $ (α_R, β_R, γ_R, β'_R) $ or a type A'' situation $ (α_R, \id_R, γ_R, β'_R) $, two integers $ m ≥ 0 $ and $ n ≥ 1 $ and an orbigon $ D $ whose interior angles are the values of $ F_1, …, F_k $ together with $ β'_R ℓ^m $. Then this data defines a new flower $ F $.
\begin{itemize}
\item The concluding situation of $ F $ is the type A situation $ (α_R, β_R, γ_R, β'_R) $ or the A'' situation $ (α_R, \id_R, γ_R, β'_R) $.
\item The value of $ F $ is a multiple of $ β_R ℓ^n $ (A) or $ \id_R ℓ^n $ (A''), depending on whether the concluding puncture of $ F $ is of type A or A''. Denote by $ q $ the puncture in the middle of the concluding situation and by $ r_D ∈ B $ the deformation parameter of the orbigon $ D $. Let $ \coeff{F_i} ∈ B $ denote the coefficient of the value of $ F_i $, stripping away the information which angle it concerns.

In case the given situation is of type A, the precise value of $ F $ is defined as
\begin{equation*}
(-1)^{‖β_R α_R β'_R‖ \sum_{i = 1}^k ‖F_i‖ + \sum_{i<j} ‖F_i‖‖F_j‖ + ‖β_R‖‖α‖ + ‖α‖} r_D r_{m+n+1, q} \coeff{F_k} … \coeff{F_1} β_R ℓ^n.
\end{equation*}
In case the given situation is instead of type A'', the precise value of $ F $ is defined as
\begin{equation*}
(-1)^{‖α_R β'_R‖ \sum_{i = 1}^k ‖F_i‖ + \sum_{i<j} ‖F_i‖‖F_j‖} r_D r_{q, m+n+1} \coeff{F_k} … \coeff{F_1} \id_R ℓ^n.
\end{equation*}
\end{itemize}
\item[ID stem flowers:] Let $ F_1 $ be a flower whose concluding situation is an ID situation $ (\id_R, α_R, γ_R) $ or an A'' situation $ (α_R, \id_R, γ_R, β'_R) $. Depending on whether the source and target strands of $ \id_R $ separate below the situation, we obtain a new A' situation $ (α_S, β_S, γ_S) $ or ID situation $ (α_S, \id_S, γ_S) $. Let $ F_2 $ be an $ α' $ flower whose value is the complement of $ γ_S $. This defines a new flower $ F $.
\begin{itemize}
\item The concluding situation of $ F $ is the A' situation $ (α_S, β_S, γ_S) $ or the ID situation $ (α_S, \id_S, γ_S) $.
\item The value of $ F $ is a multiple of $ β ℓ^m $ (A') or $ \id_S ℓ^m $ (ID), depending on whether the concluding situation is of type A' or ID. In case the situation is of A' type, the precise value of $ F $ is
\begin{equation*}
(-1)^{‖F_2‖ ‖F_1‖ + |F_1| + ‖β_S‖ ‖α_S‖ + |α_S|} F_2 F_1 α_S^{-1}.
\end{equation*}
In case the situation is of ID type, the precise value of $ F $ is
\begin{equation*}
(-1)^{‖F_2‖ ‖F_1‖ + |F_1|} F_1 F_2 α_S^{-1}.
\end{equation*}
\end{itemize}
\item[Compound flowers:] Let $ F_1, …, F_k $ be a sequence of $ k ≥ 2 $ flowers of ID stem or orbigon type with concluding situations $ (α_i, β_i, γ_i) $ (A/A'), such that (a) $ α_i $ agrees with $ γ_{i+1} $ and (b) the angle $ β'_R $ given by the union of the $ β'_i $ (or $ \id_i $) angles, including all $ α_i $ and $ γ_i $, is at most a full turn. Then this data defines a new compound flower $ F $.
\begin{itemize}
\item The concluding situation of $ F $ is the type A or A'' situation $ (α_R, β_R, γ_R, β'_R) $ or $ (α_R, \id_R, γ_R, β'_R) $ which contains the $ β'_R $ just constructed.
\item The value of $ F $ is a multiple of $ β_R ℓ^n $ (A) or $ \id_R ℓ^n $ (A''), depending on whether the concluding situation is of type A or A''. To determine the precise coefficient, we shall use a trick by applying recursion. Regard the flowers $ F_1, …, F_{k-1} $. In case $ k = 2 $, this is the single ID stem or orbigon flower $ F_{1, 1} ≔ F_1 $. In case $ k ≥ 3 $, this sequence determines a (smaller) compound flower $ F_{1, k-1} $. In either case, we can assume that we already know the value of $ F_{1, k-1} $.

We distinguish two similar cases: Assume $ F_k $ itself is an orbigon flower whose interior angles come from the flowers $ G_1, …, G_l $. Let $ m ≥ 0 $ be such that $ β'_R ℓ^m $ together with $ G_1, …, G_l $ are the interior angles of the orbigon $ F_k $. Then we define the value of $ F $ as
\begin{equation*}
(-1)^{‖F_{1, k-1}‖ \sum_{j = 1}^l ‖G_j‖ + ‖β_R‖ ‖α_R‖ + |α_R|} F_{1, k-1} (β'_R ℓ^m)^{-1}.
\end{equation*}
Assume otherwise that $ F_k $ itself is an ID stem flower. Then we define the value of $ F $ as
\begin{equation*}
(-1)^{‖F_{1, k-1}‖ + ‖β_R‖ ‖α_R‖ + |α_R|} F_{1, k-1} α_R^{-1}.
\end{equation*}
\end{itemize}
\end{description}

With this sophisticated construction of flowers, we can define $ δ_q $ simply as the sum over the values of all flowers of $ X $:

\begin{definition}
We put
\begin{equation*}
δ_q = \sum_{F \text{ flower of } X} F \quad ∈ \quad \Hom_{\Add\Gtl_r \cA}^1 (X, X).
\end{equation*}
\end{definition}

\begin{lemma}
The element $ δ_q $ is well-defined and its leading term is $ δ $.
\end{lemma}

\begin{proof}
Let us explain why $ δ_q $ is well-defined. The first observation is that there are only finitely many orbigons for a given sequence of interior angles. In conclusion, for any $ C > 0 $ there is only a finite number of flowers $ F $ for which it has taken at most $ C $ recursive steps to define $ F $. The second observation is that the value of any flower lies in a certain power of the maximal ideal $ \mathfrak{m} ⊂ B $. Every time a new flower is formed, its $ \mathfrak{m} $-adic exponent increases. Together, both observations show that the series which defines $ δ_q $ converges in the $ \mathfrak{m} $-adic topology.

It is very easy to see that the leading term of $ δ_q $ is $ δ $: The only flowers whose value does not lie in a power of the maximal ideal are the $ α $ flowers. Since $ δ $ is just the sum of all $ α $ flowers, this finishes the proof.
\end{proof}


We aim at showing that $ \sum_{k ≥ 0} μ^k_{\Add\Gtl_r \cA} (δ_q, …, δ_q) = 0 $. In order to flexibly cancel terms in this sum, we introduce an obvious notion of result components:

\begin{definition}
A \emph{result component} of $ \sum_{k ≥ 0} μ^k_{\Add\Gtl_r \cA} (δ_q, …, δ_q) $ consists of a sequence of $ k ≥ 0 $ flowers $ F_1, … F_k $ together with an additive component of $ μ^k_{\Add\Gtl_r \cA} (F_k, …, F_1) $. More precisely,
\begin{itemize}
\item for $ k = 0 $ this entails the curvature $ μ^0_X $,
\item for $ k = 2 $ this entails a nonvanishing product $ μ^k_{\Add\Gtl_r \cA} (F_k, …, F_1) $,
\item for $ k ≥ 3 $ this entails a choice of orbigon contributing to $ μ^k_{\Add\Gtl_r \cA} (F_k, …, F_1) $.
\end{itemize}
\end{definition}

\begin{figure}
\centering
\begin{subfigure}[b]{0.3\linewidth}
\centering
\begin{tikzpicture}[scale=1.5]
\path[draw] (210:1) -- (210:0.1) -- (160:1);
\path[draw] (140:1) -- (140:0.1) coordinate[pos=0.7] (beta-start) -- (120:1);
\path[draw] (100:1) -- (100:0.1) -- (80:1);
\path[draw] (60:1) -- (60:0.1) -- (40:1) coordinate[pos=0.3] (beta-end);
\path[draw] (20:1) -- (330:0.1) -- (330:1);
\path[draw, dashed] (160:1) to[bend left=90, looseness=3] (140:1);
\path[draw, dashed] (120:1) to[bend left=90, looseness=3] (100:1);
\path[draw, dashed] (80:1) to[bend left=90, looseness=3] (60:1);
\path[draw, dashed] (40:1) to[bend left=90, looseness=3] (20:1);
\path (150:1) node {$ F_4 $};
\path (110:1) node {$ F_3 $};
\path (70:1) node {$ F_2 $};
\path (30:1) node {$ F_1 $};
\end{tikzpicture}
\caption{Given compound flower}
\end{subfigure}
\begin{subfigure}[b]{0.3\linewidth}
\centering
\begin{tikzpicture}[scale=1.5]
\path[draw] (210:1) -- (210:0.1) coordinate[pos=0.5] (b-start) -- (160:1) coordinate[pos=0.5] (g-start);
\path[draw] (140:1) -- (140:0.1) coordinate[pos=0.7] (beta-start) -- (120:1);
\path[draw] (100:1) -- (100:0.1) -- (80:1);
\path[draw] (60:1) -- (60:0.1) -- (40:1);
\path[draw] (20:1) -- (330:0.1) coordinate[pos=0.6] (alpha-end) -- (330:1) coordinate[pos=0.3] (beta-end) coordinate[pos=0.5] (alpha-start);
\path[draw, dashed] (160:1) to[bend left=90, looseness=2] (140:1);
\path (150:1) node {$ F_4 $};
\path (110:1) node {$ F_3 $};
\path (70:1) node {$ F_2 $};
\path (30:1) node {$ F_1 $};
\path[draw, bend right=80, looseness=2, ->] (beta-start) to node[pos=0.7, above] {$ F_{1, 3} $} (beta-end);
\path[draw, bend right, ->] (g-start) to node[midway, left] {$ α_4 $} (b-start);
\path[draw, bend right=60, ->] (b-start) to node[midway, below] {$ F $} (alpha-start);
\end{tikzpicture}
\caption{$ μ^{≥3} (F_{1, 3}, …) + μ^2 (F, α_4) $}
\end{subfigure}
\begin{subfigure}[b]{0.3\linewidth}
\centering
\begin{tikzpicture}[scale=1.5]
\path[draw] (210:1) -- (210:0.1) coordinate[pos=0.7] (beta-start) coordinate[pos=0.5] (b-start) -- (160:1) coordinate[pos=0.3] (g-start);
\path[draw] (140:1) -- (140:0.1) -- (120:1);
\path[draw] (100:1) -- (100:0.1) -- (80:1);
\path[draw] (60:1) -- (60:0.1) -- (40:1) coordinate[pos=0.3] (beta-end);
\path[draw] (20:1) -- (330:0.1) coordinate[pos=0.6] (alpha-end) -- (330:1) coordinate[pos=0.5] (alpha-start);
\path[draw, dashed] (40:1) to[bend left=90, looseness=2] (20:1);
\path (150:1) node {$ F_4 $};
\path (110:1) node {$ F_3 $};
\path (70:1) node {$ F_2 $};
\path (30:1) node {$ F_1 $};
\path[draw, bend right=80, looseness=2, ->] (beta-start) to node[pos=0.35, above] {$ F_{2, 4} $} (beta-end);
\path[draw, bend right, ->] (alpha-start) to node[midway, right] {$ γ_1 $} (alpha-end);
\path[draw, bend right=60, ->] (b-start) to node[midway, below] {$ F $} (alpha-start);
\end{tikzpicture}
\caption{$ μ^{≥3} (…, F_{2, 4}) + μ^2 (γ_1, F) $}
\end{subfigure}
\caption{Compound cancellation}
\label{fig:trick-flower-cancellation}
\end{figure}

\begin{proposition}
We have $ \sum_{k ≥ 0} μ^k_{\Add\Gtl_r \cA} (δ_q, …, δ_q) = 0 $. Therefore $ X_q = (⊕ a_i [s_i], δ_q) $ is curvature-free.
\end{proposition}

\begin{proof}
We shall provide a list of cancellations and then check that every result component is contained in this list. We essentially distinguish two types of cancellations, namely the simple and compound cancellations. A sample compound cancellation is depicted in \autoref{fig:trick-flower-cancellation}. The precise list of cancellations reads as follows:
\begin{itemize}
\item The curvature $ μ^0_X $ and the $ μ^2 (α, α') $ and $ μ^2 (α', α) $ result components.
\item \emph{simple cancellation}: Let $ F_1, …, F_k $ be a sequence of flowers. Furthermore, let $ (α_R, β_R, γ_R, β_R') $ be an A situation or $ (α_R, β_R', γ_R) $ be an A'' situation. Then for every orbigon $ D $ with interior angles $ F_1, …, F_k, β_R' ℓ^m $, we have the cancellation
\begin{equation*}
μ^k_{\Add\Gtl_r \cA} (F_k, …, F_1, α_R' ℓ^{m+n}) + μ^k_{\Add\Gtl_r \cA} (γ_R' ℓ^{m+n}, F_k, …, F_1) = 0.
\end{equation*}
By abuse of notation, we have written $ μ(…) $ where we actually refer to the contribution of the very specific orbigon $ D $. We have also denoted by $ α_R' $ (D) and $ γ_R' $ (D) the complementary angles of $ α_R $ and $ γ_R $, in the sense that $ α_R' α_R = ℓ = γ_R' γ_R $.
\item \emph{compound cancellation}: Let $ F $ be a compound flower consisting of flowers $ F_1, …, F_k $.

Denote by $ G_1, …, G_l $ the flowers that the flower $ F_k $ is derived from, and by $ H_1, …, H_n $ the flowers that $ F_1 $ is derived from. Denote the concluding situation of $ F_i $ by $ (α_i, β_i, γ_i, β'_i) $ or $ (α_i, β_i, γ_i, \id_i) $, depending on whether it concerns an A situation or an A' situation.

Regard the compound flower $ F_{1, k-1} $ consisting of $ F_1, …, F_{k-1} $ and the compound flower $ F_{2, k} $ consisting of $ F_2, …, F_k $. In case $ k = 2 $, the flowers $ F_{1, 1} = F_1 $ and $ F_{2, 2} = F_2 $ are simply ID stem or orbigon flowers, instead of compound flowers. Either way, we have the cancellations
\begin{equation}
\label{eq:trick-flower-cancellation}
\begin{aligned}
μ^{≥3}_{\Add\Gtl_r \cA} (F_{1,k-1}, G_l, …, G_1) + μ^2_{\Add\Gtl_r \cA} (F, α_k) &= 0, \\
μ^{≥3}_{\Add\Gtl_r \cA} (H_n, …, H_1, F_{2, k}) + μ^2_{\Add\Gtl_r \cA} (γ_1, F) &= 0.
\end{aligned}
\end{equation}
By abuse of notation, $ μ(F_{1, k-1}), …) $ actually refers to the specific orbigon given by $ F_k $. Similarly, $ μ(…, F_{2, k}) $ refers to the specific orbigon given by $ F_1 $.

In \eqref{eq:trick-flower-cancellation}, we have silently assumed that $ F_1 $ and $ F_k $ are orbigon flowers. In case $ F_1 $ is instead an ID stem flower, the term $ μ^{≥3}_{\Add\Gtl_r \cA} (…, F_{2, k}) $ should read $ μ^2 (\id_1, F_{2, k}) $ instead. Similarly, in case $ F_k $ is an ID stem flower, $ μ^{≥3}_{\Add\Gtl_r \cA} (F_{1,k-1}, …) $ should read $ μ^2 (F_{1,k-1}, \id_k) $ instead.
\end{itemize}
Finally, let us explain why all possible result components of $ \sum_{k ≥ 0} μ^k_r (δ_q, …, δ_q) $ are captured in the above cancellation list. Indeed, all results of flowers naturally come with balloons. Therefore \autoref{th:trick-protocol-products} applies and any result component falls under one of the following cases:
\begin{itemize}
\item Regard a contribution $ D $ to $ μ^{k≥3} (β ℓ^m \text{(A/A')}, …) $ with $ β ℓ^m $ final-out or $ μ^{k≥3} (…, β ℓ^m \text{(A/A')}) $ with $ β ℓ^m $ first-out. By construction, $ β $ comes from the concluding situation of a flower. Whether it concerns an orbigon, ID stem or compound flower, the two terms fall under the compound cancellation.
\item Regard a contribution $ D $ to $ μ^{k≥3} (α' ℓ^m \text{(D)}, …) $ with $ α' ℓ^m $ final-out or $ μ^{k≥3} (…, α' ℓ^m \text{(D)}) $ with $ α' ℓ^m $ first-out. Then all the other interior angles of the orbigon are also the values of flowers and $ D $ falls directly under the simple cancellation.
\item A contribution $ μ^2 (β ℓ^m \text{(A/A')}, \id ℓ^n) $ or $ μ^2 (\id ℓ^m, β ℓ^n \text{(A/A')}) $ falls under the compound cancellation.
\item A contribution $ μ^2 (α' ℓ^m \text{(D)}, \id ℓ^n) $ or $ μ^2 (\id ℓ^m, α' ℓ^n \text{(D)}) $ falls under the compound cancellation.
\item Regard a contribution $ μ^2 (α \text{(D)}, β ℓ^m \text{(A/A')}) $ or $ μ^2 (β ℓ^m \text{(A/A')}, α \text{(D)}) $. Let $ F $ be the flower that $ β $ comes from. Depending on whether $ F $ is an orbigon or ID stem flower or a compound flower, the term falls under the simple or compound cancellation.
\end{itemize}
This shows that all terms have been canceled. It is a basic inspection that all terms have been canceled only once. This shows that $ \sum_{k ≥ 0} μ^k_{\Add\Gtl_r \cA} (δ_q, …, δ_q) = 0 $ and finishes the proof.
\end{proof}

\section{Classification of result components}
\label{sec:classification}
We collect here a few deferred proofs: In \autoref{sec:classification-shape}, we prove \autoref{th:subdisk-assignment} which classifies result components. In \autoref{sec:classification-subdiskshape}, we prove \autoref{th:subdisk-types-shape} which concerns subdisks of CR, ID, DS and DW result components. We spend \autoref{sec:classification-narrow} till \ref{sec:classification-IDDSDW} with a proof of \autoref{th:subdisk-minmodel-bijection} which classifies the image of $ \Subdisk $. In \autoref{sec:classification-signs}, we prove \autoref{th:subdisk-minmodel-signs} concerning the signs of result components.

\subsection{Shape of result components}
\label{sec:classification-shape}
We prove here \autoref{th:subdisk-assignment}, which claims that \autoref{tab:resultcomp-classification} is an exhausting classification of result components. Recall the situation: We are given a result component of an h- or π-tree and are supposed to analyze how it is derived. It is not necessary to find the entire tree it is derived from, but only so far that we recognize it fits the scheme of \autoref{tab:resultcomp-classification}. We will now go through \autoref{tab:resultcomp-classification} case-by-case:

\begin{description}[font=\normalfont]
\item[$ α_0 $, $ \id $ (C), $ α_4 $ from h-trees:] Note that $ α_0 $ and $ \id $ (C) belong to $ H $ and hence only appear as direct morphism or result component of a π-tree. Also, $ α_4 $ does not appear in the disk and multiplication tables as result component of an h-tree and therefore any $ α_4 $ is either direct or a result component of a π-tree in combination with $ α_3 $.

\item[$ α_3 $ from h-trees:] The angle $ α_3 $ does not appear in the disk table as result component of a disk $ h_q μ^{≥3} $. Whenever it appears as the result component of a product $ h_q μ^2 $, it must be as $ h_q μ^2 (α_0, α_4) $ and the involved $ α_0 $ and $ α_4 $ are direct.

\item[$ α_0' $ from h-trees:] The angle $ α_0' $ only appears in the disk and multiplication table as $ h_q μ^2 (α_0, α_0') $. The $ α_0 $ involved is necessarily direct, and inductively we conclude that the $ α_0' $ is the result component of one of the trees of \autoref{fig:subdisk-alpha0p}.

\item[$ β $ (C) from h-trees:] The angle $ β $ (C) is necessarily direct or the result component of $ h_q μ^2 (\id \text{(C)}, α_0') $. The $ \id $ (C) is necessarily direct and we already know $ α_0' $ is the result component of one of the trees in \autoref{fig:subdisk-alpha0p}.

\item[$ β' $ (C) from h-trees:] Same as $ β $ (C).

\item[$ β $ (A) from h-trees:] It appears as tail component in the deformed cohomology basis elements by \autoref{th:deformed-cohomology-basis}. It appears as tail component of $ α_3 $ and $ α_4 $ in the deformed version of $ (-1)^{\#α_3} α_3 + (-1)^{\#α_4 + 1} α_4 $, and as tail component of $ β $ (C) and $ β' $ (C) in the deformed version of $ \id $ (C).

\item[$ β $ (A) main result component of $ h_q μ^{≥3} $:] Its classification follows directly from the disk tables \ref{fig:components-first-out} and \ref{fig:components-final-out}. Note that an all-in disk does not produce a $ β $ (A) result component either.

\item[$ β $ (A) main result component of $ h_q μ^2 $:] According to Table \autoref{tab:components-multiplication}, this product is either of the type $ h_q μ^2 (β \text{(A)}, α_0) $, $ h_q μ^2 (β \text{(A)}, \id \text{(C)}) $ or $ h_q μ^2 (β/β' \text{(C)}, \id \text{(C)}) $. In the first two cases, we inductively trace the $ β $ (A) involved. In the third case, note that we already know the entire tree of $ β/β' $ (C). Ultimately, we end up either with $ β $ (C) or $ β' $ (C) or a $ β $ (A) that is a direct, $ h_q μ^{≥3} $ or tail $ h_q μ^2 $ result component, plus multiple compositions with $ α_0 $ or $ \id $ (C) on the right. All three kinds of trees are depicted in \autoref{fig:subdisk-betaA-mult-main}.

\item[$ β $ (A) tail result component:] If it is the tail result component of an $ h_q μ^2 $, then it is necessarily one of the tail result components marked by +E in \autoref{tab:components-multiplication}. All options come with a corresponding main result component $ β $ (A), $ β/β' $ (C) or $ α_3 $. If it is the tail result component of an $ h_q μ^{≥3} $, then it is the result component of one of the final-out disks in \autoref{fig:components-final-out} and comes with a corresponding $ β $ (A).

\item[$ \id $ (D) result component of an h-tree:] According to the disk and multiplication tables, this concerns a product $ h_q μ^2 (\id \text{(C)}, α_4) $ or $ h_q μ^2 (α_3, \id \text{(C)}) $. The $ α_4 $ and $ \id $ (C) involved are necessarily direct. Furthermore, the $ α_3 $ is direct or the result component of $ h_q μ^2 (α_0, α_4) $. In the latter case, we conclude that $ μ^2 (α_3, \id \text{(C)}) $ equals the co-identity $ α_0 $ itself, and therefore $ h_q μ^2 (α_3, \id \text{(C)}) $ vanishes. The $ \id $ (D) is therefore necessarily a result component of $ h_q μ^2 (\id \text{(C)}, α_4) $ or $ h_q μ^2 (α_3, \id \text{(C)}) $ with all involved $ α_3 $, $ α_4 $ or $ \id $ (C) being direct.

\item[$ \id $ (B) result component of an h-tree:] According to the disk and multiplication tables, this is a first-out disk or one of the trees in \autoref{fig:subdisk-idB}.

\item[$ α_3 + α_4 $ main result components of $ φ π_q μ^{≥3} $:] Obvious.

\item[$ α_3 + α_4 $ main result component of $ φ π_q μ^2 $:] It necessarily concerns a result component of a product of the type $ φ π_q μ^2 (\id \text{(C)}, β \text{(A)}) $, $ φ π_q μ^2 (β \text{(A)}, \id \text{(C)}) $, $ φ π_q μ^2 (\id \text{(D)}, α_4) $ or $ φ π_q μ^2 (α_4, \id \text{(D)}) $. The first and second case are depicted in \autoref{fig:subdisk-alpha34-mult-main}. The third and fourth case are depicted in \autoref{fig:subdisk-degenerate}.

\item[$ α_3 + α_4 $ tail result component of a π-tree:] It comes either from the G components in \autoref{fig:components-final-out} or from the G components in \autoref{tab:components-multiplication}. In all cases, this concerns a tail component of a certain $ φ π_q (βα) $. Then the corresponding $ h_q (βα) $ indeed has a corresponding $ β $ (A) main result component.

\item[$ \id $ (C) main result component of $ φ π_q μ^2 $:] It comes from a product $ μ^2 (\id \text{(C)}, \id \text{(D)}) $ or $ μ^2 (\id \text{(D)}, \id \text{(C)}) $. Since this $ \id $ (D) is necessarily the result component of the tree in \autoref{fig:subdisk-idD-1} or \ref{fig:subdisk-idD-2}, we obtain the four $ \id $ (C) trees in \autoref{fig:subdisk-degenerate}.

\item[$ \id $ (C) main result component of $ φ π_q μ^{≥3} $:] It comes from an all-in disk where the first zigzag path turns right at the concluding arc and the final zigzag path turns left at the concluding arc. An example is depicted in \autoref{fig:subdisk-idC-disk-main}.

\item[$ \id $ (C) tail result component of a π-tree:] It comes either from the H components in \autoref{fig:components-final-out} or from the H components in \autoref{tab:components-multiplication}. In all cases, this concerns a tail component of a certain $ φ π_q (βα) $. Then the corresponding $ h_q (βα) $ indeed has a corresponding $ β $ (A) main result component.

\item[$ \id $ (D) result component of a π-tree:] It is either (a) the $ φ π_q μ^{≥3} $ of an all-in disk with equal first and final zigzag path, or (b) the result component of $ φ π_q μ^2 (\id \text{(C)}, \id \text{(B)}) $ or (c) $ φ π_q μ^2 (\id \text{(B)}, \id \text{(C)}) $ or (d) $ φ π_q μ^2 (\id \text{(D)}, \id \text{(D)}) $. Option (d) is impossible, since nonvanishing of $ φ π_q μ^2 (\id \text{(D)}, \id \text{(D)}) $ implies that both identities involved are the identity $ \id_{a_0} $ at the zigzag path's identity location, which is in contradiction to the fact that this involved $ \id $ (D) lies in $ R $. Options (a), (b) and (c) are possible and depicted in \autoref{fig:subdisk-idD}. Note the $ \id $ (B) component involved in (b) and (c) can impossibly come from $ μ^2 (β/β' \text{(C)}, \id \text{(C)}) $ or $ μ^2 (β \text{(A)}, α_0) $, because the arrow directions along the disk mismatch resp.~because the arrow direction of $ α_0 $ contradicts \autoref{conv:alpha0-direction}. This is reviewed in \autoref{fig:subdisk-idD-impossible}.

\item[$ α_0 $ result component of a π-tree:] A glance at the multiplication and disk tables reveals that it comes from a product $ π_q μ^2 (\id \text{(C)}, α_4) $ or $ π_q μ^2 (α_3, \id \text{(C)}) $ or $ π_q μ^2 (α_0, \id \text{(D)}) $ or $ π_q μ^2 (\id \text{(D)}, α_0) $. In the first case, both inputs are necessarily direct. In the second case, $ \id $ (C) is definitely direct. Meanwhile, $ α_3 $ may be direct or come from $ h_q μ^2 (α_0, α_4) $ with both $ α_0 $ and $ α_4 $ direct. In the third and fourth case, $ α_0 $ is direct. Meanwhile, $ \id $ (D) may come from $ h_q μ^2 (α_3, \id \text{(C)}) $ or $ h_q μ^2 (\id \text{(C)}, α_4) $ with both $ α_3/α_4 $ and $ \id $ (C) direct. This results in 7 options in total, depicted in \autoref{fig:subdisk-alpha0}.
\end{description}

We have checked all cases of \autoref{tab:resultcomp-classification}. This finishes the proof of \autoref{th:subdisk-assignment}.

\subsection{The shape of subdisks}
\label{sec:classification-subdiskshape}
In this section, we prove \autoref{th:subdisk-types-shape}: Subdisks of CR, ID, DS and DW result components are CR, ID, DS and DW disks, respectively.

\begin{lemma}
Subdisks of a CR result components are CR disks: $ \Subdisk(\CRr) ⊂ \CRd $.
\end{lemma}

\begin{proof}
Let $ r $ be a CR result component. We show that $ \Subdisk(r) $ is a CR disk. First, the corners of $ \Subdisk(r) $ are convex by construction. Second, since the subdisk is obtained from gluing smaller subdisks, this inductively provides that $ \Subdisk(r) $ is indeed the boundary of some immersed disk. Third, stacked co-identities lie infinitesimally close to each other as in Figure \ref{fig:subdisk-betaA-mult-main}, but all other input morphisms and the output morphism lie apart.

Regard a stack of co-identities on a zigzag path $ L $ used in $ r $. We show that $ L $ is oriented clockwise with $ \Subdisk(r) $. Recall from \autoref{conv:alpha0-direction} that co-identities lie in angles with puncture to the right of the zigzag curve in its natural orientation. The claim now follows from inspection of Figure \ref{fig:subdisk-betaC}, \ref{fig:subdisk-betaCp}, \ref{fig:subdisk-alpha0p}, \ref{fig:subdisk-betaA-mult-main}. In all cases, our convention implies that the zigzag curve is oriented clockwise with the subdisk.
\end{proof}

\begin{lemma}
Subdisks of ID result components are ID disks: $ \Subdisk(\IDr) ⊂ \IDd $.
\end{lemma}

\begin{proof}
First, a subdisk of Figure \ref{fig:subdisk-idD-idC-first} or \ref{fig:subdisk-idD-idC-last} has degenerate C input and first and final zigzag paths are oriented towards the interior of the disk. This makes an ID disk.

Second, a subdisk of Figure \ref{fig:subdisk-idD-all-in} has degenerate B input. Let us inspect the situation at the concluding arc of the all-in disk involved. Assume $ α_3 $ is the first angle of the disk. Then the degenerate input directly succeeds the output mark and no further input follows at infinitesimally small distance. The final angle of the disk can impossibly be $ α_3 $ or $ α_4 $ due to orientation of the concluding arc, therefore no input precedes the output at infinitesimally small distance. We see that the first, equivalently final zigzag path $ L_1 = L_{N+1} $ is oriented counterclockwise with the subdisk. We conclude that the subdisk is an ID disk. A similar conclusion holds in case $ α_4 $ is the final angle.
\end{proof}

Before we tackle DS and DW result components, let us recall the nasty results $ h_q (α \text{(D)}) $. Denote by $ S_α $ the sequence of arcs running from the source or target of $ α $ to the source or target of $ α_0 $, whichever are closer, without hitting $ a_0 $. The codifferential $ h_q (α) $ is then equal to the signed sum of these arc identities:
\begin{equation*}
h_q (α) = \sum_{c ∈ S_α} \id_c.
\end{equation*}
We are now ready to deal with DS result components.

\begin{lemma}
\label{th:classification-subdiskshape-DS}
Subdisks of DS result components are DS disks: $ \Subdisk(\DSr) ⊂ \DSd $.
\end{lemma}

\begin{proof}
Let $ T $ be one of the 8 trees in \autoref{fig:subdisk-degenerate} and $ r $ a result component. We need to check that its subdisk falls under the condition of DS disks. Indeed, it is bounded by two arcs $ a $ and $ b $ and lies between a zigzag curve $ \smooth L $ and its Hamiltonian deformation $ \smooth L' $. We only need to check two borderline conditions: The first condition is that $ a ≠ a_0 $ if $ \smooth L $ is oriented towards the co-identity. The second condition is that the strip has positive width if $ \smooth L $ is oriented away from the co-identity.

Assume $ r $ is a result component with $ a = a_0 $. We will show by inspection that $ S_α $ runs away from $ a_0 $ in oppposite direction of the orientation of $ \smooth L $. For this, consider the two cases that the inner product is $ μ^2 (α_3, \id \text{(C)}) $ or $ μ^2 (\id \text{(C)}, α_4) $ and distinguish further regarding arrow directions. The following graphics depicts all four cases, with the horizontal zigzag path being $ L $:

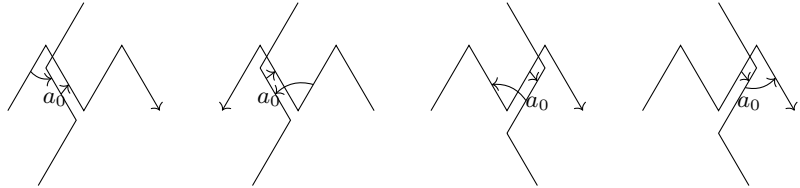
\begin{center}
\begin{tikzpicture}
\path[draw, ->] (0, 0) -- ++(60:1) coordinate (1) coordinate[pos=0.6] (alpha4-start) -- ++(300:1) coordinate (2) coordinate[pos=0.6] (idC-end) -- ++(60:1) -- ++(300:1);
\path[draw] (1)++(down:0.3) -- ++(60:1) (1)++(down:0.3) -- ++(300:0.8) node[pos=0.3, below] {\small $ a_0 $} coordinate[pos=0.2] (alpha4-end) coordinate[midway] (idC-start) -- ++(240:1);
\path[draw, ->, bend right=30] (alpha4-start) to (alpha4-end);
\path[draw, ->] (idC-start) to (idC-end);
\end{tikzpicture}
\hspace{0.5cm}
\begin{tikzpicture}
\path[draw, <-] (0, 0) -- ++(60:1) coordinate (1) -- ++(300:1) coordinate (2) coordinate[pos=0.4] (idC-end) -- ++(60:1) coordinate[pos=0.4] (alpha4-start) -- ++(300:1);
\path[draw] (1)++(down:0.3) -- ++(60:1) (1)++(down:0.3) -- ++(300:0.8) coordinate[pos=0.2] (idC-start) coordinate[pos=0.5] (alpha4-end) node[pos=0.3, below] {\small $ a_0 $} -- ++(240:1);
\path[draw, ->, bend right=30] (alpha4-start) to (alpha4-end);
\path[draw, ->] (idC-start) to (idC-end);
\end{tikzpicture}
\hspace{0.5cm}
\begin{tikzpicture}
\path[draw, ->] (0, 0) -- ++(60:1) -- ++(300:1) coordinate[pos=0.6] (alpha3-end) coordinate (1) -- ++(60:1) coordinate[pos=0.6] (idC-start) -- ++(300:1);
\path[draw] (1)++(down:0.3) -- ++(60:1) node[pos=0.7, below] {\small $ ~a_0 $} coordinate[pos=0.8] (idC-end) coordinate[pos=0.5] (alpha3-start) -- ++(120:1) (1)++(down:0.3) -- ++(300:0.8);
\path[draw, ->, bend right=30] (alpha3-start) to (alpha3-end);
\path[draw, ->] (idC-start) to (idC-end);
\end{tikzpicture}
\hspace{0.5cm}
\begin{tikzpicture}
\path[draw, ->] (0, 0) -- ++(60:1) -- ++(300:1) coordinate (1) -- ++(60:1) coordinate[pos=0.6] (idC-start) -- ++(300:1) coordinate[pos=0.5] (alpha3-end);
\path[draw] (1)++(down:0.3) -- ++(60:1) node[pos=0.7, below] {\small $ ~a_0 $} coordinate[pos=0.8] (idC-end) coordinate[pos=0.7] (alpha3-start) -- ++(120:1) (1)++(down:0.3) -- ++(300:0.8);
\path[draw, ->, bend right=30] (alpha3-start) to (alpha3-end);
\path[draw, ->] (idC-start) to (idC-end);
\end{tikzpicture}
\hspace{0.5cm}
\end{center}

In all four cases, we conclude that the arc sequence $ S_α $ runs away from $ a_0 $ against the orientation of $ \smooth L $. This means that $ \smooth L $ is oriented from $ b $ to $ a $, in other words: away from the co-identity. This proves the first condition. The second condition is checked similarly.
\end{proof}

\begin{lemma}
\label{th:classification-subdiskshape-DW}
Subdisks of DW result components are DW disks: $ \Subdisk(\DWr) ⊂ \DWd $.
\end{lemma}

\begin{proof}
This is similar to \autoref{th:classification-subdiskshape-DS}. Let us elaborate nevertheless: By definition, DW disks are a collection of three similar types of disks with $ α_0 $ output. By definition, DW result components are the $ α_0 $ result components of the 7 trees in \autoref{fig:subdisk-alpha0}. It is our task to check for every of these 7 trees that their subdisks fall under one of the three types of DW disks.

Of the 7 trees, the two trees without $ α_0 $ input fall under the triangle DW disk type. The two trees with an $ α_0 $ input at the beginning fall under the 4-gon DW disk type with $ α_0 $ succeeding the output mark. The two trees with an $ α_0 $ at the end, as well as the tree with the infinitesimally small subdisk, fall under the 4-gon DW type with $ α_0 $ preceding the output mark. The additional conditions are checked in the same way as for DS result components. This case distinction finishes the proof.
\end{proof}

\subsection{Narrow locations}
\label{sec:classification-narrow}
In this section, we begin proving \autoref{th:subdisk-minmodel-bijection} which states that all CR, ID, DS and DW disks lie in the image of $ \Subdisk $. Let us recall the situation: In \autoref{sec:classification-subdiskshape} we have already shown that $ \Subdisk $ maps only to CR, ID, DS and DW disks. Starting in the present section and ending in \autoref{sec:classification-IDDSDW}, we show that every CR, ID, DS and DW disk is actually reached by $ \Subdisk $.

\begin{figure}
\centering
\begin{tikzpicture}
\path (4, 0) node (CRDisk) {CR disk $ D $};
\path (0, -2) node[align=center] (Narrow) {Narrow locations \\ $ (l_1, m_1), (l_2, m_2) $};
\path (0.5, -5) node (NarrowTree) {Narrow tree};
\path (8, -2) node[align=center] (Subtree) {Subtree $ \Subtree(D) $ \\ Subresult $ \Subresult(D) $};
\path (7.5, -5) node (EvTree) {Evaluation tree};
\path (4, -5.3) node {\LARGE $ = $};
\begin{scope}[shift={($ (CRDisk.south) + (0, -0.5) $)}]
\path (-2, -2) coordinate (5) (2, -2) coordinate (1);
\path[draw, thick] ($ (5)!0.4!(1) $) to[out=90, in=270] ++(-0.2, 0.5) to[out=90, in=180] ($ (5)!0.5!(1) + (0, 1) $) to[out=0, in=90] ($ (5)!0.6!(1) + (0.2, 0.5) $) to[out=270, in=90] ($ (5)!0.6!(1) $);
\path[draw, thick] (5) to[out=0, in=180] coordinate[at end] (4) ++(1, 0) to[out=90, in=270] ($ (5)!0.4!(1) $);
\path[draw, thick] (1) to[out=210, in=60] coordinate[at end] (2) ++(-0.5, -0.3) to[out=180, in=0] coordinate[at end] (3) ++(-0.5, 0) to[out=120, in=270] ($ (5)!0.6!(1) $);
\path[draw, thick] (0, 0) to[out=200, in=90] (-1, -0.5) to[out=270, in=30] (-0.3, -1) to[out=210, in=90] (5);
\path[draw, thick] (0, 0) to[out=340, in=90] (1, -0.5) to[out=270, in=150] (0.3, -1) to[out=330, in=90] (1);
\foreach \i in {{0, 0}, 1, 2, 3, 4, 5} {\path[fill] (\i) circle[radius=0.05];};
\path (1) node[below] {1} (2) node[below] {2} (3) node[below] {3} (4) node[below] {4} (5) node[below] {5} (0, 0) node[above] {out};
\path[draw, gray] (-0.3, -1.8) -- (0, -0.8) -- (0.3, -1.8);
\path[draw, gray] (-0.8, -0.8) -- (-0.1, -0.8) -- (-0.8, -1.1);
\path[draw, gray] (0.8, -0.8) -- (0.1, -0.8) -- (0.8, -1.1);
\end{scope}
\begin{scope}[shift={($ (NarrowTree.east) + (0.5, 0) $)}]
\path (0, 0) node (A) {R}
node[below left of=A] {$ (l_2, m_2) $} edge (A) node[below right of=A] {$ (l_1, m_1) $} edge (A);
\end{scope}
\begin{scope}[shift={($ (Subtree.north east) + (0, 0.5) $)}]
\path (0, 0) node (5) {$ h_5 $} node[right of=5] (4) {$ h_4 $} node[right of=4] (3) {$ h_3 $} node[right of=3] (2) {$ h_2 $} node[right of=2] (1) {$ h_1 $}
node[below right of=5, align=center] (B) {$ β $ (A)} edge (5) edge (4)
node[below right of=B] {tail $ π_q μ^{≥3} = $ out} edge (B) edge (3) edge (2) edge (1);
\end{scope}
\begin{scope}[shift={($ (EvTree.west) + (-0.5, 0) $)}]
\path (0, 0) node (A) {Tail}
node[below left of=A] {Disk} edge (A) node[below right of=A] {Disk} edge (A);
\end{scope}
\path[draw] (CRDisk.south) -- (Narrow.north);
\path[draw] (Narrow.south) -- (NarrowTree.north);
\path[draw] (CRDisk.south) -- (Subtree.north);
\path[draw] (Subtree.south) -- (EvTree.north);
\end{tikzpicture}
\caption{Data structures for \autoref{sec:classification-narrow} till \ref{sec:classification-IDDSDW}}
\label{fig:classification-datastructures}
\end{figure}

Our strategy is to construct an explicit inverse map: In the present section, we analyze the shape of a given CR disk $ D ∈ \CRd $. In \autoref{sec:classification-narrowtree}, we build a candidate tree $ \Subtree(D) $. In \autoref{sec:classification-subresult}, we build a result component $ \Subresult(D) $ of $ \Subtree(D) $. In \autoref{sec:classification-verification}, we verify that its subdisk $ \Subdisk(\Subresult(D)) $ is equal to $ D $ again. We will finish our line of proof in \autoref{sec:classification-IDDSDW} by applying similar arguments to ID disks and checking the cases of DS and DW disks combinatorially. The essential data structures for the course of these sections are collected in \autoref{fig:classification-datastructures}.

Here is the observation that drives our strategy: Imagine $ r $ is a result component of a π-tree with subdisk $ D = \Subdisk(r) $. The subdisk catalog orders us to draw the subdisk stroke around all inputs of the $ μ^{≥3} $ disk and end up near its first/final arc, on both sides of the stroke. The subdisk becomes narrow there! The reader finds examples in \autoref{fig:subdisk-examples}.

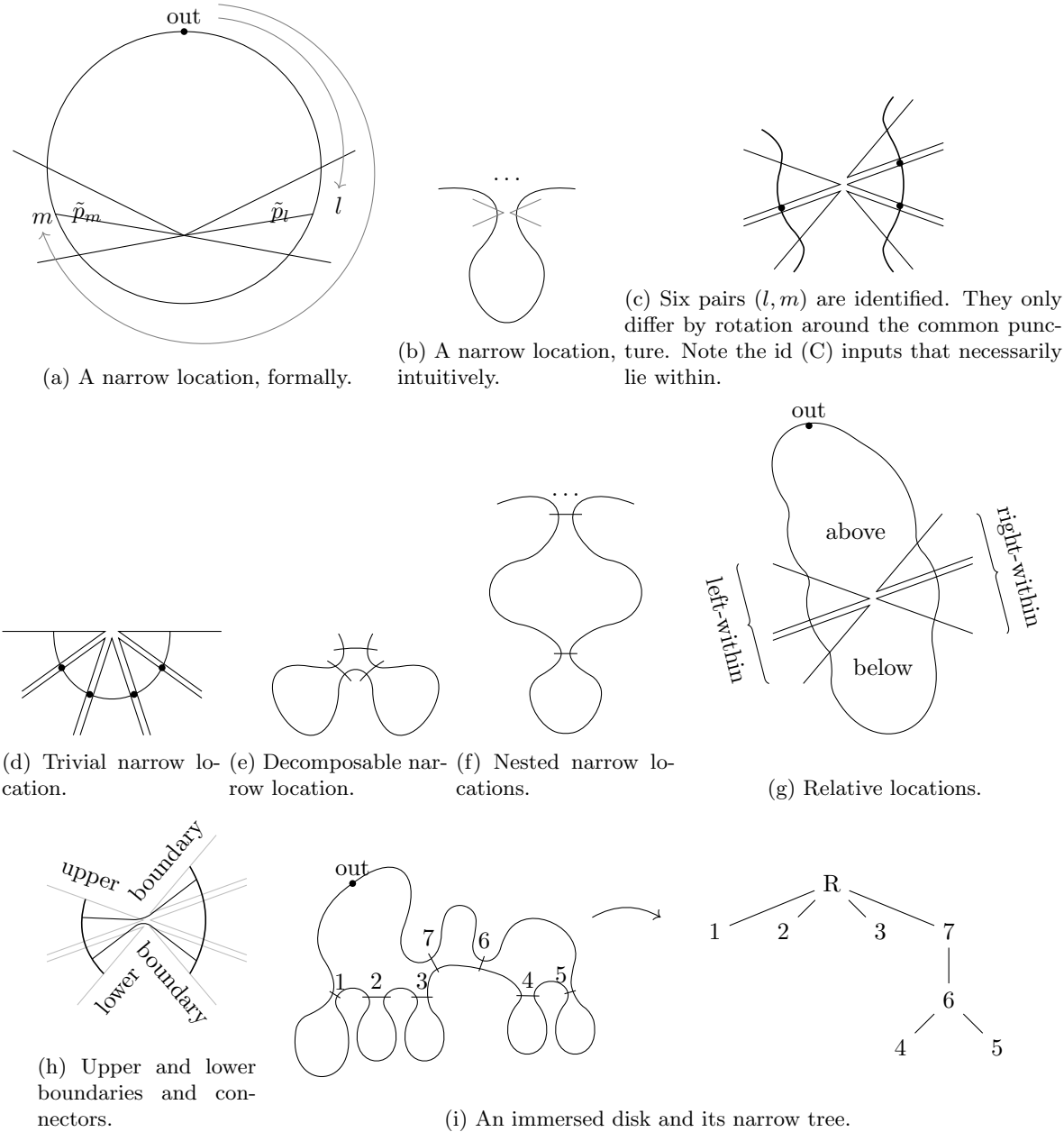
\begin{figure}
\centering
\begin{subfigure}[b]{0.35\linewidth}
\centering
\begin{tikzpicture}
\path[draw] (0, 0) circle[radius=2];
\path (320:2) coordinate (right-bottom) (340:2) coordinate (right-mid) (0:2) coordinate (right-top) (180:2) coordinate (left-top) (200:2) coordinate (left-mid) (220:2) coordinate (left-bottom);
\path[draw] (0, -1) -- ($ (0, -1)!1.25!(right-top) $);
\path[draw] (0, -1) -- ($ (0, -1)!1.4!(right-bottom) $);
\path[draw] (0, -1) -- ($ (0, -1)!1.25!(left-top) $);
\path[draw] (0, -1) -- ($ (0, -1)!1.4!(left-bottom) $);
\path[draw] (0, -1) -- ($ (right-mid) $) node[near end, shift={(0, 0.1)}] {$ \tilde p_l $};
\path[draw] (0, -1) -- ($ (left-mid) $) node[near end, shift={(0, 0.1)}] {$ \tilde p_m $};
\path[fill] (0, 2) circle[radius=0.05] node[above] {out};
\path[draw, ->, gray] (0.5, 2.2) arc(85:-15:2) coordinate[at end] (l-end);
\path (l-end) node[below] {$ l $};
\path[draw, ->, gray] (0.5, 2.4) arc(85:-160:2.5) coordinate[at end] (m-end);
\path (m-end) node[above] {$ m $};
\end{tikzpicture}
\caption{A narrow location, formally.}
\end{subfigure}
\begin{subfigure}[b]{0.2\linewidth}
\centering
\begin{tikzpicture}
\path[draw, rounded corners] (0, 0) to[bend left=70, looseness=2] (0.5, -1) to[bend right=40, looseness=1] (1, -2) to[bend right=40, looseness=1] (1.5, -1) to[bend left=70, looseness=2] (2, 0);
\path (1, 0) node[above] {…};
\path[draw, gray] (0.5, -0.15) -- (0.95, -0.35) -- (0.5, -0.55);
\path[draw, gray] (1.5, -0.15) -- (1.05, -0.35) -- (1.5, -0.55);
\end{tikzpicture}
\caption{A narrow location, intuitively.}
\end{subfigure}
\begin{subfigure}[b]{0.4\linewidth}
\centering
\begin{tikzpicture}
\path[draw] (0, 0) -- ++(340:1.5) coordinate (mid-left) coordinate[pos=0.4] (top-left) -- ++(200:1.5) ++(0, -0.1) -- ++(20:1.5) -- ++(230:1.5) coordinate[pos=0.6] (bottom-left);
\path[draw] (mid-left) ++(0.1, 0.1) -- ++(20:1.5);
\path[draw] (mid-left) ++(0.1, 0.1) -- ++(50:1.5) coordinate[pos=0.6] (top-right);
\path[draw] (mid-left) ++(0.1, 0) -- ++(20:1.5);
\path[draw] (mid-left) ++(0.1, 0) -- ++(340:1.5) ++(0, -0.1) -- ++(160:1.5) -- ++(310:1.5) coordinate[pos=0.6] (bottom-right);
\path[draw, semithick, rounded corners] ($ (bottom-left) + (-0.1, -0.5) $) to[bend right] (bottom-left) to[bend left] coordinate[pos=0.44] (id1) (top-left) to[bend right] ++(-0.3, 0.5);
\path[draw, semithick, rounded corners] (top-right) ++(0.1, 0.5) to[bend right] (top-right) to[bend left] coordinate[pos=0.3] (id2) coordinate[pos=0.7] (id3) (bottom-right) to[bend right] ($ (bottom-right) + (0.1, -0.5) $);
\path[fill] (id1) circle[radius=0.05] (id2) circle[radius=0.05] (id3) circle[radius=0.05];
\end{tikzpicture}
\caption{Six pairs $ (l, m) $ are identified. They only differ by rotation around the common puncture. Note the $ \id $ (C) inputs that necessarily lie within.}
\end{subfigure}
\begin{subfigure}[b]{0.2\linewidth}
\centering
\begin{tikzpicture}
\path[draw] (0, 0) -- ++(right:1.5) coordinate[pos=0.5] (left) -- ++(216:1.5) ++(0, -0.1) -- ++(36:1.5) -- ++(252:1.5) ++(0.1, 0) -- ++(72:1.5) -- ++(288:1.5) ++(0.1, 0) -- ++(108:1.5) -- ++(324:1.5) ++(0, 0.1) -- ++(144:1.5) -- ++(right:1.5) coordinate[pos=0.5] (right);
\path[draw, bend right=90, looseness=2] (left) to coordinate[pos=0.16] (id1) coordinate[pos=0.37] (id2) coordinate[pos=0.63] (id3) coordinate[pos=0.84] (id4) (right);
\path[fill] (id1) circle[radius=0.05];
\path[fill] (id2) circle[radius=0.05];
\path[fill] (id3) circle[radius=0.05];
\path[fill] (id4) circle[radius=0.05];
\end{tikzpicture}
\caption{Trivial narrow location.}
\end{subfigure}
\begin{subfigure}[b]{0.2\linewidth}
\centering
\begin{tikzpicture}
\path[use as bounding box] (-1.5, -1.5) -- (1.5, 0.5);
\path[draw, rounded corners] (0, 0) to[bend left=60, looseness=1.5] coordinate[pos=0.3] (top-1) coordinate[pos=0.75] (left-1) ++(-0.3, -0.5) to[bend right=130, looseness=6] ++(0.3, -0.5) to[bend left=90, looseness=3] coordinate[pos=0.3] (left-2) coordinate[pos=0.7] (right-2) ++(0.5, 0) to[bend right=130, looseness=6] ++(0.3, 0.5) to[bend left=60, looseness=1.5] coordinate[pos=0.25] (right-1) coordinate[pos=0.7] (top-2) ++(-0.3, 0.5);
\path[draw, bend left=9] ($ (top-1)!-0.4!(top-2) $) to ($ (top-2)!-0.4!(top-1) $);
\path[draw, bend left=9] ($ (left-1)!-0.4!(left-2) $) to ($ (left-2)!-0.4!(left-1) $);
\path[draw, bend right=9] ($ (right-1)!-0.4!(right-2) $) to ($ (right-2)!-0.4!(right-1) $);
\end{tikzpicture}
\caption{Decomposable narrow location.}
\end{subfigure}
\begin{subfigure}[b]{0.2\linewidth}
\centering
\begin{tikzpicture}
\path[draw, rounded corners] (0, 0) to[bend left=80, looseness=2.5] coordinate[pos=0.6] (top-left) ++(0.5, -0.8) to[bend right=80, looseness=2] ++(0, -1) to[bend left=60, looseness=2] coordinate[pos=0.5] (bottom-left) ++(0, -0.8) to[bend right=40, looseness=1] ++(0.5, -0.8) to[bend right=40, looseness=1] ++(0.5, 0.8) to[bend left=60, looseness=2] coordinate[pos=0.5] (bottom-right) ++(0, 0.8) to[bend right=80, looseness=2] ++(0, 1) to[bend left=80, looseness=2.5] coordinate[pos=0.4] (top-right) ++(0.5, 0.8);
\path (1, 0) node[above] {…};
\path[draw] ($ (top-left)!-0.4!(top-right) $) -- ($ (top-right)!-0.4!(top-left) $);
\path[draw] ($ (bottom-left)!-0.4!(bottom-right) $) -- ($ (bottom-right)!-0.4!(bottom-left) $);
\end{tikzpicture}
\caption{Nested narrow locations.}
\end{subfigure}
\begin{subfigure}[b]{0.35\linewidth}
\centering
\begin{tikzpicture}
\path[draw] (0, 0) coordinate (brace-left-top) -- ++(340:1.5) coordinate (mid-left) coordinate[pos=0.4] (top-left) -- ++(200:1.5) ++(0, -0.1) -- ++(20:1.5) -- ++(230:1.5) coordinate[pos=0.6] (bottom-left) coordinate (brace-left-bottom);
\path[draw] (mid-left) ++(0.1, 0.1) -- ++(20:1.5);
\path[draw] (mid-left) ++(0.1, 0.1) -- ++(50:1.5) coordinate[pos=0.6] (top-right) coordinate (brace-right-top);
\path[draw] (mid-left) ++(0.1, 0) -- ++(20:1.5);
\path[draw] (mid-left) ++(0.1, 0) -- ++(340:1.5) coordinate[pos=0.6] (bottom-right) coordinate (brace-right-bottom);
\path[draw, rounded corners] (bottom-left) to[bend left] (top-left) to[bend left] ++(-0.3, 1) to[bend left=80, looseness=2] coordinate[pos=0.7] (out) ++(1, 1) to[bend left] (top-right) to[bend left] (bottom-right) to[bend right=10] ++(0, -0.5) to[bend left=80, looseness=2] coordinate[pos=0.5] (bottom) ($ (bottom-left) + (0, -0.5) $) to[bend right=10] (bottom-left);
\path[fill] (out) circle[radius=0.05] node[above] {out};
\path[draw, decoration={brace}, decorate] ($ (brace-left-bottom) + (-0.5, 0) $) to node[sloped, shift={(0, -0.4)}] {left-within} ($ (brace-left-top) + (-0.5, 0) $);
\path[draw, decoration={brace}, decorate] ($ (brace-right-top) + (0.5, 0) $) to node[sloped, shift={(0, 0.4)}] {right-within} ($ (brace-right-bottom) + (0.5, 0) $);
\path (mid-left) ++(-0.2, 1) node {above};
\path (mid-left) ++(0.2, -1) node {below};
\end{tikzpicture}
\caption{Relative locations.}
\end{subfigure}
\begin{subfigure}[b]{0.2\linewidth}
\centering
\begin{tikzpicture}
\path[draw, gray!50] (0, 0) coordinate (left-1) -- ++(340:1.5) node[midway, sloped, black, shift={(-0.2, 0.3)}] {upper} coordinate (mid-left) coordinate[pos=0.4] (top-left) -- ++(200:1.5) coordinate (left-2) ++(0, -0.1) coordinate (left-3) -- ++(20:1.5) -- ++(230:1.5) coordinate[pos=0.6] (bottom-left) node[midway, sloped, black, shift={(-0.2, -0.3)}] {lower} coordinate (left-4);
\path[draw, gray!50] (mid-left) ++(0.1, 0.1) -- ++(20:1.5) coordinate (right-2);
\path[draw, gray!50] (mid-left) ++(0.1, 0.1) -- ++(50:1.5) coordinate[pos=0.6] (top-right) node[midway, sloped, black, shift={(0, 0.3)}] {boundary} coordinate (right-1);
\path[draw, gray!50] (mid-left) ++(0.1, 0) -- ++(20:1.5);
\path[draw, gray!50] (mid-left) ++(0.1, 0) -- ++(340:1.5) ++(0, -0.1) coordinate (right-3) -- ++(160:1.5) -- ++(310:1.5) coordinate[pos=0.6] (bottom-right) node[midway, sloped, black, shift={(0.1, -0.3)}] {boundary} coordinate (right-4);
\path[draw, semithick, rounded corners] (bottom-left) to[bend left] coordinate[pos=0.22] (left-6) coordinate[pos=0.78] (left-5) (top-left);
\path[draw, semithick, rounded corners] (top-right) to[bend left] coordinate[pos=0.12] (right-5) coordinate[pos=0.88] (right-6) (bottom-right);
\path[draw, rounded corners] (left-5) -- (mid-left) -- (right-5);
\path[draw, rounded corners] (left-6) -- (mid-left) -- (right-6);
\end{tikzpicture}
\caption{Upper and lower boundaries and connectors.}
\end{subfigure}
\begin{subfigure}[b]{0.7\linewidth}
\centering
\begin{tikzpicture}
\path[draw, rounded corners] (0, 0) to[bend right=30] ++(-0.5, -1) to[bend left=50, looseness=1.5] coordinate[pos=0.5] (1-l) ++(-0.3, -1) to[bend right=100, looseness=4] ++(0.7, 0) to[bend left=50, looseness=1.5] coordinate[pos=0.6] (1-r) ++(0.1, 0.5) to[bend left=50, looseness=1.5] coordinate[pos=0.4] (2-l) ++(0.1, -0.5) to[bend right=100, looseness=4] ++(0.5, 0) to[bend left=50, looseness=1.5] coordinate[pos=0.6] (2-r) ++(0.1, 0.5) to[bend left=50, looseness=1.5] coordinate[pos=0.4] (3-l) ++(0.1, -0.5) to[bend right=100, looseness=4] ++(0.5, 0) to[bend left=50, looseness=1] coordinate[pos=0.42] (3-r) coordinate[pos=0.9] (7-l) ++(0, 0.8) to[bend left=10, looseness=1] coordinate[pos=0.55] (6-l) ++(1, -0.2) to[bend left=40, looseness=1.5] coordinate[pos=0.5] (4-l) ++(0, -0.5) to[bend right=100, looseness=4] ++(0.5, 0) to[bend left=50, looseness=1.5] coordinate[pos=0.42] (4-r) ++(0.1, 0.5) to[bend left=50, looseness=1.5] coordinate[pos=0.5] (5-l) ++(0.1, -0.5) to[bend right=100, looseness=4] ++(0.5, 0) to[bend left=30, looseness=1] coordinate[pos=0.6] (5-r) ++(-0.3, 0.6) to[bend right=80, looseness=1.5] ++(-1, 0.5) to[bend left=90, looseness=3] coordinate[pos=0.5] (6-r) ++(-0.4, 0.1) to[bend right=80, looseness=2] ++(-0.5, 0.1)
to[bend left=100, looseness=4] coordinate[pos=0.3] (7-r) ++(-0.5, 0.1) to[bend right=70, looseness=2] (0, 0);
\path[fill] (0, 0) circle[radius=0.05] node[above] {out};
\path[draw] ($ (1-l)!-0.2!(1-r) $) -- ($ (1-r)!-0.2!(1-l) $) node[above] {1};
\path[draw] ($ (2-l)!-0.2!(2-r) $) -- ($ (2-r)!-0.2!(2-l) $) node[midway, above] {2};
\path[draw] ($ (3-l)!-0.2!(3-r) $) -- ($ (3-r)!-0.2!(3-l) $) node[midway, above] {3};
\path[draw] ($ (4-l)!-0.2!(4-r) $) -- ($ (4-r)!-0.2!(4-l) $) node[midway, above] {4};
\path[draw] ($ (5-l)!-0.2!(5-r) $) -- ($ (5-r)!-0.2!(5-l) $) node[above left] {5};
\path[draw] ($ (6-l)!-0.2!(6-r) $) -- ($ (6-r)!-0.2!(6-l) $) node[above] {6};
\path[draw] ($ (7-l)!-0.2!(7-r) $) -- ($ (7-r)!-0.2!(7-l) $) node[above] {7};
\path[draw, bend left, ->] (3.5, -0.5) to (4.5, -0.5);
\path (7, 0) node (R) {R}
node[below left of=R] (2) {2} edge(R) node[left of=2] {1} edge(R) node[below right of=R] (3) {3} edge(R) node[right of=3] (7) {7} edge (R)
node[below of=7] (6) {6} edge(7)
node[below left of=6] {4} edge(6) node[below right of=6] {5} edge(6);
\end{tikzpicture}
\caption{An immersed disk and its narrow tree.}
\label{fig:classification-narrow-tree}
\end{subfigure}
\caption{Illustrations of narrow locations and their terminology}
\label{fig:classification-narrow}
\end{figure}

Conversely, assume a CR disk $ D $ is given without further knowledge. In order to guess a tree $ T $, we simply need to record all the narrow locations of $ D $ We are ready for a precise definition of narrow locations. The reader may already have a glance at \autoref{fig:classification-narrow}, where all upcoming notions are depicted.

\begin{definition}
Let $ D $ be a CR disk. Index the angles that the boundary of $ D $ cuts, in clockwise order. A \emph{narrow location} of $ D $ consists of two indices $ m > l $ on the boundary such that:
\begin{itemize}
\item Both $ m $ and $ l $ lie in angles whose centers (which are punctures) lie on the inside of the disk.
\item Let $ p_m $ be the path connecting $ m $ to its puncture and $ p_l $ the path connecting $ l $ to its puncture. Both lift to paths $ \tilde p_m $ and $ \tilde p_l $ in the unit disk model. Then require that $ \tilde p_m $ and $ \tilde p_l $ actually meet.
\end{itemize}
In particular both punctures are equal. This is the \emph{connecting puncture}. The union of $ \tilde p_m $ and $ \tilde p_l $ is the \emph{connector} of $ (l, m) $. Identify two narrow locations if they only differ by rotation around the respective punctures (which are co-identities or situation C morphisms around these punctures). A narrow location $ (l, m) $ is \emph{trivial} if $ l $ and $ m $ only differ by rotation around their connecting puncture. A narrow location $ (l, m) $ is \emph{indecomposable} if it is nontrivial and there does not exist an index $ l < n < m $ such that both $ (l, n) $ and $ (n, m) $ are nontrivial narrow locations. Two narrow locations $ (l, m) $ and $ (l', m') $ are \emph{disjoint} if $ l' ≥ m $ or $ l ≥ m' $ (up to rotation around the connecting punctures). A \emph{decomposition} of a narrow location $ (l, m) $ into indecomposables consists of indecomposable disjoint narrow locations whose union is $ (l, m) $.
\end{definition}

\begin{lemma}
Let $ D $ be a CR disk. Then any two indecomposable narrow locations $ (l, m) $ and $ (l', m') $ are either nested or disjoint. Any narrow location decomposes uniquely into indecomposables.
\end{lemma}

\begin{proof}
To prove the first claim, assume $ (l, m) $ and $ (l', m') $ are neither nested nor disjoint. Without loss of generality, we have $ l < l' < m < m' $. Looking at the unit disk model, the connector of $ (l', m') $ then has to pass through the lift of the connecting puncture of $ (l, m) $. In particular the lifts of the connecting punctures of $ (l, m) $ and $ (l', m') $ are actually equal. Moreover, since the paths $ \tilde p_l $, $ \tilde p_{l'} $, $ \tilde p_m $, $ \tilde p_{m'} $ now all meet, we have that all of $ (l, l') $, $ (l', m) $, $ (m, m') $ are actually narrow locations, contradicting indecomposability.

For the second claim, existence of a decomposition is clear (one keeps decomposing until the components are indecomposable). Uniqueness follows from the first claim.
\end{proof}

\subsection{Narrow trees}
\label{sec:classification-narrowtree}
In this section, we introduce narrow trees. The idea is to capture all narrow locations of a given CR disk in a structured way. Since we have already seen that narrow locations are nested or disjoint, the most natural structure to capture them is a tree.

\begin{definition}
Let $ D $ be a CR disk with boundary of length $ |D| $. Then its \emph{narrow tree} is the ordered decorated tree defined as follows:
\begin{itemize}
\item The nodes are all indecomposable narrow locations.
\item The nodes are connected according to inclusion.
\item The nodes are ordered horizontally from high $ (l, m) $ to low $ (l, m) $.
\item The nodes are decorated with their narrow location $ (l, m) $.
\item Except in the case where $ (1, |D|) $ is a narrow location, insert a root standing for the artificial narrow location $ (1, |D|) $. The root's children are the maximal indecomposable narrow locations.
\end{itemize}
\end{definition}

An schematic example of a disk and its narrow tree is shown in \autoref{fig:classification-narrow-tree}. In this example, the root is named R. Note that the decorations $ (l, m) $ have been ignored in the narrow tree. In fact, it is not even possible to given concrete numbers for the decorations $ (l, m) ∈ ℕ × ℕ $, since we have only drawn the disk schematically.

\begin{definition}
Let $ D $ be a CR disk and $ (l, m) $ a narrow location (other than the artificial narrow location). Note that $ (l, m) $ contains multiple identified pairs $ (l', m') $. The connector of minimal $ l' $ and maximal $ m' $ is the \emph{upper connector} and the connector of maximal $ l' $ and minimal $ m' $ the \emph{lower connector}.

Let \emph{above} $ (l, m) $ refer to the portion of $ D $ minus the disk bounded by the upper connector and the corresponding $ D $ boundary segment, \emph{within} $ (l, m) $ refer to the portion of $ D $ between upper and lower connector, and \emph{below} $ (l, m) $ refer to the portion of $ D $ bounded by the lower connector and the corresponding $ D $ boundary segment. The segment of the boundary of $ D $ within $ (l, m) $ splits into two components, which are the \emph{left-within} and \emph{right-within} the narrow location. The \emph{upper boundary} of $ (l, m) $ consists of the two situation A arcs just above $ (l, m) $ and the \emph{lower boundary} of $ (l, m) $ consists of the two situation A arcs just below $ (l, m) $.

Now let $ (l, m) $ be an indecomposable narrow location of an immersed disk $ D $. The \emph{stray morphisms} are the input morphisms that are below $ (l, m) $, but above all of its children. A (direct) \emph{left (resp.~right) sibling} consists of a sibling $ (l', m') $ of $ (l, m) $ in the narrow tree such that $ l' = m $ up to rotation around the connecting puncture (resp.~$ m' = l $). A (direct) \emph{sibling} is a direct left or right sibling.

A \emph{stack} in $ (l, m) $ is one of the following: (a) a combination of a stray morphism $ α_3 + α_4 $ directly followed by an $ α_0 $ within the next angle being cut, or (b) the combination of multiple $ α_0 $ and $ \id $ (C) differing only by rotation around a puncture, or (c) a child together with all siblings and their within morphisms, or (d) any stray morphism that is not part of one of these combinations. The narrow location is \emph{2-rich} if it has at least two stacks. It is \emph{1-rich} if it has precisely one stack, and if this stack contains a child then it has a morphism within
or
a sibling. The narrow location is \emph{0-rich} if it has precisely one stack, and this stack is a child without morphisms within and without direct sibling.
\end{definition}

\begin{remark}
Due to zigzag consistency, any narrow location has at least one child or at least one stray morphism. This means that an indecomposable narrow location $ (l, m) $ is either 2-rich, 1-rich or 0-rich.
\end{remark}

\subsection{Subresults}
\label{sec:classification-subresult}
In this section, we build a result component from any given CR disk. More precisely, we associate to every CR disk $ D $ a π-tree $ \Subtree(D) $ with a result component $ \Subresult(D) $, in the hope that $ \Subdisk(\Subresult(D)) $ equals $ D $ again. We will call $ \Subtree(D) $ the subtree and $ \Subresult(D) $ the subresult associated with $ D $.

The strategy of constructing $ \Subtree(D) $ is to take the narrow tree as a starting point and keep inserting $ h_q μ^2 $ nodes to bind together morphisms that lie directly next to each other. We also have to put special attention to the relation of narrow locations: Siblings have to inserted in a specific order, sometimes irregularly.

Let us start from regarding a narrow location $ (l, m) $. Whenever $ (l, m) $ has a left direct sibling or a morphism left-within, we can easily make a final-out $ μ^{≥3} $ disk of $ (l, m) $, whose final morphism stems from the direct sibling and the further morphisms left-within $ (l, m) $. Indeed, $ (l, m) $ has at least one child or a stray morphism, therefore this really yields a $ μ_{\Tw\Gtl_q Q}^{≥2} $.

However in case $ (l, m) $ has no left direct sibling and no morphism left-within, we have to distinguish whether $ (l, m) $ is 2-rich, 1-rich or 0-rich. In the 2-rich case, we can proceed with an ordinary $ μ^{≥3} $ and obtain a main result component. In the 1-rich case, we can proceed with $ μ^2 $ and obtain a first-order tail result component. In the 0-rich case, we interpret $ (l, m) $ as a tail node of the first 1-rich or 2-rich narrow location we arrive at when tracing the tree from $ (l, m) $ towards the leaves.

This is our basic recipe of turning narrow locations into trees. Let us record a lemma affirming that this construction works.

\begin{lemma}
\label{th:subtree-connecting}
Let $ D $ be a CR disk and let $ (l, m) $ be either:
\begin{itemize}
\item an indecomposable narrow location. Then let $ E $ be the sequence of zigzag segments below $ (l, m) $ and above all children, together with the lower boundary of $ (l, m) $ and the upper boundary of all children.
\item the root $ (1, |D|) $ of the narrow tree of $ D $, and assume it is of type $ \id $ (C). Then let $ E $ denote the sequence of zigzag segments starting at the output mark, staying above the children of $ (1, |D|) $, and ending at the output mark, including the 2/5 arc at the output mark.
\item the root $ (1, |D|) $ of the narrow tree of $ D $, and assume it is of type $ α_3 + α_4 $. Then let $ E $ denote the sequence of zigzag segments starting at the output mark, staying above the children of $ (1, |D|) $, and ending at the output mark, excluding the 2/5 arc at the output mark.
\end{itemize}
Then $ E $ bounds a discrete immersed disk. Upon reversing Figure \ref{fig:subdisk-direct} and the stack figures thereafter, the sequence of stacks injects into the sequence of interior boundary angles of $ E $. The complement of the image consists of $ δ $ insertions and the $ β $ (A) morphisms at the children.

In exception to this assignment, a stack directly after or before a child (or the output mark in case 1) shall map to the entire corresponding $ β $ (A) morphism in $ F $, not only the $ β $ (A) morphism surrounding the left or right part. This means that if there is both a stack directly after and a stack directly before the output mark (or the output mark in case 1), they map in particular to the same $ β $ (A), in exception to injectivity.
\end{lemma}

\begin{proof}
Note that all $ D $ boundary segments involved can be split into small pieces cutting through indecomposable angles. We shall now construct the discrete immersed disk $ F $ as follows. At each of these angles whose puncture lies outside $ D $, flow $ E $ outwards to the puncture.

At each of the angles whose puncture lies inside $ D $, flow $ E $ inwards to the puncture. At all angles whose puncture lies outside $ D $, this procedure enlarges the disk and in particular keeps it immersed. If $ F $ becomes non-immersed as discrete disk, then this is due to two angles whose punctures lie inside $ D $ and meet. In other words, loss of immersedness constitutes a narrow location $ (l', m') $ of $ D $.

Let us check the possible locations of $ (l', m') $. Note that $ (l', m') $ cannot equal $ (1, |D|) $. Indeed, in case $ (1, |D|) $ is a narrow location at all, it is nontrivial by assumption, its lower boundary consists of merely two arcs and it does not pose an obstruction to immersedness at all. Similarly, deduce that $ (l', m') $ is not a trivial narrow location.

Since $ F $ cuts all children away and children are indecomposable, any child is either contained in $ (l', m') $ or disjoint. Splitting $ (l', m') $ into indecomposable narrow locations then necessarily yields a chain of direct siblings, which are children of $ (l, m) $. By construction of $ F $, such a chain of direct siblings children does not constitute an obstruction to immersedness at all.

The second part of the statement consists of generically checking whether any two consecutive stacks occurring along the boundary of $ D $ may fall into the same interior boundary angle. For all subdisks in Figure \ref{fig:subdisk-connecting}, this is definitely not the case. We shall therefore check that Figure \ref{fig:subdisk-connecting} actually displays all possible consecutive input stacks. We will illustrate this in case of an $ α_3 + α_4 $ input, whose 2/5 arc is oriented counterclockwise with the disk, and an arbitrary successor and predecessor.

Since $ D $ has convex corners, it cuts the two angles $ α_2 $ and $ α_3 $ before respectively after $ α_3 + α_4 $. By arrow directions, possible $ α_0 $ inputs may occur both on $ α_2 $ and $ α_3 $. Since we assumed the co-identity rule, there is at most one $ α_0 $ on both angles. By construction, we assign to the combination of $ α_3 $ and possibly one $ α_0 $ the corresponding $ α_3 $ interior boundary morphism. Regardless of its precise nature, the predecessor stack will definitely not map to $ α_3 $ at the same time.

It remains to check the the successor stack. If the successor lies directly at the target of $ α_3 $, then due to arrow direction that morphism produces the next $ α_3 $ type interior boundary angle indeed lying one arc apart. If the successor lies farther apart, it does not neighbor with $ α_3 $ anymore and will not yield this as interior boundary angle, which is also in line with injectivity. The three combinations that do occur are those in Figures \ref{fig:subdisk-alpha3-delta}, \ref{fig:subdisk-alpha3-betaA}, \ref{fig:subdisk-alpha3-betaC}.
\end{proof}

We are now ready for the complete definition of subresults.

\begin{definition}
Let $ D $ be a CR disk. Then its \emph{subresult} is the result component $ \Subresult(D) ∈ \PiTr $ on the \emph{subtree} $ \Subtree(D) $ constructed below.
\end{definition}

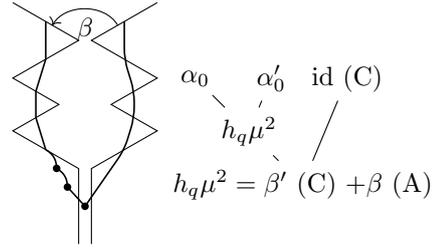
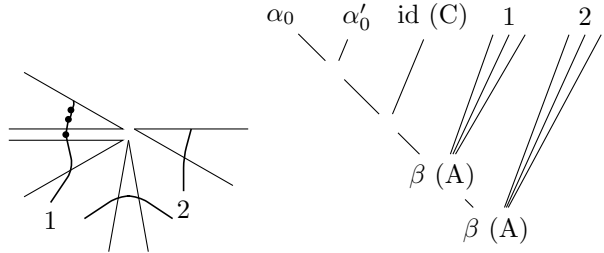
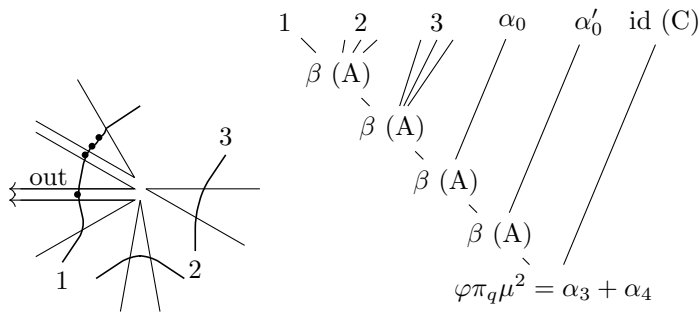
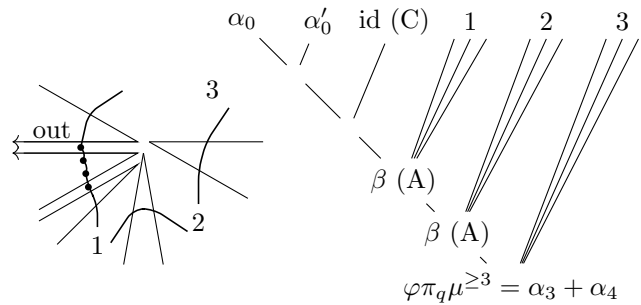
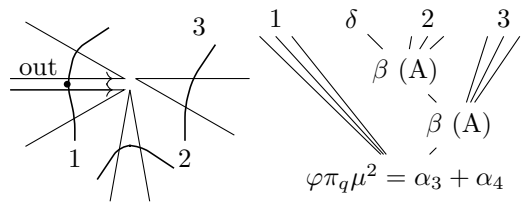
\begin{figure}
\centering
\begin{subfigure}[b]{0.45\linewidth}
\centering
\begin{tikzpicture}
\path[draw] (0, 0) -- ++(330:1) coordinate[midway] (1) coordinate[pos=0.6] (beta-end) -- ++(210:1) coordinate[midway] (2) -- ++(330:0.7) coordinate[midway] (3) -- ++(210:0.7) coordinate[midway] (4) -- ++(330:1) coordinate[midway] (5) -- ++(210:0.7) coordinate[midway] (6);
\path[draw] (1.9, 0) -- ++(210:1) coordinate[midway] (7) coordinate[pos=0.6] (beta-start) -- ++(330:1) coordinate[midway] (8) -- ++(210:0.7) coordinate[midway] (9) -- ++(330:0.7) coordinate[midway] (10) -- ++(210:1) coordinate[midway] (11) -- ++(330:0.7) coordinate[midway] (12);
\path[draw, rounded corners, semithick] (6) -- (5) -- (4) -- (3) -- (2) -- (1);
\path[draw, rounded corners, semithick] (7) -- (8) -- (9) -- (10) -- (11) -- (12) node[below left] {…};
\path[draw, ->, bend right=60] (beta-start) to node[midway, below] {$ β $} (beta-end);
\path (3, -1) node[align=center] {(tail result \\ component)};
\end{tikzpicture}
\caption{A non-root node of the narrow tree. This example has no siblings, no morphisms left-within, and is 0-rich. Whatever is below the dots, $ β $ becomes a tail result component.}
\label{fig:subtree-nonroot-0-rich}
\end{subfigure}
\hspace{0.05\linewidth}
\begin{subfigure}[b]{0.45\linewidth}
\centering
\begin{tikzpicture}
\path[draw] (0, 0) -- ++(330:1) coordinate[midway] (1) coordinate[pos=0.6] (beta-end) -- ++(210:1) coordinate[midway] (2) -- ++(330:0.7) coordinate[midway] (3) -- ++(210:0.7) coordinate[midway] (4) -- ++(330:1) coordinate[midway] (5) -- ++(down:1) coordinate[midway] (6);
\path[draw] (1.9, 0) -- ++(210:1) coordinate[midway] (7) coordinate[pos=0.6] (beta-start) -- ++(330:1) coordinate[midway] (8) -- ++(210:0.7) coordinate[midway] (9) -- ++(330:0.7) coordinate[midway] (10) -- ++(210:1) coordinate[midway] (11) -- ++(down:1) coordinate[midway] (12);
\path ($ (5)!0.33!(6) $) coordinate (switch1) ($ (5)!0.66!(6) $) coordinate (switch2);
\path[draw, semithick] ($ (6)!0.5!(12) $) -- (switch2);
\path[draw, semithick, bend right] (switch2) to (switch1);
\path[draw, semithick, bend right] (switch1) to (5);
\path[draw, semithick] (5) -- (4) -- (3) -- (2) -- (1);
\path[draw, rounded corners, semithick] (7) -- (8) -- (9) -- (10) -- (11) -- ($ (12)!0.5!(6) $);
\path[draw, ->, bend right=60] (beta-start) to node[midway, below] {$ β $} (beta-end);
\path[fill] ($ (6)!0.5!(12) $) circle[radius=0.05];
\path[fill] (switch1) circle[radius=0.05];
\path[fill] (switch2) circle[radius=0.05];
\path (2.4, -1) node (alpha0) {$ α_0 $} node[right of=alpha0] (alpha0p) {$ α_0' $} node[right of=alpha0p] (idC) {$ \id $ (C)}
node[below right of=alpha0] (m1) {$ h_q μ^2 $} edge (alpha0) edge (alpha0p)
node[below right of=m1] {$ h_q μ^2 = β' $ (C) $ + β $ (A)} edge (m1) edge (idC);
\end{tikzpicture}
\caption{A non-root node of the narrow tree. This example has no siblings, no morphisms left-within, and is 1-rich.  $ β $ becomes a tail result component of the stack of the 3 indicated inputs.}
\label{fig:subtree-nonroot-1-rich}
\end{subfigure}
\begin{subfigure}[b]{\linewidth}
\centering
\begin{tikzpicture}[scale=1.5]
\path[draw] (0, 0) -- ++(330:1) coordinate (mid) coordinate[midway] (1) -- ++(left:1) coordinate[midway] (2) ++(down:0.1) -- ++(right:1) coordinate[midway] (3) -- ++(210:1) coordinate[midway] (4);
\path[draw] (mid) ++(0.05, -0.1) -- ++(260:1) coordinate[midway] (5);
\path[draw] (mid) ++(0.05, -0.1) -- ++(280:1) coordinate[midway] (6);
\path[draw] (mid) ++(0.1, 0) -- ++(right:1) coordinate[midway] (7);
\path[draw] (mid) ++(0.1, 0) -- ++(330:1) coordinate[midway] (8);
\path (1) -- (2) coordinate[pos=0.33] (switch2) coordinate[pos=0.66] (switch1);
\path[draw, semithick, rounded corners] (7) -- (8) -- ($ (8) + (0, -0.3) $) node[below] {2};
\path[draw, semithick, rounded corners] (5) ++(-0.3, -0.2) -- (5) -- (6) -- ($ (6) + (0.3, -0.2) $);
\path[draw, semithick, rounded corners] ($ (4) + (-0.2, -0.3) $) node[below] {1} -- (4) -- ($ (2)!0.5!(3) $);
\path[draw, semithick] ($ (2)!0.5!(3) $) -- (switch1);
\path[draw, semithick, bend right] (switch1) to (switch2);
\path[draw, semithick, bend right] (switch2) to (1);
\path[fill] (switch1) circle[radius=0.03];
\path[fill] (switch2) circle[radius=0.03];
\path[fill] ($ (2)!0.5!(3) $) circle[radius=0.03];
\end{tikzpicture}
\begin{tikzpicture}
\path (2.4, 0.5) node (alpha0) {$ α_0 $} node[right of=alpha0] (alpha0p) {$ α_0' $} node[right of=alpha0p] (idC) {$ \id $ (C)} node[right of=idC] (i1) {1} node[right of=i1] (i2) {2}
node[below right of=alpha0] (m1) {} edge (alpha0) edge (alpha0p)
node[below right of=m1] (m2) {} edge (m1) edge (idC)
node[below right of=m2] (m3) {$ β $ (A)} edge (m2) edge (i1.south west) edge (i1.south) edge (i1.south east)
node[below right of=m3] {$ β $ (A)} edge (m3) edge (i2.south west) edge (i2.south) edge (i2.south east);
\end{tikzpicture}
\caption{A non-root node of the narrow tree. This example has siblings, and those are bound left to right.}
\label{fig:subtree-nonroot-2-rich}
\end{subfigure}
\begin{subfigure}[b]{\linewidth}
\centering
\begin{tikzpicture}[scale=1.5]
\path[draw] (0, 0) -- ++(330:1) coordinate (mid) coordinate[midway] (1) -- ++(left:1) coordinate[midway] (2) ++(down:0.1) -- ++(right:1) coordinate (mid-low) coordinate[midway] (3) -- ++(210:1) coordinate[midway] (4);
\path[draw] (mid) ++(0, 0.1) -- ++(150:1) coordinate[midway] (0) (mid) ++(0, 0.1) -- ++(120:1) coordinate[midway] (m1);
\path[draw] (mid) ++(0.05, -0.1) -- ++(260:1) coordinate[midway] (5);
\path[draw] (mid) ++(0.05, -0.1) -- ++(280:1) coordinate[midway] (6);
\path[draw] (mid) ++(0.1, 0) -- ++(right:1) coordinate[midway] (7);
\path[draw] (mid) ++(0.1, 0) -- ++(330:1) coordinate[midway] (8);
\path[draw, semithick, rounded corners] (7) ++(0.2, 0.3) node[above] {3} -- (7) -- (8) -- ($ (8) + (0, -0.3) $) node[below] {2};
\path[draw, semithick, rounded corners] (5) ++(-0.3, -0.2) -- (5) -- (6) -- ($ (6) + (0.3, -0.2) $);
\path ($ (2)!0.5!(3) $) coordinate (out);
\path ($ (0)!0.5!(1) $) coordinate (idC);
\path (idC) -- (m1) coordinate[pos=0.33] (switch1) coordinate[pos=0.66] (switch2);
\path[draw, semithick, rounded corners] ($ (4) + (-0.2, -0.3) $) node[below] {1} -- (4) -- (out) -- (idC) -- (switch1);
\path[draw, semithick, bend right] (switch1) to (switch2);
\path[draw, semithick, bend right] (switch2) to (m1);
\path[draw, semithick] (m1) -- ($ (m1) + (0.3, 0.2) $);
\path[fill] (out) circle[radius=0.03];
\path[fill] (idC) circle[radius=0.03];
\path[fill] (switch1) circle[radius=0.03];
\path[fill] (switch2) circle[radius=0.03];
\path (out) node[above left] {out};
\path[draw, ->] (mid) to ++(left:1.1);
\path[draw, ->] (mid-low) to ++(left:1.1);
\end{tikzpicture}
\begin{tikzpicture}
\path (2.4, 0.5) node (i1) {1} node[right of=i1] (i2) {2} node[right of=i2] (i3) {3} node[right of=i3] (alpha0) {$ α_0 $} node[right of=alpha0] (alpha0p) {$ α_0' $} node[right of=alpha0p] (idC) {$ \id $ (C)}
node[below right of=i1] (m1) {$ β $ (A)} edge (i1) edge (i2.south west) edge (i2.south) edge (i2.south east)
node[below right of=m1] (m2) {$ β $ (A)} edge (m1) edge (i3.south west) edge (i3.south) edge (i3.south east)
node[below right of=m2] (m3) {$ β $ (A)} edge (m2) edge (alpha0)
node[below right of=m3] (m4) {$ β $ (A)} edge (m3) edge (alpha0p)
node[below right of=m4] (m5) {$ φ π_q μ^2 = α_3 + α_4 $} edge (m4) edge (idC);
\end{tikzpicture}
\caption{An $ α_3 + α_4 $ output with 2/5 arc pointing away from $ D $. Siblings are bound left-to-right. In this example, sibling 1 has no morphisms left-within. Its subtree is a $ μ^{≥3} $ disk or a tail result component.}
\label{fig:subtree-alpha34-away-with-right}
\end{subfigure}
\begin{subfigure}[b]{\linewidth}
\centering
\begin{tikzpicture}[scale=1.5]
\path[use as bounding box] (0, -2) -- (1.5, 0);
\path[draw] (0, 0) -- ++(330:1) coordinate (mid) coordinate[midway] (1) -- ++(left:1) coordinate[midway] (2) ++(down:0.1) -- ++(right:1) coordinate (mid-low) coordinate[midway] (3) -- ++(210:1) coordinate[midway] (4)
++(0, -0.1) -- ++(30:1) coordinate[midway] (m1) -- ++(225:1) coordinate[midway] (m2);
\path[draw] (mid) ++(0.05, -0.1) -- ++(260:1) coordinate[midway] (5);
\path[draw] (mid) ++(0.05, -0.1) -- ++(280:1) coordinate[midway] (6);
\path[draw] (mid) ++(0.1, 0) -- ++(right:1) coordinate[midway] (7);
\path[draw] (mid) ++(0.1, 0) -- ++(330:1) coordinate[midway] (8);
\path[draw, semithick, rounded corners] (7) ++(0.2, 0.3) node[above left] {3} -- (7) -- (8) -- ($ (8) + (0, -0.3) $) node[below] {2};
\path[draw, semithick, rounded corners] (5) ++(-0.2, -0.3) -- (5) -- (6) -- ($ (6) + (0.3, -0.2) $);
\path ($ (2)!0.5!(3) $) coordinate (out);
\path ($ (4)!0.5!(m1) $) coordinate (idC);
\path (out) -- (idC) coordinate[pos=0.33] (switch1) coordinate[pos=0.66] (switch2);
\path[draw, semithick, rounded corners] ($ (m2) + (0, -0.2) $) node[below] {1} -- (m2) -- (idC) -- (switch2);
\path[draw, semithick, bend right] (switch2) to (switch1);
\path[draw, semithick, bend right] (switch1) to (out);
\path[draw, semithick, rounded corners] (out) -- (1) -- ++(0.3, 0.2);
\path[fill] (out) circle[radius=0.03];
\path[fill] (idC) circle[radius=0.03];
\path[fill] (switch1) circle[radius=0.03];
\path[fill] (switch2) circle[radius=0.03];
\path (out) node[above left] {out};
\path[draw, ->] (mid) to ++(left:1.1);
\path[draw, ->] (mid-low) to ++(left:1.1);
\end{tikzpicture}
\begin{tikzpicture}
\path (2.4, 0.5) node (alpha0) {$ α_0 $} node[right of=alpha0] (alpha0p) {$ α_0' $} node[right of=alpha0p] (idC) {$ \id $ (C)} node[right of=idC] (i1) {1} node[right of=i1] (i2) {2} node[right of=i2] (i3) {3}
node[below right of=alpha0] (m1) {} edge (alpha0) edge (alpha0p)
node[below right of=m1] (m2) {} edge (m1) edge (idC)
node[below right of=m2] (m3) {$ β $ (A)} edge (m2) edge (i1.south west) edge (i1.south) edge (i1.south east)
node[below right of=m3] (m4) {$ β $ (A)} edge (m3) edge (i2.south west) edge (i2.south) edge (i2.south east)
node[below right of=m4] (m5) {$ φ π_q μ^{≥3} = α_3 + α_4 $} edge (m4) edge (i3.south west) edge (i3.south) edge (i3.south east);
\end{tikzpicture}
\caption{An $ α_3 + α_4 $ output with 2/5 arc pointing away from $ D $. Siblings are bound left-to-right. In this example, sibling 1 is bound as a final-out disk with the stack formed by its three morphisms left-within.}
\label{fig:subtree-alpha34-away-without-right}
\end{subfigure}
\begin{subfigure}[b]{\linewidth}
\centering
\begin{tikzpicture}[scale=1.5]
\path[draw] (0, 0) -- ++(330:1) coordinate (mid) coordinate[midway] (1) -- ++(left:1) coordinate[midway] (2) ++(down:0.1) -- ++(right:1) coordinate (mid-low) coordinate[midway] (3) -- ++(210:1) coordinate[midway] (4);
\path[draw] (mid) ++(0.05, -0.1) -- ++(260:1) coordinate[midway] (5);
\path[draw] (mid) ++(0.05, -0.1) -- ++(280:1) coordinate[midway] (6);
\path[draw] (mid) ++(0.1, 0) -- ++(right:1) coordinate[midway] (7);
\path[draw] (mid) ++(0.1, 0) -- ++(330:1) coordinate[midway] (8);
\path[draw, semithick, rounded corners] (7) ++(0.2, 0.3) node[above left] {3} -- (7) -- (8) -- ($ (8) + (0, -0.3) $) node[below] {2};
\path[draw, semithick, rounded corners] (5) ++(-0.2, -0.3) -- (5) -- (6) -- ($ (6) + (0.3, -0.2) $);
\path ($ (2)!0.5!(3) $) coordinate (out);
\path ($ (4)!0.5!(m1) $) coordinate (idC);
\path[draw, semithick, rounded corners] ($ (4) + (0, -0.2) $) node[below] {1} -- (4) -- (out) -- (1) -- ++(0.3, 0.2);
\path[fill] (out) circle[radius=0.03];
\path (out) node[above left] {out};
\path[draw, ->] (mid) ++(left:1) to ++(right:0.9);
\path[draw, ->] (mid-low) ++(left:1) to ++(right:0.9);
\end{tikzpicture}
\begin{tikzpicture}
\path (2.4, 0.5) node (i1) {1} node[right of=i1] (delta) {$ δ $} node[right of=delta] (i2) {2} node[right of=i2] (i3) {3}
node[below right of=delta] (m1) {$ β $ (A)} edge (delta) edge (i2.south west) edge (i2.south) edge (i2.south east)
node[below right of=m1] (m2) {$ β $ (A)} edge (m1) edge (i3.south west) edge (i3.south) edge (i3.south east)
node[below left of=m2] {$ φ π_q μ^2 = α_3 + α_4 $} edge (m2) edge (i1.south west) edge (i1.south) edge (i1.south east);
\end{tikzpicture}
\caption{An $ α_3 + α_4 $ outputs with 2/5 arc pointing towards $ D $. Siblings 2 and higher are bound first. The result is then used as first morphism for a first-out disk of sibling~1. In this example, sibling 2 is assumed to be 2-rich and has no morphisms left-within, therefore a $ δ $ insertion is used.}
\label{fig:subtree-alpha34-towards-without-left}
\end{subfigure}
\caption{Examples of subtree construction.}
\end{figure}

\begin{description}[font=\normalfont, align=parleftalign]
\item[Basic structure of $ \Subtree(D) $ and $ \Subresult(D) $:] Inductively for every node in the narrow tree of $ D $, we construct a corresponding h-tree (and finally π-tree) and explain how it shall be inserted into $ \Subtree(D) $. We refer to any h-tree being constructed during this process as a \emph{subtree}. We select inductively for every node of $ \Subtree(D) $ a single result component. This gives a final result component $ \Subresult(D) $ of $ \Subtree(D) $. The construction of $ \Subtree(D) $ is intended to inverse the process of taking the subdisk of result components. One can keep an eye on this during the construction of $ \Subtree(D) $.

\item[Non-root 0-rich nodes $ (l, m) $ without left direct sibling and without morphisms left-within:] Interpret the narrow location $ (l, m) $ as tail node. Trace the tree down to the first narrow location that is 1-rich or 2-rich. Then this narrow location produces a $ β $ (A) tail result component. Use this as result component of the subtree. An example is depicted in Figure \ref{fig:subtree-nonroot-0-rich}.

\item[Non-root 1-rich nodes $ (l, m) $ without left direct sibling and without morphisms left-within:] Now $ (l, m) $ has precisely one stack, and it is either a stack of $ α_4 $ and $ α_0 $, or of multiple $ α_0 $ and $ \id $ (C), or a combination of at least two direct siblings or within morphisms, or just a single stray morphism. In all cases, binding the stack together yields the desired first-order $ β $ (A) first-order tail result component. Note that in the case of a single stray morphism, the $ β $ (A) is actually an additive component of the deformed cohomology input element, and constitutes a leaf node of $ \Subtree(D) $. An example is depicted in Figure \ref{fig:subtree-nonroot-1-rich}.

\item[Non-root nodes $ (l, m) $, 2-rich or with a left direct sibling or a morphism left-within:] Bind all stacks of $ (l, m) $ together. They will serve as inputs for a final-out disk $ h_q μ^{≥3} $. Now let us prepare the final-out morphism. In case $ (l, m) $ has a left direct sibling, the final-out morphism is the $ β $ (A) output of the left direct sibling. If there are morphisms left-within but no left direct sibling, bind them together in a $ β $ (A) or $ β/β' $ (C) tree and use this as final-out morphism. If there are no left direct siblings and no left-within morphisms, then use a simple $ δ $ insertion as as final-out morphism. Note that we have ensured that in terms of $ \Tw\Gtl_q $, this node is really $ μ_{\Tw\Gtl_q}^{≥2} $ and not $ μ_{\Tw\Gtl_q}^1 $. After making the $ h_q μ^{≥3} $ disk, compose it afterwards with any $ \id $ (C) or $ α_0 $ morphisms that may lie right-within $ (l, m) $. An example is depicted in Figure \ref{fig:subtree-nonroot-2-rich}.
\end{description}

\paragraph*{The $ α_3 + α_4 $ case.} If the output is of type $ α_3 + α_4 $, then $ (1, |D|) $ itself is a narrow location. We will treat the case where $ (1, |D|) $ is trivial first and the non-trivial cases afterwards.

\begin{description}[font=\normalfont, align=parleftalign]
\item[Output of type $ α_3 + α_4 $ with $ (1, |D|) $ being a trivial narrow location:] If $ (1, |D|) $ is a trivial narrow location, then $ D $ revolves around a single puncture, with only $ \id $ (C) and $ α_0 $ inputs. The corresponding π-tree is found by composing from left to right.

\item[Output of type $ α_3 + α_4 $ with 2/5 arc pointing away from $ D $, with $ (1, |D|) $ having a right-within morphism:] Due to arrow direction, there does not lie an $ α_0 $ input in the very first angle cut by $ D $, after the output mark. This means the very first input of $ D $ is an $ \id $ (C) just one cut angle after the output mark. Build the subtree of $ (1, |D|) $ as if it were a non-root node, only putting finally $ φ π_q μ^2 $ instead of $ h_q μ^2 $. The final $ μ^2 $ is a product of $ β $ (A) and $ \id $ (C). Since the 2/5 arc points away from $ D $, this $ μ^2 $ has a $ α_4 $ component and its $ φ π_q μ^2 $ has the desired $ α_3 + α_4 $ main result component. An example is depicted in Figure \ref{fig:subtree-alpha34-away-with-right}.

\item[Output of type $ α_3 + α_4 $ with 2/5 arc pointing away from $ D $, with $ (1, |D|) $ without right-within morphisms, but decomposable or with a left-within morphism:] Build the subtree of $ (1, |D|) $ as if it were a non-root node, only putting finally $ φ π_q μ^{≥3} $ instead of $ h_q μ^{≥3} $. Let us check that this indeed yields a $ α_3 + α_4 $ result component. Indeed, the final evaluation $ μ^{≥3} $ yields $ α_4 $ and $ φ π_q μ^{≥3} $ has the desired $ α_3 + α_4 $ main result component. An example is depicted in Figure \ref{fig:subtree-alpha34-away-without-right}.

\item[Output of type $ α_3 + α_4 $ with 2/5 arc pointing away from $ D $, with $ (1, |D|) $ without morphisms within, indecomposable and 0-rich:] Interpret the narrow location $ (l, m) $ as tail node, do not insert a node into $ T $, and continue with the child.

\item[Output of type $ α_3 + α_4 $ with 2/5 arc pointing away from $ D $, with $ (1, |D|) $ without morphisms within, indecomposable and 1-rich:] Then take the subtree of $ (1, |D|) $ as if it were not the root node, but put $ φ π_q μ^2 $ instead of $ h_q μ^2 $ and $ φ π_q μ^{≥3} $ instead of $ h_q μ^{≥3} $ as root. As established in the multiplication and disk tables, this indeed yields an $ α_3 + α_4 $ tail result component.

\item[Output of type $ α_3 + α_4 $ with 2/5 arc pointing towards $ D $, with $ (1, |D|) $ having a morphism left-within:] Note that in the final angle being cut before the output mark, there can not lie an $ α_0 $ input because of the arc direction. Build the subtree of $ (1, |D|) $ now as if $ (1, |D|) $ were not the root node, excluding the final $ \id $ (C) for the moment. This subtree has a $ β $ (A) result component that is in fact an indecomposable angle. Finally, compose the remaining $ \id $ (C) with this $ β $ (A) and obtain $ α_3 + α_4 $ as main result component.

\item[Output of type $ α_3 + α_4 $ with 2/5 arc pointing towards $ D $, with $ (1, |D|) $ without morphism left-within, but decomposable or with morphism right-within:] Decompose $ (1, |D|) $ into indecomposables $ C_1, …, C_k $ (note $ k ≥ 1 $). First combine $ C_2, …, C_k $ from left to right, together with their morphisms within and the morphisms right-within $ C_1 $. By assumption there are direct siblings to be combined or there is at least one morphism right-within $ C_k $. This yields a $ β $ (A) or $ β/β' $ (C) result component. Now use this as first-out morphism to bind a $ φ π_q μ^{≥3} $ disk of $ C_1 $. This disk yields a result $ μ^{≥3} = α_2 $ and hence a main result component $ φ π_q μ^{≥3} = α_3 + α_4 $. An example is depicted in Figure \ref{fig:subtree-alpha34-towards-without-left}.

\item[Output of type $ α_3 + α_4 $ with 2/5 arc pointing towards $ D $, with $ (1, |D|) $ indecomposable without morphisms within and 0-rich:] Then $ α_3 + α_4 $ actually appears as tail result component of a $ φ π_q μ^2 $ or $ φ π_q μ^{≥3} $ further downwards.

\item[Output of type $ α_3 + α_4 $ with 2/5 arc pointing towards $ D $, with $ (1, |D|) $ indecomposable without morphisms within and 1-rich:] By zigzag consistency, the stack of $ (1, |D|) $ cannot consist of a single morphism and also not of the combination of $ α_4 $ and $ α_0 $. Therefore binding this stack together yields a final $ μ^2 $ component equal to $ βα $ (A), $ β' α_2 $ or $ β α_3 $. Now its $ φ π_q μ^2 $ has the desired $ α_3 + α_4 $ first-order tail result component.

\item[Output of type $ α_3 + α_4 $ with 2/5 arc pointing towards $ D $, with $ (1, |D|) $ indecomposable without morphisms within and 2-rich:] Then make a $ φ π_q μ^{≥3} $ first-out disk of $ (1, |D|) $ with a $ δ $-insertion as first morphism. Note that due to the 2-richness this is $ μ^2_{\Tw\Gtl_q} $ or $ μ^{≥3}_{\Tw\Gtl_q} $ and not $ μ^1_{\Tw\Gtl_q} $.
\end{description}

\paragraph*{The $ \id $ (C) case.} If the output is of type $ \id $ (C), then $ (1, |D|) $ is not a narrow location itself. Building the disk is a little easier.

\begin{description}[font=\normalfont, align=parleftalign]
\item[Output of type $ \id $ (C), with $ (1, |D|) $ 0-rich:] Then trace the tree from $ (1, |D|) $ towards the leaves. Pick the subtree of the first narrow location that is 1-rich or 2-rich. Note it has a $ β $ (A) main result component. Changing its root from $ h_q μ^{≥3} $ to $ φ π_q μ^{≥3} $ or from $ h_q μ^2 $ to $ φ π_q μ^2 $ then yields the desired $ \id $ (C) tail result component.

\item[Output of type $ \id $ (C), with $ (1, |D|) $ 1-rich:] In case of a single child with morphism within, it has a subtree with $ h_q μ^2 $ main result component $ β $ (A) associated. The corresponding $ φ π_q μ^2 $ version has the desired $ \id $ (C) first-order tail result component. In case of a stack, this stack cannot consist of a single stray morphism, since this would contradict zigzag consistency. Instead, it must be a stack of $ α_4 $ and $ α_0 $ or multiple $ α_0 $ and $ \id $ (C). Such a stack comes with a subtree $ h_q μ^2 $ and main result component $ β $ (A) or $ β/β' $ (C). Upon replacing $ h_q μ^2 $ by $ φ π_q μ^2 $, we obtain the desired $ \id $ (C) main result component.

\item[Output of type $ \id $ (C), with $ (1, |D|) $ 2-rich:] We now have at least one of the following: two children, one child and one stray morphism, or no child and two stray morphisms. Pick the subtrees of all children, binding direct siblings and morphisms within together as usual. Bind all stacks of stray morphisms together. Finally, tie everything into a $ μ^{≥3} $ disk.
\end{description}

\paragraph*{The $ \id $ (D) case.}

\begin{description}[font=\normalfont, align=parleftalign]
\item[Output of type $ \id $ (D), with $ (1, |D|) $ 0-rich:] Trace the tree downwards to the leaves. Pick the subtree of the first narrow location that is 1-rich or 2-rich. Note it has a $ β $ (A) main result component. Changing its root from $ h_q μ^{≥3} $ to $ φ π_q μ^{≥3} $ or from $ h_q μ^2 $ to $ φ π_q μ^2 $ then yields the desired $ \id $ (D) tail result component.

\item[Output of type $ \id $ (D), with $ (1, |D|) $ 1-rich:] Then $ (1, |D|) $ contains a single stack, and this stack is not a single child without morphisms within. This means the stack produces a $ β $ (A) or $ β/β' $ (C) main result component and the corresponding $ φ π_q μ^2 $ or $ φ π_q μ^{≥3} $ version includes the desired $ \id $ (D) first-order tail result component.

\item[Output of type $ \id $ (D), with $ (1, |D|) $ 2-rich:] Bind together all stacks in $ (1, |D|) $ and finally take $ φ π_q μ^{≥3} $. This is possible due to \autoref{th:subtree-connecting}.
\end{description}

\subsection{Verifying the inverse}
\label{sec:classification-verification}
In this section, we verify that $ \Subdisk $ maps CR result components bijectively to CR disks. Believe it or not, we have defined $ \Subresult(D) $ in such a way that its subdisk is $ D $ again. This already shows that $ \Subdisk $ reaches all CR disks. Proving injectivity of $ \Subdisk $ on $ \CRr $ is harder and requires further constructions.

\begin{lemma}
Assigning subresults provides a map $ \Subresult: \CRd → \PiTr $. We have $ \Subdisk ∘ \Subresult = \Id \restr \CRd $.
\end{lemma}

\begin{proof}
This follows inductively from the construction of $ \Subresult(D) $.
\end{proof}

The above lemma already shows that $ \Subdisk $ reaches all CR disks. Prove injectivity of $ \Subdisk $ on $ \CRr $ is harder. It requires us to show how to reconstruct the basic structure of $ r $ from $ \Subdisk(r) $. By basic structure, we mean a very specific notion: the evaluation tree.

\begin{definition}
Let $ r $ be a result component of a Kadeishvili π-tree $ (T, h_1, …, h_N) $. Then its \emph{evaluation tree} is the decorated ordered tree defined as follows:
\begin{itemize}
\item There is a node for every tail result component of type $ h_q μ^2 $ or $ φ π_q μ^2 $ and every result component of type $ h_q μ^{≥3} $ or $ φ π_q μ^{≥3} $ used in $ r $.
\item For every node, insert as many subsequent nodes above as the used tail part is long.
\item Connect the nodes according to the tree structure of $ T $.
\item Order the nodes horizontally according to their horizontal appearance in $ T $.
\item Regarding decoration of a node $ X $, note that $ X $ determines a result component of a subtree on its own and comes with a $ β $ (A) morphism. This determines a narrow location of $ D(r) $ and is the decoration of $ X $.
\end{itemize}
\end{definition}

\begin{lemma}
\label{th:classification-verification-rough}
Let $ r $ be a result component of a π-tree. Then the narrow tree of $ \Subdisk(r) $ is equal to the evaluation tree of $ r $.
\end{lemma}

\begin{proof}
All nodes in the evaluation tree of $ r $ stand for taking immersed disks and yielding $ β $ (A) morphisms and hence determine narrow locations. The nodes of the evaluation tree are also connected according to inclusion. Now let us show by induction on the height that this inclusion of the evaluation tree in the narrow tree is actually surjective. Let $ N $ be a node in the narrow tree all of whose children appear in the evaluation tree. We will show that then also $ N $ appears in the evaluation tree. In fact, $ r $ needs to bind all children $ C_1, …, C_k $ together in some immersed disk. This determines a node $ M $ in the evaluation tree, and also a narrow location $ M $. But now the children are nested in both $ M $ and $ N $, which means one of $ M $ and $ N $ is included in the other. Now if $ M $ is strictly included in $ N $, then all children $ C_1, …, C_k $ are not direct children of $ M $, in contradiction to our assumption. If $ N $ is strictly included in $ M $, then $ M $ can impossibly be an immersed disk since a version shorter on both sides already bounds a disk. Therefore $ M = N $ as narrow locations and since $ M $ appears in the evaluation tree, we have that $ N $ appears in the evaluation tree, which was to be shown. Finally given the equal structure of the trees, their decorations are also equal.
\end{proof}

Equipped with this characterization, we are ready to prove $ \Subdisk $ injective on $ \CRr $.

\begin{lemma}
The map $ \Subdisk: \CRr → \CRd $ is injective, and hence $ \Subdisk: \CRr \isoto \CRd $.
\end{lemma}

\begin{proof}
Let $ r_1, r_2 ∈ \CRr $ denote two π-trees together with result components whose subdisks $ \Subdisk(r_1) $ and $ \Subdisk(r_2) $ are equal. By \autoref{th:classification-verification-rough}, the evaluation trees of $ r_1 $ and $ r_2 $ are then equal. This means that these trees differ only by the order in which stacks are bound and the order in which the nodes of the evaluation trees are linked together and bound together with morphisms within. It is readily checked using the disk and multiplication schemes of \autoref{tab:components-multiplication}, \ref{fig:components-first-out} and \ref{fig:components-final-out} that there is only a unique way to bind these combinations.

For example, a stack of $ α_4 $ and $ α_0 $ can only be bound in one way, the choice of result component is clear from the next higher node in the evaluation tree. A stack of $ α_0 $ and $ \id $ (C) morphisms can only be bound in one way. A stack of direct siblings and $ α_0 $ and $ \id $ (C) morphisms can also be bound only in one way. Essentially, all steps in the construction of $ \Subtree(D) $ are the unique way to obtain a result component. In other words, there is only one way to compose a result component whose evaluation tree equals the narrow tree of $ \Subdisk(r_1) = \Subdisk(r_2) $. We conclude that $ r_1 = r_2 $.
\end{proof}

\subsection{The case of ID, DS and DW disks}
\label{sec:classification-IDDSDW}
In this section, we verify that $ \Subdisk $ sends ID, DS and DW result components bijectively to ID, DS and DW disks. For ID disks, we sketch inverse constructions similar to the CR case. For DS and DW disks, the statement reduces to combinatorics.

Let us first dedicate ourselves to ID result components.

\begin{lemma}
The map $ \Subdisk: \IDr → \SLd $ is injective and its image is precisely $ \IDd $.
\end{lemma}

\begin{proof}
Injectivity is similar to the case of CR result components. We now show that all ID disks are reached by $ \Subdisk $. Let $ D $ be an ID disk. We provide a preimage through explicit construction. In all cases it can be checked that its subdisk is $ D $ again. We shall distinguish the “regular” inputs from the degenerate ones.

First, assume the degenerate input is of C type. Then the procedure is similar to the CR case. We construct a result component of Figure \ref{fig:subdisk-idD-idC-first} or \ref{fig:subdisk-idD-idC-last}. Let us apply the formalism of result components to $ D $. Regard the narrow location $ (1, |D|) $ itself.

If $ (1, |D|) $ is trivial, then evaluate from left to right all identities except the degenerate input and the final regular one. Then compose with the final regular $ \id $ (C) and note that this precisely produces an $ \id $ (B) result component of Figure \ref{fig:subdisk-idD-idC-first} or \ref{fig:subdisk-idD-idC-last}. Note that $ D $ may have stacked $ α_0 $ inputs directly after the output mark. They are welcome in our construction: They are composed one after another with the ultimate $ β $ (A).

If $ (1, |D|) $ is nontrivial, decompose it into narrow locations $ C_1, …, C_k $. As in the CR case, generally from left to right. In particular, if there is no regular identity left-within $ C_1 $, then evaluate $ C_2, …, C_k $ first and use their result for a first-out disk of $ C_1 $. If there is however a regular identity left-within $ C_1 $, then evaluate $ C_1, …, C_k $ entirely from left to right. In both cases this yields $ h_q (α_2) = \id $ (B). Finally compose with the degenerate $ \id $ (C) and obtain the desired $ \id $ (C) result component. Again, note that stacked $ α_0 $ inputs directly after the output mark are welcome, and that the orientation of the 2/5 arc at the degenerate $ \id $ (C) is relevant.

Second, assume the degenerate input is of B type. We construct an all-in disk as in Figure \ref{fig:subdisk-idD-all-in}. Denote by $ L $ the source and target zigzag path of $ D $. Denote by $ a_0 $ its identity arc. The degenerate B input just before the output mark dictates that $ L $ turns right at $ a_0 $.

If $ L $ is oriented clockwise with $ D $, then the B input preceding the output identity gives rise to an $ α_4 $. As in \autoref{th:subtree-connecting}, the stacks of $ D $ inject into the boundary angles of an all-in disk. The tree corresponding to $ D $ is obtained by binding all stacks into trees, and then evaluating their all-in disk with $ μ^{≥3} = \id $ (D) and hence $ φ π_q μ^{≥3} = \id $ (D).

If $ L $ is oriented counterclockwise with $ D $, then the B input succeeding the output identity gives rise to an $ α_3 $. Again, an all-in disk yields the desired $ \id $ (D) result component.
\end{proof}

Next are DS result components. Since both DS result components and DS disks are defined by combinatorics, this reduces to simple checks.

\begin{proposition}
The map $ \Subdisk: \DSr → \SLd $ is injective and its image is precisely $ \DSd $.
\end{proposition}

\begin{proof}
Let a DS disk of $ L $ and $ (a, b) $ be given. We recapitulate how to reconstruct its result component.

Depending on whether $ L $ turns right or left at $ a $, let the first factor of the inner multiplication be the corresponding $ \id $ (C) or $ α_3 + α_4 $. Correspondingly the second factor will be $ α_3/α_4 $ or $ \id $ (C). This indeed produces an $ \id $ (B) result component, since $ b $ lies in the arc sequence $ S_α $ starting at $ a $. Similarly depending on whether $ L $ turns right or left at $ b $, let the other factor of the outer multiplication be the corresponding $ \id $ (C) or $ α_3/α_4 $.

Now if the two inputs at $ a $ come first in $ D $, insert the inner multiplication on the right side of the tree. If $ a $ comes last, insert the inner multiplication on the left side of the tree.

Finally, note that $ a $ may equal $ a_0 $. The $ α $ (D) angle involved lies on the opposite side of where $ L $ goes. Therefore if $ S $ lies in the same direction as $ L $, then the angle lies on the opposite side and its $ h_q $ involves the identity on $ i $, hence $ a = b $ really yields a result component. If $ S $ lies in the opposite direction of where $ L $ goes, then the angle lies on the side of $ S $ and the possible set of $ b $ indices does not include $ a $.
\end{proof}

\begin{lemma}
The map $ \Subdisk: \DWr → \SLd $ is injective and its image is precisely $ \DWd $.
\end{lemma}

\begin{proof}
This is similar to the DS case. The reader may find the \autoref{th:classification-subdiskshape-DW} and its proof helpful.
\end{proof}

\subsection{Signs and $ q $-parameters}
\label{sec:classification-signs}
In this section, we prove \autoref{th:subdisk-minmodel-signs}, which claims that the sign of a result component is precisely the Abouzaid sign of its subdisk. Our strategy is to start from direct inputs and work our way up. In particular, we first compute signs of result components of h-trees. This approach requires that we define Abouzaid signs also for subdisks of h-trees. Note we treat signs additively everywhere.

\begin{definition}
\label{def:classification-signs-abouzaid}
Let $ D $ be a subdisk of a result component of an h-tree or π-tree. Then its \emph{Abouzaid sign} $ \Abouzaid(D) ∈ ℤ/2ℤ $ is the sum of all $ \# $ signs around $ D $, plus the number of odd inputs $ h_i: L_i → L_{i+1} $ where $ L_{i+1} $ is oriented counterclockwise with $ D $, plus one if it concerns a π-tree and its output $ t: L_1 → L_{N+1} $ is odd and $ L_{N+1} $ is oriented counterclockwise. In case $ r $ is a $ β $ (A) result component, the long version of the subdisk shall be taken.
\end{definition}

Recall that a result component of an h- or π-tree is not only a morphism $ r ∈ \Hom_{\H\DefZigzagCat} (L_1, L_{N+1}) $, but also remembers how it was derived from the tree. In particular, the value of any result component of an h-tree or π-tree does not carry any scalars, except signs. It is of the form
\begin{equation}
\label{eq:signs-relative}
± Q ε \text{ resp. } ± Q t,
\end{equation}
where $ ± $ is a sign, $ Q = q_1 … q_k ∈ ℂ⟦Q_0⟧ $ is a pure product of punctures, $ ε: L_1 → L_{N+1} $ is an elementary morphism resp.~$ t $ is a cohomology basis element of $ \ZigzagCat $. Note \eqref{eq:signs-relative} means we measure the sign relative to the natural signs of the cohomology basis elements.

Let us now show that any result component comes precisely with the Abouzaid sign. For a result component $ r $ of a π-tree or h-tree, we denote by $ S(r) $ the Abouzaid sign of the subdisk of $ r $. We proceed by induction.

\subsubsection*{Signs of direct inputs}
\begin{proposition}
Let $ r $ be a direct result component of an h-tree or π-tree that has a subdisk associated. Then the sign of $ r $ as in \eqref{eq:signs-relative} equals the Abouzaid sign of the subdisk of $ r $.
\end{proposition}

\begin{proof}
The sign is indeed correct for direct inputs. A co-identity $ α_0 $ (D) comes sign $ \#α_0 + 1 $. Its subdisk is by definition on the $ α_0 $ side (instead of the $ α_0' $ side) and by \autoref{conv:alpha0-direction} the zigzag paths runs counterclockwise, which makes the Abouzaid sign of the co-identity also equal to $ \#α_0 + 1 $. A direct $ α_3 $ input comes with sign $ \#α_3 + 1 $. Its subdisk consists of an input on the 2/5 arc and cutting the $ α_3 $ angle, which means the Abouzaid sign is also $ \#α_3 + 1 $. Similarly, a direct $ α_4 $ comes with sign $ \#α_4 $, a direct $ β $ (C) comes with sign $ \#α_3 + \#α_4 + S^D $, a direct $ β' $ (C) with sign $ \#α_1 + \#α_2 + S^D $, a direct $ α_0' $ with sign $ \#α_0 $. Finally, a direct $ β $ (A) as tail component of an $ α_3 + α_4 $ input comes with sign $ \#α_3 + 1 + S^D $ if it is from the $ α_3 $ part of \eqref{eq:deformed-cohomology-basis-alpha34} (and hence $ L_{i+1} $ is counterclockwise) and with sign $ \#α_4 + S^D $ if it is from the $ α_4 $ part of \eqref{eq:deformed-cohomology-basis-alpha34} (and hence $ L_{i+1} $ is clockwise). A direct $ β $ (A) as tail component of an $ \id $ (C) input comes with sign $ \#α_1 + \#α_2 + S^D $ if it is from the $ β' $ part of \eqref{eq:deformed-cohomology-basis-idC} and with sign $ \#α_3 + \#α_2 + S^D $ if it is from the $ β $ part of \eqref{eq:deformed-cohomology-basis-idC}. All these signs are equal to the Abouzaid signs.
\end{proof}

\subsubsection*{Signs of h-trees}
Next, let us check the signs for result components of h-trees.

\begin{lemma}
\label{th:classification-signs-htrees}
Let $ r $ be a result component of an h-tree or π-tree that has a subdisk associated. Then the sign of $ r $ as in \eqref{eq:signs-relative} equals the Abouzaid sign of the subdisk of $ r $.
\end{lemma}

\begin{proof}
The signs are collected in \autoref{tab:signs-h-trees}. This table informs about the applicable $ h_q $ rule, the sum of the input signs, the sum of signs due to $ μ $ applications, the sum of signs due to $ h_q $ applications and the total Kadeishvili sign. Note that by induction assumption, the inputs already have the correct Abouzaid sign.

\begin{table}
\centering
\begin{tabular}{@{}l|lll@{}}
Component & $ α_0' $ of \autoref{fig:subdisk-alpha0p} & $ β $ (C) of \autoref{fig:subdisk-betaC} & $ β' $ (C) of \autoref{fig:subdisk-betaCp} \\\hline
$ h_q $ rule & $ h_q (αα') $ & $ h_q (αα') $, $ h_q (α_4 β) $ & $ h_q (αα') $, $ h_q (β' α_2) $ \\
i-sign & $ N_{α_0} (\#α_0 + 1) - 1 $ & $ N_{α_0} (\#α_0 + 1) - 1 $ & $ N_{α_0} (\#α_0 + 1) - 1 $ \\
$ μ $-sign & $ N_{α_0} - 1 $ & $ N_{α_0} $ & $ N_{α_0} - 1 $ \\
$ h_q $-sign & $ (N_{α_0} - 1) (\#α_0 + 1) $ & $ (N_{α_0} - 1) (\#α_0 + 1) + \#α_4 $ & {$ \scriptstyle (N_{α_0} - 1) (\#α_0 + 1) + \#α_2 + 1 $} \\
K-sign & $ N_{α_0} - 1 $ & $ N_{α_0} $ & $ N_{α_0} $ \\
&&& \\
Component & $ β $ (A) of $ h_q μ^2 (\tilde{β} \text{(A)}, α_0) $ & $ β $ (A) of $ h_q μ^2 (α_0, α_4) $ & $ β $ (A) of $ h_q μ^2 (\id \text{(C)}, α_0') $ \\\hline
$ h_q $ rule & $ h_q (βα) $ & $ h_q (α_3 α_4) $ & $ h_q (α_4 β) $ \\
i-sign & $ S(\tilde{β}) + \#α_0 + 1 $ & $ \#α_0 + 1 + \#α_4 $ & $ S(α_0') $ \\
$ μ $-sign & $ 1 $ & 1 & 1 \\
$ h_q $-sign & $ S^D + \#α_0 + 1 $ & $ S^D + \#α_4 + 1 $ & $ S^D + \#α_4 $ \\
K-sign & + 1 & +1 & +1 \\
&&& \\
Component & $ β $ (A) of $ h_q μ^2 (\tilde{β} \text{(A)}, \id \text{(C)}) $ & $ β $ (A) of $ h_q μ^2 (β/β' \text{(C)}, \id \text{(C)}) $ & $ β $ (A) of $ h_q μ^2 (α_0', \id \text{(C)}) $ \\\hline
$ h_q $ rule & $ h_q (βα) $ & $ h_q (βα) $ & $ h_q (βα) $ \\
i-sign & $ S(\tilde{β}) $ & $ S(β/β') $ & $ S(α_0') $ \\
$ μ $-sign & 0 & 0 & 0 \\
$ h_q $-sign & $ S^D + \#α + 1 $ & $ S^D + \#α + 1 $ & $ S^D + \#α + 1 $ \\
K-sign & +1 & +1 & +1 \\
&&& \\
Component & $ β $ (A) of final-out $ h_q μ^{≥3} $ & $ α_3 $ (B) of \autoref{fig:subdisk-alpha3} & \\\hline
$ h_q $ rule & $ h_q (βα) $ & $ h_q (α_3 α_4) $ & \\
i-sign & $ \sum S(m_i) $ & $ \#α_0 + 1 + \#α_4 $ & \\
$ μ $-sign & 0 & $ 1 $ & \\
$ h_q $-sign & $ S^D + \#α + 1 $ & $ \#α_4 + 1 + S^D $ & \\
K-sign & +1 & +1 &
\end{tabular}
\caption{Signs of most result components of h-trees}
\label{tab:signs-h-trees}
\end{table}

In those rows of \autoref{tab:signs-h-trees} that concern a disk result component, the $ m_1, …, m_k $ refer to the inner angles of the disk. That is, the disk is supposed to be $ μ^{k≥3} (m_k, …, m_1) $. Some $ m_i $ may be $ δ $ insertions. Their $ S(m_i) $ shall stand for zero. For $ β $ (A) result components, $ S^D $ shall stand for zero if it concerns a main result component, while equal to the sum of the hash signs around the tail if it concerns a tail result component. For the $ α_0' $, $ β/β' $ (C) result components, $ N_{α_0} $ denotes the number of $ α_0 $ inputs in the result component (including the $ α_0' $ used). A single time the notation $ \tilde{β} $ was used to distinguish two different $ β $ (A) result components.

Let us examine the example of an $ β $ (A) result component of $ h_q μ^2 (α_0, α_4) $ in detail. Both inputs $ α_0 $ and $ α_4 $ are necessarily direct. The morphism $ α_0 $ comes with sign $ \#α_0 + 1 $ and $ α_4 $ comes with sign $ \#α_4 $. This gives a total input sign of $ \#α_0 + 1 + \#α_4 $. We have $ μ^2 (α_0, α_4) = - α_3 α_4 $, which gives a sign of $ 1 $ due to application of $ μ $. According to \autoref{th:tail-sum}, we have
\begin{equation*}
h_q (α_3 α_4) = (-1)^{\#α_4 + 1} \left(α_3 + \sum_{D ∈ T(α_3) \setminus ε} (-1)^{S^D} Q^D β^D\right).
\end{equation*}
This gives a sign of $ \#α_4 + 1 + S^D $ due to application of $ h_q $. A Kadeishvili sign of $ 1 $ is added. The total sign the $ β $ (A) tail result component is now
\begin{equation*}
\#α_0 + 1 + \#α_4 + 1 + \#α_4 + 1 + S^D + 1 \equiv \#α_0 + S^D.
\end{equation*}
Let us compare with the Abouzaid sign $ S(β) $. This sign consists of all $ \# $ signs around the tail, including $ \#α_0 $, plus two because both odd inputs $ α_0 $ and $ α_3 + α_4 $ are counterclockwise. This is precisely the same.

Note \autoref{tab:signs-h-trees} does not treat explicitly the $ β $ (A) result component of \autoref{fig:subdisk-betaA-mult-main}. This result component is however a combination of $ β/β' $ (C) or $ β $ (A) and $ α_0 $ and $ \id $ (C) compositions on the right, which are already present in \autoref{tab:signs-h-trees}.
\end{proof}

\subsubsection*{Signs of $ π $-trees}
\begin{lemma}
Let $ r $ be a direct result component of a π-tree that has a subdisk associated. Then the sign of $ r $ as in \eqref{eq:signs-relative} equals the Abouzaid sign of the subdisk of $ r $.
\end{lemma}

\begin{proof}
The most signs of π-trees are checked in \autoref{tab:signs-pi-trees}. Let us treat the others manually.

\begin{table}
\centering
\begin{tabular}{l|ll}
Component & $ α_3 + α_4 $ of final-out & $ α_3 + α_4 $ of first-out $ μ^{k≥3} (m_k, …, m_1) $ \\\hline
$ μ^{≥2} $ & $ α_4 $ & $ α_2 $ \\
i-sign & $ \sum S(m_i) $ & $ \sum S(m_i) $ \\
$ μ $-sign & 0 & $ 1 $ \\
$ φπ_q $-sign & $ \#α_4 + 1 $ & $ \#α_2 + 1 $ \\
&& \\
Component & $ α_3 + α_4 $ of $ μ^2 (β \text{(A)}, \id \text{(C)}) $ & $ α_3 + α_4 $ of $ μ^2 (\id \text{(C)}, β \text{(A)}) $ \\\hline
$ μ^{≥2} $ & $ α_4 $ & $ α_2 $ \\
i-sign & $ S(β) $ & $ S(β) $ \\
$ μ $-sign & 0 & $ 1 $ \\
$ φπ_q $-sign & $ \#α_4 + 1 $ & $ \#α_2 + 1 $ \\
&& \\
Component & $ \id $ (C) of all-in \\\hline
$ μ^{≥2} $ & $ \id $ (C) \\
i-sign & $ \sum S(m_i) $ \\
$ μ $-sign& 0 \\
$ φπ_q $-sign & 0
\end{tabular}
\caption{Signs of some π-trees}
\label{tab:signs-pi-trees}
\end{table}

Let us check the case where $ α_3 + α_4 $ is the G tail result component of some $ φπ_q (βα) $. We could theoretically check this by going through all the cases. It is easier to rely on what we already have. Namely $ - h_q (βα) $ has an associated main result component $ β $, which comes out of $ h_q μ^{≥2} $ with the correct sign $ S(β) $. This means the $ μ^{≥2} $ must have had sign $ S(β) + \#α + 1 + 1 $ in front of $ βα $. Then its $ α_3 + α_4 $ G tail result component comes with an additional sign of $ S^D + \#α + 1 $ in case of G1 and $ S^D + \#α $ in case of G2. In total, we get a sign of $ S(β) + S^D + 1 $ in case of a G1 tail result component and $ S(β) + S^D $ in case of a G2 result component, precisely the Abouzaid sign.

Let us check the case where $ \id $ (C) is the H tail result component of some $ φπ_q (βα) $. Then $ - h_q (βα) $ has an associated result component $ β $, which comes out of $ h_q μ^{≥2} $ with the correct sign $ S(β) $. This means the $ μ^{≥2} $ must have had sign $ S(β) + \#α + 1 + 1 $ in front of $ βα $. Then its $ \id $ (C) H tail result component comes with an additional sign of $ S^D + \#α $. In total, we get a sign of $ S(β) + S^D $, precisely the Abouzaid sign.

Checks for the $ \id $ (D) result components of \autoref{fig:subdisk-idD} are contained in \autoref{tab:signs-idD-trees}. For example, regard the case the $ \id $ (B) comes from a first-out disk $ μ^{k≥3} (m_k, …, m_1) $. According to \autoref{fig:components-first-out}, the outside part of $ m_1 $ is $ α_2 $ and hence odd. We get $ μ^{k≥3} (m_k, …, m_1) = - α_2 $ and evaluation of $ h_q (α_2) $ gives another sign of $ (-1)^{\#α_2} $. Together with the Kadeishvili sign we obtain that $ \id $ (D) has sign $ S(m_k) + … + S(m_1) + \#α_2 $ as result component, precisely the Abouzaid sign.

\begin{table}
\centering
\begin{tabular}{l|ll}
Component & all-in disk $ μ^{k≥3} (m_k, …, m_1) $ & $ φπ_q μ^2 (h_q μ^2 (\id \text{(C)}, β \text{(A)}), \id \text{(C)}) $ \\\hline
$ h_q $ & – & $ h_q (α_2) $ \\
i-sign & $ \sum S(m_i) $ & $ S(β) $ \\
$ μ/h_q/φπ_q $-sign & 0 & $ 1 + \#α_2 $ \\
K-sign & 0 & $ 1 $ \\
&& \\
Component & $ φπ_q μ^2 (h_q μ^2 (β \text{(A)}, \id \text{(C)}), \id \text{(C)}) $ & $ φπ_q μ^2 (\id \text{(C)}, h_q μ^2 (\id \text{(C)}, β \text{(A)})) $ \\\hline
$ h_q $ & $ h_q (βα) $ & $ h_q (α_2) $ \\
i-sign & $ S(β) $ & $ S(β) $ \\
$ μ/h_q/φπ_q $-sign & $ S^D + \#α + 1 $ & $ 1 + \#α_2 $ \\
K-sign & $ 1 $ & $ 1 $ \\
&& \\
Component & $ φπ_q μ^2 (\id \text{(C)}, h_q μ^2 (β \text{(A)}, \id \text{(C)})) $ & $ φπ_q (\id \text{(C)}, h_q μ^{k≥3} (m_k, …, m_1)) $ \\\hline
$ h_q $ & $ h_q (βα) $ & $ h_q (α_2) $ \\
i-sign & $ S(β) $ & $ \sum S(m_i) $ \\
$ μ/h_q/φπ_q $-sign & $ S^D + \#α + 1 $ & $ 1 + \#α_2 $ \\
K-sign & $ 1 $ & $ 1 $ \\
&& \\
Component & $ φπ_q (h_q μ^{k≥3} (m_k, …, m_1), \id \text{(C)}) $ & \\\hline
$ h_q $ & $ h_q (α_2) $ & \\
i-sign & $ \sum S(m_i) $ & \\
$ μ/h_q/φπ_q $-sign & $ 1 + \#α_2 $ & \\
K-sign & $ 1 $ &
\end{tabular}
\caption{Signs of $ \id $ (D) result components of \autoref{fig:subdisk-idD}.}
\label{tab:signs-idD-trees}
\end{table}

Finally, let us check the $ α_3 + α_4 $ and $ \id $ (C) result components of the 8 trees of \autoref{fig:subdisk-degenerate}. Recall such a result component a degenerate strip on a zigzag path $ L $ as subdisk.

Let us investigate the inner product first. In case of $ μ^2 (α_3, \id \text{(C)}) $ resp.~$ μ^2 (\id \text{(C)}, α_4) $, the inner product has sign $ \#α_3 + 1 $ resp.~$ \#α_4 + 1 $. Regard the infinitesimally short stem of the strip between the two factors of the inner product. If the angle $ α_3 $ resp.~$ α_4 $ as morphism $ L → L $ falls under case 1 of \autoref{fig:coh_splitting-identity-location}, then $ h_q $ adds a sign of $ \#α_3 + 1 $ resp.~$ \#α_4 + 1 $. Together with the Kadeishvili sign, the inner product has a total sign of $ 1 $. Indeed, the stem is counterclockwise in this case. If $ α_3 $ resp.~$ α_4 $ falls under case 2, then $ h_q $ adds a sign of $ \#α_3 $ resp.~$ \#α_4 $. Together with the Kadeishvili sign, the inner product has a total sign of $ 0 $. Indeed, the stem is clockwise in this case.

Let us now investigate the outer product. Regard the infinitesimally short stem at the output mark.

When $ \id $ (C) or $ α_3 + α_4 $ comes at the end of the strip, $ φ π_q μ $ adds no sign (in case of $ α_3 + α_4 $, the intrinsic sign $ \#α_4 $ stays correctly until the end). In case $ α_3 + α_4 $ comes at the end of the strip, this indeed constitutes an odd intersection and an odd output and both add the same sign since both refer to the orientation of the stem. In case $ \id $ (C) comes at the end of the strip, this indeed produces an even intersection and an even output mark.

When $ \id $ (C) comes at the beginning of the strip, $ φ π_q μ $ gets no sign. Indeed, the output is then also even. When $ α_3 + α_4 $ comes at the beginning of the strip, then $ φ π_q μ $ gets a single extra sign. The output is then also $ α_3 + α_4 $ and their claimed signs refer to the orientation of $ L $ and of its Hamiltonian deformation $ L' $. Since both point in the same direction, but lie on opposite sides of the strip, they contribute to the Abouzaid sign with $ 1 $.
\end{proof}

\section{The case of punctured spheres}
\label{sec:sphere}
In this section, we redo our entire minimal model computation in the case of specific punctured spheres. In particular, our treatise includes the 3-punctured sphere, also known as pair of pants. The simplest yet instructive example of mirror symmetry for punctured surfaces, it would be a shame not to know the minimal model of its deformed zigzag category. However, no dimer on a sphere is consistent and \autoref{th:subdisk-minmodel-th} therefore fails to apply. This is the reason we redo the entire calculation in the case of specific sphere dimers with $ M ≥ 3 $ punctures.

In contrast to the case of geometrically consistent dimers, the simplified Kadeishvili construction of \autoref{sec:2Bkadeishvili} does not apply. For this reason, we need to resort to the general deformed Kadeishvili theorem from \papertwoA. We shall not attempt to review the general construction here, and only offer the following minimalistic explanation: First, there is no direct sum decomposition which deserves the name of “deformed splitting”. Instead, the construction builds a \emph{deformed decomposition} with lesser properties. Second, the construction iteratively uncurves the category until it has \emph{optimal curvature}. Third, the structure on the minimal model is described in terms of Kadeishvili trees. We refer to \papertwoA\ for details.

As we conclude in \autoref{th:sphere-minmodel-odd} and \ref{th:sphere-minmodel-even}, the minimal model $ \H\DefZigzagCat $ can be described explicitly by means of CR, ID, DS and DW disks. This description also accurately captures the curvature und residual differential on $ \H\DefZigzagCat $. Explicitly, for odd $ M $ and specific choice of spin structure the category $ \H\DefZigzagCat $ is curvature-free and has a residual differential. For even $ M $ the category $ \H\DefZigzagCat $ has curvature and residual differential.

In \autoref{sec:sphere-problem}, we give an overview of what goes wrong for non-consistent dimers. In \autoref{sec:sphere-zigzagcat}, we focus on specific sphere dimers $ Q_M $ with an odd number of punctures $ M ≥ 3 $. We choose a specific type of spin structure, tailored to the use in mirror symmetry. In \autoref{sec:sphere-splitting}, we choose a homological splitting. In \autoref{sec:sphere-defzigzagcat}, we describe the deformed zigzag category $ \DefZigzagCat $. In \autoref{sec:sphere-defsplitting}, we compute the deformed decomposition of $ \DefZigzagCat $. In \autoref{sec:sphere-resultcomp}, we introduce the suitable notion of result components. In \autoref{sec:sphere-minmodel}, we assemble the minimal model $ \H\DefZigzagCat $. In \autoref{sec:sphere-even}, we comment on the case of $ Q_M $ for even $ M $.

\subsection{Absence of consistency}
\label{sec:sphere-problem}
In this section, we list consequences of the lack of geometric consistency. Within \autoref{sec:splitting} until \ref{sec:subdisk}, we have used geometric consistency heavily. Is geometric consistency actually a necessity? To find an answer, we collect the most important statements which we proved using consistency. For every statement, we explain how it depends on consistency and whether it can be partially recovered when dropping the consistency assumption.
\begin{description}
\item[Zigzag segments do not bound disks.]
We proved this statement directly using geometric consistency. Upon dropping geometric consistency, zigzag segments can easily bound disks.
\item[Deformed zigzag paths are uncurvable.]
We proved this by explicit uncurving, which succeeds because zigzag segments do not bound disks. When zigzag segments bound disks, two issues can occur: If the zigzag segment ends in an A situation, we need to adapt the uncurving procedure, but it provides no hindrance to uncurving. If the zigzag curve actually bounds a teardrop, then the zigzag path can inherently not be uncurved.
\item[A, B, C, D situations exhaust all angles.]
This statement mainly uses that every puncture has at least 4 arc incidences. In fact, a non-consistent dimer may have punctures with only 2 arc incidences. This gives rise to new types similar to B and C situations, with the difference that the head or tail of the shared arc $ 2=5 $ may have 2 arc incidences. More practically, this would push $ β $ (C) or $ β' $ (C) to be empty angles.
\item[$ \DefZigzagCat $ satisfies $ μ^1_{\DefZigzagCat} (H) ⊂ μ^1_{\DefZigzagCat} (ℂ⟦Q_0⟧ \htensor R) $.]
If $ Q $ is so nonconsistent that $ \DefZigzagCat $ has inherent curvature, then this statement is already not applicable any more. Namely, $ μ^1_{\DefZigzagCat} $ does not square to zero anymore and we cannot invoke the simplified Kadeishvili theorem of \autoref{sec:2Bkadeishvili}. If $ Q $ is only so nonconsistent that $ \DefZigzagCat $ can be uncurved (albeit by an adapted procedure), then the inclusion typically does not hold. Nevertheless, the deformed Kadeishvili theorem applies and yields a minimal model with residue differential. It is clear that the description of $ H_q $ will be very complicated.
\item[E, F, G, H disks as classification of tail terms.]
The shape of the terms in $ μ^{≥3}_q (δ, …, ε, …, δ) $ is analyzed by zooming in at the concluding puncture. In distinguishing E, F, G, H disks, we have used that the concluding puncture has at least 4 arc incidences. When dropping consistency, the resulting terms need not be of A, B, C, D type, but also of the variants explained above.
\item[Situation B/C cohomology basis elements have only type E tail.]
This is proved directly using geometric consistency: An intersection between $ L_1 $ and $ L_2 $ renders it impossible to find type F, G, H disks when tracing $ L_1 $ and $ L_2 $ away from the intersection. Upon dropping consistency, B/C cohomology basis elements acquire tail also from F, G, H disks. This raises additional complexity: The G and H disks will contribute result components, forcing us into capturing them.
\item[Every narrow locations has at least one below morphism.]
A narrow location without below morphism constitutes a zigzag segment bounding a disk. Upon dropping consistency, imagine a segment of a zigzag path $ L $ that bounds a disk. According to our explanations above, the $ δ $-matrix should already be adapted to facilitate uncurving. In fact, the new $ δ $-matrix of $ L $ will have a situation A morphism for every disk bounded by one of its segments and more situation A morphisms inserted on its tail. These situation A morphisms compensate for the lack of below morphisms. This modification allows us to construct result components for a given CR/ID disk $ D $ even if a narrow location $ (l, m) $ of $ D $ has no below morphism.
\end{description}

The calculation for the sphere dimers $ Q_M $ basically proceeds as in the geometrically consistent case. Based on the above list of issues, we can however point out a few differences: In the consistent case, $ \DefZigzagCat $ is always uncurvable. In the sphere case, for odd $ M $ and specific choice of spin structure it is uncurvable, for even $ M $ it is not uncurvable. In the consistent case, minimal models constructed by our deformed Kadeishvili construction are based on a deformed cohomology space $ H_q $ satisfying $ μ^1_{\DefZigzagCat} (H_q) = 0 $. In the sphere case, we only achieve $ μ^1_{\DefZigzagCat} (H_q) ⊂ H_q $. In the consistent case, the deformed counterpart $ φ^{-1} (h) $ of $ h = \id $ (C) morphisms includes tails of $ β/β' $ (C). In the sphere case, the deformed counterpart $ φ^{-1} (h) $ of $ h = \id $ (C) morphisms includes no tails, but a single nearby $ \id $ (B) morphism.

\subsection{The sphere and its zigzag category}
\label{sec:sphere-zigzagcat}
In this section, we define specific sphere dimers and define their category of zigzag paths. The dimers we pick are those also used as A-side for commutative mirror symmetry in \cite{MS-sphere}. The dual dimers of these spheres are consistent and therefore suited for noncommutative mirror symmetry of \cite{Bocklandt}.

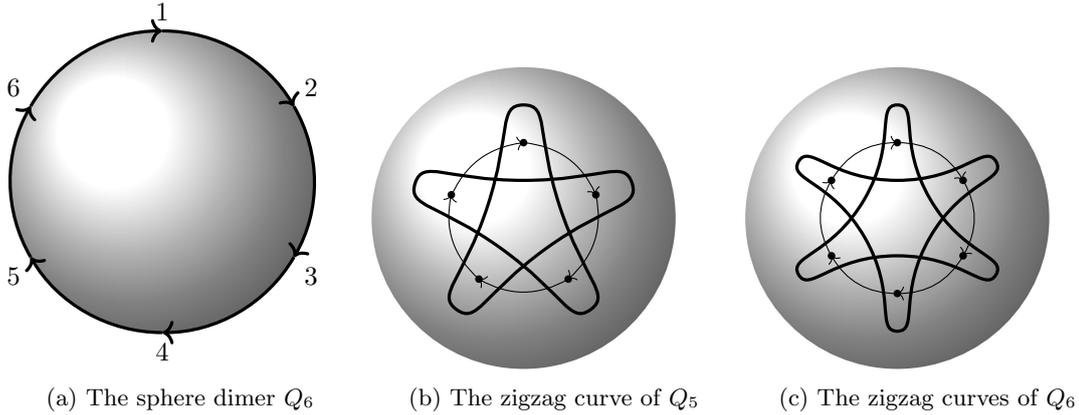
\begin{figure}
\centering
\begin{subfigure}[b]{0.3\linewidth}
\begin{tikzpicture}
\path[shade, ball color=white] (0, 0) circle[radius=2];
\path[draw, very thick, ->] (90:2) arc(90:30:2) node[at start, above] {1};
\path[draw, very thick, ->] (30:2) arc(30:-30:2) node[at start, above right] {2};
\path[draw, very thick, ->] (-30:2) arc(-30:-90:2) node[at start, below right] {3};
\path[draw, very thick, ->] (270:2) arc(270:210:2) node[at start, below] {4};
\path[draw, very thick, ->] (210:2) arc(210:150:2) node[at start, below left] {5};
\path[draw, very thick, ->] (150:2) arc(150:90:2) node[at start, above left] {6};
\end{tikzpicture}
\caption{The sphere dimer $ Q_6 $}
\label{fig:sphere-zigzagcat-dimer}
\end{subfigure}
\begin{subfigure}[b]{0.3\linewidth}
\begin{tikzpicture}
\path[shade, ball color=white] (0, 0) circle[radius=2];
\path[draw, very thick] plot[smooth cycle, tension=1.3] coordinates {(90:1.5) (18:0.5) (-54:1.5) (234:0.5) (162:1.5) (90:0.5) (18:1.5) (-54:0.5) (234:1.5) (162:0.5)};
\path[fill] (90:1) circle[radius=0.05];
\path[fill] (18:1) circle[radius=0.05];
\path[fill] (-54:1) circle[radius=0.05];
\path[fill] (234:1) circle[radius=0.05];
\path[fill] (162:1) circle[radius=0.05];
\path[draw, ->, bend left] (90:1) to (18:1);
\path[draw, ->, bend left] (18:1) to (-54:1);
\path[draw, ->, bend left] (-54:1) to (234:1);
\path[draw, ->, bend left] (234:1) to (162:1);
\path[draw, ->, bend left] (165:1) to (90:1);
\end{tikzpicture}
\caption{The zigzag curve of $ Q_5 $}
\label{fig:sphere-zigzagcat-Q5zigzag}
\end{subfigure}
\begin{subfigure}[b]{0.3\linewidth}
\begin{tikzpicture}
\path[shade, ball color=white] (0, 0) circle[radius=2];
\path[draw, very thick] plot[smooth cycle, tension=1.2] coordinates {(90:1.5) (30:0.5) (-30:1.5) (270:0.5) (210:1.5) (150:0.5)};
\path[draw, very thick] plot[smooth cycle, tension=1.2] coordinates {(30:1.5) (-30:0.5) (270:1.5) (210:0.5) (150:1.5) (90:0.5)};
\path[fill] (90:1) circle[radius=0.05];
\path[fill] (30:1) circle[radius=0.05];
\path[fill] (-30:1) circle[radius=0.05];
\path[fill] (270:1) circle[radius=0.05];
\path[fill] (210:1) circle[radius=0.05];
\path[fill] (150:1) circle[radius=0.05];
\path[draw, ->, bend left] (90:1) to (30:1);
\path[draw, ->, bend left] (30:1) to (-30:1);
\path[draw, ->, bend left] (-30:1) to (270:1);
\path[draw, ->, bend left] (270:1) to (210:1);
\path[draw, ->, bend left] (210:1) to (150:1);
\path[draw, ->, bend left] (150:1) to (90:1);
\end{tikzpicture}
\caption{The zigzag curves of $ Q_6 $}
\label{fig:sphere-zigzagcat-Q6zigzag}
\end{subfigure}
\caption{The sphere dimer and its zigzag curves}
\end{figure}

The dimer we regard is the sphere dimer $ Q_M $ for $ M ≥ 3 $ depicted in \autoref{fig:sphere-zigzagcat-dimer}. This dimer has $ M $ punctures and $ M $ arcs. It has two polygons, namely the clockwise \emph{front} side and the counterclockwise \emph{rear} side. We shall briefly discuss the differences between the cases of odd and even $ M $, and then focus on the odd case. In order to apply our findings to deformed mirror symmetry later on, we shall define one specific spin structure.

The zigzag curves of this sphere dimer $ Q_M $ are described as follows: In case $ M $ is odd, there is precisely one zigzag curve. It cycles around the arcs once, and then cycles around the arcs again with opposite $ δ $ angles. In case $ M $ is even, there are precisely two zigzag curves, each of them cycling around the arcs once. The smooth zigzag curves in both cases are depicted in \autoref{fig:sphere-zigzagcat-Q6zigzag}. In the picture, the arc system has been pulled to the front side of the sphere so that the zigzag curves become clearly visible.

Let us now focus on the case of odd $ M ≥ 3 $ and fix spin structure as follows.

\begin{convention}
The letter $ Q = Q_M $ stands for the sphere dimer with $ M ≥ 3 $ odd. The spin structure of the zigzag path is chosen by assigning $ \# α = 1 $ to an odd number of interior angles on the rear side of $ Q_M $, and $ \# α = 0 $ to all other angles. The co-identity locations $ α_0 $ are supposed to lie on the rear side and the identity locations $ a_0 $ are arbitrary indexed arcs.
\end{convention}

In \autoref{def:sphere-zigzagcat-def}, we define the category of zigzag paths $ \ZigzagCat ⊂ \Tw\Gtl Q_M $ as in the case of geometrically consistent dimers. In the case $ Q = Q_M $, the only object in the category is the single zigzag path $ L $.

\begin{definition}
\label{def:sphere-zigzagcat-def}
The \emph{category of zigzag paths} $ \ZigzagCat ⊂ \Tw\Gtl Q_M $ is the full subcategory consisting of the single zigzag path.
\end{definition}

We intend to write down the explicit twisted complex for $ L $. The main issue consists of numbering all indexed arcs of $ L $ and the angles between them. In fact, a zigzag path consists of indexed arcs, as opposed to purely arcs of $ Q_M $. We shall therefore denote the arcs in sequence by $ a_1, …, a_{2M} $, with the convention that $ h(a_1) = q_M $ and $ t(a_1) = q_1 $ and $ L $ turns right at the head of $ a_1 $.

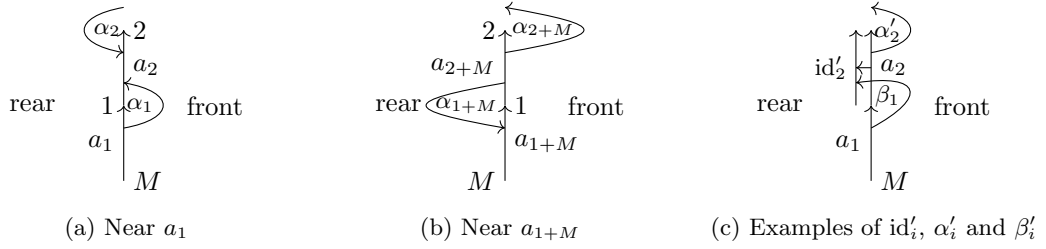
\begin{figure}
\centering
\begin{subfigure}[b]{0.3\linewidth}
\centering
\begin{tikzpicture}
\path[draw, ->] (0, 0) node[right] {$ M $} -- (0, 1) node[midway, left] {$ a_1 $} node[left] {$ 1 $} coordinate[pos=0.7] (1-start);
\path[draw, ->] (0, 1) -- (0, 2) node[right] {$ 2 $} coordinate[pos=0.3] (1-end) coordinate[pos= 0.7] (2-end) node[midway, right] {$ a_2 $};
\path[draw, ->, bend right=80, looseness=3] (1-start) to node[midway, left] {\small $ α_1 $} (1-end);
\path[draw, ->, bend right=80, looseness=3] (2-end)++(up:0.6) to node[midway, right] {\small $ α_2 $} (2-end);
\path (1.2, 1) node {front} (-1.2, 1) node {rear};
\end{tikzpicture}
\caption{Near $ a_1 $}
\label{fig:sphere-zigzagcat-numbering-1}
\end{subfigure}
\begin{subfigure}[b]{0.3\linewidth}
\centering
\begin{tikzpicture}
\path[draw, ->] (0, 0) node[left] {$ M $} -- (0, 1) node[midway, right] {$ a_{1+M} $} node[right] {$ 1 $} coordinate[pos=0.7] (1-start);
\path[draw, ->] (0, 1) -- (0, 2) node[left] {$ 2 $} coordinate[pos=0.3] (1-end) coordinate[pos= 0.7] (2-end) node[midway, left] {$ a_{2+M} $};
\path[draw, ->, bend right=80, looseness=6] (1-end) to node[midway, right] {\small $ α_{1+M} $} (1-start);
\path[draw, ->, bend right=80, looseness=6] (2-end) to node[midway, left] {\small $ α_{2+M} $} ($ (2-end) + (up:0.6) $);
\path (1.2, 1) node {front} (-1.4, 1) node {rear};
\end{tikzpicture}
\caption{Near $ a_{1+M} $}
\label{fig:sphere-zigzagcat-numbering-M}
\end{subfigure}
\begin{subfigure}[b]{0.3\linewidth}
\centering
\begin{tikzpicture}
\path[draw, ->] (0, 0) node[right] {$ M $} -- (0, 1) node[midway, left] {$ a_1 $} node[left] {} coordinate[pos=0.7] (1-start);
\path[draw, ->] (0, 1) -- (0, 2) node[left] {} coordinate[pos=0.3] (1-end) coordinate[pos= 0.7] (2-end) node[midway, right] {$ a_2 $} coordinate[midway] (idp-start);
\path[draw, ->] (-0.2, 1) -- (-0.2, 2) coordinate[pos=0.3] (beta1-end) coordinate[pos=0.5] (idp-end);
\path[draw, ->, bend right=80, looseness=3.5] (1-start) to node[pos=0.4, left] {\small $ β_1 $} (beta1-end);
\path[draw, ->, bend right=80, looseness=3] (2-end) to node[midway, left] {\small $ α_2' $} ($ (2-end) + (up:0.6) $);
\path[draw, ->] (idp-start) -- (idp-end) node[left] {\small $ \id_2' $};
\path (1.2, 1) node {front} (-1.2, 1) node {rear};
\end{tikzpicture}
\caption{Examples of $ \id_i' $, $ α_i' $ and $ β_i' $}
\label{fig:sphere-zigzagcat-numbering-exotic}
\end{subfigure}
\caption{Numbering of punctures, arcs and angles}
\end{figure}

The indexed small angles of $ L $ shall be denoted by $ α_1, …, α_{2M} $, such that $ α_1: a_1 → a_2 $ and $ α_2: a_3 → a_2 $ and so on. In other words, $ α_i $ runs at the head of $ a_i $ if $ i $ odd and at the tail of $ a_i $ if $ i $ is even. In other words, we have $ α_{2i}: a_{2i + 1} → a_{2i} $ and $ α_{2i+1}: α_{2i+1} → α_{2i+2} $. These angles are depicted in \autoref{fig:sphere-zigzagcat-numbering-1}.

We count all indices modulo $ 2M $. In contrast, an index shift of $ M $ typically turns a situation from left to right and from right to left. For example, we have $ α_{i+M} ≠ α_i $. Compare \autoref{fig:sphere-zigzagcat-numbering-1} and \ref{fig:sphere-zigzagcat-numbering-M}. We also have $ a_i = a_{i+M} $ as arcs of $ Q_M $, but not as indexed arcs of $ L $.

For $ i = 1, …, 2M $, we denote the complementary angle to $ α_i $ by $ α_i' $. For instance if $ i $ is odd, then $ α_i $ runs from $ a_{i+1} $ to $ a_i $. We denote the arc identity of $ a_i $ by $ \id_i $, and the arc identity $ a_i → a_{i+M} $ by $ \id_i' $. We set $ β_i ≔ \id_j' α_i $ and $ β_i' ≔ \id_j' α_i' $, where $ j $ is chosen so that the composition makes sense. We are now ready to define the deformed category of zigzag paths:

Generically denoting a full turn by $ ℓ $, these angles together with their multiples with $ ℓ $ powers form a basis of the hom space $ \Hom_{\ZigzagCat} (L, L) $. Examples of these angles are depicted in \autoref{fig:sphere-zigzagcat-numbering-exotic}.

\subsection{Homological splitting}
\label{sec:sphere-splitting}
In this section, we provide a homological splitting for $ \ZigzagCat $ in case of odd $ M ≥ 3 $. We first explain the analogy of all basis morphisms with the consistent case. Then, we write down an explicit choice of cohomology basis elements and an explicit choice of $ R $. We explain why it constitutes a homological splitting.

We have seen that $ \End_{\ZigzagCat} (L, L) $ has a basis given by basis morphisms of the kind $ \id_i $, $ \id_i' $, $ α_i $, $ α_i' $, $ β_i $ and $ β_i' $. Let us compare with the consistent case. The angle $ \id_i $ is an arc identity. In terms of A, B, C, D situations, we denoted it as $ \id $ (D). The angle $ \id_i' $ is comparable to an $ \id $ (C) morphism for odd $ i $, and comparable to an $ \id $ (B) morphism for even $ i $. The angle $ α_i $ is simply $ α $ (D), and similarly $ α_i' $ is $ α' $ (D). There are ambiguities of interpreting $ β_i $ and $ β_i' $ in terms of A, B, C, D situations. A possible choice is matching $ β_i $ with $ α_4 $ (B) for odd $ i $ and with $ α_1 $ (B) for even $ i $, and matching $ β_i' $ with $ α_3 $ (B) for odd $ i $ and with $ α_2 $ (B) for even $ i $. In short, the two differences are that we have explicit indices $ i $ instead of using A, B, C, D situations as enumeration tools and that we have less morphisms overall since ever puncture has only 2 arc incidences.

We now define a candidate splitting $ H ⊕ I ⊕ R $, modeled after the consistent case:

\begin{definition}
Denote by $ H ⊂ \Hom_{\ZigzagCat} (L, L) $ the space spanned by the \emph{cohomology basis elements}
\begin{itemize}
\item $ id_i' $ for $ i $ odd
\item $ (-1)^{\#(i+M) + 1} β_i' + (-1)^{\#i} β_i $ where $ i $ odd,
\item $ \sum_{i = 1}^{2M} \id_{a_i} $,
\item $ (-1)^{\#α_0 + 1} α_0 $.
\end{itemize}
Choose the space $ R ⊂ \Hom_{\ZigzagCat} (L, L) $ to be spanned by $ β_i $ for even $ i $, $ β_i' $ for odd $ i $, $ (α_i' α_i)^{k+1} $ for all $ i $ and $ k ≥ 0 $, $ β_i ℓ^{k+1} $ for even $ i $ and $ k ≥ 0 $, $ β_i' ℓ^{k+1} $ for odd $ i $ and $ k ≥ 0 $, $ \id_i $ for $ i ≠ i_0 $, $ \id_i' $ for even $ i $, $ β_i α_i' ℓ^k $ for odd $ i $ and $ k ≥ 0 $, $ α_i' ℓ^k $ for all $ i $ and $ k ≥ 0 $, $ β_i α_i' ℓ^k $ for even $ i $ and $ k ≥ 0 $. The spaces $ H $ and $ R $, together with $ I ≔ \Im(μ^1_{\ZigzagCat}) $, constitute the \emph{standard splitting} for $ \ZigzagCat $.
\end{definition}

\begin{table}
\centering
\begin{tabular}{rcp{0.8\linewidth}}
$ \id_i $ & $ = $ & $ \id_i ∈ R, \text{ if } i ≠ i_0,  $ \\
$ \id_i $ & $ = $ & $ \sum_j \id_j - \sum_{j ≠ i} \id_j, \text{ if } i = i_0,  $ \\
$ \id_i' $ & $ = $ & $ \id_i' ∈ R, \text{ even } i,  $ \\
$ \id_i' $ & $ = $ & $ \id_i' ∈ H, \text{ odd } i,  $ \\
$ α_i $ & $ = $ & $ α_i ∈ H, \text{ if } α_i = α_0,  $ \\
$ α_i $ & $ = $ & $ μ^1 (± \id_{a_j} ±  … ± \id_{a_i}) ± α_0,  $ \\
$ α_i' $ & $ = $ & $ α_i' ∈ R,  $ \\
$ β_i $ & $ = $ & $ β_i ∈ R, \text{ even } i,  $ \\
$ β_i $ & $ = $ & $ (-1)^{\#i} ((-1)^{\#(i+M) + 1} β_i' + (-1)^{\#i} β_i) + (-1)^{\#i + \#(i+M)} β_i' \text{ odd } i,  $ \\
$ β_i' $ & $ = $ & $ (-1)^{\#(i+M)} μ^1 (\id_i') + (-1)^{\#(i+M)} [(-1)^{\#(i+M-1) + 1} β_{i-1}' + (-1)^{\#(i-1)} β_{i-1}] + (-1)^{\#i + \#(i+M)} β_i, \text{ even } i,  $ \\
$ β_i' $ & $ = $ & $ β_i' ∈ R, \text{ odd } i,  $ \\
$ α_i ℓ^k $ & $ = $ & $ (-1)^{\#i} μ^1 ((α_i' α_i)^k),  $ \\
$ α_i' ℓ^k $ & $ = $ & $ α_i' ℓ^k ∈ R,  $ \\
$ β_i ℓ^k $ & $ = $ & $ β_i ℓ^k ∈ R, \text{ even } i,  $ \\
$ β_i ℓ^k $ & $ = $ & $ (-1)^{\#i + 1} μ^1 (β_i α_i' ℓ^{k-1}) + (-1)^{\#i + \#(i+M)} β_i' ℓ^k, \text{ odd } i,  $ \\
$ β_i' ℓ^k $ & $ = $ & $ β_i' ℓ^k ∈ R, \text{ odd } i,  $ \\
$ β_i' ℓ^k $ & $ = $ & $ (-1)^{\#(i+M)} μ^1 (β_i α_i' ℓ^{k-1}) + (-1)^{\#(i+M) + \#i} β_i ℓ^k, \text{ even } i,  $ \\
$ (α_i' α_i)^k $ & $ = $ & $ (α_i' α_i)^k ∈ R,  $ \\
$ (α_i α_i')^k $ & $ = $ & $ (-1)^{\#i + 1} μ^1 (α_i' ℓ^{k-1}) + (-1) (α_i' α_i)^k,  $ \\
$ β_i α_i' ℓ^k $ & $ = $ & $ β_i α_i' ℓ^k ∈ R,  $ \\
$ β_i' α_i ℓ^k $ & $ = $ & $ (-1)^{\#i + 1} μ^1 (β_i' ℓ^k), \text{ odd } i,  $ \\
$ β_i' α_i ℓ^k $ & $ = $ & $ (-1)^{\#(i+M) + 1} μ^1 (β_i ℓ^k), \text{ even } i. $ 
\end{tabular}
\caption{Decomposing arbitrary morphisms into $ H $, $ I $ and $ R $}
\label{tab:sphere-odd-decomposing}
\end{table}

In \autoref{tab:sphere-odd-decomposing} we have checked that every morphism in $ \End_{\ZigzagCat} (L, L) $ can be written in terms of $ H $, $ I $ and $ R $. The table also serves as a convenient reference for the definition of $ H $ and $ R $. In analogy to the consistent case, the sum $ H + I + R $ is in fact direct. Let us record this as follows:

\begin{lemma}
The spaces $ H $, $ I = \Im(μ^1_{\ZigzagCat}) $ and $ R $ provide a homological splitting for $ \ZigzagCat $.
\end{lemma}

\subsection{Deformed category of zigzag paths}
\label{sec:sphere-defzigzagcat}
In this section, we define the category $ \DefZigzagCat $ of deformed zigzag paths. As in the consistent case, its definition is based on the complementary angle trick. In contrast to the consistent case, we must expect that uncurving fails. Due to our specific spin structure, curvature cancels nevertheless.

\begin{definition}
\label{def:sphere-defzigzagcat-def}
Regard a sphere dimer $ Q = Q_M $ with odd $ M ≥ 3 $. Let $ \# $ denote a choice of spin structure as in \autoref{def:sphere-zigzagcat-def}. Then the \emph{deformed category of zigzag paths} is category $ \DefZigzagCat ⊂ \Tw\Gtl_q Q_M $ consisting of the single deformed zigzag path
\begin{align*}
L &= (a_1 ⊕ … ⊕ a_{2M}, δ), \\
δ &= \begin{bmatrix} 0 & (-1)^{\#1} q_1 α_1' & 0 & … & 0 \\ (-1)^{\#1} α_1 & 0 & (-1)^{\#2} α_2 & … & 0 \\ 0 & (-1)^{\#2} q_2 α_2' & 0 & \ddots & 0 \\ \vdots & \vdots & \ddots & \ddots & (-1)^{\#(2M)} α_{2M} \\ 0 & 0 & 0 & (-1)^{\#(2M)} q_{2M} α_{2M}' & 0 \end{bmatrix}.
\end{align*}
\end{definition}

Let us introduce the following shorthand notation for $ k ∈ ℤ $:
\begin{equation}
\label{eq:sphere-defzigzagcat-shorthand}
\begin{aligned}
\#\#k & ≔ \#k + … + \#(k+M-1), \\
Q_i &≔ q_i q_{i+2} … q_{i+M-3}, \\
Q^{\sphereodd} (k) &≔ \prod_{\substack{s = k \\ s \text{ odd}}}^{k+M-1} q_s, \\
Q^{\sphereeven} (k) &≔ \prod_{\substack{s = k \\ s \text{ even}}}^{k+M-1} q_s.
\end{aligned}
\end{equation}

We now come to our first meaningful calculation in the sphere case: the curvature of $ \DefZigzagCat $. We expect curvature in principle, since $ L $ is contractible when regarded in the closed surface $ |Q| $. With our specific spin structure, the curvature contributions from front and rear side however cancel each other.

\begin{lemma}
\label{th:sphere-defzigzagcat-curvature}
The curvature of $ L $ vanishes, we have $ μ^0_{\DefZigzagCat} = 0 $.
\end{lemma}

\begin{proof}
We have to evaluate $ μ^0_{\Add\Gtl_q Q} + μ^2_{\Add\Gtl_q Q} (δ, δ) + μ^{≥3}_{\Add\Gtl_q Q} (δ, …, δ) $. As in the geometrically consistent case, the first two terms cancel each other. In contrast, the term $ μ^M (δ, …, δ) $ yields two individual contribution for each index $ i = 1, …, 2M $, one from the front and one from the rear side:
\begin{equation*}
μ^M (δ, …, δ) = \sum_{i = 1}^{2M} (-1)^{\#\#i} Q_i^{\sphereeven} \id_i' + \sum_{i = 1}^{2M} (-1)^{\#\#(i-M)} Q_{i-M}^{\sphereodd} \id_i' = 0.
\end{equation*}
We have used that $ Q_i^{\sphereeven} = Q_{i-M}^{\sphereodd} $ and $ \#\#i + \#\#(i-M) ≡ 1 ∈ ℤ/2ℤ $ by assumption on the spin structure.
\end{proof}

\subsection{Deformed decomposition}
\label{sec:sphere-defsplitting}
In this section, we provide the deformed decomposition for $ \DefZigzagCat $. More precisely, the category $ \DefZigzagCat $ is curvature-free but does not satisfy $ μ^1_q (B \htensor H) ⊂ μ^1_q (B \htensor R) $. The deformed Kadeishvili theorem from \papertwoA\ nevertheless applies and defines a deformed decomposition $ H_q ⊕ μ^1_q (B \htensor R) ⊕ (B \htensor R) $ with $ μ^1_q (H_q) ⊂ H_q $. In this section, we compute $ H_q $ explicitly, together with a few values of the deformed codifferential.

As a preparation, we perform here a few calculations of $ μ^1_{\DefZigzagCat} $:
\begin{align*}
μ^1_{\DefZigzagCat} (\id_i') &= μ^2 (\id_i', (-1)^{\# α_{i-1}} α_{i-1} + (-1)^{\# α_i} α_i) + μ^2 ((-1)^{\# α_{i+M}} α_{i+M} + (-1)^{\# α_{i+M-1}} α_{i+M-1}, \id_i') \\
&= (-1)^{\#(i-1) + 1} β_{i-1} + (-1)^{\#i + 1} β_i + (-1)^{\#(i+M)} β_i' + (-1)^{\#(i+M-1)} β_{i-1}' \\
& \text{ for even } i, \\
μ^1_{\DefZigzagCat} (\id_i') &= μ^2 (\id_i', (-1)^{\#i} q_i α_i' + (-1)^{\#(i-1)} q_{i-1} α_{i-1}') \\
& \qquad + μ^2 ((-1)^{\#(i+M)} q_{i+M} α_{i+M}' + (-1)^{\#(i+M-1)} q_{i+M-1} α_{i+M-1}', \id_i') \\
&= (-1)^{\#i + 1} q_i β_i' + (-1)^{\#(i-1) + 1} q_{i-1} β_{i-1}' + (-1)^{\#(i+M)} q_i β_i + (-1)^{\#(i+M-1)} q_{i-1} β_{i-1} \\
& \text{ for odd } i, \\
μ^1_{\DefZigzagCat} (α_i ℓ^k) &= μ^2 (α_i ℓ^k, (-1)^{\#i} q_i α_i') + μ^2 ((-1)^{\#i} q_i α_i', α_i ℓ^k) + μ^M (δ, …, α_i ℓ^k, …, δ) \\
&= (-1)^{\#i + 1} q_i [(α_i α_i')^{k+1} + (α_i' α_i)^{k+1}] + μ^M (δ, …, α_i ℓ^k, …, δ), \\
μ^1_{\DefZigzagCat} (α_i' ℓ^k) &= μ^2 (α_i' ℓ^k, (-1)^{\#i} α_i) + μ^2 ((-1)^{\#i} α_i, α_i' ℓ^k) + μ^M (δ, …, α_i' ℓ^k, …, δ) \\
&= (-1)^{\#i + 1} [(α_i α_i')^{k+1} + (α_i' α_i)^{k+1}] + μ^M (δ, …, α_i' ℓ^k, …, δ), \\
μ^1_{\DefZigzagCat} (β_i ℓ^k) &= μ^2 (β_i ℓ^k, (-1)^{\#i} q_i α_i') + μ^2 ((-1)^{\#(i+M)} α_{i+M}, β_i ℓ^k) + μ^M (δ, …, β_i ℓ^k, …, δ) \\
&= (-1)^{\#i + 1} q_i β_i α_i' ℓ^k + (-1)^{\#(i+M) + 1} β_i' α_i ℓ^k + μ^M (δ, …, β_i ℓ^k, …, δ), \\
μ^1_{\DefZigzagCat} (β_i' ℓ^k) &= μ^2 (β_i' ℓ^k, (-1)^{\#i} α_i) + μ^2 ((-1)^{\#(i+M)} q_i α_{i+M}', β_i' ℓ^k) + μ^M (δ, …, β_i' ℓ^k, …, δ) \\
&= (-1)^{\#i + 1} β_i' α_i ℓ^k + (-1)^{\#(i+M) + 1} q_i β_i α_i' ℓ^k + μ^M (δ, …, β_i' ℓ^k, …, δ).
\end{align*}

\begin{remark}
\label{th:sphere-defsplitting-paramprod}
The shorthand notation \eqref{eq:sphere-defzigzagcat-shorthand} has the property that
\begin{equation*}
Q^\sphereeven_k = q_k Q^{\sphereeven}_{k+1} \text{ for even } k,
\quad Q^\sphereodd_k = q_k Q^\sphereodd_{k+1} \text{ for odd } k.
\end{equation*}
For even $ k $ we have
\begin{equation*}
Q^\sphereodd_{k-M+1} = Q^{\sphereeven}_{k+2} = q_{k+2} … q_{k+M-2},
\quad Q^\sphereeven_k = q_k … q_{k+M-1} = q_k · (q_{k+2} … q_{k+M-2}) = q_k · Q^\sphereodd_{k-M+1}.
\end{equation*}
We have used that $ \#\#(k+1) = \#\#k + \#(k+M) - \#k $.
\end{remark}

In order to apply the deformed Kadeishvili theorem, we need to compute the space $ H_q $ along the procedure of \autoref{def:2Bkadeishvili-deformed-counterpart}. In other words, we shall calculate the deformed counterparts $ h - Eh $ of the cohomology basis elements $ h ∈ H $. By abuse of notation, let us write $ \#α_0 $ for the $ \# $ sign associated with $ α_0 $ and $ q_{α_0} $ for the puncture around which $ α_0 $ winds.

\begin{lemma}
\label{th:sphere-defsplitting-basis}
The deformed cohomology basis elements are given by
\begin{equation*}
(-1)^{\#(i+M) + 1} β_i' + (-1)^{\#i} β_i \text{ for odd } i
\end{equation*}
and
\begin{equation*}
\id_i' + (-1)^{\#(i-1) + \#(i+M-1)} q_{i-1} \id_{i-1}' \text{ for odd } i
\end{equation*}
and
\begin{equation*}
\sum_{i = 1}^{2M} \id_i \quad \text{and} \quad (-1)^{\#α_0 + 1} α_0 + (-1)^{\#α_0} q_{α_0} α_0'.
\end{equation*}
The differentials $ μ^1_{\DefZigzagCat} $ for these morphisms are given by
\begin{align*}
μ^1_{\DefZigzagCat} ((-1)^{\#(i+M) + 1} β_i' + (-1)^{\#i} β_i) &= (-1)^{\#\#(i+1) + 1} Q_{i+2} \id_L ∈ H_q
\end{align*}
and
\begin{align*}
μ^1_{\DefZigzagCat} (\id_i' + (-1)^{\#(i-1) + \#(i+M-1)} q_{i-1} \id_{i-1}') &= (-1)^{\#i + \#(i+M)} q_i \bigg((-1)^{\#(i+M) + 1} β_i' + (-1)^{\#i} β_i\bigg) \\
& \hspace{-10em} + (-1)^{\#(i-1) + \#(i+M-1) + 1} q_{i-1} \bigg((-1)^{\#(i-2+M) + 1} β_{i-2}' + (-1)^{\#(i-2)} β_{i-2}\bigg) ∈ H_q
\end{align*}
and
\begin{align*}
μ^1_{\DefZigzagCat} \left(\sum_{i = 1}^{2M} \id_i\right) &= 0
\end{align*}
and
\begin{align*}
μ^1_{\DefZigzagCat} ((-1)^{\#α_0 + 1} α_0 + (-1)^{\#α_0} q_{α_0} α_0') &= \\
& \hspace{-10em} \sum_{\substack{j = 0 \\ j \text{ even}}}^{M-1} (-1)^{\#\#(i+j-M+1) + 1} Q^\sphereodd_{i+j-M+1}
[\id_{i+j+1}' + (-1)^{\#(i+j) + \#(i+j-M)} q_{i+j-M} \id_{i+j}'] \\
& \hspace{-10em} + \sum_{\substack{j = 2 \\ j \text{ even}}}^{M-1} (-1)^{\#\#(i-j+1)} Q^\sphereeven_{i-j+1}
[\id_{i-j+1}' + (-1)^{\#(i-j+M) + \#(i-j)} q_{i-j} \id_{i-j}']
 ∈ H_q.
\end{align*}
\end{lemma}

\begin{proof}
We need to check two things: First, all added infinitesimal terms lie in $ ℂ⟦Q_0⟧ \htensor R $. Second, the map $ μ^1_{\DefZigzagCat} $ sends all the deformed basis elements to $ H_q $.

The first step is an easy observation: Indeed $ \id_{i-1}' $ for odd $ i $ and $ α_0' $ lie in $ R $. For the second part, we need to evaluate $ μ^1_q $ on the deformed cohomology basis elements and check that the result belongs to $ H_q $. We execute all calculations in order:

First, we regard the morphism $ (-1)^{\#(i+M) + 1} β_i' + (-1)^{\#i} β_i $ for odd $ i $. We compute
\begin{align*}
& μ^1_{\DefZigzagCat} ((-1)^{\#(i+M) + 1} β_i' + (-1)^{\#i} β_i) = (-1)^{\#(i+M) + \#i} β_i' α_i + q_i β_i α_i' - q_i β_i α_i' + (-1)^{\#i + \#(i+M) + 1} β_i' α_i \\
& \quad + \sum_{j = 0}^{M-1}  μ^M \bigg((-1)^{\#(i+j+1)} [q_{i+j+1}] α_{i+j+1} ['], …, (-1)^{\#(i+M-1)} q_{i+M-1} α_{i+M-1}, \\
& \qquad\qquad (-1)^{\#(i+M) + 1} β_i', (-1)^{\#(i+1)} α_{i+1}, (-1)^{\#(i+2)} q_{i+2} α_{i+2}', …, (-1)^{\#(i+j)} [q_{i+j}] α_{i+j}[']\bigg) \\
& \quad + \sum_{j = 0}^{M-1}  μ^M \bigg((-1)^{\#(i+2M-1-j)} [q_{i+2M-1-j}] α_{i+2M-1-j}['], …, (-1)^{\#(i+M+1)} α_{i+M+1}, \\
& \qquad\qquad (-1)^{\#i} β_i, (-1)^{\#(i-1)} α_{i-1}', …, (-1)^{\#(i-j)} [q_{i-j}] α_{i-j}[']\bigg) \\
&= \sum_{j = 0}^{M-1} (-1)^{\#\#(i+1) + 1} Q_{i+2} \id_{i+j+1} + \sum_{j = 0}^{M-1} (-1)^{\#\#(i+M+1)} Q_{i+2} \id_{i-j} \\
&= (-1)^{\#\#(i+1) + 1} Q_{i+2} \id_L ∈ H_q.
\end{align*}
We have used that $ \#\#(k+M) + \#\#(k) $ is the total number of $ \# $ signs in the dimer, which is odd by assumption.

Second, we regard the morphism $ \id_i' + (-1)^{\#(i-1) + \#(i+M-1)} q_{i-1} \id_{i-1}' $ for odd $ i $. We compute
\begin{align*}
& μ^1_{\DefZigzagCat} (\id_i' + (-1)^{\#(i-1) + \#(i+M-1)} q_{i-1} \id_{i-1}') \\
& = (-1)^{\#i + 1} q_i β_i' + (-1)^{\#(i-1) + 1} q_{i-1} β_{i-1}' \\
& \quad + (-1)^{\#(i+M)} q_i β_i + (-1)^{\#(i+M-1)} q_{i-1} β_{i-1} \\
& \quad + (-1)^{\#(i-2) + \#(i-1) + \#(i+M-1) + 1} q_{i-1} β_{i-2} + (-1)^{\#(i+M-1) + 1} q_{i-1} β_{i-1} \\
& \quad + (-1)^{\#(i-1)} q_{i-1} β_{i-1}' + (-1)^{\#(i+M-2) + \#(i-1) + \#(i+M-1)} q_{i-1} β_{i-2}' \\
&= (-1)^{\#i + 1} q_i β_i' + (-1)^{\#(i+M)} q_i β_i \\
& \quad + (-1)^{\#(i-2) + \#(i-1) + \#(i+M-1) + 1} q_{i-1} β_{i-2} + (-1)^{\#(i+M-2) + \#(i-1) + \#(i+M-1)} q_{i-1} β_{i-2}' \\
&= (-1)^{\#i + \#(i+M)} q_i \bigg((-1)^{\#(i+M) + 1} β_i' + (-1)^{\#i} β_i\bigg) \\
& \quad + (-1)^{\#(i-1) + \#(i+M-1) + 1} q_{i-1} \bigg((-1)^{\#(i-2+M) + 1} β_{i-2}' + (-1)^{\#(i-2)} β_{i-2}\bigg) ∈ H_q.
\end{align*}
Third, we regard the identity $ \id_L = \sum_i \id_i $ and compute
\begin{equation*}
μ^1_q (\id_L) = μ^2 (\id_L, δ) + μ^2 (δ, \id_L) = 0.
\end{equation*}
Fourth, we deal with the co-identity $ (-1)^{\#α_0 + 1} α_0 + (-1)^{\#α_0} q_{α_0} α_0' $. Let $ i $ be the index where $ α_0 $ is located, such that $ α_i = α_0 $. Note that $ i $ is even, since $ α_0 $ is supposed to lie on the counterclockwise side of $ Q_M $. In evaluating $ μ^1 $ on the deformed co-identity, there appear two types of terms: four $ μ^2 $ terms and many $ μ^M $ terms. The $ μ^2 $ terms cancel each other as in the case for consistent dimers:
\begin{align*}
& μ^2 ((-1)^{\#α_0 + 1} α_0, (-1)^{\#α_0} q_{α_0} α_0') + μ^2 ((-1)^{\#α_0} q_{α_0} α_0', (-1)^{\#α_0} α_0) \\
& \quad + μ^2 ((-1)^{\#α_0} q_{α_0} α_0', (-1)^{\#α_0 + 1} α_0) + μ^2 ((-1)^{\#α_0} α_0, (-1)^{\#α_0} q_{α_0} α_0') \\
&= q_{α_0} α_0 α_0' - q_{α_0} α_0' α_0 + q_{α_0} α_0' α_0 - q_{α_0} α_0 α_0' = 0.
\end{align*}
We are now ready to calculate the $ μ^M $ terms:
\begin{align*}
& \sum_{j = 0}^{M-1} μ^M \bigg((-1)^{\#(i+j-M+1)} [q_{i+j-M+1}] α_{i+j-M+1}^{[']}, …, (-1)^{\#(i-1)} q_{i-1} α_{i-1}', \\
& \qquad\qquad (-1)^{\#α_0 + 1} α_0, (-1)^{\#(i+1)} q_{i+1} α_{i+1}', …, (-1)^{\#(i+j)} [q_{i+j}] α_{i+j}^{[']}\bigg) \\
& \qquad + \sum_{j = 0}^{M-1} μ^M \bigg((-1)^{\#(i-j+M-1)} [q_{i-j+M-1}] α_{i-j+M-1}^{[']}, …, (-1)^{\#(i+1)} α_{i+1}, \\
& \qquad\qquad (-1)^{\#α_0} q_{α_0} α_0', (-1)^{\#(i-1)} α_{i-1}, …, (-1)^{\#(i-j)} [q_{i-j}] α_{i-j}^{[']}\bigg) \\
&= \sum_{j = 0}^{M-1} (-1)^{\#\#(i+j-M+1) + 1} Q^{\sphereodd}_{i+j-M+1} \id_{i+j+1}' + \sum_{j = 0}^{M-1} (-1)^{\#\#(i-j)} Q^{\sphereeven}_{i-j} \id_{i-j}' \\
&= 
(-1)^{\#\#(i-M+1) + 1} Q^{\sphereodd}_{i-M+1} \id_{i+1}'
+ (-1)^{\#\#(i)} Q^{\sphereeven}_i \id_i' \\
& \qquad + \sum_{\substack{j = 2 \\ j \text{ even}}}^{M-1} (-1)^{\#\#(i+j-M+1) + 1} Q^{\sphereodd}_{i+j-M+1} \id_{i+j+1}'
+ \sum_{\substack{j = 2 \\ j \text{ even}}}^{M-1} (-1)^{\#\#(i+j-M) + 1} Q^{\sphereodd}_{i+j-M} \id_{i+j}' \\
& \qquad + \sum_{\substack{j = 2 \\ j \text{ even}}}^{M-1} (-1)^{\#\#(i-j)} Q^{\sphereeven}_{i-j} \id_{i-j}'
+ \sum_{\substack{j = 2 \\ j \text{ even}}}^{M-1} (-1)^{\#\#(i-j+1)} Q^{\sphereeven}_{i-j+1} \id_{i-j+1}' \\
&=
\sum_{\substack{j = 0 \\ j \text{ even}}}^{M-1} (-1)^{\#\#(i+j-M+1) + 1} Q^\sphereodd_{i+j-M+1}
[\id_{i+j+1}' + (-1)^{\#(i+j) + \#(i+j-M)} q_{i+j-M} \id_{i+j}'] \\
& \qquad + \sum_{\substack{j = 2 \\ j \text{ even}}}^{M-1} (-1)^{\#\#(i-j+1)} Q^\sphereeven_{i-j+1}
[\id_{i-j+1}' + (-1)^{\#(i-j+M) + \#(i-j)} q_{i-j} \id_{i-j}']
 ∈ H_q.
\end{align*}
We have used \autoref{th:sphere-defsplitting-paramprod} and $ \#\#(i+1) ≡ \#\#(i-M+1) + 1 ∈ ℤ/2ℤ $. These calculations show that the claimed elements are indeed the deformed cohomology basis elements, which finishes the proof.
\end{proof}

The deformed Kadeishvili theorem lays out the following procedure: We have already uncurved the category $ \DefZigzagCat $ successfully. We have the full space $ H_q $ in our hands. Next we have to compute the deformed codifferential $ h_q: \Hom_{\DefZigzagCat} (L, L) → B \htensor R $. After that, we will be able to evaluate Kadeishvili trees and derive the minimal model.

It is not necessary to compute the entire codifferential $ h_q $. Instead, the most important cohomology basis morphisms are of the form $ β_i + β_i' $ and $ α_0 + q α_0' $ and $ \id_i' + q \id_{i-1}' $. Any product $ μ^{≥3}_q $ of these can only produce an identity. Any product $ μ^2_q $ of these can only produce $ β_i α_i' $ or $ β_i' α_i $ or $ α_0 α_0' $ or $ α_0' α_0 $ or $ β_i $ or $ β_i' $ or $ α_i $ or $ α_i' $. It suffices to calculate the codifferential of these morphisms.

Let us analyze all the easy cases before we calculate the harder ones: For those morphisms lying in $ R $, the codifferential immediately vanishes. Moreover, for odd $ i $ the element $ β_i + β_i' $ lies in $ H_q $ and $ β_i' $ lies in $ R $, thus $ h_q (β_i) = 0 $ for odd $ i $. The nontrivial cases are as follows:

\begin{lemma}
\label{th:sphere-odd-codifferential}
We have the following values of the codifferential:
\begin{align*}
h_q (α_i) &= ± \id_{a_j} + … ± \id_{a_i} \text{ for } α_i ≠ α_0, \\
h_q (β_i') &= (-1)^{\#(i+M)} \id_i' \text{ for even } i, \\
h_q (β_i' α_i) &= (-1)^{\#i + 1} β_i' \text{ for odd } i, \\
h_q (β_i' α_i) &= (-1)^{\#(i+M) + 1} β_i \text{ for even } i, \\
h_q (α_0 α_0') &= (-1)^{\#i + 1} α_0'.
\end{align*}
\end{lemma}

\begin{proof}
The first two cases are simple: The value of $ μ^1_q $ on identities equals the value of $ μ^1 $ and therefore decomposition of $ α_i $ and $ β_i' $ from \autoref{tab:sphere-odd-decomposing} remains valid. We remark that for $ α_i $ the right-hand side needs to be written as $ α_0 + q α_0' - q α_0' $, but since $ α_0' ∈ R $ this is no issue.

The third case consists of checking $ β_i' α_i $ for odd $ i $:
\begin{align*}
μ^1_q (β_i') &= (-1)^{\#i + 1} β_i' α_i + (-1)^{\#(i+M) + 1} q_i β_i α_i' ℓ^k + μ^M (…).
\end{align*}
The term $ β_i α_i' $ lies in $ R $. The terms resulting from $ μ^M (…) $ are all of the form $ \id_i $. The $ h_q $ of such terms necessarily vanishes, and we deduce the above codifferential equation.

The fourth case of $ β_i' α_i $ for even $ i $ is similar. Finally, we check the fifth case of $ α_0 α_0' $:
\begin{equation*}
μ^1_q (α_0') = (-1)^{\#i + 1} [(α_0 α_0') + (α_0' α_0)] + μ^M (…).
\end{equation*}
The terms resulting from $ μ^M (…) $ are all of the form $ \id_i' $. These either lie directly in $ R $ or they lie in $ H $ when combined with additional $ \id_{i-1}' ∈ R $.
\end{proof}

In \autoref{th:sphere-odd-codifferential}, we have saved ourselves from computing the describing the correct signs of $ h_q (α_i) $. In fact, the signs are analogous to those presented in \autoref{sec:splitting-splitting}.

\subsection{Result components}
\label{sec:sphere-resultcomp}
In this section, we analyze result components of $ \H\DefZigzagCat $ and match them with CR, ID, DS and DW disks. The starting point is the category $ \DefZigzagCat $. In \autoref{sec:sphere-defsplitting} we have already computed the deformed cohomology basis elements and the deformed codifferential. Here, we regard Kadeishvili trees, analyze the shape of their outputs and introduce a suitable notion of result components. We introduce a suitable notion of CR, ID, DS and DW disks and match all result components with smooth disks of these four types.

As in the classical case, we start by computing a multiplication table for important endomorphisms of $ L ∈ \DefZigzagCat $. The multiplication table is found in \autoref{tab:sphere-odd-multiplication}.

\begin{table}
\centering
\renewcommand{\arraystretch}{1.1}
\begin{tabular}{cc|cccccc}
$ m_2 $ & \textbackslash $ m_1 $ & $ \id_i $ & $ \id_i' $ & $ β_i $ & $ β_i' $ & $ α_0 $ & $ α_i' $ \\\hline
$ \id_i $ & $ μ^2 = $ & $ \id_i $ & $ \id_i' $ & $ β_i $        & $ β_i' $               & $ α_0 $ & $ α_i' $ \\
          & $ i = $   & any       & odd        & odd            & even                   & any     & any   \\
          & $ h_q = $ & $ 0 $     & $ 0 $      & $ 0 $          & $ \id_i' $             & $ 0 $   & $ 0 $ \\
          & $ π_q = $ & $ (\id_i) $     & $ \id_i' $ & $ β_i + β_i' $ & $ β_{i-1} + β_{i-1}' $ & $ α_0 $ & $ 0 $ \\
\hline
$ \id_i' $ & $ μ^2 = $ & $ \id_i' $ & $ \id_i $ & $ α_i $ & $ α_i' $ & $ β_i $ & $ β_i' $ \\
           & $ i = $   & odd        & any       & any     & any      & even  & even   \\
           & $ h_q = $ & $ 0 $      & $ 0 $     & $ \id_i $ & $ 0 $  & $ 0 $ & $ \id_i' $ \\
           & $ π_q = $ & $ \id_i' $ & $ (\id_i) $     & $ α_0 $ & $ 0 $  & $ 0 $ & $ β_{i-1} + β_{i-1}' $ \\
\hline
$ β_i $ & $ μ^2 = $ & $ β_i $ & $ α_i' $ & $ α_i' α_i $ & imp         & imp         & $ β_i α_i' $ \\
        & $ i = $   & odd     & any      & any          & any   & any   & any \\
        & $ h_q = $ & $ 0 $   & $ 0 $    & $ 0 $        & $ 0 $ & $ 0 $ & $ 0 $ \\
        & $ π_q = $ & $ β_i + β_i' $ & $ 0 $ & $ 0 $    & $ 0 $ & $ 0 $ & $ 0 $ \\
\hline
$ β_i' $ & $ μ^2 = $ & $ β_i' $               & $ α_i $   & imp         & $ α_i α_i' $ & $ β_i' α_0 $ & imp \\
         & $ i = $   & even                   & even      & any   & any          & even         & any \\
         & $ h_q = $ & $ \id_i' $             & $ \id_i $ & $ 0 $ & $ α_i' $     & $ β_i $      & $ 0 $ \\
         & $ π_q = $ & $ β_{i-1} + β_{i-1}' $ & $ α_0 $   & $ 0 $ & $ 0 $        & $ 0 $        & $ 0 $ \\
\hline
$ α_0 $ & $ μ^2 = $ & $ α_0 $ & $ β_i' $ & $ β_i' α_i $ & imp         & imp         & $ α_i α_i' $ \\
        & $ i = $   & any     & even     & odd          & any   & any   & even  \\
        & $ h_q = $ & $ 0 $   & $ 0 $    & $ β_i' $     & $ 0 $ & $ 0 $ & $ α_i' $ \\
        & $ π_q = $ & $ α_0 $ & $ 0 $    & $ 0 $        & $ 0 $ & $ 0 $ & $ 0 $ \\
\hline
$ α_i' $ & $ μ^2 = $ & $ α_i' $ & $ β_i $        & imp         & $ β_i α_i' $ & $ α_i' α_i $ & imp \\
         & $ i = $   & any      & odd            & any   & any          & any          & any \\
         & $ h_q = $ & $ 0 $    & $ 0 $          & $ 0 $ & $ 0 $        & $ 0 $        & $ 0 $ \\
         & $ π_q = $ & $ 0 $    & $ β_i + β_i' $ & $ 0 $ & $ 0 $        & $ 0 $        & $ 0 $
\end{tabular}
\caption{Multiplication table. Whenever the parity of $ i $ is specified, this refers to the parity of the index of the $ μ^2 $ result, instead of the indices of the inputs or $ h_q $ and $ π_q $ values.}
\label{tab:sphere-odd-multiplication}
\end{table}

\begin{remark}
Most values in \autoref{tab:sphere-odd-multiplication} are checked easily using \autoref{tab:sphere-odd-decomposing} and more specifically \autoref{th:sphere-odd-codifferential}. They can be grouped essentially in three types: those multiplications which always yield a particular value (with respect to $ μ^2 $, $ h_q μ^2 $ and $ π_q μ^2 $), those which vanish if $ i $ is even or odd and yield a nonzero value if $ i $ is odd respectively even, and those which involve $ α_0 $ where only close inspection proves them to vanish. The products $ μ^2 (α_0, \id_i') $ and $ μ^2 (\id_i', α_0) $ notably fall in the latter category. Let us digest this in case of $ μ^2 (α_0, \id_i') $: From the fact that $ α_0 $ lies in the counterclockwise polygon, we deduce that the source arc of $ α_0 $ is odd and therefore $ i $ is even. The result $ μ^2 (α_0, \id_i') $ then equals $ β_{i-1}' $. Now $ i-1 $ is odd, the element $ β_{i-1}' $ lies in $ R $ and we conclude $ h_q (β_{i-1}') = π_q (β_{i-1}') = 0 $. This explains the entry of $ m_2 = α_0 $ and $ m_1 = \id_i' $ in \autoref{tab:sphere-odd-multiplication}.
\end{remark}

\begin{remark}
As in the case of consistent dimers, the multiplication table merely indicates possible products, as opposed to products that actually exist. For example, the three indices “$ i $” in a product rule like $ β_i = μ^2 (\id_i, β_i) $ are not meant to denote the same index, but rather indicate the type of morphism: The first is an indexed $ β $ morphism, the second an indexed arc identity and the third again an indexed $ β $ morphism. The table merely implies that any actually existing product is of the form $ β_i = μ^2 (\id_j, β_k) $ for some combination of indices $ (i, j, k) $ allowed. Of course, we can check which combinations actually yield nonvanishing results: Those are precisely $ β_i = μ^2 (\id_{i+M+1}, β_i) $ for odd $ i $ and $ β_i = μ^2 (\id_{i+M}, β_i) $ for even $ i $. We will refer to precise combinations of indices $ (i, j, k) $ that yield nonvanishing results as \emph{precise shape} of the product. We may also refer to precise shapes when referring to $ h_q $ or $ π_q $ evaluations like $ β_i + β_i' = π_q μ^2 (\id_i, β_i) $. In any case, the precise shape is understood to link all indices involved.
\end{remark}

\begin{definition}
Kadeishvili $ h $-trees, $ π $-trees and their result components are defined as in the consistent case. In particular, a tree is supposed to have at least two leaves. The grouping rule for result components specifically reads as follows: The $ π $-tree result components $ (-1)^{\#(i+M)+1} β_i' $ and $ (-1)^{\#i} β_i $ shall be grouped together as one result component. Also, the result components $ \id_i $ shall be grouped together as one result component.
\end{definition}

We now analyze which result components are possible. As in the consistent case, we can assume that the inputs of a Kadeishvili tree do not include the identity element $ \id_L = \sum_{i = 1}^{2M} \id_i $. We analyze all possible result components of $ h $-trees first, before proceeding to $ π $-trees. Their inputs may be deformed basis elements of type $ α_0 + q α_0' $, $ β_i + β_i' $ and $ \id_i' $.

A Kadeishvili $ h $-tree is decorated by $ h_q μ_{\DefZigzagCat} $ on all its non-leaf nodes. A Kadeishvili $ π $-tree is decorated by $ h_q μ_{\DefZigzagCat} $ on all its non-leaf non-root nodes, and $ π_q μ_{\DefZigzagCat} $ on the root. Our notation $ μ^2_q $ or $ μ^{≥3}_q $ refers to the products of $ \Add\Gtl_q Q $.

As a first clue towards our analysis, we claim that $ μ^{≥3}_q $ can only be applied at the root. Indeed a disk $ μ^{≥3}_q $ can only yield $ \id_i $ and $ \id_i' $. In both cases, their $ h_q $-values vanish. Their $ π_q $-values are given by
\begin{equation*}
π_q (\id_i) = \begin{cases} \id_L & \text{ if } i = i_0 \\ 0 & \text{ else} \end{cases} \quad \text{and } \quad π_q (\id_i') = \begin{cases} \id_i' + q_{i-1} \id_{i-1}' & \text{ if } i \text{ odd} \\ 0 & \text{ else.} \end{cases}
\end{equation*}
This shows that $ μ^{≥3}_q $ can only be applied at the root. The decoration at all other nodes necessarily concerns a $ μ^2_q $.

\begin{lemma}
\label{th:sphere-odd-hclassification}
Any result component $ α_0 $ or $ β_i $ of an $ h $-tree is direct. Any result components $ β_i' $, $ α_i' $, $ \id_i' $, $ \id_i $ of $ h $-trees are derived from one of the trees in \autoref{fig:sphere-odd-betap}, \ref{fig:sphere-odd-alphap}, \ref{fig:sphere-odd-idp}, \ref{fig:sphere-odd-id}.
\end{lemma}

\begin{proof}
We start with explaining the first statement, and then delve into the second one. Our first observation is that $ α_0 $ does not appear in the multiplication table \ref{tab:sphere-odd-multiplication} at all and therefore any result component $ α_0 $ is necessarily direct. Regard now a result component $ β_i $ and assume it is non-direct. According to the multiplication table, it must be derived from a product of the form $ μ^2 (β_i', α_0) $. Let us investigate the precise shape of this product: As $ α_0 $ is located on the rear side, the index  $ i $ of the morphism $ β_i' $ is necessarily even. Therefore $ β_i' $ cannot be direct, while a glance at the multiplication table simultaneously reveals that $ β_i' $ with even $ i $ cannot be produced as a non-direct result component either. We conclude that no single result component can be derived from a product $ μ^2 (β_i', α_0) $. Finally, this means that any result component $ β_i $ is direct. This proves the first desired statement.

Regard now a result component $ β_i' $ and assume it is non-direct. We have already seen that $ i $ is necessarily odd and $ β_i' $ is derived from a product of the form $ μ^2 (α_0, β_i) $. The precise shape of this product is $ μ^2 (α_0, β_i) $ with $ i $ being equal to the index of the co-identity angle $ α_0 $ incremented by $ M $. Finally, we also realize that both $ α_0 $ and $ β_i $ are direct. The tree is depicted in \autoref{fig:sphere-odd-betap}.

Regard a result component $ α_i' $ and assume it is non-direct. According to the multiplication table, it is derived from $ μ^2 (β_i', β_i') $ or $ μ^2 (α_0, α_i') $. Let us examine both cases separately. In the first case, the precise shape of the product is $ μ^2 (β_{i+M}', β_i') $. In particular either $ i $ or $ i+M $ is even, while there are in fact no result component $ β_i' $ with even $ i $. The first case is therefore impossible. In the second case, the precise shape is $ μ^2 (α_0, α_0') $. We recall that $ α_0 $ is necessarily direct, and $ α_0' $ may either be direct or be derived from a product $ μ^2 (α_0, α_0') $ again. This gives a recursion on how $ α_0' $ is derived. Solving this recursion gives the tree in \autoref{fig:sphere-odd-alphap}.

Regard a result component $ \id_i' $ and assume it is non-direct. According to the multiplication table, it is derived from $ μ^2 (\id_i, β_i') $ or $ μ^2 (β_i', \id_i) $ or $ μ^2 (\id_i', α_i') $. Let us examine all three cases. In the first and second case, $ β_i' $ needs even index in order to have nonvanishing $ h_q $. However, we have already seen that result components $ β_i' $ all have odd index. This means there is no result component derived from the first or second case. In the third case, the precise shape of the product is $ \id_i' = h_q μ^2 (\id_{i+1}', α_i') $ and $ i $ is even. We conclude that $ \id_{i+1}' $ is direct. Even better, the result component $ α_i' $ is necessarily $ α_0' $ and is derived from one of the trees in \autoref{fig:sphere-odd-alphap}. This gives rise to the tree in \autoref{fig:sphere-odd-idp}.

Regard a result component $ \id_i $. It is necessarily non-direct, since we excluded the identity cohomology elements from the tree inputs. According to the multiplication table, it is derived from $ μ^2 (\id_i', β_i) $ or $ μ^2 (β_i', \id_i') $. Let us explore both cases. In the first case, the precise shape is $ \id_j = h_q μ^2 (\id_{i+M+1}', β_i) $ and $ i $ is odd. Note that $ j $ is free, and in fact $ h_q μ^2 (\id_{i+M+1}', β_i) $ produces many arc identities at once. Finally, due to parity both $ \id_{i+M+1}' $ and $ β_i $ are necessarily direct. This yields one tree. In the second case, the precise shape is $ α_i = μ^2 (β_{i+M}', \id_{i+1}') $ with even $ i $. We already know that $ \id_{i+1}' $ is necessarily direct. In contrast, $ β_{i+M}' $ may be either direct or derived from $ μ^2 (α_0, β_{i+M}) $. The latter case however entails that $ α_i = α_0 $, hence $ μ^2 (β_{i+M}', \id_{i+1}') = α_0 $ and $ h_q μ^2 (β_{i+M}', \id_{i+1}') = 0 $. We conclude that $ β_{i+M}' $ is necessarily direct instead. This yields one single tree. In total, both trees producing $ \id_i $ result components are depicted in \autoref{fig:sphere-odd-id}.
\end{proof}

\begin{figure}
\centering
\begin{minipage}{0.5\textwidth}
\centering
\begin{subfigure}{0.25\linewidth}
\begin{tikzpicture}
\path node {direct};
\end{tikzpicture}
\end{subfigure}
\begin{subfigure}{0.24\linewidth}
\begin{tikzpicture}
\path node (A) {$ α_0 $} node[right of=A] (B) {$ β_i $}
node[below right of=A, align=center, anchor=north] {$ h_q = β_i' $ \\ $ i $ odd} edge (A) edge (B);
\end{tikzpicture}
\end{subfigure}
\captionof{figure}{The $ β_i' $ trees}
\label{fig:sphere-odd-betap}
\end{minipage}
\begin{minipage}{0.49\linewidth}
\centering
\begin{tikzpicture}
\path node (A) {$ α_0 $} node[right of=A] (dots) {…} node[right of=dots] (B) {$ α_0 $} node[right of=B] (C) {$ α_0' $}
node[below right of=B] (D) {$ α_0' $} edge (B) edge (C)
node[below left of=D] (E) {$ h_q = α_0' $} edge[densely dotted, thick] (D) edge (A);
\path[draw, decoration=brace, decorate] (A.north) to node[midway, above, shift={(0, 0.2)}] {$ ≥ 0 $} (B.north);
\end{tikzpicture}
\captionof{figure}{The $ α_0' $ trees}
\label{fig:sphere-odd-alphap}
\end{minipage}
\end{figure}

\begin{figure}
\centering
\begin{minipage}{0.5\textwidth}
\centering
\begin{subfigure}{0.2\linewidth}
\begin{tikzpicture}
\path node {direct};
\end{tikzpicture}
\end{subfigure}
\begin{subfigure}{0.3\linewidth}
\begin{tikzpicture}
\path node (A) {$ \id_{i+1}' $} node[right of=A] (B) {$ α_0 $} node[right of=B] (dots) {…} node[right of=dots] (C) {$ α_0 $} node[right of=C] (D) {$ α_0' $}
node[below right of=C] (E) {$ α_0' $} edge (C) edge (D)
node[below left of=E] (F) {$ α_0' $} edge[densely dotted, thick] (E) edge (B)
node[below left of=F, align=center, anchor=north] {$ h_q = \id_i' $ \\ $ i = h(α_0) $ even} edge (A) edge (F);
\path[draw, decoration=brace, decorate] (B.north) to node[midway, above, shift={(0, 0.2)}] {$ ≥ 0 $} (C.north);
\end{tikzpicture}
\end{subfigure}
\captionof{figure}{The $ \id_i' $ trees}
\label{fig:sphere-odd-idp}
\end{minipage}
\end{figure}

\begin{figure}
\centering
\begin{subfigure}{0.3\linewidth}
\begin{tikzpicture}
\path node (A) {$ \id_{i+M+1}' $} node[right of=A] (B) {$ β_i $}
node[below right of=A, align=center, anchor=north] {$ h_q = \sum \id_j $ \\ $ i $ odd} edge (A) edge (B);
\end{tikzpicture}
\end{subfigure}
\begin{subfigure}{0.3\linewidth}
\begin{tikzpicture}
\path node (A) {$ β_{i+M}' $} node[right of=A] (B) {$ \id_{i+1}' $}
node[below right of=A, align=center, anchor=north] {$ h_q = \sum \id_j $ \\ $ i $ even} edge (A) edge (B);
\end{tikzpicture}
\end{subfigure}
\caption{The $ \id_i $ trees}
\label{fig:sphere-odd-id}
\end{figure}

We are now ready to approach result components of $ π $-trees.

\begin{lemma}
All $ π_q μ^2 $ result components $ \id_i' $, $ β_i' + β_i $, $ \id_i $, $ α_0 $ are derived from one of the trees in \autoref{fig:sphere-odd-piidp}, \ref{fig:sphere-odd-pibeta}, \ref{fig:sphere-odd-piid}, \ref{fig:sphere-odd-pialpha0}.
\end{lemma}

\begin{proof}
The idea is to systematically read off from the multiplication table \ref{tab:sphere-odd-multiplication} all possible ways these result components may be be derived from result components of $ h $-trees. All result components of $ h $-trees falls under the regime of \autoref{th:sphere-odd-hclassification}, allowing us to make statements on how they are derived themselves. In each case, we acquire full knowledge of the entire $ π $-tree.

Regard an $ \id_i' $ result component of a $ π $-tree. According to the multiplication table \ref{tab:sphere-odd-multiplication}, it is necessarily derived from $ μ^2 (\id_i, \id_i') $ or $ μ^2 (\id_i', \id_i) $. In the first case, the precise shape is $ \id_i' = π_q μ^2 (\id_{i+M}, \id_i') $ with $ i $ odd. We realize that $ \id_i' $ is necessarily direct, while $ \id_{i+M} $ may come from two possible trees. In the second case, the precise shape is $ \id_i' = π_q μ^2 (\id_i', \id_i) $ with $ i $ odd. We realize that $ \id_i' $ is direct, while $ \id_i $ may again come from two possible trees. In total, the four possible trees are depicted in \autoref{fig:sphere-odd-piidp}.

Regard an $ β_i' + β_i $ result component. According to the multiplication table, it is derived from $ μ^2 (\id_i, β_i) $ or $ μ^2 (\id_i, β_i') $ or $ μ^2 (\id_i', α_0') $ or $ μ^2 (β_i, \id_i) $ or $ μ^2 (β_i', \id_i) $ or $ μ^2 (α_0', \id_i') $. Let us investigate all these six cases. In the first case, the precise shape is $ β_i' + β_i = π_q μ^2 (\id_{i+M+1}, β_i) $ with $ i $ odd. We realize that $ β_i $ is direct and $ \id_{i+M+1} $ may come from two possible trees. In the second case, the precise shape is $ β_{i-1}' + β_{i-1} = π_q μ^2 (\id_{i+M}, β_i') $ with $ i $ even. Since no result component $ β_i' $ with even $ i $ exists, this case is impossible. In the third case, the precise shape is $ β_{i-1}' + β_{i-1} = π_q μ^2 (\id_{i+1}', α_0') $ with $ i $ even. We realize that $ \id_{i+1}' $ is direct and $ α_0' $ comes from the known standard type of tree. In the fourth case, the precise shape is $ β_i' + β_i = π_q μ^2 (β_i, \id_i) $ with $ i $ odd. We realize that $ β_i $ is direct and $ \id_i $ may come from two possible trees. In the fifth case, the precise shape is $ β_{i-1}' + β_{i-1} = π_q μ^2 (β_i', \id_i) $ with $ i $ even. Since no result component $ β_i' $ with odd $ i $ exists, this case is impossible. In the sixth case, the precise shape is $ β_{i+M}' + β_{i+M} = π_q μ^2 (α_0', \id_{i+M}') $ with $ i = h(α_0) $ even. We realize that $ \id_{i+M}' $ is direct and $ α_0' $ comes from the known tree. In total, all six trees are depicted in \autoref{fig:sphere-odd-pibeta}.

Regard an $ \id_i $ result component. According to the multiplication table, it is derived from $ μ^2 (\id_i, \id_i) $ or $ μ^2 (\id_i', \id_i') $. In the first case, the precise shape is $ \id_{i_0} = π_q μ^2 (\id_{i_0}, \id_{i_0}) $. Since the arc identity $ \id_{i_0} $ never appears as result component of an $ h $-tree, this case is however vacuous. In the second case, the precise shape is $ \id_{i_0} = π_q μ^2 (\id_{i_0 + M}', \id_{i_0}') $. This tree is depicted in \autoref{fig:sphere-odd-piid}.

Regard an $ α_0 $ result component. According to the multiplication table, it is derived from $ μ^2 (\id_i, α_0) $ or $ μ^2 (α_0, \id_i) $ or $ μ^2 (\id_i', β_i) $ or $ μ^2 (β_i', \id_i') $. Let us investigate all four cases. In the first case, the precise shape is $ α_0 = π_q μ^2 (\id_i, α_0) $ with $ i = h(α_0) $ even. We realize that $ α_0 $ is direct, while $ \id_i $ may come from two possible trees. In the second case, the precise shape is $ α_0 = π_q μ^2 (α_0, \id_{i+1}) $ with $ i = h(α_0) $ even. We realize that $ α_0 $ is direct, while $ \id_{i+1} $ may come from two possible trees. In the third case, the precise shape is $ α_0 = π_q μ^2 (\id_{i+M}', β_i) $ with $ i $ even or $ α_0 = π_q μ^2 (\id_{i+M+1}', β_i) $ with $ i $ odd. The former case is impossible, since there is no result component $ β_i $ with $ i $ even. In the latter case, both $ \id_{i+M+1}' $ and $ β_i $ are direct. In the fourth case, the precise shape is $ α_0 = π_q μ^2 (β_{i+M}', \id_{i+1}') $ with $ i $ even. We realize that $ \id_{i+1}' $ is direct, while $ β_{i+M}' $ may be direct or derived from $ μ^2 (α_0, β_{i+M}) $. We recall that in the latter case it is necessary that $ i = h(α_0) $. The total collection of seven trees is depicted in \autoref{fig:sphere-odd-pialpha0}.
\end{proof}

\begin{figure}
\centering
\begin{subfigure}{0.3\linewidth}
\begin{tikzpicture}
\path node (A) {$ \id_{j+M+1}' $} node[right of=A] (B) {$ β_j $} node[right of=B] (C) {$ \id_i' $}
node[below right of=A] (D) {$ \id_{i+M} $} edge (A) edge (B)
node[below right of=D, align=center, anchor=north] {$ π_q = \id_i' $ \\ $ i $ odd, $ j $ odd} edge (D) edge (C);
\end{tikzpicture}
\end{subfigure}
\begin{subfigure}{0.3\linewidth}
\begin{tikzpicture}
\path node (A) {$ β_{j+M}' $} node[right of=A] (B) {$ \id_{i+1}' $} node[right of=B] (C) {$ \id_i' $}
node[below right of=A] (D) {$ \id_{i+M} $} edge (A) edge (B)
node[below right of=D, align=center, anchor=north] {$ π_q = \id_i' $ \\ $ i $ odd, $ j $ even} edge (D) edge (C);
\end{tikzpicture}
\end{subfigure}
\begin{subfigure}{0.3\linewidth}
\begin{tikzpicture}
\path node (A) {$ \id_i' $} node[right of=A] (B) {$ \id_{j+M+1}' $} node[right of=B] (C) {$ β_j $}
node[below right of=B] (D) {$ \id_i $} edge (B) edge (C)
node[below left of=D, align=center, anchor=north] {$ π_q = \id_i' $ \\ $ i $ odd, $ j $ odd} edge (A) edge (D);
\end{tikzpicture}
\end{subfigure}
\begin{subfigure}{0.3\linewidth}
\begin{tikzpicture}
\path node (A) {$ \id_i' $} node[right of=A] (B) {$ β_{j+M}' $} node[right of=B] (C) {$ \id_{j+1}' $}
node[below right of=B] (D) {$ \id_i $} edge (B) edge (C)
node[below left of=D, align=center, anchor=north] {$ π_q = \id_i' $ \\ $ i $ odd, $ j $ even} edge (A) edge (D);
\end{tikzpicture}
\end{subfigure}
\caption{π-trees for $ \id_i' $}
\label{fig:sphere-odd-piidp}
\end{figure}

\begin{figure}
\centering
\begin{subfigure}{0.3\linewidth}
\begin{tikzpicture}
\path node (A) {$ \id_{j+M+1}' $} node[right of=A] (B) {$ β_j $} node[right of=B] (C) {$ β_i $}
node[below right of=A] (D) {$ \id_{i+M+1} $} edge (A) edge (B)
node[below right of=D, align=center, anchor=north] {$ π_q = β_i' + β_i $ \\ $ i $ odd, $ j $ odd} edge (D) edge (C);
\end{tikzpicture}
\end{subfigure}
\begin{subfigure}{0.3\linewidth}
\begin{tikzpicture}
\path node (A) {$ β_{j+M}' $} node[right of=A] (B) {$ \id_{j+1}' $} node[right of=B] (C) {$ β_i $}
node[below right of=A] (D) {$ \id_{i+M+1} $} edge (A) edge (B)
node[below right of=D, align=center, anchor=north] {$ π_q = β_i' + β_i $ \\ $ i $ odd, $ j $ even} edge (D) edge (C);
\end{tikzpicture}
\end{subfigure}
\begin{subfigure}{0.3\linewidth}
\begin{tikzpicture}
\path node (A) {$ \id_{i+1}' $} node[right of=A] (B) {$ α_0 $} node[right of=B] (dots) {…} node[right of=dots] (C) {$ α_0 $} node[right of=C] (D) {$ α_0' $}
node[below right of=C] (E) {$ α_0' $} edge (C) edge (D)
node[below left of=E] (F) {$ α_0' $} edge[densely dotted, thick] (E) edge (B)
node[below left of=F, align=center, anchor=north] {$ π_q = β_{i-1}' + β_{i-1} $ \\ $ i = h(α_0) $ even} edge (F) edge (A);
\path[draw, decoration=brace, decorate] (B.north) to node[midway, above, shift={(0, 0.2)}] {$ ≥ 0 $} (C.north);
\end{tikzpicture}
\end{subfigure}
\begin{subfigure}{0.3\linewidth}
\begin{tikzpicture}
\path node (A) {$ β_i $} node[right of=A] (B) {$ \id_{j+M+1}' $} node[right of=B] (C) {$ β_j $}
node[below right of=B] (D) {$ \id_i $} edge (B) edge (C)
node[below left of=D, align=center, anchor=north] {$ π_q = β_i' + β_i $ \\ $ i $ odd, $ j $ odd} edge (D) edge (A);
\end{tikzpicture}
\end{subfigure}
\begin{subfigure}{0.3\linewidth}
\begin{tikzpicture}
\path node (A) {$ β_i $} node[right of=A] (B) {$ β_{j+M}' $} node[right of=B] (C) {$ \id_{j+1}' $}
node[below right of=B] (D) {$ \id_i $} edge (B) edge (C)
node[below left of=D, align=center, anchor=north] {$ π_q = β_i' + β_i $ \\ $ i $ odd, $ j $ even} edge (D) edge (A);
\end{tikzpicture}
\end{subfigure}
\begin{subfigure}{0.3\linewidth}
\begin{tikzpicture}
\path node (A) {$ α_0 $} node[right of=A] (dots) {…} node[right of=dots] (B) {$ α_0 $} node[right of=B] (C) {$ α_0' $} node[right of=C] (D) {$ \id_{i+M}' $}
node[below right of=B] (E) {$ α_0' $} edge (B) edge (C)
node[below left of=E] (F) {$ α_0' $} edge[densely dotted, thick] (E) edge (A)
node[below right of=F, align=center, anchor=north] {$ π_q = β_{i+M}' + β_{i+M} $ \\ $ i = h(α_0) $ even} edge (F) edge (D);
\path[draw, decoration=brace, decorate] (A.north) to node[midway, above, shift={(0, 0.2)}] {$ ≥ 0 $} (B.north);
\end{tikzpicture}
\end{subfigure}
\caption{π-trees for $ β_i' + β_i $}
\label{fig:sphere-odd-pibeta}
\end{figure}

\begin{figure}
\centering
\begin{subfigure}{0.3\linewidth}
\begin{tikzpicture}
\path node (A) {$ \id_{i_0 + M}' $} node[right of=A] (B) {$ \id_{i_0}' $}
node[below right of=A] {$ π_q = \id_L $} edge (A) edge (B);
\end{tikzpicture}
\end{subfigure}
\caption{π-trees for $ \id_L = \sum \id_i $}
\label{fig:sphere-odd-piid}
\end{figure}

\begin{figure}
\centering
\begin{subfigure}{0.3\linewidth}
\begin{tikzpicture}
\path node (A) {$ \id_{j+M+1}' $} node[right of=A] (B) {$ β_j $} node[right of=B] (C) {$ α_0 $}
node[below right of=A] (D) {$ \id_i $} edge (A) edge (B)
node[below right of=D, align=center, anchor=north] {$ π_q = α_0 $ \\ $ i = h(α_0) $ even, $ j $ odd} edge (D) edge (C);
\end{tikzpicture}
\end{subfigure}
\begin{subfigure}{0.3\linewidth}
\begin{tikzpicture}
\path node (A) {$ β_{j+M}' $} node[right of=A] (B) {$ \id_{j+1}' $} node[right of=B] (C) {$ α_0 $}
node[below right of=A] (D) {$ \id_i $} edge (A) edge (B)
node[below right of=D, align=center, anchor=north] {$ π_q = α_0 $ \\ $ i = h(α_0) $ even, $ j $ even} edge (D) edge (C);
\end{tikzpicture}
\end{subfigure}
\begin{subfigure}{0.3\linewidth}
\begin{tikzpicture}
\path node (A) {$ α_0 $} node[right of=A] (B) {$ \id_{j+M+1}' $} node[right of=B] (C) {$ β_j $}
node[below right of=B] (D) {$ \id_{i+1} $} edge (B) edge (C)
node[below left of=D, align=center, anchor=north] {$ π_q = α_0 $ \\ $ i = h(α_0) $ even, $ j $ odd} edge (D) edge (A);
\end{tikzpicture}
\end{subfigure}
\begin{subfigure}{0.3\linewidth}
\begin{tikzpicture}
\path node (A) {$ α_0 $} node[right of=A] (B) {$ β_{j+M}' $} node[right of=B] (C) {$ \id_{j+1}' $}
node[below right of=B] (D) {$ \id_{i+1} $} edge (B) edge (C)
node[below left of=D, align=center, anchor=north] {$ π_q = α_0 $ \\ $ i = h(α_0) $ even, $ j $ odd} edge (D) edge (A);
\end{tikzpicture}
\end{subfigure}
\begin{subfigure}{0.3\linewidth}
\begin{tikzpicture}
\path node (A) {$ \id_{i+M+1}' $} node[right of=A] (B) {$ β_i $}
node[below right of=A, align=center, anchor=north] {$ π_q = α_0 $ \\ $ i $ odd} edge (A) edge (B);
\end{tikzpicture}
\end{subfigure}
\begin{subfigure}{0.3\linewidth}
\begin{tikzpicture}
\path node (A) {$ β_{i+M}' $} node[right of=A] (B) {$ \id_{i+1}' $}
node[below right of=A, align=center, anchor=north] {$ π_q = α_0 $ \\ $ i $ even} edge (A) edge (B);
\end{tikzpicture}
\end{subfigure}
%
\begin{subfigure}{0.3\linewidth}
\begin{tikzpicture}
\path node (A) {$ α_0 $} node[right of=A] (B) {$ β_{i+M} $} node[right of=B] (C) {$ \id_{i+1}' $}
node[below right of=A] (D) {$ β_{i+M}' $} edge (A) edge (B)
node[below right of=D, align=center, anchor=north] {$ h_q = \sum \id_j $ \\ $ i = h(α_0) $ even} edge (D) edge (C);
\end{tikzpicture}
\end{subfigure}
\caption{π-trees for $ α_0 $}
\label{fig:sphere-odd-pialpha0}
\end{figure}

\newcommand{\SphereOddPiDiskDraw}{%
\path[draw, ->] (0, 0) -- ++(up:1) coordinate[pos=0.7] (7) to[bend right=90, looseness=2] ++(up:0.4) -- ++(up:1) coordinate[pos=0.3] (5) coordinate[pos=0.7] (3) to[bend left=90, looseness=2] ++(up:0.4) -- ++(up:1) coordinate[pos=0.3] (1);
\path[draw, dashed] (0.3, 0) -- ++(up:1) coordinate[pos=0.7] (8) to[bend left=90, looseness=2] ++(up:0.4) -- ++(up:1) coordinate[pos=0.3] (6) coordinate[pos=0.7] (4) to[bend right=90, looseness=2] ++(up:0.4) -- ++(up:1) coordinate[pos=0.3] (2);
\path (-0.5, 3.5) node[anchor=south] {\small rear} (0.8, 3.5) node[anchor=south] {\small front};
\path[draw] ($ (3)!0.4!(5) + (-0.2, 0) $) -- ($ (4)!0.4!(6) + (-0.1, 0) $);
\path[draw] ($ (3)!0.6!(5) + (-0.2, 0) $) -- ($ (4)!0.6!(6) + (-0.1, 0) $);
}

\newcommand{\SphereOddPiDiskDrawAlt}{%
\path[draw, ->] (0, 0) -- ++(up:1) coordinate[pos=0.7] (7) to[bend left=90, looseness=2] ++(up:0.4) -- ++(up:1) coordinate[pos=0.3] (5) coordinate[pos=0.7] (3) to[bend right=90, looseness=2] ++(up:0.4) -- ++(up:1) coordinate[pos=0.3] (1);
\path[draw, dashed] (0.3, 0) -- ++(up:1) coordinate[pos=0.7] (8) to[bend right=90, looseness=2] ++(up:0.4) -- ++(up:1) coordinate[pos=0.3] (6) coordinate[pos=0.7] (4) to[bend left=90, looseness=2] ++(up:0.4) -- ++(up:1) coordinate[pos=0.3] (2);
\path (-0.5, 3.5) node[anchor=south] {\small rear} (1, 3.5) node[anchor=south] {\small front};
\path[draw] ($ (3)!0.4!(5) + (-0.2, 0) $) -- ($ (4)!0.4!(6) + (-0.1, 0) $);
\path[draw] ($ (3)!0.6!(5) + (-0.2, 0) $) -- ($ (4)!0.6!(6) + (-0.1, 0) $);
}

\newcommand{\SphereOddPiBorderDraw}{%
\path[draw, thick] (-1, 0) to[out=90, in=210] (0, 1) to[out=30, in=270] (1, 2) to[out=90, in=330] coordinate[pos=0.7] (4) (0, 3) to[out=150, in=270] coordinate[pos=0.3] (1) (-1, 4) coordinate (5) to[out=90, in=210] (0, 5) to[out=30, in=270] (1, 6);
\path[draw, thick] (1, 0) to[out=90, in=330] (0, 1) to[out=150, in=270] (-1, 2) coordinate (6) to[out=90, in=210] coordinate[pos=0.7] (2) (0, 3) to[out=30, in=270] coordinate[pos=0.3] (3) (1, 4) to[out=90, in=330] (0, 5) to[out=150, in=270] (-1, 6);
}

\begin{figure}
\centering
\begin{subfigure}{0.31\linewidth}
\centering
\begin{tikzpicture}
\SphereOddPiDiskDraw
\path[draw, ->, bend right=80, looseness=1.5] (1) to node[midway, left] {$ α_0 / α_i $} (3);
\path[draw, ->, bend right=80, looseness=1.5] (5) to node[midway, left] {$ α_i' $} (7);
\end{tikzpicture}
\begin{tikzpicture}[scale=0.7]
\SphereOddPiBorderDraw
\path[draw, thick] (0, 0.7) to[out=30, in=270] (1.2, 2) to[out=90, in=30] (1) to[out=220, in=210] (5) to[out=40, in=210] (0, 4.7);
\path[draw, decoration={border, amplitude=4pt, angle=90, segment length=2pt}, color=gray, decorate] (0, 1) to[out=30, in=270] (1, 2) to[out=90, in=330] (0, 3) to[out=150, in=310] (1) to[out=220, in=210] (5) to[out=90, in=210] (0, 5);
\path[fill] (1) circle[radius=0.07] (5) circle[radius=0.07];
\path (1) node[above] {\small out} (5) node[right] {\small $ α_0 $};
\end{tikzpicture}
\end{subfigure}
\begin{subfigure}{0.31\linewidth}
\centering
\begin{tikzpicture}
\SphereOddPiDiskDraw
\path[draw, ->, bend right=80, looseness=1.5] (2) to node[midway, left] {$ β_i' $} (3);
\path[draw, ->, bend right=80, looseness=1.5] (5) to node[midway, left] {$ α_i' $} (7);
\end{tikzpicture}
\begin{tikzpicture}[scale=0.7]
\SphereOddPiBorderDraw
\path[draw, thick] (0, 0.7) to[out=30, in=270] (1.2, 2) to[out=90, in=0] (1) to[out=180, in=270] (-1.2, 4) to[out=90, in=210] coordinate[pos=0.83] (cross) (0, 5.3);
\path[draw, decoration={border, amplitude=4pt, angle=90, segment length=2pt}, color=gray, decorate] (1, 2) to[out=90, in=330] (0, 3) to[out=150, in=0] (1) to[out=180, in=270] (-1.2, 4) to[out=90, in=210] (cross) to[out=150, in=270] (-1, 6);
\path[fill] (1) circle[radius=0.07] (cross) circle[radius=0.07];
\path (1) node[above] {\small out} (cross) node[above] {\small $ β_i' $};
\end{tikzpicture}
\end{subfigure}
\begin{subfigure}{0.31\linewidth}
\centering
\begin{tikzpicture}
\SphereOddPiDiskDraw
\path[draw, ->, bend right=80, looseness=1.5] (1) to node[midway, left] {$ α_0 / α_i $} (3);
\path[draw, ->, bend right=80, looseness=1.5] (5) to node[midway, left] {$ β_i' $} (8);
\end{tikzpicture}
\begin{tikzpicture}[scale=0.7]
\SphereOddPiBorderDraw
\path[draw, thick] (0, 0.7) to[out=30, in=270] (1.2, 2) to[out=90, in=30] (1) to[out=220, in=210] (5) to[out=40, in=210] (0, 4.7);
\path[draw, decoration={border, amplitude=4pt, angle=90, segment length=2pt}, color=gray, decorate] (0, 1) to[out=150, in=270] (-1, 2) to[out=90, in=210] (0, 3) to[out=150, in=310] (1) to[out=220, in=210] (5) to[out=90, in=210] (0, 5);
\path[fill] (1) circle[radius=0.07] (5) circle[radius=0.07] (0, 3) circle[radius=0.07];
\path (1) node[above] {\small out} (5) node[right] {\small $ α_0 $} (0, 3) node[below] {\small $ β_i' $};
\end{tikzpicture}
\end{subfigure}
\begin{subfigure}{0.31\linewidth}
\centering
\begin{tikzpicture}
\SphereOddPiDiskDraw
\path[draw, ->, bend right=80, looseness=1.5] (2) to node[midway, left] {$ β_i' $} (3);
\path[draw, ->, bend right=80, looseness=1.5] (5) to node[midway, left] {$ β_i' $} (8);
\end{tikzpicture}
\begin{tikzpicture}[scale=0.7]
\SphereOddPiBorderDraw
\path[draw, thick] (0, 0.7) to[out=30, in=270] (1.2, 2) to[out=90, in=0] (1) to[out=180, in=270] (-1.2, 4) to[out=90, in=210] coordinate[pos=0.83] (cross) (0, 5.3);
\path[draw, decoration={border, amplitude=4pt, angle=90, segment length=2pt}, color=gray, decorate] (0, 1) to[out=150, in=270] (-1, 2) to[out=90, in=210] (0, 3) to[out=150, in=310] (1) to[out=180, in=270] (-1.2, 4) to[out=90, in=210] (cross) to[out=150, in=270] (-1, 6);
\path[fill] (1) circle[radius=0.07] (cross) circle[radius=0.07] (0, 3) circle[radius=0.07];
\path (1) node[above] {\small out} (cross) node[above] {\small $ β_i' $} (0, 3) node[below] {$ β_i' $};
\end{tikzpicture}
\end{subfigure}
\begin{subfigure}{0.31\linewidth}
\centering
\begin{tikzpicture}
\SphereOddPiDiskDrawAlt
\path[draw, ->, bend right=80, looseness=1.5] (1) to node[midway, left] {$ α_i' $} (3);
\path[draw, ->, bend right=80, looseness=1.5] (5) to node[midway, left] {$ α_0 / α_i $} (7);
\end{tikzpicture}
\begin{tikzpicture}[scale=0.7]
\SphereOddPiBorderDraw
\path[draw, thick] (0, 0.7) to[out=150, in=220] (6) to[out=40, in=300] (2) to[out=120, in=270] (0.8, 4) to[out=90, in=330] (0, 4.7);
\path[draw, decoration={border, amplitude=4pt, angle=90, segment length=2pt}, color=gray, decorate] (0, 0.7) to[out=150, in=220] (6) to[out=90, in=230] (2) to[out=120, in=270] (0.8, 4) to[out=90, in=330] (0, 4.7);
\path[fill] (2) circle[radius=0.07] (6) circle[radius=0.07];
\path (2) node[right] {\small out} (6) node[right] {$ α_0 $};
\end{tikzpicture}
\end{subfigure}
\begin{subfigure}{0.31\linewidth}
\centering
\begin{tikzpicture}
\SphereOddPiDiskDrawAlt
\path[draw, ->, bend right=80, looseness=1.5] (3) to node[midway, right] {$ α_i $} (1);
\path[draw, ->, bend right=80, looseness=1.5] (7) to node[midway, right] {$ α_0'/α_i' $} (5);
\end{tikzpicture}
\begin{tikzpicture}[scale=0.7]
\SphereOddPiBorderDraw
\path[draw, thick] (0, 0.7) to[out=150, in=220] (6) to[out=40, in=300] (2) to[out=120, in=270] (0.8, 4) to[out=90, in=330] (0, 4.7);
\path[draw, decoration={border, amplitude=4pt, angle=90, segment length=2pt, mirror}, color=gray, decorate] (0, 1) to[out=150, in=270] (6) to[out=40, in=300] (2) to[out=30, in=210] (0, 3) to[out=30, in=270] (1, 4) to[out=90, in=330] (0, 5);
\path[fill] (2) circle[radius=0.07] (6) circle[radius=0.07];
\path (2) node[left] {\small out} (6) node[left] {\small $ α_0 $};
\end{tikzpicture}
\end{subfigure}
\begin{subfigure}{0.31\linewidth}
\centering
\begin{tikzpicture}
\SphereOddPiDiskDrawAlt
\path[draw, ->, bend right=80, looseness=1.5] (3) to node[midway, right] {$ α_i $} (1);
\path[draw, ->, bend right=80, looseness=1.5] (8) to node[midway, right] {$ β_i $} (5);
\end{tikzpicture}
\begin{tikzpicture}[scale=0.7]
\SphereOddPiBorderDraw
\path[draw, thick] (0, 1.3) to[out=150, in=270] (-0.8, 2) to[out=90, in=240] coordinate[pos=0.8] (cross) (3) to[out=60, in=270] (0.8, 4) to[out=90, in=330] (0, 4.7);
\path[draw, decoration={border, amplitude=4pt, angle=90, segment length=2pt, mirror}, color=gray, decorate] (0, 1) to[out=30, in=270] (1, 2) to[out=90, in=330] (cross) to[out=50, in=240] (3) to[out=30, in=270] (1, 4) to[out=90, in=330] (0, 5);
\path[fill] (cross) circle[radius=0.07] (3) circle[radius=0.07];
\path (cross) node[below] {\small $ β_i $} (3) node[left, shift={(0, 0.1)}] {\small out};
\end{tikzpicture}
\end{subfigure}
\begin{subfigure}{0.31\linewidth}
\centering
\begin{tikzpicture}
\SphereOddPiDiskDrawAlt
\path[draw, ->, bend right=80, looseness=1.5] (3) to node[midway, right] {$ β_i $} (2);
\path[draw, ->, bend right=80, looseness=1.5] (7) to node[midway, right] {$ α_0'/α_i' $} (5);
\end{tikzpicture}
\begin{tikzpicture}[scale=0.7]
\SphereOddPiBorderDraw
\path[draw, thick] (0, 0.7) to[out=150, in=220] (6) to[out=40, in=300] (2) to[out=120, in=270] (0.8, 4) to[out=90, in=330] (0, 4.7);
\path[draw, decoration={border, amplitude=4pt, angle=90, segment length=2pt, mirror}, color=gray, decorate] (0, 1) to[out=150, in=270] (6) to[out=40, in=300] (2) to[out=30, in=210] (0, 3) to[out=30, in=270] (1, 4) to[out=90, in=330] (0, 5) to[out=30, in=270] (1, 6);
\path[fill] (2) circle[radius=0.07] (6) circle[radius=0.07] (0, 5) circle[radius=0.07];
\path (2) node[left] {\small out} (6) node[left] {\small $ α_0 $} (0, 5) node[above] {\small $ β_i $};
\end{tikzpicture}
\end{subfigure}
\begin{subfigure}{0.31\linewidth}
\centering
\begin{tikzpicture}
\SphereOddPiDiskDrawAlt
\path[draw, ->, bend right=80, looseness=1.5] (3) to node[midway, right] {$ β_i $} (2);
\path[draw, ->, bend right=80, looseness=1.5] (8) to node[midway, right] {$ β_i $} (5);
\end{tikzpicture}
\begin{tikzpicture}[scale=0.7]
\SphereOddPiBorderDraw
\path[draw, thick] (0, 1.3) to[out=150, in=270] (-0.8, 2) to[out=90, in=240] coordinate[pos=0.78] (cross) (3) to[out=60, in=270] (0.8, 4) to[out=90, in=330] (0, 4.7);
\path[draw, decoration={border, amplitude=4pt, angle=90, segment length=2pt, mirror}, color=gray, decorate] (0, 1) to[out=30, in=270] (1, 2) to[out=90, in=330] (cross) to[out=50, in=240] (3) to[out=50, in=270] (1, 4) to[out=90, in=330] (0, 5) to[out=30, in=270] (1, 6);
\path[fill] (cross) circle[radius=0.07] (3) circle[radius=0.07] (0, 5) circle[radius=0.07];
\path (cross) node[below] {\small $ β_i $} (3) node[left, shift={(0, 0.1)}] {\small out} (0, 5) node[above] {\small $ β_i $};
\end{tikzpicture}
\end{subfigure}
\begin{subfigure}{0.31\linewidth}
\centering
\begin{tikzpicture}
\SphereOddPiDiskDraw
\path[draw, ->, bend right=80, looseness=1.5] (3) to node[midway, right] {$ α_0'/α_i' $} (1);
\path[draw, ->, bend right=80, looseness=1.5] (7) to node[midway, right] {$ α_i $} (5);
\end{tikzpicture}
\begin{tikzpicture}[scale=0.7]
\SphereOddPiBorderDraw
\path[draw, thick] (0, 0.7) to[out=30, in=270] (1.2, 2) to[out=90, in=30] (1) to[out=220, in=210] (5) to[out=40, in=210] (0, 4.7);
\path[draw, decoration={border, amplitude=4pt, angle=90, segment length=2pt, mirror}, color=gray, decorate] (0, 0.7) to[out=30, in=270] (1.2, 2) to[out=90, in=30] (1) to[out=130, in=270] (5) to[out=40, in=210] (0, 4.7);
\path[fill] (1) circle[radius=0.07] (5) circle[radius=0.07];
\path (1) node[below, shift={(0, 0.1)}] {\small out} (5) node[left] {\small $ α_0' $};
\end{tikzpicture}
\end{subfigure}
\caption{$ \id_L $ disk result components and their subdisks. The first five disks lie on the rear side and the second five on the front side. The first and final arc is supposed to be the identity arc $ i_0 $ of $ L $ and is highlighted by a crossing double line. This double line also indicates the separation between the first and the final angle of the disk.}
\label{fig:sphere-resultcomp-pidiskid}
\end{figure}
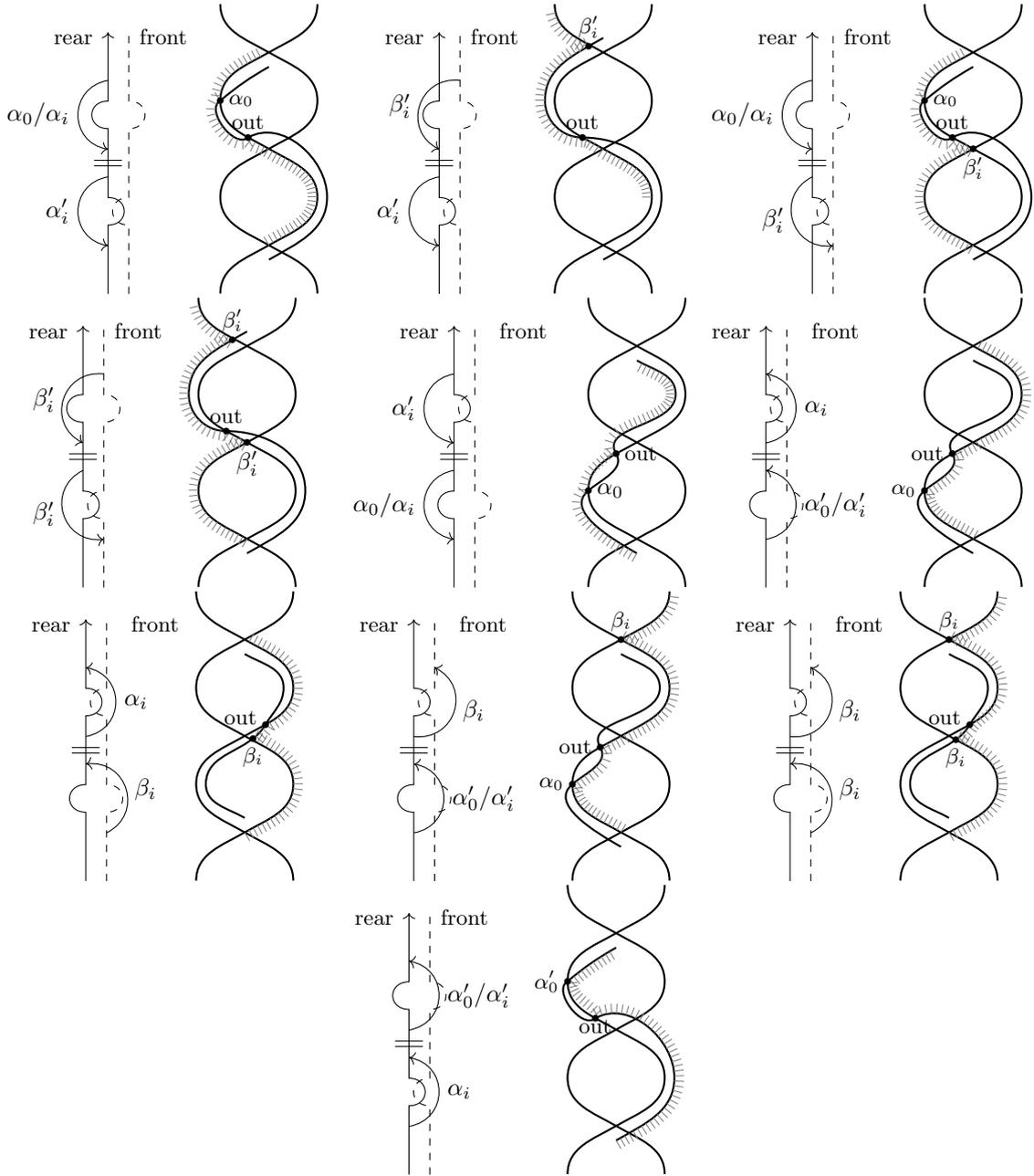

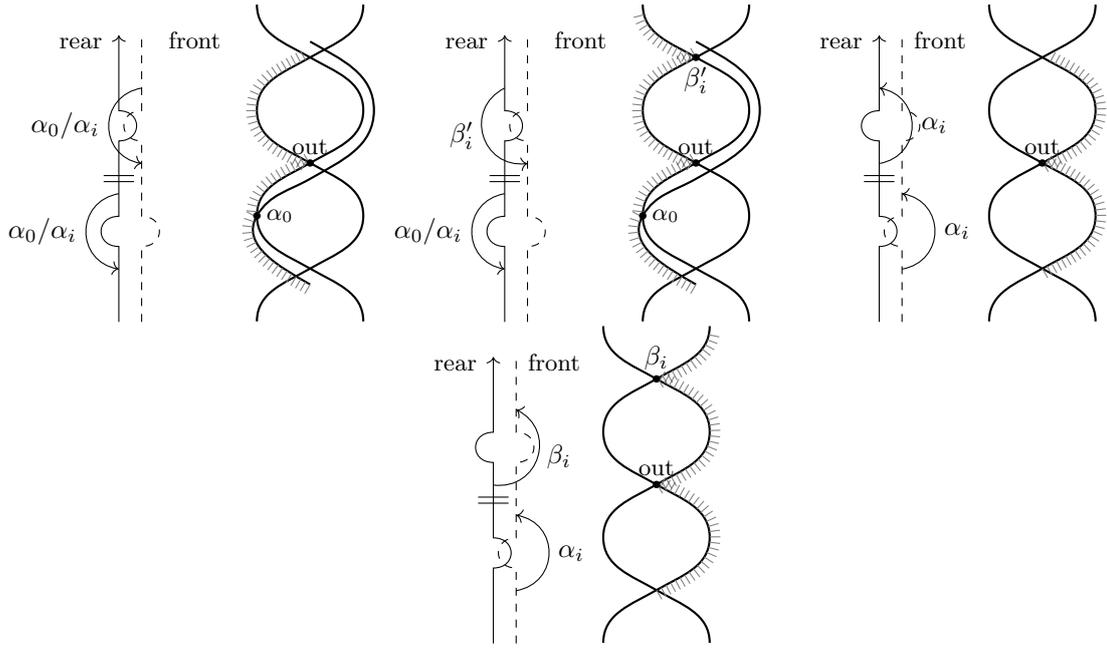
\begin{figure}
\centering
\begin{subfigure}{0.31\linewidth}
\centering
\begin{tikzpicture}
\SphereOddPiDiskDrawAlt
\path[draw, ->, bend right=80, looseness=1.5] (2) to node[midway, left] {$ α_0 / α_i $} (4);
\path[draw, ->, bend right=80, looseness=1.5] (5) to node[midway, left] {$ α_0 / α_i $} (7);
\end{tikzpicture}
\begin{tikzpicture}[scale=0.7]
\SphereOddPiBorderDraw
\path[draw, thick] (0, 0.7) to[out=150, in=240] (6) to[out=60, in=210] (0, 2.7) to[out=30, in=270] (1.2, 4) to[out=90, in=330] (0, 5.3);
\path[draw, decoration={border, amplitude=4pt, angle=90, segment length=2pt}, color=gray, decorate] (0, 0.7) to[out=150, in=240] (6) to[out=90, in=210] (0, 3) to[out=150, in=270] (-1, 4) to[out=90, in=210] (0, 5);
\path[fill] (0, 3) circle[radius=0.07] (6) circle[radius=0.07];
\path (0, 3) node[above] {\small out} (6) node[right] {\small $ α_0 $};
\end{tikzpicture}
\end{subfigure}
\begin{subfigure}{0.31\linewidth}
\centering
\begin{tikzpicture}
\SphereOddPiDiskDrawAlt
\path[draw, ->, bend right=80, looseness=1.5] (1) to node[midway, left] {$ β_i' $} (4);
\path[draw, ->, bend right=80, looseness=1.5] (5) to node[midway, left] {$ α_0 / α_i $} (7);
\end{tikzpicture}
\begin{tikzpicture}[scale=0.7]
\SphereOddPiBorderDraw
\path[draw, thick] (0, 0.7) to[out=150, in=240] (6) to[out=60, in=210] (0, 2.7) to[out=30, in=270] (1.2, 4) to[out=90, in=330] (0, 5.3);
\path[draw, decoration={border, amplitude=4pt, angle=90, segment length=2pt}, color=gray, decorate] (0, 0.7) to[out=150, in=240] (6) to[out=90, in=210] (0, 3) to[out=150, in=270] (-1, 4) to[out=90, in=210] (0, 5) to[out=150, in=270] (-1, 6);
\path[fill] (0, 3) circle[radius=0.07] (6) circle[radius=0.07] (0, 5) circle[radius=0.07];
\path (0, 3) node[above] {\small out} (6) node[right] {\small $ α_0 $} (0, 5) node[below] {$ β_i' $};
\end{tikzpicture}
\end{subfigure}
\begin{subfigure}{0.31\linewidth}
\centering
\begin{tikzpicture}
\SphereOddPiDiskDraw
\path[draw, ->, bend right=80, looseness=1.5] (3) to node[midway, right] {$ α_i $} (1);
\path[draw, ->, bend right=80, looseness=1.5] (8) to node[midway, right] {$ α_i $} (6);
\end{tikzpicture}
\begin{tikzpicture}[scale=0.7]
\SphereOddPiBorderDraw
\path[draw, decoration={border, amplitude=4pt, angle=90, segment length=2pt, mirror}, color=gray, decorate] (0, 1) to[out=30, in=270] (1, 2) to[out=90, in=330] (0, 3) to[out=30, in=270] (1, 4) to[out=90, in=330] (0, 5);
\path[fill] (0, 3) circle[radius=0.07];
\path (0, 3) node[above] {\small out};
\end{tikzpicture}
\end{subfigure}
\begin{subfigure}{0.31\linewidth}
\centering
\begin{tikzpicture}
\SphereOddPiDiskDraw
\path[draw, ->, bend right=80, looseness=1.5] (3) to node[midway, right] {$ β_i $} (2);
\path[draw, ->, bend right=80, looseness=1.5] (8) to node[midway, right] {$ α_i $} (6);
\end{tikzpicture}
\begin{tikzpicture}[scale=0.7]
\SphereOddPiBorderDraw
\path[draw, decoration={border, amplitude=4pt, angle=90, segment length=2pt, mirror}, color=gray, decorate] (0, 1) to[out=30, in=270] (1, 2) to[out=90, in=330] (0, 3) to[out=30, in=270] (1, 4) to[out=90, in=330] (0, 5) to[out=30, in=270] (1, 6);
\path[fill] (0, 3) circle[radius=0.07] (0, 5) circle[radius=0.07];
\path (0, 3) node[above] {\small out} (0, 5) node[above] {$ β_i $};
\end{tikzpicture}
\end{subfigure}
\caption{$ \id_i' $ disk result components of π-trees and how to draw their subdisks}
\label{fig:sphere-resultcomp-pidiskidp}
\end{figure}

We now introduce the suitable version of CR, ID, DS and DW disks.

\begin{definition}
CR, ID, DS and DW disks are defined as in the case of geometrically consistent dimers. More specifically, the definitions read as follows: A \emph{CR disk} is an SL disk where all inputs and the output lie apart, with the exception that stacks of $ α_0 $ inputs are allowed if the SL disk lies on the front side. An \emph{ID disk} is an SL disk where all inputs and the output lie apart, with the exception that stacks of $ α_0 $ inputs are allowed if the disk lies on the front side, and a $ β_i' + β_i $ may infinitesimally precede respectively succeed the output mark if the disk lies on the front respectively rear side. A \emph{DS disk} is one of the particular types of degenerate strips fitting between $ L $ and its Hamiltonian deformation. A \emph{DW disk} is one of the particular types of degenerate wedges fitting fitting between $ L $ and its Hamiltonian deformation, with one corner being the co-identity of $ L $. The collections of CR/ID/DS/DW disks are denoted $ \CRd $, $ \IDd $, $ \DSd $ and $ \DWd $ respectively.
\end{definition}

Subdisks of result components of π-trees are defined in the same way as in the consistent case. A few peculiarities of the subdisk construction are depicted in \autoref{fig:sphere-resultcomp-pidiskid} and \autoref{fig:sphere-resultcomp-pidiskidp}. As in the consistent case, associating subdisks provides a bijection between result components and CR, ID, DS and DW disks. Since a π-tree has at least two inputs by definition, the subdisk mapping only reaches CR, ID, DS and DW disks which have at least two inputs as well. We denote these classes of CR, ID, DS and DW disks by $ \CRd^{≥2} $, $ \IDd^{≥2} $, $ \DSd^{≥2} $ and $ \DWd^{≥2} $. We record the bijectivity statement as follows:

\begin{lemma}
\label{th:sphere-resultcomp-bijectivity}
The subdisk mapping $ \Subdisk $ is a bijection
\begin{equation*}
\Subdisk: \PiTr \isoto \CRd^{≥2} \disjoint \IDd^{≥2} \disjoint \DSd^{≥2} \disjoint \DWd^{≥2}.
\end{equation*}
\end{lemma}

\begin{proof}
Injectivity should be clear. Proving surjectivity entails recovering for every CR, ID or DS disk $ D $ a result component $ r $ whose drawing $ \Subdisk(r) $ is $ D $. We will not prove this in detail. In fact, the cases to be studied are merely a subset of the cases of the case of consistent dimers.
\end{proof}

\subsection{Minimal model}
\label{sec:sphere-minmodel}
In this section we provide our minimal model for $ \DefZigzagCat $. The assembly works as follows: In \autoref{sec:sphere-resultcomp}, we have already enumerated all result components for the products $ μ^{≥2}_{\H\DefZigzagCat} $ in terms of CR, ID, DS and DW disks. In \autoref{sec:sphere-defsplitting}, we have computed the differential $ μ^1_{\DefZigzagCat} $ on the deformed cohomology basis elements. In the present section, we assemble the minimal model $ \H\DefZigzagCat $. In particular, we show that not only the higher products $ μ^{≥2}_{\H\DefZigzagCat} $ are computed by CR, ID, DS and DW disks, but also the differential $ μ^1_{\H\DefZigzagCat} $. We offer an explicit list of the CR, ID, DS and DW disks that contribute to the differential.

In \autoref{th:sphere-minmodel-odd}, we claim that the differential $ μ^1_{\H\DefZigzagCat} $ is enumerated accurately by CR, ID, DS and DW disks. It makes sense to compile a list of these disks in advance. Recall that the an SL disk with a single input is a digon: a smooth immersed disks with two corners. In what follows, we try to spot and list all digons of which the input is a given morphism $ h $. We can already ignore DS and DW disks, since they have at least two inputs. The following is our \emph{sphere digon list:}

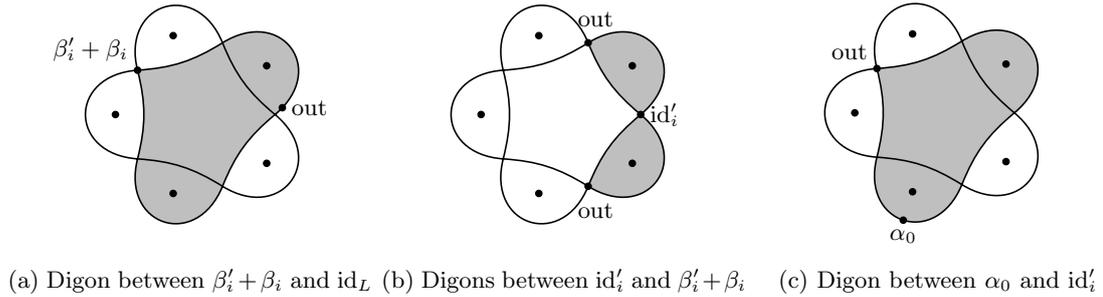
\begin{figure}
\centering
\begin{subfigure}{0.3\linewidth}
\centering
\begin{tikzpicture}
\path[fill, color=gray!50] (0:1) to[out=40, in=32, looseness=2] (72:1) to[out=212, in=-356] (144:1) to[out=284, in=-284] (216:1) to[out=-104, in=248, looseness=2] (288:1) to[out=68, in=-140] cycle;
\path[draw, semithick] (0:1) to[out=140, in=-68] (72:1) to[out=212, in=-356] (144:1) to[out=284, in=-284] (216:1) to[out=356, in=-212] (288:1) to[out=68, in=-140] cycle;
\path[draw, semithick, looseness=2] (0:1) to[out=40, in=32] coordinate[pos=0.05] (id) (72:1) to[out=-248, in=104] (144:1) to[out=-176, in=176] (216:1) to[out=-104, in=248] (288:1) to[out=-32, in=-40] cycle;
\path[fill] \foreach \i in {36, 108, 180, 252, 324} {(\i:1.1) circle[radius=0.05]};
\path[fill] (144:1) circle[radius=0.05] node[above left] {\small $ β_i' + β_i $};
\path[fill] (id) circle[radius=0.05] node[right] {\small out};
\end{tikzpicture}
\caption{Digon between $ β_i' + β_i $ and $ \id_L $}
\label{fig:sphere-minmodel-digons1}
\end{subfigure}
\begin{subfigure}{0.3\linewidth}
\centering
\begin{tikzpicture}
\path[fill, color=gray!50] (0:1) to[out=140, in=-68] (72:1) to[out=32, in=40, looseness=2] cycle;
\path[fill, color=gray!50] (0:1) to[out=-140, in=68] (-72:1) to[out=-32, in=-40, looseness=2] cycle;
\path[draw, semithick] (0:1) to[out=140, in=-68] (72:1) to[out=212, in=-356] (144:1) to[out=284, in=-284] (216:1) to[out=356, in=-212] (288:1) to[out=68, in=-140] cycle;
\path[draw, semithick, looseness=2] (0:1) to[out=40, in=32] (72:1) to[out=-248, in=104] (144:1) to[out=-176, in=176] (216:1) to[out=-104, in=248] coordinate[midway] (id) (288:1) to[out=-32, in=-40] cycle;
\path[fill] \foreach \i in {36, 108, 180, 252, 324} {(\i:1.1) circle[radius=0.05]};
\path[fill] (72:1) circle[radius=0.05] node[above, shift={(0.1, 0.1)}] {\small out};
\path[fill] (-72:1) circle[radius=0.05] node[below, shift={(0.1, -0.1)}] {\small out};
\path[fill] (0:1) circle[radius=0.05] node[right] {\small $ \id_i' $};
\end{tikzpicture}
\caption{Digons between $ \id_i' $ and $ β_i' + β_i $}
\label{fig:sphere-minmodel-digons2}
\end{subfigure}
\begin{subfigure}{0.3\linewidth}
\centering
\begin{tikzpicture}
\path[fill, color=gray!50] (0:1) to[out=40, in=32, looseness=2] (72:1) to[out=212, in=-356] (144:1) to[out=284, in=-284] (216:1) to[out=-104, in=248, looseness=2] (288:1) to[out=68, in=-140] cycle;
\path[draw, semithick] (0:1) to[out=140, in=-68] (72:1) to[out=212, in=-356] (144:1) to[out=284, in=-284] (216:1) to[out=356, in=-212] (288:1) to[out=68, in=-140] cycle;
\path[draw, semithick, looseness=2] (0:1) to[out=40, in=32] (72:1) to[out=-248, in=104] (144:1) to[out=-176, in=176] (216:1) to[out=-104, in=248] coordinate[midway] (id) (288:1) to[out=-32, in=-40] cycle;
\path[fill] \foreach \i in {36, 108, 180, 252, 324} {(\i:1.1) circle[radius=0.05]};
\path[fill] (144:1) circle[radius=0.05] node[above left] {\small out};
\path[fill] (id) circle[radius=0.05] node[below] {\small $ α_0 $};
\end{tikzpicture}
\caption{Digon between $ α_0 $ and $ \id_i' $}
\label{fig:sphere-minmodel-digons3}
\end{subfigure}
\caption{Illustration of digons in $ Q_5 $}
\label{fig:sphere-minmodel-digons}
\end{figure}

\begin{description}
\item[Digons for the odd morphism $ h = β_i' + β_i $:] There is precisely one single digon with input $ h $. It is a CR or ID disk and its output is $ \id_L $. This type of digon is depicted in \autoref{fig:sphere-minmodel-digons1}.

In general, one spots this digon as follows: The intersection point $ h $ cuts the zigzag curve $ \smooth L $ into two segments. One segment departs towards the front side at $ h $ and the other depars towards the rear side. The identity $ \id_L $ lies by choice on one of these two segments. Whether it lies on the front- or rear-bound segment determines the location of the claimed digon. Specifically, the digon lies on the front side if departing at $ h $ towards the front we hit the identity $ \id_L $ before returning to $ h $, and on the rear side if departing at $ h $ towards the rear side we hit the identity $ \id_L $ before returning to $ h $.
\item[Digons for the even morphism $ h = \id_i' $:] There are precisely two digons with input $ \id_i' $. Both are CR disks with output type $ β_i' + β_i $. They simply reach around the punctures neighboring the input $ \id_i' $. They are depicted in \autoref{fig:sphere-minmodel-digons2}.
\item[Digons for the co-identity $ h = α_0 $:] We spot $ M $ digons contributing to $ μ^1 ((-1)^{\#α_0 + 1} α_0) $, namely $ (M-1)/2 $ on the front side and $ (M+1)/2 $ on the rear side. These digons are all CR disks and have output of type $ \id_i' $. In case of $ M = 5 $, these digons are all heart-shaped and depicted in \autoref{fig:sphere-minmodel-digons3}.

In general, one spots these digons as follows: Of the $ M $-many self-intersection point $ p ∈ L ∩ L $, fix an arbitrary one. We shall construct from this data one certain digon that has corners $ h $ and $ p $. For this, note that $ p $ cuts the zigzag curve $ \smooth L $ into two segments. Only one of these two segments contains the co-identity location $ α_0 $. The digon associated with $ p $ is then the digon bounded by this segment. In other words, if the segment containing $ α_0 $ departs to the front side at $ h $, then the digon lies on the front side. If the segment containing $ α_0 $ departs to the rear side at $ h $, then the digon lies on the rear side. This determines a digon contributing to $ μ^1 (h) $ for every of the self-intersection points $ p ∈ L ∩ L $.
\end{description}

Before we reach the main theorem, we shall comment on the signs of result components. Recall that the Abouzaid sign of an SL disk is defined in \autoref{def:subdisk-minmodel-Abouzaid} and allows an arbitrary nonnegative number of inputs. The definition of the Abouzaid sign carries over without change to the case of $ Q = Q_M $. In analogy to \autoref{th:subdisk-minmodel-signs}, the sign of a result component agrees with the Abouzaid sign of its subdisk:

\begin{lemma}
\label{th:sphere-minmodel-signs}
Let $ r $ be the result component of a $ π $-tree. Then the sign of $ r $, relative to the signs of the output value, equals the Abouzaid sign of its subdisk $ \Subdisk(r) $. The $ q $-parameter $ ∈ ℂ⟦Q_0⟧ $ is equal to $ \punctures(D) $, the product of all punctures covered by $ \Subdisk(r) $ counted with multiplicities.
\end{lemma}

In contrast to the consistent case, the category $ \H\DefZigzagCat $ has a residual differential. We can in fact describe the differential by means of CR and ID digons, the sign being equal to the Abouzaid sign. In contrast to the consistent case, the definition of the Abouzaid sign rule is here also used for digons. We are now ready to formulate our freshly built interpretation as a description of the minimal model.

\begin{proposition}
\label{th:sphere-minmodel-odd}
Let $ Q_M $ be the standard sphere dimer with an odd number $ M ≥ 3 $ of punctures. Let $ h_1, …, h_N $ be a sequence of $ N ≥ 0 $ non-identity basis morphisms with $ h_i: L_i → L_{i+1} $. Then their product is given by
\begin{equation*}
μ_{\H\DefZigzagCat}^N (h_N, …, h_1) = \sum_{\substack{D ∈ \CRd \disjoint \IDd \disjoint \DSd \disjoint \DWd \\ D \text{ has inputs } h_1, …, h_N}} (-1)^{\Abouzaid(D)} \punctures(D) \disktarget(D).
\end{equation*}
More explicitly,
\begin{itemize}
\item The curvature $ μ^0_{\H\DefZigzagCat} $ vanishes.
\item The differential $ μ^1_{\H\DefZigzagCat} $ is given by the digons in the sphere digon list, with Abouzaid sign rule.
\item We have $ (-1)^{|h|} μ^2_{\H\DefZigzagCat} (\id_L, h) = μ^2_{\DefZigzagCat} (h, \id_L) = h $ and $ μ^3_{\H\DefZigzagCat} (…, \id_L, …) = 0 $.
\item Products $ μ^{≥2}_{\H\DefZigzagCat} $ on all sequences of non-identity inputs are given by the CR, ID, DS and DW disks, with Abouzaid sign rule.
\end{itemize}
\end{proposition}

\begin{proof}
The minimal model $ \H\DefZigzagCat $ is given by the deformed Kadeishvili theorem from \papertwoA. To start with, recall that we have already computed the deformed cohomology basis elements in \autoref{th:sphere-defsplitting-basis}. The deformed Kadeishvili theorem then provides curvature, differential and products for $ \H\DefZigzagCat $. We shall put focus on checks for the differential, since the structure of the products is very similar to the structure observed in the geometrically consistent case.

We have noted in \autoref{th:sphere-defzigzagcat-curvature} that the curvature $ μ^0_{\DefZigzagCat} $ already vanishes. According to our Kadeishvili theorem the curvature $ μ^0_{\H\DefZigzagCat} $ then vanishes as well. This already proves the first statement. We check the three remaining statements in order.

For the second statement, regard the differential $ μ^1_{\DefZigzagCat} $. According to the deformed Kadeishvili theorem, the differential $ μ^1_{\H\DefZigzagCat} $ is given by the composition of $ μ^1_{\DefZigzagCat} $ and the projection to deformed cohomology $ π_{H_q} $.

We have computed the differential $ μ^1_{\DefZigzagCat} $ already in \autoref{th:sphere-defsplitting-basis} and observed that $ μ^1_{\DefZigzagCat} (H_q) ⊂ H_q $. The projection is therefore without effect and we have $ μ^1_{\H\DefZigzagCat} (h) = μ^1_{\DefZigzagCat} (h) $ for deformed cohomology basis morphisms $ h $. For example, we have
\begin{equation*}
μ^1_{\H\DefZigzagCat} ((-1)^{\#(i+M) + 1} β_i' + (-1)^{\#i} β_i) = (-1)^{\#\#(i+1) + 1} Q_{i+2} \id_L.
\end{equation*}
It remains to show that for every deformed cohomology basis element $ h $, its differential $ μ^1_{\H\DefZigzagCat} (h) $ is enumerated accurately by CR, ID, DS and DW disks. We have listed all CR, ID, DS and DW disks in the sphere digon list. For any deformed cohomology basis element $ h $ of one of the three types $ β_i/β_i' $, $ \id_i' $ and $ α_0 $, it remains to interpret the result terms $ r $ of $ μ^1_{\H\DefZigzagCat} (h) $ as enumeration over the digons presented in the sphere digons list. The crucial part is to prove the sign of every term $ r $ equal to the Abouzaid sign $ \Abouzaid(D) $ of the corresponding digon $ D $.

First, regard the odd morphism $ h = β_i' + β_i $, with odd $ i $. The calculation of \autoref{th:sphere-defsplitting-basis} gives one single output term, namely the identity $ \id_L $. This is exactly the result enumerated by the single digon $ D $ from the sphere digons list. It remains to compare the sign $ (-1)^{\#\#(i+1) + 1} $ with the Abouzaid sign of $ D $. The sign of $ D $ is computed as follows: In case the digon lies on the front side, the $ \# $ signs to be summed up are $ \#i, …, \#(i-M+1) $. Their sum amounts to $ \#\#(i-M+1) $, which has equal parity with $ \#\#(i+1) + 1 $, since the total number of $ \# $ signs in $ Q_M $ is assumed to be odd. In case the digon lies on the rear side, the $ \# $ signs to be summed up are $ \#(i+M), …, \#(i+1) $. Their sum amounts to $ \#\#(i+1) $. An additional sign flip is due, since $ h $ is odd and lies counterclockwise with respect to the rear side. Ultimately, both front and rear cases give the sign $ \#\#(i+1) + 1 $. This sign agrees with our calculation of $ μ^1_{\DefZigzagCat} (h) $ in \autoref{th:sphere-defsplitting-basis}.

Second, regard the even morphism $ \id_i' $. There are two digons in our digon list contributing to the differential $ μ^1 (\id_i') $, namely the two small digons reaching around the neighboring punctures. The Abouzaid sign rule predicts a sign of $ (-1)^{\#i + \#(i+M)} $ for the upper puncture and $ (-1)^{\#(i-1) + \#(i+M-1)} $ for the lower puncture, exactly as calculated in \autoref{th:sphere-defsplitting-basis}.

Third, regard the co-identity $ (-1)^{\#α_0 + 1} α_0 $. There are $ M $ digons in our digon list contributing to the differential $ μ^1 ((-1)^{\#α_0 + 1} α_0) $, namely $ (M-1)/2 $ on the front side and $ (M+1)/2 $ on the rear side. For the front side disk with output $ \id_{i-j+1}' $, the Abouzaid sign rule predicts a sign of $ (-1)^{\#\#(i-j+1)} $. For the rear side disk with output $ \id_{i+j+1}' $, the Abouzaid sign rule predicts a sign of $ (-1)^{\#\#(i+j-M+1) + 1} $, the absolute sign flip $ (-1) $ coming from the odd co-identity whose zigzag path runs counterclockwise with respect to the rear side.

We conclude that the differential $ μ^1_{\H\DefZigzagCat} (h) $ is computed accurately by the digons from the sphere digons list for any of the three types of morphisms $ h $. This finishes the checks for the second statement of the proposition.

The third statement of the proposition is trivial, following immediately from unitality of $ \DefZigzagCat $ and the choice that $ \id_L ∈ H $. The fourth statement follows from \autoref{th:sphere-resultcomp-bijectivity} and \autoref{th:sphere-minmodel-signs}, in a way entirely analogous to the geometrically consistent case. This finishes the proof.
\end{proof}

\subsection{The case of even $ M $}
\label{sec:sphere-even}
In this section, we comment on the category of zigzag paths of $ Q_M $ for even $ M $. We define the category $ \ZigzagCat $ of zigzag paths and explain how to obtain a homological splitting. We explain how to run the curvature optimization for the corresponding subcategory of $ \Tw\Gtl_q Q $ and define the category of deformed zigzag paths $ \ZigzagCat $. Finally, we provide a minimal model for $ \H\DefZigzagCat $.

The dimer $ Q_M $ for even $ M ≥ 4 $ has two zigzag paths, each consisting of $ M $ arcs. There are $ M $ intersections between the zigzag curves. The front side of the dimer is a clockwise polygon, the rear side is a counterclockwise polygon.

\begin{convention}
The letter $ Q = Q_M $ with $ M ≥ 4 $ even denotes the standard sphere dimer with $ M $ punctures. The spin structure is chosen by assigning an even number of $ \# $ signs on the rear side of $ Q_M $ and none on the front side. The identity locations are arbitrary chosen, and the co-identity location is chosen to lie on the rear side of $ Q_M $.
\end{convention}

\begin{definition}
The category $ \ZigzagCat ⊂ \Tw\Gtl Q_M $ is the category consisting of the two zigzag paths in $ Q_M $. The standard splitting $ H ⊕ I ⊕ R $ for $ \ZigzagCat $ is defined in the analogous way as for odd $ M $. The category $ \DefZigzagCat ⊂ \Tw'\Gtl_q Q_M $ is the category consisting of the two zigzag paths with deformed twisted differential analogous to \autoref{def:sphere-defzigzagcat-def}.
\end{definition}

A priori, it is our task to compute a minimal model of the category of zigzag paths in $ \Tw\Gtl_q Q $ consisting of the same twisted complexes as $ \ZigzagCat ⊂ \Tw\Gtl Q $. As usual, we are allowed to apply gauge in order to optimize the curvature. In contrast to the case of geometrically consistent $ Q $ or $ Q_M $ for odd $ M $, the category $ \DefZigzagCat $ is not curvature-free, but its curvature is optimal nevertheless:

\begin{lemma}
\label{th:sphere-even-optimal}
The curvatures of both zigzag paths $ L_1, L_2 ∈ \DefZigzagCat $ are multiple of their respective identities $ \id_{L_1} ∈ H $ and $ \id_{L_2} ∈ H $. In particular, $ \DefZigzagCat $ has optimal curvature.
\end{lemma}

The deformed decomposition $ H_q ⊕ μ^1_q (B \htensor R) ⊕ (B \htensor R) $ of $ \DefZigzagCat $ is similar to the case of odd $ M $. The differential does not vanish and maps to $ H_q ⊕ (B \htensor R) $. According to the deformed Kadeishvili construction, we can compute $ \H\DefZigzagCat $ by setting
\begin{align*}
μ^0_{\H\DefZigzagCat} &= μ^0_{\DefZigzagCat}, \\
μ^1_{\H\DefZigzagCat} (h) &= π_{H_q} μ^1_{\DefZigzagCat} (h), \\
μ^{N≥2}_{\H\DefZigzagCat} (h_N, …, h_1) &= \sum_{T ∈ \mathcal{T}_N} (-1)^{N_T} \Res(T, h_1, …, h_N).
\end{align*}
The computation for $ μ^{N≥2}_{\H\DefZigzagCat} $ is similar to the case of odd $ M $. The computation for $ μ^1_{\H\DefZigzagCat} $ is similar to the case of odd $ M $ as well, with the difference that $ μ^1_{\H\DefZigzagCat} $ does not cancel because of the different choices of $ \# $ signs. The computation for $ μ^0_{\H\DefZigzagCat} $ is elementary. As in the case of odd $ M $, it turns out that the entire $ A_∞ $-structure of the minimal model can be described through CR, ID, DS and DW disks:

\begin{proposition}
\label{th:sphere-minmodel-even}
Let $ Q_M $ be the standard sphere dimer with an even number $ M ≥ 4 $ of punctures. Let $ h_1, …, h_N $ be a sequence of $ N ≥ 0 $ non-identity basis morphisms with $ h_i: L_i → L_{i+1} $. Then their product is given by
\begin{equation*}
μ_{\H\DefZigzagCat}^N (h_N, …, h_1) = \sum_{\substack{D ∈ \CRd \disjoint \IDd \disjoint \DSd \disjoint \DWd \\ D \text{ has inputs } h_1, …, h_N}} (-1)^{\Abouzaid(D)} \punctures(D) \disktarget(D).
\end{equation*}
\end{proposition}

\section{Calculating the mirror objects}
\label{sec:mirobjects}
The aim of this section is to perform further minimal model calculations which we need for the third paper in the series. Our calculations are based on the general deformed Kadeishvili construction of \papertwoA, because the simplified Kadeishvili construction of \autoref{sec:2Bkadeishvili} does not apply.

In \autoref{sec:mirobjects-ms}, we explain which products in the minimal model need to be computed and why. In \autoref{sec:mirobjects-goal}, we describe the input data of the minimal model construction. In \autoref{sec:mirobjects-splitting}, we construct a homological splitting. In \autoref{sec:mirobjects-defsplitting}, we compute the deformed decomposition. In \autoref{sec:mirobjects-resultcomp}, we introduce a suitable notion of result components and classify them into two types which we call MD and MT result components. In \autoref{sec:mirobjects-products}, we show how to match MD/MT result components with disks of two types which we call MD/MT disks. In \autoref{th:mirobjects-products-th} we finally describe the desired products in terms of MD and MT disks.

\subsection{Mirror symmetry for punctured surfaces}
\label{sec:mirobjects-ms}
In this section, we explain the reason we need to perform further minimal model calculations. The starting point is a brief recapitulatation of mirror symmetry for punctured surfaces. We then explain the idea of the deformed Cho-Hong-Lau construction and describe which products we need to compute.

Mirror symmetry for punctured surfaces \cite{Bocklandt} entails a quasi-isomorphism
\begin{equation*}
F: \Gtl Q → \mf(\Jac \mirQ, ℓ).
\end{equation*}
Here $ Q $ is a dimer and $ \mirQ $ is its mirror dimer which is assumed to be zigzag consistent. The vertices of $ \mirQ $ are identified with the zigzag paths in $ Q $. The algebra $ \Jac \mirQ $ is the so-called Jacobi algebra of the dimer and is explicitly defined as the quiver algebra $ ℂ\mirQ $ modulo relations. The element $ ℓ ∈ \Jac \mirQ $ is a central element known as the potential. The category $ \MF(A, ℓ) $ denotes the dg category of so-called matrix factorizations of $ (A, ℓ) $. The category $ \mf(\Jac \mirQ, ℓ) $ denotes one a certain small subcategory of $ \MF(\Jac \mirQ, ℓ) $, specific to mirror symmetry.

The deformed mirror symmetry which we prove in the third paper entails a quasi-isomorphism of deformed $ A_∞ $-categories
\begin{equation*}
F_q: \Gtl_q Q → \mf(\Jac_q \mirQ, ℓ_q).
\end{equation*}
The category $ \Gtl_q Q $ on the left-hand side has an object for every arc $ a ∈ Q_1 $. The category on the right-hand side is what we will call a deformed category of matrix factorizations. It has one object for every arc $ a ∈ Q_1 $ as well. The functor $ F_q $ matches the arc $ a ∈ \Gtl_q Q $ with an deformed matrix factorization $ F_q (a) $.

In the third paper, we compute the deformed algebra $ \Jac_q \mirQ $ and deformed potential $ ℓ_q $. The starting point is the category of deformed zigzag paths $ \DefZigzagCat ⊂ \Tw'\Gtl_q Q $. Thanks to the description of the minimal model $ \H\DefZigzagCat $ which we provided in the present paper in \autoref{sec:subdisk-minmodel}, we express in the third paper the deformed algebra $ \Jac_q \mirQ $ and the deformed potential $ ℓ_q $ explicitly in terms of combinatorical data of $ Q $. Viewed the other way around, the present paper is the technical cornerstone for the third paper.

In the third paper, we also compute the precise shape of the deformed matrix factorizations contained in $ \mf(\Jac_q \mirQ, ℓ_q) $. According to the deformed Cho-Hong-Lau construction, the precise shape is given by certain products in $ \HTw\Gtl_q Q $ which not only involve zigzag paths, but also the arc objects. The description of these products is not included in \autoref{sec:subdisk}. Therefore, the we have devoted the present \autoref{sec:mirobjects} to supplementing these products.

The objects of $ \mf(\Jac_q, ℓ_q) $ are explicitly of the form
\begin{align*}
F_q (a) &= \left(\bigoplus_{L ∈ \ZigzagCat} \Hom_{\HTw\Gtl Q} (L, a) ¤ (\Jac_q \mirQ) v_L, δ\right), \\
δ (m) &= \sum_{N ≥ 0} μ_{\HTw\Gtl_q Q} (m, \underbrace{b, …, b}_N).
\end{align*}
Here $ a $ denotes any arc in $ Q $, the letter $ v_L $ denotes the vertex of $ ℂ\mirQ $ defined by the zigzag path $ L $, and $ m $ denotes an element of $ \Hom_{\HTw\Gtl Q} (L, a) $. The element $ b $ denotes essentially a formal sum over all type B cohomology basis elements between zigzag paths in $ Q $. Geometrically, the element $ b $ includes all odd transversal intersections between zigzag curves.

The essential calculation which we shall therefore perform in the present \autoref{sec:mirobjects} consists of determining the hom space $ \Hom_{\HTw\Gtl Q} (L, a) $ and computing all possible kinds of products of the form $ μ_{\HTw\Gtl_q Q} (m, h_N, …, h_1) $. Here $ h_i: L_i → L_{i+1} $ are type B cohomology basis elements between zigzag paths and $ m: L_{N+1} → a $ is a cohomology basis element from $ L_{N+1} $ to an arc $ a ∈ Q_1 $.

\subsection{The desired products}
\label{sec:mirobjects-goal}
In this section, we examine which minimal model $ \H\Tw\Gtl_q Q $ we shall compute. In principle, we are free to choose any minimal model. When computing the mirror Jacobi algebra $ \Jac_q Q $ and potential $ ℓ_q $, we have however already made a choice for minimal model of $ \DefZigzagCat ⊂ \Tw\Gtl_q Q $. The minimal model we compute here needs to be compatible with these earlier choices.

\begin{convention}
\label{conv:mirobjects-dimer}
$ Q $ is a geometrically consistent dimer or a standard sphere dimer $ Q = Q_M $ with $ M ≥ 3 $. The dimer is equipped with choices of spin structure, identity location $ a_0 $ and co-identity location $ α_0 $ for every zigzag path. In case $ Q = Q_M $ with $ M $ odd, the spin structure is chosen by assigning $ \#α = 1 $ to an odd number of interior angles $ α $ on the rear side and $ \# α = 0 $ to all other angles. In case $ Q = Q_M $ with $ M $ even, the spin structure is chosen by assigning $ \#α = 1 $ to an even number of interior angles $ α $ on the rear side and $ \# α = 0 $ to all other angles. The co-identity $ α_0 $ shall be chosen to lie in a counterclockwise polygon.
\end{convention}

\begin{definition}
We denote by $ \ZigzagCat ⊂ \Tw\Gtl Q $ and $ \DefZigzagCat ⊂ \Tw'\Gtl_q Q $ the categories of zigzag paths and of deformed zigzag paths, defined as follows:
\begin{itemize}
\item If $ Q $ is geometrically consistent, $ \ZigzagCat $ is the category of zigzag paths as defined in \autoref{sec:prelim-zigzagcat} and $ \DefZigzagCat $ is the category of deformed zigzag paths as defined in \autoref{sec:deformed-uncurving}.
\item If $ Q = Q_M $ for odd $ M ≥ 3 $, then $ \ZigzagCat $ is the category of zigzag paths as defined in \autoref{sec:sphere-zigzagcat} and $ \DefZigzagCat $ is the category of deformed zigzag paths as defined in \autoref{sec:sphere-defzigzagcat}.
\item If $ Q = Q_M $ for even $ M ≥ 3 $, then $ \ZigzagCat $ is the category of zigzag paths as defined in \autoref{sec:sphere-even} and $ \DefZigzagCat $ is the category of deformed zigzag paths as defined in \autoref{sec:sphere-even}.
\end{itemize}
\end{definition}

To construct a deformed mirror functor $ F_q: \H\Tw\Gtl_q Q → \mf(\Jac_q \mirQ, ℓ_q) $, we need a choice of minimal model of the entire category $ \Tw\Gtl_q Q $. According to our deformed Kadeishvili theorem, we can obtain such a minimal model by optimizing curvature for all objects and performing a familiar Kadeishvili tree construction. The input data for this construction consists only of a homological splitting for every hom space in $ \Tw\Gtl Q $. All other steps are automatic.

\begin{remark}
\label{rem:mirobjects-goal-kadeishvili}
Our deformed Kadeishvili construction has the property that the restriction of a minimal model $ \H\cat C_q $ to a subcategory $ \cat D ⊂ \cat C $ is the same as the minimal model $ \H\cat D $, when the homological splitting chosen for $ \cat D $ is the restriction of the homological splitting chosen for $ \cat C $:
\begin{equation*}
\H\cat C \restr{\Ob\cat D} = \H\cat D.
\end{equation*}
\end{remark}

To construct the restriction of this functor to $ \Gtl_q Q $, we however do not need to compute the entire minimal model $ \HTw\Gtl_q Q $ explicitly. It suffices to know products of the kind $ μ(h_N, …, h_1) $ and $ μ(m, h_N, …, h_1) $, where $ h_1, …, h_N $ are morphisms between zigzag paths and $ m $ is a morphism from a zigzag path to an arc. By \autoref{rem:mirobjects-goal-kadeishvili}, it suffices to compute a minimal model of the category $ Q_1 ∪ \DefZigzagCat $, which is defined as the subcategory of $ \Tw'\Gtl_q Q $ consisting of arcs and deformed zigzag paths. We define this category precisely as follows:

\begin{definition}
The subcategory of $ \Tw\Gtl Q $ given by the union of $ \Gtl Q $ and $ \ZigzagCat $ is denoted
\begin{equation*}
Q_1 ∪ \ZigzagCat ⊂ \Tw\Gtl Q.
\end{equation*}
The subcategory of $ \Tw\Gtl_q Q $ given by the union of $ \Gtl_q Q $ and $ \DefZigzagCat $ is denoted
\begin{equation*}
Q_1 ∪ \DefZigzagCat ⊂ \Tw\Gtl_q Q.
\end{equation*}
\end{definition}

\begin{remark}
\label{rem:mirobjects-goal-extension}
Applying the Kadeishvili construction to $ Q_1 ∪ \DefZigzagCat $ involves choosing a homological splitting for $ Q_1 ∪ \ZigzagCat $. The deformed mirror symmetry construction in the third paper departs from a single minimal model model $ \HTw\Gtl_q Q $. In consequence, it is not allowed to compute the products of the two kinds $ μ_{\Tw\Gtl_q Q} (h_N, …, h_1) $ and $ μ_{\Tw\Gtl_q Q} (m, h_N, …, h_1) $ via different homological splittings of $ \ZigzagCat $. Instead, the homological splitting for $ Q_1 ∪ \ZigzagCat $ needs to extend the homological splitting already chosen for $ \ZigzagCat $.
\end{remark}

\subsection{Homological splitting}
\label{sec:mirobjects-splitting}
In this section, we construct a homological splitting for the category $ Q_1 ∪ \ZigzagCat $. The starting point is the definition of the category $ Q_1 ∪ \ZigzagCat $ in \autoref{sec:mirobjects-goal}. According to \autoref{rem:mirobjects-goal-extension}, we have to define the homlogical splitting for $ Q_1 ∪ \ZigzagCat $ in the following way:
\begin{itemize}
\item For the hom spaces $ \Hom_{\Tw\Gtl Q} (L_1, L_2) $ between two zigzag paths $ L_1, L_2 $, the homological splitting is the homological splitting already established for $ \ZigzagCat $. In case $ Q $ is geometrically consistent, this refers to the homological splitting of \autoref{sec:splitting-splitting}. In case $ Q = Q_M $ for odd $ M ≥ 3 $, this refers to the homological splitting of \autoref{sec:sphere-splitting}. In case $ Q = Q_M $ for even $ M ≥ 4 $, this refers to the analog of \autoref{sec:sphere-splitting} indicated in \autoref{sec:sphere-minmodel}.
\item For the hom spaces $ \Hom_{\Tw\Gtl Q} (L, a) $ between a zigzag path $ L $ and an arc $ a ∈ Q_1 $, we are free to choose a homological splitting.
\item For the hom spaces $ \Hom_{\Tw\Gtl Q} (a, L) $ between an arc $ a ∈ Q_1 $ and a zigzag path $ L $, we are free to choose a homological splitting. In practice, this choice is irrelevant for the calculation of the products $ μ_{\Tw\Gtl Q} (m, h_N, …, h_1) $, so we will merely assume any arbitrary splitting has been chosen.
\item For the hom spaces $ \Hom_{\Tw\Gtl Q} (a, b) $ between two arcs $ a, b ∈ Q_1 $, the homological splitting is predetermined as $ H = \Hom_{\Tw\Gtl Q} (a, b) $ and $ I = R = 0 $ by the fact that $ μ^1_{\Gtl Q} = 0 $.
\end{itemize}

According to this list, the only remaining task is to choose a homological splitting for $ \Hom_{\Tw\Gtl Q} (L, a) $ whenever $ L $ is a zigzag path and $ a $ an arc. We start by classifying morphisms $ L → a $ into three types of situations.

As in the case of $ \ZigzagCat $, let an \emph{elementary morphism} $ ε: L → a $ refer to a morphism between twisted complexes consisting of a single angle of $ Q $. We shall associate with every elementary morphism a \emph{situation}. As usual, the terminology is as follows: Every elementary morphism belongs to a unique situation, and every situation is of one given type. Running out of letters, we will denote the situation types by S1, S2, S3. Every situation is defined by a certain collection of arcs and angles. The other way around, every situation has a collection of elementary morphisms associated, constructed from its defining angles.

\begin{figure}
\centering
\begin{subfigure}{0.3\linewidth}
\centering
\begin{tikzpicture}
\path[draw, ->] (0, 0) -- ++(150:1) coordinate[midway] (beta-start) -- ++(up:1) coordinate[pos=0.4] (betap-start) coordinate[pos=0.6] (alpha-end) -- ++(150:1) coordinate[pos=0.4] (alpha-start);
\path[draw] (0.2, 0)++(150:1) -- ++(up:1) coordinate[pos=0.4] (beta-end) coordinate[pos=0.6] (alphap-start);
\path[draw, bend right=80, ->] (beta-start) to node[pos=0.9, below] {\small $ β $} (beta-end);
\path[draw, bend right=80, ->] (alpha-start) to node[pos=0.4, right] {\small $ α $} (alpha-end);
\path[draw, bend right=100, looseness=2.5, ->] (betap-start) to node[pos=0.5, right] {\small $ β' $} (beta-start);
\path[draw, bend right=100, looseness=2.5, ->] (alphap-start) to node[pos=0.7, below] {\small $ α' $} (alpha-start);
\end{tikzpicture}
\caption{S1 situation}
\label{fig:mirobject-splitting-S1}
\end{subfigure}
\begin{subfigure}{0.3\linewidth}
\centering
\begin{tikzpicture}
\path[draw, ->] (0, 0) -- ++(30:1) coordinate[midway] (beta-end) -- ++(up:1) coordinate[pos=0.4] (betap-end) coordinate[pos=0.6] (alpha-start) -- ++(30:1) coordinate[pos=0.4] (alpha-end);
\path[draw] (-0.2, 0)++(30:1) -- ++(up:1) coordinate[pos=0.4] (beta-start) coordinate[pos=0.6] (alphap-end);
\path[draw, bend right=80, ->] (beta-start) to node[pos=0.2, below] {\small $ β $} (beta-end);
\path[draw, bend right=80, ->] (alpha-start) to node[pos=0.2, above] {\small $ α $} (alpha-end);
\path[draw, bend right=100, looseness=2.5, ->] (beta-end) to node[pos=0.4, above] {\small $ β' $} (betap-end);
\path[draw, bend right=100, looseness=2.5, ->] (alpha-end) to node[pos=0.3, below] {\small $ α' $} (alphap-end);
\end{tikzpicture}
\caption{S2 situation}
\label{fig:mirobject-splitting-S2}
\end{subfigure}
\begin{subfigure}{0.3\linewidth}
\centering
\begin{tikzpicture}
\path[draw] (0, 0) -- ++(170:1) coordinate[pos=0.5] (alpha-start) coordinate (stop) -- ++(60:1) coordinate[pos=0.5] (alpha-end);
\path[draw] (stop)++(210:0.2) -- ++(210:1) coordinate[pos=0.4] (beta-end);
\path[draw, ->, bend right=30] (alpha-start) to node[pos=0.9, below] {\small $ α $} (alpha-end);
\path[draw, ->, bend right=80, looseness=2] (alpha-end) to node[pos=0.4, below] {\small $ β $} (beta-end);
\path[draw, ->, bend right=80, looseness=2] (beta-end) to node[pos=0.5, above] {\small $ γ $} (alpha-start);
\end{tikzpicture}
\caption{S3 situation}
\label{fig:mirobject-splitting-S3}
\end{subfigure}
\end{figure}

\begin{definition}
An \emph{S1} situation consists of a zigzag path $ L $ and an indexed arc $ a $ of $ L $ such that $ L $ turns left at the head of $ a $. The nearby angles of the situation are denoted $ α $, $ α' $, $ β $, $ β' $ as in \autoref{fig:mirobject-splitting-S1}. The elementary morphisms belonging to an S1 situation are the morphisms $ ε: L → a $ given by $ \id_{L → a} $, $ β ℓ^k $, $ β β' ℓ^k $, $ α ℓ^k $, $ α α' ℓ^k $.

An \emph{S2} situation consists of a zigzag path $ L $ and an indexed arc $ a $ such that $ L $ turns right at the head of $ a $. The nearby angles of the situation are denoted $ α $, $ α' $, $ β $, $ β' $ as in \autoref{fig:mirobject-splitting-S2}. The elementary morphisms belonging to an S2 situation are the morphisms $ ε: L → a $ given by $ \id_{L → a} $, $ α' ℓ^k $, $ α' α ℓ^k $, $ β' ℓ^k $, $ β' β ℓ^k $.

An \emph{S3} situation consists of two consecutive indexed arcs on a zigzag path $ L $ and an arc $ a $ such that $ a $ is incident at the common puncture of the two arcs but not equal to either of them. The nearby angles of the situation are denoted $ α $, $ β $, $ γ $ as in \autoref{fig:mirobject-splitting-S3}. The elementary morphisms belonging to an S3 situation are the morphisms $ ε: L → a $ given by $ β ℓ^k $, $ βα ℓ^k $.
\end{definition}

We have constructed these definitions so that the situations exhaust all elementary morphisms in $ \Hom_{\Tw\Gtl Q} (L, a) $. We record this as follows:

\begin{lemma}
Let $ L $ be a zigzag path and $ a $ an arc. Then any elementary angle $ ε: L → a $ belongs to precisely one S1, S2 or S3 situation.
\end{lemma}

We are now ready to construct our homological splitting for $ \Hom_{\Tw\Gtl Q} (L, a) $. This means to provide a choice of basis elements for $ H $ and $ R $.

\begin{definition}
\label{def:mirobjects-splitting-def}
Let $ L $ be a zigzag path and $ a $ an arc. We let $ R ⊂ \Hom_{\Tw\Gtl Q} (L, a) $ be the subspace spanned by:
\begin{itemize}
\item for every S1 situation the morphisms $ \id_{L → a} $, $ α α' ℓ^k $ and $ β β' ℓ^k $,
\item for every S2 situation the morphisms $ α' ℓ^k $ and $ β' ℓ^k $,
\item for every S3 situation the morphism $ β ℓ^k $.
\end{itemize}
The space $ H ⊂ \Hom_{\Tw\Gtl Q} (L, a) $ is spanned by:
\begin{itemize}
\item for every S1 situation the morphism $ (-1)^{\#β} β $,
\item for every S2 situation the morphism $ \id_{L → a} $.
\end{itemize}
Setting $ I = μ^1_{\Tw\Gtl Q} (R) $, we refer to $ H, I, R $ as the \emph{(standard) splitting} for $ \Hom_{\Tw\Gtl Q} (L, a) $.
\end{definition}

It is an elementary check that the standard splitting is indeed a homological splitting.

\begin{lemma}
Let $ a $ be an arc and $ L $ a zigzag path. Then the standard splitting indeed forms a homological splitting for $ \Hom_{\Tw\Gtl Q} (L, a) $.
\end{lemma}

\subsection{Deformed decomposition}
\label{sec:mirobjects-defsplitting}
In this section, we determine the relevant part of the deformed decomposition of $ Q_1 ∪ \DefZigzagCat $. The starting point is the category $ Q_1 ∪ \DefZigzagCat $ defined in \autoref{sec:mirobjects-goal} and the homological splitting for $ Q_1 ∪ \ZigzagCat $ defined in \autoref{sec:mirobjects-splitting}. In the present section, we show that $ Q_1 ∪ \DefZigzagCat $ has optimal curvature. We determine explicitly the deformed decomposition of the hom spaces $ \Hom_{Q_1 ∪ \DefZigzagCat} (L, a) $, where $ L ∈ \ZigzagCat $ is a zigzag path and $ a ∈ Q_1 $ an arc.

\begin{lemma}
The category $ Q_1 ∪ \DefZigzagCat $ has optimal curvature.
\end{lemma}

\begin{proof}
If $ Q $ is geometrically consistent, then $ \DefZigzagCat $ is curvature-free by \autoref{th:deformed-uncurving-curvaturefree}. If $ Q = Q_M $ with $ M $ odd, then $ \DefZigzagCat $ is curvature-free by \autoref{th:sphere-defzigzagcat-curvature}. If $ Q = Q_M $ with $ M $ even, then $ \DefZigzagCat $ has optimal curvature by \autoref{th:sphere-even-optimal}. The subcategory $ Q_1 ⊂ Q_1 ∪ \DefZigzagCat $ has optimal curvature by nature, so we conclude that $ Q_1 ∪ \DefZigzagCat $ has optimal curvature.
\end{proof}

Since $ Q_1 ∪ \DefZigzagCat $ already has optimal curvature, the products $ μ_{\HTw\Gtl_q Q} (m, h_N, …, h_1) $ can be obtained by computing the deformed decomposition of $ \Hom_{Q_1 ∪ \DefZigzagCat} (L, a) $ and evaluating Kadeishvili π-trees. As next step, we shall therefore focus on finding the deformed decomposition for $ \Hom_{Q_1 ∪ \DefZigzagCat} $. More precisely, we are interested in the deformed cohomology basis elements.

Finding the deformed decomposition entails finding for every cohomology basis element $ h ∈ H $ a deformed cohomology basis element $ φ^{-1} (h) = h + r $ such that $ r ∈ B \htensor R $ and
\begin{equation*}
μ^1_{\Tw\Gtl_q} (h + r) ∈ (B \htensor H) ⊕ (B \htensor R).
\end{equation*}
For the cohomology basis elements $ h $ of the hom space $ \Hom_{\Tw\Gtl Q} (L, a) $, we can compute $ φ^{-1} (h) $ explicitly:

\begin{lemma}
\label{th:mirobjects-defsplitting-th}
Let $ a ∈ Q_1 $ be an arc and $ L $ a zigzag path. Then the space $ H $ satisfies
\begin{equation*}
μ^1_{\Tw\Gtl_q Q} (H) ⊂ (B \htensor H) ⊕ (B \htensor R).
\end{equation*}
In particular, we have $ H_q = B \htensor H $ and the map $ φ: H_q → B \htensor H $ is the identity.
\end{lemma}

\begin{proof}
Let us start by checking for every cohomology basis element $ h $ that $ μ^1_{\Tw\Gtl_q Q} (h) $ lies in $ (B \htensor H) ⊕ (B \htensor R) $. Let $ L $ be a zigzag path and $ a $ an arc. Denote by $ q = h(a) $ the puncture at the head of $ a $ and by $ p = t(a) $ the puncture at the tail of $ a $. Denote by $ δ ∈ \Hom_{\Add\Gtl_q Q}^1 (L, L) $ the twisted differential of $ L ∈ \DefZigzagCat $.

Regard an S1 situation between $ L $ and $ a $. Denote by $ β, β' $ the angles associated with the S1 situation. We want to compute the differential of the cohomology basis element $ h = (-1)^{\#β} β $. We have
\begin{align*}
μ^1_{\Tw\Gtl_q Q} (β) &= \sum_{k ≥ 0} μ^{k+1}_{\Add\Gtl_q Q} (β, δ, …, δ) \\
&= p β β' β \quad [± \id\text{(S2)} ± \id\text{(S1)} ± β\text{(S3)}] \\
&∈ (B \htensor H) ⊕ (B \htensor R).
\end{align*}
In the first row, we have simply spelt out the definition of $ μ^1_{\Tw\Gtl_q Q} $. In the second row, we have evaluated all products. The first term $ p β β' β $ arises from $ k = 1 $. Further terms may arise from $ k ≥ 2 $, depending on the situation. If $ Q $ is geometrically consistent, then $ β $ (S3) terms may appear, stemming from first-out discrete immersed disks. If $ Q $ is not geometrically consistent, also $ \id $ (S1) and $ \id $ (S2) terms can appear. Either way, we see $ μ^1_{\Tw\Gtl_q Q} (β) $ lies in $ (B \htensor H) ⊕ (B \htensor R) $.

Regard an S2 situation between $ L $ and $ a $. Denote by $ α' $ and $ β' $ the associated angles. Then
\begin{equation*}
μ^1_{\Tw\Gtl_q Q} (\id_{L → a}) = - q α' - p β' ∈ B \htensor R.
\end{equation*}
This proves the claimed inclusion $ μ^1_{\Tw\Gtl_q Q} (H) ⊂ B \htensor R $. In the terminology of \autoref{sec:2Bkadeishvili}, this means $ E = 0 $ and we conclude
\begin{equation*}
H_q = \{h - Eh \running h ∈ B \htensor H\} = B \htensor H.
\end{equation*}
According to \autoref{def:2Bkadeishvili-deformed-counterpart}, the map $ φ: H_q → B \htensor H $ is the identity. This finishes the proof.
\end{proof}

The deformed decomposition for the hom space $ \Hom_{Q_1 ∪ \DefZigzagCat} (L_1, L_2) $ between two zigzag paths $ L_1, L_2 $ is simply the deformed decomposition described earlier. In case $ Q $ is geometrically consistent, this deformed decomposition was computed in \autoref{th:deformed-cohomology-basis}. In case $ Q = Q_M $ with odd $ M $, it was computed in \autoref{th:sphere-defsplitting-basis} and in case $ Q = Q_M $ with even $ M $, it was indicated in \autoref{sec:sphere-even}.

The deformed decomposition for the hom space $ \Hom_{Q_1 ∪ \DefZigzagCat} (a, b) $ between two arcs $ a, b ∈ Q_1 $ is trivially
\begin{equation*}
(H_q, μ^1_{Q_1 \DefZigzagCat} (B \htensor R), B \htensor R) = (B \htensor \Hom_{\Gtl Q} (a, b), 0, 0).
\end{equation*}
The deformed decomposition for the hom space $ \Hom_{Q_1 ∪ \DefZigzagCat} (a, L) $ between an arc $ a ∈ Q_1 $ and a zigzag path $ L $ depends on the choice one makes for the homological splitting of $ \Hom_{\Tw\Gtl Q} (a, L) $, but is entirely irrelevant to the present computation.

\subsection{Result components}
\label{sec:mirobjects-resultcomp}
In this section, we define and analyze result components for the products $ μ_{\HTw\Gtl_q Q} (m, h_N, …, h_1) $. The starting point is the description of the deformed decomposition from \autoref{sec:mirobjects-defsplitting}. In the present section, we introduce a notion of result components suitable for computing the products $ μ_{\HTw\Gtl_q Q} (m, h_N, …, h_1) $. We show that all result components fall into two classes which we call MD and MT result components.

According to the deformed Kadeishvili theorem of \autoref{sec:2Bkadeishvili}, the product of the morphisms $ h_1, …, h_N, m $ in the minimal model $ \HTw\Gtl_q Q $ is described in terms of Kadeishvili π-trees. Here the sequence $ h_1, …, h_N $ denotes type B cohomology basis elements $ h_i: L_i → L_{i+1} $ and $ m $ denotes a cohomology basis element $ m: L_{N+1} → a $. It is our task to evaluate all Kadeishvili π-trees $ T $ with inputs $ h_1, …, h_N, m $. For sake of convenience, we use the notation of \autoref{sec:splitting-situations} and \ref{sec:splitting-splitting} to denote angles, as opposed to the notation from \autoref{sec:sphere-zigzagcat} specific to $ Q = Q_M $. For instance, we denote the type B cohomology basis elements by $ α_3 + α_4 $. As usual, we start with a description of the possible terms that may possibly appear during evaluation of a Kadeishvili π-tree:

\begin{lemma}
\label{th:mirobjects-products-types}
Let $ (T, h_1, …, h_N, m) $ be a Kadeishvili π-tree with a type B cohomology basis elements $ h_1, …, h_N $ with $ h_i: L_i → L_{i+1} $ and $ m $ a cohomology basis element $ m: L_{N+1} → a $. Then:
\begin{itemize}
\item A proper subtree of $ T $ whose input morphisms only cover the morphisms between zigzag paths may only have result component $ β $ (A).
\item A proper subtree of $ T $ whose input morphisms also cover $ m $ vanishes.
\end{itemize}
Any nonvanishing result component $ r $ is derived either as a disk $ π_q μ^{≥3} (m, …) $, or as a product $ π_q μ^2 (\id\text{(S2)}, \allowbreak α_3/α_4 \text{(B)}) $ with direct inputs.
\end{lemma}

\begin{proof}
The statement on the subtrees that only involve morphisms between zigzag paths is familiar from the calculation of $ \H\DefZigzagCat $.

Regard now a product $ μ_{\Tw\Gtl_q Q} (m, m_k, …, m_1) $ where all $ m_1, …, m_k $ are of type $ β $ (A) or $ α_3/α_4 $ (B) and $ m $ is $ β $ (S1) or $ \id $ (S2). We claim this product lies in $ R $, apart from the case of $ μ^2_{\Add\Gtl_q Q} (\id \text{(S2)}, α_4) $ and all-in disks $ μ^{≥3}_{\Add\Gtl_q Q} (β \text{(S1)}, m_k, …, m_1) $. In these two exceptional cases, the product lies in $ H $.

The first part of checking this claim is to regard the case the product is a $ μ^2 $. The product is then of the form $ μ^2 (β \text{(S1)}, β \text{(A)}) $ or $ μ^2 (β \text{(S1)}, α_3/α_4 \text{(B)}) $ or $ μ^2 (\id\text{(S2)}, β \text{(A)}) $ or $ μ^2 (\id\text{(S2)}, α_3/α_4 \text{(B)}) $. The first case yields $ β \text{(S3)} ∈ R $, the second and third type of composition are impossible, the fourth case yields $ β \text{(S1)} ∈ H $.

The second part of checking the claim is to regard the case of a disk $ μ^{≥3} $. If it concerns an all-in disk, then the result is an arc identity $ \id_{L → a} ∈ R $ or $ \id_{L → a} ∈ H $. A final-out disk is impossible, since $ β $ (S1) is an indecomposable angle. If it concerns a first-out disk, then the first angle of the disk may be a $ δ $-morphism or $ β $ (A). In both cases, the result is of the type $ μ^{≥3} = β \text{(S3)} ∈ R $.

Finally, we draw two conclusions: Any h-tree consuming $ m $ has vanishing result. A given π-tree $ T $ with nonvanishing result must therefore consume $ m $ directly at the root. This finishes the proof.
\end{proof}

In analogy with \autoref{sec:resultcomp}, we can define result components also for Kadeishvili π-trees with inputs $ h_1, …, h_N, m $. Virtually the same definition can be applied.

\begin{definition}
Let $ (T, h_1, …, h_N, m) $ be a Kadeishvili π-tree with a type B cohomology basis elements $ h_1, …, h_N $ with $ h_i: L_i → L_{i+1} $ and $ m $ a cohomology basis element $ m: L_{N+1} → a $. Then a \emph{result component} of $ (T, h_1, …, h_N, m) $ is defined in analogy with \autoref{def:resultcomp-def}. The set of result components of all π-trees, ranging over all choices of $ h_1, …, h_N $ and $ m $ and $ T $, is denoted $ \PiTrM $.
\end{definition}

By \autoref{th:mirobjects-products-types}, result components of $ μ(m, h_N, …, h_1) $ can be split into two types: those which are derived from a disk and those which are derived from a product $ π_q μ^2 $ with direct inputs. In analogy with \autoref{sec:subdisk-types}, we shall give these two types the names mirror disks and mirror triangles, respectively.

\begin{definition}
\label{def:mirobjects-minmodel-rMDMT}
A result component $ r ∈ \PiTrM $ is a
\begin{itemize}
\item \emph{MD result component} if it is derived as $ π_q μ^{≥3} (m, …) $.
\item \emph{MT result component} if it derived as $ π_q μ^2 (\id\text{(S2)}, α_3/α_4 \text{(B)}) $.
\end{itemize}
We denote the set of MD and MT result components by $ \MDr $ and $ \MTr $, respectively.
\end{definition}

\begin{figure}
\centering
\begin{subfigure}{0.4\linewidth}
\centering
\begin{tikzpicture}
\path node (A) {$ m $} node[right of=A] (B) {$ h_N $} node[right of=B] (dots) {…} node[right of=dots] (C) {$ h_1 $}
node[below of=dots] (D) {…} edge (B) edge (C)
node[below left of=D] {$ π_q μ^{≥3} $} edge (A) edge (D);
\end{tikzpicture}
\caption{MD result component}
\label{fig:mirobjects-products-MD}
\end{subfigure}
\begin{subfigure}{0.4\linewidth}
\centering
\begin{tikzpicture}[node distance=2cm]
\path node (A) {$ \id $(S2)} node[right of=A] (B) {$ α_3/α_4 \text{(B)} $}
node[below right of=A] {$ π_q μ = β $(S1)} edge (A) edge (B);
\end{tikzpicture}
\caption{MT result component}
\end{subfigure}
\caption{How MD and MT result components are derived}
\label{fig:mirobjects-products-MDMT}
\end{figure}
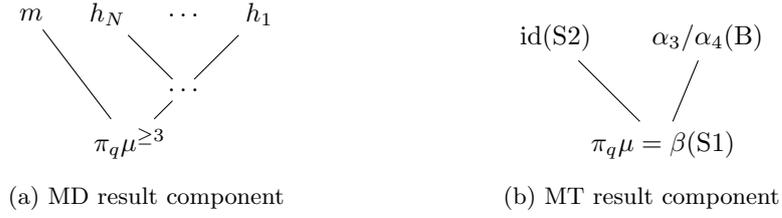

The distinction between MD and MT result components is depicted in \autoref{fig:mirobjects-products-MDMT}. Result components of these two types compute the products $ μ_{\HTw\Gtl_q Q} (m, h_N, …, h_1) $. The two names MD (mirror disk) and MT (mirror triangle) have been chosen in order to reflect their use in computing the deformed mirror in the third paper of the present series.

\subsection{The higher products}
\label{sec:mirobjects-products}
In this section, we compute the desired products of the kind $ μ_{\H\Tw\Gtl_q Q} (m, h_1, …, h_N) $. The starting point for the computation is the classification of result components from \autoref{sec:mirobjects-resultcomp}. In the present section, we introduce MD and MT disks with the aim of expressing the products as an enumeration of disks. We define a matching, the subdisk mapping, between MD/MT result components and MD/MT disks. In \autoref{th:mirobjects-products-th}, we collect the desired description of the products.

The first step of the present section is to define the notion of MD and MT disks, meant to capture MD and MT result components geometrically. The key observation is that
\begin{itemize}
\item The final input of an MD result component is an odd morphism of type $ m = β \text{(S1)}: L_{N+1} →a $ and the output is the morphism $ \id\text{(S2)}: L_1 → a $.
\item The final input of an MT result component is an even morphism of type $ m = \id\text{(S1)}: L_{N+1}: a $ and the output is the morphism $ β\text{(S1)}: L_1 → a $.
\end{itemize}
In either case, we see that the two zigzag curves $ \smooth L_1 $, $ \smooth L_{N+1} $ and the arc $ a $ have a triple intersection at the midpoint of the arc $ a $. In order to capture MD/MT result components by means of disks, the disk therefore also need to have nontransversal input sequence. This means the correct definition of MD/MT disks cannot be inferred from parallels with the relative Fukaya pre-category $ \relpreFuk Q $. Instead, the parallel needs to be drawn with the relative Fukaya category $ \relFuk Q $. The products of this category can be determined on a best-effort basis by performing Hamiltonian deformations on the involved curves. In the present context of the products $ μ(m, h_N, …, h_1) $, this means we have to choose a Hamiltonian deformation of some of the zigzag paths or arcs in order to guess the correct notion of MD/MT disks.

\begin{figure}
\centering
\begin{subfigure}{0.3\linewidth}
\centering
\begin{tikzpicture}
\path (0, 0) coordinate (A) (-1, 1) coordinate (B) (0, 2) coordinate (C) (1, 1) coordinate (D);
\path[draw, semithick] (A) to[out=300, in=100, looseness=0.5] ++(0.1, -0.2);
\path[draw, semithick, ->] (A) to[out=240, in=80, looseness=0.5] ++(-0.1, -0.2);
\path[draw, semithick, ->] (B) to[out=150, in=350, looseness=0.5] ++(-0.2, 0.1);
\path[draw, semithick] (B) to[out=210, in=10, looseness=0.5] ++(-0.2, -0.1);
\path[draw, semithick, ->] (C) to[out=60, in=260, looseness=0.5] ++(0.1, 0.2);
\path[draw, semithick] (C) to[out=120, in=280, looseness=0.5] ++(-0.1, 0.2);
\path[draw, semithick, ->] (D) to[out=330, in=170, looseness=0.5] ++(0.2, -0.1);
\path[draw, semithick] (D) to[out=30, in=190, looseness=0.5] ++(0.2, 0.1);
\path (A) to[out=120, in=330] coordinate[pos=0.3] (out) (B) to[out=30, in=240] (C) to[out=300, in=150] (D) to[out=210, in=60] coordinate[pos=0.7] (m) cycle;
\begin{scope}
\clip (out) -- (B) -- (C) -- (D) -- (m) -- cycle;
\path[fill=gray!50] (A) to[out=120, in=330] (B) to[out=30, in=240] (C) to[out=300, in=150] (D) to[out=210, in=60] cycle;
\end{scope}
\path[draw, semithick] (A) to[out=120, in=330] node[midway, above] {$ L' $} (B) to[out=30, in=240] (C) to[out=300, in=150] (D) to[out=210, in=60] node[midway, above] {$ L $} cycle;
\path[draw, semithick] ($ (out)!0.5!(m) $) ++(left:1) -- ++(right:2) node[pos=0.9, above] {$ a $};
\path[fill] \foreach \i in {A, B, C, D, out, m} {(\i) circle[radius=0.05]};
\path (out) node[below left, inner sep=0.1em] {\small out};
\path (m) node[below right, inner sep=0.1em] {\small $ m $};
\path (B) node[above] {\small $ b $};
\path (C) node[above, shift={(up:0.1)}] {\small $ b $};
\path (D) node[above] {\small $ b $};
\end{tikzpicture}
\caption{Disk}
\label{fig:mirobjects-products-disk}
\end{subfigure}
\begin{subfigure}{0.3\linewidth}
\centering
\begin{tikzpicture}
\path (0, 0) coordinate (A) (-1, 1) coordinate (B) (0, 2) coordinate (C) (1, 1) coordinate (D);
\path[draw, semithick] (A) to[out=300, in=100, looseness=0.5] ++(0.1, -0.2);
\path[draw, semithick, ->] (A) to[out=240, in=80, looseness=0.5] ++(-0.1, -0.2);
\path[draw, semithick, ->] (B) to[out=150, in=350, looseness=0.5] ++(-0.2, 0.1);
\path[draw, semithick] (B) to[out=210, in=10, looseness=0.5] ++(-0.2, -0.1);
\path[draw, semithick, ->] (C) to[out=60, in=260, looseness=0.5] ++(0.1, 0.2);
\path[draw, semithick] (C) to[out=120, in=280, looseness=0.5] ++(-0.1, 0.2);
\path[draw, semithick, ->] (D) to[out=330, in=170, looseness=0.5] ++(0.2, -0.1);
\path[draw, semithick] (D) to[out=30, in=190, looseness=0.5] ++(0.2, 0.1);
\path (A) to[out=120, in=330] coordinate[pos=0.3] (out) (B) to[out=30, in=240] (C) to[out=300, in=150] (D) to[out=210, in=60] coordinate[pos=0.7] (m) cycle;
\begin{scope}
\clip (out) -- (m) -- (A) -- cycle;
\path[fill=gray!50] (A) to[out=120, in=330] (B) to[out=30, in=240] (C) to[out=300, in=150] (D) to[out=210, in=60] cycle;
\end{scope}
\path[draw, semithick] (A) to[out=120, in=330] node[midway, above] {$ L' $} (B) to[out=30, in=240] (C) to[out=300, in=150] (D) to[out=210, in=60] node[midway, above] {$ L $} cycle;
\path[draw, semithick] ($ (out)!0.5!(m) $) ++(left:1) -- ++(right:2) node[pos=0.9, above] {$ a $};
\path[fill] \foreach \i in {A, B, C, D, out, m} {(\i) circle[radius=0.05]};
\path (out) node[below left, inner sep=0.1em] {\small out};
\path (m) node[below right, inner sep=0.1em] {\small $ m $};
\path (A) node[below, shift={(down:0.1)}] {\small $ b $};
\end{tikzpicture}
\caption{Triangle}
\label{fig:mirobjects-products-triangle}
\end{subfigure}
\caption{The two types of disks we expect to contribute to $ μ(m, h_N, …, h_1) $}
\label{fig:mirobjects-products-MDMTdisk}
\end{figure}
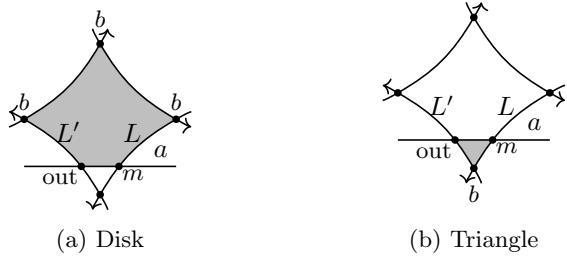

There is one particular Hamiltonian deformation of the arcs that gives the correct notion of MD/MT disks: Push every arc $ a $ a little into the neighboring clockwise polygon, leaving the zigzag curves in place. This specific Hamiltonian deformation simultaneously resolves all triple intersections between zigzag curves in $ Q $ and arcs. It predicts us to find disks of two types, depicted in \autoref{fig:mirobjects-products-MDMTdisk}. We will verify in \autoref{th:mirobjects-products-bijection} that it is the correct Hamiltonian deformation to capture the products $ μ_{\HTw\Gtl_q Q} (m, h_N, …, h_1) $. If we had chosen a different homological splitting in \autoref{def:mirobjects-splitting-def}, we would have needed a different Hamiltonian deformation.

In \autoref{def:mirobjects-minmodel-MDMT} we provide a rigorous definition of MD/MT disks. For a given arc $ a $, we have denoted by $ L $ and $ L' $ the two zigzag paths which depart from $ a $. Geometrically speaking, $ \smooth L $ and $ \smooth L' $ are the two zigzag curves which intersect at the midpoint of $ a $. This is depicted in \autoref{th:mirobjects-products-bijection}.

\begin{definition}
\label{def:mirobjects-minmodel-MDMT}
An \emph{MD disk} (mirror disk) is a CR disk with $ N ≥ 0 $ inputs $ h_1, …, h_N $ whose
\begin{itemize}
\item inputs $ h_1, …, h_N $ are all odd and do not contain co-identities,
\item output is even and not an identity,
\item zigzag segments all run clockwise,
\end{itemize}
which has undergone the following surgery: The output mark, located at a certain arc $ a $, has been cut off. The odd morphism at $ a $ is added as final input, and the even morphism at $ a $ is indicated as new output.

An \emph{MT disk} (mirror triangle) is a triangle sitting between the deformed arc $ a $ and the intersection of the two zigzag curves $ \smooth L $ and $ \smooth L' $ intersecting at $ a $.

We denote by $ \MDd $ and $ \MTd $ the sets of MD and MT disks, respectively.
\end{definition}

In the remainder of this section, we show that the product $ μ_{\HTw\Gtl_q Q} (m, h_N, …, h_1) $ is indeed given by counting MD and MT disks with inputs $ h_1, …, h_N, m $. The first step is to map a given result component to an MD or MT disk. The description of this mapping is eased by the classification of result components given in \autoref{def:mirobjects-minmodel-rMDMT}. According to this classification, result components can be categorized into what we have called MD and MT result components. An MD result component $ r ∈ \PiTrM $ is necessarily derived as $ π_q μ^{≥3}_q (β \text{(S1)}, m_l, …, m_1) $ where all $ m_1, …, m_k $ are of type $ α_3/α_4 $ (B) or $ β $ (A). In particular, every $ m_i $ is the result component of an h-tree with inputs being a subsequence of $ h_1, …, h_N $. By \autoref{sec:subdisk-subdisk}, every $ m_i $ comes with an associated subdisk. To define the subdisk of the MD result component $ r $, we essentially glue together the subdisks of the $ m_i $. The precise definition of subdisks for MD and MT result components reads as follows:

\begin{definition}
\label{def:mirobjects-products-subdisk}
Let $ r $ be a result component of a π-tree $ (T, h_1, …, h_N, m) $. Then its \emph{subdisk} $ \Subdisk(r) $ is the drawing defined as follows:
\begin{itemize}
\item If $ T $ is derived as $ π_q μ^2 (\id\text{(S2)}, α_3 + α_4) $ with both inputs being direct, then its subdisk $ \Subdisk(r) $ is the infinitesimal triangle sitting between the Hamiltonian deformed arc $ a $ and the intersection of the two input zigzag curves.
\item If $ T $ is derived as $ π_q μ^{≥3} (β \text{(S1)}, m_k, …, m_1) $, then its subdisk $ \Subdisk(r) $ is given by connecting the subdisks of $ m_1, …, m_k $ and finishing with the input $ β $ (S1). The output mark is at $ \id $ (S2) and lies infinitesimally apart from the final input $ β $ (S1).
\end{itemize}
\end{definition}

In \autoref{th:mirobjects-products-bijection}, we claim that every MD or MT disk is the subdisk of a unique single result component, in analogy with \autoref{th:subdisk-minmodel-bijection}. We have to restrict MD disks to those with at least two inputs because subdisks of result components always have at least two inputs.

\begin{lemma}
\label{th:mirobjects-products-bijection}
The subdisk of an MD result component is an MD disk. The subdisk of an MT result component is an MT disk. Denoting by $ \MDd^{≥2} $ the set of MD disks with at least two inputs, the map $ \Subdisk: \PiTrM → \MDd^{≥2} \disjoint \MTd $ is a bijection.
\end{lemma}

\begin{proof}
We divide the proof into two parts: First we comment on MT disks, then we comment on MD disks.

For MT result components and MT disks, there is not much to say: An MT disk $ D $ is a small triangle located on the counterclockwise side of an. It immediately gives rise to two morphisms $ α_3 + α_4 $ and $ m = \id\text{(S2)} $ which multiply to $ π_q μ^2 = β \text{(S1)} $. This gives an MT result component whose subdisk is $ D $ again. This shows that $ \Subdisk $ matches MT result components bijectively with MT disks.

For MD disks, it is our task to recover a tree and result component for a given MD disk $ D ∈ \MDd^{≥2} $. After our long journey in \autoref{sec:classification}, we content ourselves with merely a brief description: Find the indecomposable narrow locations of $ D $ and place them in a tree, ordered by inclusion. The root of this tree is the artifical narrow location, which we denoted $ (1, |D|) $ in \autoref{sec:classification-narrowtree}. According to \autoref{sec:classification-subresult}, the children $ C_1, …, C_k $ of the root all come with subtrees $ T_1, …, T_k $ and subresults $ r_1, …, r_k $. Each tree $ T_i $ is an h-tree which consumes part of the inputs $ h_1, …, h_N $, and $ r_i $ is a result component of $ T_i $. In fact, all $ r_i $ are of $ β $ (A) type.

With this in mind, we are ready to associate with $ D $ a tree $ T $ and a result component $ r $ of $ T $, such that $ \Subdisk(r) = D $. To construct the tree $ T $, put all trees $ T_1, …, T_k $ next to each other, insert a root $ π_q μ^{≥3} $ and connect all outputs of $ T_1, …, T_k $ together with all remaining $ b $ inputs and the input $ m $ with the root of $ T $. This gives the desired tree $ T $. The result component $ r $ is the $ π_q μ^{≥3} $ result component of $ \id $ (S2) type simply given by the data of result components $ r_1, …, r_k $ on each $ T_i $, bound together with the all remaining morphisms and $ m $ following the geometry of $ D $. The result component $ r $ defined this way satisfies $ \Subdisk(r) = D $. This shows that $ \Subdisk $ maps surjectively onto $ \MDd^{≥2} \disjoint \MTd $.
\end{proof}

We proceed by checking signs.

\begin{lemma}
\label{th:mirobjects-products-signs}
Let $ r ∈ \PiTrM $ be a result component. Then the absolute sign of $ r $ equals the Abouzaid sign $ \Abouzaid(\Subdisk(r)) $.
\end{lemma}

\begin{proof}
First we check the case of MD result components and second the case of MT result components, both focusing on the geometrically consistent case. Third we comment on the case of the sphere dimers $ Q_M $.

For the first part, regard an MD result component $ r $. Let $ T $ be the π-tree from which $ r $ stems. Then $ T $ has shape as depicted in \autoref{fig:mirobjects-products-MD}. Let $ T_1, …, T_k $ be the children of the root of $ T $, not counting the direct input $ m $. Then $ r $ is derived as $ π_q μ^{≥3} (m, r_k, …, r_1) $, with $ r_1, …, r_k $ being resukt components of $ T_1, …, T_k $. It is our task to evaluate the sign of $ r $. Our procedure is analogous to \autoref{sec:classification-signs}.

Let us compute the sign of the result component $ r $. Since it is derived as the product $ μ^{≥3} (m, r_k, …, r_1) $ and the disk is all-in, its total sign is the sum of: the $ \# $ signs of the $ δ $ insertions, the signs of the result components $ r_1, …, r_k $ and the sign of $ m = (-1)^{\#β} β $.

On the other hand, let us compute the Abouzaid sign of $ \Subdisk(r) $. By \autoref{def:mirobjects-products-subdisk}, the disk $ \Subdisk(r) $ is formed by tying together the subdisks of $ r_1, …, r_k $. Correspondingly, its Abouzaid sign is the sum of the signs coming from odd counterclockwise intersections within the subdisks of $ r_1, …, r_k $, plus $ \# $ signs of the zigzag curve segments lying between two neighboring subdisks. By \autoref{th:classification-signs-htrees} and \autoref{def:classification-signs-abouzaid}, the sign coming from odd counterclockwise intersections in the subdisk of $ r_i $ is equal to the absolute sign of $ r_i $. In other words, the total Abouzaid sign of $ \Subdisk(r) $ is the sum of the absolute signs of $ r_1, …, r_k $ and the $ \# $ signs. This $ \# $ sign already includes sign $ \#β $ of final input angle $ β $. Finally, we conclude that both signs are equal.

For the second part, regard an MT result component $ r $. It is derived as $ π_q μ^2 (\id\text{(S2)}, α_3/α_4 \text{(B)}) $, where $ \id\text{(S2)} $ comes from an S2 situation and $ α_3/α_4 $ comes from a B situation such that both angles are composable. Recall that $ α_3/α_4 $ is merely an abbreviation for the morphism $ φ^{-1} ((-1)^{\#α_3 + 1} α_3 + (-1)^{\#α_4} α_4) $. The relevant result of the product is $ μ^2 (\id\text{(S2)}, (-1)^{\#α_4} α_4) = (-1)^{\#β + 1} β \text{(S1)} $, noting that $ β $ is the same angle in $ Q $ as $ α_4 $. Relative to the sign of the cohomology basis element $ (-1)^{\#β} $, the MT result component $ r $ has a total sign of $ -1 $.

On the other hand, regard the subdisk $ \Subdisk(r) $ associated with $ r $. It is a small triangle with two inputs and one output. The first input $ α_3/α_4 $ is odd but its zigzag curves are oriented clockwise with the triangle, the second input is even, and the output is odd and its target zigzag curve is oriented counterclockwise with the triangle. The Abouzaid sign of $ \Subdisk(r) $ is therefore $ -1 $. We conclude that both signs agree.

For the third part, let us comment on the case of $ Q_M $. The choice of signs in the cohomology basis elements, $ (-1)^{\#(i+M) + 1} β_i' + (-1)^{\#i} β_i $ and $ id_i' $ is analogous to the choice of sign in the geometrically consistent case. In particular, the sign computations for MD and MT result components carry over to the case of $ Q_M $ without change.
\end{proof}

In \autoref{th:mirobjects-products-th}, we describe the desired products $ μ_{\HTw\Gtl_q Q} (m, h_N, …, h_1) $. It turns out that for any value of $ N ≥ 0 $, this product is determined by enumerating MD and MT disks with inputs $ h_1, …, h_N, m $. In case $ N ≥ 1 $, this description follows from \autoref{th:mirobjects-products-bijection} and \autoref{th:mirobjects-products-signs}. In case $ N = 0 $, the description follows from inspection of the differential $ μ^1_{\Tw\Gtl_q Q} (m) $. The notation $ Q_{i+2} ∈ ℂ⟦q_1, …, q_M⟧ $ and $ \#\#(i-M+1) $ appearing in \autoref{eq:mirobjects-products-oddm1} is taken over from \autoref{eq:sphere-defzigzagcat-shorthand}. The index $ i $ denotes the index such that $ α_i $ is the same angle in $ Q = Q_M $ as $ β $.

\begin{proposition}
\label{th:mirobjects-products-th}
Let $ Q $ be a geometrically consistent dimer or standard sphere dimer $ Q = Q_M $ with $ M ≥ 3 $, as in \autoref{conv:mirobjects-dimer}. Let $ h_1, …, h_N $ be a sequence of $ N ≥ 0 $ type B cohomology basis elements $ h_i: L_i → L_{i+1} $. Let $ m = (-1)^{\#β} β: L_{N+1} → a $ be another cohomology basis element. Then we have
\begin{equation}
\label{eq:mirobjects-products-th}
μ_{\HTw\Gtl_q Q} (m, h_N, …, h_1) = \sum_{\substack{D ∈ \MDd \disjoint \MTd \\ \text{with inputs } h_1, …, h_N, m}} (-1)^{\Abouzaid(D)} \punctures(D) \disktarget(D).
\end{equation}
The differential $ μ^1_{\HTw\Gtl_q Q} (m) $ is accurately described by this equality. More explicitly, it is given as follows:
\begin{itemize}
\item If $ Q $ is a geometrically consistent dimer, then $ μ^1_{\HTw\Gtl_q Q} (m) $ vanishes.
\item If $ Q = Q_M $ with $ M ≥ 3 $ odd, then for $ m = \id_{L → a} $ we have $ μ^1_{\HTw\Gtl_q Q} (m) = 0 $. For $ m = (-1)^{\#β} β $ (S1), let $ i $ be such that $ a_i = t(β) $. Then
\begin{equation}
\label{eq:mirobjects-products-oddm1}
μ^1_{\HTw\Gtl_q Q} ((-1)^{\#β} β) = (-1)^{\#\# (i-M+1)} Q_{i+2} \id_{L → a}.
\end{equation}
\item If $ Q = Q_M $ with $ M ≥ 4 $ even, then $ μ^1_{\HTw\Gtl_q Q} (m) $ vanishes.
\end{itemize}
\end{proposition}

\begin{proof}
The entire computation is completely analogous to \autoref{sec:subdisk-minmodel}: In \autoref{sec:subdisk}, we have shown how to match result components for $ μ_{\H\DefZigzagCat} $ with CR/ID/DS/DW disks. In the present \autoref{th:mirobjects-products-bijection}, we have shown how to match result components for $ μ(m, h_N, …, h_1) $ with MT/MD disks in case $ N ≥ 1 $. In \autoref{th:mirobjects-products-signs}, we have checked that the sign of a result component agrees with the Abouzaid sign of its associated subdisk. We conclude that the claimed product description \eqref{eq:mirobjects-products-th} in case $ N ≥ 1 $ follows as in \autoref{sec:subdisk-minmodel}.

It remains to comment on the case $ N = 0 $. This entails determining the differential $ μ^1_{\HTw\Gtl_q Q} (m) $ explicitly, where $ m $ is a cohomology basis element $ m: L → a $. We need to distinguish cases on whether $ Q $ is geometrically consistent or $ Q = Q_M $. If $ Q $ is geometrically consistent, then $ μ^1_{\Tw\Gtl_q Q} (m) $ lies in $ B \htensor R $, as we have seen in the proof of \autoref{th:mirobjects-defsplitting-th}. Therefore
\begin{equation*}
μ^1_{\HTw\Gtl_q Q} (m) = π_{H_q} μ^1_{\Tw\Gtl_q Q} (m) = 0.
\end{equation*}
If $ Q = Q_M $ with $ M ≥ 3 $ odd, then we can calculate the differential easily by looking at the $ \id $ (S2) terms appearing in $ μ^1_{\Tw\Gtl_q Q} (m) $. For $ m = \id_{L → a} $ there are no such terms, but for $ m = (-1)^{\#β} β $ we find a single type $ \id $ (S2) term, namely
\begin{multline*}
μ^M_{\Add\Gtl_q Q} ((-1)^{\#β} β, (-1)^{\#(i-1)} q_{i-1} α_{i-1}', (-1)^{\#(i-2)} α_{i-2}, …, (-1)^{\#(i-M+1)} α_{i-M+1}) \\
= (-1)^{\#\# (i-M+1)} Q_{i+2} \id_{L → a}.
\end{multline*}
The differential $ μ^1_{\HTw\Gtl_q Q} (m) $ is given by projecting this term to $ H_q $. Since the term already lies in $ H_q $, we conclude the desired formula for $ μ^1_{\HTw\Gtl_q Q} (m) $. This finishes the case distinction for $ Q $ and proves the explicit description of $ μ^1_{\HTw\Gtl_q Q} (m) $ in all cases.

It remains to check \eqref{eq:mirobjects-products-th} in case $ N = 0 $. This entails reinterpreting the explicit description of $ μ^1_{\HTw\Gtl_q Q} (m) $ in terms of MD disks. This is an easy exercise and we finish the proof here.
\end{proof}

This finishes the computation of the desired products $ μ_{\HTw\Gtl_q Q} (m, h_N, …, h_1) $.

\section{Discussion}
\label{sec:discussion}
In this section, we provide more explanation on the results of the present paper.

\subsection{Relation to the literature}
\label{sec:literature}
In this section, we list seven existing papers which have contributed to cornerstones of the present paper. For every paper, we explain its core insights, the way our paper builds on that paper, and the aftermath our paper provides to that paper. The first cornerstone of this paper is the category $ \Gtl Q $ as a model for the wrapped Fukaya category. The second cornerstone is the use of the deformation $ \Gtl_q Q $ as a candidate model for the relative Fukaya category. The third cornerstone is the observation that the objects in $ \Gtl_q Q $ and $ \Tw\Gtl_q Q $ have curvature. The fourth cornerstone is the classification of objects in $ \Tw\Gtl Q $. The fifth cornerstone is the uncurving procedure for band objects in $ \Tw\Gtl_q Q $. The sixth cornerstone is the category of zigzag paths $ \DefZigzagCat ⊂ \Tw\cat C_q $. The seventh cornerstone is the choice of homological splitting for $ \ZigzagCat $ and the subsequent lengthy calculation of the minimal model, ending up in a beautiful geometric interpretation. These seven cornerstones are related to the following seven papers:

\begin{itemize}
\item Bocklandt \cite{Bocklandt} introduces gentle algebras as model for Fukaya categories.
\item Barmeier, Schroll and Wang (based on \cite{Barmeier-Wang}) describe deformations of gentle algebras.
\item Seidel \cite{Seidel-relative} foresees the need of curvature.
\item Haiden, Katzarkov and Kontsevich \cite{HKK} classify objects of derived gentle algebras.
\item Lowen and Van den Bergh \cite{Lowen-vdB} demonstrate strategies to gauge away curvature.
\item Lekili, Perutz and Polishchuk \cite{Lekili-Perutz, Lekili-Polishchuk} calculate with zigzag paths.
\item Bocklandt \cite{Bocklandt-book} demonstrates how to perform minimal model computations.
\end{itemize}

Those papers dealing with the A-side of mirror symmetry all depart from either the geometric model $ \Fuk Q $ or the discrete model $ \Gtl Q $. Moreover, we can classify the papers according to whether they work with deformations or not. This gives the following diagram:
\begin{center}
\begin{tabular}{c|cc}
starting point & non-deformed & deformed \\\hline
geometric & \cite[Appendix B]{Bocklandt} & \cite{Seidel-relative}, \cite{Lekili-Perutz}, \cite{Lekili-Polishchuk} \\
discrete & \cite{Bocklandt}, \cite{Bocklandt-book}, \cite{HKK} & \emph{this paper}
\end{tabular}
\end{center}
The present paper fills this square by departing from a deformation of the discrete model $ \Gtl Q $.

\subsubsection{Bocklandt-Abouzaid}
\label{sec:literature-BA}
In \cite{Bocklandt}, Bocklandt introduces the $ A_∞ $-structure on $ \Gtl Q $ and proposes it as discrete model for the wrapped Fukaya category $ \wFuk Q $. In that paper's appendix, Abouzaid computes the minimal model of an arc system as part of $ \wFuk Q $ and obtains indeed the gentle algebra. This way, Bocklandt and Abouzaid contribute towards the very starting point of this paper, specifically \autoref{sec:splitting}. The paper approaches the A-side via both the non-deformed discrete and non-deformed geometric side.

As first step of the paper, Bocklandt defines an $ A_∞ $-structure on $ \Gtl Q $. Abouzaid then shows that this $ \Gtl Q $ is in fact a discrete model for $ \wFuk Q $. He regards the wrapped Fukaya category, as defined in \cite{Abouzaid-Seidel}. He discovers that one can pass with relative ease to the minimal model if one restricts to those string objects given by an arc system. More concretely, he shows that on $ A_∞ $-level one can get rid of the so-called continuation map.

With this in mind, the work of Bocklandt-Abouzaid is a non-deformed prototype for our result: If the discrete $ \Gtl Q $ provides a model for the geometric $ \wFuk Q $, then the deformation $ \Gtl_q Q $ is necessarily a model for a certain deformation of $ \wFuk Q $ (see \autoref{sec:whyshould-candidate}).

In contrast to Abouzaid's appendix, our calculations have to depart from the discrete side. Indeed, a deformed wrapped Fukaya category does not exist as of yet, so that we cannot work ourselves from geometric to discrete (see \autoref{sec:whyshould-alternative}).

Our paper provides a new proof of Abouzaid's appendix in \cite{Bocklandt}, at least on the subset of zigzag paths. Indeed, we show that $ \H\Tw\Gtl_q Q $ matches with the relative Fukaya category. Both categories are deformations over the same deformation base $ ℂ⟦Q_0⟧ $, i.e.~they have one deformation parameter per puncture in $ Q $. As soon as we restrict both $ \H\Tw\Gtl_q Q $ and $ \relFuk Q $ to the special fiber $ q = 0 $, we hold in our hands an explicit matching between $ \H\Tw\Gtl Q $ and $ \Fuk Q $, at least on the category of zigzag paths. This recovers part of the result of Bocklandt and Abouzaid.


\subsubsection{Barmeier-Schroll-Wang}
\label{sec:literature-BSW}
Intriguingly, Barmeier, Schroll and Wang are working on $ A_∞ $-deformations on Fukaya categories as well, in parallel to the present paper. The subject of their work is known to the author, so we would like to point out a few relations.

In the article \cite{Barmeier-Wang}, Barmeier and Wang investigate deformations of quiver algebras with relations. The idea behind the new work of Barmeier-Schroll-Wang is to apply their methods to topological Fukaya categories as well.

To understand their line of thought, we should look into the work \cite{HKK} of Haiden, Katzarkov and Kontsevich. They define topological Fukaya categories $ \Fuk(S, M) $ also for marked surfaces $ (S, M) $ beyond our notion of punctured surfaces. Indeed, \cite{HKK} allows the surface to have a boundary instead of punctured, and the boundary is supposed to consist alternatingly of markings and “boundary arcs”. In the simplest case without boundary arcs, their notion is equivalent to our punctured surfaces. In the case with at least one boundary arc, the topological Fukaya category $ \Fuk(S, M) $ however allows for a very explicit model: a graded algebra without differential and higher products. This is the point where the deformation theory of \cite{Barmeier-Wang} comes into play.

The work of Barmeier-Schroll-Wang yields results complementary to ours, namely deformations in the case every boundary component has at least one boundary arc. Since our case of $ \Gtl Q $ is an $ A_∞ $-localization of the case with boundary arcs, it will be interesting to speculate about the relations between our work and Barmeir-Schroll-Wang's.

\subsubsection{Seidel}
\label{sec:literature-seidel}
In \cite{Seidel-relative}, Seidel introduces the notion of relative Fukaya categories. He departs from the exact Fukaya category and explains how to work relative to a divisor. He foresees the necessity to use curvature for those Lagrangians that have teardrops intersecting with the divisor, while according to him all other Lagrangians would be free of curvature. Through his understanding of curvature in relative Fukaya categories, Seidel \cite{Seidel-relative} contributes to \autoref{sec:fukaya} and \ref{sec:uncurving}.


While Seidel does not provide anything explicit in case of punctured surfaces, his ideas carry over without difficulty: The divisor becomes a finite collection of points, which in our paper correspond to the punctures $ M ⊂ S $. Each immersed disk should be weighted with the power of a deformation parameter whose exponent is the intersection number of the disk with the divisor $ D $.

Seidel envisions those Lagrangians to be infinitesimally curved which have teardrops intersecting the divisor. This expectation has fueled our expectations towards uncurvability of objects in $ \Tw\cat C_q $: According to Seidel, we should expect that those band objects which are topologically nontrivial and do not bound a teardrop in $ S $ are uncurvable, while those with teardrop in $ S $ are inherently curved.

Seidel's definition provides a deformed Fukaya category of pre-category style: Its higher products are only defined on transversal sequences. At the time of Seidel's paper, it was not clear how to turn this definition into an actual category. This was accomplished in general only 20 years later by Sheridan and Perutz \cite{Sheridan-Perutz}. Yet, their construction relies on the Hamiltonian deformation approach, which renders the $ A_∞ $-structure on the non-transversal sequences very complicated.

The aftermath of our paper is a very down-to-earth description of the relative Fukaya category, at least on the subset of zigzag paths: We describe explicitly all the immersed disks one needs for its definition, also on all non-transversal sequences. A small caution: Technically, we cannot prove that our explicit category $ \H\DefZigzagCat $ is indeed (a subcategory of) the relative Fukaya category, but its higher products on the transversal sequences suggest so.

We confirm Seidel's expectations regarding curvature in the relative Fukaya category in \autoref{sec:uncurving}. We also extend the width of Seidel's deformation in that we use one deformation parameter per puncture. It would be interesting to reintroduce Seidel's relative Fukaya category with more deformation parameters even in the higher dimensional case.

\subsubsection{Haiden-Katzarkov-Kontsevich}
\label{sec:comparison-HKK}
In \cite{HKK}, Haiden, Katzarkov and Kontsevich famously analyze stability conditions on partially wrapped Fukaya categories. Twisted complexes of the gentle algebra $ \Gtl Q $ serve as model for their actual work. Their work contributes to \autoref{sec:uncurving}, departing from the non-deformed discrete side.

As first step, Haiden, Katzarkov and Kontsevich introduce a notion of marked surfaces. In a marked surface, each boundary component is supposed to consist alternatingly of markings and “boundary arcs”. Those marked surfaces where every $ S^1 $ boundary component is fully marked are precisely the punctured surfaces we use in the present paper.
\begin{figure}
\centering
\begin{tikzpicture}
\begin{scope}
\path[draw] (0, 0) circle[radius=0.7];
\path[draw, very thick] ++(-20:0.7) arc(-20:20:0.7);
\path[draw, very thick] ++(100:0.7) arc(100:140:0.7);
\path[draw, very thick] ++(200:0.7) arc(200:240:0.7);
\path (1.2, 0) node {\Large $ \rightsquigarrow $} (2.5, 0) node {$ \mathcal{F}_A (S) $};
\end{scope}
\begin{scope}[shift={(4.5, 0)}]
\path[draw, very thick] (0, 0) circle[radius=0.7];
\path (1.2, 0) node {\Large $ \rightsquigarrow $} (2.5, 0) node {$ \Gtl Q $};
\end{scope}
\end{tikzpicture}
\caption{Markings on a boundary component after \cite{HKK}}
\label{fig:literature-hkk-marking}
\end{figure}
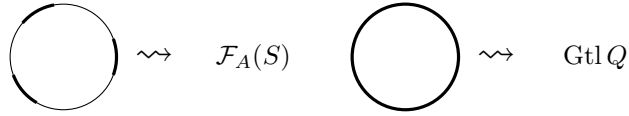

Using arc systems, they define topological Fukaya categories $ \Tw\mathcal{F}_A (S) $. If one restricts to the case of marked surfaces where each $ S^1 $ boundary component is fully marked, this is just $ \Tw\Gtl Q $ in our terminology. This is depicted in \autoref{fig:literature-hkk-marking}.

By an explicit analysis of all possible twisted complexes, they classify the objects of $ \Tw\Fuk (S, M) $ up to quasi-isomorphism. This yields two different classes, the string and band objects. We recall this classification in \autoref{sec:uncurving-stringsbands}. This classification led Bocklandt \cite{Bocklandt-book} to write down the explicit correspondence between curves and twisted complexes.

The paper \cite{HKK} then continues to classify stability conditions on a subcategory of the topological Fukaya category. The result is that these can be identified with singular flat structures on the marked surface with given poles or zeros.

In our paper, we depart from a special case of the topological Fukaya categories of \cite{HKK}. Indeed, a dimer model $ Q $ is a specific type of marked surface. Its topological Fukaya category in the sense of \cite{HKK} is simply $ \Tw\Gtl Q $.

The paper \cite{HKK} also helps us in \autoref{sec:splitting} to skip a few checks. Let us recapitulate the claims in that section: Given two zigzag paths $ L_1, L_2 ∈ \Tw\Gtl Q $, we would like to compute the definition of $ \H\Hom_{\Tw\Gtl Q} (L_1, L_2) $. While this could be checked by hand, we propose to exploit Bocklandt's equivalence \cite{Bocklandt}
\begin{equation*}
F: \H\Tw\Gtl Q \isoto \wFuk Q.
\end{equation*}
The zigzag paths $ L_1, L_2 $ live on the left-hand side, and hence
\begin{equation}
\label{eq:literature-hkk-dim}
\dim \Hom_{\H\Tw\Gtl Q} (L_1, L_2) ≅ \dim \H\Hom_{\wFuk Q} (F(L_1), F(L_2)).
\end{equation}
We are left with computing the right-hand side. For this, we need to know which curves the objects $ F(L_1) $ and $ F(L_2) $ are. Here \cite{HKK} comes into play and suggests that $ F(L_1) $ and $ F(L_2) $ are simply the smoothed-out versions of $ L_1 $ and $ L_2 $. With this assumption, the right-hand side of \eqref{eq:literature-hkk-dim} becomes simply the number of intersections between $ L_1 $ and $ L_2 $, plus two in case $ L_1 = L_2 $. This finishes the calculation of the hom space, but has cheated slightly in the identification of $ F(L_1) $ and $ F(L_2) $.

Far away on the horizon, Bridgeland \cite{Bridgeland-moduli} has suggested a conjecture regarding stability conditions versus deformations. The idea is as follows: If we reinterpret the flat structures of \cite{HKK} as deformations of the complex structure, they should constitute deformations of the derived category of coherent sheaves of the marked surface:
\begin{equation*}
\Stab\Fuk Q ≅ \Def\Coh Q.
\end{equation*}
Mirror symmetry of punctured surfaces ensures that under some conditions there is a dual dimer such that $ \Fuk Q ≅ \Coh \mirQ $ and $ \Coh Q ≅ \Fuk \mirQ $. Here $ \Coh $ is abuse of notation and mean a noncommutative version of coherent sheaves, e.g.~matrix factorizations. We then arrive at
\begin{gather}
\Stab\Coh \mirQ ≅ \Def\Fuk \mirQ, \\
\Def\Fuk Q ≅ \Def\Coh \mirQ, \label{eq:literature-hkk-def} \\
\Stab\Fuk Q ≅ \Stab\Coh \mirQ.
\end{gather}
Simply speaking, the conjecture arising from \cite{HKK} is that mirror symmetry swaps stability conditions and deformations. To prove this monster conjecture, we need a solid understanding of deformations of Fukaya categories. Our series of three papers will set up a deformed version of mirror symmetry, providing an explicit realization of the correspondence \eqref{eq:literature-hkk-def}. The present paper provides the preliminary step of equating deformations of the discrete model with those of the geometric side.

\subsubsection{Lowen-van den Bergh}
\label{sec:literature-LvdB}
In \cite{Lowen-vdB}, Lowen and Van den Bergh explain how to remove curvature from $ A_∞ $-deformations of dg categories, a cornerstone for \autoref{sec:deformed-uncurving}.

Lowen and Van den Bergh depart from a dg algebra $ A $ together with an infinitesimally curved $ A_∞ $-deformation $ A_q $ over $ ℂ⟦q⟧ $. Lowen and Van den Bergh observe that a category $ \Tw(A_q) $ of twisted complexes over $ A_q $ can be formed even with infinitesimal entries below the diagonal, just as in our definition of $ \Tw'\cat C_q $. Interpret $ A_q $ as an $ A_∞ $-deformation with a single object. Then the core observation of Lowen and Van den Bergh is that the following twisted complex has vanishing curvature:
\begin{equation}
\label{eq:literature-LvdB-trick}
X ≔ \left(A ⊕ A[1], \begin{pmatrix} 0 & μ^0_q / q \\ q \id_A & 0 \end{pmatrix}\right) ∈ \Tw(A_q).
\end{equation}
This means that $ B_q ≔ \End_{\Tw(A_q)} (X, X) $ is a curvature-free deformed $ A_∞ $-algebra. What is its special fiber $ B $? The higher products $ μ^{≥3} $ on $ B_q $ are given by embracing $ μ_{A_q} $ with the matrix entries $ μ^0_q /q $ and $ q \id_A $:
\begin{equation*}
μ^{k≥3}_{B_q} (a_k, …, a_1) = \sum μ^{≥3}_{\Add A_q} (δ, …, δ, a_k, …, δ, …, δ, a_1, δ, …, δ).
\end{equation*}
Restricting this sum to $ q = 0 $ yields only higher products $ μ^{≥3} $ of $ A $. Since $ A $ is a dg algebra, we deduce $ μ^{k≥3}_B = 0 $ and therefore $ B $ is a dg algebra as well. We conclude: $ B_q $ is a curvature-free $ A_∞ $-deformation $ B_q $ of some algebra $ B $.

Lowen and van den Bergh prove that $ A $ and $ B $ are in fact related by Morita equivalence. This costs substantial effort and uses the assumption that the curvature $ μ^0_{A_q} $ is nilpotent in the cohomology of $ A $. The result is however that $ A $ and $ B $ are Morita equivalent, and moreover that $ B_q $ is the deformation of $ B $ corresponding to the deformation $ A_q $ of $ A $ along this Morita equivalence:
\begin{center}
\begin{tikzpicture}
\path (0, 0) node {dg algebra A} (3, 0) node {\Large $ \rightsquigarrow $} (7, 0) node {dg algebra $ B $};
\path (0, -0.5) node {curved $ A_∞ $-deformation $ A_q $} (3, -0.5) node {\Large $ \rightsquigarrow $} (7, -0.5) node {uncurved $ A_∞ $-deformation $ B_q $};
\end{tikzpicture}
\end{center}

The work of Lowen and Van den Bergh helped us understand that curvature is essential in the notion of $ A_∞ $-deformations, but not an invariant on its own. While Lowen and Van den Bergh exchange the dg algebra itself to remove curvature, our \autoref{sec:uncurving} provides an example where a mere gauge transformation suffices to remove curvature.

In our recollection of $ A_∞ $-deformation theory, Lowen and Van den Bergh have greatly helped us understand how deformations can be transferred from one category to another. We have built on their understanding that the transfer should happen by means of a $ L_∞ $-quasi-equivalence between Hochschild DGLAs
while the corresponding map on Maurer-Cartan elements is always secondary:
\begin{center}
\begin{tikzpicture}
\path (0, 0) node (C) {$ \cat C $} (2.5, 0) node (HCC) {$ \HC(\cat C) $} (6, 0) node (MCC) {$ \MCb(\HC(\cat C), B) $};
\path (0, -1) node (D) {$ \cat D $} (2.5, -1) node (HCD) {$ \HC(\cat D) $} (6, -1) node (MCD) {$ \MCb(\HC(\cat D), B) $};
\path[draw, ->] (C.south) -- (D.north) node[pos=0.4, below, sloped] {$ \sim $};
\path[draw, ->] (HCC.south) -- (HCD.north) node[pos=0.4, below, sloped] {$ \sim $};
\path[draw, ->] (MCC.south) -- (MCD.north) node[pos=0.4, below, sloped] {$ \sim $};
\path (1.25, -0.5) node {\Large $ \rightsquigarrow $} (4, -0.5) node {\Large $ \rightsquigarrow $};
\end{tikzpicture}
\end{center}

As aftermath of our paper, we conclude that there is theoretically no hindrance to forming a derived category of a curved $ A_∞ $-deformation: By the definition of twisted completion and minimal models in \papertwoA, a derived category $ \H\Tw\cat C_q $ exists even for infinitesimally curved deformations. The statement of Lowen and Van den Bergh that a curved $ A_∞ $-deformation has no classical derived category remains true, but our paper contends that the study of deformations profits greatly from permitting also these “non-classical derived categories” $ \H\Tw\cat C_q $.

Our method in \autoref{sec:uncurving-gtl} seems to be both a variant and alternative to Lowen and Van den Berghs uncurving construction \eqref{eq:literature-LvdB-trick}. It is a closely related variant in that our uncurving procedure factorizes the curvature of $ (X, δ) $ into components of $ δ $ and new infinitesimal entries. By comparison, Lowen and Van den Bergh simply factorize $ μ^0 = (μ^0 / q) (q \id_A) $. Our procedure is also an alternative in that we uncurve the twisted complex itself, without passing to a different category $ \cat D $. This way there is no doubt that we have only performed a gauge equivalence, and checks for Morita equivalence are not required. Our method relies a lot on the fact that the twisted differential $ δ $ is very rich, and it would be interesting to know which other twisted complex categories have such property.

\subsubsection{Lekili-Perutz-Polishchuk}
\label{sec:comparison-LP}
In \cite{Lekili-Perutz}, Lekili and Perutz find a commutative mirror for the relative Fukaya category of the $ 1 $-punctured torus, apparently the first use of a relative Fukaya category in mirror symmetry. In \cite{Lekili-Polishchuk}, Lekili and Polishchuk generalize this result to the case of the $ n $-punctured torus. They depart from a finite collection of split-generators of the Fukaya category and compute part of their deformed products in the relative Fukaya category. This way, they contribute to \autoref{sec:subdisk}, with a viewpoint from the deformed geometric side.

Let $ \mathbb{T}_1 $ denote the $ 1 $-punctured torus. Lekili and Perutz depart from an explicit definition of the relative Fukaya pre-category $ \relpreFuk T_1 $: Working over the local ring $ ℤ⟦q⟧ $, every immersed disk is weighted by the number it covers the single puncture.

Lekili and Polishchuk regard the $ n $-punctured torus $ \mathbb{T}_n $. One might expect that they use an explicit model of the relative Fukaya category $ \relFuk \mathbb{T}_n $ and then prove it equivalent to their commutative mirror. Instead, they pick a set of $ n+1 $ curves $ L_0, …, L_n $ in $ \mathbb{T}_n $ which split-generate the wrapped Fukaya category $ \wFuk \mathbb{T}_n $.

\begin{figure}
\centering
\begin{tikzpicture}
\path[draw, lightgray] (0, 0) rectangle (3, 3);
\foreach \i in {0, 1, 2} \path[fill] (\i, 1) ++(0.2, 0) circle[radius=0.05];
\foreach \i in {1, 2, 3} \path[draw, semithick] (\i, 0) ++(-0.3, 0) -- ++(up:3) node[pos=0.6, left] {$ L_{\i} $};
\path[draw] (0, 0.5) -- (3, 0.5) node[at start, left] {$ L_0 $};
\end{tikzpicture}
\caption{The curves of Lekili and Polishchuk in case of $ n = 3 $ punctures}
\end{figure}
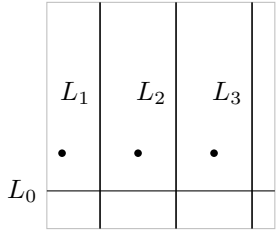

They do not attempt to compute the higher $ A_∞ $-products on this set of generators entirely, but rather show that the $ A_∞ $-structure must come from the perfect complexes of some complex curve $ T_n $ \cite[Theorem 1.1.1]{Lekili-Polishchuk}. The rest of their argument is devoted to guessing which curve $ T_n $ is the right one.

This deduction up to isomorphism yields a functor $ F: \{L_0, …, L_n\} → \Perf(T_n) $. To extend this functor to all of the relative Fukaya category, Lekili and Polishchuk view all objects of the Fukaya category as modules over these curves. More precisely, they regard a fully faithful Yoneda functor $ \wFuk(T_n) → \Mod(\{L_0, …, L_n\}) $. The right-hand side again maps to $ \Perf(T_n) $ by an extension of $ F $ to modules.

In the present paper, we have a very similar desire: to equate the $ A_∞ $-structure on $ \cat C_q ≔ \H\DefZigzagCat $ and the zigzag subcategory $ \cat D_q $ of the relative Fukaya category. If we tried to follow Lekili and Polishchuk's approach, we would start from the observation that the non-deformed versions $ \cat C $ and $ \cat D $ are isomorphic by \cite{Bocklandt}. We would then compute a few deformed higher products of $ \cat C_q $ and compare those with $ \cat D_q $, just enough to prove that $ \cat C_q ≅ \cat D_q $.

In the present paper, we do not follow the approach of Lekili and Polishchuk. In a sense, is a pity we were not able to guess the right structure like they did.

As aftermath of our paper, we recover the meaning of the curves $ L_0, …, L_n $ of Lekili and Polishchuk. Indeed, let $ Q $ denote the standard $ n $-punctured torus of \autoref{fig:prelim-std-dimers-torus}. Then the zigzag paths of $ Q $ are depicted in \autoref{fig:literature-polishchuk-zigzag}. There are precisely $ n $ diagonal, $ n $ vertical and $ 1 $ horizontal zigzag paths. Out of these, the vertical and horizontal are precisely the collection of Lekili and Polishchuk.

\begin{figure}
\centering
\begin{tikzpicture}[scale=0.8]
\begin{scope}[shift={(0, 0)}] 
\path[draw, lightgray] (0, 0) rectangle (4, 1);
\foreach \i in {1, 2, 3} \path[draw, lightgray] (\i, 0) -- (\i, 1);
\foreach \i in {0, 1, 2, 3} \path[draw, lightgray] (\i, 0) -- ++(1, 1);
\foreach \i in {0, 1, 2, 3} \path[draw, semithick, rounded corners] (\i, 0.1) ++(-0.1, 0) -- ++(up:1) -- ++(right:1) -- ++(up:1) -- ++(right:1) -- ++(up:1) -- ++(right:1);
\end{scope}
\begin{scope}[shift={(10.5, 0)}]  
\path[draw, lightgray] (0, 0) rectangle (4, 1);
\foreach \i in {1, 2, 3} \path[draw, lightgray] (\i, 0) -- (\i, 1);
\foreach \i in {0, 1, 2, 3} \path[draw, lightgray] (\i, 0) -- ++(1, 1);
\foreach \i in {0, 1, 2, 3} \path[draw, semithick, rounded corners] (\i, 0.1) -- ++(1, 1) -- ++(left:1) -- ++(1, 1) -- ++(left:1) -- ++(1, 1);
\end{scope}
\begin{scope}[shift={(5.5, 0)}] 
\path[draw, lightgray] (0, 0) rectangle (4, 1);
\foreach \i in {1, 2, 3} \path[draw, lightgray] (\i, 0) -- (\i, 1);
\foreach \i in {0, 1, 2, 3} \path[draw, lightgray] (\i, 0) -- ++(1, 1);
\path[draw, semithick, rounded corners] (4.1, 1.1) -- ++(-1, -1) -- ++(up:1) -- ++(-1, -1) -- ++(up:1) -- ++(-1, -1) -- ++(up:1) -- ++(-1, -1) -- ++(up:1) -- ++(-1, -1) -- ++(up:1);
\end{scope}
\end{tikzpicture}
\caption{Zigzag paths of the standard $ 4 $-punctured torus}
\label{fig:literature-polishchuk-zigzag}
\end{figure}

Our paper completes Lekili-Polishchuk's understanding of the deformed $ A_∞ $-structure on $ \{L_0, …, L_n\} $. Indeed, we compute an entire minimal model category $ \cat C_q = \H\DefZigzagCat $, which has the same deformed $ A_∞ $-structure on transversal sequences as the relative Fukaya category. It is technically not legitimate, but we could assume that $ \cat C_q $ indeed is a model for the relative Fukaya category. This would mean that we have computed all missing $ A_∞ $-structure that Lekili and Polishchuk were looking for.

\subsubsection{The lectures of Bocklandt}
\label{sec:literature-lectures}
A recent textbook \cite{Bocklandt-book} of Bocklandt explains gentle algebras in detail, shows how to stitch arcs together to form bands, and how to move towards the Fukaya category. Bocklandt's book contributes heavily to \autoref{sec:splitting}, departing from the discrete perspective without deformation.

In its Section 9, Bocklandt recollects the definition of the gentle algebra $ \Gtl Q $. Next, he shows how to stitch arcs together along shared angles. This procedure results in twisted complexes in $ \Tw\Gtl Q $.

For us, Bocklandt's explicit stitching procedure makes it entirely transparent how zigzag paths should be realized as twisted complexes. A zigzag path does not have a unique twisted complex representation, but there is a particularly simple one which makes direct use of the path's zigzag nature. This point of view is facilitated heavily by Bocklandt's section 9.2.

With this in mind, we can state that the twisted complex construction for gentle algebras is not the only one where the result can be identified geometrically. In fact, also twisted complexes of curves in the Fukaya category or wrapped Fukaya category can be identified as being quasi-isomorphic to curves that result from gluing together the arcs involved, see the book's Section 6.4.1.

Bocklandt's textbook contains several more hints relevant to the present paper, namely how to recognize similarity of $ \H\Tw\Gtl Q $ with the Fukaya category: In its section 9.2, the hom spaces in the minimal model of $ \Tw\Gtl Q $ are computed. Bocklandt delivers a basis of representatives of the cohomology $ \H\Hom(X, Y) $, in case $ X, Y $ are twisted complexes model transversal curves, and in case $ X = Y $ as well. He combines these ingredients into a description of some higher products of $ \H\Tw\Gtl Q $.

These calculations of Bocklandt provide a direct starting shot for the present paper: They tell us how to choose cohomology representatives for $ \H\Hom(X, Y) $ and indicate how to obtain the higher products. There are also vague indications as to how to build a homological splitting.

Our paper essentially completes the calculations of Bocklandt: First complete the cohomology basis elements of Bocklandt to an entire homological splitting, at least in the case of zigzag paths. Second, we compute the entire $ A_∞ $-structure on $ \H\DefZigzagCat $, including on non-transversal sequences, where Bocklandt's calculations are lacking. Third, we extend $ \Gtl Q $ to the deformed case and show how to obtain the relative Fukaya category. Our deformed case demonstrates how also complicated Kadeishvili trees can contribute to the higher products, in contrast to Bocklandt's non-deformed case where only the simplest Kadeishvili trees yield nonzero results. This renders our paper a powerful extrapolation of Bocklandt's method.

\subsection{Why should it work?}
\label{sec:whyshould}
This paper shows that the relative Fukaya category can be obtained from a small, discrete model. But why should such a small model exist? The question is why one expects the candidate we give indeed to be equivalent to the relative Fukaya category. In this section, we explain how one is led to believe from an a priori perspective that it should work, and explain why $ \Gtl_q Q $ is suited as a candidate.

\subsubsection{The model question}
\label{sec:whyshould-model}
In order to prove results concerning an $ A_∞ $-category $ \cat C $, one tries to switch between different models of $ \cat C $. This means, one is interested in $ A_∞ $-categories $ \cat D $ that are isomorphic, quasi-isomorphic, quasi-equivalent or derived equivalent to the $ \cat C $. If such a category $ \cat D $ satisfies certain geometric or algebraic properties or size constraints, it is called a \emph{model} (of the given kind) for $ \cat C $:

\begin{center}
\begin{tikzpicture}
\path (0, 0) node (C) {$ \cat C $} (1, 0) node {$ \cong $} (2, 0) node (D) {$ \cat D $};
\path[draw, ->, bend right=40] ($ (D.north) + (0.5, 0.3) $) coordinate (B) to ($ (D.north) + (0, 0.1) $);
\path (C.south) node[below] {original};
\path (D.south) node[below] {model};
\path (B) node[right] {better behaved};
\end{tikzpicture}
\end{center}

\noindent A standard question in symplectic geometry is then: Can we find a small model $ \cat D $ for the Fukaya category $ \cat C = \Fuk Q $? Ideally, this category $ \cat D $ would have very few objects, and still generate the whole Fukaya category. It does not work however, because cones over a small set of band objects do not yield all other bands. The question arises how to relax the task so that a small model can still be achieved. A very natural alternative is to require only that $ \cat C $ is \emph{contained} in the model $ \cat D $. Actually, one would not require $ \cat C ⊂ \cat D $, because $ \cat D $ itself is supposed to have few objects, but one would aim at:

\begin{center}
\begin{tikzpicture}
\path (0, 0) node (C) {$ \cat C $} (1, 0) node {$ ⊂ $} (2, 0) node (D) {$ \H\Tw\cat D $};
\path[draw, ->, bend right=40] ($ (D.north) + (0.5, 0.3) $) coordinate (B) to ($ (D.north) + (0, 0.1) $);
\path (C.south) node[below] {original};
\path (D.south) node[below] {model};
\path (B) node[right] {better behaved};
\end{tikzpicture}
\end{center}

\noindent Thanks to Bocklandt and Abouzaid \cite{Bocklandt}, it is now known that the Fukaya category is indeed contained in the derived category $ \H\Tw\Gtl Q $ of the gentle algebra $ \Gtl Q $. There is a quasi-fully-faithful inclusion
\begin{equation*}
\Fuk Q ⊂ \H\Tw\Gtl Q.
\end{equation*}
In fact, the category $ \H\Tw\Gtl Q $ is not all too large: It is quasi-equivalent to the wrapped Fukaya category $ \wFuk Q $. In other words, Bocklandt and Abouzaid resolve the (relaxed) model question for $ \Fuk Q $ positively.

Let us now pose the same model question for deformed $ A_∞ $-categories: Given a category $ \cat C $ with a deformation $ \cat C_q $, can we find a better behaved category $ \cat D $ with a deformation $ \cat D_q $ such that $ \cat C_q $ and $ \cat D_q $ are isomorphic, quasi-isomorphic, quasi-equivalent or derived equivalent?

\begin{center}
\begin{tikzpicture}
\path (0, 0) node (C) {$ \cat C_q $} (1, 0) node {$ ⊂ $} (2, 0) node (D) {$ \cat D_q $};
\path[draw, ->, bend right=40] ($ (D.north) + (0.5, 0.3) $) coordinate (B) to ($ (D.north) + (0, 0.1) $);
\path (C.south) node[below] {original};
\path (D.south) node[below] {model};
\path (B) node[right] {better behaved};
\end{tikzpicture}
\end{center}

\noindent Let us discuss what this means. In the above sketch, we have used $ \cong $ to indicate one of the four notions of equivalence. In either case, an equivalence on the level of deformations necessarily requires an equivalence on the non-deformed level. Conversely, one can transport deformations along equivalence of (non-deformed) categories. Let us summarize as follows:
\begin{itemize}
\item Let $ \cat D $ be a model for $ \cat C $, and let $ \cat C_q $ be a deformation of $ \cat C $. Then there exists a model $ \cat D_q $ for $ \cat C_q $, obtained as a deformation of $ \cat D $.
\item If $ \cat C $ has no good model (of a certain kind), then $ \cat C_q $ does not have a good model either.
\end{itemize}
Relative Fukaya categories were introduced by Seidel as a deformation of ordinary Fukaya categories. We may now ask: Is it possible to provide a small model for the relative Fukaya category? Unfortunately, this is not possible either. A small model for $ \relFuk Q $ would also include a small model for $ \Fuk Q $ itself, which does not exist. The right approach becomes apparent by relaxing the task again. Let us first spell this out in general:

Let $ \cat C $ be a category with a relaxed model $ \cat D $, and let $ \cat C_q $ be a deformation. Does a deformation $ \cat D_q $ exist such that it is a relaxed model for $ \cat C_q $? The answer is that this does not necessarily exist. The reason is that deformations cannot necessarily be lifted from $ \cat C $ to $ \cat D $. In fact, the restriction map $ \HC(\cat D) → \HC(\cat C) $ induced by the inclusion $ \cat C → \H\Tw\cat D $ does absolutely not have to be a quasi-isomorphism of $ L_∞ $-algebras. An easy example is the inclusion of quivers

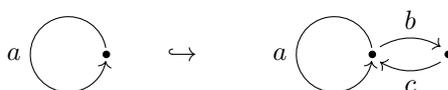
\begin{center}
\begin{tikzpicture}
\path[draw, ->] (-0.5, 0) ++(10:0.5) arc(10:350:0.5) node[midway, left] {$ a $};
\path[fill] (0, 0) circle[radius=0.05];
\path (1, 0) node {$ \embeds $};
\begin{scope}[shift={(3.5, 0)}]
\path[draw, ->] (-0.5, 0) ++(10:0.5)  arc(10:350:0.5) node[midway, left] {$ a $};
\path[fill] (0, 0) circle[radius=0.05];
\path[draw, ->] (0.1, 0.1) to[out=30, in=150] node[midway, above] {$ b $} (0.9, 0.1);
\path[draw, ->] (0.9, -0.1) to[out=210, in=330] node[midway, below] {$ c $} (0.1, -0.1);
\path[fill] (1, 0) circle[radius=0.05];
\end{scope}
\end{tikzpicture}
\end{center}

The center of the quiver algebra $ ℂQ_1 $ on the left is of course $ ℂ[a] $, while the center of the algebra $ ℂQ_2 $ on the right is just $ ℂ^2 $, spanned by the two idempotents. We deduce that the map $ \HH^0 (ℂQ_2) → \HH^0 (ℂQ_1) $ is not surjective. Here $ \HH^0 $ denotes classical Hochschild cohomology, which is the same as $ \HH^{-1} $ in the $ A_∞ $-grading. In short, restriction maps between Hochschild cohomologies are far from surjective due to “global” phenomena.

Let us tie this back to the question of finding a small model for the relative Fukaya category. We have already discussed that $ \Gtl Q $ provides a small relaxed model for the Fukaya category. As we have just seen, this does however not imply the existence of a relaxed model for $ \relFuk Q $ in the form of a deformation of $ \Gtl Q $. One starting point for understanding the present paper is therefore:

\begin{center}
To find an $ A_∞ $-category $ \cat D $ together with a deformation 
$ \cat D_q $ \par such that $ \relFuk Q $ embeds quasi-fully-faithfully into $ \H\Tw\cat D_q $.
\end{center}

A priori it is not clear that such a category $ \cat D $ and deformation $ \cat D_q $ should exist. The reason is that the ordinary and wrapped Fukaya categories are not equivalent and have different deformation theory. For the same reason, such a pair is not uniquely determined. There are however several ways of trying to find such a pair:
\begin{itemize}
\item[A1] Guessing,
\item[A2] Trying out the candidate deformation $ \cat D_q ≔ \Gtl_q Q $ of $ \cat D ≔ \Gtl Q $.
\item[A3] Extending the relative Fukaya category to a deformation $ \wFuk_q $ of the wrapped Fukaya category.
\end{itemize}

In this paper, we succeed in approach A2: We show that $ \Gtl_q Q $ is a relaxed small model for $ \relFuk Q $, at least on the subcategory of zigzag paths. In \autoref{sec:whyshould-candidate} we explain why approach A2 is plausible and in \autoref{sec:whyshould-alternative} we explain why approach A3 is promising for mathematicians who can handle wrapped symplectic geometry.

There are three reasons why the author picked approach A2 instead of A3. First, we already have a concrete $ \Gtl_q Q $ available from \paperone. Second, approach A2 comes only with combinatorial calculations, as opposed to deforming and working with the wrapped Fukaya category in approach A3. The expertise in symplectic geometry on the side of the author was simply not enough. Third, this paper was originally written not in order to find a small model, but to compute the subcategory of zigzag paths in $ \H\Tw\Gtl_q Q $. The interpretation as a small model for $ \relFuk Q $ has come out as a useful byproduct.

\subsubsection{The candidate $ \Gtl_q Q $}
\label{sec:whyshould-candidate}
The goal of this section is to describe why our candidate $ \Gtl_q Q $ is plausible as a (relaxed) model for $ \relFuk Q $. The category $ \Gtl_q Q $ itself is a seemingly arbitrary choice defined in \paperone, so that it is a priori not clear why it should be a model for $ \relFuk Q $. There are however reasons why one should expect $ \Gtl_q Q $ to be a model, even before performing any calculations. In this section, we explain those reasons.

Bocklandt-Abouzaid showed that the gentle algebra $ \Gtl Q $ is equivalent to the wrapped Fukaya category. In the words of \autoref{sec:whyshould-model}, this implies that $ \Gtl Q $ is a relaxed model for $ \Fuk Q $. As we have seen in \autoref{sec:whyshould-model}, it is however far from clear that a relaxed model for a deformation can be obtained as a deformation of a relaxed model. In other words, if $ \Gtl Q $ is a model for $ \Fuk Q $, why should the deformation $ \Gtl_q Q $ be a model for $ \relFuk Q $?

There are three reasons why one might expect $ \Gtl_q Q $ to be a relaxed model for $ \relFuk Q $:
\begin{itemize}
\item The derived category $ \H\Tw\Gtl_q Q $ exists by construction, and it is a deformation of $ \H\Tw\Gtl Q $. In particular, it is equivalent to a deformation of $ \wFuk Q $ and has a restriction to $ \Fuk Q $. In other words, it contains some deformation of $ \Fuk Q $. One may now speculate which deformation of $ \Fuk Q $ it concerns.
\item A glance at the deformed higher products shows that $ \Gtl_q Q $ closely resembles $ \relFuk Q $: Although the objects of both categories are completely disjoint, every disk containing one puncture gets multiplied by that puncture. Every disk containing two punctures gets multiplied by both, etc. One easily becomes suspicious that the deformation of $ \Fuk Q $ contained in $ \H\Tw\Gtl_q Q $ is actually $ \relFuk Q $.
\item Reasoning with the beauty of mathematics, one should expect that $ \relFuk Q $ is such a reasonable deformation that is extends to $ \wFuk Q $. By the Bocklandt-Abouzaid equivalence, it then induces a deformation on $ \Gtl Q $, and one may now guess which one this is: probably isomorphic to $ \Gtl_q Q $.
\end{itemize}

Against the second reason, one might object that similarity of deformations is not the same as equality. It might be possible that the deformation of $ \Fuk Q $ contained in $ \H\Tw\Gtl_q Q $ is slightly off, even though the products of $ \relFuk Q $ and $ \Gtl_q Q $ look so similar. For example, $ \Gtl_q Q $ intrinsically multiplies disks by $ q $. The process of deriving $ \Gtl_q Q $ may change this factor however to $ q + q^2 $ instead. This would imply that $ \H\Tw\Gtl_q Q $ does not have the same higher products as $ \relFuk $.

The fact that this $ \Gtl_q Q $ actually is a relaxed model for $ \relFuk Q $ and the higher products on $ \H\Tw\Gtl_q Q $ are identical to those of $ \relFuk Q $ is therefore out of pure luck. We will comment on this fact in \autoref{sec:whydoes-luck}.

\subsubsection{Alternative via the wrapped Fukaya category}
\label{sec:whyshould-alternative}
In this section, we explain another approach to obtain a small (relaxed) model for $ \relFuk Q $. Namely, we comment on the idea to deform the wrapped Fukaya category, labeled A3 in \autoref{sec:whyshould-model}. We will see why it is realistic, and what the difficulties are.

Let us recall approach A3 as follows: One tries to lift the deformation of $ \Fuk Q $ given by $ \relFuk Q $ to a deformation $ \wFuk_q Q $ of $ \wFuk Q $. If one succeeds at this approach, then one immediately has $ \relFuk Q $ as a subcategory of $ \wFuk_q Q $. Pick a generating set $ X ⊂ \wFuk Q $, typically a collection of arcs that split the surface. Denote by $ X_q ⊂ \wFuk_q Q $ the restriction of the deformation $ \wFuk Q $ to the generating set $ X $. Since $ X $ is a generating set for $ \wFuk Q $, we have a quasi-equivalence
\begin{equation*}
\H\Tw X → \wFuk Q
\end{equation*}
induced from the inclusion $ X ⊂ \wFuk Q $. The deformation $ \wFuk_q Q $ is therefore already determined by the deformation $ X_q $. In other words, we have a quasi-equivalence
\begin{equation*}
\H\Tw X_q → \wFuk_q Q.
\end{equation*}
Since the right-hand side $ \wFuk_q Q $ contains the relative Fukaya category, we conclude that

\begin{center}
\begin{tikzpicture}
\path (0, 0) node (C) {$ \H\relFuk Q $} (1, 0) node {$ ⊂ $} (2, 0) node (D) {$ \H\Tw X_q $};
\path[draw, ->, bend right=40] ($ (D.north east) + (0.5, 0.3) $) coordinate (B) to ($ (D.north) + (0, 0.1) $);
\path (C.south) node[below] {original};
\path (D.south) node[below] {model};
\path (B) node[right] {small};
\end{tikzpicture}
\end{center}

\noindent In other words, $ X_q $ is a small model for $ \relFuk $. We conclude: A lift of the deformation $ \relFuk Q $ to $ \wFuk Q $ solves the (relaxed) model question for $ \relFuk Q $. Such a lift does not need to exist a priori and it is not unique.

Let us explain how one may obtain a candidate deformation $ \wFuk_q Q $ such that $ \relFuk Q ⊂ \wFuk_q Q $. We can already guess several of its properties:
\begin{itemize}
\item On band objects, the higher products are just given by disks multiplied by $ q $-parameters, as in $ \relFuk Q $.
\item String objects need to have curvature. There is no technical necessity for this, but it is likely from the point of view that our combinatorial model $ \Gtl_q Q $ also has curvature.
\item The definition of higher products through Hamiltonian deformations needs to be completely revised to be compatible with the curvature. Due to the new higher products, there now exist infinitesimal results of disks even on teardrops. The obstruction theory in the definition of the Fukaya category needs to be completely revised.
\end{itemize}
This list already highlights some of the difficulties. The author has no clue how to properly define such an extension.

Let us assume for a moment that the paper's result extends beyond zigzag paths. From this a posteriori perspective we can deduce that a lift from $ \relFuk Q $ to a deformation $ \wFuk_q Q $ exists: Regard the Bocklandt-Abouzaid quasi-equivalence
\begin{equation*}
\H\Tw\Gtl Q \isoto \wFuk Q.
\end{equation*}
Then the deformation $ \Gtl_q Q $ of $ \Gtl Q $ induces a deformation $ \H\Tw\Gtl_q Q $ of $ \H\Tw\Gtl Q $, and by transport through the quasi-equivalence also a deformation $ \wFuk_q Q $. Since $ \H\Tw\Gtl_q Q $ contains the relative Fukaya category, we deduce that the same holds for $ \wFuk_q Q $:
\begin{equation*}
\relFuk ⊂ \H\Tw\Gtl_q Q \isoto \wFuk_q Q.
\end{equation*}
In other words, if one believes for a moment that the result of this paper extends to all band objects, then a lift from $ \relFuk Q $ to $ \wFuk Q $ necessarily exists. Approach A3 does therefore have a solution, although it is unclear how to construct it explicitly.

\subsection{Why does it work?}
\label{sec:whydoes}
This paper shows that the relative Fukaya category can be obtained from a small, discrete model. But why does the calculation work out? What are the ingredients that make it work? In contrast to the a priori discussion in \autoref{sec:whyshould}, we explain in the present section why it works from an a posteriori perspective. In particular, we discuss the role of choices and luck.

Let us paraphrase the methods of this paper. The starting point is the deformed gentle algebra $ \Gtl_q Q $. The task is to prove that its derived category $ \H\Tw\Gtl_q Q $ contains the relative Fukaya category. To achieve this, we need to realize all Lagrangians in the Fukaya category as specific twisted complexes over $ \Gtl Q $, and show that the subcategory of these twisted complexes equals the relative Fukaya category up to quasi-equivalence of deformations.

How would we achieve an equivalence between this subcategory of $ \H\Tw\Gtl_q Q $ and the relative Fukaya category? The relative Fukaya category $ \relFuk Q $ has mostly vanishing differential $ μ^1 $, while the category $ \Tw\Gtl_q Q $ has large hom spaces and non-vanishing differential. They are clearly far away, but the category $ \H\Tw\Gtl_q Q $ already comes closer to the relative Fukaya category. In the present paper, we show how to actually match them. During the calculations, four facilitating factors have come into play:
\begin{itemize}
\item Zigzags: Instead of proving the whole relative Fukaya category to lie inside $ \H\Tw\Gtl_q Q $, we only prove this for the subcategory $ \DefZigzagCat $ of zigzag paths.
\item Choices: We choose a “natural” homological splitting of $ \DefZigzagCat $.
\item Luck: During the calculation of the minimal model structure of $ \DefZigzagCat $, our choice of homological splitting proves to be right one both for efficient calculation and to obtain exactly the relative Fukaya category.
\item Fearless calculations: Performing the model computation for $ \DefZigzagCat $ emits enormous amounts of data and requires us to construct a tower of data structures as depicted in \autoref{fig:intro-data}. Binding the discrete data structures together to form smooth disk requires us to work through hundreds of calculations and special cases in order to bring order into the chaos. Practically, lots of trees need to be classified and large multiplication tables need to be filled. This paper performs the calculation until the bitter end.
\end{itemize}
We explain these four facilitating factors in more detail in \autoref{sec:whydoes-zigzag}, \ref{sec:whydoes-choice} and \ref{sec:whydoes-luck}.

\subsubsection{Restriction to zigzag paths}
\label{sec:whydoes-zigzag}
The result presented in this paper is restricted to zigzag paths. In this section we explain how this restriction eases the calculations and how the general case may be obtained later on.

Recall that zigzag paths are paths in a dimer that alternatingly turn left and right. When we say “zigzag path”, we frequently refer to their realization as twisted complex in $ \Tw\Gtl Q $ or as a band object in $ \Fuk Q $. Zigzag paths are a small class out of a large set of objects in both categories. Three factors distinguish zigzag paths from other band objects in $ \Fuk Q $:
\begin{itemize}
\item The higher structure on zigzag paths is necessary to compute a mirror for $ \Gtl_q Q $, according to Cho-Hong-Lau.
\item The arcs in the twisted complex representation of zigzag paths have only small angles between each other, i.e.~no full turns or larger angles. This makes it easy to get grip on the disks between zigzag paths.
\item If one assumes that $ Q $ is geometrically consistent, a mild requirement, then all zigzag paths in $ Q $ bound neither discrete nor smooth immersed disks. This is very useful.
\end{itemize}

It appears possible that the restriction to zigzag paths be overcome in the future, even without redoing the calculations. Let us sketch how this will work. The first step is to prove mirror symmetry for $ \Gtl_q Q $, and the second step is to realize that the mirror depends only on the higher structure on zigzag paths.

Indeed, both $ \relFuk $ and $ \H\Tw\Gtl_q Q $ produce mirror functors
\begin{equation*}
\Mod\relFuk → \mf(A_q, ℓ_q) \quad \text{and} \quad \H\Tw\Gtl_q Q → \mf(A_q, ℓ_q).
\end{equation*}
Both mirrors $ \mf(A_q, ℓ_q) $ are equal, since the Cho-Hong-Lau construction only depends on the structure on the zigzag paths. The module category $ \Mod\relFuk $ contains quasi-fully-faithfully some deformed copy $ (\Gtl Q)_q' $ of $ \Gtl_q Q $ and so does $ \H\Tw\Gtl_q Q $ contain the deformation $ \Gtl_q Q $. Both are mapped quasi-equivalently to the mirror. It seems that we can deduce this way that $ (\Gtl Q)_q' \cong \Gtl_q Q $ as deformations of $ \Gtl Q $. Together with $ \relFuk ⊂ \H\Tw(\Gtl Q)_q' $, we should be able to deduce that $ \relFuk $ is simply a subcategory of $ \H\Tw\Gtl_q Q $. In other words, this should imply that $ \Gtl_q Q $ is a small model for $ \relFuk $.

\subsubsection{Choice}
\label{sec:whydoes-choice}
This paper presents a minimal model for (part of) $ \H\Tw\Gtl_q Q $. Such a minimal model is by no means unique. In this section, we explain why our specific choice of homological splitting works so well.

Let $ \cat C $ be an $ A_∞ $-category. Recall that by a minimal model for $ \cat C $ one means any other $ A_∞ $-category $ \cat D $ such that $ \cat D $ is minimal and $ \cat C $ and $ \cat D $ are quasi-isomorphic:
\begin{equation*}
μ^1_{\cat D} = 0 \text{ and } \cat C \cong \cat D.
\end{equation*}
Given a category $ \cat C $, one may look for minimal models simply by guessing. Such a guess involves
\begin{itemize}
\item Possibly identifying the cohomology $ \H\Hom(X, Y) $ for every $ X, Y ∈ \cat C $ with some explicit graded vector space $ \cat D(X, Y) $.
\item Guessing an $ A_∞ $-structure on these spaces $ \cat D(X, Y) $, turning them into an $ A_∞ $-category $ \cat D $.
\item Finding an $ A_∞ $-quasi-isomorphism $ \cat C → \cat D $ or $ \cat D → \cat C $.
\end{itemize}
Guessing minimal models requires an enormous imagination.

There are also systemic ways of finding minimal models. In fact, the Kadeishvili theorem grants the existence of minimal models and provides an explicit way to construct them. The formula for the minimal model depends on the choice of a so-called homological splitting $ \Hom_{\cat C} = H ⊕ I ⊕ R $.


Assume we have chosen a homological splitting $ \Hom_{\cat C} = H ⊕ I ⊕ R $. Then the map $ μ^1: R → I $ is bijective. One then defines the so-called codifferential $ h: I → R $ as the inverse of $ μ^1: R → I $. The Kadeishvili construction then describes the $ \H\cat C $ as follows: The objects are the same as in $ \cat C $. The hom spaces are the chosen cohomology representatives $ H $. The differential is defined as $ μ^1_{\H\cat C} ≔ 0 $. The interesting part in the definition are the (higher) products. They are defined as sums over trees of the form

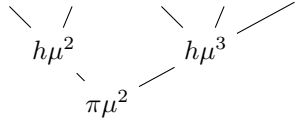
\begin{center}
\begin{tikzpicture}
\path node (A) {} node[right of=A] (B) {} node[right of=B] (C) {} node[right of=C] (D) {} node[right of=D] (E) {}
node[below right of=A] (A1) {$ h μ^2 $} edge (A) edge (B)
node[below right of=C] (B1) {$ h μ^3 $} edge (C) edge (D) edge (E)
node[below right of=A1] {$ π μ^2 $} edge (A1) edge (B1);
\end{tikzpicture}
\end{center}

For two inputs, there is precisely 1 such tree. For three inputs, there are 3 such trees. For four inputs, there are 11 such trees. The result of each tree shall be multiplied by a sign. The sign is given by $ (-1)^{N_T} $, where $ N_T $ is the number of nodes in the tree, excluding the root. In other words, $ s $ is the number of nodes in the tree labeled $ h μ $. For instance, the product $ μ^2 (a, b) $ for $ a, b ∈ H $ is simply given by
\begin{equation*}
μ^2_{\H\cat C} (a, b) = π μ^2 (a, b).
\end{equation*}
The higher product $ μ^3 (a, b, c) $ for $ a, b, c ∈ H $ is given by
\begin{equation*}
μ^3_{\H\cat C} (a, b, c) = π μ^3 (a, b, c) - π μ^2 (h μ^2 (a, b), c) - π μ^2 (a, h μ^2 (b, c)).
\end{equation*}
Observing these formulas, we conclude that the minimal model does depend on the choice of $ H $ and $ R $. One may also say: The minimal model depends on the choice of codifferential.

In this paper, we select one concrete choice of a homological splitting for the category $ \ZigzagCat $ of zigzag paths in $ \Tw\Gtl Q $. The choice looks arbitrary, but has some sophistication behind it. Let us explain the philosophy behind the cohomology representatives $ H $ in our choice:
\begin{itemize}
\item We know how many representatives we have to choose: as many as $ \H\Hom(L_1, L_2) $ has dimension.
\item The dimensions of $ \H\Hom(L_1, L_2) $ and the dimension of the hom space in the Fukaya category are equal (either by calculation or by using Bocklandt-Abouzaid). Hom spaces in the Fukaya category are spanned by intersection points, therefore we should try to find one representatives of $ \H\Hom(L_1, L_2) $ for every intersection point.
\item For each intersection point $ p ∈ L_1 ∩ L_2 $, choose the representative in $ H $ such that we have the best chance of obtaining the Fukaya category as minimal model. For example, a disk existing in the Fukaya category should be realizable as a product $ μ^{≥3}_{\ZigzagCat} $ of the corresponding basis elements in $ H $.
\item The signs of the elements in $ H $ should be chosen such that in the minimal model we obtain exactly the Abouzaid sign rule, without further sign conversion.
\end{itemize}

Regard an endomorphism space $ \End(L, L) $ of a zigzag path $ L ∈ \ZigzagCat $. Our choice for $ H $ consists of two morphism of $ \End(L, L) $: the identity and a co-identity. While the identity element of $ \End(L, L) $ naturally stems from the unitality of $ \Gtl Q $, the choice of co-identity involves a choice. We namely define the co-identity to be any of the angles involved in the $ δ $-matrix of $ L $. In other words, we choose the connecting angle between an arbitrary pair of consecutive arcs in $ L $.

Why is this a sensible choice? One of the reasons to use the identity for $ H $ is that it is very natural and it provides a strict unit in the minimal model $ \H\ZigzagCat $. This strict unit is simultaneously necessary to exist if we want to make $ \H\ZigzagCat $ equal to the zigzag paths in the Fukaya category.

A reason why we choose the other basis element of $ H $ to be a small angle between two consecutive arcs of $ L $ is that this angle is easily seen to lie in the kernel of $ μ^1: \End(L, L) → \End(L, L) $. Moreover, we want to obtain the Fukaya category as minimal model, which means that we have to reflect the arbitrary location of the co-identity morphism of Fukaya categories an closely as possible by means of the combinatorical datum of an angle.

\subsubsection{Luck}
\label{sec:whydoes-luck}
A decent amount of luck has been involved in the functioning of the present paper. In this section, we present five specific occasions where luck is decisive. The reader instead interested in a technical explanation why our choice of $ \Gtl_q Q $ and the homological splitting are wise choices is referred to \autoref{sec:whyshould-candidate} and \autoref{sec:whydoes-choice}.

\paragraph*{Transparency of deformed cohomology basis}
After building the homological splitting in the non-deformed case, we prepare in \autoref{sec:deformed-tails} the calculation of the deformed differential $ μ^1_q $ on $ \DefZigzagCat $. It turns out that the differential $ μ^1_q $ of any morphisms falls apart in contributions of certain types E, F, G and H.

Cohomology basis elements come from type B and C situations which restrict the tail to type E disks. The entire tail of a cohomology basis element then becomes relatively simple: It depends only on type E disks, and its tail terms are all of the form $ β $ (A). The description of the deformed cohomology basis elements becomes not only explicit this way, but also very homogeneous.

\paragraph*{Requirements for deformed Kadeishvili theorem}
We are lucky that the deformed Kadeishvili theorem can be established in the full generality. From a technical point of view, the Kadeishvili theorem is the only part of the paper that is not straightforward. It form a bottleneck for the minimal model computation and without its working we could not have pursued the calculation.

\paragraph*{Transparency of the deformed codifferential}
Luck comes into play in our computation of the deformed codifferential $ h_q $ in \autoref{sec:deformed-codif}. As always in this paper, this computation is rather an enumeration in terms of disks than a calculation with a concrete output. The deformed codifferential that illustrates the impact of luck best is $ h_q (βα) $, where the angles $ α, β $ are from an A situation. In this case, we have to find a sum of angles in $ R $ whose differential totals to $ βα $ plus possibly terms of $ R $. The first-order guess is $ β $ itself, however $ μ^1_q (β) $ may also contain disk terms from E, F, G and H disks.

We are double lucky. First, the F disks only produce $ β $ angles from A situations, the G2 disks only produce $ α_3 $ and $ α_4 $ angles, and the H disks only produce $ β $ and $ β' $ angles from C situations. All of these angles lie in the kernel of $ h_q $. In other words, those angles are in fact irrelevant in order to compute $ h_q (βα) $. We conclude that only the type E and type G1 disks are relevant for computing $ h_q (βα) $, which greatly reduces complexity.

As for $ α_4 $, it can be written as a signed sum $ α_4 = ± h ± α_3 $ of the cohomology basis element $ h = (-1)^{\#α_3 + 1} α_3 + (-1)^{\#α_4} α_4 $ and the angle $ α_3 $ lying in $ R $. As a cohomology basis element, $ h $ in turn can be written as a sum $ h = h' + r $ of a deformed cohomology basis element $ h' $ and an remainder $ r ∈ R $. All of $ h' $, $ r $ and $ α_3 $ have vanishing codifferentials $ h_q $, so that we conclude $ h_q (α_4) = 0 $.

Second, the G1 disks yields result of the form $ α_1 ± α_2 $, where $ α_1 $ and $ α_2 $ are from a B situation. The angle $ α_1 $ again lies in $ R $, while $ α_2 $ equals $ d(\id_{2→5}) $ modulo kernel of $ h_q $. Since $ μ^1_q (\id_{2→5}) = d(\id_{2→5}) $, we can simply add $ \id_{2→5} ∈ R $ to $ βα $ and $ μ^1_q (β ± \id_{2→5}) $ will eliminate the $ α_2 $ term. Ultimately, every G1 disk only adds in a simple B situation identity into the $ h_q (βα) $. This is the reason we obtain the comparatively neat formula in \autoref{th:deformed-codif-ba}.

\paragraph*{The chaos and order of result components}
In \autoref{sec:resultcomp}, we introduce the notion of result components. The subsequent classification of result components, its matching with immersed disks and the analysis of the immersed disks obtained this way is a roller coaster ride of case distinctions. Despite the intermediate chaos, the result collapses into a manageable description: four types of immersed disks (CR, ID, DS), following more or less the same rules. This collapse is a very fortunate turn.

\paragraph*{Just the right products in $ \H\DefZigzagCat $}
Even with a slightly different homological splitting, we might already have obtained a minimal model $ \H\DefZigzagCat $ that looks entirely different from the relative Fukaya category. It would be hopeless to compare even a slightly different result to the relative Fukaya category. We are very fortunate that we obtain the higher products of the relative Fukaya category up to strict isomorphism.

\subsection{Which calculations can be reused?}
\label{sec:reuse}
The heart of the present paper is a long and very specific calculation. In fact, the starting point consists of a very concrete deformation $ \Gtl_q Q $ of the gentle algebra and the specific subcategory of $ \Tw\Gtl_q Q $ given by the zigzag paths. This raises the question how the calculations and the result presented here can ever be used by other mathematicians for their own purposes.

In this section, we would like to answer this question. We explain how one can use the gentle algebra, the specific deformation $ \Gtl_q Q $, the homological splitting and the notion of result components in a modular way as standard tools in computations.

We are convinced that while the precise calculations apply only to the specific situation of $ \Gtl_q Q $, the versatility lies in the manner of performing the calculations and matching their result with the expected outcome. We contend that the mathematical value of the present paper mainly lies in making Kadeishvili trees computable.

\subsubsection{The gentle algebra}
The use of the gentle algebra to perform calculations in mirror symmetry of punctured surfaces is not yet standard, as of writing. For example, in \cite{MS-sphere} the notion is still implicit. Some newer work \cite{HKK} however uses the notion actively. In this section, we would like to highlight how easy the gentle algebra makes it to describe intersection theory.

Second, the twisted complexes of $ \Gtl Q $ can be classified up to quasi-isomorphism. Recall from \autoref{sec:uncurving-stringsbands} that the twisted complexes of $ \Gtl Q $ can be classified as string and band objects. Formulated the other way around, every twisted complex of $ \Gtl Q $ can be obtained up to quasi-isomorphism by stitching together arcs along angles. Regard two twisted complexes $ X, Y ∈ \Tw\Gtl Q $ stitched together from arcs. Then the hom space $ \Hom(X, Y) $ is spanned by all angles from arcs of $ X $ to arcs of $ Y $.

Given a whole sequence $ X_1, …, X_{k+1} $ of twisted complexes in $ \Tw\Gtl Q $ and angles $ α_i: X_i → X_{i+1} $, how to evaluate the higher product $ μ^k (α_k, …, α_1) $? By definition, this product is taken in the $ A_∞ $-category $ \Tw\Gtl Q $ and as such is made up of $ δ $-insertions. For each $ X_i $, the possible $ δ $-insertions are insertions of arbitrary angles used to stitch together the arcs of $ X_i $. In total, this higher product $ μ^k (α_k, …, α_1) $ gives a result if the angle sequence $ α_1, …, α_k $ can be filled up with $ δ $-insertions to form an immersed disk.

We see that even the twisted completion $ \Tw\Gtl Q $ is an utterly geometric object and can be used for geometric proofs.

\subsubsection{The deformation $ \Gtl_q Q $}
In \paperone\ we introduced the deformed gentle algebra $ \Gtl_q Q $. In fact, we provided even broader deformations and proved that they exhaust all deformations of $ \Gtl Q $ up to gauge equivalence. In this section, we would like to explain what makes $ \Gtl_q Q $ so versatile for studying deformations of Fukaya categories and mirror symmetry.

First, $ \Gtl Q $ itself is a small category itself and such ideally suited for computations. Its deformation $ \Gtl_q Q $ can be described fairly easily. Already the crude insight that $ \Gtl_q Q $ is a deformation capturing behavior similar to the relative Fukaya category makes $ \Gtl_q Q $ an interesting A-side of mirror symmetry. For comments on the use of $ \Gtl_q Q $ as model for the relative Fukaya category, see \autoref{sec:whyshould-model}.

As we recall in \autoref{sec:uncurving-stringsbands}, the twisted complexes of $ \Gtl Q $ can be classified as string and band objects. As we show in \autoref{sec:uncurving-gtl}, most band objects can be uncurved. The uncurving procedure adds in infinitesimal connecting angles into the $ δ $-matrix. Let us explain the effect of this procedure. Regard a sequence of uncurved twisted complexes $ X_1, …, X_{k+1} $ and angles $ α_i: X_i → X_{i+1} $. Then the higher product $ μ^k (α_k, …, α_1) $, taken in $ \Tw\Gtl_q Q $ now includes $ δ $-insertions of the additional infinitesimal angles in the $ δ $-matrices of the $ X_i $. This makes that also immersed disks count that are bounded by whole segments of the curves $ X_i $, instead of only a single arc as is the case without deformation. In particular, immersed disks between the curves $ X_i $ that also cover an arbitrary number of punctures now contribute to the product.

In the present paper, we match $ \H\Tw\Gtl_q Q $ with $ \relFuk Q $. In other words, the deformation $ \Gtl_q Q $ of $ \Gtl Q $ induces a deformation $ \H\DefZigzagCat $ of $ \H\ZigzagCat $ that looks like the relative Fukaya category. It is interesting to speculate what happens if we start with other deformations on the $ \Gtl Q $ side. More precisely, recall from \paperone\ that $ \Gtl Q $ also permits deformations where “orbifold disks” contribute to the higher products. Such a deformation of course also induces a deformation on $ \H\ZigzagCat $. Since $ \H\ZigzagCat $ is a full subcategory of the Fukaya category, this makes a plausible case for new deformations of the entire Fukaya category. For sure, such deformations of the Fukaya category have not been discovered yet. Future readers may therefore find joy in experimenting with other deformations of $ \Gtl Q $ and for example obtain deformed “Fukaya categories with orbifold points”.

\subsubsection{The homological splitting}
Whenever one wants to compute a minimal model of an $ A_∞ $-category explicitly, one needs a homological splitting $ R ⊕ I ⊕ H $ of the $ A_∞ $-category. A homological splitting is by no means unique, and different homological splittings result in different but quasi-equivalent minimal models. The present paper deploys a specific choice of homological splitting for the category $ \ZigzagCat ⊂ \Tw\Gtl Q $ of zigzag paths. In this section, we explain why this homological splitting should be established as the standard splitting for $ \ZigzagCat $. We also comment on how to extend it to curves other than zigzag paths.

The homological splitting we choose in this paper is very well suited for the category $ \ZigzagCat $. This splitting is chosen under the expectation that $ \H\ZigzagCat $ is a full subcategory of the Fukaya category. The reader finds the definition of the homological splitting in \autoref{sec:splitting-splitting}, and comments on why this particular splitting is suited in \autoref{sec:whydoes-choice}. In fact, the homological splitting is both the right splitting to simplify the calculations, and the right one to prove $ \H\DefZigzagCat $ equal to the relative Fukaya category without further hassle with gauge equivalence. There is no doubt that the homological splitting is the best one for $ \ZigzagCat $.

It is clear that minor modifications to the homological splitting are possible. Most obviously, in the splitting we present the author is free to choose where to put identity and co-identity morphisms of each zigzag path (the choices of $ a_0 $ and $ α_0 $). A few actual changes are also possible: For instance, regard a transversal odd crossing between two zigzag paths. In the words of \autoref{sec:splitting-splitting}, this corresponds to a B situation. Our choice of cohomology basis elements consists of the angle sum $ (-1)^{\#α_3 + 1} α_3 + (-1)^{\#α_4} α_4 $. Choosing $ (-1)^{\#α_1 + 1} α_1 + (-1)^{\#α_2} α_2 $ instead is however possible just as well.

While basis morphisms in the Fukaya category have a unique “location” in the surface, cohomology basis morphisms of $ \Tw\Gtl Q $ can only imitate this behavior. Basis morphism in the Fukaya category lie on arcs of $ Q $, while odd cohomology basis morphisms of $ \Tw\Gtl Q $ can only lie around punctures of $ Q $. The quality of this imitation determines whether the minimal model computation yields a result in the desired shape or not.

\begin{figure}
\centering
\begin{subfigure}{0.4\linewidth}
\centering
\begin{tikzpicture}
\path[draw] (0, -0.1) -- ++(right:1.2) coordinate (B) coordinate[midway] (1) coordinate[pos=0.6] (2-end) -- ++(60:0.95) coordinate[midway] (2) -- ++(330:1) coordinate[midway] (3);
\path[draw] (B) ++(left:0.2) ++(up:0.1) ++(300:1) -- ++(120:1) coordinate[midway] (4) -- ++(60:1) coordinate[pos=0.3] (2-start) coordinate[pos=0.7] (1-end) coordinate[midway] (5) -- ++(120:1) coordinate (A) coordinate[midway] (6) -- ++(60:1) coordinate[midway] (7);
\path[draw] (A) ++(left:1.2) -- ++(right:1) coordinate[midway] (8) -- ++(300:0.95) coordinate[pos=0.7] (1-start) coordinate[midway] (9) -- ++(30:1) coordinate[midway] (10);
\path[draw, bend right=40, ->] (1-start) to node[midway, left] {$ α_3 $} (1-end);
\path[draw, bend right, ->] (2-start) to node[midway, left] {$ α_3 $} (2-end);
\path (2) -- (5) coordinate[midway] (m1);
\path (6) -- (9) coordinate[midway] (m2);
\path[draw, thick, rounded corners] ($ (1)!-0.5!(2) $) -- (1) -- (m1) -- (3) -- ($ (3)!-0.5!(m1) $);
\path[draw, thick, rounded corners] ($ (4)!-0.5!(5) $) -- (4) -- (m1) -- (m2) -- (7) -- ($ (7)!-0.5!(m2) $);
\path[draw, thick, rounded corners] ($ (8)!-0.5!(9) $) -- (8) -- (m2) -- (10) -- ($ (10)!-0.5!(m2) $);
\path[fill] (m1) circle[radius=0.05];
\path[fill] (m2) circle[radius=0.05];
\end{tikzpicture}
\caption{The simplest tree detects the disk.}
\label{fig:reuse-splitting-rapid-OK}
\end{subfigure}
\hspace{0.5cm}
\begin{subfigure}{0.4\linewidth}
\centering
\begin{tikzpicture}
\path[draw] (0, -0.1) -- ++(right:1.2) coordinate (B) coordinate[midway] (1) coordinate[pos=0.6] (2-end) -- ++(60:0.95) coordinate[midway] (2) -- ++(330:1) coordinate[midway] (3) coordinate[pos=0.3] (3-end);
\path[draw] (B) ++(left:0.2) ++(up:0.1) ++(300:1) -- ++(120:1) coordinate[midway] (4) -- ++(60:1) coordinate[pos=0.3] (2-start) coordinate[pos=0.7] (1-end) coordinate[midway] (5) -- ++(120:1) coordinate (A) coordinate[midway] (6) -- ++(60:1) coordinate[midway] (7);
\path[draw] (A) ++(left:1.2) -- ++(right:1) coordinate[midway] (8) -- ++(300:0.95) coordinate[pos=0.7] (1-start) coordinate[midway] (9) -- ++(30:1) coordinate[midway] (10);
\path[draw, bend right, ->] (1-start) to node[midway, left] {$ α_3 $} (1-end);
\path[draw, bend right, ->] (1-end) to node[midway, below] {$ α_1 $} (3-end);
\path (2) -- (5) coordinate[midway] (m1);
\path (6) -- (9) coordinate[midway] (m2);
\path[draw, thick, rounded corners] ($ (1)!-0.5!(2) $) -- (1) -- (m1) -- (3) -- ($ (3)!-0.5!(m1) $);
\path[draw, thick, rounded corners] ($ (4)!-0.5!(5) $) -- (4) -- (m1) -- (m2) -- (7) -- ($ (7)!-0.5!(m2) $);
\path[draw, thick, rounded corners] ($ (8)!-0.5!(9) $) -- (8) -- (m2) -- (10) -- ($ (10)!-0.5!(m2) $);
\path[fill] (m1) circle[radius=0.05];
\path[fill] (m2) circle[radius=0.05];
\end{tikzpicture}
\caption{The simplest tree detects no disk.}
\label{fig:reuse-splitting-rapid-broken}
\end{subfigure}
\caption{Detecting disks with Kadeishvili trees}
\end{figure}
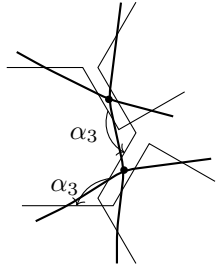
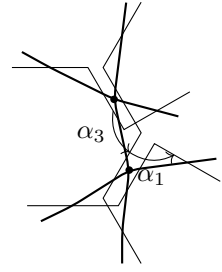

In our choice of homological splitting, we consistently choose $ α_3 + α_4 $ for every single B situation. This has the advantage that many immersed disks with intersection points $ h_1, …, h_N $ between zigzag paths can be imitated by the simplest possible Kadeishvili tree $ πμ(β_N, …, β_1) $, where $ β_1, …, β_N $ denote the corresponding B situation cohomology morphisms of type $ α_3 + α_4 $. More specifically, the simplest Kadeishvili tree is capable of capturing immersed disks where two situation B crossings follow each other within one arc distance. An illustration is shown in \autoref{fig:reuse-splitting-rapid-OK}.

If we were to choose $ α_3 + α_4 $ for some B situations and $ α_1 + α_2 $ for other B situations, the simplest Kadeishvili tree would not recognize disks where B situations follow each other rapidly. An example is shown in \autoref{fig:reuse-splitting-rapid-broken}. That figure depicts three curves and a piece of an immersed disk between them. For the upper B situation the morphism $ h_1 = α_3 + α_4 $ was chosen as cohomology basis representative, while for the lower B situation the morphism $ h_2 = α_1 + α_2 $ was chosen. It is impossible to form a disk $ μ^{\geq 3} (…, h_2, h_1, …) $ in $ \Tw\Gtl_q Q $. We conclude that a random choice of cohomology basis morphisms makes the minimal model calculation much less tractable.

Let us put the versatility of our homological splitting in the context of result components. Any choice of homological splitting provides an automatic notion of result components. To exploit result components for a minimal model calculation, one however needs to enumerate all possible result components by some target set, see \autoref{sec:reuse-resultcomp}. This enumeration by a target set is not automatic and depends on situational insight.

In our case of computing $ \H\ZigzagCat $, the notion of result components only needs to be tweaked minimally in order to map bijectively to the target set of immersed disks. Upon choice of a very different homological splitting for $ \ZigzagCat $, a notion of result components is still automatic, but the collection of result components does not biject to immersed disks anymore. Instead, it will biject to a complicated set of disk-like objects that requires far more detailed analysis. In other words, our choice of homological splitting has the advantage that its result components have a very simple target set.

It seems possible to find a homological splitting also objects in $ \Tw\Gtl Q $ which are not zigzag paths. The idea is still to sort elementary morphisms into different kinds of situations and to define the spaces $ H $ and $ R $ explicitly. The difficulty is however that general string and band objects have no limit with regards to the kind of angles they involve between two arcs. This means it is hard to find explicit cohomology representatives and to check that it concerns a homological splitting.

\subsubsection{The notion of result components}
\label{sec:reuse-resultcomp}
Result components are a technical tool serving as the main carrier of information in this paper, see \autoref{sec:resultcomp}. The idea is easy: A term like $ (3x + 5y) (2x + 3y) $ has the result components $ 6x^2 $, $ 9xy $, $ 10xy $, $ 15y^2 $. In other words, there are four distinct result components, even though the result can be abbreviated to only three terms. Result components provide maximum insight into part of the result instead of the whole, and how it is obtained instead of what is obtained. In this section, we argue that result components provide a means to analyze complicated Kadeishvili trees.

Regard a Kadeishvili tree $ T $ with $ N $ leaves and let $ h_1, …, h_N $ be inputs for the tree that lie in cohomology. Then for every node $ N ∈ T $, there is attached a set of result components. The set of result components is determined from choice of result components of all children of $ N $. In other words, result components are an inductive notion.

Let us paraphrase how we use result components in the present paper. We map a set of result components to “open” smooth immersed disks, which are called subdisks in \autoref{sec:subdisk}. This map is defined inductively: Given a result component at a node $ N $, it is analyzed how the result component was obtained from result components of the node's children. By induction hypothesis, every of the node's children already has a subdisk assigned. The subdisk associated with the result component at $ N $ is then obtained by gluing together the subdisks of the children in a way specific to the type of result component. This provides an inductively defined map from the set of result components of a Kadeishvili tree to the set of immersed disks.

For some categories, result components are better suited than for others. If the reader suspects that its minimal model has limited higher products, result components will not provide any use since most Kadeishvili trees result in zero anyway. If he however suspects that the minimal model calculation will result in a certain infinite “hierarchy” of higher products, then result components capture the higher products effectively.

For the reader who wishes to calculate the minimal model of some $ A_∞ $-category via result components, we suggest the following roadmap:
\begin{enumerate}
\item Find a homological splitting of the category. Typically, cohomology representatives must be found at the beginning and the rest space $ R $ can be accumulated on the go. The next step is to perform a few test calculations of products $ μ^k (h_k, …, h_1) $, where $ h_i $ are cohomology basis elements. The typical node in a Kadeishvili tree has output covering one or multiple basis elements of $ R $. This is the time to start accumulating basis elements into $ R $. The reader would then try to evaluate some products $ μ^k (…) $ where the inputs are mixed from both cohomology and $ R $. Which inputs from $ R $ multiply to a nonzero product and how does the product depend on its inputs?
\item Construct a notion of result components. The exact way to do this depends on the situation. In the easiest case, a result component would simply be defined inductively as an output term of the evaluation of $ hμ $ at each node, or $ πμ $ at the root. For other calculations like ours, it makes sense to distinguish or identify some output terms of $ hμ $ or $ πμ $ at every node (for example $ α_3 $ and the corresponding $ α_4 $ output are always collected as a combined result component $ α_3 + α_4 $).
\item Analyze how result components are derived. It is by no means necessary to classify all result components directly. Rather, it is important to classify result components into different types and understand which result components of which type can be derived from result components of which other types.
\item Determine a “target structure” or “target set”. The idea is to match result components with instances of some kind of better understood structure. For example, we have identified immersed disks as the correct target structure for result components $ \DefZigzagCat $. Upon commencing this step, a vague idea of what the target structure or target set will be may help. In either case, the target structure becomes clearer as the application of result components proceeds.
\item Matching result components with target objects. This step is hardest. But when performed successfully, this step ensures that the correspondence between result components and target objects can be written down explicitly and in a recursive manner.
\item Perform an inverse construction. The idea is to classify which instances of the target structure have been obtained via the identification. By constructing an explicit inverse mapping, it becomes clear which target objects have been reached and which not.
\end{enumerate}
The hard part always lies in identifying the correct target structure and the right identification of result components with target objects. Depending on what is expected from the particular minimal model, it might be possible to interpret the structure of a given Kadeishvili tree in a geometric way, so as to guess what the correct target object is.

The field of homological algebra requires us to perform a lot of minimal model calculations. Many minimal model calculations can be simplified vastly by choosing a clever homological splitting. However, minimal models are often not computed in their entirety. An example is Bocklandt's partial computation of $ \H\mf(\Jac Q, ℓ) $ in \cite{Bocklandt}, which is nevertheless sufficient to prove mirror symmetry for punctured surfaces. We are convinced that result components can facilitate the execution of complete minimal model calculations wherever a geometric outcome is expected.

\section{Notation}
\label{sec:notation}
The following is a list of heavily used notation specific to this paper:

\begin{center}
\begin{longtable}{lll}
Notation & Meaning & Reference \\\hline
\endhead
$ \cat C $ & $ A_∞ $-category & \autoref{def:2Binfty-ainfty-def} \\
$ \cat C_q $ & $ A_∞ $-deformation of $ \cat C $ & \autoref{def:2Bprelim-ainfty-defo} \\
$ μ_q $ & $ μ_{\cat C_q} $, more specifically $ μ_{\Gtl_q Q} $ or $ μ_{\Add \Gtl_q Q} $ & \autoref{def:2Bprelim-ainfty-defo} \\
$ B $ & deformation base & \autoref{def:2Bainfty-defo-base} \\
$ \mathfrak{m} $ & maximal ideal of $ B $ & \autoref{def:2Bainfty-defo-base} \\
$ B \htensor V $ & completed tensor product with vector space $ V $ & \autoref{sec:2Bainfty-defo} \\
$ \mathfrak{m} V $ & shorthand for $ \mathfrak{m} \htensor V ⊂ B \htensor V $ & \autoref{sec:2Bainfty-defo} \\
$ \Tw\cat C_q $ & twisted completion of $ \cat C_q $ & \autoref{def:2B-ainfty-defo-twisted} \\
$ \Tw'\cat C_q $ & liberal twisted completion of $ \cat C_q $ & \autoref{rem:2B-ainfty-defo-liberaltwisted} \\
$ \H\cat C_q $ & minimal model of $ \cat C_q $ & \autoref{def:2Bkadeishvili-minmodel} \\
$ \mathcal{T} $ & set of Kadeishvili tree shapes & \autoref{def:2Bkadeishvili-classical-treeshape} \\
$ N_T $ & number of internal nodes of a tree & \autoref{def:2Bkadeishvili-classical-treeshape} \\
$ φ $ & bijection $ H_q → B \htensor H $ & \autoref{def:2Bkadeishvili-deformed-counterpart} \\
$ F: \cat C → \cat D $ & $ A_∞ $-functor & \autoref{def:2Bprelim-ainfty-functor} \\
$ F_q: \cat C_q → \cat D_q $ & functor of deformed $ A_∞ $-categories & \autoref{def:2Bprelim-defo-functor} \\
$ (S, M) $ & punctured surface & \autoref{def:prelim-punctured-def} \\
$ \cA $ & arc system & \autoref{def:prelim-arcsys-def} \\
$ a $ & arc in $ \cA $ & \autoref{def:prelim-arcsys-def} \\
$ h(a), t(a) $ & puncture at head/tail of arc $ a $ & \autoref{def:prelim-arcsys-def} \\
$ α $ & angle in $ \cA $ & \autoref{def:prelim-gtl-gtl-def} \\
$ h(α), t(α) $ & arc at head/tail of angle $ α $ & \autoref{def:prelim-arcsys-def} \\
$ Q $ & dimer, typically geometrically consistent & \autoref{def:prelim-gtl-dimers-def} \\
$ Q_M $ & standard sphere dimer & \autoref{sec:sphere-zigzagcat} \\
$ \id_a $ & arc identity & \autoref{sec:prelim-gtl-gtl} \\
$ L $ & zigzag path & \autoref{def:prelim-gtl-zigzag-def} \\
$ \ZigzagCat $ & zigzag category & \autoref{def:splitting-zigzagcat-def} \\
$ \DefZigzagCat $ & deformed zigzag category & \autoref{def:deformed-zigzagcat-def} \\
$ a_0 $ & identity location on zigzag path & \autoref{conv:alpha0-direction} \\
$ α_0 $ & co-identity location on zigzag path & \autoref{conv:alpha0-direction} \\
$ P_k $ & standard $ k $-gon & \autoref{sec:prelim-gtl-infty} \\
$ ε $ & elementary morphism $ ε: L_1 → L_2 $ & \autoref{sec:prelim-terminology} \\
$ T $ & Kadeishvili tree shape & \autoref{def:2Bkadeishvili-classical-treeshape} \\
$ T $ & tail of a morphism $ ε: L_1 → L_2 $ & \autoref{def:deformed-tails-tail} \\
$ \PiTr $ & class of result components of π-trees & \autoref{def:subdisk-shapeless-def} \\
$ \SLd $ & class of shapeless disks & \autoref{def:subdisk-shapeless-def} \\
$ \CRr $ & class of CR result components & \autoref{def:subdisk-types-rescomp} \\
$ \IDr $ & class of ID result components & \autoref{def:subdisk-types-rescomp} \\
$ \DSr $ & class of DS result components & \autoref{def:subdisk-types-rescomp} \\
$ \DWr $ & class of DW result components & \autoref{def:subdisk-types-rescomp} \\
$ \CRd $ & class of CR disks & \autoref{def:subdisk-types-CRdisk} \\
$ \IDd $ & class of ID disks & \autoref{def:subdisk-types-IDdisk} \\
$ \DSd $ & class of DS disks & \autoref{def:subdisk-types-DSdisk} \\
$ \DWd $ & class of DW disks & \autoref{def:subdisk-types-DWdisk} \\
$ \Subdisk $ & subdisk mapping $ \Subdisk: \PiTr → \SLd $ & \autoref{th:subdisk-subdisk-Dconstruction} \\
$ \disktarget(D) $ & target/output morphism of disk $ D $ & \autoref{def:subdisk-shapeless-def} \\
$ \Abouzaid(D) $ & Abouzaid sign of disk $ D $ & \autoref{def:subdisk-minmodel-Abouzaid} \\
$ \punctures(D) $ & product of punctures covered by $ D $ & \autoref{def:subdisk-minmodel-Abouzaid}
\end{longtable}
\end{center}

\printbibliography

\end{document}